\documentclass[openright]{cambridge7A}

\usepackage[scaled=0.9]{couriers}

\usepackage{amsthm}
\theoremstyle{plain}
  \newtheorem{theorem}{Theorem}[section]
  \newtheorem{lemma}[theorem]{Lemma}
  \newtheorem{proposition}[theorem]{Proposition}
  \newtheorem{corollary}[theorem]{Corollary}
  
  \newtheorem{warning}[theorem]{Warning}
  \newtheorem{notat}[theorem]{Notation}

  \newtheorem*{theorem*}{Theorem}
  \newtheorem*{lemma*}{Lemma}
  \newtheorem*{proposition*}{Proposition}
  \newtheorem*{corollary*}{Corollary}
  \newtheorem*{conjecture*}{Conjecture}
  \newtheorem*{notat*}{Notation}

\theoremstyle{definition}
  \newtheorem{definition}[theorem]{Definition}
  \newtheorem{example}[theorem]{Example}
  
  \newtheorem{remark}[theorem]{Remark}
  \newtheorem{notation}[theorem]{Notation}
  \newtheorem{axiom}[theorem]{Axiom}

  \newtheorem*{definition*}{Definition}
  \newtheorem*{example*}{Example}
  \newtheorem*{prob*}{Problem}
  \newtheorem*{remark*}{Remark}
  \newtheorem*{notation*}{Notation}

\def\makeRRlabeldot#1{\hss\llap{#1}}
\renewcommand\theenumi{{\rm (\roman{enumi})}}
\renewcommand\theenumii{{\rm (\alph{enumii})}}
\renewcommand\theenumiii{{\rm (\arabic{enumiii})}}
\renewcommand\theenumiv{{\rm (\Alph{enumiv})}}

\usepackage[usenames,dvipsnames]{xcolor}

\usepackage[export]{adjustbox}
\usepackage{
    amssymb
  , amsfonts
  , amsmath
  , amsthm
  , booktabs
  , caption
  , datetime
  , epigraph
  , etex
  , etoolbox
  , graphicx
  , hyphenat
  , lettrine
  , lipsum
  , lmodern
  , manfnt
  , marginnote
  , mathtools
  , microtype
  , multirow
  , newclude
  , rotating
  , stmaryrd
  , tcolorbox
  , tikz
  , tikz-cd
  , url
  , xparse
  , xspace
  }
\usepackage{esint}
\usepackage{inconsolata}

\usepackage[normalem]{ulem}

\usetikzlibrary{ 
  matrix
  , calc
  , babel
  , fit
  , shapes
  , positioning
  , arrows
  }

\usepackage[ 
  spanish
  , greek
  , english
  ]{babel}

\usepackage[utf8]{inputenc}
\usepackage[T1]{fontenc}

\usepackage[all,2cell]{xy}\UseAllTwocells
\usepackage{cjhebrew}

\newcommand{\eend}{\underline{\mathrm{end}}}
\newcommand{\coend}{\underline{\mathrm{coend}}}

\def\mho{\rotatebox[origin=c]{180}{$\Omega$}}

\def\proP{\mathfrak{p}}
\def\proQ{\mathfrak{q}}
\def\proH{\mathfrak{r}}
\def\proL{\mathfrak{l}}

\def\proG{\mathfrak{g}}

\def\psh#1{[{#1}^\opp,\Set]}
\newcommand{\celtag}[2][dr]{\ar[#1,white, "#2"{black,description}]}

\def\Fin{\BF{Fin}}
\def\Map{\BF{Map}}

\newcommand{\po}[1][dr]{\save*!/#1-3pc/#1:(1,-1)@^{|-}\restore}
\newcommand{\pb}[1][dr]{\save*!/#1-2pc/#1:(-1,1)@^{|-}\restore}

\newcommand{\bigast}{\mathop{\scalebox{1.5}{\raisebox{-0.2ex}{$\ast$}}}}

\def\tr{\mathrm{tr}}

\providecommand{\abbrv}[1]{#1.\@\xspace}
  \providecommand{\ie}	 {\abbrv{i.e}}

\usepackage{CJKutf8}

\newcommand{\yon}{\text{\begin{CJK}{UTF8}{min}よ\end{CJK}}}
\newcommand{\coyon}{\rotatebox[origin=c]{180}{$\yon$}\!}
\newcommand{\wk}{\textsc{wk}}
\newcommand{\cof}{\textsc{cof}}
\newcommand{\fib}{\textsc{fib}}

\def\itsnonsense{\upeyes}

\newcommand{\bDelta}{\boldsymbol{\Delta}}
\newcommand{\xto}[2][]{\xrightarrow[#1]{#2}}
\newcommand{\xot}[2][]{\xleftarrow[#1]{#2}}

\newcommand{\bsmat}{\left[\begin{smallmatrix}}
\newcommand{\esmat}{\end{smallmatrix}\right]}

\newenvironment{xsmallmatrix}[1]
  {\renewcommand\thickspace{\kern#1}\smallmatrix}
  {\endsmallmatrix}

\renewcommand{\phi}{\varphi}
\newcommand{\var}[3][]{
  \left[\begin{smallmatrix} #2 \\
  #1\downarrow \\ #3
  \end{smallmatrix}\right]}
\newcommand{\cvar}[3]{
  \begin{xsmallmatrix}{0em}
  & #1 \\ #2 & \downarrow \\ & #3
  \end{xsmallmatrix}}

\newcommand{\laxto}{\dashrightarrow}
\newcommand{\ccirc}{\diamond}

\DeclareMathOperator{\id}{id}
\DeclareMathOperator{\Lan}{Lan}
\DeclareMathOperator{\lan}{lan}
\DeclareMathOperator{\Ran}{Ran}
\DeclareMathOperator{\ran}{ran}
\DeclareMathOperator{\hoLan}{hoLan}
\DeclareMathOperator{\hoRan}{hoRan}
\DeclareMathOperator{\Rift}{Rift}
\DeclareMathOperator{\Lift}{Lift}
\DeclareMathOperator{\rift}{rift}
\DeclareMathOperator{\leeft}{lift}

\DeclareMathOperator{\eq}{eq}

\newcommand{\opp}{\mathrm{op}}
\newcommand{\co}{\mathrm{co}}
\newcommand{\coop}{\mathrm{coop}}
\newcommand{\oo}{\text{oo}}
\newcommand{\Ex}{\text{Ex}}

\newcommand{\cate}[1]{\mathcal{#1}}
  \newcommand{\A}   {\cate{A}}
  \newcommand{\B}   {\cate{B}}
  \newcommand{\C}   {\cate{C}}
  \newcommand{\D}   {\cate{D}}
  \newcommand{\E}   {\cate{E}}
  \newcommand{\F}   {\mathfrak{F}}
  \newcommand{\M}   {\cate{M}}
  \newcommand{\J}   {\cate{J}}
  \renewcommand{\I} {\cate{I}}
  
  \newcommand{\K}   {\cate{K}}
  \renewcommand{\P} {\cate{P}}
  \newcommand{\LL}  {\cate{L}}
  \newcommand{\V}   {\cate{V}}
  \newcommand{\cX}  {\cate{X}}
  \newcommand{\cY}  {\cate{Y}}
  \newcommand{\cZ}  {\cate{Z}}
  \newcommand{\cT}  {\cate{T}}
  
  \newcommand{\cN}  {\cate{N}}
  \newcommand{\cW}  {\cate{W}}
  \newcommand{\catS}{\cate{S}}

  \newcommand{\nw}{\protect{\xymatrix@C=1.5cm{\A \ruppertwocell^F{\alpha}
\rlowertwocell_H{\beta} & \B \ruppertwocell^K{\gamma}
\rlowertwocell_M{\delta} & \C }}}
\def\sw{\protect{\xymatrix@C=1.5cm{\A\rtwocell^F_H{\quad\beta\circ\alpha} & \B\rtwocell^K_M{\quad\delta\circ\gamma} & \C}}}
\def\northeast{\protect{\xymatrix@C=2cm{
\A \ruppertwocell^{KF}{\quad\gamma\boxminus\alpha} 
\rlowertwocell_{MH}{\quad\delta\boxminus\beta} & \C 
}}}

  \newcommand{\tCat} {\mathsf{Cat}}
  \newcommand{\caat} {\mathsf{cat}}
  \newcommand{\sfK} {\mathsf{K}}
  \newcommand{\sfA} {\mathsf{A}}
  \newcommand{\adm}{\mathsf{Ads}}

  \newcommand{\BF}[1]{\mathrm{#1}}
  \newcommand{\Cat} {\BF{Cat}}
  \newcommand{\Ch}  {\BF{Ch}}
  
  \newcommand{\Ab}  {\BF{Ab}}
  \newcommand{\Mod} {\BF{Mod}}
  
  \newcommand{\Sets}{\BF{Set}}
  \newcommand{\Set} {\BF{Set}}
  \newcommand{\sSet}{\BF{sSet}}
  \newcommand{\Spc} {\BF{Spc}}

  \newcommand{\N} {\mathbb{N}}
  \newcommand{\Z} {\mathbb{Z}}
  
  \newcommand{\bR}{\mathbb{R}}
  \newcommand{\bC}{\mathbf{C}}
  \newcommand{\sD}{\mathbb{D}}

  \newcommand{\fko}{\mathfrak{o}}
  \newcommand{\fki}{\mathfrak{i}}
  \newcommand{\fkj}{\mathfrak{j}}
  \newcommand{\fkP}{\mathfrak{P}}
  \newcommand{\fkQ}{\mathfrak{Q}}

  \newcommand{\Y}     {\boldsymbol{Y}}
  \newcommand{\W}     {\boldsymbol{W}}
  \newcommand{\VCat}  {\V\text{-}\Cat}
  \newcommand{\twoCat}{2\text{-}\Cat}

\def\bemo{\flat}
\def\diesis{\sharp}
\DeclareMathOperator{\colim}{colim}
\DeclareMathOperator*{\llim}{llim}
\DeclareMathOperator*{\lcolim}{lcolim}
\DeclareMathOperator*{\hocolim}{hcolim}
\DeclareMathOperator*{\holim}{hlim}
\DeclareMathOperator{\lLan}{lLan}
\let\varinjlim\colim
\let\varprojlim\lim
\newcommand{\pto}{\leadsto}
\newcommand{\otp}{\mathrel{\reflectbox{$\pto$}}}
\newcommand{\wlim}[1]{{\textstyle\lim^{#1}}}
\newcommand{\wcolim}[1]{{\textstyle\colim^{#1}}}

\newcommand{\elts}[2]{
  \mathchoice{#1 \rotatebox[origin=c]{15}{$\int$} #2}
  {#1 \rotatebox[origin=c]{15}{$\int$} #2}
  {#1 \rotatebox[origin=c]{15}{\scriptsize$\int$} #2}
  {#1 \rotatebox[origin=c]{15}{\tiny$\int$} #2}
}
\newcommand{\tw}{\textsc{tw}}

\newcommand{\twoint}{\sqint}
\newcommand{\infint}{\oint}

\def\To{\Rightarrow}

\newcommand{\rotArrow}[1]{\rotatebox[origin=c]{#1}{$\Rightarrow$}}
  \providecommand{\Nearrow}{\rotArrow{45}}
  \providecommand{\Nwarrow}{\rotArrow{135}}
  \providecommand{\Searrow}{\rotArrow{-45}}
  \providecommand{\Swarrow}{\rotArrow{225}}
  \providecommand{\Sarrow} {\rotArrow{-90}}

\newadjustboxenv*{noverfullmath}{max width=\textwidth, vspace={\the\abovedisplayskip} {\the\belowdisplayskip},innercode={$\displaystyle}{$}}

\let\to\rightarrow
\def\diag{\mathsf{d}}
\def\tot{\mathsf{t}}
\def\Dist{\mathsf{Prof}}
\def\Kl{\BF{Kl}}
\def\Alg{\BF{Alg}}
\def\Ab{\BF{Ab}}

\def\sfSd{\mathsf{Sd}}

\def\sfSpli{\mathsf{Spli}}
\def\sfA{\mathsf{A}}
\def\sfDer{\mathsf{Der}}
\def\clX{\mathcal{X}}
\def\clY{\mathcal{Y}}
\def\clO{\mathcal{O}}
\def\DFib{\BF{DFib}}

\DeclareMathOperator{\coker}{coker}

\newenvironment{abstract}{%
 \begin{center}
 \begin{minipage}{.9\textwidth}
  \setlength{\parindent}{1em}
 \textsc{Summary.}
 \small\linespread{1.1}
}
{
 \end{minipage}
 \end{center}
 \normalsize
 \linespread{1}
}

\def\ceil#1{\lceil #1 \rceil}

\def\ULbullet{{}^\bullet\kern-.1em}
\newcommand{\isbO}{\mathsf{O}}
\def\Spec{\mathsf{Spec}}

\newcommand{\doboxminus}[1]{\mathbin{
  \begin{tikzpicture}[scale=#1,line width=.5pt]
    \draw[rounded corners=.25pt]
      (0,0) rectangle (1,1);
    \draw[shorten <=.75pt, shorten >=.75pt, cap=round]
      (0,.5) -- (1,.5);
    \end{tikzpicture}
}}

\renewcommand{\boxminus}{
\mathchoice{\doboxminus{.165}}
  {\doboxminus{.165}}
  {\,\mathbin{\doboxminus{.1325}}\,}
  {\mathbin{\doboxminus{.1}}}
}

\usepackage{turnstile}
\newcommand{\adjunct}[2]%
  {\nsststile{#2}{#1}}

\DeclareDocumentEnvironment{enumtag}{m}
{\begin{enumerate}[ label=\textsc{#1}\oldstylenums{\arabic*})
                  , ref=\textsc{#1}\oldstylenums{\arabic*}%
	                ]}%
{\end{enumerate}}

\newcommand{\Cn}{\BF{Cn}}
\newcommand{\Seq}{\BF{Seq}}
\newcommand{\Ccn}{\BF{Ccn}}
\newcommand{\wed}{\BF{Wd}}
\newcommand{\cwed}{\BF{Cwd}}
\newcommand{\Tamb}{\BF{Tamb}}

\newcommand{\din}
  {\mathrel{%
    \begin{tikzpicture}[baseline=-0.58ex]
    \begin{scope}[line cap=round, thin]
    \draw[yshift=-.25mm] (0,0) -- (-.275,0);
    \draw[yshift= .25mm] (0,0) -- (-.275,0);
    \end{scope}
    \filldraw[ rounded corners=.25pt
             , thin
             , rotate around={45:(0,0)}
             , fill=white
             ] (-.45mm,-.45mm) rectangle (.45mm,.45mm);
    \end{tikzpicture}%
  }}

\newcommand{\deduction}[4]{%
  \begin{array}{c}
    #1 \To #2 \\ \hline
    #3 \To #4
    \end{array}%
  }

\def\bsP{\boldsymbol{P}}
\def\bsT{\boldsymbol{T}}
\def\fh{\mathfrak{h}}

\newcommand{\uno}{\BF{1}}

\newcommand{\comma}[1]
  {\left[\begin{smallmatrix}
  	#1
  \end{smallmatrix}
  \right]}

\makeatletter
  \def\@cite#1#2{[\textsf{#1}\if@tempswa , #2\fi]}
  \def\@biblabel#1{[\textsf{#1}]}
\makeatother

\newlength{\seplen}
\newlength{\sepwid}
\def\firstblank{\,\rule{\seplen}{\sepwid}\,}
\def\secondblank{\firstblank\llap{\raisebox{2pt}{\firstblank}}}

\renewcommand{\textbf}[1]{\text{\fontseries{b}\selectfont{\upshape #1}}}

\usepackage{enumitem}

\usepackage[titletoc]{appendix}

\allowdisplaybreaks

\newcommand{\standard}{
\draw[ultra thin] (0.25,0) -- (1.75,0) -- (1.75,1) -- (0.25,1) -- cycle;
\draw (.5,0) -- (.5,1);
\draw (1,0) .. controls (1,.5) and (1.5,.5) .. (1.5,0);
\draw (1,1) .. controls  (1,.5) and (1.5,.5) .. (1.5,1);
}

\newcommand{\standardbis}[3]{
\draw[ultra thin] (0.25,0) -- (1.75,0) -- (1.75,1) -- (0.25,1) -- cycle;
\draw[#1] (.5,0) -- (.5,1);
\draw[#2] (1,0) .. controls (1,.5) and (1.5,.5) .. (1.5,0);
\draw[#3] (1,1) .. controls  (1,.5) and (1.5,.5) .. (1.5,1);
}

\newcommand{\fst}[1]{\draw[ultra thin] (0.25,0) -- (1.75,0) -- (1.75,1) -- (0.25,1) -- cycle;
\draw (.5,0) -- (.5,1);
\draw[xshift=.5cm] (.5,0)  -- (.5,1);
\draw[xshift=1cm] (.5,0)-- (.5,1);
\filldraw[lightgray!70] (.5,.5) circle (6pt) node[black] {\footnotesize $ #1 $};}
\newcommand{\snd}[1]{\draw[ultra thin] (0.25,0) -- (1.75,0) -- (1.75,1) -- (0.25,1) -- cycle;
\draw (.5,0) -- (.5,1);
\draw[xshift=.5cm] (.5,0)  -- (.5,1);
\draw[xshift=1cm] (.5,0)-- (.5,1);
\filldraw[lightgray!70] (1,.5) circle (6pt) node[black] {\footnotesize $ #1 $};}
\newcommand{\trd}[1]{\draw[ultra thin] (0.25,0) -- (1.75,0) -- (1.75,1) -- (0.25,1) -- cycle;
\draw (.5,0) -- (.5,1);
\draw[xshift=.5cm] (.5,0)  -- (.5,1);
\draw[xshift=1cm] (.5,0)-- (.5,1);
\filldraw[lightgray!70] (1.5,.5) circle (6pt) node[black] {\footnotesize $ #1 $};}

\newcommand{\upeyes}{%
  \begin{tikzpicture}
  \draw circle (3pt);
  \draw[fill,yshift=1.5pt] circle (1.5pt);
  \begin{scope}[xshift=7pt]
  \draw circle (3pt);
  \draw[fill,yshift=1.5pt] circle (1.5pt);
  \end{scope}
  \end{tikzpicture}\xspace%
}

\newcommand{\righteyes}{%
  \begin{tikzpicture}
  \draw circle (3pt);
  \draw[fill,xshift=1.5pt] circle (1.5pt);
  \begin{scope}[xshift=7pt]
  \draw circle (3pt);
  \draw[fill,xshift=1.5pt] circle (1.5pt);
  \end{scope}
  \end{tikzpicture}\xspace%
}

\newcommand{\awful}{%
   \begin{tikzpicture}
   \draw circle (3pt);
   \draw[fill] circle (1.5pt);
   \begin{scope}[xshift=7pt]
   \draw circle (3pt);
   \draw[fill] circle (1.5pt);
   \end{scope}
   \end{tikzpicture}\xspace%
}

\newcommand{\sierpo}{%
\'{S}}
\newcommand{\dis}{%
D}
\newcommand{\codis}{%
C}

\usepackage{hyperref}
\hypersetup{%
    pdftoolbar=true
  , pdfmenubar=true
  , pdffitwindow=true
  , pdftitle={This is the co/end}
  , pdfauthor={F. Loregian}
  , colorlinks=true
  , linkcolor=black
  , citecolor=black}

\def\[{\begin{equation}}
\def\]{\end{equation}}

\usepackage{makeidx}
\makeindex
\begin{document}
\title%
  [the book formerly known as `This is the co/end']%
  {Coend calculus}

\author{Fosco Loregian}

\bookabstract{}
\bookkeywords{coend calculus}

\frontmatter

\maketitle
\tableofcontents

\mainmatter

\cleardoublepage
\thispagestyle{empty}
\begin{flushright}
	\thispagestyle{empty}
	\vspace*{\fill}
	A \textgreek{`Ek'atera} \\
	ambrosia pura \\
	vestale di fuoco.
\end{flushright}
\vspace*{4cm}
\cleardoublepage

\pagenumbering{roman}
\chapter*{Preface}
\epigraph{
	No podemos esperarnos que ningún aspecto de la realidad cambie si seguimos usando los medios [\dots\unkern] de un lenguaje que lleva el peso de toda la negatividad del pasado. El lenguaje está ahí, pero tenemos que limpiarlo, revisarlo y, sobre todo, debemos desconfiar de él.
}{J. Cort\'azar}
\paragraph{What is co\fshyp{}end `calculus'.}
Coend calculus rules the behaviour of suitable universal objects associated to functors of two variables $T : \C^\opp\times \C \to \D$.

The intuition behind the process of attaching a special invariant to such a functor $T$ can be motivated in many ways.

It is well\hyp{}known that a measurable scalar function $f :X \to \bR$ from a measurable space $(X, \Omega)$ can be integrated `against' a measure $\mu$ defined on $\Omega$ to yield a real number
\[\notag
	\int_X f(x)d\mu
\]
(for example, when $X$ is a smooth space, the measure can be legitimately thought to depend `contravariantly' on $x$, as $d\mu$ is a volume form living in the top\hyp{}degree exterior algebra of $X$). In a similar fashion, the evaluation map $V^\lor\otimes V \to k$ for a vector space $V$ is a pairing $\langle \zeta,v\rangle = \zeta(v)$ between a vector $v$ and a co\hyp{}vector $\zeta : V \to k$, that becomes the sum $\sum_i \zeta_i v_i$ once a basis for $V$, and its dual basis, is chosen and the vector $v$ has coordinate $(v_1,\dots,v_d)$, whereas $\zeta$ has coordinates $(\zeta_1,\dots,\zeta_d)$.

At the cost of pushing this analogy further than permitted, a functor $T : \C^\opp\times \C \to \D$ can be thought as a generalised form of evaluation of an object of $\C$ against another; the `quantity' $T(C,C')$ can then be `integrated' to yield two distinct objects having dual universal properties:
\begin{enumtag}{c}
	\item a \emph{coend}, resulting by the symmetrisation along the diagonal of $T$, \ie by modding out the coproduct $\coprod_{C\in \C}T(C,C)$ by the equivalence relation generated by the arrow functions $T(\firstblank,C') : \C^\opp(X,Y) \to \D(T(X,C'), T(Y,C'))$ and $T(C,\firstblank) : \C(X,Y) \to \D(T(C,X), T(C,Y))$;
	\item an \emph{end}, \ie an object $\int_C T(C,C)$ arising as an `object of invariants' of `fixed points' for the same action of $T$ on arrows; by dualisation, if a coend is a quotient of $\coprod_{C\in \C}T(C,C)$, an end is a subobject of the product $\prod_{C\in \C} T(C,C)$.
\end{enumtag}
This also suggests a fruitful analogy with modules over a ring: if a functor $T : \C^\opp\times \C \to \Set$ is a `bimodule', that lets $\C$ act once on the left and once on the right on the sets $T(C,C')$, the end $\int_C T(C,C)$ is the subspace of invariants for the action of $\C$, whereas the coend $\int^C T(C,C)$ is the space of orbits (or ``coinvariants'') of said action.

In fact, a rather common way to employ coends is the following: consider a functor $F : \C \to \Set$ (a `left module') and a functor $G : \C^\opp\to \Set$ (a `right module'), and tensor them together into a functor $(C,C')\mapsto GC \times FC'$; the symmetrisation of $F\times G$ yields a \emph{functor tensor product} of $F,G$ as the set
\[\notag
	F\boxtimes G := \int^C FC \times GC.
\]
Note that in this light, the analogy is meaningful: if $\C$ is a single\hyp{}object category (so a monoid or a group $G$), such a pair of modules constitutes a pair $(X,Y)$ of a left and a right $G$\hyp{}set, and their functor tensor product can be characterized as the product $X\times_G Y$ obtained as the quotient of $X\times G$ for the equivalence relation $(g.x,y)\sim (x, g.y)$, so that $X\times_G Y$ is the universal \emph{$G$\hyp{}bilinear} product of sets, in that the ``scalar'' $g\in G$ can pass left\hyp{}to\hyp{}right from $(g.x,y)$ to $(x,g.y)$ in the quotient. Of course, the terminology works better when $X,Y$ are vector spaces carrying a \emph{linear} representation of $G$.

Theorems involving ends and coends $\int^C T(C,C)$ and $\int_C T(C,C)$ can now be proved by means of the universal properties that define them; it is easily seen that given $T : \C^\opp\times \C \to \D$ there exist a category $\bar{\C}$ and a functor $\bar T : \bar{\C}\to \D$ such that $\int^C T(C,C)\cong \colim_{\bar{\C}} \bar T$ and $\int_C T(C,C)\cong \lim_{\bar{\C}} \bar T$. In the example above, the tensor product of a left $G$\hyp{}module and a right $G$\hyp{}module (here `module' means `$k$\hyp{}vector space') can be characterized as the coequalizer
\[
	\xymatrix{
		\bigoplus_{g\in G} X \otimes_k Y \ar@<4pt>[r]^-\alpha\ar@<-4pt>[r]_-\beta & X \otimes_k Y \ar[r] & X\otimes_G Y
	}
\] where $\alpha(g,(x,y)) = (g.x,y)$ and $\beta(g, (x,y)) = (x, g.y)$.

So, all co\fshyp{}ends can be characterised as co\fshyp{}limits; but they provide a richer set of computational rules than mere co\fshyp{}limits. Often, establishing that an object has a certain universal property is a difficult task, because a direct argument tangles the reader to use elements. A general tenet of modern category theory is that cleaner, more conceptual arguments shall be preferred against element\hyp{}wise proofs that are evil in spirit, if not in shape.

Co/end calculus provides such conceptualisation for many classical arguments of category theory: it is in fact possible to prove that two objects of a category, at least one of which is defined as a coend, are isomorphic by means of a chain of `deduction rules'.

These rules are described in the first half of the book, but here we glimpse at what they look like.\footnote{The somewhat far\hyp{}fetched conjecture that permeates all the book is that coend calculus provides an higher\hyp{}dimensional version of a deductive system, suited for category theory (see \cite{winskel} for some preliminary steps in this direction), having deduction rules similar to those of Gentzen's sequent calculus. We will never attempt to turn this enticing conjecture into a theorem, or even to make a precise claim; the interested reader is thus warned that their curiosity will not get satisfaction --not in the present book, at least. We record that the idea that coends categorify logical calculus comes from William Lawvere, and it was first proposed in \cite{LawvereFW:metsgl}.}

In order to make clear what this paragraph is about, let us consider the statement that \emph{right adjoints preserve limits}; it is certainly possible to prove it by hand. Nevertheless, using little more than the Yoneda lemma it is possible to prove that if $R : \D\to\C$ is right adjoint to $L : \C \to \D$, there is a natural isomorphism of hom\hyp{}sets $\C(C, R(\lim_\J D_J))\cong \C(C, \lim_\J RD_J)$ for every object $C$, and every diagram $D : \J \to \D$ by arguing as follows:
\begin{center}
	\begin{tabular}{c}
		$\C(C, R(\lim_\J D_J))$ \\[1mm]\midrule
		$\C(LC, \lim_\J D_J)$   \\[1mm]\midrule
		$\lim_\J \C(LC, D_J)$   \\[1mm]\midrule
		$\lim_\J \C(C, RD_J)$   \\[1mm]\midrule
		$\C(C, \lim_\J RD_J)$
	\end{tabular}
\end{center}
where each step of this ``deduction'' is motivated either by the fact that $L\dashv R$ are adjoint functors, or by the fact that all functors $\C(X,\firstblank)$ preserve limits. Once this is proved, Yoneda lemma (see \ref{lem:the-real-yoda} for the statement) entails that there is an isomorphism $R(\lim D_J)\cong \lim RD_J$.

A similar argument is a standard way to prove that a certain object, defined (say) as the left adjoint to a certain functor, must admit an `integral expansion' to which it is canonically isomorphic. For example, in the proof of what we called \emph{ninja Yoneda lemma} in \ref{ninjayo}, we carry on the following computation:
\begin{center}
	\begin{tabular}{c}
		$\Set\Big( \int^{C\in\C} KC\times \C(X,C),Y\Big)$            \\[2mm]\midrule
		$\int_{C\in\C} \Set\big( KC\times \C(X,C),Y \big) $          \\[2mm]\midrule
		$\int_{C\in\C}\Set(\C(X,C),\Set(KC,Y)) $                     \\[2mm]\midrule
		$[\C, \Set]\big(\C(X,\firstblank),\Set(K\firstblank,Y)\big)$ \\[2mm]\midrule
		$\Set(KX,Y)$
	\end{tabular}
\end{center}
where each step has to be interpreted as an application of a certain deduction rule that interchanges coends with ends, places them in and out of a hom functor, etc.

The reduction of proofs to a series of deduction steps embodies some sort of `logical calculus', whose introduction rules resemble formulas as
\[\notag
	\Cat(\C,\D)(F,G) \rightsquigarrow \int_C \D(FC,GC)
\]
where the object $\Cat(\C,\D)(F,G)$ is decomposed into an integral like $\int_C \D(FC,GC)$, and elimination rules look like
\[\notag
	\int^C FC\times \C(C,X) \rightsquigarrow FX
\]
where an integral is packaged into the object $FX$ (of course, the symmetric nature of the canonical isomorphism relation makes all elimination rules reversible into introductions, and vice versa). Altogether, this allows to derive the validity of a canonical isomorphism as a result of a chain of deductions, in a `categorified' fashion.

The reader shall of course not concentrate now on the meaning of these derivations at all; all notation will be duly introduced at the right time; when the statement of our \ref{ninjayo} will be introduced, the chain of deductions above will look almost tautological, and rightly so.

It is clear that done in this way, category theory acquires an alluring algorithmic nature, and becomes (if not easy, at least) \emph{easier to understand}.

Thus, the `calculus' arising from these theorems encodes many, if not all, elementary constructions in category theory (we shall see that it subsumes the theory of co\fshyp{}limits, it allows for a reformulation of the Yoneda lemma, it provides an explicit formula to compute pointwise Kan extensions, and it is a cornerstone of the `calculus of bimodules', encoding the compositional nature of \emph{profunctors}, the natural categorification of relational composition).

As it stands, co\fshyp{}end calculus describes pieces of abstract and universal algebra \cite{curiennone,gambo-joy}, algebraic topology \cite{markl2007operads,may1972geometry,getzler2009operads}, representation theory \cite{loday2012algebraic}, logic, computer science \cite{kmett}, as well as pure category theory. The present book wishes to explore in detail such a theory and its applications.

\smallskip
So far, the motivations for the \emph{topic} of this book. What about the motivation \emph{for the book itself}? It shall be noted that co\fshyp{}ends are not absent from the already existing literature on category theory: the topic is covered in \cite{McL}, a statutory reading for every categorephile, and Mac Lane himself used coends to characterise a construction in algebraic topology as a `tensor product' operation between functors in his \cite{mac1970milgram}; coends are mentioned (but not used as widely as they deserve) in Borceux's \emph{Handbook} (in its first two tomes \cite{Bor1,Bor2}); however, the topic lacks a treatment that it is at the same time systematic, easy to read, and monographic.

As a result, co\fshyp{}end calculus still lies just beyond the grasp of many people, and even of a few category theorists, because the literature that could teach its simple rules is a vast constellation of scattered papers, drawing from a large number of diverse disciplines.

This situation is all the more an issue because nowadays category theory produces fruitful contamination with applied sciences: in the opinion of the author, it is of the utmost importance to provide his growing community with a single reference that accounts for the simplicity and unitary nature of category theory through co\fshyp{}end calculus, thereby providing proof for the plethora of its different applications, and popularising this `secret weapon' of category theorists, making it available to novices and non\hyp{}mathematicians. The present endeavour is but an humble attempt to address this issue.
\paragraph{A brief history of co\fshyp{}ends.}
Like many other pieces of mathematics, co\fshyp{}end calculus was developed as a tool for homological algebra: the first definition of co\fshyp{}end was given in a paper studying the $\text{Ext}$ functors, and the father of co\fshyp{}end calculus is none other than Nobuo Yoneda. In his \cite{Yoneda} he singled out most of the definitions we will introduce along the first five chapters of this book.

Having read Yoneda's original paper to write the present introduction, we find no better way than to quote the original text, untouched, just occasionally adding a few details here and there to frame Yoneda's words in a modern perspective, (but also in order to adapt them to our choice of notation).

This strategy has multiple purposes: \cite{Yoneda} is a mathematical gem, an enticing prelude of all the theory developed in the subsequent decades, and a perfect prelude to the story this book tries to tell; even more so, in reporting Yoneda's words we believe we are also doing a service to the mathematical community, since the integral text of \cite{Yoneda} is somewhat difficult to find.

\smallskip
Our sincere hope is that this introduction, together with the whole book the readers are about to read, credits the visionary genius of Yoneda: category theory has few theorems, and one of them is a lemma. The \emph{Yoneda lemma}, in its myriad of incarnations, is certainly a cornerstone of structural thinking, way before than of category theory: if anything more was needed to revere Nobuo Yoneda, let this be co\fshyp{}end calculus.

The paper \cite{Yoneda} starts introducing co\fshyp{}ends in the following way:
\begin{quote}
	Let $\C$ be a category. By a \emph{left $\C$\hyp{}group} we mean a covariant functor $M$ of $\C$ with values in the category $\Ab$ of abelian groups and homomorphisms. [\dots\unkern] Also by a $\C^*$\hyp{}group (or a \emph{right $\C$\hyp{}group}) we mean a contravariant functor $K : \C \to \Ab$ [\dots\unkern]. Functors of several variables with values in $\Ab$ will accordingly be called $\B$-$\C$\hyp{}groups, $\B^*$-$\C$\hyp{}groups, etc.

	Let $H$ be a $\C^*$-$\C$\hyp{}group, and $G$ an additive group. By a \emph{balanced homomorphism} $\mu : G \din H$ we mean a system of homomorphisms $\mu(C) : G \to H(C,C)$ defined for all objects $C\in\C$ such that for every map $\gamma : C \to C'$ in $\C$ commutativity holds in the diagram
	\[\notag
		\xymatrix{
			G\ar[d]_{\mu(C')}\ar[r]^{\mu(C)} & H(C,C) \ar[d]^{H(C,\gamma)}\\
			H(C', C') \ar[r]_{H(\gamma,C')}& H(C,C').
		}
	\]
	Also by a balanced homomorphism $\lambda : H \din G$ we mean a system of homomorphisms $\lambda(C) : H(C,C) \to G$ defined for all objects $C\in\C$ such that for every map $\gamma : C \to C'$ in $\C$ commutativity holds in the diagram
	\[\notag
		\xymatrix{
			H(C',C) \ar[r]^{H(C',\gamma)}\ar[d]_{H(\gamma,C)}& H(C',C') \ar[d]^{\lambda(C')}\\
			H(C,C) \ar[r]_{\lambda(C)}& G
		}
	\]
\end{quote}
Of course, here a left/right `$\C$\hyp{}group' is merely an $\Ab$\hyp{}enriched presheaf (covariant or contravariant) with domain $\C$. The above paragraphs define the fundamental notions we will use throughout the entire book, \emph{wedges} and \emph{cowedges}. These are exactly natural maps from/to a constant, that vary taking into account the fact that $H(C,C)$ depends both covariantly and contravariantly on $C$. There is a category of such co\fshyp{}wedges, and a process dubbed \emph{co\fshyp{}integration} picks the initial and terminal objects of such categories:
\begin{quote}
	An additive group $I$ together with a balanced $\theta : H \din I$ is called \emph{integration} of a $\C^*$-$\C$\hyp{}group $H$ if it is universal among balanced homomorphisms from $H$, \ie if for any other balanced homomorphism $\lambda : H \din G$ there is a unique morphism $\zeta : I \to G$ such that $\zeta\circ \theta(C) = \lambda(C)$ for every object $C\in\C$.

		[Integrations and cointegrations] are [\dots\unkern] given as follows: for a map $\gamma :C \to C'$ in $\C$ we put $H(\gamma) = H(C',C)$, $H(\gamma^*) = H(C,C')$, and define homomorphisms
	\index{Integration}
	\index{Cointegration|see{Integration}}
	\begin{gather*}
		\partial_\gamma : H(\gamma) \to H(C,C) \oplus H(C',C') \\
		\delta_\gamma : H(C,C)\oplus H(C',C') \to H(\gamma^*)
	\end{gather*}
	by
	\begin{gather*}
		\partial_\gamma(h') = h'\circ\gamma \oplus (-\gamma\circ h') \\
		\delta_\gamma(h\oplus h'') = \gamma\circ h - h''\circ \gamma
	\end{gather*}
	Denote by $\Sigma_0$ and $\Pi^0$ the direct sum $\sum_{C\in\C} H(C,C)$ and the direct product $\prod_{C\in\C} H(C,C)$ respectively. Also, denote by $\Sigma_1$ and $\Pi^1$ the direct sum $\sum_{\gamma\in\hom(\C)} H(\gamma)$ and $\prod_{\gamma\in\hom(\C)} H(\gamma^*)$ respectively. Then $\partial_\gamma, \delta_\gamma$ are extended to homomorphisms
	\[\notag\partial : \Sigma_1 \to \Sigma_0 \qquad \delta_{CC'} : \Pi^0 \to \Pi^1.\]
\end{quote}
Now \cite{Yoneda} proves that the `integration' $\int_C H$ and the `cointegration' $\int^*_C H$ of a $\C^*$-$\C$\hyp{}group are respectively given by the cokernel of $\partial_{CC'}$, and by the kernel of $\delta_{CC'}$, for suitably defined maps $\partial_{CC'}$ and $\delta_{CC'}$.

To avoid confusion, we stress that in this definition Yoneda employed the \emph{opposite} choice of terminology that we will introduce later on: an \emph{integration} for $H$ is a coend $\int^C H(C,C)$, but Yoneda denotes it $\int_C H(C,C)$; a \emph{cointegration} is an end $\int_C H(C,C)$ but Yoneda denotes it $\int_C^* H(C,C)$. Always be careful if you consult \cite{Yoneda}, or references from the same age.

\smallskip
The coend $\int^C H$ of a functor $H : \C^\opp\times \C \to \D$ is a colimit (in the particular case of $\Ab$\hyp{}enriched functors, a cokernel) built out of $H$, and the end $\int_C H$ is a limit; precisely the kernel of a certain group homomorphism. Once this terminology has been set up, given two functors $F,G : \C \to \D$, if we let $H$ be the functor $(C,C')\mapsto \D(FC,GC')$, a family of arrows $\alpha_C : FC \to GC$ forms the components of a natural transformation $\alpha : F \To G$ if and only if the components $\alpha_C : FC \to GC$ lie in the kernel of a `differential' $\bar\delta : \prod_{C\in\C} H(C,C) \to \prod_{\gamma \in\hom(\C)} H(\gamma^*)$, obtained in the obvious way `gluing' all the $\delta_{CC'}$ together.

\begin{quote}
	In dealing with functors of more variables, we shall often inscribe $x$ (or $y,z$) to indicate the two entries to be considered in the (co)integration, namely we write
	\[\int_{X\in \C} H(\dots,X,\dots, X,\dots). \]
	This is based on the following fact: let $H, H'$ be $\C$-$\C^*$\hyp{}groups, and let $\theta : H \to \int_\C H$, $\theta' : H' \to \int_\C H'$ be the integrations. Then a natural transformation $\eta : H \To H'$ induces a unique homomorphism $\int_\C \eta : \int_\C H \to \int_\C H'$ such that $\left(\int_\C \eta\right)\circ \theta(C) = \theta'(C)\circ \eta(C,C)$. Thus if $H$ is a $\B$-$\C^*$-$\C$\hyp{}group, then $\int_{X\in\C}^{(*)} H(B,X,X)$ is a $\B$\hyp{}group. On this account, for an $\A$-$\B^*$-$\B$-$\C$-$\C^*$\hyp{}group $H$ we have
	\[\notag
		\int_{Y\in\B}\int_{X\in\C} H(A,Y,Y,X,X) =\int_{X\in\C}\int_{Y\in\B} H(A,Y,Y,X,X)
	\]
\end{quote}
Here Yoneda introduces one of the pillars of coend calculus, the \emph{Fubini rule}, \ie the fact that the result of a co\fshyp{}integration is the same irregardless of the order of integration; this is ultimately just a consequence of the functoriality of the assignment sending a $\C$-$\C^*$\hyp{}group $H$ into its co\fshyp{}integration. Of course, there is nothing special about the codomain of $H$ being the category of abelian groups: any sufficiently co\fshyp{}complete category $\D$ will do, as long as the co/integrations involved exist.

In modern terms, the Fubini rule can be obtained as a consequence of a much deeper, and hopefully more enlightening, result: we prove it in our \ref{fubozzo}.
\begin{quote}
	Next for a $\B$\hyp{}group $M$ and a $\C$\hyp{}group $N$, $M\otimes N : N(C)\otimes M(B)$ is a $\B$-$\C$\hyp{}group, and $\underline{\hom}(M,N) = \hom(MB, NC)$ is a $\B^*$-$\C$\hyp{}group. For an $\A$\hyp{}group $M$ and a $\B$-$\C$-$\C^*$\hyp{}group $H$ we have:
	\begin{gather*}
		\int_{X\in\C} MA\otimes H(B,X,X) = MA \otimes \int_{X\in\C} H(B,X,X) \\
		\int^*_{X\in\C} \hom(MA, H(B,X,X)) = \hom\left(MA, \int^*_{X\in\C}H(B,X,X)\right) \\
		\int^*_{X\in\C} \hom(H(B,X,X), MA) = \hom\left( \int^*_{X\in\C} H(B,X,X), MA \right)
	\end{gather*}
\end{quote}
As an immediate consequence of these statements we get another fundamental building\hyp{}block of a coend `calculus': given two functors $M : \B \to \Ab$ and $N : \B^\opp\to \Ab$, they can be \emph{tensored} by the integration $M \boxtimes N := \int_{B\in\B} N(B)\otimes_{\Z} M(B)$.

A rather interesting perspective on this construction is the following: the result remains true when the $\Ab$-category $\B$ has a single object, so it is merely a ring $B$: in such a case, a functor $M : B \to \Ab$ is a left module, and a functor $N : B^\opp\to \Ab$ is a right module; the integration (or in modern terms, the coend) $\int^B N\otimes_{\Z} M$ in this case is exactly the tensor product of $B$\hyp{}modules: it has the universal property of the cokernel of the map
\[\notag
	\bigoplus_{b\in B} M\otimes_{\Z} N \xto{\varrho} M\otimes_{\Z} N
\]
defined by $\varrho (b,m,n) = b.m - n.b$.

\smallskip
As the reader might now suspect, few analogies of them are more fruitful than the one between modules over which a monoid object acts, and presheaves $\C \to \Set$.

\paragraph{Structure of the book.}
We shall now briefly review the structure of the book: in the first three chapters we outline the basic rules of co\fshyp{}end calculus; after having defined a co\fshyp{}end as a universal object and having proved that it can be characterised as a co\fshyp{}limit, we start denoting such object as an integral $\int_C T(C,C)$ or $\int^C T(C,C)$. This notation is motivated by the fact that co\fshyp{}ends `behave like integrals' in that a Fubini rule of exchange holds: see \ref{fubozzo}.

Then we introduce the first rules of the calculus: the Yoneda lemma \ref{lem:the-real-yoda} can be restated in terms of a certain coend computation, and pointwise Kan extensions can be computed by means of a co\fshyp{}end.

After this, we study the single case of left Kan extensions along the Yoneda embedding: in some sense, the theory of such extensions alone embodies ``all'' category theory.

The subsequent chapter begins to introduce more modern topics described by means of co\fshyp{}end calculus. The theory of \emph{weighted co\fshyp{}limits}, of \emph{profunctors} and \emph{operads} is a cornerstone of `formal' approaches to category theory. Weighted co\fshyp{}limits are the correct notion of co\fshyp{}limit in an enriched or formal\hyp{}categorical (see \ref{def:yoneda-struc} and \cite{Graya}) setting; profunctors are a bicategory where one can re\hyp{}enact all category theory, and are deeply linked to categorical algebra and representation; operads, initially introduced as a technical mean to solve an open problem in homotopy theory, constitute now the common ground where universal algebra and algebraic topology meet. The final point of the chapter will be \ref{curien_main}, where we draw a tight link between profunctors and operads.

The subsequent chapter studies higher\hyp{}dimensional analogues of co\fshyp{}ends; first, we study co\fshyp{}ends in 2\hyp{}categories; then we move up to infinity and study homotopy\hyp{}coherent analogues of co\fshyp{}ends, in simplicial categories \cite{bergner2007model}, quasicategories \cite{HTT}, model categories \cite{Hov} and derivators.

Appendix A serves as a short introduction to category theory: it fixes the notation we have employed in the previous chapters. Elementary mathematics is a prerequisite to appreciate it, but we will introduce most categorical jargon from scratch.

\smallskip
Each chapter has a short introduction to its content, in the form of a small abstract; this allows the interested reader to get a glimpse into the content and fundamental results of each chapter (often, one or two main theorems). We believe this format is easier to consult than a comprehensive survey of each chapter given all at once in the introduction, so we felt free to keep this short introductory paragraph on the content of the book pretty terse.

Several exercises follow each chapter of the book; there are questions of every level, sometimes easy, sometimes more difficult; some of them make the reader rapidly acquainted with the computational approach to category theory offered by co\fshyp{}end calculus; some others shed a new light on old notions. In approaching them, we advise you to avoid element-wise reasoning; instead, find either an abstract argument, or a `deduction\hyp{}style' one.

Some of the exercises are marked with an $\awful$ symbol (eyes wide with fear): this means they are more difficult and less well\hyp{}posed questions than the others. This can happen on purpose (and thus part of the exercise is understanding what the question is) or not (and thus the question and its answer are not completely clear even to the author). In this second case, it is likely that a complete answer might result in new Mathematics that the solvers are encouraged to develop.

Other kinds of `small eyes' are present along the book: the paragraphs decorated with a \upeyes contain material that can be skipped at first reading, or material that deepens a prior topic in a not\hyp{}so\hyp{}interesting detour; \righteyes is used to signal key remarks and more generally important material that we ask the reader to digest properly and analyse in full detail.
\paragraph{Notation.}
Having to deal with many different sources along the exposition, the hope to maintain a coherent choice of notation throughout the whole book is wishful thinking; however, the author did his best to provide a coherent enough one, striving to make it at the same time expressive and simple.

In general, 1-dimensional category-like structures will be denoted as calligraphic letters $\C,\D,\dots$; objects of $\C$ are denoted $C,C'\dots \in\C$. Instead, 2-categories are often denoted with a sans-serif case $\tCat, \sfK, \sfA,\dots$; in this case, an object of the 2-category of small categories is denoted $\C\in\tCat$, but an object of an abstract 2-category is denoted $A\in\sfK$.

Functors between categories are denoted as capital Latin letters like $F,G,H,K$ and suchlike (although there can be little deviations to this rule); the category of functors $\C\to \D$ between two categories is almost always denoted as $\Cat(\C,\D)$ (or less often $[\C,\D]$; this will be done especially when $[\C,\D]$ is regarded as the internal hom of the closed structure in $\tCat$, or when it is necessary to save some space); the symbols $\firstblank$, $\secondblank$ are used as placeholders for the ``generic argument'' of a functor or bifunctor (they mark temporal precedence of saturation of a variable); morphisms in the category $\Cat(\C,\D)$ (\ie natural transformations between functors) are often written in lowercase Greek, or lowercase Latin alphabet, and collected in the set $\Cat(\C,\D)(F,G)$.

The simplex category $\bDelta$ is the \emph{topologist's delta} (opposed to the \emph{algebraist's delta} $\bDelta_+$ which has an additional initial object $[-1]:= \varnothing$), having objects \emph{nonempty} finite ordinals $[n]:=\{0<1\dots<n\}$; we denote $\Delta[n]$ the representable presheaf on $[n]\in\bDelta$, \ie the image of $[n]$ under the Yoneda embedding of $\bDelta$ in the category $\mathsf{sSet} = \widehat{\bDelta}$ of simplicial sets. More generally, we indicate the Yoneda embedding of a category $\C$ into its presheaf category with $\yon_\C$ --or simply $\yon$--, \ie with the hiragana symbol for ``yo''.
\paragraph{Acknowledgments.}
Writing this book spanned the last five years of my life; it would probably take five more years to do it in the form the topic deserves to be told. I did my best, and if I managed to be proud of the result, it's because I've been helped a lot.

The draft of this book accompanied me in all sorts of places: this has made the list of people who helped and guided me through it quite long. After all this time, it remains all the more true today that `in some sense, I am not the only author of this note', as I used to write in the early stages of this project.

Several people participated and extensively commented such draft(s); each of them helped to make the present book a little better than it could have been if I would have been alone.

I would like to thank
	{\sf S. Ariotta} (for providing me more help than I deserved, in the proof of the $\infty$-Fubini rule in \ref{fubi_for_infty})%
, {\sf J. D. Christensen}%
, {\sf I. Di Liberti}%
, {\sf E. de Oliveira Santos} (for an attentive proofreading when my eyes were bleeding)%
, {\sf D. Fiorenza} (a friend and an advisor I will never be able to refund for his constant support)%
, {\sf A. Gagna}%
, {\sf N. Gambino} (for offering me me the opportunity to discuss the content of this note in front of his students in Leeds; in just a few days I realised years of meditation were still insufficient to teach co\fshyp{}end calculus)%
, L. (for accompanying me there, and for being the best \emph{mate} for a much longer trip)%
, {\sf F. Genovese} (for believing in co\fshyp{}end calculus with all his heart, and for making dg\hyp{}categories more understandable through it)%
, {\sf A. Joyal} (for letting me sketch the statement of \ref{elts-as-coend} on the back of a napkin in a cafe in Paris)%
, {\sf G. Marchetti}%
, {\sf A. Mazel\hyp{}Gee}%
, {\sf B. Milewski}%
, {\sf G. Mossa}%
, {\sf D. Palombi} (for sharing with me his truffles, his friendship, his monoidal categories, and his typographical kinks)%
, {\sf D. Ravenel} (for believing in the title of this draft more than I will ever do)%
, {\sf M. Román}%
, {\sf G. Ronchi} (for accepting my harsh ways while I was learning how to be a mentor, and because he Believes that Mathematics is Beautiful)%
, {\sf J. Rosick\'{y}} (for he trusted in me and gave me pearls of wisdom)%
, {\sf E. Riehl}%
, {\sf E. Rivas}%
, {\sf T. Trimble}%
, {\sf M. Vergura}%
, {\sf C. Williams}%
, {\sf K. Wright}%
, {\sf N. Yanofsky}, and a countless number of other people who pointed out mistakes, faulty arguments, typos.%

\smallskip
A separate thanks goes to
\textsf{S\@.}
, \textsf{P\@.}
, \textsf{G\@.}
, \textsf{C\@.}: 
they opened me their doors when I was frail and broken\hyp{}hearted. Writing means, always, opening a wound, exposing your tendons to the blade of a sharp knife. During these years I have been protected by the warmth of your smile, of your houses, of your brilliance, of your friendship.

\smallskip
Finally, a special thanks goes to: Jelly Roll Morton, Robert Musil, Daisetzu Teitar\={o} Suzuki, Raymond Roussel, Kurt Schwitters, Vieira da Silva, Akutagawa, Anton Webern, Greta Garbo, Jos\'e Lezama Lima, Bu\~{n}uel, Louis Armstrong, Borges, Michaux, Dino Buzzati, Max Ernst, Pevsner, Gilgamesh (?), Garcilaso, Arcimboldo, René Clair, Piero di Cosimo, Wallace Stevens, Izak Dinesen, \sout{Rimbaud, Picasso, Chaplin, Alban Berg}.

Various institutions supported me through the many years necessary to finish this book. I would like to thank the University of Western Ontario, Masaryk University, The Max Planck Institute for Mathematics in Bonn for providing me time and money to keep writing, and the University of Coimbra, that added to the deal a \emph{completely} distraction\hyp{}free environment. Additional help came from Taltech, where I am working while this book approaches its final form, and from all Estonia: \emph{aitäh!}

\smallskip
Ringrazio per ultimo mio padre, che non comprenderebbe questa dedica se non la scrivessi in italiano, e che non leggerà mai questo libro `assorto com'è in altra e lontana scienza'. A sei anni mi fece sedere sulle sue ginocchia, e scrisse su un foglio il numero 1 seguito da nove zeri, per farmi capire quanto è grande un miliardo: si può dire tutto questo sia iniziato lì.

\smallskip
Grazie a tutti.
\cleardoublepage
\vspace*{\fill}

\begin{figure}[h]
	\begin{adjustbox}{width=\textwidth}
		\thispagestyle{empty}
		\begin{tikzpicture}[
				x={ 3cm*sin(60)},
				y={-3cm*cos(60)},
				sephirot/.style={
						draw
						, fill=white
						, line width=.8pt
						, shape=circle
						, inner sep=0sp
						, minimum size=1.75cm
						, text width=1.4cm
						, font=\fontsize{6pt}{8pt}\selectfont\scshape
						, align=center
					},
				sibihil/.style={
						line width=.4pt
						, double=lightgray
						, double distance=8pt
					},
				mosaic/.style={
						line width=.2pt
						, gray
					},
				frame/.style={
						line width=.4pt
					}
			]
			\coordinate (S0) at ( 0,0);
			\coordinate (S1) at (-1,1);
			\coordinate (S2) at (+1,1);
			\coordinate (S3) at (-1,3);
			\coordinate (S4) at (+1,3);
			\coordinate (S5) at ( 0,4);
			\coordinate (S6) at (-1,5);
			\coordinate (S7) at (+1,5);
			\coordinate (S8) at ( 0,6);
			\coordinate (S9) at ( 0,8);
			\draw[frame] (0,4) circle (3*3cm);
			\draw[frame] (0,4) circle (3*3cm+.25cm);
			\draw[mosaic] (-2,2) circle (3cm);
			\draw[mosaic] (-2,4) circle (3cm);
			\draw[mosaic] (-2,6) circle (3cm);
			\draw[mosaic] (-1,1) circle (3cm);
			\draw[mosaic] (-1,3) circle (3cm);
			\draw[mosaic] (-1,5) circle (3cm);
			\draw[mosaic] (-1,7) circle (3cm);
			\draw[mosaic] ( 0,0) circle (3cm);
			\draw[mosaic] ( 0,2) circle (3cm);
			\draw[mosaic] ( 0,4) circle (3cm);
			\draw[mosaic] ( 0,6) circle (3cm);
			\draw[mosaic] ( 0,8) circle (3cm);
			\draw[mosaic] (+1,1) circle (3cm);
			\draw[mosaic] (+1,3) circle (3cm);
			\draw[mosaic] (+1,5) circle (3cm);
			\draw[mosaic] (+1,7) circle (3cm);
			\draw[mosaic] (+2,2) circle (3cm);
			\draw[mosaic] (+2,4) circle (3cm);
			\draw[mosaic] (+2,6) circle (3cm);
			\draw[mosaic] (0,-2)
			arc (-150:-90:3cm)
			arc (-150:-90:3cm)
			arc (-150:-90:3cm)
			arc (150:210:3cm)
			arc (150:210:3cm)
			arc (150:210:3cm)
			arc (90:150:3cm)
			arc (90:150:3cm)
			arc (90:150:3cm)
			arc (30:90:3cm)
			arc (30:90:3cm)
			arc (30:90:3cm)
			arc (-30:30:3cm)
			arc (-30:30:3cm)
			arc (-30:30:3cm)
			arc (-90:-30:3cm)
			arc (-90:-30:3cm)
			arc (-90:-30:3cm)
			;
			\draw [mosaic] (0,-2)
			arc (-210:-30:3cm)
			arc (-210:-90:3cm)
			arc (-270:-90:3cm)
			arc ( 30:210:3cm)
			arc ( 30:150:3cm)
			arc (-30:150:3cm)
			arc (-90:90:3cm)
			arc (-90:30:3cm)
			arc (-150:30:3cm)
			;
			\draw [mosaic] (0,10)
			arc (-30:150:3cm)
			arc (-30:90:3cm)
			arc (-90:90:3cm)
			arc (-150:30:3cm)
			arc (-150:-30:3cm)
			arc (-210:-30:3cm)
			arc (-270:-90:3cm)
			arc (-270:-150:3cm)
			arc (30:210:3cm)
			;
			\draw[sibihil] (S9) -- (S8) node[midway] {\textcjheb{t|}};
			\draw[sibihil] (S8) -- (S7) node[midway] {\textcjheb{n|}};
			\draw[sibihil] (S7) -- (S4) node[midway] {\textcjheb{k|}};
			\draw[sibihil] (S4) -- (S2) node[midway] {\textcjheb{b|}};
			\draw[sibihil] (S2) -- (S0) node[midway] {\textcjheb{h|}};
			\draw[sibihil] (S0) -- (S1) node[midway] {\textcjheb{w|}};
			\draw[sibihil] (S1) -- (S3) node[midway] {\textcjheb{g|}};
			\draw[sibihil] (S3) -- (S6) node[midway] {\textcjheb{p|}};
			\draw[sibihil] (S6) -- (S8) node[midway] {\textcjheb{l|}};
			\draw[sibihil] (S5) -- (S8) node[pos=.65] {\textcjheb{r|}};
			\draw[sibihil] (S5) -- (S7) node[midway] {\textcjheb{y|}};
			\draw[sibihil] (S5) -- (S4) node[midway] {\textcjheb{.h|}};
			\draw[sibihil] (S5) -- (S2) node[pos=.7] {\textcjheb{.t|}};
			\draw[sibihil] (S5) -- (S1) node[pos=.7] {\textcjheb{`|}};
			\draw[sibihil] (S5) -- (S3) node[midway] {\textcjheb{.s|}};
			\draw[sibihil] (S5) -- (S6) node[midway] {\textcjheb{s|}};
			\draw[sibihil] (S1) -- (S4) node[pos=.3] {\textcjheb{q|}};
			\draw[sibihil] (S2) -- (S3) node[pos=.3] {\textcjheb{z|}};
			\draw[sibihil] (S5) -- (S0) node[midway] {\textcjheb{d|}};
			\draw[sibihil] (S1) -- (S2) node[midway] {\textcjheb{/s|}};
			\draw[sibihil] (S3) -- (S4) node[midway] {\textcjheb{'|}};
			\draw[sibihil] (S6) -- (S7) node[midway] {\textcjheb{m|}};
			\node[sephirot] at (S0) {Wedges and co\fshyp{}ends};
			\node[sephirot] at (S1) {Yoneda and Kan};
			\node[sephirot] at (S2) {Basic category theory};
			\node[sephirot] at (S3) {Profunctors};
			\node[sephirot] at (S4) {Nerves and realisations};
			\node[sephirot] at (S5) {Weighted co\fshyp{}limits};
			\node[sephirot] at (S6) {Addenda};
			\node[sephirot] at (S7) {Operads};
			\node[sephirot] at (S8) {Higher co\fshyp{}ends};
			\node[sephirot] at (S9) {Preface};
		\end{tikzpicture}
	\end{adjustbox}
	\caption*{Conceptual dependencies between the chapters, embedded in a sephirotic tree.}
\end{figure}
\vspace*{\fill}
\cleardoublepage

\pagenumbering{arabic}

\chapter{Dinaturality and co\fshyp{}ends}\label{section:due}
\begin{abstract}
	Naturality of a family of morphisms $\alpha_C : FC \to GC$ defines the correct notion of map between functors $F,G$; yet, it is not capable to describe more subtle interactions that can occur between $F$ and $G$, for example when both functors have a product category like $\C^\opp\times \C$ as domain. A transformation that takes into account the fact that $F,G$ act on morphisms once covariantly and once contravariantly is called \emph{dinatural}.

	As ill-behaved as it may seem (in general, dinatural transformations can't be composed), this notion leads to the definition of a \emph{co\fshyp{}wedge} and \emph{co\fshyp{}end} for a functor $T : \C^\opp\times \C\to \D$: a dinatural transformation having constant co\fshyp{}domain, and a suitable universal property. This is in perfect analogy with the theory of co\fshyp{}limits: universal natural transformations from/to a constant functor. Unlike colimits, however, co\fshyp{}ends support a \emph{calculus}, that is a set of inference rules allowing to mechanically prove nontrivial statements as initial and terminal points of a chain of deductions.

	The purpose of this chapter, and indeed of the entire book, is to familiarise its readers with the rules of calculus.
\end{abstract}
\epigraph{Los idealistas arguyen que las salas hexagonales son una forma necesaria del espacio absoluto o, por lo menos, de nuestra intuición del espacio.}{J.L. Borges --- \emph{La biblioteca de Babel}}
\section{Supernaturality}
We chose to let the name ``supernaturality'' describe the two sorts of generalisations of naturality for functors that we will investigate throughout the book: \emph{di}naturality, in \ref{dinatra}, and \emph{extra}naturality, in \ref{extranatural}.
\subsection{Dinaturality}
This first section starts with a simple example. We denote $\Sets$ the category of sets and functions, considered with its natural cartesian closed structure (see \ref{ex:adjoints}.\ref{ad:cin}): this means that we have a bijection of sets
\[\label{adjisos}
	\Sets(A\times B,C)\cong \Sets(A, C^B)
\]
natural in all three arguments, if we let $C^B$ denote the set of functions $f : B\to C$.
\index{Adjunction!unit and counit}
\index{Counit!--- of $\firstblank\times A \dashv (\firstblank)^A$}
The bijection above is defined by the maps
\begin{align*}
	(f : A\times B \to C) & \mapsto A \xto{\eta_{A,(B)}}(A\times B)^B \xto{f^B} C^B                \\
	(g : A \to C^B)       & \mapsto A\times B \xto{g\times B} C^B\times B \xto{\epsilon_{C,(B)}} C
\end{align*}
by means of suitable \emph{unit} and \emph{counit} maps $\eta$ and $\epsilon$ (see \ref{unit_counit}) witnessing the adjunction. Let us concentrate on the counit map alone (a dual reasoning will yield similar conclusions for the unit): it is a natural transformations having components
\[
	\{\epsilon_{X,(B)} : X^B\times B\to X \mid X \in\Set\}.
\]
This family of functions sends a pair $(f,b)\in X^B\times B$ to the element $fb\in X$, and thus deserves the name of \emph{evaluation}.

For the purpose of our discussion, we shall consider this family of morphisms not only natural in $X$ (as every counit morphism), but also mutely depending on the variable $B$ in its codomain. This means that $X^B\times B$ is the image of the pair $(B,B)$ under the functor $(U,V)\mapsto X^U\times V$, and $X$ can be regarded similarly as the image of $(B,B)$ under the constant functor in $X$. Both functors have thus `type' $\Set^\opp\times\Set \to \Set$.

The evaluation maps $\epsilon_{X,(B)}$ however do not vary naturally in the variable $B$; the most we can say is that for each function $f\in \Sets(B,B')$ the following square is commutative:
\[
	\vcenter{\xymatrix{
	X^{B'}\times B \ar[r]^{X^f\times B}\ar[d]_{X^{B'}\times f}& X^B\times B \ar[d]^{\epsilon}\\
	X^{B'}\times B' \ar[r]_{\epsilon}& X.
	}}
\]
This relation doesn't remind naturality so much, but it can be easily deduced from the request that the adjunction isomorphisms \eqref{adjisos} are natural in the variable $B$; in fact, such naturality imposes the commutativity of the square
\[
	\vcenter{\xymatrix{
			\Set(A, X^{B'})\ar[d]_{\Set(A, X^f)} \ar[r] & \Set(A\times B', X)\ar[d]^{\Set(A\times f, X)}\\
			\Set(A, X^B) \ar[r] & \Set(A\times B, X)
		}}
\]
for an arrow $f : B\to B'$ (the horizontal maps are the adjunction isomorphisms $\firstblank\mapsto \epsilon\circ (\firstblank\times B)$), and this in turn entails that we have an equation
\begin{align}
	\epsilon_{X,(B')}\circ (u\times B')\circ (A\times f)     & = \epsilon_{X,(B)}\circ (X^f\circ u)\times B \label{chista}              \\
	\epsilon_{X,(B')}\circ (X^{B'}\times f)\circ (u\times B) & = \epsilon_{X,(B)}\circ (X^f\times B)\circ (u\times B)\label{chillautra}
\end{align}
for every $u : A\to X^{B'}$. But since this is an equality for every such $u$, then the functions $\epsilon_{X(B')}\circ (X^{B'}\times f)$ and $\epsilon_{X(B)}\circ (X^f\times B)$ must also be equal.

So, it would seem that there's no way to frame the diagram above in the usual context of naturality for a transformation of functors. Fortunately, a suitable generalisation of naturality (a `super\hyp{}naturality' condition), encoding the above commutativity, is available to describe this and other similar phenomena.

As already said, the correspondence $(B,B')\mapsto C^B\times B'$ is a functor with domain $\Set^\opp\times \Set$; it turns out that these functors, where the domain is a product of a category with its opposite, supports a notion of \emph{di}naturality besides the classical naturality; this notion is more suited to capture the phenomenon we just described: in fact, most of the transformations that are canonical, depending on two variables $(C,C')\in\C^\opp\times \C$, but not natural, can be seen as dinatural.
\begin{definition}[Dinatural transformation]\label{dinatra}
	\index{Transformation!dinatural ---}
	\index{_aaa_dinat@$\din$}
	Let $\C,\D$ be two categories. Given two functors $P,Q : \C^\opp\times \C\to \D$ a \emph{dinatural transformation} $\alpha : P\din  Q$ consists of a family of arrows
	\[
		\alpha_C : P(C,C)\to Q(C,C)
	\]
	indexed by the objects of $\C$ and such that for any $f : C \to C'$ the following diagram commutes
	\[\label{dinana}
		\vcenter{\xymatrix{
		P(C',C) \ar[r]^{P(f,C)}\ar[d]_{P(C',f)} & P(C,C) \ar[r]^{\alpha_C} & Q(C,C) \ar[d]^{Q(C,f)}\\
		P(C',C') \ar[r]_{\alpha_{C'}} & Q(C',C') \ar[r]_{Q(f, C')} & Q(C,C').
		}}
	\]
\end{definition}
\begin{remark}
	The notion of dinaturality takes into account the fact that a functor $P : \C^\opp\times \C \to \D$ maps at the same time two `terms' of the same `type' $\C$, once covariantly in the second component, and once contravariantly in the first: on arrows $f : C\to C'$ the functor $P$ acts in fact as follows:
	\[\label{gnogni}
		\vcenter{\xymatrix{
				& P(C',C) \ar[dr]^{P(C',f)} \ar[dl]_{P(C,f)}& \\
				P(C,C) && P(C',C')
			}}
	\]
	Given two such functors, say $P,Q : \C^\opp\times \C \to \D$, we can consider the two diagrams \eqref{gnogni} and
	\[
		\vcenter{\xymatrix{
				Q(C,C)\ar[dr]_{Q(C,f)} && Q(C',C')\ar[dl]^{Q(f,c)} \\
				& Q(C,C') &
			}}
	\]
	In the same way a natural transformation $F \To G$ can be seen as a family of maps that `fill the gap' between $F(f)$ and $G(f)$ in a commutative square, a \emph{dinatural} one between $P$ and $Q$ can be seen as a way to close the hexagonal diagram connecting the action on arrows of $P$ to the action on arrows of $Q$:
	\[
		\vcenter{\xymatrix{
				& P(C',C) \ar[dr]^{P(C',f)} \ar[dl]_{P(C,f)}& \\
				P(C,C) \ar@{.>}[d]&& P(C',C') \ar@{.>}[d]\\
				Q(C,C)\ar[dr]_{Q(C,f)} && Q(C',C')\ar[dl]^{Q(f,C')} \\
				& Q(C,C') &
			}}
	\]
	This is precisely the diagram drawn in \eqref{dinana}.
\end{remark}
\begin{remark}
	If we let $P_C$ be the functor $(U,V)\mapsto C^U\times V$, the counit components $\epsilon_{C(B)} : P_C(B,B)\to C$ of the cartesian closed adjunction form a dinatural transformation $\epsilon : P_C \din \Delta_C$, where $\Delta_C$ is the constant functor at $C$.
\end{remark}
Such dinatural transformations, having constant codomain, deserve a special name:
\begin{definition}[Co/wedge]\index{Functor!Wedge for a ---}\index{Co/wedge}
	Let $P : \C^\opp\times\C\to\D$ be a functor;
	\begin{enumtag}{wc}
		\item A \emph{wedge} for $P$ is a dinatural transformation $\Delta_D\din P$ from the constant functor on the object $D\in\D$ (we often denote such constant functor simply by the name of the constant, $D : \C^\opp\times\C\to\D$), defined by the rules $(C,C')\mapsto D$, $(f,f')\mapsto \id_D$.
		\item Dually, a \emph{cowedge} for $P$ as above is a dinatural transformation $P \din \Delta_D$ having codomain the constant functor on the object $D\in\D$.
	\end{enumtag}
\end{definition}
\begin{remark}
	Wedges for a fixed functor $P$ as above form the class of objects of a category $\wed(P)$, where a morphism of wedges is a morphism between their domains that makes an obvious triangle commute; given two wedges $\alpha : D\din P$ and $\alpha' : D'\din P$ a morphism $u : $ consists of an arrow $u : D\to D'$ such that the triangle
	\[\label{chistu_cowedge}
		\vcenter{\xymatrix{
		D \ar[rr]^u \ar[dr]_{\alpha_{CC}}&& D'\ar[dl]^{\alpha'_{CC}} \\
		&P(C,C)&
		}}
	\]
	is commutative for every component $\alpha_{CC}$ and $\alpha'_{CC}$. (Note the role of quantifiers: the same $u$ makes \eqref{chistu_cowedge} commute for \emph{every} component of the wedges.)

	Dually, there is a category $\cwed(P)$ of cowedges for $P$, where morphisms of cowedges are morphisms between codomains (of course there is a relation between the two categories: cowedges for $P$ coincide with the opposite category of wedges for the opposite functor).
\end{remark}
We now define the end of $P$ as a terminal object in $\wed(P)$, and the coend as an initial object in $\cwed(P)$.
\index{Co/end}
\begin{definition}[Co/end]
	Let $P : \C^\opp\times\C\to \D$ be a functor;
	\begin{itemize}
		\item The \emph{end} of $P$ consists of a terminal wedge $\omega : \eend(P) \din P$; the object $\eend(P) \in \D$ itself is often called \emph{the end} of the functor.
		\item Dually, the \emph{coend} of $P$ as above consists of an initial cowedge $\alpha : P \din \coend(P)$; similarly, the object $\coend(P)$ itself is often called the coend of $P$.
	\end{itemize}
\end{definition}
Spelled out explicitly, the universality requirement means that for any other wedge $\beta : D\din P$ the diagram
\[
	\vcenter{\xymatrix{
	D\ar@{.>}[dr]^h \ar@/^1pc/[drr]^{\beta_C}\ar@/_1pc/[ddr]_{\beta_{C'}} \\
	& \eend(P) \ar[r]^{\omega_C}\ar[d]_{\omega_{C'}} & P(C,C) \ar[d]^{P(1,f)}\\
	&P(C',C') \ar[r]_{P(f,1)}& P(C,C')
	}}
\]
commutes for a unique arrow $h : D\to \eend(P)$, for every arrow $f : C \to C'$. Note again the role of quantifiers: the arrow $h$ is the same for every component of the wedge. A dual diagram can be depicted for the coend of $P$.
\begin{remark}[Functoriality of ends]\label{fun_4_ends}
	Given a natural transformation $\eta : P\To P'$ between functors $P,P' : \C^\opp\times \C\to \D$ there is an induced arrow $\eend(\eta) : \eend(P)\to \eend(P')$ between their ends, as depicted in the diagram
	\[
		\vcenter{\xymatrix@R=5mm@C=5mm{
		&\eend(P') \ar[rr]^{\omega'_{C'}}\ar[dd]|\hole&& P'(C',C')\ar[dd]^{P'(f,C')}\\
		\eend(P) \ar@{.>}[ur]^{\eend(\eta)}\ar[rr]^(.3){\omega_{C'}}\ar[dd]_{\omega_C}&& P(C',C')\ar[dd]\ar[ur]_{\eta_{C'C'}}\\
		&P'(C,C) \ar[rr]|(.49)\hole&& P'(C,C')\\
		P(C,C) \ar[ur]_{\eta_{CC}}\ar[rr]_{P(C,f)}&& P(C,C')\ar[ur]_{\eta_{CC'}}
		}}
	\]
	When all ends exist, sending a functor $P$ into its end $\eend(P)$ is a (covariant) functor $\eend : \Cat(\C^\opp\times \C,\D)\to \D$: the usual argument applies, as the arrow $\eend(\eta)\circ\eend(\eta')$ must coincide with $\eend(\eta\circ\eta')$ in a suitable pasting of cubes. Similarly, the unique arrow induced by $\id_P : P\To P$ must be the identity of $\eend(P)$.
\end{remark}
\subsection{Extranaturality}
A slightly less general, but better behaved\footnote{We say \emph{better behaved} since extranaturality admits a graphical calculus translating commutativity\hyp{}checking into checking that certain string diagrams can be deformed one into the other.} notion of super\hyp{}naturality, that allows again to define co\fshyp{}wedges and thus co\fshyp{}ends, is available: the notion is called \emph{extra\hyp{}naturality} and it was introduced in \cite{eilenberg1966generalization}.
\begin{definition}[Extranatural transformation]\label{extranatural}\index{Transformation!extranatural ---}
	Let $\A,\B,\C,\D$ be categories, and $P,Q$ be functors
	\begin{gather*}
		P : \A\times \B^\opp\times\B\to \D,\\
		Q : \A\times \C^\opp\times\C \to \D.
	\end{gather*}
	An \emph{extranatural transformation} $\alpha : P\din Q$ consist of a collection of arrows
	\[
		\alpha_{ABC} : P(A,B,B) \longrightarrow Q(A, C,C)
	\]
	indexed by triples of object in $\A\times\B\times\C$ such that the following hexagonal diagram commutes for every triple of arrows $f : A\to A'$, $g : B\to B'$, $h : C\to C'$, all taken in their suitable domains:
	\[
		\vcenter{\xymatrix@C=1.6cm{
		P(A,B',B) \ar[r]^{P(f,B',g)}\ar[d]_{P(A,g,B)} & P(A', B', B') \ar[r]^{\alpha_{A'B'C}} & Q(A', C,C) \ar[d]^{Q(A', C,h)}\\
		P(A,B,B) \ar[r]_{\alpha_{ABC'}} & Q(A, C', C') \ar[r]_{Q(f,h,C')} & Q(A', C, C');
		}}
	\]
\end{definition}
Notice how this commutative hexagon can be equivalently described as the juxtaposition of three distinguished commutative squares, depicted in \cite{eilenberg1966generalization}: the three can be obtained letting respectively $f$ and $h$, $f$ and $g$, or $g$ and $h$ be identities in the former diagram, which thus collapses to
\begin{gather}
	\xymatrix@C=1.5cm{
	P(A, B,B) \ar[r]^{P(f,B,B)}\ar[d]_{\alpha_{ABC}}&\ar[d]^{\alpha_{A'BC}} P(A', B,B) \\
	Q(A,C,C) \ar[r]_{Q(f,C,C)}& Q(A', C,C)
	}\quad
	\xymatrix@C=1.5cm{
	P(A,B',B) \ar[r]^{P(A,B',g)}\ar[d]_{P(A,g,B)}&\ar[d]^{\alpha_{AB'C}} P(A, B', B') \\
	P(A,B,B) \ar[r]_{\alpha_{ABC}}& Q(A, C,C)
	}\notag\\
	\vcenter{\xymatrix@C=1.5cm{
	P(A,B,B) \ar[r]^{\alpha_{ABC}}\ar[d]_{\alpha_{ABC'}}&\ar[d]^{Q(A,C,h)} Q(A, C,C) \\
	Q(A, C', C') \ar[r]_{Q(A,h,C')}& Q(A, C, C')
	}}
	\label{extrana}
\end{gather}
\begin{remark}
	We can again define co\fshyp{}wedges in this setting: if $\B=\C$ and in $P(A,B,B)\to Q(A,C,C)$ the functor $P$ is the constant functor on $D\in\D$, and $Q(A,C,C)=\bar Q(C,C)$ is mute in $A$, we get a wedge condition for $D\din Q$; dually we obtain a cowedge condition for $P(B,B)\to Q(A,B,B)\equiv D'$ for all $A,B,C$.

	It's worth to mention that a extranatural transformation contains more information than a dinatural, since in \ref{extranatural} we are given arrows
	\[F(B,B) \xto{\alpha_{BB'}} G(B',B')\]
	that are simultaneously a cowedge in $B$ for each $B'$, and a wedge in $B'$ for all $B,B'\in\B$. We shall see in a while that extranaturality can be obtained as special case of dinaturality.
\end{remark}
Both dinatural and extranatural transformations give rise to the same notion of co\fshyp{}end, defined as a universal co\fshyp{}wedge for a bifunctor $F : \C^\opp\times\C\to \D$. (More formally: the notion of dinatural co\fshyp{}wedge is indistinguishable from the notion of extranatural co\fshyp{}wedge, and thus the two give rise to the same notion of co\fshyp{}end.)

We should prefer extranaturality for a variety of reasons:
\begin{itemize}
	\item it is less general (see \ref{extraisdi}), but it still makes co\fshyp{}ends available;
	\item it gives rise to a fairly intuitive \emph{graphical calculus} (see \ref{gracal}); moreover, it behaves better under composition (see Exercise \ref{ex1:compoextra});
	\item extranaturality is the correct notion in the enriched setting (see \ref{enridina} and the caveat right after).
\end{itemize}
\begin{definition}[\upeyes Graphical calculus for extranaturality]\label{gracal}
	The graphical calculus for extranatural transformations depicts the components $\alpha_{ABC}$, and arrows $f : A\to A'$, $g : B\to B'$, $h : C \to C'$, respectively as planar diagrams like
	\begin{center}
		\begin{tikzpicture}[scale=1.5]
			\begin{scope}[thick]
				\draw[ultra thin] (0.25,0) -- (1.75,0) -- (1.75,1) -- (0.25,1) -- cycle;
				\draw (.5,0) node[below] {$G(A,$} -- (.5,1) node[above] {$F(A$} ;
				\draw (1,0) node[below] {$\phantom{(}C,$} .. controls (1,.5) and (1.5,.5) .. (1.5,0) node[below] {$C)$};
				\draw (1,1) node[above] {$\phantom{(}B,$} .. controls (1,.5) and (1.5,.5) .. (1.5,1) node[above] {$B)$};
			\end{scope}
			\begin{scope}[thick, xshift=3cm]
				\draw[ultra thin] (0.25,0) -- (1.75,0) -- (1.75,1) -- (0.25,1) -- cycle;
				\draw (.5,0) node[below] {$A'$} -- (.5,1) node[above] {$A$};
				\draw[xshift=.5cm] (.5,0)node[below] {$B'$} -- (.5,1) node[above] {$B$};
				\draw[xshift=1cm] (.5,0) node[below] {$C'$}-- (.5,1) node[above] {$C$};
				\filldraw[lightgray!70] (.5,.5) circle (4pt) node[black] {$\scriptstyle f$};
				\filldraw[lightgray!70] (1,.5) circle (4pt) node[black] {$\scriptstyle g$};
				\filldraw[lightgray!70] (1.5,.5) circle (4pt) node[black] {$\scriptstyle h$};
			\end{scope}
		\end{tikzpicture}
	\end{center}
	where wires are labeled by objects and must be thought oriented from top to bottom. The commutative squares of (\ref{extrana}) become, in this representation, the following three string diagrams, whose equivalence is graphically obvious (the labels $f,g,h$ are allowed to `slide' along the wire they live in):
	\begin{center}
		\begin{tikzpicture}[thick, scale=.9]
			\standard
			\begin{scope}[yshift=-1cm]
				\fst{f}
			\end{scope}
			\begin{scope}[xshift=2.5cm]
				\fst{f}
			\end{scope}
			\begin{scope}[xshift=2.5cm, yshift=-1cm]
				\standard
			\end{scope}
			\draw (2.25,0) node {\scalebox{2}{=}};
		\end{tikzpicture}
		\hspace{\fill}
		\begin{tikzpicture}[thick, scale=.9]
			\snd{g}
			\begin{scope}[yshift=-1cm]
				\standard
			\end{scope}
			\begin{scope}[xshift=2.5cm]
				\trd{g}
			\end{scope}
			\begin{scope}[xshift=2.5cm, yshift=-1cm]
				\standard
			\end{scope}
			\draw (2.25,0) node {\scalebox{2}{=}};
		\end{tikzpicture}
		\hspace{\fill}
		\begin{tikzpicture}[thick, scale=.9]
			\standard
			\begin{scope}[yshift=-1cm]
				\snd{h}
			\end{scope}
			\begin{scope}[xshift=2.5cm]
				\standard
			\end{scope}
			\begin{scope}[xshift=2.5cm, yshift=-1cm]
				\trd{h}
			\end{scope}
			\draw (2.25,0) node {\scalebox{2}{=}};
		\end{tikzpicture}
	\end{center}
\end{definition}
\begin{remark}
	The notion of extranatural transformation can be specialised to encompass various other constructions: simple old naturality arises when $F, G$ are both constant in their co\fshyp{}wedge components, so the cap and cup in $\alpha_{ABC}$ vanish:
	\begin{center}\begin{tikzpicture}[thick, scale=1.4]
			\standardbis{black}{dashed,black!40}{dashed,black!40}
		\end{tikzpicture}\end{center}
	The wedge and cowedge conditions arise when either $F,G$ are constant, so that the straight line and one among the cup and the cap in $\alpha_{ABC}$ vanish:
	\begin{center}\begin{tikzpicture}[thick, scale=1.4]
			\standardbis{dashed, black!40}{black}{dashed, black!40}
			\begin{scope}[xshift=4cm]
				\standardbis{dashed, black!40}{dashed, black!40}{black}
			\end{scope}
		\end{tikzpicture}\end{center}
	All the others mixed situations (a wedge\hyp{}cowedge condition, naturality and a wedge, etc. that do not have a specified name) admit a graphical representation of the same sort, and follow similar graphical rules of juxtaposition, when the boundaries of their associated cells agree in shape in the obvious sense.
\end{remark}
All extranatural transformations can be obtained as particular cases of dinatural; on the contrary, there are dinatural transformations which are not extranatural: an example is given in Exercise \ref{ex1:dinatarentextra}.
\begin{proposition} \upeyes \label{extraisdi}
	Extranatural transformations are particular kinds of dinatural transformations.
\end{proposition}
\begin{proof}[Proof (due to T\@.~Trimble)]
	Given functors $F : \C^\opp \times \C \times \C \to \D$ and $G : \C \times \C \times \C^\opp \to \D$, set $\A = \C \times \C^\opp \times \C^\opp$, and form two new functors $F', G': \A^\opp \times \A \to \D$ by taking the composites
	\begin{gather*}
		F' = (\C^\opp \times \C \times \C) \times (\C \times \C^\opp \times \C^\opp) \xto{\text{proj}} \C^\opp \times \C \times \C \xto{F} \D \\
		(X', Y', Z'; X, Y, Z) \longmapsto (X', X, Y') \stackrel{F}{\mapsto} F(X', X, Y')
		\\
		G' = (\C^\opp \times \C \times \C) \times (\C \times \C^\opp \times \C^\opp) \xto{\text{proj}'} \C \times \C \times \C^\opp \xto{G} \D\\
		(X', Y', Z'; X, Y, Z) \mapsto (Y', Z', Z) \stackrel{G}{\mapsto} G(Y', Z', Z)
	\end{gather*}
	Now let's put $A' = (X', Y', Z')$ and $A = (X, Y, Z)$, considered as objects in $\A$. An arrow $\phi : A' \to A$ in $\A$ thus amounts to a triple of arrows $f : X' \to X$, $g : Y \to Y'$, $h : Z \to Z'$ all in $\C$.
	Following the instructions above, we have $F'(A', A) = F(X', X, Y')$ and $G(A', A) = G(Y', Z', Z)$.
	Now if we write down a dinaturality hexagon for $\alpha: F' \din G'$, we get a diagram of shape
	\[
		\vcenter{\xymatrix{
		F'(A, A') \ar[d]_{F(\phi, 1)}\ar[r]^{F'(1, \phi)} & F'(A, A) \ar[r]^{\alpha_A} & G'(A, A) \ar[d]^{G'(\phi, 1)} \\
		F'(A', A') \ar[r]_{\alpha_{A'}} & G'(A', A') \ar[r]_{G'(1, \phi)} & G(A', A)
		}}\]
	which translates to a hexagon of shape
	\[
		\vcenter{\xymatrix{
				F(X, X', Y) \ar[d]_{F(f, 1, g)} \ar[r]^{F(1, f, 1)} & F(X, X, Y) \ar[r] & G(Y, Z, Z) \ar[d]^{G(g, h, 1)}\\
				F(X', X', Y') \ar[r] & G(Y', Z', Z') \ar[r]_{G(1, h, 1)} & G(Y', Z', Z)
			}}
	\]
	where the unlabeled arrows are the extranatural components. This is the extranaturality hexagon of \ref{extranatural}.
\end{proof}
\subsection{The integral notation for co\fshyp{}ends}\label{the_integral_not}
\index{Yoneda!---'s notation for co\fshyp{}ends}
\index{_aaa_int@$\int$}
A suggestive notation for co\fshyp{}ends, alternative to the eponymous one `$\underline{\text{co}}/\eend(F)$', is due to N\@.~Yoneda, which in \cite{Yoneda} introduces most of the notions we are dealing with, in the setting of $\Ab$\hyp{}enriched functors $\C^\opp\times \C\to \Ab$:
\begin{notation}
	The \emph{integral notation} denotes the end of a functor $F\in \Cat(\C^\opp\times \C,\D)$ as a `subscripted\hyp{}integral' $\int_CF(C,C)$, and the coend $\coend(F)$ as the `superscripted\hyp{}integral' $\int^CF(C,C)$.
\end{notation}
From now on we will systematically adopt this notation to denote the universal co\fshyp{}wedge $\underline{\text{co}}/\eend(F)$ or, following a well\hyp{}established abuse of notation, the object itself; when the domain of $F$ has to be made explicit, we will also employ more pedantic variants of $\int_C F$ and $\int^C F$ like
\[
	\int_{C\in \C} F(C,C), \qquad \int^{C\in \C} F(C,C).
\]
\begin{remark}
	In reading \cite{Yoneda}, one should be aware that Yoneda employs a reversed notation to denote ends and coends: he calls \emph{integration} what we call a coend, which he denotes as $\int_{C\in\C} F(C,C)$, \ie in the way we denote an end; and he calls \emph{cointegrations} our ends, which he denotes $\int_{C\in\C}^*F(C,C)$.

	No trace of this ambiguity survived in the current literature, so we will not mention the Yoneda convention ever again.
\end{remark}
\index{Co/end!functoriality of ---} Functoriality of co\fshyp{}ends acquires a particularly suggestive flavor when written in integral notation: the dream of every freshman learning calculus is that the integral of a product of functions is just the product of the integrals of the two functions. In category theory this is true, provided the integral of a function is the map induced between two co\fshyp{}ends, and that product is composition of arrows.
\begin{notation}
	The unique arrow $\eend(\eta)$ induced by a natural transformation $\eta : F\To G$ between $F,G\in \Cat(\C^\opp \times \C ,\D)$ can be written as $\int_C\eta : \int_C F\to \int_C G$, and uniqueness of this induced arrow entails functoriality, \ie $\int_C (\eta\circ\sigma)=\int_C \eta\circ\int_C \sigma$ and $\int_C \id_F=\id_{\int F}$.

	A similar convention holds for coends.
\end{notation}
\begin{remark}
	As \cite[IX.5]{McL} puts it,
	\begin{quote}
		[\dots] the `variable of integration' $C$ [in $\int_C F$] appears twice under the integral sign (once contravariant, once covariant) and is `bound' by the integral sign, in that the result no longer depends on $C$ and so is unchanged if $C$ is replaced by any other letter standing for an object of the category $\C$.
	\end{quote}
	This somehow motivates the integral notation for co\fshyp{}ends, and yet the analogy between integral calculus and co\fshyp{}ends seems to be too elusive to justify.

	There seems to be no chance to give a formal explanation of the similarities between integrals and coends, but it is nevertheless very suggestive to employ an informal justification for such an analogy to exist. See Exercises \ref{ex1:vector-of-coends} or \ref{dirac}, or even \ref{stokkio}.
\end{remark}
\section{Co/ends as co\fshyp{}limits}\label{coends_as_colims}\index{Co/end!---s as colimits}
A general tenet of elementary category theory is that universal objects (\ie objects having the property of being initial or terminal in some category) can be equivalently characterized
\begin{enumtag}{u}
	\item as limits (so the existence and uniqueness is simply translated in a category of \emph{diagrams}),
	\item as adjoints (so the uniqueness follows from uniqueness of adjoints, see \ref{adjbasta}),
	\item as the representing object of a certain functor (so the uniqueness follows from the fact that the Yoneda embedding $\yon_\C$ is fully faithful, \ref{yon_is_ff}).
\end{enumtag}
The language of co\fshyp{}ends makes no exception: the scope of the following subsection is to characterise the co\fshyp{}end of a functor $F : \C^\opp\times\C\to \D$ as a co\fshyp{}limit over a suitable diagram $\bar F$ (obtained from a canonically chosen correspondence $F\mapsto \bar F$), and (see \ref{lims_iff_pis_and_eqs}) consequently as the co\fshyp{}equaliser of a single pair of arrows.
\begin{remark}
	Given $F : \C^\opp\times \C\to \D$ and a wedge $\tau : D\din F$, we can build the following commutative diagram
	\[
		\vcenter{\xymatrix@R=.4cm@C=.4cm{
		& F(C,C) \ar[rr]^{F(C, f)} & & F(C,C ')\ar[dd]^{F(C, g)}\\
		D\ar[ur]^{\tau_C}\ar@/_1pc/[rrrd]|(.67)\hole_{\tau_{g\circ f}}\ar[rr]^{\tau_{C'}}\ar[dd]_{\tau_{C''}} && F(C', C')\ar[dd]\ar[ur]_{F(f, C')} \\
		& && F(C, C'')\\
		F(C'', C'') \ar[rr]_{F(g, C'')} && F(C', C'')\ar[ur]_{F(f, C'')}
		}}
	\]
	where
	$C\xto{f}C'\xto{g}C''$ are morphisms in $\C$. From this commutativity we deduce the following relations:
	\begin{align*}
		\tau_{g\circ f} & =F(g\circ f,C'')\circ \tau_{C''} = F(C,g\circ f)\circ \tau_C \\
		                & =F(f,C'')\circ F(g,C'')\circ \tau_{C''}=F(f,C'')\circ \tau_g \\
		                & =F(C,g)\circ F(C,f)\circ\tau_C = F(C,g)\circ \tau_f.
	\end{align*}
	where $\tau_f$, $\tau_g$ are the common values $F(f, C')\circ\tau_{C'} = F(C,f)\circ\tau_C$ and $F(C',g)\circ\tau_{C'} = F(g, C'')\circ\tau_{C''}$ respectively, and $\tau_{g\circ f}$ is the common value $F(C,g)\circ\tau_f = F(f, C'') \circ\tau_g$.
\end{remark}
These relations imply that there is a link between co\fshyp{}wedges and co\fshyp{}cones, encoded in the following definition.
\begin{definition}[\righteyes The twisted arrow category of $\C$]\label{twisted}\index{Category!twisted arrows ---}\index{Twisted arrow category}
	For every category $\C$ we define $\tw(\C)$, the category of \emph{twisted arrows} in $\C$ as follows:
	\begin{itemize}
		\item $(\tw(\C))_o =\hom(\C)$ (of course, this will not be $\mho$\hyp{}small if $\C$ was only $\mho^+$\hyp{}small);
		\item Given $f : A \to A'$, $g : B\to B'$ a morphism $f\to g$ is given by a pair of arrows $(h : B\to A,k : A'\to B')$, such that the square
		      \[\vcenter{\xymatrix{
					      A \ar[d]_f & \ar[l]_h B\ar[d]^g \\
					      A' \ar[r]_k & B'
				      }}\]
		      commutes (asking that the arrow between domains is reversed is \emph{not} a mistake), \ie that $g=k\circ f\circ h$.
	\end{itemize}
\end{definition}
Endowed with the obvious rules for composition and identity, $\tw(\C)$ is easily seen to be a category, and now we can find a functor
\[
	\vcenter{\xymatrix{
			\Cat(\C^\opp\times\C,\D)\ar[rr]&&\Cat(\tw(\C),\D)
		}}
\]
defined sending $F : \C^\opp\times\C \to \D$ to the functor $\bar F : \tw(\C)\to \D : \var{C}{C'}\mapsto F(C,C')$; it is extremely easy now to check that bifunctoriality for $F$ corresponds to functoriality for $\bar F$, but there is more to this remark.
\begin{remark}[Co/ends are co\fshyp{}limits, I]\label{is.a.colim}\righteyes
	The family of arrows
	\[ \{\tau_f \mid f\in\hom(\C)\} \]
	constructed above is a cone for the functor $\bar F$, and conversely any such cone determines a wedge for $F$, obtained setting $\{\tau_C = \tau_{\id_c}\}_{C\in \C}$.

	A morphism of cones maps to a morphism between the corresponding wedges, and conversely every morphism between wedges induces a morphism between the corresponding cones; these operations are mutually inverse and form an equivalence between the category $\Cn(\bar F)$ of cones for $\bar F$ and the category $\wed(F)$ of wedges for $F$. (We leave this to the reader to check and to properly dualise.)
\end{remark}
Equivalences of categories obviously preserve initial and terminal objects, thus we have isomorphisms\footnote{Notice that the colimit is taken over the category $\tw^\opp(\C)$, the \emph{opposite} of $\tw(\C^\opp)$: an \emph{object} of $\tw^\opp(\C)$ is an arrow $f : C'\to C$ in $\C^\opp$, and a \emph{morphism} from $f : C \to C'$ to $g : D \to D'$ is a commutative square $(u,v)$ such that $vgu = f$.}
\[
	\int_C F(C,C) \cong {\textstyle \lim_{\tw(\C)}} \bar F; \qquad
	\int^C F(C,C) \cong \colim_{\tw(\C^\opp)^\opp} \bar F
\]
\begin{remark}\label{endsareeq}\righteyes
	According to \ref{lims_iff_pis_and_eqs}, co\fshyp{}limits in a category exist as soon as it has co\fshyp{}products and co\fshyp{}equalisers. So we would expect a characterisation of co\fshyp{}ends in terms of these simpler pieces as well; such a characterisation exists, and it turns out to be extremely useful in explicit computations.

	It is rather easy to extract from the bare universal property that there must be an isomorphism
	\[\label{ends_are_lims}
		\vcenter{\xymatrix{
				\int_C F(C,C) \cong \text{eq}\left( \displaystyle \prod_{C\in\C} F(C,C) \ar@<3pt>[r]^(.6){F^*}\ar@<-3pt>[r]_(.6){F_*}\right.& \left.\displaystyle \prod_{\varphi : C \to C'} F(C,C')\right)
			}}
	\]
	where the product over morphisms $ C \to C'$ can be expressed as a double product (over the objects $C,C'\in \C$, and over the arrows $f$ between these two fixed objects), and the arrows $F^*, F_*$ are easily obtained from the arrows whose $(f; C,C')$\hyp{}components are (respectively) $F(f, C')$ and $F(C, f)$. This is a consequence of the fact that an `element' in $\prod_C F(C,C)$, regarded as a family $(x_C\mid C\in\C)$, shall equalise both actions of $F$ on arrows at the same time in order to belong to the end $\int_C F(C,C)$.
	The dual statement of \ref{endsareeq}, expressing
	\[\label{coends_are_coeqs}
		\vcenter{\xymatrix{
		\int^C F(C,C) \cong \text{coeq}\left( \displaystyle \coprod_{C\in\C} F(C,C) \ar@{<-}@<3pt>[r]^-(.6){F^*}\ar@{<-}@<-3pt>[r]_-(.6){F_*}\right.& \left.\displaystyle \coprod_{\varphi : C \to C'} F(C',C)\right)
		}}
	\]
	is left as an exercise for the reader to formalise, in \ref{ex1:iscoeq}.
\end{remark}

The following remark is elementary but extremely useful: it asserts that the co\fshyp{}limit of a functor has the same universal property of the co\fshyp{}end of the same functor, when it is `promoted' as mute in its remaining variable.
\begin{remark}\label{coends_are_colims_are_coends}
	Let $F : \C \to \D$ be a functor; we can always regard it as a functor $F' : \C^\opp\times \C\to \D$ mute in its first variable: this means that we can extend the action of $F$ defining a new functor $F'$ such that $F'(C,C')=FC'$ for every $C,C'\in\C$, and $F'(C,f)=Ff$, $F'(f,C')=\id_{C'}$ for every $f : C\to C'$; from this, and from \ref{is.a.colim} above, since all mute functors can be regarded as arising this way, it follows that the co\fshyp{}end of a functor that is mute in one of its variables coincides with its co\fshyp{}limit.
\end{remark}
\begin{definition}
	There is an obvious definition of \emph{preservation} of co\fshyp{}ends from their description as co\fshyp{}limits, which reduces to the preservation of the particular kind of co\fshyp{}limit involved in the definition of $\eend(T)$ and $\coend(T)$: let $F : \D \to \E$ be a functor, and let $T : \C^\opp\times \C \to \D$ be a functor; we say that
	\begin{enumtag}{pr}
		\item $F$ \emph{preserves} the end of $T$ if the family of maps
		\[\xymatrix{F\Big(\int_CT(C,C)\Big) \ar[r]^-{F\omega_C} & FT(C,C)}\]
		exhibits the universal property of the end of the composed functor $\C^\opp\times \C \xto{T} \D \xto{F} \E$;
		\item $F$ \emph{reflects} the end of $T$ if, whenever $F\omega_C$ is the end of the composition $F\circ T$, then $\omega_C$ exhibits the end of $F$.
	\end{enumtag}
\end{definition}
A dual definition defines the concept of preserving and reflecting coends; since ends are limits, it is clear that a functor that preserves all limits preserves all ends, and since limits are ends of mute functors, also the converse is true; dually for coends.

An alternative, equivalent way to put this result is the following:
\begin{theorem}\label{coconti}
	Every co\fshyp{}continuous functor $F : \D\to\E$ preserves the co\fshyp{}ends that exist in $\D$:
	\begin{itemize}
		\item if $T : \C^\opp\times\C\to \D$ has an end $\int_C T(C,C)$, and $F : \D\to\E$ commutes with all limits, then
		      \[
			      F\Big(\textstyle \int_C T(C,C) \Big) \cong \int_C FT(C,C)
		      \]
		      meaning that the image of the terminal wedge of $T$ under $F$ is a terminal wedge for the composite functor $F\circ T : \C^\opp\times\C \to \E$, and thus the two terminal objects are canonically isomorphic;
		\item Dually, if $T : \C^\opp\times\C\to \D$ has a coend $\int^C T(C,C)$, and $F : \D\to\E$ commutes with all colimits, then
		      \[
			      F\Big(\textstyle \int^C T(C,C) \Big) \cong \int^C FT(C,C)
		      \]
		      meaning that the image of the initial wedge of $T$ under $F$ is a initial cowedge for the composite functor $F\circ T : \C^\opp\times\C \to \E$, and thus the two initial objects are canonically isomorphic.	\end{itemize}
\end{theorem}
Similarly, a functor that reflects all co\fshyp{}limits reflects all co\fshyp{}ends.

As a capital example of a functor that preserves all ends, we have the hom functors:
\begin{corollary}[The $\hom$ functor commutes with integrals]\label{commuhom}
	From the fact that the $\hom$ bifunctor $\C(\firstblank,\secondblank) : \C^\opp\times\C\to\Set$ is such that
	\begin{align}
		\C\big(\colim F, C\big) & \cong \lim \C\big( F, C\big) \notag \\
		\C\big(C, \lim F \big)  & \cong \lim \C\big(C, F\big)
	\end{align}
	and from \eqref{ends_are_lims}, \eqref{coends_are_coeqs} we deduce that for every $D\in \D$ and every functor $F : \C^\opp\times\C \to \D$ we have canonical isomorphisms
	\begin{align}
		\D\Big( \int^C F(C,C),D \Big) & \cong \int_C \D(F(C,C),D)\notag \\
		\D\Big(D, \int_C F(C,C)\Big)  & \cong \int_C\D(D,F(C,C))
	\end{align}
\end{corollary}
We close the section recording a notable but somewhat technical result.
\begin{remark}[\protect{\cite[IX.5.1]{McL}} Co/ends are co\fshyp{}limits, II]\upeyes
	Define the \emph{subdivision graph} $(\wp \C)^\S$ of a category $\C$ as the directed graph having a vertex $C^\S$ for each object $C\in \C$, and a vertex $f^\S$ for each morphism $f : C\to C'$ in $\C$, and edges all the arrows $C^\S \to f^\S$ and $C'^\S \to f^\S$, as $f : C \to C'$ runs over morphisms of $\C$.

	The \emph{subdivision category} $\C^\S$ is obtained from $(\wp \C)^\S$ formally adding identities and giving to the resulting category the trivial composition law (composition is defined only if one of the arrows is the identity).
\end{remark}
Every functor $F : \C^\opp\times \C\to \D$ induces a functor $F^\S : \C^\S\to \D$, whose limit (provided it exists) is isomorphic to the end of $F$. More precisely,
\begin{proposition}\upeyes In the above notation,
	\begin{itemize}
		\item There is a final functor $\varsigma$ from $\C^\S$ to $\tw(\C)$; thus, the limit of $\bar F$ computed over $\tw(\C)$ is the same as the limit of $\bar F \circ \varsigma$ computed over $\C^\S$;
		\item There is a functor $F^\S : \C^\S \to \D$ such that wedges for $F$ correspond bijectively to cones for $F^\S$;
	\end{itemize}
	the assignment $F\mapsto F^\S$ sets up an isomorphism of categories between wedges for $F$ and ones for $F^\S$, and thus the terminal wedge, \ie the end of $F$, must correspond to the terminal cone, the limit of $F^\S$.
	Dually,
	\begin{itemize}
		\item There is a cofinal functor $\vartheta$ from $\C^\S$ to $\tw(\C^\opp)^\opp$; thus, the colimit of $\bar F$ done over $\tw(\C^\opp)^\opp$ is the same as the colimit of $\bar F \circ \vartheta$ done over $\C^\S$;
		\item There is a functor $F^\S : \C^\S \to \D$ such that wedges for $F$ correspond bijectively to cocones for $F^\S$;
	\end{itemize}
	this assignment sets up an equivalence of categories between cowedges for $F$ and ones for $F^\S$, and thus the initial cowedge, \ie the coend of $F$, must correspond to the initial cocone, the colimit of $F^\S$.
\end{proposition}
We just address the case of ends, leaving dualisation process to the reader. Recall that a functor $F : \C \to \D$ is final if every comma category $(C\downarrow F)$ is nonempty and connected.

To sum up, given a functor $F : \C^\opp\times\C \to \D$ one can equivalently compute the value of the end $\int_C F$ as
\begin{itemize}
	\item The terminal wedge $\alpha : \int_C F(C,C) \din F$;
	\item The terminal cone $\alpha^\S : \lim_{\C^\S}F^\S \To F^\S$;
	\item The terminal cone $\bar \alpha : \lim_{\tw(\C)}\bar F \to \bar F$.
\end{itemize}
\begin{proof}
	Let us first define the functor $\varsigma$: on objects,
	\[ \varsigma(C^\S) = \id_C \qquad\qquad \varsigma(f^\S) = f \]
	while on morphisms $f :C \to C'$
	\[ \varsigma\var{C^\S}{f^\S} =
		\begin{smallmatrix}
			C &=& C \\ \rotatebox[origin=c]{90}{$\scriptstyle =$} & & \downarrow \\ C' & \to & C'
		\end{smallmatrix}
		\qquad\qquad
		\varsigma\var{C'^\S}{f^\S} =
		\begin{smallmatrix}
			C' & \leftarrow & C \\ \rotatebox[origin=c]{90}{$\scriptstyle =$} && \downarrow \\ C' &=& C'
		\end{smallmatrix} \]
	This easily shows how $\varsigma$ is surjective on objects, and thus each comma category $(f\downarrow\varsigma)$ is nonempty. The image of $\varsigma$ visibly contains very few morphisms; yet it is still final, because the comma category $(f \downarrow \varsigma)$ is easily seen to be connected.

	To prove the second point, we shall show that every functor $F : \C^\opp\times \C \to \D$ induces a functor $F^\S : \C^\S \to \D$ from the subdivision category; we can define
	\begin{itemize}
		\item an object function posing $F^\S(C^\S) = F(C,C)$ and $F^\S\var{C}{C'}^\S = F(C,C')$;
		\item a morphism function, posing $F^\S\var{C^\S}{f^\S} = F(C, f)$ and $F^\S\var{C'^\S}{f^\S}=F(f,C')$.
	\end{itemize}
	Let now $\alpha : D \To F^\S $ be a cone for $F^\S$; this means that for every morphism of $\C^\S$ one among the diagrams
	\[
		\vcenter{\xymatrix{
		& D\ar[dr]^{\alpha_{CC'}}\ar[dl]_{\alpha_{C}} &  & & D\ar[dr]^{\alpha_{CC'}}\ar[dl]_{\alpha_{C'}} & \\
		F(C,C)\ar[rr]_{F(1,f)}&&F(C,C') & F(C',C') \ar[rr]_{F(f,1)} && F(C,C')
		}}
	\]
	(chosen accordingly to the shape of $f^\S$) is commutative. But then, the restriction of $\alpha$ to its diagonal component forms a wedge for $F$ with domain $D$, because the square
	\[\label{ueggio}
		\vcenter{\xymatrix{
		D \ar[r]^{\alpha_C}\ar[d]_{\alpha_{C'}}& F(C,C)\ar[d]^{F(1,f)}\\
		F(C',C') \ar[r]_{F(f,1)} & F(C, C')
		}}
	\]
	has $\alpha_{CC'}$ as diagonal. Vice versa, given a wedge $\alpha : D \din F$ the square \eqref{ueggio} defines the components of a cone $\alpha : D \To F^\S$ by $\alpha_{CC'} = F(1,f)\circ \alpha_{CC} = F(f,1)\circ \alpha_{C'C'}$. It is evident that this sets up a bijection between $\wed(F)$ and $\Cn(F^\S)$, and that this can in fact be promoted to an isomorphism between the two categories.
\end{proof}
From this it ultimately follows that both characterisations of co\fshyp{}ends as co\fshyp{}limits, either as diagrams with domain $\C^\S$, or with domain $\tw(\C)$, lead to the same theory.
\section{The Fubini rule}
\index{Fubini rule!--- for coends}
An absolutely central theorem for coend calculus is the `exchange rule' for integrals known as the \emph{Fubini rule}: informally, it says that the result of an integration on more than one variable does not depend on the order in which the operation is performed: in mathematical analysis, if $f : X \times Y \to \bR$ is a function such that $\int_{X\times Y} |f|d\mu_{X\times Y}$ exists finite, the three integrals
\[
	\int_X\left(\int_Y f(x,y)\,dy\right)\,dx=\int_Y\left(\int_X f(x,y)\,dx\right)\,dy=\int_{X\times Y} f(x,y)\,d(x,y).
\]
are equal. In category theory, if a functor $F$ defined on $\C^\opp\times \C\times \E^\opp\times \E$ admits a co\fshyp{}end, then so do the two functors obtained integrating over $\C$ first and over $\E$ second, or in the opposite order.
\begin{theorem}[\righteyes Fubini theorem for co\fshyp{}ends]\label{fubozzo}
	Given a functor
	\[
		F : \C^\opp\times\C\times \E^\opp\times \cate E\to \D,
	\]
	we can form the end $\int_CF(C,C,\firstblank,\secondblank)$ obtaining a functor $\E^\opp\times \cate E\to \D$ whose end is $\int_E\int_CF(C,C,E,E)\in \D$; we can also form the ends $\int_C\int_E F(C,C,E,E)\in \D$ and $\int_{(C,E)}F(C,C,E,E)$ identifying $\C^\opp\times\C\times \E^\opp\times \cate E$ with $(\C\times\cate E)^\opp\times (\C\times\cate E)$.

	Then, there are canonical isomorphisms between the three objects:
	\[
		\int_{(C,E)}F(C,C,E,E)\cong
		\int_E\int_CF(C,C,E,E)\cong
		\int_C\int_E F(C,C,E,E)
	\]
	Dually, there are canonical isomorphisms between the iterated coends
	\[\label{fub_for_coends}
		\int^{(C,E)}F(C,C,E,E)\cong
		\int^E\kern-.5em\int^CF(C,C,E,E)\cong
		\int^C\kern-.5em\int^E F(C,C,E,E)
	\]
\end{theorem}
As usual, this is an existence and uniqueness result: one of the three objects above exists if and only if so do the other two, and there are canonical isomorphisms between them.

We shall only prove the statement for ends, \eqref{fub_for_coends} for coends being the exact dual statement.\footnote{The reader will have enough care to dualise properly the proof that follows: the coend functor $\int^C :\Cat(\C^\opp\times\C, \D)\to\D $ will turn out to be a left adjoint, with right adjoint $H^\C$, etc.}

One possible strategy to prove Fubini theorem would be to find suitable canonical map between the end $\int_{(C,E)}F(C,C,E,E)$ and (say) the end $\int_E\int_CF(C,C,E,E)$ (the same argument will then produce a canonical map between the end $\int_{(C,E)}F(C,C,E,E)$ and the end $\int_E\int_CF(C,C,E,E)$), and then prove its invertibility. This is of course a viable option, but has little conceptual content, and can result in a long, unenlightening proof.

Instead, we prefer to offer a more elegant argument, especially because the strategy of our proof will be recycled in \ref{fubi_for_infty} when we prove the Fubini rule for $\infty$\hyp{}coends (the analogue of coends in $(\infty,1)$\hyp{}category theory, and more specifically in the setting of Joyal-Lurie quasicategories \cite{Joyal, HTT}).

Our strategy goes as follows: we assume that all ends exist. Thus, as stated in \ref{fun_4_ends}, the correspondence $T\mapsto\int_C T$ is a functor. We shall show that $\int_C$ has a left adjoint $H_\C : \D \to \Cat(\C^\opp\times\C, \D)$, and that $H_{\C\times \E}\cong H_\C \circ H_\E\cong H_\E\circ H_\C$; the Fubini rule then follows from the uniqueness of adjoints (see \ref{adjbasta}) because every such `interchange isomorphism' between the left adjoints must induce a similar isomorphism between the right adjoints.

Note that this argument yields for free that the isomorphisms witnessing the Fubini rule above are natural in $F$, \ie that $\int_{(C,E)}, \int_C\int_E, \int_E\int_C$ are isomorphic \emph{as functors}.

Let us then define the functor $H_\C :  \D \to \Cat(\C^\opp\times\C, \D)$ as follows: to the object $D\in\D$ we associate the functor $\hom_\C\otimes D : (C,C')\mapsto \C(C,C')\otimes D$, where for a set $X$ and an object $D\in\D$ we denote $X\otimes D$ the coproduct $\coprod_{x\in X}D$ of $X$ copies of $D$, also called the \emph{copower} or \emph{tensor} of $D$ by $X$. See \ref{tenscotens} for the whole definition; note that this yields a chain of isomorphisms
\index{_aaa_pitch@$\pitchfork$}
\[\label{local_thc}\D(X\otimes D,D')\cong \D(D, X\pitchfork D')\cong \Set(X, \D(D,D')),\]
where $X\pitchfork D$ is defined dually as $\prod_{x\in X}D$. We will freely employ \eqref{local_thc} when needed.

We shall prove that there is a bijection
\[\label{deadj}\textstyle \Cat(\C^\opp\times\C,\D)(\hom_\C\otimes D, F)\cong \D(D, \int_CF)\]
If $\alpha : \hom_\C\otimes D\To F$ is a natural transformation it has components
\[ \alpha_{CC'} : \C(C,C')\otimes D \to F(C,C'); \notag\]
by \eqref{local_thc}, these maps mate to
\[\tilde{\alpha}_{CC'} : D \to \C(C,C')\pitchfork F(C,C')\notag\]
It is easy to show that these mates form a wedge \emph{in the pair $(C,C')$}, and thus there is an induced morphism
\[\label{dat}\hat\alpha : D \to \int_{(C,C')} \C(C,C')\pitchfork F(C,C').\]
\begin{lemma}\label{usefullemma}
	The latter integral $\int_{(C,C')} \C(C,C')\pitchfork F(C,C')$ is in fact isomorphic to $\int_C F(C,C)$.
\end{lemma}
\begin{proof}[Proof of \ref{usefullemma}.]
	We shall first find a candidate wedge
	\[
		\omega_{(CC')} : \int_C F(C,C) \to \C(C,C')\pitchfork F(C,C');
	\]
	for the sake of exposition, we treat the case $\D=\Set$; the core of the argument is the same in the general case, using generalised elements.

	Thanks to \eqref{ends_are_lims} we now know that an `element' of $\int_C F(C,C)$ is a coherent sequence $(a_C \mid C\in \C)$ in the product $\prod_C F(C,C)$, such that $F(u, \id_C)(a_{C'}) = F(\id_{C'}, u)(a_C)$ for every $u : C \to C'$ in $\C$. This means that we can map a coherent sequence $(a_C\mid C\in\C)$ into (say) a sequence $(F(u,C')(a_{C'}) \mid u : C \to C')$, and this sets up a map
	\[\int_C F(C,C) \xto{\varpi_{(CC')}} \C(C,C')\pitchfork F(C,C') = \prod_{C,C'}\prod_{u : C\to C'} F(C,C')\]
	This in turn defines a wedge in the pair $(C,C')$, thus inducing a unique morphism
	\[\xymatrix{
		\int_C F(C,C) \ar[r]^-{\bar\varpi} & \int_{(C,C')} \C(C,C')\pitchfork F(C,C')
		}\]
	between the ends. Now that we put all notation in place, we leave as an exercise for the reader to show that $\varpi_{(CC')}$ is a terminal wedge, and thus that $\bar\varpi$ is in fact an isomorphism.
\end{proof}
\begin{proof}[Proof of \ref{fubozzo}.]
	A natural transformation $\alpha : \hom_\C\otimes D\To F$ induces an arrow like in \eqref{dat}, and this induces a unique arrow $\bar\varpi^{-1} \circ \hat\alpha : D \to \int_C F(C,C)$. The correspondence $\alpha\mapsto \bar\varpi^{-1} \circ \hat\alpha$ sets up the desired adjunction like in \eqref{deadj}.

	Unwinding its definition, we easily see that the functor $H_\C : \D \to \Cat(\C^\opp\times\C, \D)$ has the property  that $H_{\C\times \E}\cong H_\C \circ H_\E\cong H_\E\circ H_\C$ (it easily follows from a commutativity of the involved coproducts: we invite the reader to fill in the details). This, in turn, shows the Fubini rule, since every isomorphism between left adjoint functors $L \To L'$ induces an isomorphism of right adjoint functors $R'\To R$.
\end{proof}
\begin{remark}
	From the binary case shown above an easy induction shows that the iterated co\fshyp{}end of a functor $F : \prod_{i=1}^n \C_i^\opp\times \C_i \to \D$ gives the same result with respect to their integration variables, taken in whatever order. More formally, if $\sigma : \{1,\dots,n\}\to \{1,\dots,n\}$ is any permutation of an $n$ elements set, then there is a canonical isomorphism
	\[
		\int_{C_{\sigma 1}}\kern-1em\dots\int_{C_{\sigma n}} F(C_1, C_1,\dots, C_n, C_n) \cong \int_{(C_1,\dots, C_n)} F(C_1, C_1,\dots, C_n, C_n)
	\] and similarly for the coend of $F$.
\end{remark}
\section{First instances of co\fshyp{}ends}
\index{Co/end!natural transformations as ---s}\index{Natural transformation!---s as ends}
A basic example exploiting the whole machinery introduced so far is the proof that the set of natural transformations between two functors $F,G : \C\to \D$ can be characterised as an end:
\begin{theorem}\label{naturalu}
	Given functors $F,G : \C\to \D$ whose domain is a small category, and whose codomain is locally small, we have the canonical isomorphism of sets
	\[
		\Cat(\C,\D)(F,G)\cong \int_C \D(FC,GC).
	\]
\end{theorem}
\begin{proof}
	A wedge $\tau_C : Y\to \D(FC,GC)$ consists of a function $y\mapsto (\tau_{C,y} : FC \to GC\mid C\in\C)$, which is natural in $C\in \C$ (this is simply a rephrasing of the wedge condition): the equation
	\[
		G(f)\circ \tau_{C,y} = \tau_{C',y}\circ F(f)
	\]
	valid for any $f : C \to C'$, means that for a fixed $y\in Y$ the arrows $\tau_C$ form the components of a natural transformation $F \To G$; thus, 	 there exists a unique way to close the diagram
	\[
		\vcenter{\xymatrix{
		Y \ar[r]^{\tau_C}\ar@{.>}@/_1pc/[dr]_h & \D(FC,GC) \\
		& \Cat(\C,\D)(F,G) \ar[u]
		}}
	\]
	with a function sending $y\mapsto \tau_{\firstblank,y}\in \prod_{C\in\C} \D(FC,GC)$, and where $\Cat(\C,\D)(F,G)\to\D(FC,GC)$ is the wedge sending a natural transformation to its $c$\hyp{}component; the diagram commutes for a single $h : Y\to \Cat(\C,\D)(F,G)$, and this is precisely the desired universal property for $\Cat(\C,\D)(F,G)$ to be $\int_C\D(FC,GC)$.
\end{proof}
\begin{remark}\label{keryoneda}
	A suggestive way to express naturality as a `closure' condition is given in \cite[4.1.1]{Yoneda}, where for an $\Ab$\hyp{}enriched functor (see \ref{enrifun}) $F : \C^\opp\times \C\to \Ab$ from a complete $\Ab$\hyp{}category $\C^\opp\times\C$, one can prove that natural transformations $F \To G$ form the kernel $\ker\delta$ of a map
	\[
		\delta : \bigoplus_{Y\in \C} \Ab(FY, GY) \to \bigoplus_{f : X \to Y} \Ab(FX, GY).
	\]
\end{remark}
\begin{remark}
	If we let $F=G=\id_\C$ in the isomorphism above, we get that the end of the hom functor $\hom : \C^\opp\times\C\to\Set$ is the monoid $M(\C)$ of endo-natural transformations of the identity functor $\id_\C$. This monoid is of great importance in homological algebra and algebraic geometry, as it constitutes a precious source of information on $\C$ (morally, $M(\C)$ can be thought of as the 0\hyp{}th term of a sequence of cohomology groups $H^n(\C)$ associated to $\C$). Similarly, the presence of a `nice' functor $F : \C \to \Set$ (for example, a faithful and conservative one) yields a fairly rich monoid of endomorphisms $M(F) = \int_C \Set(FC,FC)$.

	`Reconstruction theory' is the branch of category theory that asserts that under suitable conditions on $F$, its domain is a category of $M(F)$\hyp{}modules, \ie there is an equivalence (or a full embedding) between $\C$ and $\Mod(M(F))$.
\end{remark}
The set of dinatural transformations between two functors can also be characterised as an end, where the hom functor has been ``completely symmetrised''. The result was first proved in \cite{Dubuc_1970}.
\begin{example}[Dinatural transformations as an end]
	Let $F,G : \C^\opp\times\C \to \D$ be two functors; define a new functor \[\D_d(F,G)[\firstblank,\secondblank] : \C^\opp\times \C\to \Set\]
	\begin{itemize}
		\item on objects, sending the pair $(A,B)$ to the set $\D(F(B,A),G(A,B))$;
		\item on morphisms, sending a pair of arrows $
			      \left[
				      \begin{smallmatrix}
					      A & X\\
					      f\downarrow & \downarrow g \\
					      B & Y
				      \end{smallmatrix}
				      \right]
		      $ to the diagonal of the commutative square
		      \[\vcenter{\xymatrix{
					      \D_d(F,G)[B,X]\ar[r]\ar[d] & \D_d(F,G)[B,Y] \ar[d]\\
					      \D_d(F,G)[A,X] \ar[r]& \D_d(F,G)[A,Y]
				      }}\]
		      whose horizontal and vertical arrows are defined by the action of $F,G$ on morphisms: for example the left vertical arrow is defined as
		      \begin{align}
			      \D(F(X,B),G(B,X)) & \to \D(F(X,A), G(A,X))\notag        \\
			      u                 & \mapsto G(f,1)\circ u \circ F(1,f).
		      \end{align}
	\end{itemize}
	It is a natural question to wonder what is the end of $\D_d(F,G)$: \cite{Dubuc_1970} was the first to observe that such end is in fact isomorphic to the set of dinatural transformations $\alpha : F \din G$.

	In order to prove this, observe that if we consider the parallel morphisms
	\[\prod_{C\in\C}\D(F(C,C),G(C,C)) \rightrightarrows \prod_{f : C \to C'} \D(F(C',C),G(C,C'))\]
	induced by the conjoint action of $\D_d(F,G)$ on morphism as in \eqref{ends_are_lims}, then  $u,v$ are defined sending a sequence $(x_C \mid C \in \C)$ respectively in
	\[\begin{cases}
			(x_C\mid C\in\C) & \mapsto (G(1,f)\circ x_A\circ F(f,1)\mid f : A \to B) \\
			(x_C\mid C\in\C) & \mapsto(G(f,1)\circ x_B\circ F(1,f)\mid f : A \to B)
		\end{cases}\]
	The equaliser of this pair of maps is evidently selecting the set of ``coherent'' sequences $(x_C)$ such that these two actions are equal, \ie such that the family $x_C : F(C,C) \to G(C,C)$ specifies a dinatural transformation $F \din G$  in the sense of \ref{dinana}.
\end{example}
\begin{example}[\upeyes Stokes' theorem is about co\fshyp{}ends]\label{stokkio}
	\index{Co/end!Stokes' theorem using ---s}\index{Stokes' theorem}
	This remark was first observed by \cite{1229249} and it follows their exposition almost word\hyp{}by\hyp{}word. Let $\N$ be the poset of natural numbers in the usual ordering, and $\Mod(\bR)$ be the category of real vector spaces.

	Fix a manifold $X$ (or some other sort of smooth space). Then we have functors
	\begin{itemize}
		\item $C: \N^\opp \to \Mod(\bR)$ where $C_n$ is the vector space freely generated by smooth maps $Y \to X$ where $Y$ is a compact, $n$\hyp{}dimensional, oriented manifold with boundary, and the induced map $\partial: C_{n+1} \to C_n$ is the boundary map. This is a chain complex, since $\partial\partial=\varnothing$.
		\item $\Omega: \N \to \Mod(\bR)$ is the de Rham complex; $\Omega_n = \Omega_n(X)$ is the space of $n$\hyp{}forms on $X$ and the induced map $d: \Omega_n \to \Omega_{n+1}$ is the exterior derivative. This is the usual de Rham complex of $X$.
	\end{itemize}
	Consider the usual tensor product on $\Mod(\bR)$; taking the object-wise tensor product of $C$ nd $\Omega$, we obtain a functor $C \otimes \Omega: \N^\opp \times \N \to \Mod(\bR)$, while there is also the constant functor $\bR: \N^\opp \times \N \to \Mod(\bR)$

	Then Stokes' theorem asserts that we have a cowedge
	\[\notag
		\rbag : C \otimes \Omega \to \bR
	\]
	which, given a map $Y \to X$ and a differential form $\omega$ on $X$, pulls the form back to $Y$ and integrates it (returning $0$ if it's the wrong dimension): $\rbag_n : C_n\otimes\Omega_n \to \bR$ is defined by $(\var{Y}{X},\omega)\mapsto \int_Y \varphi^*\omega$ if $\varphi : Y\to X$. (The integral $\int_Y$, here, is not a co\fshyp{}end).

	In fact, the cowedge condition for $\rbag$ amounts to the commutativity of the square
	\[
		\vcenter{\xymatrix{
		C_{n+1}\otimes\Omega_n \ar[r]^{\partial\otimes \id}\ar[d]_{\id\otimes d} & C_n\otimes\Omega_n \ar[d]^{\rbag_n}\\
		C_{n+1}\otimes \Omega_{n+1} \ar[r]_-{\rbag_{n+1}} & \bR
		}}
	\]
	which according to the definition of $\rbag$ is like saying that
	\[\textstyle \int_{\partial Y} \varphi|_{\partial Y}^*\omega = \int_Y \varphi^* d\omega.\]
	This is precisely Stokes' theorem.

	Exercise for the reader: what is the natural induced map $\int^{n \in \N} C_n \otimes \Omega_n \to \bR$? (hint: express the coequaliser defining $\int^{n \in \N} C_n \otimes \Omega_n$ as the degree\hyp{}zero cohomology of a suitable bicomplex $(C_\bullet\otimes \Omega_\bullet, D)$).
\end{example}
\begin{exercises}
\item \index{Wedge}
Prove equations \eqref{chista} and \eqref{chillautra}. Which commutative diagrams do you need to derive them?
\item \label{ex1:donotcomp} Show with an example that dinatural transformations $\alpha : P\din  Q, \beta : Q\din  R$ cannot be composed in general. Nevertheless, there exists a composition rule of a dinatural $\alpha : P\din Q$ with a natural $\eta : P'\To P$ which is again dinatural $P'\din Q$, as well as a composition $P\din Q\To Q'$ (hint: the appropriate diagram results as the pasting of a dinaturality hexagon and two naturality squares).
\item \label{ex1:endofconstant} What is the end of the constant functor $\Delta_D : \C^\opp\times\C\to \D$ at the object $D\in \D$? What is its coend?
\item \label{ex1:compoextra} \upeyes Show that extranatural transformations compose accordingly to these rules:
\begin{itemize}
	\item (stalactites) Let $F, G$ be functors of the form $\C^\opp \times \C \to \D$. If $\alpha_{X, Y} : F(X, Y) \to G(X,Y)$ is natural in $X, Y$ and $\beta_X : G(X, X) \to H$ is extranatural in $X$ (for some object $H$ of $\D$), then
	      \[\beta_X \circ \alpha_{X, X}: F(X, X) \to H\notag\]
	      is extranatural in $X$.
	\item (stalagmites) Let $G, H$ be functors of the form $\C^\opp \times \C \to \D$. If $\alpha_X : F \to G(X, X)$ is extranatural in $X$ (for some object $F$ of $\D$) and $\beta_{X, Y} : G(X, Y) \to H(X, Y)$ is natural in $X, Y$, then
	      \[
		      \beta_{X, X} \circ \alpha_X: F \to H(X, X)
		      \notag\]
	      is extranatural in $X$.
	\item (yanking) Let $F, H$ be functors of the form $\C \to \D$, and let $G : \C \times \C^\opp \times \C \to \D$ be a functor. If $\alpha_{X, Y} : F(y) \to G(X, X, Y)$ is natural in $y$ and extranatural in $X$, and if $\beta_{X, Y} : G(X, Y, Y) \to H(X)$ is natural in $X$ and extranatural in $Y$, then
	      \[\beta_{X, X}\circ \alpha_{X, X} : F(X) \to H(X)\notag\]
	      is natural in $X$.
\end{itemize}
Express these laws as equalities between suitable string diagrams (explaining also the genesis of the names `stalactite' and `stalagmite').
\item \label{ex1:dinatarentextra} \righteyes Prove that dinaturality is strictly more general than extranaturality, following this plot.

Let $\Delta[1]=\{ 0 \to 1 \}$ be the `generic arrow' category, and $S,T : \Delta[1]^\opp\times \Delta[1] \to \Set$ the functors respectively defined by
\[
	\vcenter{
		\xymatrix{
			\{1\}\ar@{}[dr]|S\ar[r]^{c_1}\ar[d] & \{1,2\}\ar[d]^\sigma & (0,0)\ar[r]\ar[d] & (0,1)\ar[d]& \{1\} \ar[r]\ar[d]\ar@{}[dr]|T & \{1\}\ar[d]^{c_2}\\
			\{1\} \ar[r]_{c_2} & \{1,2\} 	 		 & (1,0) \ar[r] & (1,1) &		 \{1\} \ar[r]_{c_2} & \{1,2\}
		}
	}
	\notag\]
where $c_i$ chooses element $i\in\{1,2\}$, and $\sigma$ permutes the two elements. Show that there exists a dinatural transformation $T \din  S$, whose components are identities, which is not extranatural when in \ref{extranatural} we choose $\A=*$ and $\B =\C^\opp = \Delta[1]$.
\item \label{ex1:vector-of-coends} If $(X, \Omega, \mu)$ is a measure space, the integral of a vector function $\vec F : X \to \bR^n$ such that each $F_i = \pi_i \circ F : X\to \bR$ is measurable and has finite integral, is the vector whose entries are $\Big(\int_X F_1 d\mu, \cdots, \int_X F_n d\mu \Big)$.

Prove that category theory possesses a similar formula, \ie that if $F : \C^\opp\times \C \to \prod_{i=1}^n\A_i$ is a functor towards a product category, such that
\begin{itemize}
	\item each $\A_i$ has both an initial and a terminal object;
	\item each co\fshyp{}end $\int (\pi_i\circ F)$ exists;
\end{itemize}
then the `vector' of all these co\fshyp{}ends, as an object $\big( \int F_1 , \dots, \int F_n\big) \in \prod \A_i$, is the (base of a universal co\fshyp{}wedge forming the) co\fshyp{}end of $F$. Where did you use the assumption that each $\A_i$ has an initial and a terminal object? \index{Wedge}
\item \label{ispull} Let $\D$ be a category. Show that the end of a functor $T : \Delta[1]^\opp\times \Delta[1]\to \D$ is the pullback of the morphisms $$T(0,0)\xto{T(0,d_0)} T(0,1)\xot{T(d_0,1)}T(1,1),$$ \ie that the following square is a pullback in $\D$:
\[
	\vcenter{\xymatrix{
			\int_{i\in\Delta[1]} T(i,i) \pb \ar[r]\ar[d] & T(0,0)\ar[d]\\
			T(1,1) \ar[r] & T(0,1)
		}}
	\notag\]
Dualise to the coend being a pushout.
\item \label{ex1:forgroups} \awful Let $G$ be a topological group, and $\text{Sub}(G)$ the set of its subgroups partially ordered by inclusion; let $X$ be a $G$\hyp{}space, \ie a topological space with a continuous action $G\times X \to X$ (the product $G\times X$ has the product topology).

We can define two functors $\text{Sub}(G) \to \Spc$, sending $(H\le G) \mapsto G/H$ (this is a covariant functor, and $G/H$ has the induced quotient topology as a space; there is no need for $H$ to be normal) and $(H\le G)\mapsto X^H$ (the subset of $H$\hyp{}fixed points for the action; this is a contravariant functor).
\begin{itemize}
	\item Compute the coend
	      \[
		      \fko_G(X) = \int^{H\le G} X^H \times G/H
		      \notag\]
	      in the category $\Spc$ of spaces, if $G = \Z/2$ has the discrete topology;
	\item Give a general rule for computing $\fko_G(X)$ when $G = \Z/n\Z$ is cyclic with $n$ elements;
	\item Let instead $\text{Orb}(G)$ be the \emph{orbit category} of subgroups of $G$, whose objects are subgroups but $\hom(H,K)$ contains $G$-\emph{equivariant} maps $G/H \to G/K$. Let again $X^{\firstblank}$ and $G/\firstblank$ define the same functors, now with different action on arrows. Prove that \[\notag \int^{H\in\text{Orb}(G)} X^H \times G/H \cong X\] (\emph{Elmendorf reconstruction}, \cite{Elmendorf1983}).
	\item Let $E|F$ be a field extension, and $\{H \le \text{Gal}(E|F)\}$ the partially ordered set of subgroups of the Galois group of the extension. Compute (in the category of \emph{rings}) the coend
	      \[
		      \int^{H} E^H \times \text{Gal}(H|F)
		      \notag\]
\end{itemize}
\item \index{_aaa_int@$\int$} \label{ex1:iscoeq} Dualise the construction in \ref{coends_as_colims}, to obtain a characterisation for the coend $\int^C F(C,C)$, characterised as the coequaliser of a pair $(F^*, F_*)$ as in
\[
	\coprod_{C\to C'} F(C',C) \rightrightarrows \coprod_{C\in\C} F(C,C)
	\notag\]
\item \label{ex1:niftynat} Find an alternative proof that natural transformations can be written as an end (see \ref{naturalu}), using the characterisation of $\int_C\D(FC,GC)$ as an equaliser in \ref{endsareeq}: as a subset of $\prod_{C\in\C} \D(FC,GC)$, is precisely the subset of natural transformations $\{\tau_C : FC \to GC\mid Gf\circ \tau_C = \tau_{C'}\circ Ff,\;\forall f : C \to C'\}$.
\item \label{ex1:coendofid} What is the co\fshyp{}end of the identity functor $1_{\C^\opp\times \C} : \C^\opp\times \C\to \C^\opp\times \C$? Use the bare definition; use the characterisation of co\fshyp{}ends as co\fshyp{}limits; feel free to invoke Exercise~\ref{ex1:vector-of-coends}.
\item A set of objects $\catS\subset\C$, regarded as a full subcategory, \emph{finitely generates} a category $\C$ if for each object $X\in\C$ and each arrow $f : S\to C$ from $S\in \catS$ there is a factorisation
\[
	S \xrightarrow{g} \coprod_{i=1}^n S_i\xrightarrow{h_C} C
	\notag
\]
where $h_C$ is an epimorphism and $\{S_1, \dots, S_n\}\subset \catS$ ($n$ depends on $C$ and $f$).

Suppose $T : \C^\opp\times \C\to \Sets$ is a functor, finitely continuous in both variables, and $\C$ is finitely generated by $S$. Then if we denote $T|_{\catS} : \catS^\opp\times\catS \to \Sets$ the restriction, we have an isomorphism
\[
	\int^{C\in \catS} T|_{\catS}(C,C) \cong \int^{C\in \C} T(C,C)
	\notag
\]
induced by a canonical arrow $\int^{C\in \catS} T'(C,C) \to \int^{C\in \C} T(C,C)$.
\item \label{ex1:acriterion} Let $F\dashv U : \C\leftrightarrows \D$ be an adjunction, and $G : \D^\opp\times \C\to \E$ a functor; then there is an isomorphism
\[
	\int^C G(FC,C) \cong \int^D G(D,UD).
	\notag\]
Show that the converse of this result is true: if the above isomorphism is true for any $G$ and natural therein, then there is an adjunction $F\dashv U$.
\item Give a different proof of Fubini theorem \ref{fubozzo} as follows: define again $H_\C : D\mapsto \hom_\C \otimes D$, and find suitable unit and counit maps (see \ref{unit_counit}).
\begin{itemize}
	\item The unit of the adjunction $\eta : \id_\D \To \int_C \circ H_\C$ is determined as a morphism
	      \[\notag\eta_D : D \to \int_C \big(\C(C,C)\otimes D\big)\]
	      (parenthesisation is important: $\int_C \big(\C(C,C)\otimes D\big)$ is \emph{not} isomorphic to $\big(\int_C \C(C,C)\big)\otimes D$) which in turn corresponds to a wedge
	      \[\notag\eta_{D,(u)} : D \to \coprod_{v : C \to C} D\]
	      defined on the component $u : C\to C$ by the inclusion $i_u$ in the coproduct.
	\item The counit of the adjunction $\epsilon : H_\C\circ \big(\int_C \firstblank ) \To \id_{\Cat(\C^\opp\times\C,\D)}$ is determined as a natural transformation having components
	      \[\notag\epsilon_{CC'} : \coprod_{u : C \to C'} \int_C F(C,C) \to F(C,C')\]
	      which in turn have components
	      \[\notag
		      \epsilon_{CC'}^u : \int_C F(C,C) \to F(C,C')
	      \]
	      by the universal property of a coproduct, determined by sending a coherent family $(a_C\mid C\in \C)$ into $F(C, u)(a_C)$ (o equivalently by the wedge condition, $F(u, C')(a_{C'})$).
\end{itemize}
Show that these two definitions set up the desired adjunction (prove the triangle identities for $\eta$ and $\epsilon$, as in \ref{zigozago}).

Hint: to show that $\epsilon^u_{CC'}$ is indeed natural in $C,C'$, take two arrows $\cvar{C_0}{\alpha}{C_1}$ and $ \cvar{C_0'}{\beta}{C_1'}$ and split the naturality square into the pasting of two smaller squares
\[\notag
	\vcenter{\xymatrix@C=2cm{
	\coprod_{C_1 \to C_0'} \int_C F\ar[r]^{\epsilon^u_{C_1C_0'}}\ar[d]_{\firstblank{\circ\beta}} & F(C_1, C_0')\ar[d]^{F(1,\beta)} \\
	\coprod_{C_1 \to C_1'} \int_C F\ar[r]^{\epsilon^u_{C_1C_1'}}\ar[d]_{\alpha\circ\firstblank} & F(C_1, C_1')\ar[d]^{F(\alpha,1)} \\
	\coprod_{C_0 \to C_1'} \int_C F \ar[r]_{\epsilon^u_{C_0C_1'}}& F(C_0, C_1')
	}}
\]
each of which commutes by evident reasons.
\end{exercises}

\chapter{Yoneda and Kan}\label{sec:tre}
\begin{abstract}
	In this chapter we begin to learn the rules of co\fshyp{}end calculus.	First, we re-enact the well known Yoneda lemma as an isomorphism of co\fshyp{}ends. Yoneda lemma is one of the deepest results that can be examined with the technology built so far: every functor $F: \C^\opp\to\Set$ can be decomposed as a coend
	\[F \cong \int^CFC \times \C(\firstblank,C).\notag\]
	This isomorphism has plenty of consequences; it is in fact equivalent to the assertion that every presheaf is a colimit of a certain diagram $D : \C \to [\C^\opp,\Set]$ of representables, and that there is a canonical choice for such $D$.	Then we move to investigate \emph{Kan extensions}; the famous tenet that `every thing is a Kan extension' is, here, translated into the more appealing (for us!) statement that every sufficiently nice thing is a co\fshyp{}end.	We then glimpse to \emph{formal category theory}; this latter part draws on \cite{Graya} and other classical sources of the Australian school of category theory.
\end{abstract}
\epigraph{The reason why this technique is so fast is that you are not trying to cut them; you are throwing your sword into them.}{s\=oke M\@. Hatsumi}
Along the present chapter, an abstract 2\hyp{}category will be denoted in sans\hyp{}serif, $\sfA, \sfK\dots$; an object in $\sfK$ will be denoted in uppercase roman $A,B,C\dots$; an object of a 0\hyp{}cell $A\in\tCat$ will be denoted in roman lowercase $a,b,a',b'\dots$

This complies with the common notation for papers in 2\hyp{}dimensional algebra, and the different notation employed elsewhere should not cause any confusion. Note in particular that we make the distinction between the 1\hyp{}category $\Cat$ (categories and functors) and the 2\hyp{}category $\tCat$ (categories, functors and natural transformations).
\section{The Yoneda lemma and Kan extensions}
For the ease of exposition, we recall the statement of the \emph{Yoneda lemma}:
\begin{lemma*}
	Let $\yon : \C \to \Cat(\C^\opp,\Sets)$ be the functor sending $X\in\C$ to the presheaf associated to $X$; then, for every $F\in\Cat(\C^\opp,\Sets)$, there exists a bijection between the set of natural transformations $\yon X\To F$ and the set $FX$. This bijection is moreover natural in the object $X$.
\end{lemma*}
We carry on the proof in full detail in our \ref{lem:the-real-yoda} to convince even the most skeptical reader that this result is a tautology. Few tautologies are, however, richer of meaning.
\begin{remark}[The Yoneda-Grothendieck philosopy]\label{yogrophilo}
	Yoneda lemma entails that there is a copy of a category $\C$ in its category of presheaves $\Cat(\C^\opp,\Set)$; indeed, if we let $F=\yon B$ to be representable (see \ref{def:reprepsh}) in the isomorphism above, we see that
	\[ \tCat(\C^\opp,\Set)(\yon A, \yon B)\cong \C(A,B) \]
	This means precisely that the Yoneda embedding functor $\yon_\C : \C \to \Cat(\C^\opp,\Set)$ is fully faithful. To all practical purposes then, every question about $\C$ is a question about the copy of $\C$ inside its category of presheaves; the latter category is, however, always very well\hyp{}behaved under certain aspects; for example, it is always complete and cocomplete (see \ref{being-complete}) even when $\C$ isn't.

	Thus, every diagram in $\C$, say $D: \J \to \C$ admits a limit \emph{in $[\C^\opp,\Set]$}, as we can build the limit
	\[\lim \yon(D_I) = \lim \hom_{\C}(\firstblank, D_I).\]
	If this functor is representable, say by an object $L$, it is easy to see that $L$ realises the universal property of $\lim D_I$ \emph{in $\C$}; in fact, this is a necessary and sufficient condition: thus, a completeness request on $\C$ has a 1-1 translation in terms of a request of representability of certain functors.

	The idea that properties of $\C$ can (and should) be translated into representability properties, and more extensively that representability conditions play a major role in category theory, algebra and geometry, was first advocated by A. Grothendieck, and it makes heavy use of Yoneda lemma; thus we colloquially refer to it as the \emph{Yoneda-Grothendieck philosophy}.

	\smallskip
	The next step of this introduction to the chapter is to single out certain universal 2-cells in a 2-category, and the fundamental properties thereof.
\end{remark}
\begin{notation}[Extensions and lifts]\label{notat:liftext}
	Let $f,g$ be 1\hyp{}cells of a 2\hyp{}category $\sfK$, respectively $f : A\to B$ (the extendable arrow) and $g : A \to C$ (the extendant arrow). We say that a pair $\langle u,\eta\rangle$ \emph{exhibits} the left extension $\lan_gf$ of $f$ along $g$ if the 1\hyp{}cell $u : C \to B$ and the 2\hyp{}cell $\eta : f \To ug$ can be arranged in a triangle
	\[\vcenter{\xymatrix{
		&A\ar[dr]^g\ar@{}[d]|{\overset{\eta}\To} \ar[dl]_f & \\
		B &&\ar[ll]^u C
		}}\]
	initial among all such. This means that every pair $\langle v, \alpha\rangle$ of a 1\hyp{}cell $\var[v]{C}{B}$ and a 2\hyp{}cell $\alpha : f \To vg$ factors uniquely through $\eta$ as a composition $\alpha = (\bar\alpha * g)\circ\eta$ (the cell $\bar\alpha * g$ is the \emph{whiskering} of \ref{comporizz}).
\end{notation}
\begin{remark}\label{specchietto_dei_lift}
	The notion of left extension is subject to dualisation, in a way that it is worth recording explicitly: it is useful to have a diagram illustrating at once all these universal constructions. The diagram below has to be parsed as follows: moving horizontally reverses the direction of 1\hyp{}cells, but not of 2\hyp{}cells; moving vertically but not horizontally reverses the direction of 2\hyp{}cells, but not of 1\hyp{}cells. The universality request acquires the same shape, and when we write (for example) $\scriptstyle\deduction{\leeft_gf}{h}{f}{gh}$ we mean that there is a bijection between 2\hyp{}cells $\leeft_gf\To h$ and 2\hyp{}cells $f \To gh$.
	\[\label{left-and-right-lift-and-ext}
		\begin{array}{cc}
			\xymatrix{
			A \ar@{}[dr]|(.3){\Swarrow\eta}\ar[d]_g \ar[r]^f        & B                                           \\
			C \ar@{.>}[ur]_{\lan_gf}                                & {\scriptstyle\deduction{\lan_gf}{h}{f}{hg}}
			}
			                                                        &
			\xymatrix{
			{\scriptstyle\deduction{\leeft_gf}{h}{f}{gh}}           & C\ar[d]^g                                   \\
			B\ar[r]_f \ar@{.>}[ur]^{\leeft_gf}                      & \ar@{}[ul]|(.3){\Nearrow\eta} A
			}                                                                                                     \\
			\xymatrix{
			A \ar@{}[dr]|(.3){\Nearrow\varepsilon}\ar[d]_g \ar[r]^f & B                                           \\ C
			\ar@{.>}[ur]_{\ran_gf}                                  & {\scriptstyle\deduction{hg}{f}{h}{\ran_gf}}
			}
			                                                        &
			\xymatrix{
			{\scriptstyle\deduction{h}{\rift_gf}{gh}{f}}            & C\ar[d]^g                                   \\ B\ar[r]_f
			\ar@{.>}[ur]^{\rift_gf}                                 & \ar@{}[ul]|(.3){\Swarrow\varepsilon} A
			}                                                                                                     \\
		\end{array}
	\]
	So, a right extension is a left extension in $\sfK^\co$, the 2\hyp{}category where all 2\hyp{}cells have been reversed; a left lifting is a left extension in $\sfK^\opp$, the 2\hyp{}category where all 1\hyp{}cells have been reversed; and a right lifting is a left extension in $\sfK^\coop$, where both 1- and 2\hyp{}cells have been reversed.
\end{remark}
\begin{definition}[Pointwise and absolute extension]\label{ptwise_ext}
	We say that a left (resp., right) extension is \emph{pointwise} if for every object $C$ and $k : X\to C$ the diagram obtained pasting at $g$ the comma object (see \ref{def:comma}) $(g/k)$ (resp., $(k/g)$) is again a left extension \cite[5.2]{street1981conspectus}; this means that in every diagram of 2\hyp{}cells
	\[
		\vcenter{\xymatrix{
				(g/k) \ar@{}[dr]|\Swarrow \ar[r]\ar[d]& A \ar@{}[dr]|(.35)\Swarrow \ar[r]^f\ar[d]_g  & B \\
				X \ar[r]_k & C \ar[ur] &
			}}
	\]
	if the left triangle is a left extension and the square is a comma object, then the whole pasting diagram is a left extension; dually (in $\sfK^\co$) for a right extension.

	We say that a left extension is \emph{absolute} if it is preserved by all 1\hyp{}cells, in the same sense a functor preserves co\fshyp{}limits: in the left extension diagram $\langle k,\eta\rangle$
	\[
		\vcenter{\xymatrix{
				A \ar@{}[dr]|(.35)\Swarrow \ar[r]^f\ar[d]_g  & B \ar[r]^k & Y\\
				C \ar[ur]_h &
			}}
	\]
	the 1\hyp{}cell $k$ \emph{preserves} the extension if the whiskering $k * \eta$ exhibits $k\circ h$ as the left extension of $k\circ f$ along $g$: there is an isomorphism $k(\lan_gf)\cong \lan_g(kf)$, canonically.

	Of course, the same nomenclature applies to define pointwise and absolute left and right liftings.
\end{definition}
\begin{remark}
	Later in the chapter we will see that extensions in $\tCat$ exist, whenever sufficient co\fshyp{}limits exist in the codomain of the extendable arrow, and they carry a very expressive theory (to the point that Mac Lane states that `everything is a Kan extension' in \cite{McL}). Instead, few lifts (both left and right) do exist in $\tCat$.

	The deep reason why this is true is that the Yoneda lemma does not hold in the category $\tCat^\opp$, as the opposite of a functor is not a functor in general.

	In fact, one tenet of the present chapter is that left extensions in $\tCat$ (where they are called \emph{Kan} extensions, in honour of D.M. Kan (see \ref{why_lans_are_lans}): this justifies the notation as $\lan$ is a portmanteau for $\text{l(-eft)} + \text{(K-)an}$) can be computed as certain explicit co\fshyp{}limits (if they are left extensions, they are colimits; right extension, they are limits). More precisely, pointwise Kan extensions can be expressed as co\fshyp{}ends.

	It is not possible to do the same for lifts, due to the peculiar behavior of the universality requests in \ref{left-and-right-lift-and-ext}; the reason is somewhat trivial --- the dual of a functor is generally not a functor any more.
\end{remark}
\begin{example}
	To see that there is no way to compute a lift as a co\fshyp{}limits, let us consider the following example: let us consider the discrete category $\{0,1\}$, and let $\delta_0 \colon \C \to \{0,1\}$, $\delta_1 \colon \D \to \{0,1\}$ be the obvious constant functors (with disjoint images) from any two non\hyp{}trivial categories $\C,\D$.

	Then there are no liftings of $\delta_0$ through $\delta_1$, no matter how complete and cocomplete $\C$ and $\D$ are. In fact, for any functor $F \colon \C \to \D$, there are no natural transformations $\delta_1 \circ F \to \delta_0$ nor $\delta_0 \to \delta_1 \circ F$. The same is true, if we substitute $\{0,1\}$ with its free cocompletion $\Cat(\{0,1\}^\opp,\Set)\cong \Set/\{0,1\}$.
\end{example}

The deep reason why we do not have nice formulae for Kan liftings, is that $\tCat^\opp$ lacks an internal concept of (co)end; in other words the internal language of $\tCat^\opp$ is not expressive enough.

On the contrary, the theory of extensions behaves much better, and in fact it is possible to give a neat characterisation of pointwise extensions:
\begin{proposition}
	Let $\A,\B,\C$ be categories, $f,g$ be functors as in the triangle below. Then, the following conditions are equivalent:
	\begin{itemize}
		\item The triangle
		      \[	\vcenter{\xymatrix{
					      \A \ar@{}[dr]|(.35)\Nearrow \ar[r]^f\ar[d]_g  & \B \\
					      \C \ar[ur]_h &
				      }}\]
		      is a pointwise right Kan extension;
		\item The triangle is a right extension, and it is preserved by every representable $\yon^\lor(B) = \B(B,\firstblank) : \B \to \Set$, meaning that
		      \[\ran_g\hom_\B(B,f) \cong \hom_\B(B,\ran_gf)\]
	\end{itemize}
\end{proposition}
\begin{proof}
	See \cite[5.4]{lehner}.
\end{proof}
\section{Yoneda lemma using co\fshyp{}ends}
Tightly linked to the Yoneda lemma is the \emph{density theorem} of \ref{thm:yoda-is-dense}: every presheaf $F : \C^\opp\to \Set$ is, canonically, the colimit of a diagram of representables, \ie functors lying in the image of $\yon_\C : \C \to \tCat(\C^\opp, \Sets)$; we prove this result in full detail in \ref{lem:the-real-yoda}, without using coends; however, the scope of the present chapter is to convince the reader that  co\fshyp{}end calculus allows to concisely and effectively rephrase both the Yoneda lemma and this result.

In our chapter \ref{section:weight}, thanks to the machinery of \emph{weighted co\fshyp{}limits}, the fact that every presheaf is a colimit of representables acquires the more alluring form
\begin{quote}
	Every presheaf is the weighted colimit of the Yoneda embedding, weighted by the presheaf itself.
\end{quote}
\begin{proposition}[ninja Yoneda Lemma]\label{ninjayo}\index{Yoneda lemma!ninja ---}
	For every functor $K : \C^\opp\to\Set$ and $H : \C\to \Set$, we have the following isomorphisms (natural equivalences of functors):
	\begin{align}
		K & \cong \int^CKC\times \C(\firstblank,C) \label{nyo_1}  \\
		K & \cong \int_C \Set(\C(C,\firstblank),KC)\label{nyo_2}  \\
		H & \cong \int^CHC\times \C(C,\firstblank) \label{nyo_3}  \\
		H & \cong \int_C \Set(\C(\firstblank,C),HC) \label{nyo_4}
	\end{align}
	where the functor $\C\times\C^\opp\times\C^\opp$ in \eqref{nyo_2} is defined by
	\[\lambda AXC.\Set(\C(A,X),KC),\notag \]
	and similarly for $\Set(\C(\firstblank,C),HC)$.
\end{proposition}
\begin{remark}
	The name \emph{ninja Yoneda lemma} is chosen to instill in the reader the sense that this result \emph{is} the Yoneda lemma in disguise; the name comes from \cite{20451} where T\@. Leinster offers the same argument we are about to give as a proof of \ref{ninjayo}:
	\begin{quote}
		Th[e above one is] often called the \emph{Density Formula}, or (by Australian ninja category theorists) simply the Yoneda Lemma. (but Australian ninja category theorists call \emph{everything} the Yoneda Lemma\dots).
	\end{quote}
\end{remark}
Note that the isomorphisms \eqref{nyo_2} and \eqref{nyo_4} follow directly from \ref{naturalu}. We prove only isomorphism \eqref{nyo_1}, as the remaining one can be easily obtained by dualisation.

We put a particular emphasis on giving a detailed proof once, as it is the first argument that truly exploits co\fshyp{}end \emph{calculus}; it is a paradigmatic example of the style of proofs we're using from now on.
\begin{proof}
	Consider the chain of isomorphisms
	\begin{align*}
		\Set\Big(\displaystyle \int^{C\in\C} KC\times \C(X,C),Y\Big) & \cong \int_{C\in\C} \Set\big( KC\times \C(X,C),Y \big)           \\
		                                                             & \cong \int_{C\in\C}\Set(\C(X,C),\Set(KC,Y))                      \\
		                                                             & \cong [\C, \Set]\big(\C(X,\firstblank),\Set(K\firstblank,Y)\big) \\
		                                                             & \cong\Set(KX,Y)
	\end{align*}
	where the first step is motivated by the coend\hyp{}preservation property of the $\hom$ functor \ref{commuhom}, the second follows from the fact that $\Sets$ is a cartesian closed category, where
	\[
		\Sets(A\times B, C)\cong \Sets(A, \Set(B,C))
	\]
	for all three sets $A,B,C$, naturally in all arguments, and where $C^B = \Set(B,C)$ is the set of all functions $B\to C$, and the final step exploits Theorem \ref{naturalu} and the classical Yoneda Lemma that says that the set of natural transformations $\yon^\lor(C)\To F$ is in bijection with the set $FC$ for every presheaf $F : \C\to \Set$. Here, $F = \Set(K\firstblank,Y)$.

	Every step of this chain of isomorphisms is natural in the object $Y$; since the Yoneda embedding $\yon_\C : \C^\opp \to [\C,\Set]$ is fully faithful, the isomorphism of functors
	\[
		\Sets\Big(\textstyle \int^C KC\times \C(X,C),Y\Big)\cong \Sets(KX,Y)
	\]
	ensures in turn that there exists an isomorphism between the represented objects $\int^C KC\times \C(X,C)\cong KX$.

	This is moreover natural in the object $X$. Accepting the validity of \eqref{nyo_1} and re-doing each step backwards, we can get back the Yoneda lemma in the form expressed in \ref{lem:the-real-yoda}.
\end{proof}
From now one we will make frequent use of the notion of \emph{tensor} and \emph{cotensor} in an enriched category; the definitions for a generic $\V$\hyp{}category can be found for example in \cite[Ch. 6]{Bor2}, and in particular its Definition 6.5.1: we have already introduced this definition in the proof of our \ref{fubozzo}.
\begin{definition}[Tensor and cotensor in a $\V$\hyp{}category]\label{tenscotens}\index{Co/tensors}\index{_aaa_pitch@$\pitchfork$}\index{Co/tensors}\index{_aaa_otimes@$\otimes$}
	Let $\C$ be a $\V$\hyp{}enriched category (see \cite[6.2.1]{Bor2} or our \ref{enrichcat}), the \emph{tensor} $\otimes : \V\times \C\to \C$ (when it exists) is a functor $(V, C)\mapsto V\otimes C$ such that there is the isomorphism \[ \label{ten} \C(V\otimes C, C')\cong \V(V, \C(C,C')),\] natural in all components; dually, the \emph{cotensor} in an enriched category $\C$ (when it exists) is a functor $(V, C)\mapsto V\pitchfork C$ (contravariant in $V$) such that there is the isomorphism \[ \label{coten} \C(C', V\pitchfork C)\cong \V(V, \C(C',C)),\] natural in all components.
\end{definition}
\begin{example}
	Every co\fshyp{}complete, locally small category $\C$ is naturally $\Sets$\hyp{}co\fshyp{}tensored by choosing $V\pitchfork C\cong \prod_{v\in V}C$ and $V\otimes C\cong \coprod_{v\in V}C$. We have employed this construction in the proof of \ref{fubozzo}.
\end{example}
\begin{remark}\label{cocotens}
	The tensor, hom and cotensor functors are the prototype of what's called a `two variable adjunction' (see \cite[§1.1]{Gray}); given the $\hom$\hyp{}objects of a $\V$\hyp{}category $\C$, the tensor
	\[\firstblank\otimes\secondblank : \V\times \C\to \C\]
	and the cotensor
	\[\firstblank\pitchfork\secondblank : \V^\opp\times \C\to \C\]
	can be characterised as adjoint functors to the hom functors, saturated in their first or second component: usual co\fshyp{}continuity properties of the co\fshyp{}tensor functors are implicitly derived from this characterisation.
\end{remark}
\begin{remark}[\itsnonsense, The Yoneda embedding is a Dirac delta]\label{dirac}\index{Functor!representable}
	In functional analysis, the Dirac delta appears in the following convenient abuse of notation:
	\[
		\int_{-\infty}^\infty f(x) \delta(x-y)dx = f(y)
	\]
	(the integral sign is not a co\fshyp{}end). Here $\delta(x-y):=\delta_y(x)$ is the \emph{$y$\hyp{}centered delta distribution}, and $f : \bR\to \bR$ is a continuous, compactly supported function on~$\bR$.

	It is really tempting to draw a parallel between this relation and the ninja Yoneda lemma, conveying the intuition that representable functors on an object $c\in\C$ play the r\^ole of $c$\hyp{}centered delta distributions.

	If the relation above is written as $\langle f,\delta_y\rangle = f(y)$, interpreting integration as an inner product between functions (or more precisely as a universal pairing between a space and its dual), then the ninja Yoneda lemma says the same thing for categories: each presheaf $F : \C^\opp \to \Sets$, can be paired with a distribution concentrated on the point $C$ in an `inner product' $\langle \yon_C, F\rangle = \int^X \yon_C(X)\times FC$; the latter object is now isomorphic to $FC$, in the same way the integral of a smooth function $f$ against a $y$\hyp{}centered delta equals $f(y)$.
	\begin{center}
		\begin{figure}[h!]
			\label{fig:profuncs}
			\begin{tikzpicture}
				\matrix (m) [
					matrix of math nodes,
					every node/.append style={
							minimum height=36pt,
							inner sep=8pt,
							execute at begin node=\strut\displaystyle,
						},
				]{
					\int_{ x \in\bR} & f(x) & \circ  & \delta(x-y)dx & =     & f(y) \\
					\int^X                  & FX   & \times & \C(Y,X)   & \cong & FY   \\
				};
				\foreach \i in {1,...,6}
				\node[inner sep=-4pt,draw,dashed,gray,fit=(m-1-\i)(m-2-\i)]{};
			\end{tikzpicture}
			\caption{The analogy between the pairing of a function and a delta distribution, and the ninja Yoneda lemma.}
		\end{figure}
	\end{center}
\end{remark}
\section{Kan extensions using co\fshyp{}ends}
In the following series of remarks, $G : \C \to \E$ and $F : \C \to \D$ are functors; for the sake of exposition we assume that $\Lan_GF$ exists for all such $F$'s; we shall see in the future that a sufficient condition for this to be true is that $\D$ is a cocomplete category (see \ref{being-complete}), thus we assume that all colimits exist in $\D$.
\begin{remark}
	The correspondence $F\mapsto \Lan_GF$ is a functor
	\[\tCat(\C,\D) \to \tCat(\E,\D);\]
	to every natural transformation $\alpha : F \To F'$ is associated a natural transformation $\Lan_G\alpha : \Lan_GF\To \Lan_GF'$ defined as the unique 2\hyp{}cell $\xi$ such that
	\[
		\vcenter{\xymatrix{
		\C \ar@{}[dr]|(.3){\Swarrow\eta'}\ruppertwocell^F{\alpha}\ar[r]\ar[d]_G & \D \\
		\E\ar[ur]_{\Lan_GF'} &
		}}
		\qquad
		=
		\qquad
		\vcenter{\xymatrix{
				\C \ar@{}[dr]|(.3){\Swarrow\eta}\ar[r]^F\ar[d]_G & \D \\
				\E\ar[ur]\urlowertwocell_{\qquad\Lan_GF'}{\xi} &
			}}
	\]
\end{remark}
\begin{remark}
	The correspondence $G\mapsto \Lan_G$ is a functor
	\[ \tCat(\C,\E)^\opp \to \tCat(\Cat(\C,\D),\Cat(\E,\D)); \]
  this means that to every natural transformation $\alpha : G\To G'$ is associated a natural transformation $\Lan_\alpha : \Lan_{G'} \To \Lan_G$, defined as the mate of the composite 2\hyp{}cell
	\[\label{Lan_is_a_functor} F \overset{\eta_F}\To \Lan_GF\circ G \overset{\Lan_GF * \alpha}\To \Lan_GF\circ G'.\]
	In order to show that the correspondence $G\mapsto \Lan_G$ is a pseudofunctor, we have to find suitable coherence isomorphisms $\Lan_{GG'} \cong\Lan_G\circ\Lan_{G'}$ and $\Lan_{\id}\cong \id$; this is the content of Exercise \ref{ex:kan_is_functor}.
\end{remark}
\begin{remark}
	The functor $\Lan_G$ is (the unique up to isomorphism) adjoint to the \emph{inverse image} functor $G^* = \firstblank\circ G$; this follows directly from the characterisation of $\Lan_GF$ as the functor $\E\to\D$ such that
	\[\tCat(\E,\D)(\Lan_GF,H)\cong \tCat(\C,\D)(F, HG).\]
\end{remark}
\begin{remark}
	Dually, the correspondence $G\mapsto \Ran_G$ is a functor
	\[ \tCat(\C,\E)^\opp \to \tCat(\Cat(\C,\D),\Cat(\E,\D)) \]
	that works as right adjoint to the precomposition functor $\firstblank\circ G$, and it is functorial, contravariant in the extendant component. Given $\alpha : G\To G'$ the components of $\Ran_\alpha$ are the mates of
	\[\label{Ran_is_a_functor} \Ran_{G'}F\circ G \overset{\Ran_{G'}*\alpha}\To \Ran_{G'}F\circ G' \overset{\epsilon}\To F. \]
\end{remark}
\begin{remark}
	From the universal property of $\Lan_G\firstblank$ we can derive the unit and the counit of the adjunctions $\Lan_G \adjunct{\epsilon^L}{\eta^L} G^*$ and $G^* \adjunct{\epsilon^R}{\eta^R} \Ran_G$: we leave this as Exercise \ref{ex:counit_of_lans} for the reader to spell out explicitly.
\end{remark}
\begin{proposition}\label{kan_are_coends}\index{Kan extension}\index{Functor!Kan extension ---}
	Let again $G : \C \to \E$ be a functor and $F : \C \to \D$ is a functor whose domain is a cocomplete category. Since both the co\fshyp{}tensors (see \ref{tenscotens}) and the co\fshyp{}ends involved in the equations below exist, then the left/right Kan extensions of $G : \C\to\E$ along $F : \C\to \D$ exist and there are isomorphisms (natural in $F$ and $G$)
	\[\label{kanend}
		\Lan_FG\cong \int^C \D(FC,\firstblank)\otimes GC\qquad\qquad
		\Ran_FG\cong \int_C {\D(\firstblank,FC)}\pitchfork GC.
	\]
	These Kan extensions are pointwise in the sense of \ref{ptwise_ext}.
\end{proposition}
Since this is our first instance of derivation in co\fshyp{}end calculus, we apply a certain pedantry to the explanation, duly recording the results that allow each step. Such verbosity will be soon abandoned, to let each reader profit from the instructive meditation following the chains of isomorphisms we write.
\begin{proof}
	The proof consists of a string of canonical isomorphisms, exploiting simple remarks and the results established so far: the same argument is offered in \cite[X.4.1-2]{McL}.
	\begin{align*}
		\tCat(\D,\E)\Big( \int^C \D(FC,\firstblank)\otimes GC,H\Big) & \cong \int_X\D\Big( \int^C \D(FC,X)\otimes GC,HX \Big) \\
		\ref{commuhom}                                               & \cong \int_{CX}\D(\D(FC,X)\otimes GC,HX)               \\
		\eqref{ten}                                                  & \cong \int_{CX}\Set(\D(FC,X),\E(GC,HX))                \\
		\ref{naturalu}                                               & \cong\int_C [\D(FC,\firstblank),\E(GC,H-)]             \\
		\eqref{nyo_4}                                                & \cong \int_C \E(GC,HFC)\cong \tCat(\C,\E)(G,HF).
	\end{align*}
	The case of $\Ran_FG$ is dually analogous and we leave it to the reader.
\end{proof}
\index{Functor!adjoint ---s}
\begin{corollary} Let $D : \A \to \B$ be a functor;
	\begin{enumtag}{cc}
		\item if $D$ is a left adjoint, then $D$ preserves all pointwise left Kan extensions that exist in $\A$; dually,
		\item if $D$ is a right adjoint, then $D$ preserves all pointwise right Kan extensions that exist in $\A$.
	\end{enumtag}
\end{corollary}
\begin{proof}
	Left adjoints commute with tensors, \ie $D(X\otimes A)\cong X\otimes DA$ for any $(X,A)\in \Sets\times \C$, and with colimits (see \ref{adj-preserve}). The result follows, and can easily be dualised.
\end{proof}
\begin{example}\index{Category!Kleisli ---}\index{Kleisli category}
	Let $T : \C\to \C$ be a monad (see \ref{def:monad}) on $\C$; the \emph{Kleisli category} $\Kl(T)$ of $T$ is defined having the same objects of $\C$ and morphisms $\Kl(T)(A,B) := \C(A, TB)$.

	Given any functor $F : \A\to \C$, the right Kan extension $T_F=\Ran_FF$, when it exists, is a monad on $\C$, that we call the \emph{codensity monad} of $F$; using the end expression for $\Ran_FF$ we get that on objects $T_F$ is defined as
	\[T_F(C) \cong \int_A \C(C,FA)\pitchfork FA.\]
	(See \ref{tenscotens} for the $\pitchfork$ operation.) Hom-sets in the Kleisli category $\Kl(T_F)$ can be characterised as
	\[\label{the_kleislona}
		\Kl(T_F)(C,C') \cong \int_A \Sets(\C(C', FA), \C(C, FA)).
	\]
	The multiplication and unit of $T_F$ can be found using the universal property of $T_F$:
	\begin{itemize}
		\item The multiplication is obtained as the mate of
		      \[
			      \label{codensity_mult}
			      \Ran_FF\circ \Ran_FF\circ F \xto{\Ran_F * \epsilon * F} \Ran_FF\circ F \xto{\epsilon_F} F
		      \]
		      under the adjunction isomorphism $[\C,\C](H,\Ran_FF)\cong [\A,\C](HF, F)$;
		\item the unit $\eta : \id \To \Ran_FF$ is obtained as the component of the unit of the adjunction $F^* \dashv \Ran_F$ at the identity.
	\end{itemize}
	The difficult part is to show that the two maps $\mu : T_F \circ T_F \To T_F$ and $\eta : \id_{\C} \To T_F$ indeed form a monad; we relegate this proof to an exercise in \ref{codensity_is_monad}; the interested reader shall try to translate in coend calculus the proof in \cite[pp. 67---71]{dubuc1970kan}.
\end{example}
\begin{remark}
	There is a dual theory of \emph{density comonads}: given a functor $F : \A \to \C$ the left Kan extension of $F$ along itself, when it exists, has the structure of a comonad $S^F$ (see \ref{comona}), where
	\begin{itemize}
		\item the co\hyp{}multiplication is obtained as the mate of
		      \[\label{density_comult}
			      F \xto{\eta_F} \Lan_FF\circ F \xto{\Lan_FF * \eta_F * F} \Lan_FF\circ \Lan_FF \circ F
		      \]
		      under the adjunction isomorphism
		\item the counit $\sigma : \Lan_FF \To \id$ is obtained as the component of the counit of the adjunction $\Lan_F \dashv F^*$ at the identity.
	\end{itemize}
	Find a similar expression for the hom sets of the co-Kleisli category $\text{coKl}(S^F)$ (and see Exercise \ref{ex:density}).
\end{remark}

\begin{example}[Stalks of a sheaf {(\cite[6.8 and {\S}7.1]{Grothendieck1972})}] (Credits to E. de Oliveira Santos)
	Let $(X,\tau)$ be a topological space, and $i_{p}\colon\{p\}\hookrightarrow X$ the inclusion of a singleton into $X$. From this, we get an induced functor
	\[
		\vcenter{\xymatrix@R=0cm{
		\mathcal{O}(i_{p}) : \mathcal{O}(X)
		\ar[r]
		&
		\mathcal{O}(\{p\})
		\\
		U
		\ar@{|->}[r]
		&
		i_{p}^{-1}(U)
		}}
	\]%
	Considering now left Kan extensions along the opposite of $\mathcal{O}(i_{p})$,
	\[
		\vcenter{\xymatrix{
		& \mathcal{O}(\{p\})^\opp
		\ar@{.>}[d]^{\Lan_{\mathcal{O}(i_{p})^\opp}\mathcal{F}}
		\\
		\mathcal{O}(X)^\opp
		\ar[ur]^{\mathcal{O}(i_{p})^\opp}
		\ar[r]_{\mathcal{F}}
		&
		\Sets,\ar@{}[ul]|(.3)\Nearrow
		}}
	\]
	we obtain a functor $\Lan_{\mathcal{O}(i_{p})^\opp}\colon \Cat(\tau^\opp,\Set)\to \Cat(\{p\},\Set)=\Set$, whose image at $\mathcal{F}$ is written $\lceil \mathcal{F}_{p}\rceil$.

	The restriction of this functor to the category of sheaves on $X$ can be identified with the \emph{stalk} functor $(-)_{p}$: we have $\mathcal{O}(\{p\})=\{\varnothing\le\{p\}\}$ and computing the images of $\varnothing$ and $\{p\}$ under $\lceil \mathcal{F}_{p}\rceil$ via the colimit formula for left Kan extensions gives
	\begin{align*}
		\lceil \mathcal{F}_{p}\rceil(\{p\}) & \cong \colim\left(\left(\mathcal{O}(\lceil p\rceil)\downarrow\underline{\{p\}}\right)^\opp\twoheadrightarrow\mathcal{O}(X)^\opp\xto {\mathcal{F}}\Set\right),  \\
		                                    & \cong \colim_{U\ni p}(\mathcal{F}(U)),                                                                                                                         \\
		                                    & \cong \mathcal{F}_{p}                                                                                                                                          \\
		\ceil{\mathcal{F}_{p}}(\varnothing) & \cong \colim\left(\left(\mathcal{O}(\ceil{p})\downarrow\underline{\varnothing}\right)^\opp\twoheadrightarrow \mathcal{O}(X)^\opp\xto {\mathcal{F}}\Set\right), \\
		                                    & \cong \colim_{U\hookrightarrow\varnothing}(\mathcal{F}(U)),                                                                                                    \\
		                                    & \cong \mathcal{F}(\varnothing).
	\end{align*}
\end{example}

\begin{example}[Analytic functors]\label{analuo}
	\index{Combinatorial species}
	A functor $F : \Set\to \Set$ is said to be \emph{analytic} if it consists of the left Kan extension of a functor $f : \B(\N)\to \Set$ (the `species' of $F$) along the functor $j : \B(\N)\to \Set$; $\B(\N)$ is the category having objects natural numbers and such that $\B(\N)(m,n)$ are the bijective functions $\{1,\dots,m\}\to \{1,\dots,n\}$ (so this set is empty if $n\neq m$).

	In other words, $\B(\N)$ is the groupoid arising as disjoint union of all symmetric groups $\coprod_{n\ge 0} \text{Sym}(n)$.

	Representing the left Kan extension $\Lan_jf$ as a coend we have
	\[
		F(T)\cong \int^n T^n\times f(n);
	\]
	a functor is `analytic' if it can be expressed as a \emph{Taylor series}, and the coend is in a suitable sense that Taylor series). The theory of analytic functors, besides having an intrinsic interest, is capable to categorify many phenomena in classical combinatorics. See for example the seminal \cite{Joyal1986foncteurs}, but also \cite{adamek2008analytic,gambo-joy}.
\end{example}
\subsection{Tannaka duality using coends}
\begin{example}\label{dualvecspace}\index{Dual space}
	Let $V$ be a finite dimensional vector space over the field $K$; let $V^\lor$ denote the dual vector space of linear maps $V\to K$. Then there is a canonical isomorphism
	\[
		\int^V V^\lor \otimes_K V\cong K.
	\]
	The fastest way to see this is to notice that
	\[
		\int^V \hom(V,\firstblank)\otimes V\cong \Lan_{\id}(\id)\cong \id_{\Mod(K)}
	\]
	(compare this argument with any proof trying to explicitly evaluate the coend from its bare definition).
\end{example}
\index{Tannaka theory|(}
\begin{remark}
	Let again $V$ be a finite dimensional vector space over the field $K$. The universal cowedge $\hom(V,V)\xto{\alpha_V} K$ sends an endomorphism $f : V\to V$ to its \emph{trace} $\tau(f)\in K$ (which in this way acquires a universal property).

	The above argument holds in fact in fair generality, adapting to the case where $V$ is an object of a compact closed monoidal category, and it is linked to the theory of \emph{Tannaka reconstruction}.
\end{remark}
Let $G$ be an internal group in a suitable category of spaces (it can be a Lie group or an affine group, \ie a group in the category of algebraic varieties or schemes). The category $\BF{Rep}(G)$ of its finite dimensional representations $\varrho : G \to \Mod(K)$ carries many important properties: it has a monoidal structure and an involution $(\firstblank)^\lor$ turning it into a \emph{rigid} monoidal category (see \cite{selinger2010survey} for a glimpse on the vast zoology of monoidal categories).

Tannaka theory tries to find sufficient conditions on a nice monoidal category $\A$ ensuring that it is equivalent to $\BF{Rep}(G)$ for some space $G$, and retrieves sufficient information to \emph{reconstruct} such $G$ (or equivalently by Gel'fand duality, its algebra of functions) from $\A$. Such algebra of functions is built out of $\A$ alone in a canonical fashion, as long as it comes equipped with a nice functor $\A \to \Mod(K)$ for some ring $K$.
\begin{remark}\upeyes
	For the purposes of this remark, unveiling the coend\hyp{}y nature of the argument in \cite[10.2.2]{schappi2013formal},
	\begin{itemize}
		\item We fix a ground ring $R$, and we let $K$ be a commutative $R$\hyp{}algebra and $\A$ an additive autonomous symmetric monoidal $R$\hyp{}linear category; this means that every object in $\A$ is \emph{dualisable}.
		\item We let $w : \A \to \Mod(K)$ be a $R$\hyp{}linear functor that is strong monoidal; since $\A$ is autonomous, this entails that the essential image of $w$ is contained in the subcategory of dualisable (\ie finitely generated projective) $B$\hyp{}modules.
		\item We also assume that $w$ is comonadic.
	\end{itemize}
	Under these assumptions, we can consider the \emph{density comonad} of $w$, \ie the left Kan extension of $w$ along itself (of course if $w$ is comonadic, and has a right adjoint $r$, then $\Lan_ww\cong w\circ r$ because $\Lan_w \cong \firstblank \circ r$): according to \ref{kan_are_coends}, the left Kan extension in study can be computed as the coend
	\[
		M \mapsto \int^{A\in\A}	\Mod(K)(wA, M)\otimes wA
	\]
	We now claim that the object $H=\Lan_ww(K)$ carries a natural structure of Hopf algebra in $\Mod(K)$: indeed $H$ results as
	\[ \Lan_ww(K)\cong \int^A \Mod(K)(wA, K)\otimes wA = \int^A (wA)^\lor\otimes wA\]
	(this is exactly how the object $H$ is introduced in \cite{schappi2013formal,bakke2007hopf,ulbrich1990hopf}); moreover, $\Lan_ww$ is strong monoidal by doctrinal adjunction \cite{Kelly1974} and thus $H = \Lan_ww(K)$ must carry a bi\hyp{}algebra structure.
\end{remark}
One form of the \emph{Tannaka reconstruction theorem} now asserts that $\A$ is monoidally equivalent to the category of representations of the space $\Spec(H)$.

Here, we prove a slightly less sophisticated version of the theorem using coends. The proof is obtained in various steps; we just sketch its backbone to let the reader appreciate how, although there's still a lot of (non\hyp{}formal) work left to do, the use of coends makes the essential idea crystal\hyp{}clear.
\begin{theorem}
	Let $K$ be a ring, $F : \A \to \text{mod}(K)$ a $K$\hyp{}linear, faithful, strong monoidal functor with domain a $K$\hyp{}linear rigid monoidal category. The codomain $\text{mod}(K)$ is the category of finitely\hyp{}generated projective $K$\hyp{}modules.

	Then there is a bialgebra $B\in \Mod(K)$ (now, not necessarily finitely generated) such that $\A$ is monoidally equivalent to the category of $B$\hyp{}modules $\Mod(B)$.
\end{theorem}
\begin{proof}
	First, consider the codensity monad $\Ran_FF : \Mod(K) \to \Mod(K)$, and the density comonad $\Lan_FF$ of $F$; the image of $K$ under the monad is the object $\Ran_FF(K)$ and corresponds to the module
	\[
		\int_A \hom_K(\hom_K(k,FA),FA) = \int_A \hom_K(FA,FA)
	\]
	\ie to the monoid of endo\hyp{}transformations $\tCat(\A,\Mod(K))(F,F)$ (see \ref{naturalu}; the monoid operation here is vertical composition, which is of course bilinear). We claim this is the algebra $B$ we are looking for. Indeed, the object $\Lan_FF(K)$ is another $K$\hyp{}module, and not very far from $B$:
	\begin{align*}
		Lan_FF(k) & \cong \int^A \hom_K(FA, k)\otimes FA             \\
		          & \cong \int^A (FA)^*\otimes FA                    \\
		          & \cong \int^A (FA\otimes FA^*)^*                  \\
		          & \cong \left(\int_A FA\otimes FA^*\right)^* = B^*
	\end{align*}
	Now, `every object of $\A$ is a $B$\hyp{}module' in the following sense: the universal wedge of the end $\int_A \hom_K(FA,FA)$ is made by maps
	\[
		\epsilon_A : B \to \hom_K(FA,FA) = End(FA)
	\]
	and this is a ring map giving $FA$ a structure of $B$\hyp{}module; a morphism $f : A\to A'$ in $\A$ now fits in the commutative square
	\[
		\vcenter{\xymatrix{
				B\ar[r]\ar[d] & \hom_K(FA,FA) \ar[d]^{Ff_*}\\
				\hom_K(FA',FA')\ar[r]_{Ff^*} & \hom_K(FA,FA')
			}}
	\]
	which means that for every $b\in B$, $Ff(b.x)=b.Ff(x)$ where $b.\firstblank = \epsilon(b)$; thus it is a homomorphism of $B$\hyp{}modules.

	This is enough to define a functor $\tilde F : \A \to \Mod(B)$ (just corestrict $F$) in such a way that it is an equivalence of categories; it is indeed full and strictly surjective on objects, and strong monoidal and faithful by the initial assumption.

	Moreover, the multiplication of $B$ given by vertical composition of natural transformation is compatible with a comultiplication on $B^*$, and precisely (if we shortly denote $[F,F] = \tCat(\A,\Mod(K))(F,F)$)
	\[
		\begin{array}{c}
			([F,F]\otimes [F,F] \to [F,F])^* \\ \hline
			[F,F]^*\otimes [F,F]^* \leftarrow [F,F]^*
		\end{array}
	\]
	(and $B\cong B^*$ because as $B$\hyp{}module it is of course 1\hyp{}dimensional).
\end{proof}
\index{Tannaka theory|)}
\section{A  Yoneda structure on $\tCat$}\label{def:yoneda-struc}
\index{Yoneda!--- structure}
\subsection{Formal category theory: a crash course}\index{Category!formal --- theory}\index{Formal category theory}
The language of category theory is built upon a certain number of fundamental notions: among these we find the universal characterisation of co\fshyp{}limits, the definition of adjunction, (pointwise) Kan extension, and the theory of monads.

It is often possible to `axiomatise' these definitions, pretending that they refer to the 1- and 2\hyp{}cells of a generic 2\hyp{}category other than $\tCat$.

In some sense, category theory arises when the way in which abstract patterns interact becomes itself an object of study, and when it is generalised to several different contexts. In a few words, the aim of \emph{formal category theory} is to provide a framework in which this process of conceptualisation can be outlined mathematically. Quoting \cite{Gray},
\begin{quote}
	The purpose of category theory is to try to describe certain general aspects of the structure of mathematics. Since category theory is also part of mathematics, this categorical type of description should apply to it as well as to other parts of mathematics.
\end{quote}
The basic idea is that the category of small categories, $\Cat$, can be promoted to a 2\hyp{}category $\tCat$ with `formal' properties in the same way $\Set$ is a category with `formal' properties (leading to the definition of a topos). The aim of formal category theory is to outline these properties, and the assumptions needed to ensure that a certain 2\hyp{}category behaves like $\tCat$ for all practical purposes.

Unfortunately, being too na\"ive when performing this process doesn't always
give the `right' answer (because it doesn't always build an object
with the right universal property).

This is ultimately due to the fact that, when moving to the setting of
$\V$\hyp{}enriched categories (which is the adjacent step of abstraction from
$\tCat$, the category $\Set$-$\Cat$ of $\Set$-enriched categories) the theory `behaves differently' in various ways, and some of these
differences prevent $\V$\hyp{}categories to be as expressive as one would have
liked them to be (a paradigmatic example of this minor expressiveness is the lack of
a \emph{Grothendieck construction} for generic $\V$\hyp{}presheaves: studying the way in which the
Grothendieck construction of \ref{eltsf} ultimately pertains to formal category theory has been
addressed in the early literature on formal category theory (see \cite{StreetFibreYoneda1974,street1978yoneda,street1980fibrations}).

Formal category theory can be thought as a way to encode the same amount of
information carried from $\tCat$ in other contexts: even though it is always possible to do some
constructions by mimicking definitions from $\tCat$ (adjunctions and
adjoint equivalences, extensions by universal 2\hyp{}cells, etc.), things get a
little hairy when we want to provide the theory with an analogue of the Yoneda lemma.

In the 2\hyp{}category $\tCat$, we can use the above mentioned Grothendieck construction to `revert'
set\hyp{}valued functors on an object $\B$ into arrows `over' $\B$;\footnote{We basically glue
	together a bunch of fibers $\coprod_B \E_B$ projecting onto $\B$, in the same
	manner we build the \emph{étale space} of a presheaf $F : \B^\opp \to \Set$; the reader might have noticed that this is secretly the same construction that gives the \emph{étale space} of a presheaf on a topological space (see \ref{etalage}).} in the
2\hyp{}category $\tCat$ the comma object of $B : 1 \to \B$ to $\id_\B : \B \to \B$ together
with its projection $B/\B \to \B$ is a good stand in for the covariant functor
represented by $B$ (more generally, \emph{discrete left fibrations} over $\B$
stand in for general functors $\B \to \Set$).

In the 2\hyp{}category $\V\text{-}\tCat$, we care about $\V$-\emph{valued}
$\V$\hyp{}functors and we would like to do the same construction. But for an
object $B$ in a $\V$\hyp{}enriched category $\B$, the comma $B/\B$ is more naturally
an \emph{internal} category (whose object of objects is $\coprod_{X \in \B}
	\B(B,X)$) rather than an \emph{enriched} one (whose objects are morphisms $p : B \to
	X$ in the underlying category of $\B$). The skewness between these two presentations of category, one `inside' a universe, and the other `with respect to' a universe, generates all sorts of subtleties and problems.

Now, we are left with the question:
\begin{quote}
	Which additional structure on a 2\hyp{}category $\sfK$ allows to recognise arrows of $\sfK$ playing the same r\^ole of discrete (op)fibrations in
	$\tCat$, thus providing with a meaningful notion of (fibrational) Yoneda lemma internal to $\sfK$ (see §\ref{alt_yoneda})?
\end{quote}
The axioms of \emph{Yoneda structure} provide a possible answer to this question.

Our aim here is to present them not for an arbitrary 2\hyp{}category, as they appear in \cite{street1978yoneda} but for the 2\hyp{}category $\tCat$: we prove the validity of each axiom as it appears in \cite{street1978yoneda}; to do this, we will extensively use the co\fshyp{}end calculus we know.

\smallskip
To start, we establish the following notation:
\begin{enumtag}{yd}
	\item $\sfK$ is a 2\hyp{}category, fixed once and for all;
	\item $\adm(A,B) \subseteq \sfK(A,B)$ is a full subcategory of
	`admissible' 1\hyp{}cells, which is moreover a \emph{right ideal}, meaning that the
	composition map restricted to admissible 1\hyp{}cells restricts as a family of maps
	\[ c_{XAB}|_{\adm} : \adm(A,B) \times \sfK(X,A) \to \adm(X,B).
	\] This means that if $\var[f]{X}{A}$ and $\var[g]{A}{B}$ are 1-cells and $g$ is admissible, then $g\circ f$ is again admissible. We call admissible an \emph{object} $A$ such that $\id_A \in \adm(A,A)$.
\end{enumtag}
Now we assume that the following structure can be found on $\sfK$:
\index{_aaa_P@$\bsP$}
\begin{enumtag}{ys}
	\item for each admissible object $A\in\sfK$ we can find an admissible 1\hyp{}cell
	$\yon_A : A \to \bsP A$ called a \emph{Yoneda arrow};\footnote{It is desirable for the correspondence $A\mapsto \bsP A$ to be a functor; it will be an axiom. For the moment, the assignment $\bsP A$ is just any object depending on $A$.}
	\item for each $f : A \to B$ admissible 1\hyp{}cell with admissible domain, we
	can find a 2\hyp{}cell
	\[\label{sopra}
		\vcenter{\xymatrix{
		A \ar[d]_f
		\drtwocell<\omit>{<2>\chi^f}\ar[dr]^{\yon_A}& \\ B \ar[r]_-{B(f,1)}& \bsP A
		}}
	\]
\end{enumtag}
We say that a 2\hyp{}category $\sfK$ \emph{has a Yoneda structure} if it is endowed with the data above, and if the following axioms are satisfied.
\begin{axiom}\label{ysax_uno}
	In \eqref{sopra}, the pair $\langle B(f,1),\chi^f\rangle$ exhibits $\lan_f \yon_A$.
\end{axiom}
\begin{proof*}
	The proof that this axiom holds in $\tCat$ consists of the following derivation in coend calculus: let $f : \A \to \B$ be a functor, then
	\begin{align*}
		\lan_f \yon_\A(B) & \cong \int^A \B(fA,B)\times \yon_\A(A)        \\
		                  & \cong \int^A \B(fA,B)\times \A(\firstblank,A) \\
		\ref{ninjayo}     & \cong \B(f\firstblank,B).
	\end{align*}
\end{proof*}
\begin{axiom}\label{ysax_due}
	In \eqref{sopra}, the pair $\langle f,\chi^f\rangle$ exhibits the absolute left lifting	$\Lift_{B(f,1)}\yon_A$.
\end{axiom}
\begin{proof}
	The validity of this axiom in $\tCat$ is again a game of coend
	calculus: if we call $N_f = \lan_f\yon_\A=\B(f,1)$ for short, we have
	$\leeft_{N_f}\dashv N_{f,*}$, where $N_{f,*}\colon g\mapsto N_f\circ g$ is the
	`direct image' (or \emph{post\fshyp{}composition}) functor; then we have an isomorphism between sets of 2\hyp{}cells
	\begin{align*}
		\tCat(\A,\bsP \A)\big( \yon_\A, N_f\circ g \big) & \cong \int_{A}[\A^\opp,\Set]\big(\yon_\A{A}, N_f\circ g(A)\big)        \\
		                                                 & \cong \int_{A}[\A^\opp,\Set]\big( \yon_\A{A}, \B(f\firstblank,gA)\big) \\
		                                                 & \cong \int_{A}\B(fA,gA)                                                \\
		                                                 & \cong \tCat(\A,\B)(f,g)
	\end{align*}
	\index{Adjunction!relative ---}
	\index{Relative adjunction}
	We leave to the reader the proof that this lifting is preserved by every functor; a terse equivalent formulation of this axiom is that there is a relative adjunction (see Exercise \ref{ex:reladjs}) $f \adjunct{[\yon_A]}{} B(f,1)$ with \emph{relative unit} the Yoneda map $\yon_A$.
\end{proof}
\begin{axiom}\label{ysax_tre}
	Given a pair of composable 1\hyp{}cells $A \xto{f} B\xto{g} C$, the
	pasting of 2\hyp{}cells
	\[
		\vcenter{\xymatrix{
		A\ar@{}[dr]|(.6){\Swarrow\chi^{\yon_B f}}\ar[rr]^{\yon_A}\ar[d]_f && \bsP A\\
		B\ar@{}[dr]|(.35){\Swarrow\chi^g} \ar[d]_g \ar[r]|-{\yon_B}& \bsP X\ar[ur]_{\bsP f} \\
		C\ar[ur]_{C(g,1)} &&
		}}
	\]
	exhibits the extension $\lan_{gf}\yon_A = C(gf,1)$, and the pair $\langle \id_{\bsP A}, \id_{\yon_A}\rangle$ exhibits $\lan_{\yon_A}\yon_A$.
\end{axiom}
\begin{remark}
	The hidden meaning of this axiom is that $\bsP$ is a pseudofunctor $\sfK^\coop \to \sfK$ (whose domain is the sub-class of admissible objects).

	Let's make this evident: given a pair of composable 1\hyp{}cells $A \xto{f} B\xto{g} C$, the universal property of
	$\chi^{gf}$ entails that there is a unique 2\hyp{}cell $\theta^{gf}$ filling the
	diagram
	\[
		\vcenter{\xymatrix@C=2cm{
		A \drtwocell<\omit>{\chi^{gf}} \ar[r]^{\yon_A} \ar[d]_f & \bsP A \\
		B \ar[d]_g & \bsP B \ar[u]_{\bsP f}\\
		C \ar[ur]\ar[uur]
		\uurtwocell<\omit>{<2>\kern1em\theta^{gf}}
		\urtwocell<\omit>{<2>}
		\ar[r]_{\yon_C} & \bsP C \ar[u]_{\bsP g}
		}}
	\] Axiom \ref{ysax_tre} is equivalent to the request that this arrow is invertible (exercise:
	draw the right diagram), and this yields that the above diagram has the same
	universal property of the square
	\[ \vcenter{\xymatrix{ A \ar@{}[dr]|\Sarrow \ar[r]^{\yon_A}\ar[d]_{gf}& \bsP A\\ C
				\ar[r]_{\yon_C}& \bsP C\ar[u]_{\bsP(gf)} }}
	\] which in turn entails that there is a unique, and invertible, 2\hyp{}cell $\bsP(gf)
		\To \bsP f \circ \bsP g$. This is of course half of the structure of
	pseudofunctor on $\bsP$; the remaining structure is given by the request that
	$\langle \id_{\bsP A}, \id_{\yon_A}\rangle$ exhibits $\lan_{\yon_A}\yon_A$, because this provides a unique, invertible 2-cell $\id_{\bsP A} \To \bsP(\id_A)$ for every admissible object $A$. This result has a variety of different interpretation: it follows from the Yoneda lemma as stated in \ref{thm:yoda-is-dense}, and from the coend expression for the density comonad, or by a direct check of the universal property in study.
\end{remark}
All the remarks in \ref{ysax_tre} are evidently true in $\tCat$, since the (pseudo) functoriality of
the correspondence $\A \mapsto [\A^\opp,\Set]$ can be proved directly without great effort. Nevertheless,
axiom \ref{ysax_tre} is still telling us something about a `reduction rule' for composition
of Kan extensions: indeed, it is possible to prove that (in the same notation of
axiom \ref{ysax_tre}) that there is a (canonical!) isomorphism
\[ \theta_{gf} : \lan_{gf}\yon_\A \cong \lan_{\yon_\B f}\yon_\A \circ \lan_g \yon_\B.
\]
The proof of this statement is another game of coends, using \ref{kan_are_coends}: try to do it as an exercise (hint: there will be many coends involved; fix an explicative notation and maintain it clear).
\subsection{Addendum: co\fshyp{}ends inside a Yoneda structure}
\begin{warning}
	The present section relies on material presented in chapter 4 and 5; it can (and should) be skipped at first reading, but can be used as reference for additional examples of weighted co\fshyp{}limits coming back from said chapters.
\end{warning}

Along the following remark we fix a 2\hyp{}category $\sfK$ and we assume that $\sfK$ has all finite limits, that it is cartesian closed, and that the functor $\bsP$ is \emph{quasi\hyp{}representable}, \ie it arises as $A \mapsto [A^\lor, \Omega]$ for some object $\Omega \in \sfK$ and a \emph{duality involution} $(\firstblank)^\lor$ on $\sfK$ (see \cite{shulmano2016contra}); it is possible to prove that there is an isomorphism $\Omega \cong \bsP 1$, where $1$ is the terminal object of $\sfK$.

The purpose of this section is to establish how the Yoneda structure having such $\bsP$ as presheaf construction possesses a formal analogue of co\fshyp{}end calculus.
\begin{remark}\label{for-enrich}\upeyes
	In a quasi-representable Yoneda structure $\bsP : \adm \to \sfK$, Yoneda embeddings are maps in $\sfK$ of the form $A \to [A^\lor, \bsP 1]$; if the product $A^\lor\times A$ is admissible in the Yoneda structure generated by $\bsP$ we can consider the admissible maps
	\[
		\fh_A : A^\lor\times A\to \bsP 1
	\]
	as the mates of the Yoneda embeddings $\yon_A$. These maps play the r\^ole of \emph{internal homs}, so that admissible objects can be thought as canonically $\bsP1$\emph{\hyp{}enriched}: in order to define co\fshyp{}ends we shall write them as certain weighted co\fshyp{}limits using the maps $\fh_A$ as weights (see \ref{weicolims}.\ref{coends_are_hom_weighted}, co\fshyp{}ends are precisely hom\hyp{}weighted co\fshyp{}limits).\footnote{It is worth to remind that fr every small category $A\in\tCat$ the composition maps are given by a family of functions $\text{c}_{abc} : A(a,b) \times A(b,c) \to A(a,c)$ such that $c_{a,\firstblank,c}$ is a cowedge, and in fact an initial one.}

	Let's stick for a moment to the case where the ambient 2-category is $\tCat$; by definition of what is a weighted colimit in a Yoneda structure \cite[§4]{street1978yoneda}, the presence of an isomorphism
	\[\int^X \A(A,X)\times \A(X,A')\cong \A(A,A')\]
	valid for a category $\A$ and objects $A,A'\in A$ means that the left lifting of $\yon_{\A^\opp\times \A}$ along $\bsP 1(\fh_\A,1)$, is $\fh_\A$, and such lifting is absolute; this yields a relative adjunction $\fh_\A \adjunct{[\yon_{\A^\opp\times \A}]}{} \bsP 1(\fh_\A,1)$; this final request means that the triangle
	\[
		\vcenter{\xymatrix{
		& \A^\opp\times \A
		\ar@{}[d]|{\overset{\chi^{\fh_\A}}\To}\ar[ld]_{\yon_{\A^\opp\times \A}} \ar[rd]^{\fh_\A} &  \\
		\bsP(\A^\opp\times \A) & {} & \bsP 1 \ar[ll]^{\bsP 1(\fh_\A,1)}
		}}
	\]
	is an absolute left lifting; this is precisely axiom \ref{ysax_due} applied to $\fh_\A$.

	It is worth to investigate how far this formalisation of co\fshyp{}end calculus can go:
\end{remark}
\begin{remark}[Co\fshyp{}ends in a Yoneda structure, \upeyes]
	The 1\hyp{}cell $\bsP 1(\fh_A,1)$ admits a left adjoint $\Lan_{\yon_{A^\lor\times A}} \fh_A$; computing this Kan extension in $\tCat$ we get
	\begin{align*}
		\Lan_{\yon_{\A^\opp\times \A}} \fh_\A(F) & \cong \int^{(A,A')} [\A^\opp\times \A,\Set](\yon_{\A^\opp\times \A}(A,A'), F)\times \fh_\A(A,A') \\
		                                         & \cong \int^{(A,A')} F(A,A')\times \A(A,A')                                                       \\
		                                         & \cong \int^a F(A,A).
	\end{align*}
	Thus, the adjunction $\Lan_{\yon_{A^\lor\times A}} \fh_A\dashv \bsP 1(\fh_A,1)$ is an abstract analogue of the adjunction $\int^A : [\A^\opp\times \A, \Set] \leftrightarrows \Set : \hom\pitchfork\firstblank$, where the left adjoint $\int^A$ takes the coend of a functor $F : A^\opp\times A\to \Set$, and the right adjoint is the hom\hyp{}weighted limit functor $\hom\pitchfork X : (A,A')\mapsto X^{\hom(A,A')}$, the same that helped prove the Fubini theorem in \ref{fubozzo} (see \ref{weicolims}.\ref{coends_are_hom_weighted} for the precise result).

	\smallskip
	It is now straightforward to go further: axiom \ref{ysax_due} entails that in a Yoneda structure an admissible object $A$ for which $A\times A^\opp$ is still admissible (this translates into the property that the domain of a \index{Monad!relative ---}\index{Relative monad} relative monad $\bsP : \sfA \to \sfK$ is closed under product and duality involution), thus there are absolute left liftings
	\[
		\vcenter{\xymatrix{
		& A^\opp\times A \ar[ld]_{\yon_{A^\opp\times A}} \ar[rd]^{\fh_A} &  \\
		\bsP(A^\opp\times A) &\utwocell<\omit>{} & \bsP 1 \ar[ll]^{\bsP 1(\fh_A,1)}
		}}
	\]
	in which $\bsP 1(\fh_A, 1)$ has a left adjoint, precisely $\Lan_{\yon_{A^\lor\times A}}\fh_A$. Such a left adjoint is the functor taking the `coend' of an internal endo\hyp{}profunctor $F : A^\lor\times A\to \bsP 1$.
\end{remark}
\section{Addendum: relative monads}\label{app:relmo}
\begin{notation}
	Along the present section, categories of functors $\cX \to \cY$ are often denoted as $[\cX,\cY]$; this complies with the notation for cartesian closed categories and serves to avoid cumbersome accumulation of symbols when iterated functor categories such as
	\[
		\Cat(\Cat(\Cat(A^\opp,\Set)_s^\opp, \Set)_s^\opp,\Set)
	\]
	are considered. Context always allows to determine where the hom\hyp{}category $[\cX, \cY]$ lies.

	Another notational simplification is the following: every functor $J : \clX \to \clY$ induces a left extension functor $\Lan_J = J_! : [\clX, \clY] \to [\clY, \clY]$. This notation draws from algebraic geometry and it is chosen since we have to iterate several left extensions and to iteratively apply the functors $J_!\dashv J^*$ and compose their unit and counit maps.
\end{notation}
In this subsection address a deeper study of the presheaf construction: we shall prove that
\begin{itemize}
	\item albeit not an endofunctor, the presheaf construction functor $\bsP$ of the Yoneda structure on $\Cat$ is a \emph{relative monad} on the category of functors $\caat \to \tCat$; this means that it is an internal monoid in the skew\hyp{}monoidal category of those functors; in simple terms, a skew\hyp{}monoidal category $(\M, \triangleleft)$ is like a monoidal category, but the associator $\alpha_{XYZ} : (X\triangleleft Y)\triangleleft Z \to X\triangleleft (Y\triangleleft Z)$ and left and right unitors are not invertible maps: thus, aggregations of objects
	      \[X_0 \triangleleft X_1\triangleleft\cdots\triangleleft X_n\]
	      are not well\hyp{}defined without specifying their parenthesisation. We introduce the notion in full generality, and we prove that, under suitable cocompleteness assumptions on $\cY$, the category $\Cat(\cX, \cY)$ becomes skew\hyp{}monoidal.
	\item We then prove that $\bsP$ is a particularly well\hyp{}behaved relative monad, that in \cite{yosegi} is called a \emph{yosegi box}. In simple terms, the action of $\bsP$ on 1\hyp{}cells is such that
	      \begin{enumtag}{yb}
		      \item \label{yb:uno} every $\bsP_! f$ fits into an adjunction $\bsP_! f\dashv \bsP^* f$;
		      \item \label{yb:due} $\bsP_!$ is a relative monad with respect to the inclusion $j : \caat\subset\tCat$;
		      \item \label{yb:tre} the monad $\bsP_!$ is lax idempotent.
	      \end{enumtag}
	      As shown in \cite{yosegi}, it turns out that these three properties characterise uniquely the presheaf construction of a Yoneda structure.
\end{itemize}
The upshot of the present section is thus the following: provided the left extension along a given functor $J : \cX\to \cY$ exists, the functor category $[\cX,\cY]$ becomes a \emph{skew\hyp{}monoidal} \cite{szlachanyi2012skew} with respect to the functor $\triangleleft : [\cX,\cY]\times[\cX,\cY]\to [\cX,\cY]$ sending $F,G$ to $J_!F\circ G$.

In the present section we offer a characterisation of relative monads as $\triangleleft$\hyp{}monoids: we can provide a formal analogue for a similar statement than the one given for $\sfK = \tCat$ in \cite{altenkirch2010monads}; while extremely useful a reference, the proof in \cite{altenkirch2010monads} does apparently work only on $\tCat$, whereas the argument appearing here can be easily adapted to an abstract 2\hyp{}category.

Although a relatively elementary argument, spelling out a complete proof of this fact turns out to be a rather tedious task.

\begin{definition}\label{it-is-skiu}
	\index{Monad!relative ---}\index{Relative monad}
	Let $J : \cX \to \cY$ be a functor such that the left extension along $J$ exists, defining a functor $J_! : [\cX, \cY]\to [\cY,\cY]$; then the category $[\cX,\cY]$ becomes a \emph{left skew\hyp{}monoidal category} under the skew multiplication defined by
	\[
		(F,G) \mapsto F\triangleleft G = {J_!}F\circ G;
	\]
	there are natural maps of \emph{association}, \emph{left} and \emph{right unit}
	\[\label{coherence-data}
		\begin{tikzcd}[column sep=large, row sep=0mm]
			(F\triangleleft G)\triangleleft H \ar[r, "\gamma_{FGH}"] & F\triangleleft(G\triangleleft H)\\
			J\triangleleft F \ar[r, "\lambda_F"] & F\\
			F \ar[r, "\varrho_F"] & F\triangleleft J \\
		\end{tikzcd}
	\]
	such that the following diagrams are commutative:
	\begin{enumtag}{skm}
		\item \label{skm:uno}\emph{skew associativity}:
		\[
			\begin{tikzcd}[column sep=2cm]
				((F\triangleleft G)\triangleleft H)\triangleleft K \ar[r,"\gamma_{FG,H,K}"]\ar[d, "\gamma_{F,G,H}\triangleleft K"']& (F\triangleleft G)\triangleleft (H\triangleleft K) \ar[dd, "\gamma_{F,G,HK}"]\\
				(F\triangleleft (G\triangleleft H))\triangleleft K \ar[d, "\gamma_{F,GH,K}"']& {}\\
				F\triangleleft ((G\triangleleft H)\triangleleft K) \ar[r, "F \triangleleft \gamma_{G,H,K}"'] & F\triangleleft (G\triangleleft (H\triangleleft K))
			\end{tikzcd}
		\]
		\item \label{skm:due}\emph{skew left and right unit}:
		\[
			\begin{tikzcd}[column sep=2cm]
				(F\triangleleft G)\triangleleft J\ar[r, "\gamma_{F,G,J}"] \celtag[pos=.1,dr]{\textsc{2r}}& F\triangleleft (G\triangleleft J)\\
				F\triangleleft G\ar[u, "\varrho_{F\triangleleft G}"]\ar[ur, "F\triangleleft \varrho_G"'] \celtag[pos=.9,dr]{\textsc{2l}} & F\triangleleft G\\
				(J\triangleleft F)\triangleleft G \ar[r, "\gamma_{J,F,G}"']\ar[ur, "\lambda_F\triangleleft G"]& J\triangleleft(F\triangleleft G)\ar[u, "\lambda_{F\triangleleft G}"']
			\end{tikzcd}
		\]
		\item \label{skm:tre}\emph{zig\hyp{}zag identity}:
		\[
			\begin{tikzcd}[column sep=small]
				& J\triangleleft J \ar[dr, "\lambda_J"]& \\
				J \ar[ur, "\varrho_J"]\ar[rr,equal]&& J
			\end{tikzcd}
		\]
		\item \label{skm:qtr}\emph{interpolated zig\hyp{}zag identity}:
		\[
			\begin{tikzcd}[column sep=large]
				(F\triangleleft J)\triangleleft G \ar[rr, "\gamma_{F,J,G}"]&& F\triangleleft (J\triangleleft G) \ar[d, "F \triangleleft \lambda_G"]\\
				F\triangleleft G \ar[rr,equal] \ar[u, "\varrho_F\triangleleft G"]&& F \triangleleft G
			\end{tikzcd}
		\]
	\end{enumtag}
\end{definition}
The natural maps $\gamma,\lambda,\varrho$ are defined as follows in this specific case:
\begin{enumtag}{s}
	\item \label{s:uno} Given $F,G,H\in [\cX, \cY]$, the cell $\gamma_{F,G,H}$ is defined by $\tilde \gamma_{F,G} * H$, where $\tilde\gamma_{F,G}$ is obtained as the mate of $J_!F * \eta_G$ under the adjunction $J_!\adjunct{\epsilon}{\eta} J^*$: the arrow
	\[\label{associatore} J_!F\circ G \xrightarrow{J_!F * \eta_G} J_!F\circ J^*J_!G  = J^*(J_!F\circ J_!G)\]
	mates to a map $J_!(J_!F\circ G) \longrightarrow J_!F\circ J_!G$.
	\item \label{s:due} Given $F\in [\cX, \cY]$, the cell $\lambda_F : J\triangleleft F \To F$ is defined by the whiskering $\sigma * F$, where $\sigma = \sigma_{\id_{\cY}}$ is the counit of the density comonad of $J$;
	\item \label{s:tre} Given $G\in [\cX,\cY]$, the cell $\varrho_G : G \to G\triangleleft J$ is the $G$\hyp{}component $\eta_G$ of the unit of the adjunction ${J_!} \dashv J^*$
\end{enumtag}
\begin{remark}\label{we-re-better-than-uustalu}
	A complete proof of \ref{it-is-skiu} appears as Theorem 3.1 in \cite{altenkirch2010monads}; the main argument is however heavily relying on the assumption that $\sfK = \tCat$; it was our desire to produce a formal proof \emph{ex novo}, while explicitly recording some useful equations satisfied by the skew monoidal structure in study; as a rule of thumb, the structure maps of this skew monoidal structure are entirely induced by the $J_!\adjunct{\eta}{\varepsilon} J^*$ adjunction and from (the equations satisfied by) the unit and counit thereof.

	Until the end of the section, we adopt the following notation, and employ the following equations wherever needed:
	\begin{enumtag}{e}
		\item \label{e:zer} we denote $\ULbullet{\varpi} : J_!U \to V$ the \emph{mate} of $\varpi : U \to J^*V$ under the adjunction $J_! \adjunct{\eta}{\varepsilon} J^*$; similarly, we denote $\chi^\bullet : U \to J^*V$ the mate of $\chi : J_!U \to V$; in this notation, the bijection
		\[[\cX, \cY](U, J^*V)\cong [\cY, \cY](J_!U, V)\]
		reads as $(\ULbullet{\varpi})^\bullet = \varpi$ and ${}^\bullet(\kern-.1em\chi^\bullet)=\chi$.
		\item \label{e:uno} the first zig\hyp{}zag identity for the adjunction $J_! \adjunct{\eta}{\varepsilon} J^*$ is $(\varepsilon_B * J)\circ \eta_{BJ} = \id_{BJ}$; in particular, if $B = \id_{\cY}$ we have $(\sigma * J)\circ \eta_J = \id_J$;
		\item \label{e:due} the second zig\hyp{}zag identity for the adjunction $J_! \adjunct{\eta}{\varepsilon} J^*$ is $\varepsilon_{J_!F} \circ J_!(\eta_F) =\id_{J_!F}$;
		\item \label{e:ter} an irreducible expansion for the associator, expliciting all its components, is
		\[\label{associator}\gamma_{FGH} = \big(\varepsilon_{J_!F\circ J_!G}\circ J_!(J_!F * \eta_G)\big) * H\]
		\item \label{e:qtr} the density comonad of $J$, whose comultiplication $\nu$ and counit $\sigma$ satisfy the comonad axioms is defined by the maps
		\begin{gather}\sigma = \varepsilon_{\id_B} : J_!(J)\to 1 \notag \\ \nu : J_!(J) \xto{J_!(\eta_J)} J_!(J_!(J)\circ J) \xto{\tilde\gamma_{JJ}} J_!(J)\circ J_!(J)
		\end{gather}
		in particular, the comultiplication and the associator of the skew\hyp{}monoidal structure determine each other.
	\end{enumtag}
\end{remark}
\begin{proof}
	First of all, axiom \ref{skm:tre} is one of the two triangle identities of the adjunction ${J_!} \dashv J^*$. It remains to prove the other coherence conditions.
	\begin{enumtag}{skm}
		\item It is easy to see that one can prove commutativity before precomposing with $K$; one has then to prove the commutativity of the diagram
		\[
			\begin{tikzcd}[column sep=4cm]
				J_!(J_!(J_!F\circ G)\circ H) \ar[r, "\tilde\gamma_{F\triangleleft G,H}"]\ar[d, "J_!(\tilde\gamma_{FG} * H)"'] & J_!(J_!F\circ G)\circ J_!H \ar[dd, "\tilde\gamma_{FG} * J_!H"]\\
				J_!(J_!F\circ (J_!G\circ H)) \ar[d, "\tilde\gamma_{F, G\triangleleft H}"'] \ar[dr,dashed] \celtag[ur]{\text{\textsc{c}\oldstylenums{1}}}& {}\\
				J_!F \circ J_!(J_!G\circ H) \celtag[near start,ur]{\text{\textsc{c}\oldstylenums{2}}}\ar[r, "J_!F * \tilde\gamma_{GH}"']& J_!F\circ J_!G\circ J_!H
			\end{tikzcd}
		\]
		note that the dashed arrow exists: it is the mate $\ULbullet\alpha$ of $\alpha = (J_!F \circ J_!G) * \eta_H$; the plan is to prove that the two sub\hyp{}diagrams in which this arrow splits the whole diagram commute separately. In order to do this, we start from the square diagram \textsc{c}\oldstylenums{1}: recall that \ref{e:uno} entails that one of the squares
		\[
			\begin{tikzcd}
				J_!A \ar[d, "J_!f"'] \ar[r,"\ULbullet a"]& X\ar[d, "g"] & A \ar[d, "f"'] \ar[r,"a"]& J^* X\ar[d, "J^*g"] \\
				J_!B \ar[r, "\ULbullet b"'] & Y & J_!B \ar[r, "b"'] & J^*Y
			\end{tikzcd}
		\]
		commutes if and only if the other does. Hence, if we denote $f=\tilde\gamma_{FG} * H,g=\tilde\gamma_{FG} * J_!H,a=(J_!F\circ G) * \eta_H,b=\alpha$, \textsc{c}\oldstylenums{1} commutes if and only if the square
		\[
			\begin{tikzcd}[column sep=4cm]
				J_!F \circ G\circ H \ar[d, "\ULbullet(J_!F * \eta_G) * H"']\ar[r, "(J_!F\circ G) * \eta_H"]& J_!(J_!F \circ G)\circ J_!H\circ J\ar[d, "\ULbullet(J_!F * \eta_G) * (J_!H \circ J)"]\\
				J_!F \circ J_!G \circ H \ar[r, "(J_!F\circ J_!G) * \eta_H"']& J_!F\circ J_!G \circ J_!H \circ J
			\end{tikzcd}
		\]
		commutes. It does, since their common value at the diagonal is simply the horizontal composition $\ULbullet{(J_!F * \eta_G)} \boxminus \eta_H$.

		A similar argument shows that \textsc{c}\oldstylenums{2} commutes: we have to establish the commutativity of
		\[\begin{tikzcd}
				J_!(J_!F\circ J_!G\circ H) \ar[dr, "\ULbullet(J_!F * \eta_{G\triangleleft H})"']\ar[rr, "\ULbullet(J_!F * J_!G * \eta_H)"]&& J_!F\circ J_!G\circ J_!H \\
				& J_!F\circ J_!(J_!G\circ H)\ar[ur, "J_!F * \ULbullet(J_!G * \eta_H)"'] &
			\end{tikzcd}\]
		a diagram that can be `straightened' to the left one below:
		\[
			\begin{tikzcd}[column sep=1.75cm]
				A \ar[r, "\ULbullet(J_!F * \eta_{G\triangleleft H})"] \ar[d,equal] & A \ar[d, "J_!F * \ULbullet(J_!G * \eta_H)"] && A \ar[d, equal]\ar[r, "J_!F * \eta_{G\triangleleft H}"] & A\ar[d, "J_!F * \ULbullet(J_!G * \eta_H) * J"]\\
				A \ar[r, "\ULbullet(J_!F * J_!G * \eta_H)"'] & A && A\ar[r, "J_!F * J_!G * \eta_H"'] & A
			\end{tikzcd}
		\] but now, the left diagram commutes if and only if the right one does; and the latter commutativity follows from the definition of $\tilde\gamma$.
		\item A separate argument works for diagrams \textsc{2r} and \textsc{2l}: unwinding the definitions, the commutativity of \textsc{2r} amounts to the commutativity of
		\[
			\begin{tikzcd}
				J_!(J_!F\circ G)\circ J \ar[r, "\gamma_{FGJ}"]\ar[d, "\eta_{J_!F\circ G}"']& J_!F\circ J_!G \circ J \ar[d,equal]\\
				J_!F\circ G \ar[r, "J_!F * \eta_G"]& J_!F \circ J_!G\circ J.
			\end{tikzcd}
		\]
		If we denote for short $J_!F * \eta_G = \varpi$, this commutativity is equivalent to the fact that $J^*(\ULbullet{\varpi})\circ \eta_{J_!F\circ G} = \varpi$, which is true since the left hand side of this equation is $(\ULbullet{\varpi})^\bullet$. For axiom \textsc{2l}, the commutativity of
		\[
			\begin{tikzcd}
				J_!(J)\circ F \ar[r, ""] & J^*J_!F\\
				J_!(J)\circ F \ar[u,equal]\ar[r, "J_!F * \eta_F"']& J_!(J)\circ J^*J_!F \ar[u, "\varepsilon * J_!F"']
			\end{tikzcd}
		\] follows from the fact that the upper horizontal row coincides with the horizontal composition $\sigma \boxminus \eta_F = \eta_F \circ (\sigma * F)$; this means that \textsc{2l} is true if and only if the square
		\[
			\begin{tikzcd}
				F \ar[r,"\eta_F"]& J^*J_!F \\
				J_!J \circ F \ar[u,"\sigma * F"] \ar[r, "J_!F * \eta_F"'] & J_!J \circ J^*J_!F\ar[u, "\varepsilon * J_!F"']
			\end{tikzcd}\]
		commutes; this is obvious by the naturality property of $\eta$.
		\item[\textsc{skm}\oldstylenums{4})] Unwinding the definition, axiom \ref{skm:qtr} asks the diagram
		\[
			\begin{tikzcd}
				J_!F \circ J_!(J) \ar[r, "J_!F * \sigma"]& J_!F \\
				J_!(J_!F\circ J)\ar[u, "\gamma_{FJ}"] & J_!F\ar[u,equal]\ar[l, "J_!(\eta_F)"]
			\end{tikzcd}\]
		to commute. In order to see that it does, we observe that the chain of equivalences
		\begin{align*}
			\underline{(J_!F * \sigma) \circ \varepsilon_{J_!F \circ J_!J}} \circ J_!(J_!F * \eta_J) & = \varepsilon_{J_!F}\circ \underline{J_!(J_!F * \sigma * J) \circ J_!(J_!F * \eta_J)}  \\
			                                                                                         & = \varepsilon_{J_!F} \circ \underline{J_!\big(J_!F * ((\sigma * J) \circ \eta_J)\big)} \\
			                                                                                         & = \varepsilon_{J_!F} \circ \id_{J_!F \circ J}
		\end{align*}
		holds, where the last step is motivated by \ref{e:uno}. So, the composition $(J_!F * \sigma)\circ \tilde\gamma_{FJ}$ equals $\varepsilon_{J_!F}$; this, together with \ref{e:due}, concludes.\qedhere
	\end{enumtag}
\end{proof}
\begin{remark}[A nice $J$ gives a nice structure]\label{adjoints-are-exts}
	In favourable cases the skew monoidal structure simplifies until it collapses to a straight monoidal structure:
	\begin{itemize}
		\item if $J$ is fully faithful, the unit $\eta_F : F \to J_!F \circ J$ is invertible for every $F$, so that $F \cong F \triangleleft J$;
		\item if $J$ is dense, the density comonad $\sigma$ is isomorphic to the identity functor, so $\lambda_F : F \cong J\triangleleft F$ is an isomorphism;
		\item if each extension $J_!F$ is absolute, then each $\tilde\gamma_{F,G}$ is invertible.
	\end{itemize}
\end{remark}
\begin{remark}
	The bifunctoriality of composition has been implicitly employed in the above proof; we spell it out explicitly for future reference.

	In order to prove such bifunctoriality, we have to show the commutativity of the square
	\[\begin{tikzcd}
			F\triangleleft G\ar[r, "F\triangleleft g"]\ar[d, "f\triangleleft G"'] & F \triangleleft K\ar[d, "f\triangleleft K"] \\
			H \triangleleft G \ar[r, "H\triangleleft g"']& H\triangleleft K
		\end{tikzcd}\]
	given $f : F\to H$ and $g : G\to K$. In fact, unwinding the definition it's easy to realize that this diagram commutes since its diagonal coincides with the horizontal composition $J_!f\boxminus g$.
\end{remark}
\begin{remark}[A note on $\triangleleft$\hyp{}whiskering]\label{on-uisge}
	This begs the question of how the whiskering of a 1\hyp{}cell $F$ with a 2\hyp{}cell $\alpha$ works, on the left and on the right; given the shape of the $\triangleleft$-skew\hyp{}monoidal structure, it turns out that
	\begin{gather}
		\alpha \triangleleft F := J_!\alpha * F\notag \\
		F \triangleleft \alpha := J_!F * \alpha
	\end{gather}
\end{remark}
\begin{definition}[Monads need not be endofunctors, but are always skew monoids]\label{def:relmo}
	We will be interested in the notion of a $J$-\emph{relative monad}, or simply a relative monad when $J$ is understood from the context; $J$\hyp{}relative monads are defined as internal monoids in the skew\hyp{}monoidal category $([\cX,\cY], \triangleleft)$, thus they come equipped with a \emph{unit} $\eta : J \To T$ and a multiplication $\mu : T\triangleleft T\To T$. However, since the monoidal structure $\triangleleft$ satisfies pretty asymmetrical coherence conditions, the commutativity conditions satisfied by a relative monad get altered accordingly. In particular,
	\begin{enumtag}{rm}
		\item the unit axiom amounts to the commutativity of
		\[\label{relmo:unit}
			\begin{tikzcd}
				J\triangleleft T \ar[r, "\eta\triangleleft T"] \ar[rd, "\lambda"'] & T\triangleleft T \ar[d, "\mu"] & T\triangleleft J \ar[l, "T\triangleleft \eta"'] \\
				& T \ar[ru, "\varrho"'] &
			\end{tikzcd}
		\]
		where the right triangle `commutes' in the sense that the composition $\mu\circ (T\triangleleft \eta)\circ \varrho$ makes the identity of $T$.
		\item the multiplication $\mu$ is `skew associative':
		\[\label{relmo:mult}
			\begin{tikzcd}
				(T\triangleleft T)\triangleleft T \ar[d, "\mu\triangleleft T"'] \ar[r, "\gamma"] & T\triangleleft(T\triangleleft T) \ar[r, "T\triangleleft \mu"] & T\triangleleft T \ar[d, "\mu"] \\
				T\triangleleft T \ar[rr, "\mu"'] &  & T
			\end{tikzcd}
		\]
	\end{enumtag}
\end{definition}
We are particularly interested in \emph{lax idempotent} monads, \ie as those that satisfy one of the following equivalent properties:
\begin{definition}\label{lax-equivs}
	Let $\bsT : \mathsf{A} \to \mathsf{B}$ be a 2\hyp{}monad between 2-categories; we say that $\bsT$ is \emph{lax idempotent} if one of the following equivalent conditions is satisfied:
	\begin{enumtag}{l}
		\item \label{l:uno} for every pair of $  \bsT $\hyp{}algebras $(a,A),(b,B)$ and morphism $f :A \to B$, the square
		\[\notag
			\begin{tikzcd}[row sep=8mm,column sep=8mm]
				\bsT  A \ar[d, "a"']\ar[dr, phantom, "\Swarrow"] \ar[r," \bsT f"] &   \bsT  B \ar[d, "b"]\\
				A \ar[r, "f"'] & B
			\end{tikzcd}
		\]
		is filled by a unique 2\hyp{}cell $\bar f : b \circ \bsT f\To f\circ a $ which is a lax morphism of algebras;
		\item \label{l:due} there exists an adjunction $a\dashv \eta_A$ with invertible counit;
		\item \label{l:ter} there exists an adjunction $\mu\dashv \eta *  \bsT $ with invertible counit;
		\item there is a modification $\Delta :  \bsT  * \eta \To \eta *  \bsT $ such that $\Delta * \eta = 1$ and $\mu * \Delta = 1$.
	\end{enumtag}
	The conditions for colax algebras are of course obtained replacing lax algebra morphisms with \emph{colax} ones.
	\begin{enumtag}{cx}
		\item \label{cx:uno} for every pair of $  \bsT $\hyp{}algebras $a,b$ and morphism $f :A \to B$, the square
		\[\notag
			\begin{tikzcd}[row sep=8mm,column sep=8mm]
				\bsT  A \ar[d, "a"']\ar[dr, phantom, "\Nearrow"] \ar[r," \bsT f"] &   \bsT  B \ar[d, "b"]\\
				A \ar[r, "f"'] & B
			\end{tikzcd}
		\]
		is filled by a unique 2\hyp{}cell $\bar f : f\circ a \To b \circ \bsT f$ which is a colax morphism of algebras;
		\item \label{cx:due} there exists an adjunction $\eta_A \dashv a$ with invertible unit;
		\item \label{cx:ter} there exists an adjunction $\eta  *   \bsT \dashv \mu$ with invertible unit;
		\item there is a modification $\Upsilon : \eta *  \bsT  \To  \bsT  * \eta$ such that $\Upsilon * \eta = 1$ and $\mu * \Upsilon = 1$.
	\end{enumtag}
\end{definition}
The scope of the next section is to prove that the presheaf construction satisfies the axioms of a lax idempotent monad.
\subsection{Relative monads and presheaves}
\index{_aaa_P@$\bsP$}
\begin{warning}\label{thewarning}
	Throughout the section we will often blur the distinction between two choices of notation: if the typeface for a 2-category is $\mathsf{A}$, objects, \ie 0-cells in $\mathsf{A}$, are denoted as Roman letters $A,B,\dots,X,Y,\dots$, 1-cell are lowercase Roman $f,g,\dots$ this applies in particular to the case of $\tCat$ and its objects, categories $A,B,C,\dots$; we invite the reader to keep in mind this small clash of notation when comparing the section with the rest of the book; a similar choice has been made -without further mention- discussing particular shapes of 2-limits in chapter 4 (see \ref{inserters} and \ref{commae_obj}).
\end{warning}
In the present section we study the pair $(\tCat, \bsP)$, where $\bsP : \caat\to \tCat$ is the presheaf construction sending $A\mapsto\psh{A}$; this is defined having domain the category $\caat$ of small categories, and codomain the locally small ones; the embedding of $\caat$ into $\tCat$ will always be denoted as $j : \caat\subset\tCat$: it is an inclusion at the level of all cells.

We fix a notation that can be easily generalized to the case of a pair $(\sfK,\bsP)$, where $\bsP : \sfA\to \sfK$ is the presheaf construction of a Yoneda structure on $\sfK$.
\begin{notation}[Presheaves]\label{presh-notat}
	We consider the functor $\bsP : A\mapsto \psh{A}$ as a \emph{covariant} correspondence on functors and natural transformations; more formally, $\bsP$ acts as a correspondence $\caat \to \tCat$ sending functors $f : A\to B$ to \emph{adjoint pairs} $\bsP_!f : \bsP A\leftrightarrows \bsP B : \bsP^* f$, and its action on 2\hyp{}cells is determined by our desire to privilege the left adjoint, inducing a 2\hyp{}cell $\alpha_! : \bsP_! f \To \bsP_!g$ for each $\alpha : f\To g$. Given $f : A\to B$, the functor $\bsP^*f := \bsP B(\yon_B\circ f,1)$ acts as pre\hyp{}composition with $f$, whereas $\bsP_! f$ is the operation of left extension along $f$. The situation is conveniently depicted in the diagram
	\[
		\begin{tikzcd}
			A \ar[r, "f"] \ar[d, "\yon_A"'] & B \ar[d, "\yon_B"] \\
			\bsP A \ar[shift right, r, "\bsP_! f"'] & \bsP B\ar[shift right, l, "\bsP^* f"']
		\end{tikzcd}
	\]
	Such diagram is filled by an isomorphism when it is closed by the 1\hyp{}cell $\bsP_! f$ and by the cell $\chi^{\yon_B\circ f}$ when it is closed by $\bsP^*f$.
\end{notation}
\begin{definition}[Small functor and small presheaf]
	Let $X$ be a category; we call a functor $F : X^\opp\to \Set$ a \emph{small presheaf} if it results as a \emph{small} colimit of representables; equivalently, $F$ is small if it exist a small subcategory $i : A\subset X$ and a legitimate presheaf $\bar F : A^\opp\to \Set$ of which the functor $F$ is the left Kan extension along $i$.
\end{definition}
\begin{remark}
	The same definition applies, of course, to a functor $F : X \to Y$ between two large categories.
\end{remark}
\begin{notation}[Small presheaves]
	It turns out that the category of small presheaves on $X$ is legitimate in the same universe of $X$ (while the category of all functors $X^\opp \to\Set$ isn't). Given a locally small category there is a Yoneda embedding $X \to [X^\opp, \Set]_s$, having the universal property of free cocompletion of $X$ (see \cite[§3]{il-vecchio}). We explicitly record how the functor $\psh{\firstblank}_s:  \tCat \to \tCat$ acts on 1\hyp{} and 2\hyp{}cells: the universal property of $\psh{\firstblank}_s$ proved in \cite{il-vecchio} implies that in the square
	\[
		\begin{tikzcd}
			A \ar[r, "f"] \ar[d, "\yon_A"'] & B \ar[d, "\yon_B"] \\
			\psh{A}_s \ar[dotted, r] & \psh{B}_s
		\end{tikzcd}
	\]
	the dotted arrow exists (it is the Yoneda extension of $\yon_B\circ f$). The action of $\psh{\firstblank}_s$ on 2\hyp{}cells is uniquely determined as a consequence of this definition.
\end{notation}
\begin{lemma}\label{this_sopra}
	Let $j : \caat\subset\tCat$ be the embedding of small categories in locally small ones. For every large category $X$ there is a natural isomorphism $[X^\opp, \Set]_s\cong \Lan_j\bsP(X)$ (a convenient shorthand is to denote the left extension of $\bsP$ along $j$ as $j_!\bsP$; this is compatible with the notation in \ref{app:relmo} and we will adopt it without further mention). More in particular, there is a canonical isomorphism between $\bsP\triangleleft\bsP = j_!\bsP \circ \bsP$ and the functor $[(\bsP \firstblank)^\opp, \Set]_s$, where for a large category $X$, the category $\psh{X}_s$ designates small presheaves $X^\opp \to \Set$.
\end{lemma}
\begin{proof}
	To show that the universal property of the Kan extension is fulfilled by $[X^\opp, \Set]_s$ we employ a density argument: given a functor $H : \tCat\to \tCat$, every natural transformation $\alpha : \bsP \To H\circ j$ can be extended to a natural transformation $\bar\alpha$ from small presheaves to $H$, using the fact that each $F\in \psh{X}_s$ can be presented as a small colimit of representables: the components of $\bar\alpha_X$ are defined, if $F \cong i_! \tilde F \in\psh{X}_s$ for $i : A\to X$, as
	\[\notag
		\psh{X}_s \xto{i^*} \psh{A}_s = \psh{A} \xto{\alpha_A} HjA = HA \xto{Hi} HX
	\]
\end{proof}
\begin{corollary}
	From this it follows that there is a canonical isomorphism
	\[[X^\opp, \Set]_s \cong \int^{A\in\caat} [X,A]^\opp\times[A,\Set]\]
	(in particular, this specific coend exists even if it is indexed over a non\hyp{}small category); this will turn out to be useful in the proof of \ref{sotto}.
\end{corollary}
\begin{remark}\label{not-a-skiu}
	It is reasonable to expect $j_!\bsP$ to be the small\hyp{}presheaf construction; this construction is the legitimate version of the Yoneda embedding associated to a (possibly large) category $X$. It is important to stress our desire to exploit the results in \cite{fiore2016relative}, but minding that $\bsP$ has additional structure (their approach qualifies $\bsP$ as a monad, but only in the sense that it is a $j$\hyp{}pointed functor, endowed with a unit $\eta : j \to \bsP$ and with a `Kleisli extension' map --instead of a monad multiplication-- sending each $f : jA\to \bsP B$ to $f^\star : \bsP A\to \bsP B$); in view of \ref{it-is-skiu}, now we would like to say that $[\caat, \tCat]$ is a skew\hyp{}monoidal category with skew unit $j$, and composition $(F, G)\mapsto j_!F\circ  G$: this would yield `iterated presheaf constructions' $\bsP \triangleleft\bsP$, $\bsP\triangleleft(\bsP\triangleleft\bsP), (\bsP \triangleleft \bsP)\triangleleft \bsP,\dots$, all seen as functors $\caat\to\tCat$. Unfortunately, given a functor $F: \caat \to \tCat$ it is impossible to ensure that $j_!F$ exists in general ($\caat$ is a ${\mho}^+$\hyp{}category, and $\tCat$ can't be ${\mho}^+$\hyp{}cocomplete), thus --if anything-- the skew\hyp{}monoidal structure of \ref{it-is-skiu} does not exist globally. Fortunately \emph{some} left extensions --precisely those we need-- exist, so we can still employ the `local' existence of $j_!\bsP$ and its iterates to work as if it was part of a full monoidal structure.

	It also turns out (and this is by no means immediate, see \ref{coirens}) that the unitors $\lambda_{\bsP} : j\triangleleft\bsP \to \bsP$ and $\varrho_{\bsP} : \bsP \to \bsP\triangleleft j$ and the associator $\gamma_{\bsP} : (\bsP\triangleleft\bsP)\triangleleft\bsP \to \bsP\triangleleft(\bsP\triangleleft\bsP)$ are all invertible.

	The preliminary results in this section legitimate the practice to na\"i\-ve\-ly consider iterated presheaf constructions: once we consider small functors, the category $[X^\opp,\Set]_s$ lives in the same universe of $X$, and so do all categories $\psh{\psh{A}}_s, \psh{\psh{\psh{A}}_s}_s$.
\end{remark}
\begin{lemma}[$j_!\bsP$ preserves itself]\label{sotto}
	There is a canonical isomorphism
	\[
		\tilde\gamma_{\bsP\bsP} : j_!(j_!\bsP \circ \bsP) \cong j_!\bsP\circ j_!\bsP.
	\]
\end{lemma}
\begin{proof}
	Relying on the previous lemma, we compute the coend
	\[
		j_!(j_!\bsP\circ \bsP) \cong \int^{A\in \caat} [A,X] \times [(\bsP A)^\opp, \Set]_s
	\]
	which is now isomorphic to $[[X^\opp,\Set]_s^\opp,\Set]_s$ in view of the Yoneda reduction and of \ref{this_sopra} (note that the definition of $\lambda A.\bsP A$ and $\lambda A.j_!\bsP A$ entail that $\lambda A.\psh{\bsP A}$ is covariant in $A$).
\end{proof}
\begin{lemma}\label{coirens}
	There are canonical isomorphisms $\lambda_{\bsP} : j\triangleleft \bsP \cong\bsP$, $\varrho_{\bsP} : \bsP \cong \bsP\triangleleft j$ and $\gamma_{\bsP\bsP\bsP} : (\bsP\triangleleft\bsP)\triangleleft\bsP\cong \bsP\triangleleft(\bsP\triangleleft\bsP)$, determined as in \ref{s:uno}--\ref{s:tre}.
\end{lemma}
\begin{proof}
	The isomorphisms come from the unit and associativity constraints of \ref{s:uno}-\ref{s:tre}; as noted in \ref{adjoints-are-exts}, $\varrho$ is invertible in every component because $j$ is fully faithful, and similarly $\lambda$ is invertible in every component if we show that $j$ is a dense functor. Once we have shown this, \ref{sotto} above will conclude, since from \eqref{associator} the associator $\gamma_{\bsP\bsP\bsP}$ coincides with the composition $\tilde\gamma_{\bsP\bsP} * \bsP$.

	Now, the functor $j$ is dense, because the full subcategory of $\tCat$ on the generic commutative triangle $[2] = \{0\to 1\to 2\}$ is dense. So $\caat$, being a full supercategory of a dense category, is dense.
\end{proof}
\begin{remark}
	We often write $\eta\bsP$, $\bsP \eta$, $\mu\bsP$\dots to denote what should be written as $\eta \triangleleft\bsP$, $\bsP \triangleleft\eta$, $\mu\triangleleft\bsP\dots$ Keeping in mind \ref{on-uisge}, that clarifies what is the formal definition for $\triangleleft$\hyp{}whiskerings, there is non chance of confusion: $\eta \triangleleft \bsP = j_!\eta * \bsP$, $\bsP \triangleleft \eta = j_!\bsP * \eta$, and similarly for every other whiskering.
\end{remark}
\begin{proposition}\label{set-is-a-ianua}
	The presheaf construction $\bsP = \psh{\firstblank}$ is a lax idempotent $j$\hyp{}relative monad.
\end{proposition}
In simple terms, a relative monad $T : \cX\to \cY$ is the formal equivalent of a monad in the skew\hyp{}monoidal structure $([\cX,\cY],\triangleleft)$; on $[\caat,\tCat]$, such structure exists `locally' in the components we need allowing us to  work as if $\bsP$ really was a $\triangleleft$\hyp{}monoid. As already said, \cite{fiore2016relative} does not assume the existence of a multiplication map $\mu : \bsP\triangleleft\bsP \to \bsP$, replacing it with coherently assigned `Kleisli extensions' to maps $f : jA\to \bsP B$.

A relative monad is now \emph{lax idempotent} (or a \emph{KZ\hyp{}doctrine}) if it satisfies the 2\hyp{}dimensional analogue of the notion of idempotency; in short, when the algebra structure on an object $A$ is unique up to isomorphism as soon as it exists. Following \cite[2.2]{GARNER20121372}, $\bsP$\hyp{}algebras as cocomplete categories, thus, when an object is a $\bsP$\hyp{}algebra it is so in a unique way (this is of course a behaviour of all formal cocompletion monads).
\begin{proposition} \label{set-is-ian:uno}
	$\bsP$ is a $j$\hyp{}relative monad, if we denote $j : \caat \subset \tCat$ the obvious inclusion.
\end{proposition}
\begin{proof}
	As already noted, the relevant left extensions involved in this proof exist; now, the diagrams we have to check the commutativity of are the following, once we define the Yoneda embedding $\yon_A : A \to \bsP A$ as unit, and $\bsP^* \yon_A : \bsP\bsP A\to \bsP A$ as multiplication of the desired monad.
	\[
		\begin{tikzcd}
			{[[(\bsP A)^\opp, \Set]}_s^\opp,\Set]_s\ar[d, "\bsP_! \mu_A"]  \ar[r,"\mu_{\bsP A}"]& {[(\bsP A)^\opp,\Set]}_s \ar[d, "\mu_A"]\\
			{[(\bsP A)^\opp,\Set]}_s \ar[r, "\mu_A"'] & \bsP A
		\end{tikzcd}
	\]
	\[
		\begin{tikzcd}
			\bsP A \ar[r, "\yon_{\bsP A}"]\ar[dr,equal]& \ar[d, "\mu_A"] \psh{(\bsP A)^\opp} & \bsP A \ar[l, "\bsP_! \yon_A"']\ar[dl,equal]\\
			& \bsP A & {}
		\end{tikzcd}
	\]
	In order to show that they commute, we will exploit the adjunctions $\bsP_!\yon_A\dashv \bsP^* \yon_A$ and $\bsP_!\mu_A \dashv \bsP^*\mu_A$ (the functor $\bsP^*\mu_A$ exists because it must coincide with the Yoneda embedding of $\psh{\bsP A}_s$ into the iterated presheaf category $\psh{\psh{\bsP A}_s}_s$; it must act as the Yoneda embedding $Q\mapsto \psh{\bsP A}_s(\firstblank,Q)$, and this evidently lands into the category of small presheaves on $\psh{\bsP A}_s$ when restricted to small functors).

	For what concerns the unit axiom, the commutativity of the left triangle can be deduced from the chain of isomorphisms
	\begin{align*}
		\mu_{A} \circ \yon_{\bsP A}     & \cong   \Lan_{\yon_{\bsP A} \yon_A}(\yon_A) \circ \yon_{\bsP A}       \\
		                                & \cong \Lan_{\yon_{\bsP A}}(\Lan_{\yon_A}(\yon_A)) \circ \yon_{\bsP A} \\
		(\yon_{\bsP A} \text{ is f.f.}) & \cong  \Lan_{\yon_A}(\yon_A)                                          \\
		(\yon_A \text{ is dense})       & \cong \text{id}_{\bsP A}.
	\end{align*}
	The right triangle corresponds to the composition
	\[\notag
		\begin{tikzcd}[row sep=0]
			\psh{A} \ar[r] & \psh{\psh{A}}_s \ar[r] & \psh{A}\\
			P \ar[r,mapsto] & (Q\mapsto [Q,P]) \ar[r,mapsto] & (a\mapsto [\yon_A(a),P]\cong Pa)
		\end{tikzcd}
	\]
	which is again isomorphic to the identity of $\bsP A$ thanks to the Yoneda lemma.

	In order to show that the multiplication is associative, we prove that $\bsP_!\mu_A\cong \mu_{\bsP A}$ as a consequence of the fact that	there is an adjunction $\mu_{\bsP}\adjunct{1}{} \bsP \mu$, having moreover invertible counit. The argument will be fairly explicit, building unit and counit from suitable universal properties of the presheaf construction and from the definition for $\bsP^*\mu_A$ and $\mu_{\bsP A}$:
	\begin{itemize}
		\item $\bsP^* \mu_A$ sends $\zeta : \bsP A^\opp \to \Set$ into $\psh{\bsP A}_s(\firstblank,\zeta)$ (it plays the exact same r\^ole of a large Yoneda embedding; this will entail that the counit of the adjunction $\mu_{\bsP}\dashv \bsP^*\mu$ is invertible);
		\item $\mu_{\bsP A}$ acts sending $\lambda F.\Theta(F) \in \psh{\psh{\bsP A}_s}_s$ to $\lambda a.\Theta(\hom(\firstblank,a))$.
	\end{itemize}
	With these definitions, the composition $\mu_{\bsP A}\circ \bsP^* \mu_A$ is in fact isomorphic to the identity of $\psh{\bsP A}_s$; we not find the unit map: we refrain from showing the zig\hyp{}zag identities as they follow right away from the explicit description of the co\fshyp{}unit.

	The unit will have as components morphisms $\Theta \To \bsP^* \mu_A(\mu_{\bsP A}(\Theta))$ natural in $\Theta\in \psh{\psh{\bsP A}_s}_s$: given one of these components, its codomain can be rewritten as the coend
	\begin{align*}
		\bsP^* \mu_A(\mu_{\bsP A}(\Theta)) & = \lambda\chi.\psh{\bsP A}_s(\chi,\mu_{\bsP A}(\Theta))                        \\
		                                   & \cong \lambda\chi.\int_{F\in\bsP A}\Set(\chi(F),\mu_{\bsP A}(\Theta)(F))       \\
		                                   & \cong \lambda\chi.\int_{F\in\bsP A}\Set(\chi(F),\Theta(\bsP A(\firstblank,F))) \\
		                                   & \leftarrow \lambda\chi.\Theta(\chi)
	\end{align*}
	and we obtain the candidate morphism in the last line as follows: its component at $\chi \in \psh{\bsP A}_s$ must be an arrow
	\[
		\Theta(\chi) \longrightarrow  \int_{F\in\bsP A}\Set(\chi(F),\Theta(\bsP A(\firstblank,F)))
	\]
	which is induced by a wedge
	\[\Theta(\chi)\to \Set(\chi(F),\Theta(\bsP A(\firstblank,F)));
	\]
	such wedge comes from (the mate of) $\Theta^\text{ar}$, $\Theta$'s function on arrows: the Yoneda lemma now entails that such action induces a map
	\[
		\chi(F)\cong \psh{\bsP A}_s(\bsP A(\firstblank,F),\chi) \xto{\;\Theta^\text{ar}\;} \Set(\Theta(\chi),\Theta(\bsP A(\firstblank,F)))
	\] that by cartesian closure can be reported to
	\[
		\Theta(\chi)\longrightarrow \Set(\chi(F),\Theta(\bsP A(\firstblank,F)))
	\] Of course, this is a wedge in $F$, and we conclude.
\end{proof}
\begin{proposition} \label{set-is-ian:due}
	The monad $\bsP$ is lax idempotent in the sense of \ref{lax-equivs}: there exist an adjunction $\mu_A\adjunct{1}{}\eta_{\bsP A}$.
\end{proposition}
\begin{proof}
	The existence of an adjunction $\mu_A\adjunct{1}{}\eta_{\bsP A}$. will imply all the equivalent conditions in \ref{lax-equivs}, that we nevertheless recall in \ref{lax-equivs-for-P} below for the convenience of the reader. From the definition of these maps, there is a natural candidate to be the counit, and this will be invertible as a consequence of the Yoneda lemma. Indeed, the isomorphism $1 \cong \mu_A\circ \eta_{\bsP A}$ corresponds to the map
	\[
		\lambda a.Fa\mapsto \lambda G.\hom(G,F)\mapsto \lambda a.\hom(A(\firstblank,a),F)\cong Fa.
	\]
	The unit is instead given by the action of a certain functor on arrows, in a similar way as above: given $\chi\in \psh{\bsP A}_s$, there is a canonical map
	\begin{align*}
		\chi(F) & \to \bsP A(F, \lambda a.\chi(\hom(\firstblank,a)))    \\
		        & \cong \int_{a\in A}\Set(Fa,\chi(\hom(\firstblank,a)))
	\end{align*}
	coming from (the mate of) a wedge
	\[
		Fa\cong \bsP A(\hom(\firstblank,a),F) \to \Set(\chi(F), \chi(\hom(\firstblank,a)))
	\]
	defined by the action on arrows of $\chi$.
\end{proof}
The following remark, which is a particular case of \ref{lax-equivs}, characterizes $\bsP$\hyp{}algebras as categories whose Yoneda embedding $\yon_A$ has a left adjoint $\alpha$. These are the \emph{cocomplete} categories \cite[§2]{GARNER20121372}; it is rather easy to see that one of the axioms of $\bsP$\hyp{}algebras asserts that $\alpha(A(\firstblank,a))\cong a$, and since $\alpha$ is a left adjoint it is uniquely determined by sending a colimit of representables into the colimit in $A$ of representing objects (all such colimits, in particular, exist).
\begin{remark}\label{lax-equivs-for-P}
	It turns out from the general theory of lax idempotent monads that $\bsP$ satisfies the following equivalent conditions:
	\begin{enumtag}{pl}
		\item \label{pl:uno} for every pair of $\bsP$\hyp{}algebras $a,b$ and morphism $f :A \to B$, the square
		\[
			\begin{tikzcd}
				\bsP A \ar[d,"a"']\ar[r,"\bsP f"] &  \bsP B \ar[d, "b"] \\
				A \ar[r,"f"'] & B
			\end{tikzcd}
		\]
		\item \label{pl:due} if $\alpha : \bsP A\to A$ is a $\bsP$\hyp{}algebra, there is an adjunction $\alpha\dashv \eta_A$ with invertible counit;
		\item \label{pl:ter} there is an adjunction $\mu_A\dashv \eta_{\bsP A}$ with invertible counit.
	\end{enumtag}
	In particular, we have shown condition \ref{pl:ter}; since \ref{pl:due} holds, there is only a possible choice up to isomorphism for a $\bsP$\hyp{}algebra structure on an object $A$, namely the arrow playing the r\^ole of left adjoint to the Yoneda embedding.
\end{remark}
\begin{exercises}
\item Show that the mate of \eqref{Lan_is_a_functor} (with the same notation therein) is a natural transformation
\[\notag\Lan_{\alpha}F : \Lan_{G'}F \To \Lan_GF\]
which is the component at $F$ of a transformation $\Lan_\alpha : \Lan_{G'} \To \Lan_G$; in other words show that the square
\[\notag
	\vcenter{\xymatrix{
	\Lan_{G'}F \ar@{=>}[d]_{\Lan_\alpha F}\ar@{=>}[r]^{\Lan_{G'}\tau}  & \Lan_{G'}H \ar@{=>}[d]^{\Lan_\alpha H}\\
	\Lan_GF \ar@{=>}[r]_{\Lan_{G}\tau} & \Lan_GH
	}}
\]
commutes for every natural transformation $\tau : F \To H$.
\item \label{ex:counit_of_lans} Write explicitly the unit and the counit of the adjunctions $\Lan_G \adjunct{\epsilon^L}{\eta^L} G^*$ and $G^* \adjunct{\epsilon^R}{\eta^R} \Ran_G$.
\item \label{closed.via.coends} Show that presheaf categories are cartesian closed, via coends: if $\Cat(\C^\opp,\Sets)$ is the category of presheaves on a small $\C$, then there exists an adjunction
\[\notag
	\tCat(\C^\opp,\Sets)(P\times Q, R)\cong \tCat(\C^\opp,\Sets)(P, R^Q)
\]
by showing that $R^Q(c) = \tCat(\C^\opp,\Sets)(\yon_c\times Q, R)$ does the job (use the ninja Yoneda lemma, as well as \ref{naturalu}).
\item Let $\sfK$ be a 2\hyp{}category, and $e : X \to X$ an endo-1\hyp{}cell; show that the following conditions are equivalent:
\begin{itemize}
	\item $e$ is the identity 1\hyp{}cell of $X$;
	\item for every $f :X \to A$ there is a triangle
	      \[\notag
		      \xymatrix{
			      &X\ar[dr]^f &\\
			      X\ar[rr]_f \ar[ur]^e &&A
		      }
	      \]
	      that exhibits $f$ as the right extension of itself along $e$.
\end{itemize}
How can this statement be dualised?
\item \label{ex:kan_is_functor} Use equations (\ref{kanend}) and the ninja Yoneda lemma that $\Lan_\text{id}$ and $\Ran_\text{id}$ are the identity functors, as expected. Use again (\ref{kanend}) and the ninja Yoneda lemma to complete the proof that $F\mapsto \Lan_F$ is a pseudofunctor, by showing that for $\A \xto{F}\B$, $\A\xto{G}\C\xto{H}\D$ there is a uniquely determined laxity cell for composition
\[\notag
	\Lan_H(\Lan_G(F))\cong \Lan_{HG}(F)
\]
(hint: play with the coend $\Lan_H(\Lan_G(F))D$ until you get
\[\notag\int^{XY}(\D(HX, D)\times \C(GY, X))\otimes FY;\]
now use the ninja Yoneda lemma plus co\hyp{}continuity of the tensor, as suggested in Remark \ref{cocotens}).
\item Show that the unit $\eta : F \To \Lan_GF\circ G$ is a wedge in the component $G$; dually, the counit $\epsilon : \Lan_G(H\circ G)\To H$ of the same adjunction is a cowedge. Are the two dinatural transformations the end and the coend of the respective functors?
\item Let $\C$ be a small compact closed monoidal category \cite{MR0470024}; Show that the functor $Y\mapsto \int^X X^\lor \otimes Y \otimes X$ carries the structure of a monad on $\C$.
\item Use coend calculus to prove that if a functor $G$ is fully faithful, then so is $\Lan_G(\firstblank)$. Use coend calculus to prove that if a functor $i : \C \to \D$ has small domain and is dense, then the left Kan extension of the Yoneda embedding $\yon_\C$ along $i$ is a fully faithful functor $\D \to [\C^\opp,\Set]$. Chapter 3 will extensively study this kind of phenomena.
\item \label{codensity_is_monad} Show that if $T_F = \Ran_FF \in \Cat(\D,\D)$ is the codensity monad of a functor $F : \C \to \D$, the two maps $\mu : T_F \circ T_F \To T_F$ and $\eta : \id_{\C} \To T_F$ indeed form a monad, according to \ref{def:monad}.
\item \label{ex:density} Use the end expression in \eqref{the_kleislona} for the hom set in the Kleisli category of the codensity monad of $F : \A \to \C$, and define the \emph{Kleisli composition} in $\Kl(T_F)$ by means of the universal property of \eqref{the_kleislona}. You shall obtain that the Kleisli composition of $f : A \to T_FB$ and $g : B \to T_FC$ is
\[\notag
	\mu_C \circ Tg \circ f : A \to TC
\]
Do the same for the co-Kleisli composition in $\text{coKl}(S^F)$, the density comonad of $F$.
\item 	\index{Adjunction!relative ---}
\index{Relative adjunction}
\label{ex:reladjs} \awful The following exercise draws a piece of the theory of relative adjunctions: we recall that if $J : \A \to \C$, $F : \A \to \B$, and $G : \B \to \C$ are functors, an adjunction between $F$ and $G$, relative to $J$, consists of a family of natural isomorphisms
\[\notag
	\B(FA,B)\cong \C(JA, GB).
\]
If this is the case, we say that $G$ has a \emph{$J$\hyp{}relative left adjoint}. This is denoted $F \adjunct{[J]}{} G$ Dually, (in the same notation) we say that $G$ has a \emph{$J$\hyp{}relative right adjoint} $F$ if there exists a natural isomorphism
\[\notag
	\C(GB,JA)\cong \B(B,FA).
\]
This is denoted $G\adjunct{}{[J]} F$.
\begin{itemize}
	\item Show that $F \adjunct{[J]}{} G$ if $F\cong \leeft_GJ$ and this lifting is absolute; the unit $\eta : J \To G\circ \leeft_GJ = GF$ is the \emph{relative unit} of the adjunction; dually, $G\adjunct{}{[J]} F$ if $F\cong \rift_GJ$ and this lifting is absolute; the counit $\epsilon : G\circ \rift_GJ = GF \To J$ is the \emph{relative counit} of the adjunction.

	      Note in particular that if $F \adjunct{[J]}{} G$ there is no counit, and if $F \adjunct{}{[J]} G$ there is no unit.
	\item Is this criterion also necessary? Namely, is it true that if $F \adjunct{[J]}{} G$ then $F\cong \leeft_GJ$ is absolute? If not, does this mean that there can be two non\hyp{}isomorphic $G,G'$ such that $F \adjunct{[J]}{} G, G'$? What is the structure, if any, of the class of functors $\{G\mid F \adjunct{[J]}{} G \}$?
	\item Assume that
	      \[\notag
		      F \adjunct{[J_1]} G_1 \text{ and } F \adjunct{}{[J_2]} G_2
	      \]
	      are relative adjunctions; it is then possible to build a 2\hyp{}cell
	      \[\notag
		      J_1G_2 \overset{\eta * G_2}\To G_1FG_2 \overset{G_1 * \epsilon}\To G_1 J_2
	      \]
	      pasting the relative unit of $F \adjunct{[J_1]}{} G_1$ with the relative counit of $F \adjunct{}{[J_2]} G_2$. Under which conditions is this 2\hyp{}cell an isomorphism?
	\item What does the invertibility of the relative unit $\eta : J\To GF$ of a relative adjunction $F \adjunct{[J]}{} G$ imply for $G,F$ ($F$ is `relatively fully faithful': what does it mean?). Dual question for the relative counit of $F \adjunct{}{[J]} G$.
\end{itemize}
\item \index{Category!ring as ---} Let $R$ be a ring regarded as a one\hyp{}object category enriched over $\Ab$. A presheaf on $R$ is a right $R$\hyp{}module; let $\A$ be a cocomplete $\Ab$\hyp{}enriched category and $X : R \to \A$ an enriched functor;
\begin{itemize}
	\item Show that the class of such functors $X$ corresponds to the class of objects of $\A$ with an action of $R$;
	\item compute the left Yoneda extension $X\otimes\firstblank$ of $M : R \to \Ab$ and its right adjoint; why the notation $X\otimes \firstblank$?
\end{itemize}
\item Let $F : \C \to \D$ be a functor, and let $\C^\rhd$ be the category $\C$ with an additional terminal object adjoined (see \ref{conecompletion} for the precise notation). There is an obvious embedding functor $\C \hookrightarrow \C^\rhd$: show that the triangle
\[\notag\vcenter{\xymatrix{
			&\C\ar[dr]^F\ar[dl] & \\
			\C^\rhd \ar[rr]_{\bar F}&& \D
		}}\]
is a right extension; who is the 2\hyp{}cell $\eta : \bar F \To F \circ \iota$? Why?
\item Prove that in the Yoneda structure on $\Cat$ the following additional axiom holds:
\begin{quote}
	Let $\xymatrix{B \rtwocell^{B(f,1)}_g{\sigma}& \bsP A}$ be a 2\hyp{}cell;
	if it has the property that the pasting
	\[\notag
		\xymatrix{ A \ar@{}[dr]|(.3)\Swarrow\ar[r]^{\yon_A}\ar[d]_f & \bsP A \\
			B\ar[ur] \urlowertwocell_g{\sigma} &
		}
	\]
	exhibits $\leeft_g\yon_A$, then $\sigma$ is invertible.
\end{quote}
\item Prove that the 2\hyp{}cell
\[\notag
	\xymatrix{
		A\ar@{}[dr]|\To \ar[r]^f \ar[d]_{\yon_A}& B\ar[d]^{\yon_B} \\
		\bsP A & \bsP B\ar[l]_{\bsP f}
	}
\]
witnesses the fact that there is a \emph{lax dinatural} transformation\index{Transformation!lax dinatural ---} $\yon : \id \din \bsP$. Take the occasion to define the notion of lax dinaturality; don't despair if you can't: wait for \ref{laxwedge}. What is the definition of lax dinaturality in the 2\hyp{}category $\tCat$? How is the cell $\upsilon : \yon_A \To \bsP f \circ \yon_B \circ f$ defined in the canonical Yoneda structure on $\tCat$?
\end{exercises}

\chapter{Nerves and realisations}\label{section:nr}
\begin{abstract}
	The present chapter studies a single kind of Kan extensions: the ones where the extendant functor is the Yoneda embedding $\yon_C : \C \to [\C^\opp,\Set]$. Far from being too narrow, the resulting theory of \emph{Yoneda extensions} is astoundingly rich and pervasive; the central object of study of the whole chapter are extensions of the form
	\[\Lan_\yon F \dashv \Lan_F \yon\notag\]
	induced by a functor $F : \C \to \D$; these are called \emph{nerve-realisation} adjunctions. Exploiting the results of the previous chapter, many theorems turn out to be purely formal consequences of the basic rules of coend calculus.

	In order to motivate such pervasivity, we propose a collection of examples of nerve-realisation adjunctions drawing from algebra, topology, geometry, logic, and more category theory.
\end{abstract}
\setlength{\epigraphwidth}{0.8\textwidth}
\epigraph{
Form itself is Void and Void itself is Form.\\
Form is not other than Void and Void is not other than Form.\\
The same is true of Feelings, Perceptions, Mind, and Consciousness.}
{Heart s\={u}tra}
\setlength{\epigraphwidth}{0.75\textwidth}
\section{The classical nerve and realisation}\label{classnerve}
\subsection{Overture: the universal property of $[\C^\opp,\Set]$}
The case in which a Kan extension is done along the Yoneda embedding acquires particular significance in the general theory of Kan extensions: since in a cocomplete category $\D$ all tensors and colimits in \eqref{kanend} exist, every functor $F : \C \to \D$ having small domain and cocomplete codomain has an \emph{extension} to the category of presheaves on $\C$, in a diagram
\[
	\vcenter{\xymatrix{
	&\C  \ar[dr]^F \ar[dl]_{\yon_\C}& \\
	[\C^\opp,\Set] \ar[rr]_{\bar F} & \utwocell<\omit>{\eta} & \D
	}}
\]
filled by an invertible 2\hyp{}cell $\eta$. The scope of this introductory section is to outline this as a universal property. The rest of the chapter is devoted to show how pervasive and expressive the theory of `Yoneda extensions' is.

Let $\C$ be a small category, $\D$ a cocomplete category; then, precomposition with the Yoneda embedding $\yon_{\C} : \C \to [\C^\opp, \Set]$ determines a functor
\[\Cat([\C^\opp, \Set], \D)\xto{\firstblank\circ \yon_{\C}} \Cat(\C,\D),\]
that restricts a functor $G : [\C^\opp, \Set]\to \D$ to act only on representable functors, confused with objects of $\C$, thanks to the fact that $\yon_\C$ is fully faithful. We then have that
\begin{theorem}\label{yext_are_good}\leavevmode
	\begin{enumtag}{ye}
		\item The universal property of the category $[\C^\opp, \Set]$ amounts to the existence of a left adjoint $\Lan_{\yon_{\C}}$ to precomposition, that has invertible unit (so, the left adjoint is fully faithful).
	\end{enumtag}
	This means that $\Cat(\C,\D)$ is a full subcategory of $\Cat([\C^\opp, \Set], \D)$. Moreover
	\begin{enumtag}{yi}
		\item The essential image of $\Lan_{\yon_{\C}}$ consists of those $F : [\C^\opp, \Set] \to \D$ that preserve all colimits.
		\item If $\D = [\E^\opp, \Set]$, this essential image is equivalent to the subcategory of left adjoints $F : [\C^\opp, \Set] \to [\E^\opp, \Set]$.
	\end{enumtag}
\end{theorem}
\begin{proof*}
	The first claim asserts that the following are equivalent:
	\begin{itemize}
		\item every $F : \C \to\D$ extends uniquely to $\hat F : [\C^\opp, \Set]\to \D$ (this in turn means that there is such an $\hat F$ and that if $\hat F', \hat F$ both extend $F$, they are canonically isomorphic);
		\item there exist an adjunction $\Lan_{\yon_\C} \adjunct{\epsilon}{\id} \firstblank\circ\yon_{\C}$ that has invertible unit.
	\end{itemize}
	Proving these claims amounts more or less to playing with definitions, and to notice that that since $\hat F$ is unique up to isomorphism, the correspondence $F\mapsto \hat F$ is a functor determined up to isomorphism. Given the universal property of $\hat F = \Lan_\yon F$, the counit of the adjunction is determined by the fact that there is a unique $\epsilon : \widehat{G \circ\yon} = \Lan_\yon(G\yon)\to G$ induced by the pasting diagram of 2-cells
	\[
		\vcenter{\xymatrix{
				&\C \ar[dl]_{\yon_\C}\ar[dr]^{Gy} & \\
				**[r] [\C^\opp, \Set] \ar[rr] \rrlowertwocell_G{}&
				\utwocell<\omit>{\eta} & \A
			}}
	\]
	half of the triangle identities is simply the fact that the above triangle factors the identity of $G\circ\yon_\C$. The unit of the adjunction is invertible, as a consequence of a general fact: if the extendant arrow $y$ is a fully faithful functor, then there is a canonical isomorphism $\Lan_y(F)\circ y\cong F$ (in case $y=\yon$ is the Yoneda embedding, there's an easy coend proof for this: find it using the ninja Yoneda lemma).

	Now, we shall prove that the counit $\epsilon_G : \Lan_{\yon_\C}(G\circ \yon_{\C}) \To G$ of the above adjunction is invertible if and only if the functor $G$ is cocontinuous.

	Let $G$ be cocontinuous; then since the counit is obtained as the canonical map
	\begin{align*}
		\Lan_\yon(G\yon)(P) & \cong \int^X PX \times G(\hom(\firstblank,X))          \\
		                    & \to G\left(\int^X PX \times \hom(\firstblank,X)\right) \\
		                    & \cong GP
	\end{align*}
	it's easy to see that $G$ is cocontinuous if this is an isomorphism. Vice versa, if the counit gives an isomorphism $G \cong \Lan_\yon(G\yon)$ then
	\begin{align*}
		G(\colim_J P_j) & \cong \int^X \colim_J P_j X \times G(\hom(\firstblank,X)) \\
		                & \cong  \colim_J\int^X P_j X \times G(\hom(\firstblank,X)) \\
		                & \cong \colim_J G(P_j)
	\end{align*}
	(this is in fact a consequence of a more general result: \emph{every} functor that is of the form $\Lan_\yon F$ is cocontinuous; prove this more general fact using the commutation of coends and colimits).

	To complete the proof, we just observe that a cocontinuous functor $L : [\C^\opp, \Set] \to [\E^\opp, \Set]$ between presheaf categories has a right adjoint given by the coend
	\[
		\int^{C \in \C} [\E^\opp, \Set](LC, X) \otimes C,
	\] where $\otimes$ is a tensor in the sense of \ref{tenscotens}.
\end{proof*}
\subsection{Realisations of simplicial sets}\label{realizia}
We start the present section running two examples in parallel. This way, it will seem evident how they are particular instances of the same general construction.
\index{Category!simplex ---}\index{Simplex category}\index{_aaa_delta@$\bDelta$} Consider the category $\bDelta$ made of finite nonempty ordinals and monotone functions, and let us consider the Yoneda embedding $\yon_{\bDelta} : \bDelta \to [\bDelta^\opp,\Sets]$; the presheaf category over $\bDelta$ is usually called the category $\sSet$ of \emph{simplicial sets}; if $\Spc$ now denotes a nice category of topological spaces\footnote{Everyone having their hands a bit dirty of algebraic topology will understand what we mean by `nice' here; the rest of our readers shall take for granted that the whole category of all topological spaces is not nice; it is in fact full of pathological objects we don't want to consider. We shall then restrict to a smaller subcategory $\Spc \subset \BF{Top}$, like the one of \emph{compactly generated} spaces or the one of spaces with the homotopy type of a CW-complex.} we can define two functors $\rho : \bDelta\to \Spc $ and $i : \bDelta\to \Cat$ which `represent' every object $[n]\in\bDelta$ either as a topological space or as a small category:
\begin{itemize}
	\item The category $i[n]$ is obtained regarding a totally ordered set $\{0\le  1\le \dots\le  n\}$ as a category (every poset is a category, in the usual sense: a morphism $x \to y$ between the elements of a poset is the judgment that $x \le y$);
	\item The topological space $\rho[n]$ is defined as the \emph{geometric $n$\hyp{}simplex} $\Delta^n$ embedded in $\bR^{n+1}$ as the subset
	      \[
		      \Big\{(x_0, \dots, x_n) \in \bR^{n+1} \mid 0\leq x_i \leq 1, \textstyle \sum_{i=0}^n x_i = 1 \Big\}.
	      \]
	      endowed with the subspace topology. (What is the affine function $\Delta^f : \Delta^m \to \Delta^n$ induced by a monotone function $f : [m] \to [n]$? Draw pictures of $\rho[n]$ for small values of $n$ and define the \emph{face} maps $\Delta^n \to \Delta^{n-1}$ and \emph{degeneracy} maps $\Delta^n \to \Delta^{n+1}$.
\end{itemize}
Looking for the left Kan extension of these two functors along the Yoneda embedding, we are in the situation depicted by the following diagrams:
\[
	\vcenter{\xymatrix{
	\bDelta \ar[r]^-i\ar[d]_{\yon_{\bDelta}} & \Cat \\
	\sSet\ar@{.>}[ur]_L
	}}
	\qquad
	\vcenter{\xymatrix{
	\bDelta \ar[r]^-\rho\ar[d]_{\yon_{\bDelta}}& \Spc \\
	\sSet\ar@{.>}[ur]_{L'}
	}}
\]
We want to consider the left extensions $L, L'$ of the two functors $i,\rho$ along the Yoneda embedding $\yon_{\bDelta} : \bDelta \to \sSet$; according to our \ref{yext_are_good} above these extensions are left adjoints.

We denote the adjunctions so determined as
\[
	\Lan_\yon i\dashv N_i \quad\text{ and }\quad \Lan_\yon \rho\dashv N_\rho;
\]
\index{Nerve!simplicial ---}
\index{_aaa_N@$N$}
the two right adjoint functors are called the \emph{nerves} associated to functors $i$ and $\rho$ respectively, and are defined as follows:
\begin{itemize}
	\item the \emph{categorical nerve} sends a category $\C$ to the simplicial set $N_i(\C) : [n]\mapsto \Cat(i[n],\C)$; the set of $n$\hyp{}simplices $N_i(\C)_n$ coincides with the set of composable tuples of arrows
	      \[
		      C_0 \xot{f_1} C_1 \xot{f_2} C_2 \leftarrow\cdots\leftarrow C_{n-1} \xot{f_n}C_n.
	      \]
	      In particular, $N(\C)_0$ is the set of objects of $\C$, and $N(\C)_1$ its set of morphisms. Thus the category $\C$ can be reconstructed from its nerve.
	\item the \emph{geometric nerve} sends a topological space to the simplicial set $N_\rho(X) : [n]\mapsto \Spc(\rho[n], X) = \Spc(\Delta^n, X)$ (the \emph{singular complex} of a space $X$).\footnote{The name is motivated by the fact that if we consider the free-abelian group on $N_\rho(X)_n$, the various $C_n = \Z\cdot N_\rho(X)_n = \coprod_{ N_\rho(X)_n}\Z$ organise as a chain complex where differentials are determined as alternating sums of face maps, and whose homology is precisely the singular homology of $X$.}\index{Singular complex}
\end{itemize}
The left adjoints to $N_\rho$ and $N_i$ must be thought as `realisations' of a simplicial set as an object of $\Spc$ or $\Cat$:
\begin{itemize}
	\item The left Kan extension $\Lan_\yon \rho$ is called the \emph{geometric} realisation $|X_\bullet|$ of a simplicial set $X_\bullet$, and it can be characterised as the coend
	      \[
		      \int^{n\in\bDelta} \Delta^n \times X_n
	      \]
	      which in turn coincides to a suitable coequaliser in $\Spc$ in view of our \ref{coends_as_colims} and \ref{endsareeq} (and their duals). The product $\Delta^n\times X_n$ is indeed a product of topological spaces when $X_n$ is thought as discrete: it is the space $\coprod_{s\in X_n}\Delta^n \subseteq \bR^{n+1}$ with the subspace topology on a disjoint union.

	      The shape of this object is determined in light of our description of the coend as a suitable coequaliser (in \ref{coends_as_colims} and \ref{endsareeq} --and their duals), \ie as a quotient space of $\coprod_n \Delta^n \times X_n$ and it agrees with the more classical description presented in almost all books in algebraic topology: the topological space $|X_\bullet|$ is obtained choosing a $n$\hyp{}dimensional disk $\Delta^n$ for each $n$\hyp{}simplex $x \in X_n$ and gluing these disks along the boundaries $\partial\Delta^n$ according to the degeneracy maps of $X_\bullet$.

	      The resulting space is, almost by definition, a CW-complex, because each standard $n$\hyp{}simplex is homeomorphic to a closed disk: this means that $|X_\bullet|$ has the topology induced by a sequential colimit of pushouts of spaces $X_{(0)} \to X_{(1)} \to\dots$ all obtained starting from a discrete space of 0-simplices $X_{(0)}\cong \Delta^0\times X_0$ ($\Delta^0\subset \mathbb{R}$ is a single point).
	\item
	      \index{Category!simplex ---}
	      \index{Simplex category}
	      \index{_aaa_delta@$\bDelta$}
	      The left Kan extension $\Lan_\yon i$ is the \emph{categorical} realisation $\tau_1(X_\bullet)$ of a simplicial set $X_\bullet$, resulting as the coend
	      \[
		      \int^{n\in\bDelta} i[n] \times X_n.
	      \]
	      The $\Set$\hyp{}tensor $i[n] \times X_n$ is interpreted as a product in $\Cat$, where the set $X_n$ is thought as a discrete category; faces and degeneracies of $X$ prescribe how the set $\coprod_{s\in X_n} \{0\to 1\cdots\to n\}$ glue together in the quotient $\left(\coprod_{n\in\bDelta} i[n] \times X_n\right)/{\simeq}$.

	      Note that $\tau_1(X_\bullet)$ \emph{realises} $X_\bullet$ in the sense that the set of objects (resp., of morphisms) of $tau_1(X_\bullet)$ is precisely the set of 0-simplices (resp., of 1-simplices) of $X_\bullet$. This is the category whose objects are 0-simplices of $X_\bullet$, arrows are 1-simplices, and where composition is defined asking that $f,g\in X_1$ compose if there exists a 2-simplex $\sigma$ having 0\hyp{}th face $g$ and 2\hyp{}nd face $f$; identities are witnessed by degenerate simplices.

	      In modern terms (like the ones in \cite{JoyS,joyal2002quasi,HTT}), when $X_\bullet$ is an $\infty$\hyp{}category, we call $\tau_1(X_\bullet)$ the \emph{homotopy category} of $X_\bullet$.
\end{itemize}
\begin{remark}\label{why_lans_are_lans}
	It should now sound as a reasonable conjecture that these examples arise as particular instances of a general theorem. This is indeed the case: the general pattern unifying both these constructions, that we will call it the `nerve and realisation paradigm', was first suggested by D.M. Kan's works in algebraic topology, and in particular on the eponymous \emph{Dold-Kan} correspondence; in this perspective we can recover the classical/singular nerve as particular instance of the paradigm, and many more examples embodied from time to time in different settings; describing this pattern by means of co\fshyp{}end calculus is the scope of the following sections.
\end{remark}
\section{Abstract realisations and nerves}
The upshot of the present section is that whenever $F : \C \to \D $ has small domain and cocomplete codomain, the universal property of the Yoneda embedding in \ref{yext_are_good} determines a functor $\Lan_\yon F$; this functor is always a left adjoint (compare this to the axioms of a Yoneda structure in \ref{ysax_uno}, \ref{ysax_due}).

Algebraic topology, representation theory, geometry and logic constitute natural factories for examples of `nerve and realisations'; more or less everywhere there is an interesting cocomplete category $\D$, there lies an interesting example of nerve-realisation adjunction, induced by a functor $F :\C\to \D$ with small domain.

We now want to lay down the foundations and the terminology allowing to collect a series of readable and enlightening examples.
\begin{definition}[Nerve and realisation contexts]\label{nr-para}\index{Nerve!--- context}
	Any functor $F : \C\to \D$ from a small category $\C$ to a (locally small) \emph{cocomplete} category $\D$ is called a \emph{nerve\hyp{}realisation context} (a NR \emph{context} for short).
\end{definition}
Given a NR context $F$, we can prove the following result:
\begin{proposition}[Nerve-realisation paradigm]\label{nervereal}
	The left Kan extension of $F$ along the Yoneda embedding $\yon : \C\to [\C^\opp, \Sets]$, \ie the functor
	\[L_F=\Lan_\yon F : [\C^\opp, \Sets]\to \D\]
	is a left adjoint, $L_F\dashv N_F$. $L_F$ is called the $\D$-\emph{realisation functor} or the \emph{Yoneda extension} of $F$, and its right adjoint the $\D$-\emph{coherent nerve}.
\end{proposition}
\begin{proof}
	The cocomplete category $\D$ is $\Sets$\hyp{}tensored, and hence $\Lan_\yon F$ can be written as the coend in \eqref{kanend}; so the claim follows from the chain of isomorphisms
	\begin{align*}
		\D\big( \Lan_\yon F(P),D \big) & \cong \D\Big(\int^C [\C^\opp,\Set](\yon_C,P)\otimes F C,D \Big) \\
		                               & \cong \int_C\D\big( [\C^\opp,\Set](\yon_C,P)\otimes F C,D \big) \\
		                               & \cong \int_C\Set\big( [\C^\opp,\Set](\yon_C,P),\D(F C,D) \big)  \\
		                               & \cong \int_C\Sets\big(PC, \D(F C,D)\big).
	\end{align*}
	If we define $N_F(D)$ to be $C\mapsto \D(F C,D)$, this last set becomes canonically isomorphic to $[\C^\opp,\Set](P,N_F(D))$.
\end{proof}
It is straightforward to recognise the choice of $F$ leading to the nerves $N_\rho$ and $N_i$. Also, in light of \ref{yext_are_good} the previous result can be rewritten as follows:
\begin{remark}
	There is an equivalence of categories, induced by the universal property of the Yoneda embedding $\yon_\C : \C \to [\C^\opp,\Sets]$,
	\[
		\firstblank\circ \yon_\C : \Cat(\C, \D)\cong \Cat_!(\Cat(\C^\opp,\Set), \D)
	\]
	whenever $\C$ is small and $\D$ is cocomplete locally small. The left hand side is the category of \emph{cocontinuous} functors $[\C^\opp,\Set] \to \D$.
\end{remark}
\begin{remark}\label{real_commutes_with_lims}
	It is well-known to algebraic topologists that the geometric realisation functor $L_\rho = |\firstblank| : \sSet\to\Spc $ commutes with finite products: coend calculus simplifies this result a lot reducing it to a direct check on representables; it is unfortunately not powerful enough to provide an additional simplification.

	In fact, we can only define a \emph{bijection} between the sets $|\Delta[n]\times \Delta[m]|\cong \Delta^n \times \Delta^m$; after that, a certain amount of dirty work is necessary to show that this bijection is also a homeomorphism with respect to the natural topologies on the two sets. The formal proof that $L_\rho$ commutes with finite products is left as an exercise at the end of the chapter.
\end{remark}
\subsection{Examples of nerves and realisations}
A natural factory of NR contexts is \emph{homotopical algebra}, as such nerve and realisation functors are often used to build Quillen equivalences between model categories. The existence of such Quillen equivalences is somewhat related to the fact that `transfer theorem' for model structures often apply to well-behaved nerve functor: if $\D$ is a model category and $F : \C \to \D$ is a NR functor, oftentimes there is a Quillen adjunction $\Set \leftrightarrows \D$, and sometimes a Quillen \emph{equivalence} induced by the right adjoint $N_F$.

Quillen adjunctions between model categories are certainly not the only examples of NR paradigms. The following list attempts to gather other important examples from different areas of algebra, geometry and logic: for the sake of completeness, we repeat the description of the two above-mentioned examples of the topological and categorical realisations.
\index{Nerve!examples of ---|(}
\begin{example}[Categorical nerve and realisation]\label{catnerve}\index{Category!simplex ---}\index{Simplex category}\index{Nerve!simplicial ---}
	\index{_aaa_N@$N$}\index{Simplicial!--- set}
	If our NR context is $i : \bDelta\to \Cat$, we obtain the \emph{classical nerve} $N_i$ of a (small) category $\C$, whose left adjoint is the \emph{categorical realisation} (the \emph{fundamental category} $\tau_1 X$ of $X$ described in \cite{joyal2002quasi}). The NR adjunction
	\[
		\tau_1 : \sSet\leftrightarrows \Cat : N_i
	\]
	gives a Quillen adjunction between the Joyal model structure on $\sSet$ (see \cite{joyal2002quasi}) and the folk model structure on $\Cat$. This adjunction yields a Quillen adjunction between the category $\Cat$ and the category of simplicial sets; the fibrant objects in $\sSet$, \emph{$\infty$\hyp{}categories}, constitute a model to study category theory in a homotopy coherent fashion.
\end{example}
\begin{example}[Geometric nerve and realisation]\label{topnerve}\index{Simplicial!--- nerve}\index{Simplicial!--- set}
	If our NR context is $\rho : \bDelta\to \Spc$ is the realisation of a representable $[n]$ in the standard topological simplex, we obtain the adjunction between the \emph{geometric realisation} $|X_\bullet|$ of a simplicial set $X_\bullet$ and the \emph{singular complex} of a topological space $Y$, \ie the simplicial set $Y_\bullet$ having as set of $n$\hyp{}simplices the continuous functions $\Delta^n\to Y$.

	If we apply object-wise the free abelian group functor $\Z[\firstblank] : \Set \to \Ab$ to this simplicial set we obtain the simplicial abelian group $\Z Y_\bullet$, which under the Dold-Kan correspondence \ref{doldekanne} gives rise to a (positive degree) chain complex, the \emph{singular complex} of $Y$. The homology of this chain complex coincides with the singular homology of $Y$.
\end{example}
\begin{example}[$\sSet$\hyp{}coherent nerve and realisation]\index{Nerve!simplicially coherent ---}\index{Simplicial!--- set}
	Let $F : \bDelta\to \Cat_{\bDelta}$ be the functor that realises every simplex $[n]$ as a simplicial category having objects the same set $[n]=\{0,1,\dots, n\}$ and as $\hom(i,j)$ the simplicial set obtained as the nerve of the poset $P(i,j)$ of subsets of the interval $[i,j]$ which contain both $i$ and $j$.\footnote{In particular if $i > j$ then $P(i,j)$ is empty and hence so its nerve is the constant simplicial set on $\varnothing$.}

	This sets up a NR context, and if we consider $\Lan_\yon F$ we obtain the \emph{(Cordier) simplicially coherent nerve and realisation}, defined as follows:
	\begin{itemize}
		\item the left adjoint sends a simplicial set into the simplicial category
		      \[ \int^n F[n]\times X_n \]
		      obtained in a similar fashion of $\tau_1(X_\bullet)$ by taking simplicial categories $F[n]$ as shapes and the simplices $s\in X_n$ as gluing instructions. Of course, colimits of (enriched) categories tend to be wildly complicated, but it is an instructive exercise to try to understand how the functor $X_\bullet\mapsto \|X_\bullet\|$ behave on simple examples. We address the reader to \cite{dugger2011,riehl2010coherent}.
		\item the right adjoint $N_{\bDelta} : \Cat_{\bDelta}\to\sSet$ sends a simplicial category $\C$ into a simplicial set constructed remembering that $\C$ carries a simplicial structure. Intuitively, simplicial functors $F[n] \to \C$ carry more information than plain set-enriched functors $[n] \to \C$.
	\end{itemize}
	This adjunction establishes another Quillen adjunction $\sSet\leftrightarrows \Cat_{\bDelta}$ which restricts to an equivalence between quasicategories (fibrant objects in the Joyal model structure on $\sSet$ \cite{Joy}) and fibrant simplicial categories (with respect to the Bergner model structure on $\Cat_{\bDelta}$ \cite{bergner2007model}).
\end{example}
\begin{example}[Moerdijk generalised intervals]\label{mordicchio}
	The construction giving the topological realisation of $\Delta[n]$ extends to the case of any `interval' in the sense of \cite[\S III.1]{Moe}, \ie any ordered topological space $J$ having `endpoints' $0,1$; indeed every such space $J$ defines a generalised topological $n$\hyp{}simplex $\Delta^n{(J)}$ as follows:
	\[
		\Delta^n(J) := \{(x_1,\dots,x_n)\mid x_i \in J, x_0\le \dots\le x_n\} \subseteq J^{n+1}
	\]
	endowed with the subspace and product topology. These data assemble into a NR context $\Delta^\bullet(J) : \bDelta\to \Spc$ that gives rise to an adjunction $\Lan_\yon \Delta^\bullet(J) \dashv N_{\Delta^\bullet(J)}$. Instead of going deep into the technicalities, we address the reader to \cite[\S III.1]{Moe} for more information.
\end{example}
\begin{example}[Toposophic nerve and realisation]\label{toposophic}\index{Topos!Grothendieck ---}
	The correspondence $D : [n]\mapsto {\rm Sh}(\Delta^n)$ defines a \emph{cosimplicial topos}, \ie a functor from $\bDelta$ to the category of  (spatial) toposes, which serves as a NR context. Some geometric properties of this nerve/realisation are studied in \cite[\S III]{Moe}: we outline an instance of a problem where this adjunction naturally arises: let $\cX = \BF{Sh}(X), \cY=\BF{Sh}(Y)$ be the categories of sheaves over topological spaces $X,Y$.

	Let $\cX\star \cY$ be the \emph{join} (see \ref{def:join}) of the two toposes seen as categories: this blatantly fails to be a topos, but there is a canonical `replacement' procedure
	\[
		\vcenter{\xymatrix@R=0mm{
		\Cat \times \Cat \ar[r]^-\star & \Cat \ar[r]^{N_i} & \sSet \ar[r]^{\Lan_\yon D}& \BF{Topos}\\
		(\cX , \cY) \ar@{|->}[r] & \cX \star \cY \ar@{|->}[r] & \Cat(\Delta^\bullet, \cX\star\cY) \ar@{|->}[r] & \cX \odot \cY
		}}
	\]
	that builds a topos out of $\cX$ and $\cY$; various questions about this construction are left as an exercise  in \ref{ex:toposofic}
\end{example}
\begin{example}[The Dold-Kan correspondence] \label{doldekanne}\index{Adjunction!Dold-Kan equivalence}
	The Dold-Kan correspondence \cite{dold1958homology} asserts that there is an equivalence of categories between simplicial abelian groups $[\bDelta^\opp,\BF{Ab}]$ and chain complexes $\Ch^+(\BF{Ab})$ concentrated in positive degree, and it can be seen as an instance of the NR paradigm.

	In this case, the functor $dk : \bDelta\to \Ch^+(\BF{Ab})$ sending $[n]$ to $\Z^{\Delta[n]}$ (the free abelian group on $\Delta[n]$) and then to the \emph{Moore complex} $M(\Z^{\Delta[n]})$ determined by any simplicial group $A\in [\bDelta^\opp,\BF{Ab}]$ as in \cite{GoJ} is the NR context. The resulting adjunction
	\[
		\Lan_\yon(dk) = DK : [\bDelta^\opp,\Ab] \leftrightarrows \Ch^+(\Ab) : N_{dk}
	\]
	sets up an equivalence of categories.
\end{example}
\begin{example}[Étale spaces as Kan extensions]\label{etalage}
	(The present example reworks \cite[1.5]{Bredon1997}) Let $X$ be a topological space, and $o(X)$ its poset of open subsets; let $\Spc/X$ be the slice category of morphisms with codomain $X$ and commutative triangles as morphisms.

	There exists a tautological functor
	\[
		\iota : o(X) \to \Spc/X
	\]
	sending $U\subseteq X$ to itself, regarded as an object $\var{U}{X}$; this works as a NR context, yielding a pair of adjoint functors
	\[
		\Lan_\yon \iota \dashv N_\iota
	\]
	where $N_A$ is defined taking the (pre)sheaf of sections of $p\in \Spc/X$. The resulting left adjoint is precisely the functor sending a presheaf $F : o(X) \to \Sets$ to the disjoint union of stalks $\tilde F = \coprod_{x\in X} F_x$, endowed with the final topology that makes continuous all maps of the form $\tilde s : U\to \tilde{F}$ sending $x\in X$ to the equivalence class $[s]_x \in F_x$.

	In order to see this, let's unpack the definition of the left Kan extension in study: we have to compute the coequaliser 
  \[ 
		\vcenter{\xymatrix{
			\coprod_{V\subseteq U} FU \otimes \var{V}{X} \ar@<.5em>[r]^{l_{UV}} \ar@<-.5em>[r]_{r_{UV}} & \coprod_{U\in o(X)} FU\otimes\var{U}{X} \ar[r] & \Lan_y\iota(F)
		}}
  \]
  where the parallel maps are defined on components by $l_{VU}(s,V) = (s|_V,V)$ and $r_{VU}(s,V)=(s,U)$, and $\otimes$ is the canonical $\Set$-tensor of $\Spc/X$. Such coequaliser imposes on $\coprod_{U\in o(X)} FU\otimes\var{U}{X}$ a relation identifying $(s,U)$ with all its restrictions $(s|_V,V)$ to smaller open sets. 
  
  This means that in the equivalence relation generated by these pairs, a section $(s\in FU)$ is identified with a section $(t\in FV)$ if they coincide on at least an open $W\subseteq U\cap V$. This is very near to the universal property defining the set of all germs $\coprod_{x\in X} F_x$, and in fact it is exactly what is needed to define a natural map $q_U : FU \otimes \var{U}{X} \to \coprod_{x\in X} F_x$ for each $U\in o(X)$ (the cocone condition for these $q_U$ entails that they are ``defined on germs of sections'', in a precise sense that can be easily spelled out). 
  
  In fact, $q_U : FU \otimes \var{U}{X} \to \coprod_{x\in X} F_x$ corresponds, by the universal property of the tensor functor on $\Spc/X$, to a natural family of functions
  \[ \textstyle \bar q_U : FU \to \Spc/X(\iota[U], \coprod_x F_x). \]
  This is to say, every abstract section $s\in FU$ shall give rise to a ``true section'' $\dot s : U \to \coprod_x F_x$.

  A routine argument now shows that the family $(q_U\mid U\in o(X))$ is also initial: every other cocone for the same diagram, say $\zeta_U : FU\otimes\var{U}{X} \to Z$ for some space $Z$, must determine a unique map $\bar \zeta : \coprod_x F_x\to Z$, whose components are $\bar\zeta_x : F_x\to Z$, sending $[s]_x\in F_x$ to $\zeta_u(s)$ for some (the cocone condition makes this choice well-defined) $(s,U)\in FU\times \mathcal U[x]$ ($\mathcal U[x]$ is the filter of neighbourhoods of $x\in X$) having germ $[s]_x$ at $x$.

	This adjunction restricts to an equivalence of categories between the subcategory $\BF{Sh}(X)$ of sheaves on $X$ and the subcategory $\BF{Ét}(X)$ of \emph{étale spaces} over $X$, giving a formal method to prove \cite[II.6.2]{mac1992sheaves}.
\end{example}
\begin{example}[The tensor product of modules as a coend]\index{Product!tensor ---}\label{Tensor}
	Any ring $R$ can be regarded as an $\BF{Ab}$\hyp{}enriched category (see \ref{enrichcat}) having a single object, whose set of endomorphisms is the ring $R$ itself; once noticed this, we obtain natural identifications for the categories of modules over $R$ as covariant and contravariant enriched presheaves on $R$:
	\begin{align}
		\BF{Mod}_R   & \cong \Cat(R^\opp,\BF{Ab})\notag \\
		{}_R\BF{Mod} & \cong\Cat(R,\BF{Ab}).
	\end{align}
	Given $A\in \BF{Mod}_R, B\in {}_R\BF{Mod}$, we can define a functor $T_{AB} : R^\opp\times R\to \BF{Ab}$ which sends the unique object to the tensor product of abelian groups $A\otimes_\Z B$. The coend of this functor can be computed as the coequaliser
	\[
		\coker\Big(
		\xymatrix{
			\coprod_{r\in R}A\otimes_\Z B \ar@<3pt>[r]^-{r\otimes 1}\ar@<-3pt>[r]_-{1\otimes r}& A\otimes_\Z B
		}\Big),
	\]
	that quotients the object $A\otimes_\Z B$ for the submodule generated by the sums $ra\otimes b - a\otimes rb$.

	In other words, there is a canonical isomorphism $\displaystyle\int^{*\in R}T_{AB}\cong A\otimes_R B$.
\end{example}
\begin{remark}\label{Module}\index{Module}\index{Module!bicategory of ---s}\index{Bicategory}
	The previous construction is in fact part of a richer structure: we can define a bicategory $\Mod$ (see \ref{bicat}) having
	\begin{itemize}
		\item 0\hyp{}cells the rings $R,S,\dots$;
		\item 1\hyp{}cells $R \to S$ the modules ${}_RM_S$, regarded as functors $M : R\times S^\opp \to \Ab$;
		\item 2\hyp{}cells $f : {}_RM_S \To {}_RN_S$ are the module homomorphisms $f : M \to N$.
	\end{itemize}
	This bicategory $\Mod$ has a fairly rich structure induced by the one of ${}_R\Mod_S$: for example, bifunctoriality of the tensor product $\otimes$ amounts to the interchange law in the bicategory $\Mod$.

	Chapter \ref{sec:profunctors} on profunctors will extensively generalise this point of view, extending it to the case of multi-object categories, enriched over a generic base.
\end{remark}
\begin{remark}\index{Functor!tensor product of ---}\label{tenso_pro_offunc}\index{_aaa_boxtimes@$\boxtimes_\C$}
	The previous point of view on tensor products can be generalised further (see \cite[\S IX.6]{McL}, but more on this has been written in \cite[\S 4]{Yoneda}): given functors $F,G : \C^\opp,\C\to \V$ having values in a cocomplete (see \ref{being-complete}) monoidal category, we can define the \emph{tensor product} of $F,G$ as the coend
	\[
		F\boxtimes_\C G:= \int^C FC\otimes_{\V}GC.
	\]
	Chapter \ref{section:weight} will generalise in some ways this point of view, regarding this example as an instance of a \emph{weighted colimit} of $F$ with weight $G$, in case $\C$ is a $\V$\hyp{}enriched category.
\end{remark}
\index{Nerve!examples of ---|)}
To appreciate the next example we need to recall the following
\begin{proposition}\label{cofiltness}
	The following properties for a functor $F : \C\to \Sets$ are equivalent:
	\begin{enumtag}{fl}
		\item \label{fl:uno} $F$ commutes with finite limits;
		\item \label{fl:due} $\Lan_\yon F$ commutes with finite limits;
		\item \label{fl:qtr} The \emph{category of elements} $\elts{\C}{F}$ of $F$ (see Def\@. \ref{eltsf} and Prop\@. \ref{elementi}) is cofiltered.
	\end{enumtag}
\end{proposition}
\begin{proof*}
	We provide a proof using as much coend calculus as we can: first of all, from \ref{thm:yoda-is-dense} we know that every presheaf is the colimit of the Yoneda embedding over its own category of elements, thus if $F$ commutes with finite limits its category of elements $\elts{\C}{F}$ has the property that
	\[\colim_{\elts{\C}{F}}\lim_J \yon(A)(D_i) \cong \lim_J \colim_{\elts{\C}{F}} \yon(A)(D_i) \]
	for every object $X\in\C$ and diagram $D : J \to \C$ with finite domain. This is true only if $\elts{\C}{F}$ is filtered.

	Condition \ref{fl:due} now implies condition \ref{fl:uno} because $\Lan_\yon F \circ \yon \cong F$ and the left hand side is a composition of finitely continuous functors.

	Last, condition \ref{fl:qtr} implies condition \ref{fl:due}, because
	\begin{align*}
		\Lan_\yon F(P) & \cong \int^C PC \otimes FC                              \\
		               & \cong \int^C PC \otimes \big(\colim_J \yon(D_j) C\big)  \\
		               & \cong \int^C \colim_J \big(PC \otimes \yon(D_j) C\big)  \\
		               & \cong \colim_J  \int^C \big(PC \otimes \yon(D_j) C\big) \\
		               & \cong \colim_J  P(D_j)
	\end{align*}
	Now, assume $P$ results as a finite limit $\lim_{A\in \A} P_A$; the category $J$ is filtered, thus we can complete the step as
	\begin{align*}
		\textstyle \Lan_\yon F(\lim_\A P_A) & \textstyle \cong \colim_J  \lim_\A P_A(D_j) \\
		                                    & \textstyle \cong \lim_\A \colim_J P_A(D_j)  \\
		                                    & \textstyle \cong \lim_\A \Lan_\yon F(P_A).
	\end{align*}
\end{proof*}
\begin{example}[Giraud theorem using coends]\label{giraudo}\index{Giraud theorem}\index{Topos!Grothendieck ---}
	\emph{Giraud theorem} asserts that every Grothendieck topos is equivalent to a left exact localisations of a presheaf category $[\C^\opp, \Sets]$ (a Grothendieck topos is a category of \emph{sheaves} $\BF{Sh}(\C, \fkj)$ with respect to a \emph{Grothendieck topology} $\fkj$).

	Such a classical `representation' theorem is deeply intertwined with the theory of locally presentable categories (see \cite{Adamek1994,Claudi-2006}), is contained at the end of \cite{mac1992sheaves}.

	We now try to outline an argument that employs coend calculus and Yoneda extensions; we will realise a Grothendieck topos $\E$ as a full subcategory of $[\C^\opp, \Sets]$ for $\C=\E_{<\omega}\subset \E$ the subcategory of \emph{compact} (or \emph{finitely presentable}) objects of $\E$.\footnote{Recall that an object $X\in\E$ is \emph{compact} or \emph{finitely presentable} if the functor $\hom(X,\firstblank)$ commutes with filtered colimits. This essentially means that $X$ can be presented with finitely many generators and relations in the `theory' of $\E$. Also, an object is finitely \emph{generated} if it commutes with filtered colimits of \emph{monomorphisms}.} A finite limit of finitely presentable objects is again finitely presentable, and thus the inclusion $\iota : \E_{<\omega}\subseteq\E$ preserves finite limits.

	The full embedding $\iota : \C \subset \E$ works as a NR context as in \ref{nr-para}, and moreover it is a dense functor; this enables us to use coend calculus to prove that
	\begin{enumerate}
		\item The $\iota$\hyp{}nerve $N_\iota$ is full and faithful; it will turn out to be the inclusion of sheaves into presheaves $[\C^\opp,\Sets]$;
		\item $\Lan_\yon\iota$ is a (left exact, thanks to \ref{cofiltness} and to the fact that $\iota$ is left exact) reflection of $[\C^\opp,\Sets]$ onto $\E$.
	\end{enumerate}
	Hence, the NR adjunction has the following form:
	\[
		\xymatrix@C=2cm{
		\Lan_\yon(\iota) : [\C^\opp, \Sets] \ar@<5pt>[r] & \ar@<5pt>[l]\E : N_\iota
		}
	\]
	As said, since $\iota$ is dense, the nerve $N_\iota$ is fully faithful (see Exercise \ref{ex:nerveisfufai}): it only remains to prove that the functor $\Lan_\yon(\iota)$ behaves like sheafification. In view of the shape of unit and counit of $\Lan_\yon(\iota) \adjunct{\epsilon}{\eta} N_\iota$ (see exercise \ref{unit-and-counit}) means that we have to manipulate the following chain of (iso)morphisms:
	\begin{align*}
		\Lan_\yon(\iota)(P)(C) \cong & \E(\iota C, \Lan_\yon\iota(P))               \\
		\cong                        & \E\big(\iota C, \int^A PA\times \iota A\big) \\
		\leftarrow                   & \int^A \E\big(\iota C, PA\times \iota A\big) \\
		\cong                        & \int^A PA\times \E\big(\iota C, \iota A\big) \\
		\cong                        & \int^A PA\times \C\big(C,A\big)              \\
		\cong                        & PC
	\end{align*}
	This gives an arrow $PC \to LPC$ which is easily seen to be the component of a natural transformation, and in fact of a reflection (to show that $\eta_P$ has the correct universal property, it suffices to show that every morphism $P \To F$ where $F$ is an object in $\E$ uniquely extends to $LP\To F$).

	It remains to prove that this functor is left exact. To do this we invoke \ref{cofiltness}, since $\E_{<\omega}$ is closed under finite limit. It also remains to characterise \emph{sheaves} as those $P$ such that $\eta_P$ is invertible, but this is an equivalent characterisation of orthogonal classes, addressed in Exercise \ref{ex:orthoequivs}.
\end{example}
\begin{example}[Simplicial subdivision functor]\label{ex.sub}\index{Simplicial!--- subdivision}\index{_aaa_Exinfty@$\Ex^\infty$}
	Let again $\bDelta$ be the category of nonempty finite ordinals. The \emph{Kan $\Ex^\infty$ functor} is an endofunctor of $\sSet$ turning every simplicial set $X_\bullet$ into a \emph{Kan complex}.\footnote{A Kan complex is a simplicial set $Y$ such that the functor $\hom(\firstblank,X)$ turns each \emph{horn inclusion} $\Lambda^k[n]\to \Delta[n]$ (see \cite{GoJ}) into an epimorphism.} This construction is of fundamental importance in simplicial homotopy theory, and we now want to re-enact the classical argument given by Kan in the modern terms of a NR paradigm on $\bDelta$, following \cite{GoJ}.

	\index{_aaa_sn@$\text{s}[n]$}
	\index{_aaa_sdn@$\text{sd}[n]$}
	First of all, we note from \cite{GoJ} that the non-degenerate $m$\hyp{}simplices of $\Delta[n]$ are in bijective correspondence with the subsets of $\{0,\dots,n\}$ of cardinality $m+1$; this entails that the set of non-degenerate simplices of $\Delta[n]$ becomes a poset $\text{s}[n]$ ordered by inclusion.

	We can then consider the nerve $N_\rho(\text{s}[n])\in\sSet$ (see Example \ref{topnerve}). This organises into a functor $\text{sd} : \bDelta \to \sSet$, that works as a NR paradigm: using \ref{nervereal} we obtain the adjunction
	\[\xymatrix{
			\sSet \ar@<5pt>[r]^{\Ex}\ar@{}[r]|\top &\ar@<5pt>[l]^{\sfSd} \sSet
		}\]
	where $\Ex$ is the nerve $N_\text{sd}$ associated to the NR paradigm $\text{sd}$ (so a right adjoint to $\sfSd = \Lan_\yon\text{sd}$): the set of $m$\hyp{}simplices $\Ex(X)_n$ is $\sSet(\text{sd}(\Delta[n]), X)$.

	There is a canonical map $\text{sd}(\Delta[n]) \to \Delta[n]$ which by precomposition and by the Yoneda lemma, induces a map $X_n=\sSet(\Delta[n],X)\xto{j^*} \sSet(\text{sd}(\Delta[n]),X)=\Ex(X)_n$, natural in $X\in\sSet$. This gives to $\Ex(\firstblank)$ the structure of a pointed functor, and in fact a \emph{well\hyp{}pointed} functor in the sense of \cite{kelly1980unified}; this, finally, means that we can define
	\[
		\Ex^\infty(X) \cong \varinjlim\Big( X \xto{\eta} \Ex(X) \xto{\eta * \Ex} \Ex^2(X) \xto{\eta * \Ex^2} \cdots\Big)
	\]
	as an endofunctor on $\sSet$. The functor $\Ex^\infty$ enjoys a great deal of formal properties useful in the study of simplicial homotopy theory (the most important of which is that $\Ex^\infty(X)$ is a Kan complex for each $X\in\sSet$, see \cite{GoJ}). A more intrinsic characterisation of this construction is contained in \cite{ehlers2008ordinal}, and defines not only $\text{Sd} = \Lan_\yon\text{sd}$ as a Left Kan extension, but also $\text{sd}$: the authors consider the diagram of 2\hyp{}cells
	\[\label{sd_formaldef}\vcenter{\xymatrix@C=1.5cm{
		\bDelta\times\bDelta \drtwocell<\omit>{\eta}\ar[d]_\oplus \ar[r]^-{\yon_{\bDelta}\times\yon_{\bDelta}}& \sSet \times \sSet \ar[r]^-\times & \sSet \\
		\bDelta \ar[d]_{\yon_{\bDelta}}\ar@/_1pc/[urr]_{\text{sd}} \drtwocell<\omit>{\eta'} &&\\
		\sSet \ar@/_2pc/[uurr]_{\sfSd} &&
		}}\]
	where $\oplus : \bDelta\times\bDelta\to \bDelta$ is the \emph{ordinal sum} defined by $[m]\oplus[n] = [m+n+1]$. Exercise \ref{sd_a_la_cordier} expands on this. In the notation of our Chapter 6, the functor sd is the \emph{convolution} $\yon_{\bDelta} * \yon_{\bDelta}$.
\end{example}
\begin{example}[Isbell duality]\label{isbella-duella}\index{Isbell duality}
	Let $\V$ be a co\fshyp{}complete symmetric monoidal category $\V$ (this will be called a \emph{Bénabou cosmos} in the subsequent sections), and $\C \in \VCat$ a $\V$\hyp{}enriched category (see \ref{enrichcat}); then, we have an adjunction
	\[\vcenter{\xymatrix@C=1.5cm{
		\VCat(\C,\V)^\opp\ar@{}[r]|\perp \ar@<5pt>[r]^{\isbO} &\ar@<5pt>[l]^{\Spec} \VCat(\C^\opp, \V)
		}}\]
	This means that we find a bijection of hom-sets
	\begin{align}
		\VCat(\C, \V)^\opp\big( \isbO(X), Y\big) & = \VCat(\C, \V)\big(Y, \isbO(X)\big) \notag \\
		                                         & \cong \VCat(\C^\opp, \V)(X, \Spec(Y))
	\end{align}
	induced by the functors
	\begin{gather*}
		\isbO : X \mapsto \Big( C\mapsto \VCat(\C, \V)\big(X, \yon_{\C^\opp}(C)\big)\Big),\\
		\Spec : Y\mapsto \Big( C \mapsto \VCat(\C,\V)^\opp\big( \yon_{\C^\opp}(C),Y \big) \Big).
	\end{gather*}
	The adjunction property is a simple derivation in coend calculus:
	\begin{align*}
		\VCat(\C, \V)\big(Y, \isbO(X)\big) & = \int_D \V\Big(YD, \int_A \V(XA, \C(A,D))\Big)  \\
		                                   & \cong \int_{DA}\V\big(YD,\V(XA, \C(A,D))\big)    \\
		                                   & \cong \int_A\V\Big(XA,\int_D\V(YD, \C(A,D))\Big) \\
		                                   & = \VCat(\C^\opp,\V)\big(X, \Spec(Y)\big).
	\end{align*}
\end{example}
\begin{exercises}
\item Use coend calculus to show that the geometric realisation $L_\rho$ commutes with finite products of simplicial sets, assuming that it does commute with finite products \emph{of representables}.
\item Compute the $J$\hyp{}realisation (see Example \ref{mordicchio}) of $X\in \sSet$ in the case $J$ is the Sierpi\'nski space \index{Sierpi\'nski space} $\{0 < 1\}$ with topology $\{\varnothing, J, \{1\}\}$. (hint: the colimit $\int^n \Delta^n(J)\times X_n$ stops after two steps.)
\item \label{unit-and-counit} Let $F : \C \to \D$ be a NR context. Write explicitly the unit and counit of the nerve and realisation adjunction $\text{Lan}_\yon {F} \dashv N_{F}$.
\item \righteyes Show that the nerve functor $N_{F}$ is canonically isomorphic to $\Lan_{F}\yon$, so that there is an adjunction
\[\notag
	\Lan_\yon{F} \dashv \Lan_{F}\yon,
\]
without using coend calculus.
\item \label{ex:nerveisfufai} Show that the following two conditions are equivalent, for a NR context $F : \C \to \D$:
\begin{itemize}
	\item The nerve $N_F : \D \to \Cat(\C^\opp,\Set)$ is fully faithful;
	\item The functor $F$ is \emph{dense}, \ie the density comonad $\Lan_FF$ is the identity functor of $\D$.
\end{itemize}
\item \righteyes \label{tiaccaci}\index{thc situation}
A \emph{tensor-hom-cotensor} situation (\emph{thc situation} for short) consists of a triple $(\otimes , \wedge, [\firstblank,\secondblank])$ of functors between three categories $\catS, \A, \B$, whose covariance type is defined by the adjunctions
\[\notag
	\B(S\otimes A, B)\cong \catS(S, [A,B])\cong \A(A, S\land B).
\]
More precisely, if $\otimes : \catS\times \A\to \B$, then $\land : \catS^\opp\times \B\to \A$, and $[\firstblank,\secondblank] : \A^\opp\times \B\to \catS$. The aim of this exercise is to show that given a thc situation $(\otimes , \wedge, [\firstblank,\secondblank])$, we can induce a new one $(\boxtimes, \curlywedge, \langle\firstblank,\secondblank\rangle )$, on the categories $\catS^{I^\opp\times J}, \A^I, \B^J$, for any $I, J\in\Cat$;
\begin{itemize}
	\item Define $F\boxtimes G\in \B^J$ out of $F\in \catS^{I^\opp\times J}, G\in \A^I$, as the coend
	      \[\notag
		      \int^i F(i, \firstblank)\otimes Gi
	      \]
	      and show that there is an adjunction
	      \[\notag
		      {\B^J}(F\boxtimes G, H)\cong {\catS^{I^\opp\times J}}(F, \langle G,H\rangle )\cong {\A^I}(G, F\curlywedge H)
	      \]
	      for suitable functors $\langle\firstblank,\secondblank\rangle$ and $\firstblank\curlywedge\secondblank$, developing ${\B^J}(F\boxtimes G, H) = \dots$ in two ways, with coend calculus.
	\item Assuming that the relevant structure exists, is it true that composition of the 2\hyp{}category $\tCat$ is the $\otimes$ functor of a thc situation $\otimes : \tCat(\B,\C) \times \tCat(\A,\B) \to \tCat(\A,\C)$? What are the parametric adjoints of $\firstblank\circ F$ and $G\circ\firstblank$?
\end{itemize}
\item \label{ex:toposofic} \awful \index{Topos!Grothendieck ---} Example \ref{toposophic} can be expanded and studied more deeply:
\begin{itemize}
	\item Is $\odot$ a monoidal structure on $\BF{Topos}$?
	\item Under which conditions on $X,Y$ is $\cate X \odot \cate Y$ equivalent to a topos of sheaves on a topological space $X \odot Y$?
	\item What are the properties of the bifunctor $(X,Y)\mapsto X \odot Y$? Does this operation resemble or extend the topological join?
\end{itemize}
\item \awful In this exercise $\Spc$ is a nice category for algebraic topology. Define the category $\Gamma$ having objects the power-sets of finite sets, and morphisms the functions $f : 2^n\to 2^m$ preserving unions and set-theoretical differences.
\begin{enumerate}
	\item Show that there is a functor $\bDelta\hookrightarrow \Gamma$, sending the chain $\{0<1<\cdots<n\}$ in $\bDelta$ to $\{\varnothing\subset \{0\}\subset\cdots\subset \{0,\dots,n\}\}$ in $\Gamma$.
	\item \index{Gamma-space@$\Gamma$-space}
	      \index{_aaa_Gammaspace@$\Gamma$-space} The category of presheaves of spaces $A : \Gamma^\opp\to \Spc$ is called the category of \emph{$\Gamma$\hyp{}spaces}; a $\Gamma$\hyp{}space is said to satisfy the \emph{Segal condition} (or to be Segal) if it turns pushout in $\Gamma$ into homotopy pullback in $\Spc$. Describe pushouts in $\Gamma$; show that a $\Gamma$\hyp{}space is Segal if and only if the following two conditions are satisfied: (i) $A(0)$ is contractible; (ii) the canonical map $A(n)\to \prod_{i=1}^n A(1)$ is a homotopy equivalence in $\Spc$.
	\item Let $X\in\Spc$ and $A : \Gamma^\opp\to \Spc$; define $X\otimes A$ to be the coend (in $\Spc$)
	      \[\notag
		      \int^{n\in\Gamma} X^n\times A(n)
	      \]
	      Show that this defines a bifunctor $\Gamma^\opp\times\Gamma\to\Spc$. Find a canonical map connecting the tensor product of $S^1$ and $A$ with the geometric realisation of the simplicial space $\bDelta^\opp\to \Gamma^\opp\xrightarrow{A}\Spc$.
	\item Let $\C : \Gamma^\opp\to{\Spc}\text{-}\Cat$; let $X\otimes_\Gamma \C$ be the coend (in the category of topological categories)
	      \[\notag
		      \int^{n\in\Gamma} X^n\times \C(n).
	      \]
	      Show that $X\otimes_\Gamma (\firstblank) : {\Spc}\text{-}\Cat\to \Cat$ commutes with finite products, namely if $\C, \D$ are topological categories, then
	      \[\notag
		      X\otimes_\Gamma({\C\times \D})\cong (X\otimes_\Gamma\C)\times (X\otimes_\Gamma {\bf D}).
	      \]
\end{enumerate}
\item \label{sd_a_la_cordier} Write suitable coends for the Kan extensions that define sd and Sd in \eqref{sd_formaldef}, and for their right adjoints.
\item Generalise the NR paradigm to the setting of \emph{separately cocontinuous} (also called \emph{multi-linear}) functors. Given ${F} : \C_1 \times \dots \times \C_n \to \D$, where each $\C_i$ is small and $\D$ is cocomplete, show that there exists an equivalence of categories
\[\notag
	\Cat(\C_1 \times \dots \times \C_n , \D) \cong
	\text{Mult}([\C_1^\opp, \Set]\times \dots \times [\C_n^\opp, \Set], \D)
\]
where $\text{Mult}(\firstblank,\secondblank)$ is the category of all functors that are cocontinuous in each variable once all the others have been fixed (hint: show it `by induction' composing multiple Kan extensions). Given $\theta \in \Cat(\C_1 \times \dots \times \C_n , \D)$, describe the right adjoint of each $\theta^{(i)} : [\C_i^\opp, \Set] \to \D$ (it fixes all components on objects $C_j\in\C_j$ for $j\neq i$, and the $i^\text{th}$ component runs free). All these functors assemble to a `vector-nerve' $N : \D \to [\C_1^\opp, \Set]\times \dots \times [\C_n^\opp, \Set]$.
\index{_aaa_yoncontra@$\yon$}
\index{_aaa_yoncov@$\coyon$}
\item Let $\yon_\C : \C \to [\C^\opp,\Sets]$ be the Yoneda embedding, and ${\coyon}_{\C} : \C^\opp\to \Cat(\C,\Sets)$ its contravariant counterpart. Show that in Example \ref{isbella-duella} we can characterise $\isbO$ as $\Lan_{\yon_\C}({\coyon}_{\C})$.
\item Regard a combinatorial species $f : \B(\N) \to \Set$ as a NR context. Show that there is an isomorphism
\[\notag
	\Lan_\yon f(N_f(X)) \cong \Lan_jf(X)
\]
in the notation of \ref{analuo}, if $j : \B(\N) \to \Set$ is the natural functor. What can you say about the composition
\[\notag
	N_f \circ \Lan_jf : \Set \to \Set \to [\B(\N)^\opp, \Set]
\] for the same species $f$?
\item \awful Make the world a better place by providing it with a deeper study of Isbell duality:
\index{Isbell duality}
\begin{itemize}
	\item An object of a category $\A$ is called \emph{Isbell autodual} if it is a fixed point of the comonad $\isbO\circ\Spec$; is there a general way to characterize all Isbell autodual objects? (Hint: start simple: what if $\A$ is a monoid?)
	\item \index{Isbell duality!--- envelope} The \emph{Isbell envelope} of a category $\A$ consists of the category having
	      \begin{itemize}
		      \item Objects the triples $(A\in\A, \xi : \hom \To P\times Q$ where $P : \A^\opp \to \Set, Q : \A \to \Set$ and $\xi$ is a natural transformation with components $\xi_A : \hom(A,A') \to PA\times QA'$;
		      \item morphisms $(\alpha,\beta) : (A,\xi)\to (A',\eta)$ the pairs $\alpha : P\To P'$ and $\beta : Q'\To Q$ such that $(Q * \alpha)\circ\xi = (\beta * P')\circ\eta$.
	      \end{itemize}
	      What are the properties of this category? Does it have an initial object? A terminal object? Is there a functor between $\text{I}(\A)$ and the category $\A^\daleth \subset\A$ of Isbell autoduals?
	\item Explicitly determine the unit and counit of the Isbell adjunction in \ref{isbella-duella}.
	\item Let $\Gamma\subset\caat$ be a subcategory of $\caat$; let $\Cat(\A^\opp,\Set)_\Gamma \subset \Cat(\A^\opp,\Set)$ be the subcategory of presheaves commuting with limits of functors having domain in $\Gamma$; does the $\clO$ functor land in the subcategory $\Cat(\A,\Set)^\opp_\Gamma$? A particular case of this is when $\A$ is a Grothendieck site, and we want to know whether $\clO(F)$ is a sheaf if $F$ is.
	\item What does Isbell duality look like, when $\Cat(A^\opp,\Set)$ is identified with the category of discrete opfibrations over $\A$, using \ref{thm:equconfib}?
\end{itemize}
\end{exercises}

\chapter{Weighted co\fshyp{}limits}\label{section:weight}
\begin{abstract}
	The present chapter introduces the theory of \emph{weighted co\fshyp{}limits}. Such universal objects constitute a cornerstone of enriched category theory, that can be easily formulated and understood in terms of co\fshyp{}end calculus.

	After having introduced the main definition of weighted limit and colimit, we show that in presence of co\fshyp{}tensors in a $\V$\hyp{}category $\C$ the limit $\wlim{W}F$ of $F : \J \to \C$ weighted by a functor $W : \J \to \V$ can be written as an end
	\[\notag
		\int_J WJ \pitchfork FJ
	\]
	and dually the colimit $\wcolim{W}F$ is a coend
	\[\notag
		\int^J WJ \otimes FJ
	\]
	for $W : \J^\opp\to \V$. This allows to re-read many of the results we already know (for example, Kan extension are weighted co\fshyp{}limits), and to find that many constructions in category theory (comma objects, laxified version of co\fshyp{}limits,\dots) can all be expressed as weighted co\fshyp{}limits for suitable weights $W$. It turns out that weighted co\fshyp{}limits are the correct notion of such universal object in enriched category theory: thus, we conclude the section discussing the theory of \emph{enriched co\fshyp{}ends}.
\end{abstract}
\epigraph{
	No, Time, thou shalt not boast that I do change:\\
	Thy pyramids built up with newer might\\
	To me are nothing novel, nothing strange.
}{W. Shakespeare --- Sonnet CXXIII}

We recall the fundamentals of enriched category theory in a section of our Appendix, \ref{enrichcat}, but the material therein is by no means sufficient to provide a self-contained introduction to the topic; instead, the reader unfamiliar with the basic notions can (and should) consult classical references as \cite{kelly1982basic,Eilenberg1966} or \cite[§6.2]{Bor2}.

Limits and colimits constitute a cornerstone of elementary category theory, because of their ubiquity in describing universal construction. Nevertheless, the notion soon becomes too strict when one moves to the world of \emph{enriched} categories; the `conical' shape of a classical co\fshyp{}limit is not general enough to encompass the fairly rich variety of shapes in which universal objects in $\V$\hyp{}categories arise.

This feature of the enriched\hyp{}categorical world can be justified in many ways: na\"ively speaking, the notion of \emph{cone} for a functor $F : \J\to \C$ is based on the notion of constant functor; yet, in many cases, there is no such a thing as a `constant' $\V$\hyp{}enriched functor $F : \J \to \C$ (we shall notice in \ref{no_enri_dinat} that for the exact same reason $\V$\hyp{}enriched categories do not support a sensible notion of dinaturality). Weighted co\fshyp{}limits circumvent this issue by declaring that a limit depends on \emph{two} arguments: the diagram $F$ \emph{of which} we want to compute the co\fshyp{}limit, and a functor $W \in \Cat(\J ,\V)$ \emph{along which}
we `weight' the co\fshyp{}limit: the span $\V \xot{W} \J \xto{F} \C$ gives rise to an object of $W$\hyp{}shaped $F$\hyp{}diagrams, and the terminal such diagram constitutes the $W$\hyp{}weighted limit of $F$ (and dually, the initial object is the $W$\hyp{}weighted colimit of $F$).

In presence of a terminal functor (for example, when $\V=\Set$), \emph{conical} co\fshyp{}limits arise when we choose $W : \J \to \V$ to be the \emph{terminal} presheaf, that sends every object of $\J$ to the singleton of $\Set$ (and every morphism to the identity function of the singleton).

The theory of weighted co\fshyp{}limits is fairly rich and spans through several chapters of category theory. We cannot touch but the surface of this intricate topic: the interested reader can consult \cite{2catlimits}, a presentation of unmatched clarity filled with enlightening examples.
\begin{remark*}
	As a consequence of the vastness of the topic, our approach is a compromise between full generality and usability, for we are more interested in translating the fundamentals about weighted co\fshyp{}limits into results about suitable co\fshyp{}ends, than in drawing a fully general theory. The focus is not on proofs; a few arguments, not directly linked to our main topic, are only briefly sketched.
	The underlying structure of the present chapters draws a lot from \cite[II.7]{riehl2014categorical}. Quite often, the choice of notation appears to be similar.
\end{remark*}
\begin{notation*}
	Throughout this section, a \emph{weight} is a $\V$\hyp{}enriched functor $W : \J \to \V$, or more generally $W : \I \times \J^\opp\to \V$; we call \emph{ordinary} a category which is enriched over $\Sets$ with its obvious cartesian closed structure. All bases of enrichment $\V$ are \emph{B\'enabou cosmoi}, \ie symmetric monoidal closed, complete and cocomplete categories. These are the `good places' to do enriched category theory.
\end{notation*}
Throughout the section we make heavy use of the \emph{category of elements} of a weight $W : \C^\opp\to \Set$; the reader is invited to follow the present section having thoroughly meditated on \ref{eltsf:char} and \ref{fibelem}.\index{Category!--- of elements}\index{Category of elements}
\section{Weighted limits and colimits}
\begin{remark}[A sophisticated look at classical co\fshyp{}limits]
	Let $F : \J\to\C$ be a functor between small ordinary categories.
	\begin{itemize}
		\item The \emph{limit} $\varprojlim F$ of $F$ can be characterised as the representing object of a suitable presheaf: to define $\lim F$ up to isomorphism we have the natural isomorphism
		      \[
			      \C(C, \varprojlim F)\cong \Cat(\J , \Set)(1, \C(C,F\firstblank))
		      \]
		      where $1$ is a shorthand to denote the terminal functor $\C\to \Sets : C\mapsto 1$ sending every object to the singleton set, and $\C(C,F\firstblank)$ is the functor $\J\to \Sets$ sending $J$ to $\J(C, FJ)$.
		\item Dually, the colimit $\varinjlim F$ of $F : \J\to \C$ can be characterised, in the same notation, as the representing object in the natural isomorphism
		      \[
			      \C(\varinjlim F,C)\cong \Cat(\J^\opp , \Set)(1, \J(F\firstblank,C)).
		      \]
	\end{itemize}
\end{remark}
So, $\Cat(\J^\opp , \Set)(1, \C(F\firstblank,C))$ is a set (of natural transformations), for every $C$, $\J(F\firstblank,C) : \J^\opp\to\Set$ is a presheaf, as well as $C\mapsto\Cat(\J^\opp , \Set)(1, \J(F\firstblank,C)) : \C \to \Set$; if each functor of this sort is representable, we say that $F$ \emph{admits a limit}, precisely the representing object for $\Cat(\J , \Set)(1, \J(F,\firstblank))$: clearly, this is nothing but an instance of the Yoneda-Grothendieck philosophy introduced in \ref{yogrophilo}.

The leading idea behind the definition of weighted co\fshyp{}limit is to generalise this construction to admit shapes other than the terminal presheaf for the domain functor in $1\to \J(F\firstblank,C)$. We can package this rough idea in the following definition:
\begin{definition}[Weighted limit and colimit]\label{weightdef}\index{Colimit!weighted ---}\index{Weighted co\fshyp{}limit}
	Given a diagram of functors
	\[ \xymatrix{
			\C & \ar[l]_F \J \ar[r]^W & \Set
		} \]
	we define the \emph{weighted limit} of $F$ by $W$ (or, equally often, the limit of $F$, weighted by $W$) as a representing object for the presheaf sending $C\in \C$ to $\Cat(\J , \Set)(W, \C(C, F\firstblank))$.

	In other words the weighted limit of $F$ by $W$ is an object $\wlim{W}F\in \C$ such that the isomorphism
	\[\textstyle \C\big(C,\wlim{W}F\big)\cong \Cat(\J , \Set)(W, \C(C,F\firstblank))\]
	holds naturally in $C\in \C$. Dually we define the \emph{colimit} of $F : \J\to \C$ weighted by $W : \J^\opp\to\Set$ to be an object $\wcolim{W}F\in \C$ such that the natural isomorphism
	\[\textstyle \C\big(\wcolim{W}F, C\big)\cong \Cat(\J^\opp, \Set)(W, \C(F\firstblank,C))\]
	holds naturally in $C\in \C$.
\end{definition}
\begin{notat}
	A common alternative notation for the object $\wlim{W}F$ is $\{W,F\}$; a common alternative notation for $\wcolim{W}F$ is $W\otimes F$; this is meant to evoke a THC situation, through the isomorphisms
	\begin{align}
		\C\big(C,\{W,F\}\big)      & \cong \Cat(\J , \Set)(W, \C(C,F\firstblank))     \\
		\C\big(W \otimes F, C\big) & \cong \Cat(\J^\opp, \Set)(W, \C(F\firstblank,C))
	\end{align}
	although this is not properly the same arrangement of functors of \ref{tiaccaci}, the intuition is fruitful.
\end{notat}
\begin{example}\label{frecc}
	Let $[1]$ be the ``generic arrow'' category $\{0\to 1\}$, and let $\lceil f\rceil : [1]\to\C$ be the functor choosing an arrow $f : X\to Y$ in $\C$, and $W : [1]\to \Set$ the functor choosing the arrow $W0=\{0,1\}\to \{0\}=W1$; then a natural transformation $W\Rightarrow \C(C,f)$ consists of arrows $W0\to \C(C,X), W1\to \C(C,Y)$, namely on the choice of two arrows $h,k : C\to X$ such that $fh=fk$: the universal property for $\wlim{W}f$ implies that this is the \emph{kernel pair} of the arrow $f$, namely that $h,k$ fill in the pullback
	\[
		\vcenter{\xymatrix@R=8mm@C=8mm{
		C\ar@/^1pC/[drr]^h\ar@/_1pC/[ddr]_k\ar@{.>}[dr] & \\
		& \varprojlim^Wf\ar[r]\ar[d] & X \ar[d]^f \\
		& X\ar[r]_f & Y
		}}
	\]
	of the arrow $f$ along itself. (The reader is invited to compute the colimit of the same diagram as a preparatory exercise.)
\end{example}
In order to characterize weighted co\fshyp{}limits as co\fshyp{}ends, we employ the same notation of \ref{weightdef}:
\begin{proposition}[Weighted limits as ends]\label{wlimcoends}\index{Co/tensors}
	When the end below and the $\Sets$\hyp{}cotensor (see \ref{tenscotens}) $(X, A)\mapsto X\pitchfork A$ exist, we can express the weighted limit $\wlim{W}F$ for $F : \J\to \C$ as an end, as follows from the chain of computations
	\begin{align*}
		\Cat(\J , \Set)(W,\C(C,F))
		 & \cong \int_{J\in\J}\Set(WJ, \C(C, FJ))              \\
		 & \cong \int_{J\in\J}\C\big(C, WJ\pitchfork FJ \big)  \\
		 & \cong \C\Big(C,\int_{J\in\J} WJ\pitchfork FJ \Big).
	\end{align*}
\end{proposition}
The above derivation implies that there is a canonical isomorphism
\[
	\wlim{W}F\cong \int_{J\in\J} WJ \pitchfork FJ.
\]
The reader might have noticed that we didn't provide a dual statement to \ref{wlimcoends}; this dualisation process is left as an easy exercise, spelled out explicitly in \ref{ex4:wcolims}: the weighted colimit of $F$ by $W$ is a coend, precisely
\[\wcolim{W}F \cong \int^{J\in\J} WJ \otimes FJ.\]
\begin{example}
	Consider the particular case of two parallel functors $W,F : \C\to\Sets$; then we can easily see that $\wlim{W}F$ coincides with the set of natural transformations $W\Rightarrow F$, since the cotensor $WC \pitchfork FC$ is precisely the set $\Sets(WC, FC)$. So: the limit of a presheaf $F$ weighted by a parallel presheaf $W$ is the set of natural transformations $W \To F$.
\end{example}
\begin{example}[Kan extensions]\index{Kan extension!--- as weighted co\fshyp{}limit}\label{kan_are_wei}
	The ninja Yoneda lemma \ref{ninjayo}, rewritten in this notation, says that $\wlim{\C(C,\firstblank)}F\cong FC$ (or, in case $F$ is contravariant, $\wlim{\C(\firstblank,C)}F\cong FC$).

	A short slogan to remember this fact is that
	\begin{quote}
		Representably\hyp{}weighted co\fshyp{}limits are evaluations on the representing object.
	\end{quote}
	(this hides a slightly more general fact: the category of element \ref{eltsf} of $W$ has a terminal object if and only if the weight $W$ is representable, say $\yon A\cong W$; then, $\wlim{W}F\cong FA$. See also our \ref{dirac} about representables playing the role of Dirac delta functions) and suggests that \emph{Kan extensions} can be expressed as suitable weighted co\fshyp{}limits, and more precisely that they can be characterised as those weighted co\fshyp{}limits where the weight is a representable functor (possibly ``twisted'' by a functor $K : \C \to \D$):
	\begin{align}
		\Ran_KF\firstblank & \cong \int_{C\in\C} {\D(\firstblank, KC)}\pitchfork FC\cong \wlim{\D(\firstblank, K\secondblank)}F\label{ext_as_weight_1} \\
		\Lan_KF\firstblank & \cong \int^{C\in\C} {\D(KC,\firstblank)}\otimes FC\cong \wcolim{\D(K\firstblank, \secondblank)}F.\label{ext_as_weight_2}
	\end{align}
	More precisely, $\D(\firstblank,K\secondblank) : \D^\opp\times \C \to \Set$ is a functor, and for every $D\in\D$ the functor $W_D =\D(C,K\secondblank) : \C \to \Set$ works as a weight; the weighted limit of $F$ along $W_D$ is the value of $\Ran_KF$ on the object $D$.

	Note that if $K$ is the identity functor, we obtain the various forms of ninja Yoneda lemma \ref{nyo_1}--\ref{nyo_4} as special cases.
\end{example}
\begin{example}\righteyes\label{ends-are-weighted}\index{Co/end!--- as weighted co\fshyp{}limit}
	Ends themselves can be computed as weighted limits: given $T : \C^\opp\times\C\to\D$ we can take the hom functor $\C(\firstblank,\secondblank) : \C^\opp\times\C\to\Sets$ as a weight, and if the weighted limit exists, we have the chain of isomorphisms
	\begin{align*}
		\wlim{\C(\firstblank,\secondblank)}T & \cong \int_{(C,C')}  \C(C,C') \pitchfork T(C,C')              \\
		                                     & \cong \int_C \Big( \int_{C'} \C(C,C') \pitchfork T(C,C')\Big) \\
		\ref{ninjayo}                        & \cong \int_C T(C,C).
	\end{align*}
\end{example}
\begin{remark}
	It is particularly instructive to unwind the statement above and directly compute the end of $T : \C^\opp\times\C\to\D$ as the equaliser of a pair of maps
	\[
		\vcenter{\xymatrix{
				\prod_{C,C'\in \C} \C(C,C') \pitchfork T(C,C') \ar@<-3pt>[d]\ar@<3pt>[d] \\ \prod_{(f,g) : (C,C') \to (C'', C''')} \C(C,C') \pitchfork T(C'',C''')
			}}
	\]
	determined by the action of $T$ on arrows. In fact, this is exactly what we did in \ref{usefullemma} to prove the Fubini rule \ref{fubozzo}.
\end{remark}
The following Remark and Proposition constitute a central observation.
\begin{remark}[\righteyes The Grothendieck construction absorbs weights]
	Our \ref{weightdef} can be extended to the case $F : \J\to\C$ is a $\V$\hyp{}enriched functor between $\V$\hyp{}categories, and $W : \J\to\V$ is a $\V$\hyp{}presheaf; this is the setting where the notion of a weighted limit proves itself to be the correct one over the `conical' one (where the weight is the terminal presheaf).
	\index{Grothendieck!construction}
	When $\V=\Sets$, indeed, the \emph{Grothendieck construction} sending a presheaf into its category of elements turns out to simplify the theory of $\Sets$\hyp{}weighted limits, reducing every weighted limit to a conical one.
\end{remark}
The last statement is clarified by the following proposition: recall from \ref{eltsf} the definition of the \emph{category of elements} of a functor $W : \J \to \Set$.
\index{_aaa_CintF@$\elts{\C}{F}$}
\begin{proposition}[$\Sets$\hyp{}weighted limits are limits]\label{elementi}
	As shown in \ref{fibelem}, the category of elements $\elts{\J}{W}$ comes equipped with a discrete fibration $\Sigma : \elts{\J}{W}\to\J$; such a fibration is universal, in the sense that for any functor $F : \J\to \C$ one has
	\[\wlim{W}F\cong \varprojlim_{\elts{\J}{W}} \Big(\elts{\J}{W} \xto{\Sigma} \J \xto{F} \C \Big).\]
\end{proposition}
\begin{proof}
	Using \ref{ends-are-weighted}, the characterisation of the end $\int_{J\in \J}WJ \pitchfork FJ $ as an equaliser (as in \ref{endsareeq}), and the characterisation of $\Sets$\hyp{}cotensors as iterated products, we can see that
	\begin{align*}
		\int_{J\in \J}WJ \pitchfork FJ & \cong \eq\Big( \prod_{J\in \J}WJ \pitchfork FJ\rightrightarrows \prod_{J\to J'}FWJ \pitchfork J'\Big)               \\
		                               & \cong \eq\Big( \prod_{J\in\J}\prod_{x\in WJ} FJ\rightrightarrows \prod_{J\to J'}\prod_{x\in WJ}FJ' \Big)            \\
		(\star)                        & \cong \eq\Big( \prod_{(J,x)\in \elts{\J}{W}} FJ\rightrightarrows \prod_{(J,x)\to (J',x')\in \elts{\J}{W}} FJ' \Big) \\
		                               & \cong \varprojlim_{(J,x)\in\elts{\J}{W}} F\circ \Sigma
	\end{align*}
	(step $(\star)$ is motivated by the fact that, thanks to the discrete fibration property for $\Sigma$, every arrow $f : \Sigma(J,x)\to J'$ has a unique lift $(J,x)\to (J', x')$ since $W(f)(x) = x'$).
\end{proof}
\begin{remark}
	If the weight has the form $W=\D(D,K\firstblank)$ for an object $D\in\D$, and a functor $K : \J\to\D$, then the category of elements $\elts{\J}{W}$ is precisely the \emph{comma category} $(D\downarrow K)$: thus the right Kan extension of $F$ along $K$ can be computed as the conical limit of the functor $F\circ\Sigma$, where $\Sigma : (D\downarrow K)\to \J$ is the obvious forgetful functor. We just rediscovered \eqref{ext_as_weight_1} and  \eqref{ext_as_weight_2}.
\end{remark}
When \emph{every} weighted limit exists in $\C$, we can prove that the correspondence $(W,F) \mapsto\wlim{W}F$ is a bifunctor:
\[
	\wlim{(\firstblank)}\,(\secondblank) : \Cat(\J , \Set)^\opp \times \Cat(\J , \C)\to \C.
\] From this, it follows at once that
\begin{itemize}
	\item Functoriality in the $W$ component of $\wlim{W}F$ entails that the terminal morphism $W\To 1$ induces a \emph{comparison arrow} between the $W$\hyp{}weighted limit of any $F : \J\to\C$ and the classical (conical) limit: every weighted limit is a ``fattened up'' version of the conical limit, and there is a comparison arrow $\lim F\to \wlim{W}F$. This intuition has some connection with homotopy theory: it will become clearer in \ref{busfokane}.

	      As an example, consider that the conical limit of the functor $f : [1]\to \C$ described in Example \ref{frecc} consists of the object $\text{src}(f)$; hence the comparison arrow consists of the unique factorisation of two copies of the identity of $\text{src}(f)$ along the kernel pair of $f$.
	\item Using \ref{ends-are-weighted} one can prove that the functor $\wlim{(\firstblank)}F$ is \emph{continuous}, \ie the isomorphism
	      \[\textstyle\label{wlimcocont}
		      \wlim{\big(\varinjlim_{\cate{I}} W_i\big)}F\cong \varprojlim_{\cate{I}} \big(\wlim{W_i}F\big),
	      \]
	      holds for every small diagram of weights $\I\to \Cat(\C, \Sets) : i\mapsto W_i$. Indeed, for a generic object $X\in\C$ we have
	      \begin{align*}
		      \C\Big(X, \wlim{\big(\varinjlim_{\cate{I}} W_i\big)}F\Big) & \textstyle \cong \C\Big(X, \int_A \big(\varinjlim_{\cate{I}} W_iA\big) \pitchfork FA\Big) \\
		                                                                 & \textstyle \cong \C\Big(X, \int_A \varprojlim_{\cate{I}} ( W_iA \pitchfork FA)\Big)       \\
		                                                                 & \textstyle \cong \C(X, \lim_{\cate{I}} \wlim{W_i}F)
	      \end{align*}
	      so the two objects $\wlim{\big(\varinjlim_{\cate{I}} W_i\big)}F$ and $\varprojlim_{\cate{I}} \big(\wlim{W_i}F\big)$ must be canonically isomorphic.
\end{itemize}
The above observation will turn out to be useful during our discussion of \emph{simplicially coherent} co\fshyp{}ends in \ref{cohcoend}.

Weighted \emph{colimits} can be obtained as a straightforward dualisation of the above arguments. The boring technicalities are left to the reader to expand in Exercise \ref{ex4:wcolims}; for the record, we state the proper dualisation in a separate remark:
\begin{remark}\leavevmode\label{weicolims}
	\begin{enumtag}{wc}
		\item \label{are_coends} (weighted colimits as coends) Let
		\[\xymatrix{
				\C & \ar[l]_F \J \ar[r]^W & \Set
			} \] be two functors; if $\C$ admits the coend below, we can express the weight\-ed colimit $\wcolim{W}F$ as
		\[
			\wcolim{W}F\cong \int^{J\in\J} WJ\otimes FJ
		\]
		where we used, like we always do, the $\Sets$\hyp{}tensoring structure of $\C$.
		\item (left Kan extensions as weighted colimits) Let $F : \J\to \C$ and $K : \J\to \D$ be functors; then
		\[
			\Lan_KF\firstblank\cong \int^{J\in\J} \D(KJ,\firstblank)\otimes FJ\cong \wcolim{\D(K\secondblank,\firstblank)}F
		\]
		\item \label{coends_are_hom_weighted} (coends as $\hom$\hyp{}weighted colimits) The coend of $T : \C^\opp\times\C\to \D$ can be written as $\wcolim{\C(\firstblank,\secondblank)}T$, regarding $\hom_\C : \C^\opp\times \C \to \Set$ as a weight:
		\[
			\int^C 	T(C,C) \cong \int^{C,C'} \C(C,C')\otimes T(C,C')
		\]
		\item If the weight $W$ is $\Sets$\hyp{}valued, the colimit of $F : \J \to \C$ weighted by $W : \J^\opp \to \Set$ can be written as a conical colimit over $\elts{\J^\opp}{W}$ using $\Sigma : \elts{\J^\opp}{W} \to \J^\opp$:
		\[
			\wcolim{W}F\cong \varinjlim_{(J,x)\in (\elts{\J^\opp}{W})^\opp}(F \circ \Sigma^\opp)
		\]
		\item (functoriality) If the $W$\hyp{}colimit of $F : \J\to\C$ always exists, then the correspondence $(W,F)\mapsto \wcolim{W}F$ is a functor, cocontinuous in its first variable:
		\begin{gather}
			\wcolim{(\firstblank)}(=) : \Cat(\J^\opp,\Sets) \times \Cat(\J,\C)\to \C, \notag \\
			\textstyle
			\wcolim{\big(\varinjlim_\cate{I} W_i\big)}F\cong \varinjlim_\cate{I} \big(\wcolim{W_i}F\big)
		\end{gather}
		\item (comparison) There is a canonical natural transformation $W\to 1$, inducing a canonical \emph{comparison arrow} from the $W$\hyp{}colimit of any $F : \J\to \C$ to the conical colimit.
	\end{enumtag}
\end{remark}
From the description \ref{ends-are-weighted} and \ref{weicolims}.\ref{are_coends} above of $\wlim{W}F$ and $\wcolim{W}F$ as co\fshyp{}ends, it is immediate that hom functors preserve weighted limits:
\begin{remark}\label{homcommuteswei}
	For every $C\in \C$, every functor $F : \J \to \C$ and weight $W : \J^\opp\to \V$ we have a canonical isomorphism
	\[
		\C\Big( \wcolim{W}F,C \Big) \cong \wlim{W}{\C(F,C)}.
	\]
	Dually, for every functor $F : \J \to \C$ and weight $W : \J \to \V$ we have a canonical isomorphism
	\[
		\C\Big( C,\wlim{W}F \Big) \cong \wlim{W}{\C(C,F)}.
	\]
\end{remark}
As an immediate consequence, every functor that preserves co\fshyp{}ends preserves weighted co\fshyp{}limits; for example, a left adjoint must preserve all weighted colimits, and a right adjoint must preserve all weighted limits.
\section{Examples of weighted colimits}
Due to their deep connections with enriched category theory, homological algebra and algebraic topology are a natural factory of examples of weighted colimits:
\begin{example}[The cone construction as a weighted colimit]\label{cone-is-a-colim}\index{Cone!mapping}\index{Chain complex}
	Let $K$ be a ring, and $\V = \BF{Ch}(K)$ the category of chain complexes of $K$\hyp{}modules. Considering $\V$ as self-enriched, suitably defining the chain complex of maps $C_* \to D_*$, we aim to prove that the \emph{mapping cone} $\text{C}(f) = X_*[1]\oplus Y_*$ of a chain map $f : X_*\to Y_*$ \cite[1.5.1]{Weibel1994} in $\V$ can be characterised as $\wcolim{W}f$, where $f : [1]\to \V$ is the arrow $f$, and $W : [1]^\opp\to \V$ is the functor which chooses the map $S^1(K)_*\to D^2(K)_*$, where $S^n(K)_* = K[n]_*$ is the chain complex with the only nonzero term $K$ concentrated in degree $\hyp{}n$, and $D^n(K)_*$ is the complex
	\begin{equation*}
		\xymatrix{
			\cdots \ar[r]& 0 \ar[r] & K \ar@{=}[r] & K \ar[r] & 0 \ar[r] & \cdots,
		}
	\end{equation*}
	where the first nonzero term is in degree $\hyp{}n$. There is an obvious inclusion $S^n_* \hookrightarrow D^{n+1}_*$:
	\begin{equation*}
		\xymatrix{
			\cdots \ar[r]& 0 \ar[r] & 0 \ar[d] \ar[r] & K \ar[r] \ar@{=}[d] & 0 \ar[r] & \cdots \\
			\cdots \ar[r]& 0 \ar[r] & K \ar@{=}[r] & K \ar[r] & 0 \ar[r] & \cdots.
		}
	\end{equation*}
	We now aim to prove that
	\begin{equation} \label{eq:cap2_cone_coend}
		\text{C}(f) \cong \int^i W(i) \otimes f(i).
	\end{equation}
	In view of (the dual of) Exercise \ref{ispull}, it is enough to show that there is a pushout square
	\begin{equation*}
		\xymatrix{
			W(1) \otimes f(0)\ar@{}[dr]|\ulcorner \ar[r] \ar[d] & W(1) \otimes f(1) \ar@{.>}[d]\\
			W(0) \otimes f(0) \ar@{.>}[r]& \text{C}(f)
		}
	\end{equation*}
	This is a rather simple exercise in universality, given the maps
	\[
		\xymatrix@C=2cm{
		B \ar[r]^{\left(\begin{smallmatrix} 0\\1 \end{smallmatrix}\right)} & \text{C}(f) & \ar[l]_{\left(\begin{smallmatrix} 0 & 1 \\ f & 0 \end{smallmatrix}\right)} A\oplus A[1].
		}
	\]
\end{example}
The following example is more juicy. In addition to those of \ref{eltsf:char} there is a fourth characterisation for the category of elements of a presheaf as a suitable coend over $\Cat$.
\index{_aaa_CintF@$\elts{\C}{F}$}
\begin{example}[\righteyes The category of elements of a presheaf]\label{elts-as-coend}
	The category of elements of a functor $F : \C\to \Sets$ introduced in \ref{eltsf} can be characterised as a $\Cat$\hyp{}weighted colimit: it results as the colimit
	\[
		{\textstyle \elts{\C}{W}}\cong \int^{C\in\C} C/\C \times WC
	\]
	where $WC$ is a set, regarded as a discrete category; it is, in other words, isomorphic to the weighted colimit $\wcolim{S}{W}$, where $S : \C^\opp\to \Cat$ ($S$ as ``slice'') is the functor $C\mapsto C/\C$ (the `coslice' category of arrows $C\to X$ and commutative triangles under $C$).

	To prove this statement, we verify that $\elts{\C}{W}$ has the universal property of the coequaliser of the pair
	\[
		\xymatrix{
			\displaystyle \coprod_{f : A\to B} B/\C \times WA \ar@<4pt>[r]^\alpha \ar@<-4pt>[r]_\beta & \displaystyle \coprod_{C\in\C} C/\C \times WC
		}
	\]
	where $\alpha$ has components $\alpha_f : B/\C \times WA \xto{1\times Ff} B/\C \times WB$ sending $\left(\var{B}{X}, u\right)\mapsto \left( \var{B}{X}, F(f)u \right)$ and $\beta$ has components $\beta_f : B/\C \times WA\xto{f^*\times WA} A/\C \times WA$ sending $\left(\var{B}{X}, u\right)\mapsto \left( \left[\begin{smallmatrix} A &\xto{f} &B \\ &&\downarrow \\ && X \end{smallmatrix}\right], u \right)$. Of course, these maps are composed with the suitable coproduct injections.

	It's rather easy now to define a functor
	\[
		\theta : \coprod_{A\in\C} A/\C\times WA \to \textstyle \elts{\C}{W}
	\]
	having components $\theta_A : A/\C\times WA \to \elts{\C}{W}$ defined by
	\[\left( \var[\!\! f]{A}{B}, u\in FA\right)\mapsto \left( b , F(f)(u)\in Fb\right),\]
	which coequalises $\alpha$ and $\beta$. This functor $\theta$ has the universal property of the coequaliser: given any other $\zeta : \coprod_{A\in\C} A/\C\times WA \to \K$ we can define a functor $\overline{\zeta} : \elts{\C}{W}\to \K$ such that
	\[ \overline{\zeta}(A, u\in FA) = \zeta(\id_A, u). \]
	Now notice that every map $\zeta'$ that coequalises $(\alpha,\beta)$ has the property that
	\[
		\zeta'\left( \var{B}{X}, F(f)u \right) = \zeta'\left( \left[\begin{smallmatrix} A &\xto{f} & B \\ &&\downarrow \\ && X \end{smallmatrix}\right] , u \right)
	\]
	It is now a routine verification to see that $\overline{\zeta}\circ\theta_A = \zeta_A$, and every other functor with this property must coincide with our $\overline{\zeta}$. This concludes the proof.
\end{example}
\begin{remark}[Again an alternative characterisation of the category of elements]\label{its.another.nerve}
	The reader may have noticed that all the above discussion gives a fifth presentation for the category of elements $\elts{\C}{W}$, as the image of $W$ under the Kan extension $\Lan_\yon J$: in the language of \ref{section:nr}, $S : \C^\opp\to \Cat$ is the NR context of the paradigm
	\[
		\xymatrix{
			\elts{\C}{\firstblank} : \Cat(\C, \Cat) \ar@<4pt>[rr]&& \ar@<4pt>[ll] \Cat : N_S
		}
	\]
	where $N_S : \Cat \to \Cat(\C, \Cat) $ is the `nerve' functor sending $\D$ to the functor $C\mapsto \Cat(C/\C, \D)$.
\end{remark}
\begin{remark}
	An alternative approach to characterise $\elts{\C}{W}$ is the following: the category $\elts{\C}{W}$ is precisely the lax limit of $W$ regarded as a $\Cat$\hyp{}valued presheaf \cite[\S4]{2catlimits}, \cite{Graya,Street19}.
\end{remark}
\index{Colimit!homotopy ---}
\index{Bousfield\hyp{}Kan construction}
We can express the Bousfield\hyp{}Kan construction for the homotopy co\fshyp{}limit functor using co\fshyp{}end calculus (see \ref{coends-in-model} for a crash course on what's an homotopy co\fshyp{}limit). We condense Bousfield\hyp{}Kan construction in the following series of examples.
\begin{theorem}[\upeyes The Bousfield-Kan formula for homotopy co\fshyp{}limits]\label{busfokane}
	Let $F : \J\to \M$ be a diagram in a model category $\M$; let moreover $\M$ be equipped with a thc situation (see \ref{tiaccaci}) by functors $\firstblank\pitchfork\secondblank : \sSet^\opp\times \M \to \M$, $[\firstblank,\secondblank] : \M^\opp\times \M \to \sSet$ (so, in particular, $\M$ is $\sSet$\hyp{}enriched), and $\firstblank\otimes\secondblank : \sSet \times \M \to \M$.

	Let us consider the nerve functor of \ref{catnerve}: $N : \Cat \to \sSet$ sends a category $\C$ to the simplicial set of its $n$\hyp{}tuples of composable arrows, for each $[n]\in\Delta$.

	Then the homotopy limit $\holim F$ of $F$ can be computed as the end
	\[
		\int_J N(\J/J)\pitchfork FJ,
	\]
	and the homotopy colimit $\hocolim F$ of $F$ can be computed as the coend
	\[
		\int^J N(J/\J) \otimes FJ.
	\]
\end{theorem}
\begin{remark}
	These two universal objects are weighted co\fshyp{}limits in an evident way: it is possible to rewrite $\holim F\cong \wlim{N(J/\firstblank)}F$ and $\hocolim F\cong \wcolim{N(\firstblank/J)}F$.

	The idea behind this characterisation is to replace the terminal weight with an homotopy equivalent, but fibrant one (in the case of limits; cofibrant, in the case of colimits).

	Bousfield\hyp{}Kan formula arises precisely when we replace the terminal weight with a fibrant one: for every object $J$, both $N(J/\J)$ and $N(\J/J)$ are contractible categories, and they are linked to $N(1)$ by an homotopy equivalence induced by the terminal functor.

	Then, the categories $N(J/\firstblank),N(\firstblank/J)$ must be thought as proper \emph{replacements} for the co\fshyp{}limit functor that correct its failure to preserve weak equivalences (see \cite[Ch. 6]{strom2011modern} for an extremely hands-on account of the theory of homotopy co\fshyp{}limits in algebraic topology, and \cite{Hov} for a standard, easy reference on model categories).
\end{remark}
\index{Co/tensors}
A fairly large class of interesting examples of weighted co\fshyp{}limits comes from the theory of 2\hyp{}categories; many 2\hyp{}dimensional constructions are captured by the above formalism. Moreover, co\fshyp{}end calculus expresses very concretely the shape of the universal object $\wlim{W}F$, as well as its 1- and 2\hyp{}dimensional universal properties. Like in the previous section, we do not seek utter generality, but instead clarity of exposition. Thus, we restrict our attention to the 2\hyp{}category $\tCat$ of categories, strict functors and natural transformations, avoiding to study the case of pseudofunctors, pseudonatural transformations, etc.

We first study a simple example of $\Cat$\hyp{}enriched limit, and its dual; the description of other shapes of weighted co\fshyp{}limits is way more instructive when it is left as an exercise for the reader (see \cite{2catlimits}). Recall Warning \ref{thewarning}.
\begin{example}[\righteyes Inserters in $\tCat$]\label{inserters}
	\index{Weighted co\fshyp{}limit!inserter}\index{Colimit!weighted ---}\index{Inserter}
	Let $\C$ be a 2\hyp{}category, and $f,g : X\rightrightarrows Y$ two parallel 1\hyp{}cells in $\C$; the \emph{inserter} $I(f,g)$ is a pair $(p,\lambda)$ where $p : I(f,g)\to X$ is a 1\hyp{}cell and $\lambda : fp \Rightarrow gp$ is a 2\hyp{}cell, universal with respect to the property of connecting $fp, gp$: this means that the pair $(p,\lambda)$ enjoys
	\begin{enumtag}{in}
		\item A 1\hyp{}dimensional universal property: given a diagram
		\[\vcenter{\xymatrix@R=5mm{
			&X \ar[dr]^f & \\
			B\ar[ur]\ar[dr] \rrtwocell<\omit>{\mu}&&Y \\
			&X\ar[ur]_g&
			}}\]
		this can be split as the whiskering
		\[\vcenter{\xymatrix@R=5mm{
			&&X \ar[dr]^f & \\
			B\ar@{.>}[r]_-h\ar@/^1.5pc/[urr] \ar@/_1.5pc/[drr] & I(f,g)\ar[ur]^p\ar[dr]_p \rrtwocell<\omit>{\lambda}&&Y \\
			&&X\ar[ur]_g&
			}}\]
		for a unique 1\hyp{}cell $h : b\to I(f,g)$ in $\C$: this means, again, that $ph=q$ and $\lambda * h = \mu$.
		\item A 2\hyp{}dimensional universal property: given parallel 1\hyp{}cells $h,k : A \to I(f,g)$ and a 2\hyp{}cell $\beta : ph \To pk$ such that
		\[\label{robaccia}
			\vcenter{\xymatrix{
					A\ar[d]_k\ar[r]^h \drtwocell<\omit>{\beta} & I(f,g) \ar[d]^p \\
					I(f,g) \ar[r]^p\ar[d]_p \drtwocell<\omit>{\lambda} & B \ar[d]^f \\
					B \ar[r]_g & C
				}}\qquad = \qquad
			\vcenter{\xymatrix{
					A \ar[d]_h \ar[r]^k& I(f,g) \ar[d]^p \\
					I(f,g) \ar[r]_p \ar[d]_p & B \ultwocell<\omit>{\beta} \ar[d]^g \\
					B \ar[r]_f & C\ultwocell<\omit>{\lambda}
				}}
		\]
		there is a unique $\bar\beta : h \To k$ such that
		\[
			\vcenter{\xymatrix@C=5mm@R=5mm{
			& I(f,g) \ar[dr] &   & &I(f,g) \ar[dr]\ar@{=}[dd] & \\
			A \rrtwocell<\omit>{\beta}\ar[ur]\ar[dr]&& B  \ar@{}[r]|= & A \rtwocell<\omit>{\bar\beta} \ar[ur]\ar[dr]&  & B \\
			& I(f,g) \ar[ur] & & &I(f,g) \ar[ur] &
			}}
		\]
	\end{enumtag}
	Now, if $\J $ is the category $\{0 \rightrightarrows 1\}$ with two objects and two parallel non-identity arrows, the inserter $I(f,g)$ is the limit of the functor $F : \J \to \C$ choosing the two 1\hyp{}cells $f,g$, weighted by the weight $W : \J \to \Cat$ choosing the parallel `source' and `target' functors $s,t : [0] \rightrightarrows [1]$.

	We shall deduce the shape of the inserter when $\C =\tCat$ is the 2\hyp{}category of categories. In such case, the end
	\[
		\int_{J \in \J} WJ \pitchfork FJ
	\]
	that according to \ref{weightdef} defines the weighted limit boils down to the object of $\Cat$\hyp{}natural transformations $W \To F$ (see also \ref{naturalu}): it is indeed the case that such a natural transformation is determined as a pair $(b, u : c \to c') \in B \times C^{[1]}$ such that $u : fb \to gb$. In fact, naturality corresponds to the commutativity of the following two squares:
	\[ \vcenter{\xymatrix{
		{[0]} \ar[r]^b\ar[d]_s & B\ar[d]^f  & {[0]} \ar[r]^b\ar[d]_t
		& B\ar[d]^g  \\
		{[1]} \ar[r]_u & C & {[1]} \ar[r]_u & C
		}} \]
	and this means precisely that $u$ has $fb$ as domain, and $gb$ as codomain.

	Let's show that this object, as a subobject $p : (B \times C^{[1]})^\circ \subseteq B \times C^{[1]}$ has the desired universal property: first of all, $(B \times C^{[1]})^\circ$ clearly is a subcategory of the product category $B \times C^{[1]}$. There is an obvious natural transformation $\lambda : fp \To gp$ defined on components as $\lambda_{(b,u)} : fb \xto{u} gb$. We leave to the reader to check that this is indeed the component of a natural transformation.
	\begin{enumtag}{in}
		\item Every $\alpha : fq\To gq$ with components $\alpha_a : fqa \to gqa$ is such that $(qa, \alpha_a) \in (B \times C^{[1]})^\circ$, and the map $h : A \to (B \times C^{[1]})^\circ$ that sends $a$ into $(qa, \alpha_a)$ is a functor because $\alpha$ was natural.
		\item A similar argument shows that every $\beta : ph \To pk$ satisfying \eqref{robaccia} factors through $(B \times C^{[1]})^\circ$.
	\end{enumtag}
	Altogether, these two properties shows that there is a unique isomorphism between $(B \times C^{[1]})^\circ$ and $I(f,g)$. This concludes the proof.
\end{example}
\begin{example}[Comma objects]\label{commae_obj}
	\index{Comma object}
	\index{Comma category}\index{Category!comma ---}
	Let $\C$ be a 2\hyp{}category, and $f,g$ a cospan of 1-cells like
	\[\vcenter{
			\xymatrix{
				& C \ar[d]^g \\
				B \ar[r]_f & X
			}
		}\]
	This can be regarded as the image of a functor $F : \Lambda^2_2 \to \C$, where $\Lambda^2_1$ is the ``generic cospan'' $\{0 \to 2 \leftarrow 1\}$. Let us consider the weight $W : \Lambda^2_2 \to \Cat$ whose image is
	\[
		W\left[
			\begin{smallmatrix}
				&& 1 \\
				&&\downarrow \\
				0 &\to& 2
			\end{smallmatrix}
			\right] =
		\vcenter{\xymatrix@!=4mm{
		& [0] \ar[d]^{d_1} \\
		[0] \ar[r]_{d_0} & [1]
		}}
	\]
	where $d_i : \{i\} \to \{0\to 1\}$ chooses the object $i$. Let us prove that the limit of $F$ weighted by $W$ is the \emph{comma object} of $f,g$; evidently, in the special case of $\C=\Cat$, the limit $\wlim{W}F$ is the \emph{comma category} of \ref{def:comma}.

	Let us fix an object $A$ of $\C$; a natural transformation $W \To \C(A,F\firstblank)$ consists of the following data:
	\begin{enumtag}{c}
		\item A 1-cell $u : A \to B$;
		\item A 1-cell $v : A \to C$;
		\item A 2-cell $\lambda : [1] \to \C(A,X)$, whose source is forced by the naturality condition to be $fu$, and whose target is $gv$.
	\end{enumtag}
	More explicitly, natural transformations $W \To \C(A,F\firstblank)$ correspond to squares
	\[\vcenter{\xymatrix{
				A \ar[r]^v \ar[d]_u & C\ar[d]^g \\
				B \ar[r]_f & X \ultwocell<\omit>{\lambda}
			}}\]
	filled by a 2-cell $\lambda : fu\To gv$. The terminal such 2-cell is then $\wlim{W}F$.
\end{example}
One can now routinely dualise the above construction to define \emph{cocomma objects}: the following example freely uses some notions from the subsequent chapter.
\begin{example}[Cocomma objects]
	\index{Cocomma object|see{Comma object}}
	Let $\C$ be a 2\hyp{}category, and $f,g$ a span of 1-cells like
	\[\vcenter{
			\xymatrix{
				X \ar[r]^g \ar[d]_f & C\\
				B &
			}
		}\]
	This can be regarded as the image of a functor $F : \Lambda^2_0 \to \C$, where $\Lambda^2_0$ is the ``generic span'' $\{0 \to 2 \leftarrow 1\}$. Let us consider again the weight $W : (\Lambda^2_0)^\opp \cong \Lambda^2_2 \to \Cat$ whose image is again
	\[
		\vcenter{\xymatrix@!=4mm{
		& [0] \ar[d]^{d_1} \\
		[0] \ar[r]_{d_0} & [1]
		}}
	\]
	where $d_i : \{i\} \to \{0\to 1\}$ chooses the object $i$. The colimit of $F$ weighted by $W$ is the \emph{cocomma object} of the pair $(f,g)$.
\end{example}
Exercise \ref{ex5:cocomma-as-colage} provides you with another proof that the cocomma object in $\tCat$ of two functors $f,g : B \leftarrow X \to A$ can be described as the category $\comma{f\\g}$ having objects those of $B\sqcup C$, and morphisms $x\to y$ defined as follows:
\begin{itemize}
	\item $B(b,b')$ if $b=x,b'=y$ are both objects of $B$;
	\item $C(c,c')$ is $c=x,c'=y$ are both objects of $C$;
	\item if $b=x\in B_o$ and $c=y\in C_o$, the set of morphisms in $\comma{f\\g}(b,c)$ is the coend
	      \[\int^x C(gx,c)\times B(b, fx);\]
	\item otherwise, the hom set is empty.
\end{itemize}
We shall now present a proof of the universality of $\comma{f\\g}$ based solely on the construction of the coend that defines $\comma{f\\g}(b,c)$.

The intuition that shall guide the reader is that $\comma{f\\g}(b,c)$ is a set of ``fake arrows'', \ie of triples
\[
	b \xto{\varphi} fa \overset{(\xi)}{\dashrightarrow} ga \xto{\psi}c
\]
where the arrow $(\xi) : fa \dashrightarrow ga$ is in some suitable sense freely adjoined in the disjoint union $B\sqcup C$.

Inside the coend, we identify two fake arrows $b \xrightarrow{\varphi} fa \overset{(\xi)}\dashrightarrow ga \xrightarrow{\psi}c$ and $b \xrightarrow{\varphi} fa' \overset{(\xi')}\dashrightarrow ga' \xrightarrow{\psi}c$ precisely when there is a ``hammock'' diagram between $a, a'$ of the following form:
\[\scriptscriptstyle
	\vcenter{\xymatrix@!=4mm{
	& fa \ar@{.>}[r]^{(\xi)}             & ga \ar[rdd]^\psi                      &   \\
	& fa_1 \ar@{.>}[r]\ar[u] \ar[d] & ga_1 \ar[u] \ar[d] \ar[rd] &   \\
	b \ar[ruu]^\varphi \ar[rdd]_{\varphi'} \ar[ru] \ar[rd] & \vdots                   & \vdots                              & c \\
	& fa_n\ar@{.>}[r] \ar[u] \ar[d] & ga_n \ar[u] \ar[d] \ar[ru] &   \\
	& fa' \ar@{.>}[r]_{(\xi')}            & ga' \ar[ruu]_{\psi'}                     &
	}}
\]
Now we have enough material to discover the universal property of such an object: first of all, there are canonical functors $i_B : B \to \comma{f\\g}$ and $i_C : C \to \comma{f\\g}$, and a canonical natural transformation $\zeta\colon i_B f \Rightarrow i_C g$ comes from taking equivalence class of identity arrows $(1_{ga}, 1_{fa}) \in Y(ga,ga)\times X(fa,fa)$ under the composition
\[
	\vcenter{\xymatrix@R=4mm@C=1mm{
	\hom(ga,ga)\times \hom(fa,fa)\ar[d] &\ni(1_{ga}, 1_{fa}) \ar@{|->}[dd]\\
	\coprod_{a\in A}\hom(c,ga)\times \hom(fa, b) \ar[d]\\
	{\comma{f\\g}}(fa, ga) &\ni [(1_{ga}, 1_{fa})]
	}}
\]
as a natural candidate for $\zeta_a \colon fa\to ga$.

Now, suppose we are given a commutative square
\[
	\vcenter{\xymatrix{
			X \ar[r]^g \ar[d]_f & C\ar[d]^v \\
			B \ar[r]_w& Y \ultwocell<\omit>{\theta}
		}}
\]
filled by a 2-cell $\theta : wf \Rightarrow vg$. Then define a unique functor $u : \comma{f\\g} \to Y$ on objects and true arrows in $B$ or $C$ acting as $v\colon C\to Y, w\colon B\to Y$; the action of $u$ on a fake arrow
$$
	b \xrightarrow{\varphi} fa \overset{(\xi)}\dashrightarrow ga \xrightarrow{\psi}c
$$
is induced by the composition
$$
	wb \xrightarrow{w*\varphi} wfa \xrightarrow{\theta_a} vga \xrightarrow{v*\psi}vc
$$
(all arrows exist now!). At this point, all remaining checks are pure routine: $u$ is unique due to the tautological definition of $\zeta$; the cell $\zeta$ has a 2-dimensional universal property as well; $\zeta\circ \firstblank$ reflects isomorphisms.

Once this is done, the willing reader can embark on all sort of instructive computations: exchanging the r\^ole of $B,C$ in the above construction, one obtains a 2-cell in the opposite direction; in the terminology of the next chapter (see \ref{embare_i_profi} and \eqref{proco}), $\comma{f\\g}$ is the category of elements of the composite profunctor $B \overset{\proP_f}\pto X \overset{\proP^g}\pto C$ regarded as a presheaf $B^\opp\times C \to \Set$;  what if the functor $f$ or $g$ is the identity?
\section{Enriched co\fshyp{}ends}
\subsection{Preliminaries on enriched categories}
In the setting of enriched category theory, the property of being complete is stated in terms of an existence result for every weighted limit $\wlim{W}{F}$: in short, the reason for this choice is that the categories admitting only co\fshyp{}limits weighted by terminal presheaves contain too few objects and shall not be considered complete  by the internal language of $\underline{\VCat}$.

This can be made precise in the following way: there is a canonical choice of a Yoneda structure (see \ref{def:yoneda-struc}) on the 2\hyp{}category of $\V$\hyp{}enriched categories, where Yoneda maps are given by $\V$\hyp{}enriched Yoneda embeddings; every 2\hyp{}category with a Yoneda structure has a notion of a \emph{cocomplete object}, and in that Yoneda structure co\fshyp{}completeness coincides with having all weighted co\fshyp{}limits (see \cite{street1978yoneda}).

The present section serves the purpose to introduce a sensible definition of \emph{enriched co\fshyp{}end} in \ref{enriend:uno}, \ref{enriend:due}: in short, using \ref{weicolims}.\ref{coends_are_hom_weighted}, given a $\V$\hyp{}functor $H : \C^\opp\boxtimes\C\to\V$ we can define the end $\int_C H(C,C)$ to be the limit of $H$ weighted by the enriched hom $\V$\hyp{}functor $\C^\opp\boxtimes\C\to\V$ (we will introduce the $\firstblank\boxtimes\secondblank$ product in \ref{cosmuclosed} below).

In case all co\fshyp{}ends exist and the codomain category of functors $T : \C^\opp\times \C \to \A$ has co\fshyp{}tensors, we can moreover easily characterise the co\fshyp{}end functor as an adjoint to a `co\fshyp{}tensor with hom' functor.
\begin{lemma}
	\label{importantlemma} Let $\C$ be a tensored and cotensored $\V$\hyp{}category and $W$ a functor to be treated as a weight; then
	\begin{enumtag}{wa}
		\item if $W : \J \to \V$, the $\wlim{W}{\firstblank}$ functor has a left adjoint, given by tensoring with the weight $W$; in other words, there is an adjunction
		\[\vcenter{\xymatrix@C=2cm{
			\C \ar@{}[r]|-\perp\ar@<5pt>[r]^-{W\otimes\firstblank}&\ar@<5pt>[l]^-{\wlim{W}{\firstblank}} \Cat(\J, \C)
			}}\]
		where $W\otimes C : \lambda J. WJ\otimes C$.
		\item Dually, if $W : \J^\opp\to \V$, the $\wcolim{W}{\firstblank}$ functor has a right adjoint, given by cotensoring with the weight $W$; in other words, there is an adjunction
		\[\vcenter{\xymatrix@C=2cm{
			\Cat(\J, \C) \ar@{}[r]|-\perp\ar@<5pt>[r]^-{\wcolim{W}{\firstblank}}&\ar@<5pt>[l]^-{W\pitchfork\firstblank} \C
			}}\]
		where $W\pitchfork C : \lambda J. WJ\pitchfork C$.
	\end{enumtag}
\end{lemma}
\begin{proof*}
	Both arguments consist of an easy computation in coend calculus:
	\begin{align*}
		\C(C, \wlim{W}{F}) & \cong \C\Big(C, \int_J WJ\pitchfork FJ\Big) \\
		                   & \cong \int_J \C(WJ\otimes C, FJ)            \\
		                   & \cong \Cat(\J, \C)(W\otimes X, F).
	\end{align*}
	Dually, we have
	\begin{align*}
		\C(\wcolim{W}{F}, C) & \cong \C\Big(\int^J WJ\otimes FJ, C\Big) \\
		                     & \cong \int_J \C(WJ\otimes FJ, C)         \\
		                     & \cong \int_J \C(FJ, WJ\pitchfork C)      \\
		                     & \cong \Cat(\J, \C)(F, W\pitchfork X).
	\end{align*}
\end{proof*}
\begin{remark}
	This finally sheds a light on our proof of Fubini theorem in \ref{fubozzo}: given that co\fshyp{}ends are hom-weighted co\fshyp{}limits (see \ref{weicolims}.\ref{coends_are_hom_weighted}), the co\fshyp{}end functor has a right/left adjoint given by tensoring with the weight hom: this is exactly the way in which we proved that $\int^A : \Cat(\C^\opp\times \C, \D)$ had a right adjoint, only without explicitly mentioning the technology of weighted co\fshyp{}limits.
\end{remark}
\subsection{The theory of enriched co/ends}
A fundamental step in laying the foundation of weighted limit theory is the natural isomorphism
\[\label{iso_per_wlim}
	\C\big(C,\wlim{W}F\big) \cong \VCat(\J,\V)(W, \C(C,F\firstblank))
\]
valid for functors
\[ \xymatrix{
		\C & \ar[l]_F \J \ar[r]^W & \Set.
	} \]
In order to export this isomorphism to the enriched setting, we shall make \eqref{iso_per_wlim} take place in the base cosmos $\V$; in short, this means that we have to find a way to promote the category $\VCat(\J,\V)$ of $\V$-functors and $\V$\hyp{}natural transformations as an enriched category $\underline{\VCat}(\C , \V)$: this means that every $\VCat(\J,\V)(F,G)$, for functors $F,G$, must become an \emph{object of enriched natural transformations} in $\V$.

To do this, we will endow $\VCat$ with a closed symmetric monoidal structure, such that the natural isomorphism
\[\label{closurenrich}
	\V\text{-}\Cat(\C\boxtimes \E, \D) \cong
	\V\text{-}\Cat(\E, \underline{\VCat}(\C , \D))
\]
holds for $\V$\hyp{}categories $\C, \D,\E$. The $\V$-category $\underline{\VCat}(\C , \D)$ will thus be the internal hom for the closed monoidal structure given by $\boxtimes$.
\begin{definition}[Tensor product of $\V$\hyp{}categories]\label{cosmuclosed}
	Given two $\V$\hyp{}categories $\C, \D$ we define the $\V$\hyp{}category $\C\boxtimes\D$ having
	\begin{itemize}
		\item as objects the set $\C\times\D$, and
		\item as $\V$\hyp{}object of arrows $(C,D)\to (C',D')$ the object
		      \[
			      \C(C,C')\otimes \D(D,D')\in\V.
		      \]
	\end{itemize}
	The free $\V$\hyp{}category $\I$ associated to the terminal category, having a single object $*$ and where $\I(*,*)=I$, the monoidal unit of $\V$, is the unit object for this monoidal structure.
\end{definition}
Various checks are now in order:
\begin{itemize}
	\item $\C\boxtimes\D$ really is a $\V$-category;
	\item $\firstblank\boxtimes\secondblank : \VCat \times\VCat \to\VCat$ is a bifunctor, hat endows $\VCat$ with a monoidal structure.
\end{itemize}
None of these constitutes a conceptual challenge.
\begin{proposition}\label{prop:enricchio}
	The monoidal category $(\VCat,\boxtimes)$ can be promoted to a closed monoidal category, with internal hom denoted $\underline{\VCat}(\firstblank,\secondblank) : \VCat^\opp\times \VCat \to \VCat$.
\end{proposition}
\begin{proof}
	Given $\C,\D\in \VCat$ we define a $\V$\hyp{}category whose objects are $\V$\hyp{}functors $F,G : \C\to\D$ and where (with a little abstraction from \ref{naturalu} to the enriched setting) the $\V$\hyp{}object of natural transformations $F\To G$ is defined via the end
	\[
		\VCat(\C,\D)(F,G) := \int_{C\in\C} \D(FC,GC).
	\]
	In the unenriched case, the end was better understood as the equaliser of a pair of arrows:
	\[
		\int_{C\in\C}\D(FC,GC)\cong \eq\Big( \prod_{C\in \C}\D(FC,GC)\rightrightarrows \prod_{C,C'}\prod_{c\to c'}\D(FC,GC') \Big)
	\]
	In the enriched case, we can consider the same symbol, and re-interpret the product $\prod_{\C(C,C')}$ as a suitable \emph{power} in $\V$:
	\[
		\int_{C\in\C}\D(FC,GC)\cong \eq\Big( \prod_{C\in \C}\D(FC,GC)\rightrightarrows \prod_{C,C'} {\C(C,C')}\pitchfork \D(FC,GC')\Big)
	\]
	(see \cite[\S 2.3]{Gray} for more on this definition). It is now a matter of unwinding the definition to show that a $\V$\hyp{}natural transformation corresponds to a generalised element of $\int_{C\in\C}\D(FC,GC)$; we leave the proof to the reader in Exercise \ref{ex:enrichnat}

	It remains to prove, now, that the isomorphism (\ref{closurenrich}) holds: this is rather easy, since in the above notations, any functor $F : \C\boxtimes \E \to \D$ defines a unique functor $\hat F : \E \to \underline{\VCat}(\C,\D)$: for any two objects $E,E'\in\E$, the collection of arrows $\E(E,E') \to \V(\C(C,C'),\D(F(C,E),F(C',E')) = \C(C,C')\pitchfork\D(F(C,E),F(C',E')$ given by the action of $F$ on hom-objects is a wedge in the pair $(C,C')$. Thus, since
	\[\label{madonna_che_proof}
		\int_{CC'} \C(C,C')\pitchfork\D(F(C,E),F(C',E') \cong \int_C \D(F(C,E),F(C,E')
	\]
	by \ref{ends-are-weighted}, we get a correspondence on objects $\hat F : E \mapsto (\lambda C.F(C,E))$, and a correspondence on hom-objects in the form of maps $\hat F_{E,E'} : \E(E,E') \to \int_C \D(\hat FE(C),\hat FE'(C)) = \underline{\VCat}(\C,\D)(\hat FE, \hat FE')$. The fact that each $\hat F(E)$ is a $\V$-functor $\C\to \D$, and that $\hat F$, so defined, is a $\V$-functor $\E \to \underline{\VCat}(\C,\D)$, are both necessary but tedious checks. As a proof of why we choose the word ``tedious'' to describe the process, let us show the reader the argument proving that the triangle
	\[\label{zio_bubu_la_pruff}
		\vcenter{\xymatrix{
		&I\ar[dr]^{\iota'}\ar[dl]_\iota\ar@{}[d]|{\tiny \fbox{1}} & \\
		\E(E,E) \ar[rr]_-{\hat F_{E,E}} && \int_C \D(F(C,E), F(C,E))
		}}
	\]
	commutes, so that $\hat F$ preserves the identity arrows $\iota : I \to \E(E;E)$ and $\iota' : I \to \int_C \D(F(C,E), F(C,E))$. First, we have to define $\iota'$: the component at $(C,E)$ is obtained from the universal property of the end at codomain, starting from the wedge induced by $F$. Second, it is evident that diagram \eqref{zio_bubu_la_pruff} commutes if and only if the whiskered diagram
	\[
		\vcenter{\xymatrix{
		&I\ar[dr]^{\iota'}\ar[dl]_\iota & \\
		\E(E,E) \ar[rr]_-{\hat F_{E,E}} && \int_C \D(F(C,E), F(C,E)) \ar[d]^{\omega_{CC'}}\\
		&& \C(C,C')\pitchfork\D(F(C,E),F(C',E')
		}}
	\]
	commutes for every choice of $C,C'\in\C$ (we implicitly use \eqref{madonna_che_proof} to describe conveniently the universal wedge). In order to show this last commutativity, consider the diagram
	{\footnotesize \[
		\vcenter{\xymatrix{
		&&\C(C,C')\pitchfork \big(\C(C,C')\otimes \E(E,E)\big)\ar@/^1.5pc/[dr]^-{\;\C(C,C')\pitchfork F_{CC';EE}}\\
		&I\ar@{}[d]^{\tiny \fbox{1}}\ar[dl]^{\iota_E} \ar[r]^-{\eta_I}\ar[dr]^{\iota'_{F(C,E)}} & \ar@{}[ul]|{\tiny \fbox{2}}\ar@{}[r]|-{\tiny \fbox{3}}\ar[u]_{1\pitchfork(1\otimes\iota_E)} \C(C,C')\pitchfork \big(\C(C,C')\otimes I\big) & \C(C,C')\pitchfork \D(F(C,E), F(C',E))\\
		\E(E,E) \ar@/^3pc/[uurr]^{\eta_{\E(E,E)}}\ar[rr]_-{\hat F_{E,E}} && \int_C \D(F(C,E), F(C,E)) \ar@/_1pc/[ur]_{\omega_{CC'}}
		}}
	\]}
	The sub-diagrams \fbox{\tiny 2} and \fbox{\tiny 3} commute, respectively because the unit of the adjunction $\C(C,C)\otimes \firstblank \dashv \C(C,C)\pitchfork \firstblank$ is natural, and because we assumed $F$ was a $\V$-functor. From this, we deduce the desired commutativity.

	The fact that $\hat F$ preserves composition translates into the commutativity of the square
	{\small \[\label{hom_is_a_hom}
		\vcenter{\xymatrix{
		\E(E,E') \otimes \E(E', E'') \ar[d]\ar[r]^{c^\E} & \E(E, E'') \ar[d]^{\hat F_{E,E''}}\\
		\int_A \D(FAE, FAE') \otimes \int_B\D (FBE', FBE'') \ar[r]_-k & \int_C \D(FCE, FCE'')
		}}
	\]}
	where the lower horizontal map $k$ arises from the composition of maps
	\[\label{the_assoc}
		\vcenter{\xymatrix{
				\int_A \D(FAE, FAE') \otimes \int_B\D (FBE', FBE'') \ar[d]^{\omega\otimes\omega}\\
				\D(FAE, FAE') \otimes \D (FAE', FAE'')\ar[d]^\gamma\\
				\int_A \D(FAE, FAE'')
			}}
	\]
	(The last map is the unique induced by the composition law of $\D$.) The proof that diagram \eqref{hom_is_a_hom} commutes is relegated to an exercise in \ref{hom_is_a_cat}.
\end{proof}
The given definition for the enriched end allows us to state an elegant form of the $\V$\hyp{}enriched Yoneda lemma:
\begin{remark}[$\V$-Yoneda lemma]\label{yoneda_enrichio}
	\index{Yoneda lemma!enriched ---}
	Let $\D$ be a small $\V$\hyp{}category, $D\in\D$ an object, and $F : \D\to\V$ a $\V$\hyp{}functor. Then the canonical map
	\[\label{yonedarrow}
		FD\to \underline{\VCat}(\D , \V)(\D(D,\firstblank),F)
	\]
	induced by the universal property of the involved end\footnote{Notice that this is an alternative point of view on the proof of the ninja Yoneda lemma \ref{ninjayo}: the morphism in \eqref{yonedarrow} is induced by a wedge $FD\to \V(\D(D,D'), FD')$ in $D'$, whose members are the mates of the various $\D(D,D')\to \V(FD,FD')$ giving the action of $F$ on arrows.} is a $\V$\hyp{}isomorphism.
\end{remark}
Enriched co\fshyp{}ends can now be defined in the setting of enriched categories, by re-inventing all the initial definitions given in our Chapter 1, and adapting them to the enriched setting. The present section is nothing more than a graphical embellishment of \cite{dubuc1970kan}, where we make a few blanket assumptions for the sake of simplicity of exposition.

The interested reader is warmly invited to look at said text for more information and more general statements.
\begin{notation}
	Our blanket assumption throughout the section is the following: categories are $\V$\hyp{}cotensored (see \cite{dubuc1970kan}: in the absence of cotensors, the enriched counterpart of a co\fshyp{}end is not well-behaved enough to be interesting; in such cases, one loses the equivalent description of a co\fshyp{}end as weighted co\fshyp{}limit, because our \ref{importantlemma} fails to be true).

	Moreover, we sometimes blur the distinction between (di)natural families $D \to T(C,C)$, for a $\V$\hyp{}functor $T : \C^\opp\boxtimes \C \to \D$ and $\V$\hyp{}arrows $I \to \D(D, T(C,C))$. We do this quite liberally especially when drawing commutative diagrams or referring to components of $\V$\hyp{}natural transformations.
\end{notation}
The enriched analogue of extranaturality can be defined as follows:
\index{Extranaturality!enriched ---}
\begin{definition}[Enriched extranaturality]\label{enridina}
	Let $P : \A \boxtimes \B^\opp\boxtimes \B \to \E$ and $Q : \A \boxtimes \C^\opp\boxtimes \C \to \E$ be $\V$\hyp{}functors; an \emph{extranatural transformation} $\alpha P \din Q$ consists of a family of morphisms
	\[
		\alpha_{ABC} : P(A,B,B) \to Q(A,C,C)
	\]
	in $\E$, indexed by the objects of $\A,\B,\C$, such that the following three diagrams made by the action on morphisms of $P,Q$ commute:
	\begin{gather}
		\xymatrix@C=2cm{
		\A(A,A') \ar[d]_{Q(\firstblank,C,C)}\ar[r]^-{P(\firstblank,B,B)}& \E(P(A,B,B), P(A', B,B)) \ar[d]^{\E(1,\alpha_{A'BC})}\\
		\E(Q(A,C,C), Q(A', C,C)) \ar[r]_{\E(\alpha_{ABC}, 1)} & \E(P(A,B,B), Q(A', C,C))} \notag\\
		\vcenter{\xymatrix@C=2cm{
		\B(B,B') \ar[d]_{P(A,\firstblank,B)}\ar[r]^-{P(A,B',\firstblank)}& \E(P(A,B', B), P(A,B',B')) \ar[d]^{\E(1,\alpha_{AB'C})}\\
		\E(P(A,B',B), P(A,B,B)) \ar[r]_{\E(1,\alpha_{ABC})}& \E(P(A,B',B), Q(A,C,C))
		}} \\
		\xymatrix@C=2cm{
		\C(C,C') \ar[d]_{Q(A,\firstblank,C')}\ar[r]^-{A(A,C,\firstblank)}& \E(Q(A,C,C), Q(A,C,C')) \ar[d]^{\E(\alpha_{ABC},1)}\\
		\E(Q(A,C',C'), Q(A,C,C')) \ar[r]_{\E(\alpha_{ABC'}, 1)} & \E(P(A,B,B), Q(A,C,C'))
		}\notag
	\end{gather}
	($A,A',B,B', C,C'$ are objects of the respective categories, and $\E(u,1)$ is the image of $u$ under the functor $\E(\firstblank,E) : \E^\opp\to \V$ as $E$ runs over the objects of $\E$).
\end{definition}
We collect the $\V$\hyp{}extranatural transformations $\alpha : D \din T$ into the object \[\underline{\VCat}_e(\C^\opp\boxtimes \C , \D)(\Delta D, T).\]
\begin{remark}\label{no_enri_dinat}
	\index{Dinaturality!enriched ---}
	\emph{Enriched dinaturality} does not seem to appear in the literature. The scope of the present remark is to show why such a notion is almost often useless. In short, for a generic base of enrichment there is no notion of `constant' $\V$\hyp{}functor, and thus there is no way to define co\fshyp{}wedges as dinatural transformations to/from a constant.
\end{remark}
First, we define the enriched end of a functor taking value in the base of enrichment.
\begin{definition}\label{enriend:uno}
	\cite[I.3.1]{dubuc1970kan} Given a $\V$\hyp{}category $\C$ and a $\V$\hyp{}functor $ : \C^\opp \boxtimes \C \to \V$, the \emph{end} of $T$ is an object of $\V$, denoted $\int_C T(C,C)$, and a $\V$\hyp{}natural family of morphisms $\{\int_C T(C,C) \xto{p_C} T(C,C) \mid C\in \C\}$ such that given any other $\V$\hyp{}natural family $\{u_C : V \to T(C,C) \mid C \in \C\}$ there exists a unique $V \to \int_C T(C,C)$ such that $p_C \circ \bar u = u_C$ in the diagram
	\[
		\vcenter{\xymatrix{
				V\ar[rr]^{\bar u}\ar[dr]_{u_C} && \int_C T(C,C) \ar[dl]^{p_C}\\
				& T(C,C)
			}}
	\]
\end{definition}
The definition of end for a generic codomain is now given representably, following the enriched version of Yoneda-Grothendieck philosophy (see  \ref{yogrophilo}). The enriched Yoneda lemma \ref{yoneda_enrichio} draws the connection between the two: we now define
\begin{definition}\label{enriend:due}
	Let $T : \C^\opp\boxtimes \C \to \D$ be a $\V$\hyp{}functor; the \emph{end} of $T$ is an object $\int_C T(C,C)$ of $\D$ endowed with a $\V$\hyp{}natural family of morphisms $\{\int_C T(C,C) \xto{p_C} T(C,C) \mid C\in \C\}$ such that given any $D\in\D$ the family of arrows
	\[
		\xymatrix@C=2cm{
		\D(D, \int_C T(C,C)) \ar[r]^-{\D(D,p_C)}& \D(D, T(C,C))
		}
	\]
	exhibits the end of $\D(D, T(\firstblank,\secondblank)) : \C^\opp\boxtimes\C \to \D \to \V$.
\end{definition}
\begin{remark} Equivalent to the universal property above is the fact that there is a natural bijection between the $\V$\hyp{}wedges $D \xto{u_C} T(C,C)$ and the underlying set of $\D(D, \int_C T(C,C))$:
	\[ \underline{\VCat}_e(\C^\opp\boxtimes \C , \D)(\Delta D, T) \cong \V\Big(I, \D\Big(D, \int_C T(C,C)\Big)\Big) \]
\end{remark}
\begin{definition}
	Let $T : \C^\opp\boxtimes\C \to \D$ be a $\V$\hyp{}functor; the universal property of $\int_C T(C,C)$ yields a unique morphism $G\big(\int_C T(C,C)\big) \to \int_C GT(C,C)$ for every $\V$\hyp{}functor $G : \D \to \E$; this is the unique morphism closing the diagram
	\[
		\vcenter{\xymatrix{
				G\big(\int_C T(C,C)\big) \ar[dr]_{Gp_C}\ar@{.>}[rr]^\zeta && \int_C GT(C,C) \ar[dl]^{p'_C} \\
				& GT(C,C)
			}}
	\]
	to a commutative one, where $p'_C$ is the terminal wedge of $GT$. We say that $G$ \emph{preserves} the end of $T$ if whenever $\int_C T$ exists, so does $\int_C GT$, and the above comparison morphism $\zeta_G$ is invertible.
\end{definition}
\begin{remark}[Parametric ends of $\V$\hyp{}functors]
	Let $T : \C^\opp\boxtimes\C \boxtimes \E \to \D$ be a $\V$\hyp{}functor; the monoidal closed structure of $\VCat$ gives $T$ a mate $\hat T : \C^\opp\boxtimes\C \to [\E,\D]$. The \emph{parametric end} of $T$, provided it exists, consists of the end of $\hat T$, promoted to a $\V$\hyp{}functor $\E \to \D$.

	More in detail, this means that the parametric end of $T$ exists if for every $E\in\E$ the end of $T(\firstblank,\secondblank;E) : \C^\opp\boxtimes\C \to \D$ exists and there is a unique morphism
	\[\textstyle
		\E(E,E') \to \D\Big(\int_C T(C,C;E), \int_C T(C,C;E')\Big)
	\]
	giving to $\int_C T(C,C;\firstblank)$ the structure of a $\V$\hyp{}functor $\E \to \D$, in such a way that $p_{(C),E} : \int_C T(C,C;E)\to T(C,C;E)$ is a $\V$\hyp{}wedge for $T(\firstblank,\firstblank;E)$, and it is $\V$\hyp{}natural in $E$.
\end{remark}
The above remark is based on the fact that if $\D$ is a complete $\V$-category, namely a $\V$-category such that every limit of $F : \J \to \E$ weighted by $W : \J \to \V$ exists, then $[\E,\D]$ is complete as well, because limits can be computed pointwise: spelled out precisely, this means that given a $\V$-functor $F : \J \to [\E,\D]$, the correspondence
\[E\mapsto \int_J WJ \pitchfork F(J,E)\]
defines a $\V$-functor, and such $V$-functor has the universal property of the limit $\wlim{W}F$.

This statement follows at once from a simple computation with the involved end: it just remains to see that the above definition gives a well-defined $V$-functor.

Let now $\E = \B^\opp\boxtimes \B$; then an enriched Fubini rule holds for functors
\[\notag
	T : \C^\opp\boxtimes\C\boxtimes\B^\opp\boxtimes \B \to \E.
\]
\begin{theorem}[Enriched Fubini rule]\label{enrifubi}
	Whenever both inner parametric ends
	\[
		\int_C \int_B T(C,C;B,B) \qquad\qquad\int_B \int_C T(C,C;B,B)
	\]
	exist as functors $T_\B : \B^\opp\boxtimes \B \to \D$ and $T_\C : \C^\opp\boxtimes \C \to \D$, the outer ends exist if and only if  either one of them exists, and they are canonically isomorphic objects, in turn isomorphic to the end of the rearranged functor $(\B\boxtimes \C)^\opp\boxtimes (\B\boxtimes\C) \to \D$.
\end{theorem}
\begin{exercises}
\item \label{ex4:wcolims} Prove all the statements in \ref{weicolims}.
\item This is a corollary to \ref{importantlemma}, where we take the tensored and cotensored $\V$-category $\A = \V$. Prove that the adjunction $ $ reduces to the adjunction $\Lan_yW\dashv N_W$ of \ref{nervereal}; do it in two ways: first, through coend calculus, and then exhibiting unit and counit of an adjunction $W\otimes \firstblank \dashv \V(W\firstblank,\secondblank) : [\C,\V] \leftrightarrows \V$.
\item \cite{2catlimits} Mimic the argument in \ref{inserters} to give a characterisation based on coend calculus of the following weighted co\fshyp{}limits:
\begin{enumerate}
	\item the \emph{equifier} of a pair of functors $f,g : B \rightrightarrows C$ and two natural transformations $\alpha,\beta : f \To g$; it is defined as the weighted limit $\wlim{W}{F}$ where $W : \J \to \Cat$ sends the 2\hyp{}category
	      \[\notag \J =
		      \xymatrix@C=5mm{
		      0\ar@/^1pc/[rr]\ar@/_1pc/[rr]
		      \rtwocell<\omit>{}&\rtwocell<\omit>{}&1
		      }
	      \]
	      to the diagram of 2\hyp{}cells
	      $
		      \xymatrix@C=3.3mm{
		      [0]\ar@/^1pc/[rrrr]^{s}\ar@/_1pc/[rrrr]_{t} &\rtwocell<\omit>{\delta_0}&\rtwocell<\omit>{\delta_1}&& [1]
		      }
	      $
	      and $F : \J \to \Cat$ sends $\J$ to the diagram of 2\hyp{}cells
	      $
		      \xymatrix@C=3.3mm{
		      B\ar@/^1pc/[rrrr]^f\ar@/_1pc/[rrrr]_g &\rtwocell<\omit>{\alpha}&\rtwocell<\omit>{\beta}&& C
		      }
	      $ and enjoys the following universal property (see \cite[4.5]{2catlimits}): there is a universal diagram
	      \[\notag
		      \xymatrix@C=3.3mm{
		      \wlim{W}{F} \ar[rr]^-p && B\ar@/^1pc/[rrrr]^f\ar@/_1pc/[rrrr]_g &\rtwocell<\omit>{\alpha}&\rtwocell<\omit>{\beta}&& C
		      }
	      \]
	      such that $\alpha * p = \beta * p$ and given any other $q : A \to B$ such that $\alpha * q = \beta * q$ there is a unique $\bar q : A \to \wlim{W}F$ such that $q = p\circ \bar q$ in the diagram
	      \[\notag
		      \vcenter{\xymatrix@C=3.3mm{
		      \wlim{W}{F} \ar[rr]^-p && B\ar@/^1pc/[rrrr]^f\ar@/_1pc/[rrrr]_g &\rtwocell<\omit>{\alpha}&\rtwocell<\omit>{\beta}&& C\\
		      A \ar[u]^{\bar q}\ar[urr]_q
		      }}
	      \]
	      Moreover, given $h,k : A \to \wlim{W}F$ and a 2\hyp{}cell $\mu : ph \To pk$, there is a unique 2\hyp{}cell $\bar\mu$ such that $p *\bar\mu = \mu$.
	\item the \emph{co}equifier of a pair of functors $f,g : A \rightrightarrows B$; it is defined as the weighted colimit $\wcolim{W}{F}$ for the same  $W,F$ above, and enjoys the dual universal property (write it down in detail).
	\item the \emph{lax limit} of a functor $f : A \to B$; it is defined as the weighted limit $\wlim{W}F$ where $W : \{0<1\} \to \Cat$ chooses the functor $\{0\} \xto{0} \{0 \to 1\}$ and $F : \J \to \Cat$ chooses the functor $f : B \to C$. The object $\wlim{W}F$ has the following universal property: there exists a pair $(u,v)$ of 1\hyp{}cells and a 2\hyp{}cell $\lambda : fu\To v$ in a diagram
	      \[\notag
		      \xymatrix{
			      &\dtwocell<\omit>{^\lambda} \wlim{W}F\ar[dr]^v \ar[dl]_u & \\
			      B \ar[rr]_f && C
		      }
	      \]
	      terminal with respect to this property (write down the universal property in detail).
	\item the \emph{pseudo-limit} of a functor $f : A \to B$, where $\J,F$ are the same, and $W$ is instead the embedding of the domain in the generic isomorphism, \ie the functor $\{0\} \xto{0} \{0 \cong 1\}$.
\end{enumerate}
\item Let $\J = \{0 \rightrightarrows 1\}$, let $[n]$ denote the category $\{0 < 1 < \dots < n\}$ and
\begin{gather*}
	W : \J \to \Cat : \{0 \rightrightarrows 1\} \longmapsto \{[1] \underset{d_2}{\overset{d_0}\rightrightarrows} [2]\} \\
	F : \J \to \Cat : \{0 \rightrightarrows 1\} \longmapsto \{ B \underset{g}{\overset{f}\rightrightarrows} C\}
\end{gather*}
where $d_i : [n] \to [n+1]$ avoids the $i^\text{th}$ element. What is the universal property of $\wlim{W}F$?
\item \label{ex:enrichnat} Complete the proof of \ref{prop:enricchio}.
\item \label{ex:cone-again} Let $W : S^0 \hookrightarrow D^1$ be the canonical inclusion of the endpoints $\{0,1\}$ into the interval $[0,1]\subset \bR$ with the usual topology; prove that the mapping cone of a continuous map $f : X\to Y$ regarded as a functor to the category $\Spc$ is precisely the weighted colimit $\wcolim{W}f$.
\item Show that there are canonical isomorphisms $\wlim{W}FJ\cong \wlim{\Lan_JW}F$, in the diagram
\[\notag\xymatrix{
		\A\ar[d]_W\ar[r]^J & \B \ar[r]^F \ar[dl]^{\Lan_JW}& \C \\
		\V
	}\]
and dually $\wcolim{WJ}F\cong \wcolim{W}\Lan_JF$, in the diagram
\[\notag\xymatrix{
		\A\ar[d]_F\ar[r]^J & \B \ar[r]^W \ar[dl]^{\Lan_JF}& \V \\
		\C
	}\]
\item Is there a Fubini rule for weighted co\fshyp{}limits?
\item Use the universal property of \ref{enriend:due} to show that every $\V$\hyp{}natural transformation $\alpha : T \To T' : \C^\opp\boxtimes \C \to \D$ induces an arrow $\int_C \alpha : \int_C T(C,C) \to \int_C T'(C,C)$. Show that if $\alpha$ is a component-wise monomorphism\footnote{We say that $u : D \to D'$ is a $\V$-\emph{monomorphism} if the image of $u$ under the enriched Yoneda embedding $\yon(u) : \yon D \To \yon D'$ is a monomorphism in $\VCat(\D^\opp, \V)$.} in $\D$, then so is $\int_C \alpha$.
\item \label{hom_is_a_cat} Prove that the $\V$-category of $\V$-functors $\C \to \D$ is indeed a $\V$-category, in the sense that the axioms of \ref{enrichcat} hold. More precisely, the claim here is that given $\V$-functors $F,G : \C \to \D$, the $\V$-category with hom-objects $\int_C \D(FC,GC)$ is well-defined and indeed a $\V$-category. Verifying the axioms once and for all is instructive, but things get a little bit hairy.

For example, this is the commutative diagram witnessing that composition has a right identity induced by $\iota^{GX} : I \to \D(GX,GX)$:
\[\notag\vcenter{\xymatrix@C=2cm{
	\big(\int_C \D(FC,GC)\big)\otimes I \ar[d]_{\omega^{FG}_C\otimes I}\ar[r]^{1\otimes \iota^G} & \int_C \D(FC,GC) \otimes \int_B\D(GB,GB) \ar[d]^{\omega^{FG}_C\otimes\omega^{GG}_C}\\
	\D(FC,GC)\otimes I \ar@/_1pc/[dr]_\cong\ar[r]_{1\otimes \iota^{GC}}& \D(FC,GC)\otimes \D(GC,GC) \ar[d]^{\text{comp}}\\
	& \D(FC,GC)
	}}\]
where $\iota^G$ is the unique map such that $\omega^{GG}_C \circ \iota^G = \iota^{GC}$; associativity gets a bit worse, but the appropriate diagram can be reshaped and broken according to the following scheme:
\[\notag\xymatrix{
	& & & {}  \ar[rrd]^{1\otimes \gamma} & & \\
	& {}  \ar@{.>}[rd]  \ar[ldd]_{\omega\otimes\omega\otimes 1}  \ar@{.>}[rr]  \ar@{.>}[rrd]  \ar[rru]^{1\otimes\omega\otimes\omega} & & {} \ar[d] \ar[rr] & & {} \ar[ld]  \ar[dd]^{\omega\otimes\omega} \\
	& & {} \ar[d] \ar[rrd] & {} \ar[rd] \ar[r] & {} \ar[rd] & \\
	{}  \ar[rd]_{\gamma\otimes 1} & & {} \ar[rd] & & {}  \ar[ld]_{\gamma\otimes 1}  \ar[r]^{1\otimes\gamma} & {}  \ar[d]^{\gamma} \\
	& {} \ar[ru]  \ar[rr]_{\omega\otimes\omega} & & {}  \ar[rr]_{\gamma} & & {}
	}\]
where $\gamma$ are suitable composition maps, $\omega$ are suitable cowedge maps, and the dotted arrows are induced as canonical maps, say,
\[\Big(\int_C \D(FC,GC) \Big)\otimes X \to \int_C \big(\D(FC,GC) \otimes X\big).\]

Starting from this, prove that diagram \ref{hom_is_a_hom} commutes.
\item Prove \ref{enrifubi} with the aid of \cite[I.3.4]{dubuc1970kan}.
\item Prove the following: if $T : \C^\opp\boxtimes \C \boxtimes \E \to \D$ is a $\V$\hyp{}functor, we consider $T(C,C;\firstblank) : \E \to \D$ as a $\V$\hyp{}functor, and we let $W : \C \to \V$ be a weight; then
\[\notag
	\int_C \wlim{W}T(C,C;E) \cong \wlim{W} \int_C T(C,C;E).
\]
Dualise for coends.
\end{exercises}

\chapter{Profunctors}\label{sec:profunctors}
\section{The 2\hyp{}category of profunctors}
\begin{abstract}
	The present chapter introduces the theory of \emph{profunctors}; regarded as a generalisation of presheaves and modules over rings, profunctors have a pride of place in 2\hyp{}dimensional algebra. We explore the main features of the bicategory $\Dist$ that they form, heavily employing coend calculus as a mean to handle computations.

	The bicategory of profunctors has plenty of nice features: for example, it is monoidal biclosed, and in fact also \emph{compact} closed. There are canonical embeddings of the strict 2\hyp{}category $\tCat$ in the bicategory $\Dist$, that in a suitable sense `preserve' the expressiveness of the category theory done in both environments.

	The theory of profunctors has several applications in algebra \cite{gambo-joy}, algebraic geometry \cite{kuznetsov2014categorical} algebraic topology, representation theory \cite{tambara2006distributors,pastro2008doubles} and computer science \cite{kmett}.
\end{abstract}
\epigraph{On peut regarder une pièce d'un puzzle pendant trois jours et croire tout savoir de sa configuration et de sa couleur sans avoir le moins du monde avancé : seule compte la possibilité de relier cette pièce à d'autres pièces}{G. Perec --- \emph{La vie, mode d'emploi}}
The lucid presentation in the notes \cite{benabou2000distributors} and in \cite[\S 4]{Cordier2008}, \cite{Bor2} are standard references to follow this chapter. Co/end calculus provides a deeper glance and a unified description for the material appearing in \cite{benabou2000distributors}; here some arguments become neater, some other are made more explicit or computationally evident.

First, recall from \ref{Tensor} and \ref{Module} that we can define tensor product of modules as a coend:
\begin{example*}[The tensor product of modules as a coend]\index{Product!tensor ---}
	Any ring $R$ can be regarded as an $\Ab$\hyp{}enriched category having a single object: under this identification, the category of left modules over $R$ is but the category of functors $R \to \Ab$, and dually, the category of right $R$-modules is the category of contravariant functors, $R^\opp\to \Ab$
	\begin{align}
		\BF{Mod}_R   & \cong \Cat(R^\opp,\Ab)\notag \\
		{}_R\BF{Mod} & \cong\Cat(R,\Ab).
	\end{align}
	Moreover, given functors $ A : R^\opp\to\Ab, B : R\to \Ab$, there is a canonical isomorphism between the functor tensor product $A\otimes B$ defined as the coend of $A\otimes_{\mathbb{Z}} B$, and $A\otimes_R B$: in short, there is a coequaliser diagram
	\[ \vcenter{\xymatrix{
		\bigoplus_{r\in R} A\otimes_{\mathbb{Z}} B \ar@<-.275em>[r]_-{1\otimes r}\ar@<.275em>[r]^-{r\otimes 1} & A\otimes_{\mathbb{Z}} B \ar[r] & A\otimes_R B
		}} \]
\end{example*}
\begin{remark}
	\index{Module}
	\index{Module!bicategory of ---s}
	\index{Bicategory}
	We can define a bicategory (see \ref{bicat})  $\Mod$  having
	\begin{itemize}
		\item 0\hyp{}cells the rings $R,S,\dots$;
		\item 1\hyp{}cells $R \to S$ the modules ${}_RM_S$, regarded as functors $M : R\times S^\opp \to \Ab$;
		\item 2\hyp{}cells $f : {}_RM_S \To {}_RN_S$ are the module homomorphisms $f : M \to N$.
	\end{itemize}
\end{remark}
The notion of profunctor arises from a massive, but fairly natural, generalisation of this construction.

The parallel here is motivated by the fact that categories are certain monoid objects, and the features of such monoids are captured by their \emph{categories of action}: in this perspective a (left) action of a category on a set is merely a functor $\C \to \Set$. The analogy with a group action is once again evident: in fact, a (left) group action of $G$ on a set $X$ is merely a presheaf $G \to \Set$, when $G$ is regarded as a category.

This allows to state the following definition:
\begin{definition}[The bicategory of profunctors]\label{profdef}
	\index{Bicategory!--- of profunctors}\index{_aaa_Prof@$\Dist$}
	There exists a bicategory $\Dist$ having
	\begin{itemize}
		\item 0\hyp{}cells  are (small) categories $\A,\B,\dots$;
		\item 1\hyp{}cells $\proP,\proQ\dots$, denoted as arrows $\proP : \A\pto \B$, are functors
		      \[\A^\opp\times\B\to \Sets;\]
		\item 2\hyp{}cells $\alpha : \proP\Rightarrow\proQ$ are natural transformations.
	\end{itemize}
	\index{Profunctor!composition of ---s}
	Given two contiguous 1\hyp{}cells $\A\stackrel{\proP}{\pto}\B\stackrel{\proQ}{\pto}\C$ we define their composition $\proQ \diamond \proP$ as the coend
	\[\label{proco}
		\proQ \diamond \proP(A,C) := \int^{B\in\B} \proP(A,B)\times\proQ(B,C).
	\]
\end{definition}
\begin{definition}\label{prof_over_anybase}
	\index{Cosmos}
	This definition works well also with $\Sets$ replaced by an arbitrary \emph{Bénabou cosmos} $\V$, \ie in any symmetric monoidal closed, complete and cocomplete category: in this case we speak of $\V$\hyp{}profunctors in the bicategory $\Dist(\V)$.
\end{definition}
\begin{remark}[Naming a category]
	\index{Distributor|see{Profunctor}}
	\index{_aaa_Prof@$\Dist$}
	Since their first introduction, profunctors have been called many other names, depending on the leading perspective that guided their definition:
	\begin{enumtag}{pn}
		\item the 1\hyp{}cells of $\Dist$ are called \emph{profunctors}, because they generalise functors: we will see that some profunctors are \emph{representable}; they are the ones of the form $\B(B, FA)$ or $\B(B,FA)$ for some functor $F : \A \to \B$ between categories. A \emph{pro}functor thus works `on behalf' of a functor (this is one of the meanings of the prefix \emph{pro-}).
		\index{Relator|see{Profunctor}}
		\item In the same vein, \emph{relations} are generalised functions too: this is why some people (among which, A. Joyal) prefer to call the 1\hyp{}cells of $\Dist$ \emph{relators}: just as a func\hyp{}tion is a special kind of rela\hyp{}tion, a func\hyp{}tor is a special kind of rela\hyp{}tor. This analogy is deeper than it seems: we will see in \ref{profs_are_rels} that there is a bicategory of relations, and this is a $\Dist(\V)$ for a certain $\V$.\index{Topos!Grothendieck ---}
		\item Following the idea that \emph{distributions} are generalised functions in functional analysis, the 1\hyp{}cells of $\Dist$ are also called \emph{distributors}, when we follow the intuition that funct\emph{ions} are to funct\emph{ors} as distribut\emph{ions} are to distribut\emph{ors}. As the nLab says,
		\begin{quote}
			Jean Bénabou, who invented the term and originally used “profunctor”, now prefers “distributor”, which is supposed to carry the intuition that a distribut\emph{or} generalises a funct\emph{or} in a similar way to how a distribut\emph{ion} generalises a funct\emph{ion}.
		\end{quote}
		\index{Lawvere}
		\index{W. Lawvere|see{Lawvere}}
		Again, this intuition is deeper that it seems: in \cite{lawvere07cohesion} Lawvere defined a notion of `distribution between toposes', in such a way that the points of a topos $p : \Set \to \E$ behave like Dirac delta functions, and such that distributions between presheaf toposes are exactly profunctors. We discuss the analogy between Dirac deltas and representable functors in \ref{dirac}.
		\index{Correspondence|see{Profunctor}}
		\item The 1\hyp{}cells of $\Dist$ are sometimes called \emph{correspondences}: consider the case when $\V=\{0,1\}$ \ie where $\A, \B$ are sets regarded as discrete categories, and see \ref{profs_are_rels} below.
		\index{Bimodule|see{Profunctor}}
		\item Drawing from \ref{Module}, the 1\hyp{}cells of $\Dist$ are sometimes called \emph{bimodules}: indeed, the bicategory $\Mod$ is precisely the subcategory of $\Dist(\Ab)$ made by one\hyp{}object $\Ab$\hyp{}categories. See \cite{kuznetsov2014categorical} for applications in algebraic geometry, ad \cite{nashphd} for applications of bimodules in homological algebra.
	\end{enumtag}
	\index{Ring}
	Category theorists know well that an elephant can have different names according to the angle it is observed from; we accept this situation, and we tacitly stick to call `profunctor' the 1\hyp{}cells of $\Dist$ without making further mention of this variety of names. However, we invite the reader to maintain clear the intuition conveyed by all these names, in order to appreciate the variety of contexts in which the notion of profunctor naturally arises.
\end{remark}
\index{_aaa_Prof@$\Dist$}
\begin{example}[Profunctors as generalised relations]\label{profs_are_rels}
	We consider categories enriched over the monoidal category $\V =  \{0<1\}$.

	By definition, a profunctor between $\{0,1\}$\hyp{}categories is a function between sets $A^\opp\times B\to \{0,1\}$, or more precisely a function $A\times B \to \{0,1\}$ (a $\{0,1\}$\hyp{}enriched is merely a set, and dualisation on a discrete category has no effect), that is to say a \emph{relation} regarded as a subset $R\subseteq A\times B$.
\end{example}
\index{Lawvere}
The standpoint regarding profunctors as generalised relations is what Lawvere \cite[\S4,5]{LawvereFW:metsgl} calls \emph{generalised logic}, and it regards the coend in \ref{profdef}, as well as the product $\proP(A,X)\times \proQ(X,B)$ therein, as a \emph{generalised existential quantification} and a \emph{generalised conjunction} respectively, giving a composition rule for generalised relations: the coend stands on the same ground of the composition rule for relations (even more: the two constructions have the same universal property), in such a way that the composition of relations $R \subseteq Y\times Z, S \subseteq X\times Y$ is given by the rule
\[
	(x,z)\in R\circ S \iff \exists y\in Y : \big((x,y)\in S\big)\land \big((y,z)\in R\big).
\]
\begin{center}
	\begin{figure}[h!]
		\begin{tikzpicture}
			\matrix (m) [
				matrix of math nodes,
				every node/.append style={
						minimum height=36pt,
						inner sep=6pt,
						execute at begin node=\strut\displaystyle,
					},
			]{
				(x,z)\in S\circ R          & \iff & \exists y\in Y         & \big((x,y)\in S\big)   & \land  & \big((y,z)\in R\big)   \\
				(\proQ\diamond \proP)(X,Z) & =    & \smash{\int^{Y\in\cate Y}} & \proP(Y,Z) & \times & \proQ(X,Y) \\
			};
			\foreach \i in {1,...,6}
			\node[inner sep=-4pt,draw,dashed,gray,fit=(m-1-\i)(m-2-\i)]{};
		\end{tikzpicture}
		\caption{The analogy between the composition of profunctors between categories and the composition of relations between sets gives rise to what Lawvere calls \emph{generalised logic}.}
	\end{figure}
\end{center}
\index{Profunctor!composition of ---s}
\begin{example}
	Let again $A,B$ be sets, considered as discrete categories. A profunctor $\proP : A\pto B$ is then simply a collection of sets $P_{ab}$, one for each $a\in A,b\in B$. Profunctor composition then results in a `categorified' matrix multiplication, since the coend in \eqref{proco} boils down to be a mere coproduct of sets: given $\proP : A\to B$, $\proQ : B\to C$ we have
	\[
		(\proP\diamond \proQ)_{ac} = \coprod_{b\in B}P_{ab}\times Q_{bc}
	\]
	if $P_{ab}=\proP(A,B)$ and $Q_{bc}= \proQ(b,c)$.
\end{example}
\begin{remark}\label{alternative}
	There is an alternative, but equivalent definition for the composite profunctor $\proQ \diamond \proP$ which exploits the universal property of $\Cat(\C^\opp,\Set)$ as a free cocompletion: by definition, the category of profunctors $\proP : \A\pto \B$ and natural transformations fits into an isomorphism of categories
	\[
		\Cat(\A^\opp\times \B ,\Sets)\cong \Cat(\B, \Cat(\A^\opp,\Sets));
	\]
	under this isomorphism, $\proP$ corresponds to a functor $\widehat\proP : \B\to \Cat(\A^\opp,\Set)$ obtained as $B\mapsto \proP(\firstblank,B)$.

	We can thus define the composition $\A\stackrel{\proP}{\pto}\B\stackrel{\proQ}{\pto}\C$ to be $\Lan_\yon \widehat\proP\circ\widehat\proQ$ ($\circ$ is the usual composition of functors):
	\[\label{promate}
		\vcenter{\xymatrix@C=2cm{
		& \B \ar[d]_{\yon }\ar[r]^{\widehat\proP}& \Cat(\A^\opp,\Set) \\
		\C \ar[r]_-{\widehat\proQ}& \Cat(\B^\opp,\Set) \ar@/_1pc/[ur]_{\Lan_\yon \widehat\proP}&
		}}
	\]
\end{remark}
Note that this is in line with the fact that given a category $\C$, presheaves on $\C$ correspond to profunctors $\C\pto\uno$; covariant functors $\C \to \Set$ instead correspond to profunctors $\uno\pto\C$, where $\uno$ is the terminal category.

Equation \eqref{promate} above is equivalent to the previous definition of profunctor composition, in view of the characterisation of a left Kan extension as a coend in $\Cat(\A^\opp,\Set)$ that we have given in \ref{kanend}:
\[
	\Lan_\yon \widehat\proP\cong \int^B\Cat(\B^\opp,\Set)(\yon(B),\firstblank)\otimes \widehat\proP(B).
\]
We have
\begin{align*}
	\Lan_\yon \widehat{\proP}(\widehat{\proQ}(C)) & \cong \int^B \Cat(\B^\opp,\Set)(\yon(B),\widehat{\proQ}(C))\otimes \widehat{\proP}(B) \\
	                                              & \cong \int^B \widehat{\proQ}(C)(B)\otimes \widehat{\proP}(B)                          \\
	                                              & \cong \int^B \proP(\firstblank,B)\times \proQ(B,C).
\end{align*}
\begin{remark}[$\Dist$ is a bicategory]\label{profundefs}
	The properties of (strong) associativity and unitality for the composition of profunctors follow directly from the associativity of cartesian product, its co\hyp{}continuity as a functor of a fixed variable, and from the ninja Yoneda lemma \ref{ninjayo}, as shown by the following computations:
	\begin{itemize}
		\item Composition of profunctors is associative up to isomorphism: given three profunctors $\B\overset{\proH}\pto \cX \overset{\proQ}\pto \cY \overset{\proP}\pto \A$, giving the \emph{associator} of a bicategory structure:
		      \begin{align*}
			      \proP \diamond (\proQ \diamond \proH) & =\int^Y\proP(Y,A)\times (\proQ \diamond \proH)(B,Y)                  \\
			                                            & =\int^Y \proP(Y,A)\times\Big( \int^X \proQ(X,Y)\times\proH(B,X)\Big) \\
			                                            & \cong \int^{XY}\proP(Y,A)\times\Big( \proQ(X,Y)\times\proH(B,X)\Big) \\
			      (\proP \diamond \proQ) \diamond \proH & = \int^X(\proP \diamond \proQ)(X,A)\times\proH(B,X)                  \\
			                                            & \cong \int^{XY}\Big(\proP(Y,A)\times\proQ(X,Y)\Big)\times\proH(B,X)
		      \end{align*}
		      (we freely employ most of the rules of co\fshyp{}end calculus we learned to far, mostly the Fubini rule \ref{fubozzo}) these results are clearly isomorphic, and naturally so.
		\item Every object $\A$ has an identity arrow, given by the hom profunctor $\A(\firstblank,\secondblank)  : \A^\opp\times \A\to \Sets$: the fact that $\proP \diamond \hom\cong\proP$, $\hom \diamond\, \proQ\cong\proQ$ simply rewrites the ninja Yoneda lemma \ref{ninjayo}.
	\end{itemize}
\end{remark}
\begin{remark}\label{isitpentagon}
	The isomorphism above is part of the data turning $\Dist$ into a bicategory; the \emph{associator} realises the identification between different parenthesisations of 1\hyp{}cells, and the \emph{unitor} realises the identification between $\proP\diamond \hom\cong \proP$.

	In order to get a bicategory, some coherence conditions have to be imposed, precisely those of \ref{bicat}.

	It's immediate to observe that the validity of the pentagon identity in the case of the cartesian monoidal structure of $\Sets$, and the naturality thereof, ensure that the associator (whose components are) $(\proP\diamond\proQ)\diamond\proH \Rightarrow \proP\diamond(\proQ\diamond\proH)$ satisfies the pentagon identity; a similar argument shows that the unitor satisfies similar (left and right) triangular identities, as a consequence of the naturality of the ninja Yoneda lemma \ref{ninjayo}.\index{Yoneda lemma!ninja ---}

	The reader might find instructive the exercise of showing by means of the bare universal property of the product that the isomorphisms $A\times (B\times C)\cong (A\times B)\times C$ in a cartesian category satisfy the pentagon equation.
\end{remark}
\begin{notat}[Einstein notation]\label{einstein}
	\index{Einstein notation}
	The reader might have noticed, at this point, that computations with coends might become rather heavy\hyp{}handed, and that sometimes multiple indices have to be considered at the same time; we shall introduce here a useful notation to shorten them a bit, which is particularly evocative when dealing with profunctors; we choose to call it \emph{Einstein convention} for evident reasons.\footnote{To the knowledge of the author, this notation has first been adopted also in \cite{emilyreedy}, a valuable reading in itself for more than one reason. Last but not least, the fact that it provides a description of `Reedy calculus' in homotopical algebra that heavily --albeit without too an explicit mention-- employs coends.}

	Let $\proP : \A\pto \B$, $\proQ : \B\pto\C$ be two composable profunctors. If we adopt the notation $\proP^A_B, \proQ^B_C$ to denote the images $\proP(A,B), \proQ(B,C)\in\Sets$ (keeping track that superscripts are contravariant and subscripts are covariant components), then composition of profunctors acquires the form of a product of tensor components:
	\[\textstyle
		\proP\diamond \proQ(A,C) = \int^B \proP^A_B\times \proQ^B_C = \int^B \proP^A_B \proQ^B_C.
	\]
	The convention is then defined as follows: we denote as super- and subscripts the objects a bifunctor depends on; whenever an object $B\in\B$ appears once covariantly, say in $\proP_B^X$ and once contravariantly, say in $\proQ^B_Y$, the result can be integrated with a co\fshyp{}end into an object $\int^B \proP_B^X \proQ^B_Y$ or $\int_B \proP_B^X \proQ^B_Y$.

	From now on, we freely employ the Einstein summation convention when typesetting long computations.
\end{notat}
\section{Embeddings and adjoints}\label{embare_i_profi}
\index{_aaa_Prof@$\Dist$}
In the present section we shall define two canonical identity\hyp{}on\hyp{}objects embeddings $\tCat\to\Dist$: the \emph{covariant} one reverses the direction of 1\hyp{}cells, but maintains the direction of 2\hyp{}cells; the \emph{contravariant} one does the opposite.\footnote{Recall from chapter 2 that given a 2\hyp{}category $\K$ there is a 2\hyp{}category $\K^\opp$, where 1\hyp{}cells are reversed and 2\hyp{}cells are untouched, and a 2\hyp{}category $\K^\co$, where 2\hyp{}cells are reversed and 1\hyp{}cells are untouched; of course, $\K^\coop = (\K^\opp)^\co = (\K^\co)^\opp$.} The two correspondences are pictorially represented by the assignments
\[\label{direction}
	\xymatrix@R=1.5cm@C=1.5cm{
	\C\ar@{}[r]^{\Cat^\opp\to \Dist}\dtwocell^G_F{^\alpha} & \C && \C\dtwocell^G_F{^\alpha}\ar@{}[r]^{\Cat^\co\to \Dist} & \C \dtwocell^{\,\proP^G}_{\proP^F}{\proP^\alpha} \\
	\D & \D \utwocell^{\proP_F}_{\,\proP_G}{\proP_\alpha}&& \D & \D
	}
\]
and they will be usually nameless 2-functors along the discussion, just called ``the embedding'' of $\tCat$ in $\Dist$.

Given a functor $F : \A \to \B$, we define
\begin{enumtag}{e}
	\item $\proP_F$ to be the profunctor $\B(1,F) : \B \pto \A$ defined sending $(B,A)\mapsto \B(B,FA)$;
	\item $\proP^F$ to be the profunctor $\B(F,1) : \A \pto \B$ defined sending $(A,B)\mapsto \B(FA,B)$.
\end{enumtag}
The 1\hyp{}cells $\proP_F,\proP^F$ will usually be called \emph{representable} profunctors induced by $F : \A \to \B$. If the context forces us to distinguish between $\proP^F$ and $\proP_F$ we call the first the \emph{co}representable one (mnemonic: sending $F$ to its corepresentable profunctor is a 2-functor $\tCat^\co \to\Dist$); this is however uncommon, and which of the two ``representable profunctor induced by $F$'' is usually clear from the context.

Both choices define (pseudo) functors of the appropriate variance, since the isomorphisms
\begin{itemize}
	\item $\proP_{FG}\cong \proP_F \diamond \proP_G$, and $\proP^{FG}\cong \proP^G \diamond \proP^F$ and
	\item $\proP_{\id_{\A}}=\proP^{\id_{\A}}=\A(\firstblank,\secondblank)$
\end{itemize}
can be easily established using elementary coend calculus, and they hold for every composable pair of functors $F,G$ and every category $\A$: for example,
\begin{align*}
	(\proP_F \diamond \proP_G)(C,A) & \cong \int^B \B(B,FA)\times \C(C,GB) \\
	                                & \cong \C(C,GFA) = \proP_{FG}(C,A)
\end{align*}
by the ninja Yoneda lemma.

Natural transformations $\alpha : F\Rightarrow G$ are obviously sent to 2\hyp{}cells in $\Dist$, and the covariance of this assignment is evident from the definition of $\proP^F$ and $\proP_F$.

Finally, the Yoneda lemma entails that both functors $F\mapsto \proP^F$ and $F\mapsto \proP_F$ are locally fully faithful: this means that there are bijections
\[\tCat(\A,\B)(F,G)\cong \Dist(\B,\A)(\proP_F,\proP_G)\]
for every $F,G: \A \to\B$. If $\proP^F\cong \proP^G$, it means we have a natural bijection $\B(FA,B)\cong \B(GA,B)$, that entails the existence of a natural isomorphism $\alpha : F \To G$.

\smallskip
The most important fact about representable and corepresentable profunctors associated to the same functor is that they are \emph{adjoint} 1-cells in $\Dist$:
\begin{remark}\label{theyre.adjoints}
	There is a tight link between the 1\hyp{}cells $\proP_F$, $\proP^F$: they are \emph{adjoint} 1\hyp{}cells in the bicategory $\Dist$. Indeed, for every $F\in \Cat(\A,\B)$ we can define 2\hyp{}cells
	\[
		\vcenter{\xymatrix@R=0cm{
		\epsilon_F : \proP_F \diamond \proP^F \ar@{=>}[r] & \B(\firstblank,\secondblank)\\
		\eta_F : \A(\firstblank,\secondblank)\ar@{=>}[r] & \proP^F \diamond \proP_F
		}}
	\]
	\index{Counit!--- of $\proP_F \dashv \proP^F$}
	\index{Unit!--- of $\proP_F \dashv \proP^F$}
	(\emph{counit} and \emph{unit} of the adjunction), satisfying the zig\hyp{}zag identities: we choose a very explicit, set\hyp{}theoretic proof for this statement, leaving to the reader the task to generalise the argument to a generic base of enrichment.
	\begin{itemize}
		\item The components of the counit $\epsilon$ are determined by inspecting the coend $\proP_F \diamond \proP^F$ as the quotient set
		      \[
			      \int^X \B(A, FX)\times \B(FX,B) = \Big(\coprod_{x\in \A} \B(A, FX)\times \B(FX,B)\Big)/\!\simeq
		      \]
		      where $\simeq$ is the equivalence relation generated by $\big(A\xto{u}FX,FX\xto{v}B\big)\simeq\big(A\xto{u'}FY,FY\xto{v'}B\big)$ if there is $t : X\to Y$ such that $v'=Ft\circ v$ e $Ft\circ u=u'$. This can be visualised as the commutativity of the square
		      \[
			      \xymatrix@R=3mm@C=8mm{
			      & FX \ar@{.>}[dd]^{Ft}\ar[dr]^v& \\
			      A \ar[dr]_{u'}\ar[ur]^u&& B \\
			      & FY\ar[ur]_{v'} &
			      }
		      \] Now it's easily seen that sending $\big(A\xto{u}FX,FX\xto{v} B \big)$ in the composition $v\circ u$ descend to the quotient with respect to $\simeq$, hence there is a well defined map $\epsilon : \proP_F \diamond \proP^F \to \B(\firstblank,\secondblank)$. All boils down to the fact that the composition
		      \[
			      c : \B(A, FX)\times \B(FX,B)\to \B(A,B)
		      \]
		      defines a cowedge in the variable $X$.
		\item The components of the unit $\eta$ are determined as functions
		      \[
			      \xymatrix{
				      \eta : \A(A,B)\ar@{=>}[r] & \proP^F \diamond \proP_F(A,B).
			      }
		      \]
		      Since $\proP^F\ccirc\proP_F(A,B) = \int^X \B(FA,X)\times\B(X,FB)\cong \B(FA,FB)$ (as a consequence of the ninja Yoneda lemma), these are simply determined by the action of $F$ on arrows, $\A(A,B)\to \B(FA,FB)$.
	\end{itemize}
	We now have to verify that the triangle identities (see \cite[Thm\@.3.1.5.(2)]{Bor1}) hold:
	\begin{gather*}
		(\proP^F\ccirc\epsilon)\circ (\eta\ccirc \proP^F) = \id_{\proP^F}\\
		(\epsilon\ccirc \proP_F)\circ(\proP_F\ccirc \eta) = \id_{\proP_F}
	\end{gather*}
	As for the first, we must verify that the diagram
	\[
		\xymatrix@C=1.5cm{
		\proP^F \ar[r]^\sim\ar@{=}[d] & \B(\firstblank,\secondblank)\ccirc\proP^F \ar[r]^{\eta\ccirc\proP^F}& (\proP^F\ccirc\proP_F)\ccirc\proP^F \ar[d]^\cong\\
		\proP^F & \ar[l]^\sim \proP^F\ccirc\B(\firstblank,\secondblank) & \ar[l]^{\proP^F\ccirc\epsilon} \proP^F\ccirc(\proP_F\ccirc\proP^F)
		}
	\]
	commutes. One has to send $h\in \proP^F(u,v)=\B(Fu,v)$ in the class $[(\id_u,h)]\in \int^X\B(U,X)\times \B(FX,V)$, which must go under $\eta\ccirc \proP^F$ in the class $[(F(\id_u),h)]\in \int^{XY}\B(FA,X)\times\B(X,FY)\times \B(FY,B)$, canonically identified with $\int^Y\B(FA,FY)\times \B(FY,B)$. Now $\proP^F\ccirc\epsilon$ acts composing the two arrows, and one obtains $F(\id_A)\circ h=h$ back.

	Similarly, to prove the second identity, the diagram
	\[
		\xymatrix@C=1.5cm{
		\proP_F \ar[r]^\sim\ar@{=}[d]& \proP_F\ccirc\B(\firstblank,\secondblank) \ar[r]^{\proP_F\ccirc\eta}& \proP_F\ccirc(\proP^F\ccirc\proP_F) \ar[d]^\cong\\
		\proP_F & \ar[l]^\sim\B(\firstblank,\secondblank)\ccirc\proP_F & \ar[l]^{\epsilon\ccirc\proP_F}(\proP_F\ccirc\proP^F)\ccirc\proP_F
		}
	\]
	must commute (all the unlabeled isomorphisms are the canonical ones). This translates into
	\[
		\xymatrix{
			\Big(a\xto{u}FB\Big) \ar@{|->}[r]& (u, \id_B)_\sim \ar@{|->}[r]& (u, F(\id_B)) \ar@{|->}[r]& u\circ F(\id_B)=u,
		}
	\]
	which is what we want; hence $\proP_F\dashv \proP^F$. \qed
\end{remark}
\begin{remark}
	Two functors $F : \A\leftrightarrows \B : G$ are adjoints if and only if $\proP_F\cong \proP^G$ (and therewith $G\dashv F$) or $\proP_G\cong \proP^F$ (and therewith $F\dashv G$).
\end{remark}
The covariant embedding of $\tCat$ in $\Dist$ gives every 1\hyp{}cell $F :\A\to \B$ a left adjoint. It turns out that many properties of functors can be translated into properties of the associated representable profunctor. As an elementary example, we take the fact that a functor is full and faithful:
\begin{remark}
	It is a well\hyp{}known fact (see \cite[dual of Prop. 3.4.1]{Bor1}) that if $F\dashv G$, then $F$ is fully faithful if and only if the unit of the adjunction $\eta : 1\to GF$ is an isomorphism.

	This criterion can be extended also to functors which do not admit a `real' right adjoint, once noticed that $F$ is fully faithful if and only if $\A(A,B)\cong \B(FA,FB)$ for any two $A,B\in \A$, \ie if and only if the unit $\eta : \hom_{\A}\Rightarrow \proP^F\ccirc\proP_F$ is an isomorphism in the bicategory $\Dist$.
\end{remark}
It is possible to define an \emph{action} of profunctors on functors; more precisely,
\begin{definition}
	Given a profunctor $\proP : \A\pto \B$ and a functor $F : \B\to \D$ we can define $\proP\otimes F$ to be the functor $\A\to \D$ given by $\Lan_\yon F\circ\widehat{\proP}$, where $\widehat\proP : \B\to\Cat(\A^\opp,\Set)$ is the mate of $\proP$, as in \eqref{promate}.

	More explicitly,
	\[
		\proP\otimes F(A) = \int^B [\yon_B,\proP(\firstblank,A)]\otimes FB\cong \int^B \proP^B_A\otimes F_B
	\]
	Note that a conceptual definition for this operation is the following: $\proP\otimes F$, evaluated on $A\in\A$, is the weighted colimit $\wcolim{\proP_A}F$. Exploiting this definition, it is possible to prove using coend calculus that this operation can be rightly called a left action $\Dist(\A,\B)\times\tCat(\B,\D) \to \tCat(\A,\D)$.
	\begin{itemize}
		\item $\hom_{\B}\otimes F\cong F$ as a consequence of the ninja Yoneda lemma;
		\item If $\C\stackrel{\proQ}{\pto}\A\stackrel{\proP}{\pto} \B\xto{F}\cX$, then $(\proP\diamond \proQ)\otimes F\cong \proQ\otimes(\proP\otimes F)$: indeed
		      \begin{align*}
			      [(\proP\diamond \proQ)\otimes F] & = \int^B (\proP\ccirc\proQ)^B\times F_B                                       \\
			                                       & \cong \int^{BX}\proP^B_X\times \proQ^X\times F_B                              \\
			                                       & \cong \int^X\proQ^X \times\Big( \int^B\proP^B_X \times F_B \Big)              \\
			                                       & \cong \int^X\proQ^X \times (F\otimes\proP)_X = [\proQ\otimes(\proP\otimes F)]
		      \end{align*}
	\end{itemize}
\end{definition}
The composition of profunctors is a `closed structure' in the following sense: both functors $\proP\diamond \firstblank$ and $\firstblank\diamond \proQ$ have right adjoints, respectively given by a right lifting and a right extension operation (see \ref{specchietto_dei_lift}).
\begin{example}[Kan extensions in $\Dist$]\index{Kan extension!--- in $\Dist$}\label{rans}
	Every profunctor $\proP$ has a right Kan extension $\Ran_\proP$ in the sense that the notion has in any bicategory, where composition of functors or natural transformations is replaced by composition of 1- or 2\hyp{}cells.

	One has the following chain of isomorphisms in $\Dist$:
	\begin{align*}
		\Dist(\proG\ccirc\proP,\proH) & \cong \int_{AB}\Sets\big( (\proG\ccirc \proP)(A,B),\proH(A,B) \big)             \\
		                              & \cong \int_{AB} \Sets\Big( \int^X \proG(A,X) \times \proP(X,B) ,\proH(A,B)\Big) \\
		                              & \cong \int_{ABX}\Sets\big( \proG(A,X),\Sets(\proP(X,B),\proH(A,B)) \big)        \\
		                              & \cong \int_{AX}\Sets\Big( \proG(A,X),\int_B\Sets(\proP(X,B),\proH(A,B))\Big)    \\
		                              & \cong \int_{AX}\Sets\Big( \proG(A,X),\langle\proP/\proH\rangle(A,X)\Big)        \\
		                              & \cong \Dist(\proG,\langle\proP/\proH\rangle)
	\end{align*}
	when we define $\langle\proP/\proH\rangle(A,X)$ to be the set of natural transformations $\proP(X,\firstblank)\To\proH(A,\firstblank)$. This yields that $\langle\proP/\proH\rangle$ has the universal property of the right extension of $\proH$ along $\proP$.

	Similarly, we can prove that a suitable end $\langle \proG\backslash \proH\rangle$ defines the right lifting of $\proH$ along $\proG$. This is the content of Exercise \ref{riftuzzo}.
\end{example}
\section{The structure of $\Dist$}\label{struct_of_Prof}
The bicategory of profunctors has an extremely rich structure; we briefly account about its main properties, drawing from \cite{concurren}.

As \ref{profs_are_rels} shows, a relation $R \subseteq A \times B$, regarded as a function $A\times B \to \{0,1\}$, is a 1\hyp{}cell of $\Dist(\{0,1\})$; but what are domain and codomain of such a $R$? It is evident that the symmetry of the product $\sigma_{AB} : A\times B \cong B\times A$ induces an identity\hyp{}on\hyp{}objects identification
\[ \Dist(\{0,1\}) \cong \Dist(\{0,1\})^\opp \]
(see \ref{prof_over_anybase} for the notation $\Dist(\V)$) sending a relation $R \subseteq A\times B$, seen as a 1\hyp{}cell $R : A \pto B$ to the relation $\sigma_{AB}\circ R \subseteq B\times A$, regarded as a 1\hyp{}cell $R : B\pto A$.

This is a general phenomenon: every $\Dist(\V)$ admits a `tautological' bi\hyp{}equivalence with $\Dist(\V)^\opp$, determined as follows (see \cite{concurren}):
\begin{proposition}[The canonical dualiser of $\Dist$]\label{dualiseur}
	The identification $\Cat(\A^\opp\times \B,\Set)\cong \Cat(\B\times\A^\opp,\Set)$ determines an equivalence $(\firstblank)^\circ : \Dist^\opp\to \Dist$ determined by the following correspondences:
	\begin{itemize}
		\item On objects, a 0\hyp{}cell $\A$ of $\Dist$ goes to $\A^\circ=\A^\opp$;
		\item The 1\hyp{}cell $\proP : \A \pto \B$ goes to `itself', $\proP^\circ : \B^\opp\pto \A^\opp$, under the identification (isomorphism of categories)
		      \[
			      \Dist(\A,\B) = \Cat(A^\opp\times\B,\Set) \cong \Cat(\B\times\A^\opp,\Set) = \Dist(\B^\opp,\A^\opp)
		      \]
		\item The 2\hyp{}cell $\alpha : \proP \To \proQ$ goes to `itself', again under the identification
		      \[
			      \Dist(\A,\B)(\proP,\proQ)\cong \Dist(\B^\opp,\A^\opp)(\proP^\opp,\proQ^\opp)
		      \]
	\end{itemize}
\end{proposition}
We pair the existence of a \emph{dualiser} on $\Dist$ with the following result:
\begin{proposition}\label{prof_is_monoidal}
	\cite[4.3.7]{concurren} The bifunctor
	\[ \firstblank\times\secondblank : \Dist \times \Dist \to \Dist \]
	defined as the product functor in each dimension, \ie the correspondence sending
	\begin{itemize}
		\item two 0\hyp{}cells $\A,\B \in\Dist$ into the product category $\A\times \B$;
		\item two 1\hyp{}cells $\proP : \A \pto \B, \proQ : \C \pto \D$ to
		      \[ \proP\times\proQ : \A \times \C \pto \B \times \D : (A,C;B,D) \mapsto \proP(A,B)\times \proQ(C,D),\]
		\item two 2\hyp{}cells $\alpha : \proP \To \proQ$ and $\beta : \proL \To \proH$ to their product components
		      \[ (\alpha\times\beta)_{A,C;B,D} : \proP \times\proL(A,C;B,D)\To \proQ\times\proH(A,C;B,D)\]
		      equips $\Dist$ with a monoidal structure.
	\end{itemize}
\end{proposition}
\begin{remark} (\cite[4.3.3]{concurren})
	A somewhat odd feature of $\Dist$ is that (pseudo\hyp{})products and (pseudo\hyp{})coproducts both exist, they coincide on objects, and exchange each other's place on 1\hyp{}cells:
	\begin{itemize}
		\item the span
		      \[\A \overset{p_\A}\otp \A \coprod\B\overset{p_\B}\pto \B \]
		      exhibits the universal property of the pseudo\hyp{}product $\A\mathbin{\&}\B$ of $\A$ and $\B$ in $\Dist$, if the `projection' maps are defined as $\bsmat \hom_\A \\ \varnothing\esmat : (\A^\opp\times\A)\amalg (\B^\opp\times\A) \to \Set$ and $\bsmat \varnothing \\ \hom_\B \esmat : (\A^\opp\times\B)\amalg (\B^\opp\times \B) \to \Set$;
		\item the cospan
		      \[\A \overset{i_\A}\pto \A \coprod\B\overset{i_\B}\otp \B \]
		      exhibits the universal property of the pseudo\hyp{}coproduct $\A\oplus\B$ of $\A$ and $\B$ in $\Dist$, if the `injection' maps are defined as $\bsmat \hom_\A \\ \varnothing\esmat : (\A^\opp\times\A)\amalg (\A^\opp\times\B) \to \Set$ and $\bsmat \varnothing \\ \hom_\B \esmat : (\B^\opp\times\A)\amalg (\B^\opp\times \B) \to \Set$;
	\end{itemize}
	Note in particular that the monoidal structure in \ref{prof_is_monoidal} is not the cartesian one. Since the empty diagram carries no information on 1- and 2\hyp{}cells, it is thus possible to use this construction to show that $\Dist$ has a pseudo\hyp{}zero object, namely the initial category $\boldsymbol{\varnothing}$. This appears as \cite[4.3.5]{concurren}.
\end{remark}
\index{Category!compact closed ---}
\index{Compact closed category}
\index{_aaa_Prof@$\Dist$}
\begin{remark}
	Taken together, these results allow to prove that the bicategory $\Dist$ is \emph{compact closed}: given that the monoidal unit for the structure in \ref{prof_is_monoidal} is the terminal category, we shall find suitable profunctors
	\[
		\eta : 1 \pto \A^\opp\times \A  \qquad
		\epsilon : \A \times \A^\opp \pto 1
	\]
	such that the following two \emph{yanking equations} are satisfied:
	\begin{gather}
		A\cong A\times I\xto{\A\times\eta}\A\times (\A^\circ\times \A)\cong (\A\times \A^\circ)\times \A\xto{\epsilon\times \A} I\times \A\cong \A\\
		\A^\circ\cong I\times \A^\circ\xto{\eta\times \A^\circ}(\A^\circ\times \A)\times \A^\circ\cong \A^\circ\times (\A\times \A^\circ)\xto{\A^\circ \times\epsilon} \A^\circ\times I\cong \A^\circ
	\end{gather}
	It is evident that suitable hom functors (in fact, the same hom\hyp{}functor of $\A$ and its dualised copy $\hom^\circ$) do the job.
\end{remark}
\index{Category!Traced ---}\index{Traced category}
\begin{remark}\label{def:traced}
	The coend operation $\int^\C : \Cat(\C^\opp\times \C , \Set)$ endows $\Dist$ with a \emph{traced structure}: this means that the category $\Dist$ is monoidal with respect to the product defined in \ref{prof_is_monoidal}, and that there exists a family of functions
	\[ \tr^\C_{\A,\B}:\Dist(\A\times \C,\B\times \C)\to\Dist(\A,\B)\]
	satisfying suitable coherence conditions, expressed in \cite{verty1996traced}: we record these condition in Exercise \ref{ex:traced} and leave the proof to the reader.

	The map $\tr^\C_{\A,\B}$ is defined as the composition
	\begin{align}
		\Cat((\A\times \C)^\opp\times \B\times \C,\Set) & \cong
		\Cat(\A^\opp\times\B, \Cat(\C^\opp\times\C,\Set)) \notag                                   \\
		                                                & \xto{\int^\C} \Cat(\A^\opp\times\B,\Set)
	\end{align}
	(a few completely straightforward identifications have been omitted).
	\section{A more abstract look at $\Dist$}
\end{remark}
\index{_aaa_Prof@$\Dist$}
\index{Bicategory}
\begin{remark}
	The bicategory of profunctors can be promoted to a \emph{multi\hyp{}bicategory} in the sense of \cite[1.4]{cockett2003morphisms}; this means that we exploit the composition operation to specify a class of multi\hyp{}morphisms $\eta : \proP_1,\dots, \proP_n \pto \proQ$, depicted as diagrams
	\[
		\vcenter{\xymatrix{
		X_0 \ar@/_4pc/[rrr]|-@{|}_{\proQ}
		\ar[r]|-@{|}^{\proP_1} & X_1 \ar[r]|-@{|}^{\proP_2} & \dots \ar[r]|-@{|}^{\proP_n} & X_n
		\ar@{=>}(22.5,-5);(22.5,-10)_{\eta}
		}}
	\]
	whose composition, associativity, and unitality follow at once from pasting laws for 2\hyp{}cells in 2\hyp{}categories \cite{kelly1982basic} (try to outline them explicitly as a straightforward exercise). More generally, the bicategory of profunctors is a toy example of a \emph{fc\hyp{}multicategory} \cite{leinster1999fc} or, in more modern terms, a \index{Virtual double category}\index{Category!virtual double ---}\emph{virtual double category} (vdc) \cite{cruttwell2010unified}.

	A vdc $\bC$ as defined in \cite{cruttwell2010unified} is a category\hyp{}like structure whose `cells' $\alpha$ have the form
	\[
		\xymatrix{
		X_0\ar[d] \ar[r]|-@{|}^{p_1}\ar@{}[drrr]|{\Downarrow \alpha}& X_1\ar[r]|-@{|}^{p_2} & \cdots \ar[r]|-@{|}^{p_n} & X_n \ar[d] \\
		Y_0 \ar[rrr]|-@{|}_{q} & & & Y_1
		}
	\]
	where the vertical arrows are the morphisms of a category $\bC_v$ called the \emph{vertical category} of the vdc, and the cell $\alpha$ has a $n$\hyp{}tuple (for $n\ge 0$, the non-negative integer $n$ is called the \emph{arity} of the cell) of \emph{horizontal} morphisms $(p_1, p_2,\dots, p_n)$ as horizontal domain, and a single $q$ as horizontal codomain.

	These cells are subject to certain straightforward coherence conditions of associativity: first of all, vertical morphisms compose like they do in a category; second, the cells can be \emph{grafted} as follows: every
	\[\vcenter{\xymatrix{
		\ar[r]|-@{|}^{p_{11}}\ar@{}[drr]|{\Downarrow\alpha_1}\ar[d]& \cdots\ar[r]|-@{|}^{p_{1n_1}} &\ar@{}[r]|{\cdots}\ar[d] & \ar[r]|-@{|}^{p_{m1}}\ar[d]\ar@{}[drr]|{\Downarrow\alpha_m}& \cdots\ar[r]|-@{|}^{p_{m,n_m}}  & \ar[d]\\
		\ar[d]\ar[rr]|-@{|}_{q_1} \ar@{}[drrrrr]|{\Downarrow\beta}& &\ar@{}[r]|{\cdots} & \ar[rr]|-@{|}_{q_m} &&  \ar[d]\\
		\ar[rrrrr]|-@{|}_r &&&& &
		}}
	\]
	can be composed into a single cell
	\[ \vcenter{\xymatrix{
		\ar[r]|-@{|}^{p_{11}}\ar[dd]\ar@{}[ddrrrrr]|{\Downarrow \beta\odot (\alpha_1,\dots,\alpha_m)}& \cdots\ar[r]|-@{|}^{p_{1n_1}} &\ar@{}[r]|{\cdots} & \ar[r]|-@{|}^{p_{m1}}& \cdots\ar[r]|-@{|}^{p_{m,n_m}}  & \ar[dd]\\
		&&&&& \\
		\ar[rrrrr]|-@{|}_r &&&& &
		}} \]
	so that this operation is compatible with the rest of the data. The \emph{0\hyp{}ary} cells are determined by diagrams like $\scriptstyle \vcenter{\xymatrix@C=2mm@R=2mm{ & X\ar[dr]\ar[dl] & \\Y_0 \ar[rr]|-@{|} && Y_1}}$.

	The standard choice to make $\Dist(\V)$ a vdc is to take
	\begin{itemize}
		\item as objects the small categories;
		\item as vertical category the category $\Cat$;
		\item as $n$\hyp{}ary cells $\alpha : (p_1,\dots,p_n)\To q $ having vertical domain $F : X_0 \to Y_0$ and vertical codomain $G : X_n \to Y_1$ the natural transformations
		      \[ \alpha : \proP_G \diamond \proP_n \diamond\cdots\diamond\proP_1 \To \proQ \diamond \proP_F \]
	\end{itemize}
	The reader will routinely check that all coherences stated in \cite{cruttwell2010unified} hold.

	Of course, whenever we consider composable profunctors, every $n$\hyp{}ary cell can be reduced to a 1\hyp{}ary cell by composing the $n$\hyp{}tuple of horizontal arrows that form its horizontal domain; but in situations where the composition does not exist, for example for categories enriched in a non\hyp{}cocomplete base, it is still possible to define the vdc of these generalised profunctors.

	The notion of vdc can thus be considered as the `correct' generalisation of categories $\Dist(\V)$ enriched over a generic base, and exhibits many of its nice features even without horizontal compositions. Examples of vdcs abund inside and outside category theory: the reader is invited to consult \cite{cruttwell2010unified}.
\end{remark}
The following definition is modeled on the behavior of the canonical embedding of $\tCat$ into $\Dist$: it appears in \cite{wood1982abstract}.
\begin{definition}[Proarrow equipment]\label{def-di-equi}\index{Proarrow equipment}
	Let $p^* : \A \to \M$ be a 2\hyp{}functor between bicategories; $p^*$ is said to \emph{equip $\A$ with proarrows}, or to be a \emph{proarrow equipment for $\A$} if
	\begin{enumtag}{pe}
		\item \label{pe:due} $p^*$ is locally fully faithful;
		\item \label{pe:ter} for every arrow $f\in\A$, $p^*f$ has a right adjoint in $\M$.
	\end{enumtag}
\end{definition}
It is clear how the embedding $\tCat^\opp \to \Dist$ equips $\tCat$ with proarrows (see \eqref{direction}; our choice for the direction of 1\hyp{}cells in $\Dist$ forces us to treat the embedding as contravariant).\footnote{A rather unexpected tight connection between vdcs and proarrow equipments is the following: equipments arise precisely as those vdcs where horizontal arrows can all be composed, and have local horizontal identities for every object. This is \cite[§7]{cruttwell2010unified} and in particular its Definition 7.6.}
\begin{remark}[\upeyes Exact squares and profunctors]
	\index{Exact square}
	\index{Carré exact|see{Exact square}}
	Let us consider $\Cat$ as a 2\hyp{}category; a square
	\[
		\vcenter{\xymatrix{
		\A \ar[r]^T \ar[d]_S \ar@{}[dr]|{\Nearrow\alpha} & \cY \ar[d]^V \\
		\cX \ar[r]_U & \B
		}}
	\]
	filled by a 2\hyp{}cell $\alpha$ will be called a \emph{carré} (see \cite{guitart1980relations}). Any carré induces a 2\hyp{}cell
	\[
		\xymatrix{
		\A \ar@{}[dr]|{\Searrow\alpha^\flat}\ar[d]|-@{|}_{\proP_S} & \ar[l]|-@{|}_-{\proP^T} \cY\ar[d]|-@{|}^{\proP_V} \\
		\cX & \B\ar[l]|-@{|}^{\proP^U}
		}
	\]
	defined by the universal property of profunctor composition via the cowedge
	\[
		\alpha^\flat_{(A),XY} : \cX(X,SA) \times \cY(TA,Y) \xto{\qquad} \B(UX, VY)
	\]
	sending the pair $(f,g)=\left(\var{X}{SA},\var{TA}{Y}\right)$ into $Vg \circ \alpha_A \circ Uf$.

	We say that a carré is \emph{exact} if $\alpha^\flat$ is invertible, \ie if $\int^A \cX(X,SA) \times \cY(TA,Y) \cong \B(UX, VY)$.

	\cite{guitart1980relations} observes that there is a criterion for a carré to be exact that does not involve profunctors: let us consider a square as above, and the induced diagram
	\[
		\xymatrix{
		\A\ar@/^1pc/[drr]\ar@/_1pc/[ddr] \ar[dr]^W && \\
		&(U\downarrow V) \ar[r]^{d_1} \ar@{}[dr]|{\Nearrow a} \ar[d]_{d_0}& \cY \ar[d]^V \\
		&\cX \ar[r]_U & \B
		}
	\]
	where $(U\downarrow V)$ is the \emph{comma category} of the cospan $\cX \xto{U} \B \xot{V}\cY$ and $W$ is the unique induced functor to $(U\downarrow V)$. For each pair of functors $\cX \xto{P} \cZ \xot{Q} \cY $ we have the following identifications of sets:
	\[\scriptstyle
		\vcenter{\xymatrix@C=5mm{
		\Cat(U\downarrow V, \cZ)(Pd_0, Qd_1) \ar[r]\ar@{=}[d] & \Cat(\A,\cZ)(PS, QT)\ar@{=}[dd]\\
		\Dist(\cY,\cX)(\proP_{d_0}\ccirc\proP^{d_1}, \proP^P\ccirc \proP_Q) \ar@{=}[d]& \\
		\Dist(\cY,\cX)(\proP^U\ccirc\proP_V, \proP^P\ccirc \proP_Q)\ar[r] &  \Dist(\cY,\cX)(\proP_A\ccirc \proP^T, \proP^P\ccirc \proP_Q)\\
		}}
	\]
	where the horizontal arrows are induced by whiskering. The carré is exact if and only if for each $P,Q$ there is a bijection
	\[
		\Cat(\A,\cZ)(PS,QT)\cong\Cat(U\downarrow V, \cZ)(Pd_0, Qd_1) .
	\]
\end{remark}
\begin{remark}[Displayed category]\label{dispcat}\index{Category!displayed ---}\index{Displayed category}
	The bicategory of profunctors appears in a generalised form of Grothendieck construction, as described in of \ref{thm:equconfib}. There, we draw an equivalence of categories
	\[
		p_{\firstblank} : \DFib(\C) \leftrightarrows \Cat(\C^\opp,\Set) : \elts{\C}{\firstblank}
	\]\index{Category!--- of elements}
	between discrete opfibrations over $\C$ and presheaves over $\C$ (see \ref{def:dfib}). In recent years it has become popular to call the pair of adjuncts that realises this equivalence respectively `straightening' and `unstraightening'.\index{Straightening}\index{Unstraightening|see{Straightening}}

	It is a good question to ask, if there's some kind of un/straightening for generic functors $q : \E \to \C$ over a base, that do not satisfy \ref{def:dfib}. Sure, since the fibration condition on $p : \E \to \C$ is meant exactly to ensure that each morphism $C \to C'$ in the base induces a function\fshyp{}functor between the fibers $p^\leftarrow C' \to p^\leftarrow C$, for a generic functor $q : \E \to \C$ the straightened $\Gamma_q : \C^\opp \to \K$ won't be strictly functorial, nor it will ensure functorial correspondences among the fibers.

	We can however show precisely how much regularity is lost when passing from discrete fibrations over $\C$ to generic categories lying over $\C$, proving the following result (to the best of the author's knowledge, this was first observed by J. Bénabou \cite{benabou2000distributors}).\index{_aaa_gamma@$\Gamma$}
\end{remark}
\begin{theorem}
	There is an equivalence of categories
	\[
		\Pi(\firstblank) : \Cat_{l,1}(\C,\Dist) \to \Cat/\C : \Gamma
	\]
	\index{Functor!lax ---}\index{Lax functor}\index{_aaa_Catl1@$\Cat_{l,1}$}
	where $\Cat_{l,1}(\cate X,\cate Y)$ denotes the category of \emph{normal lax functors} (see \ref{pseudocolax}) between bicategories $\cate X,\cate Y$. The correspondence is defined by a `generalised Grothendieck construction', in that  the correspondence $\Pi$ is defined sending $F\in\Cat_{l,1}(\C^\opp,\Dist)$ to the (strict) pullback
	\[
		\vcenter{\xymatrix{
				\Pi(F)\ar[r]\ar[d]\pb & \Dist_*\ar[d]^U \\
				C \ar[r]_F & \Dist
			}}
	\]
	where $U : \Dist_* \to \Dist$ is the forgetful functor from \emph{pointed profunctors} that forgets the basepoint.\footnote{A pointed profunctor $\proP : (\C, C) \pto (\D,D)$ of pointed categories is a functor $\proP : (\C^\opp\times \D, (C,D)) \to \Set$ with a specified element of $\proP(C,D)$. The functor $U$ forgets this specified element and keeps only the functor.}
\end{theorem}
\begin{notat}
	The category $\Pi(F)$ corresponding to a normal lax functor $F : \C \to \Dist$ is not any more \emph{fibered}; we call it a \emph{displayed} category.
\end{notat}
\begin{proof}
	As for the $\Pi$ correspondence, there is nothing to check apart its functoriality; this is easy: given a 2\hyp{}cell $\alpha :F \To G : \C \to \Dist$ between normal lax functors we can easily induce a morphism $\Pi(F) \to \Pi(G)$ of categories over $\C$ using the defining universal property of $\Pi$.

	We shall now define the functor $\Gamma$ and prove that it lands on the declared domain: we shall define
	\begin{itemize}
		\item A correspondence on objects, $p : \E \to \C$, that sends such $p$'s into normal lax functors $\C \to \Dist$;
		      \begin{enumtag}{gp}
			      \item an object $C$ goes to the category $p^\leftarrow(C) = \bsmat X\in\E \colon pX=C,\\ f : C \to C' \colon pf = \id_C\esmat $, \ie to the \emph{fiber} of $p$ over $C\in\C$;
			      \item a morphism $f : C \to C'$ goes to a profunctor $\Gamma_p(f) : p^\leftarrow(C)^\opp\times p^\leftarrow(C') \to \Set$, defined sending the pair $(X,Y) \in p^\leftarrow(C)^\opp\times p^\leftarrow(C')$ to the set of all $u : X \to Y$ such that $pu=f$.
		      \end{enumtag}
		\item A correspondence on morphisms over $\C$, that sends a functor $h : \E \to \E'$ of categories over $\C$, \ie such that $p'\circ h = p$ for projections $p,p'$, into a morphism of normal lax functors
		      \[\Gamma_h :  \xymatrix@C=1.5cm{\C \rtwocell^{\Gamma_p}_{\Gamma_{p'}}{}& \Dist.} \]
		      It is evident how a functor $h$ as above induces a morphism of this kind simply because it respects the fibers of $p$ and $p'$ over the same object.
	\end{itemize}
	Now, we shall show that $\Gamma_p : \C \to \Dist$ is indeed a normal lax functor: for the moment we have nothing but the bare definition.

	It is however quite easy to prove that the correspondence $\Gamma_p(f)$ if functorial and contravariant in the first component: if $\alpha : X \to X'$, there is an obvious function $\Gamma_p(X', Y) \to \Gamma_p(X,Y)$. Similarly, $\Gamma_p(f)$ is covariant in the second component.

	We shall now show that $\Gamma_p(\firstblank)$ is a normal lax functor: in order to do so we shall show the following.
	\begin{itemize}
		\item The functor is indeed normal: by definition $\Gamma_p(\id_C)$ is a functor $p^\leftarrow(C)^\opp\times p^\leftarrow(C)\to \Set$, that sends a pair $(X,Y)$ into the set of all $u : X \to Y$ such that $pu=\id_C$; but this is no less than the set of \emph{all} morphisms $X \to Y$ in the fiber $p^\leftarrow(C)$, so that $\Gamma_p(\id_C)$ is the hom functor of the fiber of $p$ over $C$, \ie the identity profunctor $p^\leftarrow(C)\pto p^\leftarrow(C)$.
		\item Given a pair of composable morphisms $C \xto{g} C' \xto{f} C''$, we shall find a 2\hyp{}cell filling the diagram
		      \[
			      \vcenter{\xymatrix@R=3mm{
			      & p^\leftarrow(C') \ar[dr]|-@{|}^{\Gamma_p(f)} & \\
			      p^\leftarrow(C) \ar[rr]|-@{|}_{\Gamma_p(fg)} \rrtwocell<\omit>{<-2>}\ar[ur]|-@{|}^{\Gamma_p(g)} && p^\leftarrow(C'')
			      }}
		      \]
		      The correspondences $\Gamma_p(f)$ and $\Gamma_p(g)$ are respectively defined by
		      \begin{gather*}
			      p^\leftarrow(C')^\opp \times p^\leftarrow(C'') \xto{\Gamma_p(f)} \Set \\
			      (X,Y) \overset{\Gamma_p(f)}\longmapsto \{u : X \to Y \mid pu = f\}; \\
			      p^\leftarrow(C)^\opp \times p^\leftarrow(C') \xto{\Gamma_p(g)} \Set \\
			      (Z,X) \overset{\Gamma_p(g)}\longmapsto \{v : Z\to X \mid pu = g\}.
		      \end{gather*}
		      Composition in the categories $p^\leftarrow(C), p^\leftarrow(C'), p^\leftarrow(C'')$ now forms a cowedge that induces a unique morphism
		      \[
			      \vcenter{\xymatrix{
					      \int^{X\in p^\leftarrow(C')} \Gamma_p(g)(Z,X)\times \Gamma_p(f)(X,Y) \ar[d]^c \\
					      \Gamma_p(fg)(Z,Y) = \{w : Z \to Y \mid pw = fg\}
				      }}
		      \]
		      by the universal property of the coend involved in the definition of $\Gamma_p(f)\diamond\Gamma_p(g)$. We leave to the reader the routine verification that this is indeed part of the laxity constraint of a lax functor $\Gamma_p : \C \to \Dist$.
	\end{itemize}
	This concludes the proof, up to some fine details that we leave to the reader (Additional exercise: try to find sufficient conditions fror the laxity cell $\Gamma_p(f)\diamond \Gamma_p(g)\To \Gamma_p(fg)$ to be invertible).
\end{proof}
\section{Addendum: Fourier theory}\label{sec:promono}\index{Category!promonoidal ---}
According to our §\ref{struct_of_Prof}, the bicategory $\Dist$ is monoidal with respect to the pseudo\hyp{}cartesian structure (similarly, every $\Dist(\V)$ inherits a symmetric monoidal structure from a symmetric monoidal structure on $\VCat$).

This means that we can consider internal monoids in $\Dist$: objects endowed with maps
\[
	\xymatrix{
	\M\times\M \ar[r]|-@{|}^-{\mathfrak{m}} & \M & 1 \ar[r]|-@{|}^-{\mathfrak{i}} & \M
	}
\]
in $\Dist$ that witness the fact that $\M$ is an internal (pseudo)monoid.

Such internal monoids take the name of \emph{promonoidal categories}.

Informally speaking, a promonoidal category is what we obtain if we replace every occurrence of the word \emph{functor} with the word \emph{profunctor} in the definition of monoidal category (of course, taking care of the coherence conditions imposed by the weak 2\hyp{}category structure of $\Dist$).

More precisely, we can give the following definition:
\begin{definition}[Promonoidal structure]\index{Promonoidal category}
	Let $\C$ be a category. A \emph{promonoidal structure} consists of a tuple
	\[\fkP = (\C, P, J,\alpha,\lambda,\rho)\]
	where
	\begin{enumtag}{pm}
		\item $\C$ is a category endowed with
		\item a bi\hyp{}profunctor $P : \C\times \C \pto \C$ (the monoidal \emph{multiplication}) and
		\item a profunctor $J : 1\pto \C$ (the monoidal unit), such that the following two diagrams
		\[\notag
			\xymatrix@C=2cm{
			\C \drtwocell<\omit>{\alpha} \times \C \times \C \ar[d]|-@{|}_{\hom\times P}\ar[r]|-@{|}^{P\times\hom}& \C \times\C \ar[d]|-@{|}^{P}\\
			\C \times \C \ar[r]|-@{|}_{P}& \C
			}\qquad
			\xymatrix@C=2cm{
			\C \drtwocell<\omit>{\rho} \ar@/_1pc/@{=}[dr]_\hom\ar[r]|-@{|}^{J\times \hom} & \C \times \C \ar[d]|-@{|}^{P} & \C \ar[l]|-@{|}_{\hom\times J} \ar@/^1pc/@{=}[dl]^\hom \\
			& \C \urtwocell<\omit>{\lambda}&
			}
		\]
		are filled by the indicated 2\hyp{}cells,
		\item 2-cells, respectively called the \emph{associator}
		\[\label{pms:uno}\alpha : P \diamond (P\times \hom) \cong P\diamond(\hom\times P)\]
		and the \emph{left} and \emph{right unitors}
		\[\label{pms:due}\lambda : P\diamond(\hom\times J)\cong \hom \quad \rho : P\diamond(J\times\hom)\cong\hom\]
	\end{enumtag}
\end{definition}
\begin{remark}
	Coend calculus allows to turn \eqref{pms:uno} and \eqref{pms:due} into diagrammatic relations:
	\index{Associator!promonoidal ---}
	\begin{itemize}
		\item The associator amounts to an isomorphism linking the two sets below (note that each component $\alpha^{ABC}_{D}$ has four arguments, three contravariant and one covariant, whereas $P$ has components $P^{AB}_C = P(A,B;C)$ as a functor $\C^\opp\times\C^\opp\times \C \to \Set$).
		      \begin{align*}
			      (P\diamond (\hom\times P))_{ABC;D} & =\int^{XY} P_D^{XY}H_A^X P_Y^{BC}                     \\
			                                         & \cong \int^YZ\Big(\int^X P_D^{XY}H_A^X \Big) P_Y^{BC} \\
			                                         & \cong \int^Z P_D^{AY}P_Y^{BC}                         \\
			      (P\diamond(P\times\hom))_{ABC;D}   & \cong \int^{XY} P_D^{XY} H_Y^C P_X^{AB}               \\
			                                         & \cong \int^Z P_X^{AB} P_D^{XC}.
		      \end{align*}
		\item The left unit axiom is equivalent to the isomorphism between the functor
		      \begin{align*}
			      (A,B) & \mapsto \int^{YZ} J_Z H^A_Y P^{YZ}_B              \\
			            & \cong \int^Z J_Z\Big( \int^Y H^A_Y P^{YZ}_B \Big) \\
			            & \cong \int^Z J_Z P^{AZ}_B
		      \end{align*}
		      and the hom functor $(A,B)\mapsto \C(A,B)$.
	\end{itemize}
\end{remark}
The most interesting feature of promonoidal structure in categories is that they correspond bijectively with monoidal structures on the category of functors $\Cat(\C, \Sets)$, framing the construction of \emph{Day convolution} give in \ref{day} in its maximal generality.
\index{_aaa_ast@$\ast$}
\begin{proposition}\label{promonoshit}\index{Convolution product}
	Let $\fkP = (P, J,\alpha, \rho,\lambda)$ be a promonoidal structure on the category $\C$; then we can define a $\fkP$-\emph{convolution} monoidal structure on the category $\Cat(\C, \Set)$, via
	\begin{align}
		[F\ast_{\fkP} G]C & = \int^{AB} P(A,B;C)\times FA\times GB \\
		J_{\fkP}          & = J
	\end{align}
	and this turns out to be a monoidal structure on $\Cat(\C, \Set)$. We denote the monoidal structure $(\Cat(\C, \Set), \ast_{\fkP}, J_{\fkP})$ shortly as $[\C, \Set]_{\fkP}$.
\end{proposition}
\begin{remark}
	The same definition, changing the cartesian structure with the monoidal structure of $\V$, yields a notion of $\fkP$\hyp{}convolution on $\VCat(\C,\V)$ for a $\V$\hyp{}category $\C$.
\end{remark}
\begin{definition}
	A functor $\Phi : [\A, \Sets]_{\fkP} \to [\B, \Sets]_{\fkQ}$ is said to \emph{preserve the convolution product} if the obvious isomorphisms hold in $[\B, \Sets]_{\fkQ}$:
	\begin{itemize}
		\item $\Phi(F\ast_{\fkP} G) \cong \Phi(F)\ast_{\fkQ}\Phi(G)$;
		\item $\Phi(J_{\fkP}) = J_{\fkQ}$;
	\end{itemize}
	in other words $\Phi$ is a strong monoidal functor with respect to convolution product. When $\Phi$ is a colimit\hyp{}preserving functor, this condition is equivalent to the request that $\Phi$ defines a \emph{multiplicative kernel} between $\A,\B$ regarded as objects of $\Dist$.
\end{definition}
\begin{remark}
	It is observed in \cite{imkelly} that for a monoidal $\A$ the category of presheaves $[\A^\opp,\V]$ endowed with the convolution monoidal structure is the \emph{free monoidal cocompletion} of $\A$, having in $\BF{Mon}$ (monoidal categories, monoidal functors and monoidal natural transformations) the same universal property that $[\A^\opp,\V]$ has in $\Cat$.
\end{remark}
There is a bijection between the promonoidal structures on $\C$, and monoidal structure on $\Cat(\C,\Set)$; this is the content of Exercise \ref{es8_2}.
\subsection{Fourier transforms via coends}
\index{Multiplicative kernel|see{Fourier transform}}
\index{Fourier transform}
For the rest of the section, $\V$ is assumed to be a complete and cocomplete *\hyp{}autonomous category.
\begin{definition}\label{mulker}
	Let $\A ,\C$ be two promonoidal categories (thus implicitly regarded as objects of $\Dist$) with promonoidal structures $\fkP$ and $\fkQ$ respectively; a \emph{multiplicative kernel} from $\A$ to $\C$ consists of a profunctor $K : \A\pto \C$ endowed with two natural isomorphisms
	\begin{align}
		\int^{YZ}K^A_Y K^B_Z P^{YZ}_X & \cong \int^C K^C_X P^{AB}_C\label{k_1} \\
		\int^C K^C_X J_C              & \cong J_X\label{k_2}
	\end{align}
	These isomorphisms say that $K$ mimics the behaviour of the hom functor (in fact, the hom functor $\hom_\A$ is the \emph{identity} multiplicative kernel $\A\pto\A$: the isomorphisms above follow from \ref{ninjayo}).

	We define a \emph{multiplicative} natural transformation $\alpha : K\to H$ between two kernels as a 2\hyp{}cell in $\Dist$ commuting with the structural isomorphisms given in \ref{mulker}. This, together with the fact that multiplicative kernels compose, yields a category of kernels $\ker(\A,\C)$.
\end{definition}
\begin{definition}
	Let $K : \A\pto\C$ be a multiplicative kernel between promonoidal categories; we define the $K$-\emph{Fourier transform} $f\mapsto \hat K(f) : \C\to \Sets$, obtained as the image of $f : \A\to\Sets$ under the left Kan extension $\Lan_\yon K : [\A, \Sets]\to [\C, \Sets]$. Explicitly, this is the coend
	\[
		\F_K(f) : X\mapsto \int^A K(A,X)\otimes fA.
	\]
\end{definition}
We can also define the dual Fourier transform:
\[
	\F^\lor(g) : Y\mapsto \int_A [K(A,X), gA]
\]
and find the relation $\F^\lor_K(g) \cong \F_K(g^*)^*$.

The following results are easily proved using standard co\fshyp{}end calculus:
\begin{proposition}
	Let $K : \A \pto \cate X$ be a multiplicative kernel, and let $\A$ be a promonoidal category; then
	\begin{enumtag}{mk}
		\item $\F_K$ preserves the upper $\mathfrak P_\A$\hyp{}convolution of presheaves $f,g$, defined as
		\[f\mathrel{\overline{\ast}} g = \int^{AA'} fA\otimes gA'\otimes P(A,A',\firstblank);\] dually,
		\item $\F_K^\lor$ preserves the lower $\mathfrak P_\A$\hyp{}convolution of presheaves $f,g$, defined as
		\[f\mathrel{\underline{\ast}} g = \int_{AA'} \Big(fA^*\otimes (gA')^*\otimes P(A,A',\firstblank)\Big)^*;\]
	\end{enumtag}
\end{proposition}
Observe that (as stated in \cite{day2011monoidal}), both the upper and lower convolution product yield associative and unital monoidal structures on the functor category $\VCat(\A,\V)$; the upper product preserves $\V$\hyp{}colimits in each variable, while the lower product preserves $\V$\hyp{}limits
in each variable.

\index{Convolution!--- product}
\index{Convolution!upper and lower ---}
The  lower and upper convolution transform into each other under the equivalence of $\V$\hyp{}categories
\[
	\VCat(\A,\V)^\opp\cong \VCat(\A^\opp,\V^\opp);
\]
this means that under the above equivalence $(f \mathrel{\overline{\ast}} g)^* \cong f^* \mathrel{\underline{\ast}} g^*$.
\begin{theorem}
	Let $\V$ be a *\hyp{}autonomous monoidal base; then we can define the pairing $\VCat(\A,\V)\times \VCat(\A,\V) \to \V$ as the twisted form of functor tensor product (as defined in \ref{tenso_pro_offunc})
	\[
		\langle f,g\rangle = \int^A fA^*\otimes gA
	\]
	\index{Parseval formula|see{Fourier transform}}
	\index{Fourier transform}
	If $K$ is a kernel such that the Fourier transform $\Lan_\yon K$ is fully faithful, we have an analogue of \emph{Parseval formula}:
	\[
		\langle f,g\rangle \cong \langle \F_K(f),\F_K(g)\rangle.
	\]
\end{theorem}
\index{Combinatorial species}
Fourier theory is linked to the theory of Joyal's \emph{combinatorial species}
(see \ref{analuo}): let $E : \A\pto \cate X$ be a multiplicative kernel; if we write
$E(A,X):= X^A$ (without any reference to a tensor operation between $A$ and
$X$), the $E$-Fourier transform can be expressed as an $\A$\hyp{}indexed formal power
series as follows:
\begin{align}
	\F_E(f) & = \int^A fA\otimes X^A \notag        \\
	        & \cong \sum_{A\in\A} fA\otimes_\A X^A
\end{align}
(it is understood that $f : \A \to \V$ is a fixed combinatorial species.) This can be made precise as follows: let $\A$ be the (free $\V$\hyp{}category on the) permutation category of
\ref{prelim_notata}; then $E(n,X) := X^{\otimes n}=X\otimes\dots\otimes X$ is a
multiplicative kernel and an $E$\hyp{}analytic functor results as the left Kan
extension
\begin{align}
	FX & =\int^n f(n)\otimes X^{\otimes n} \notag                       \\
	   & \cong \sum_{n\in\P} f(n)\otimes_{\text{Sym}(n)} X^{\otimes n}.
\end{align}
Given two combinatorial species $f,g : \A \to \V$ and the associated analytic
functors $F,G$, the convolution $F \ast G$ is again an analytic functor, and its
generating combinatorial species is the upper convolution product of $f,g$ (with
respect to the implicit promonoidal structure of $\A$):
\begin{align}
	F \ast G(X) & \cong \int^{AB} FA\otimes GB\otimes p(A,B;X) \notag                              \\
	            & \cong \int^{ABUV} fU\otimes E(U,A)\otimes gV\otimes E(V,B)\otimes p(A,B;X)\notag \\
	            & \cong \int^{UVC} fU\otimes gV\otimes p(U,V;C)\otimes E(C,X) \tag{k}              \\
	            & \cong \int^C (f\mathrel{\overline{\ast}} g)(C)\otimes E(C,X).
\end{align}
(note that in (k) we used the fact that $E$ is a multiplicative kernel.)
\section{Addendum: Tambara theory}
\index{Tambara module}
\index{Profunctor}
\begin{definition}
	Let $\C$ be a monoidal category with monoidal unit $I$. A (\emph{left}) \emph{Tambara module} on $\C$ consists of:
	\begin{itemize}
		\item a profunctor
		      $P : \C^\opp \times \C \rightarrow \Set$;
		\item a family of functions
		      $\tau_{A}(X, Y) : P(X, Y) \longrightarrow P(A \otimes X, A \otimes Y)$
		      natural in $X, Y$ and a wedge in $A$, satisfying the two equations:
		      \begin{gather}
			      \vcenter{\xymatrix{
			      P(X, Y) \ar[rr]^{\tau_I(X,Y)}\ar@{=}[dr]&& P(I\otimes X, I\otimes Y) \ar[dl]^{P(l_X^{-1},l_Y)}\\
			      & P(X,Y)
			      }}
			      \notag\\
			      \vcenter{\xymatrix{
			      P(X,Y) \ar[rr]^{\tau_{A'}(X, Y)}
			      \ar[dr]_{\tau_{A \otimes A'}(X, Y)\quad}
			      && P(A'\otimes X, A'\otimes Y) \ar[dl]^{\quad\tau_{A}(A' \otimes X, A' \otimes Y)} \\
			      &  P(A \otimes A' \otimes X, A \otimes A' \otimes Y)&
			      }}
		      \end{gather}
	\end{itemize}
	The notion of \emph{right} Tambara module is given in a similar fashion, using maps $\nu_A(X,Y) : P(X, Y) \longrightarrow P(X \otimes A, Y \otimes A)$ satisfying the relations
	\begin{enumtag}{lt}
		\item $P(r_X, r_X^{-1})\circ \nu_A(X,Y) = \id_{P(X,Y)}$;
		\item $\nu_A(X\otimes A',Y\otimes A')\circ \nu_{A'}(X,Y) = \nu_{A'\otimes A}(X,Y)$.
	\end{enumtag}
	The definition can be given for profunctors enriched over any other monoidal base different from $\Set$; for example, in the original work by Tambara, \cite{tambara2006distributors}, categories are assumed to be enriched over vector spaces. In the paper \cite{pastro2008doubles} the enrichment base is completely arbitrary (\ie it is just a symmetric monoidal closed category).
\end{definition}
\begin{definition}
	Define the category $\Tamb(\C)$ whose:
	\begin{itemize}
		\item objects are Tambara modules $(P, \tau)$ consisting of a
		      profunctor $P : \C^\opp \times \C \rightarrow \Set$
		      and Tambara structures $\tau_{A}(X, Y)$.
		\item morphisms $(P, \tau) \rightarrow (Q, \sigma)$
		      are natural transformations
		      $\gamma : P \Rightarrow Q$ such that for all $A, X, Y$ the following
		      diagram commutes:
		      \[ \vcenter{\xymatrix@C=1.5cm{
			      P(X, Y)
			      \ar[r]^-{\tau_{A}(X, Y)}
			      \ar[d]_{\gamma_{(X, Y)}}
			      & P(A \otimes X, A \otimes Y)
			      \ar[d]^{\gamma_{(A \otimes X, A \otimes Y)}}
			      \\
			      Q(X, Y)
			      \ar[r]_-{\sigma_{A}(X, Y)}
			      & Q(A \otimes X, A \otimes Y)
			      }}
		      \]
	\end{itemize}
	There is a functor to the category of endo\hyp{}profunctors on $\C$,
	\[
		\iota : \Tamb(\C) \longrightarrow \Dist(\C,\C)
	\]
	which forgets the Tambara structure.
\end{definition}
The codomain of $\iota$ is monoidal with respect to composition of 1\hyp{}cells, as every hom\hyp{}category of endomorphisms: it turns out that Tambara modules can be composed, and that $\iota$ is strong monoidal with respect to this monoidal structure.
\begin{remark}
	The category $\Tamb(\C)$ has a monoidal structure whose:
	\begin{itemize}
		\item \emph{unit} is the hom\hyp{}functor
		      $\hom_\C : \C^\opp \times \C \rightarrow \Set$
		      which has a canonically associated Tambara structure:
		      \[
			      \C(X, Y) \longrightarrow \C(A \otimes X, A \otimes Y)
		      \]
		\item The profunctor composition of
		      $(P, \tau)$ and $(Q, \sigma)$ given by the coend
		      \[
			      (P \diamond Q)(X, Y) = \int^{Z} P(X, Z) \times Q(Z, Y)
		      \]
		      has a Tambara structure
		      $(P \diamond Q)(X, Y) \rightarrow (P \diamond Q)(A \otimes X, A \otimes Y)$
		      induced by the maps
		      \[
			      \tau_{A} \times \sigma_{A} : P(X, Z) \times Q(Z, Y)
			      \rightarrow P(A \otimes X, A \otimes Z) \times Q(A \otimes Z, A \otimes Y)
		      \]
		      using the universal property of the coend.
	\end{itemize}
	This makes the functor $\iota : \Tamb(\C) \rightarrow \Dist(\C,\C)$
	strong monoidal.
\end{remark}
\begin{proposition}
	The forgetful functor $\iota : \Tamb(\C) \rightarrow \Dist(\C,\C)$
	forms part of an adjoint triple:
	\[
		\vcenter{\xymatrix{
				\Tamb(\C) \ar[r]|\iota & \ar@<6pt>[l]^\varphi  \ar@<-6pt>[l]_\theta \Dist(\C,\C)
			}}
	\]
\end{proposition}
\begin{itemize}
	\item The left adjoint $\phi : \Dist(\C,\C) \rightarrow \Tamb(\C)$ constructs
	      the \emph{free Tambara module} from a profunctor. This is given by the formula
	      \[
		      \phi_{P}(X, Y)
		      = \int^{C, U, V} \C(X, C \otimes U) \times \C(C \otimes V, Y) \times P(U, V)
	      \]
	      with Tambara module structure given by
	      \[
		      \vcenter{\xymatrix{
				      \C(X, C \otimes U) \times \C(C \otimes V, Y) \times P(U, V)
				      \ar[d] \\
				      \C(A \otimes X, A \otimes C \otimes U)
				      \times \C(A \otimes C \otimes V, A \otimes Y) \times P(U, V)
			      }}
	      \]
	      together with the coprojection $q_{A \otimes C}$,
	      using the universal property of the coend.
	\item The right adjoint $\theta : \Dist(\C,\C) \rightarrow \Tamb(\C)$
	      constructs the \emph{cofree Tambara module} from a profunctor. This is given by the formula
	      \[
		      \theta_{P}(X, Y) = \int_{C} P(C \otimes X, C \otimes Y)
	      \]
	      with Tambara module structure given by
	      $\theta_{P}(X, Y) \rightarrow \theta_{P}(A \otimes X, A \otimes Y)$ is induced
	      by the projection functions,
	      \begin{gather}
		      p_{C \otimes A} : \int_{C} P(C \otimes X, C \otimes Y)
		      \rightarrow P(C \otimes A \otimes X, C \otimes A \otimes Y)\notag\\
		      p_{C} : \int_{C} P(C \otimes A \otimes X, C \otimes A \otimes Y)
		      \rightarrow P(C \otimes A \otimes X, C \otimes A \otimes Y)
	      \end{gather}
	      using the universal property of the end.
\end{itemize}
The proof of the following proposition (a `recognition principle' for Tambara modules) goes by inspection
using the definition of coalgebra: such a map is determined by
\begin{itemize}
	\item An object $P$ of $\Dist(\C,\C)$ given by
	      $P : \C^\opp \times \C \rightarrow \Set$;
	\item A structure map given by a natural transformation
	      $\tau : P \Rightarrow \theta_{P}$ in $\Dist(\C,\C)(P,\theta_P)$ whose components
	      $\tau(X, Y) : P(X, Y) \rightarrow \theta_{P}(X, Y)$ are given by:
	      \[
		      \tau(X, Y) : P(X, Y) \longrightarrow \int_{A} P(A \otimes X, A \otimes Y)
	      \]
	\item By the universal property of the end at codomain, the structure map is
	      determined by a wedge in $A$,
	      \[
		      \tau_{A}(X, Y) :
		      P(X, Y) \longrightarrow P(A \otimes X, A \otimes Y)
	      \]
	      that is moreover natural in $X,Y$.
\end{itemize}
\begin{proposition}
	The adjunction $\iota \dashv \theta$ yields a comonad
	\[
		\Theta : \Dist(\C,\C) \rightarrow \Dist(\C,\C)
	\]
	whose category of coalgebras is
	isomorphic to $\Tamb(\C)$. Dually, the adjunction $\phi \dashv \iota$ yields a monad
	$\Phi : \Dist(\C,\C) \rightarrow \Dist(\C,\C)$ whose category of algebras is
	isomorphic to $\Tamb(\C)$. Moreover, there is an
	adjunction
	\[
		\xymatrix{
			\Dist(\C,\C) \ar@<4pt>[r]^\Phi\ar@{}[r]|\perp & \ar@<4pt>[l]^\Theta\Dist(\C,\C)
		}
	\]
	between the resulting monad $\Phi = \iota\circ\varphi$ and comonad $\Theta = \iota\circ\theta$ on $[\C^\opp \times \C, \Set]$.
\end{proposition}
Profunctors $\A \pto \B$ can equivalently be described as left adjoints $\Cat(\B^\opp,\Set) \to \Cat(\A^\opp,\Set)$; thus we obtain that
\begin{corollary}
	The left adjoint
	\[
		\Phi : \Cat(\C^\opp \times \C, \Set)
		\to \Cat(\C^\opp \times \C, \Set)
	\]
	is equivalent to the \emph{endo\hyp{}profunctor} $\check\Phi : \C^\opp\times\C\pto\C^\opp\times\C$ whose action on objects is given by the coend:
	\[
		\check{\Phi}(X, Y, U, V)
		= \int^{C} \C(X, C \otimes U) \times \C(C \otimes V, Y)
	\]
\end{corollary}
This endo\hyp{}profunctor
$\check{\Phi} : \C^\opp \times \C \pto \C^\opp \times \C$
actually underlies a \emph{promonad} (see Exercise \ref{ex:promonad}) in the bicategory $\Dist$.
From the formal theory of monads \cite{Street1972} it is known that the bicategory $\Dist$ admits the Kleisli
construction for promonads, so we can ask what is the Kleisli category of $\check\Phi$: such category
is called the (\emph{left}) \emph{double} of the monoidal category $\C$ and it is denoted $\BF{Db}(\C)$; it has the same objects as
$\C^\opp \times \C$, and hom\hyp{}sets defined by the coend
\[
	\BF{Db}(\C)\big( (X, Y), (U, V) \big)
	= \int^{C} \C(X, C \otimes U) \times \C(C \otimes V, Y)
\]
This formula provides the foundation for all of \emph{profunctor optics} (see \cite{barto,noi,pickering2017profunctor}).
\cite{pastro2008doubles} proves that there is an equivalence of categories:
\[
	\Tamb(\C) \simeq \Cat(\BF{Db}(\C), \Set).
\]
\begin{exercises}
\item Describe the bicategory of profunctors between monoids, regarded as one\hyp{}object categories; describe the bicategory of profunctors between posets regarded as thin categories; similarly, the bicategory $\Mod$ of modules.
\item \label{riftuzzo} Dualise \ref{rans}: given $\proH : \D\pto \A$ and $\proL : \cate E\pto \A$ we can define
\[\notag
	\proL\rtimes \proH : \E \pto \D : (D,E) \mapsto \int_A [\proH(D,A), \proL(E,A)].
\]
Show that this second operation is a right \emph{Kan lifting} (we spell out explicitly the definition that can be evinced from \ref{specchietto_dei_lift}): given 1\hyp{}cells $p : B \to C$, $f : A \to C$ in a 2\hyp{}category $\sfK$, a \emph{right Kan lift} of $f$ through $p$, denoted $\Rift_p f$, is a 1\hyp{}cell $\Rift_p(f) : A \to B$ equipped with a 2\hyp{}cell
\[\notag
	\varepsilon: p \circ \Rift_p(f) \Rightarrow f
\]
satisfying the following universal property: given any pair $(g : A \to B, \alpha : p \circ g \Rightarrow f)$, there exists a unique 2\hyp{}cell
\[\notag
	\zeta: g \Rightarrow \Rift_p(f)
\]
such that the following diagram of 2\hyp{}cells commutes for a unique $\zeta : g\Rightarrow \Rift_p(f)$
\[\notag
	\xymatrix{
	& B \ar[dr]^p & \\
	A \ar@/^1.5pc/[ur]^g \ar[rr]_f && C
	\ar@{=>}(13,-4);(13,-11)^\alpha
	} \qquad \text{\large =}\qquad
	\xymatrix{
	& B \ar[dr]^p & \\
	A \ar@/^2pc/[ur]^g \ar[ur]|{\Rift_p(f)} \ar[rr]_f && C
	\ar@{=>}(13,-4);(13,-11)^\epsilon
	\ar@{:>}(3,-1);(7,-5)^\zeta
	}
\]
\ie there is a unique factorisation $\varepsilon \circ (p \ast \zeta) = \alpha$.
\item Show that the structure on $\Dist$ given by $\diamond$ is \emph{biclosed} (\ie, $\diamond$ is a bifunctor $\Dist(\A,\B)\times \Dist(\B,\C) \to \Dist(\A,\C)$ and each $\proP\diamond\firstblank$, as well as each $\firstblank\diamond \proQ$ have right adjoints).
\item \label{collage} The \emph{collage} of two categories $\A,\B$ along a profunctor $\proP : \A \pto \B$ is defined as the category $\A \uplus_\proP \B$ with the same objects as $\A \amalg \B$ and morphisms given by the rule
\[\notag
	\A\uplus_\proP \B(X,Y) =
	\begin{cases}
		\A(X,Y)    & \text{ if } X,Y\in\A        \\
		\B(X,Y)    & \text{ if } X,Y\in \B       \\
		\proP(X,Y) & \text{ if } X\in \A, Y\in\B
	\end{cases}
\]
and empty in every other case. Show that $\A\uplus_\proP \B$ has the universal property of the category of elements of $\proP$, regarded as a presheaf. Find an isomorphism between the category $\A \uplus_\proP \B$ so defined and the coend
\[\notag
	\int^{(A,B)} (\A^\opp\times \B) /(A,B)\otimes \proP(A,B)
\]
of \ref{elts-as-coend} ($\otimes$ is the $\Set$\hyp{}tensor of $\Cat$).
\item Show that the composition laws $\proP(A,B)\times \B(B,B')\to \proP(A,B')$, $\A(A,A')\times \proP(A',B)\to \proP(A, B)$ of arrows in the collage $\A\uplus_\proP \B$ are universal cowedges of a coend.
\item \label{ex5:cocomma-as-colage}The \emph{cocomma object} $\bsmat F \\ G \esmat$ of two functors $\cX \xot{F} \A \xto{G} \cY$ is defined to be the pushout of
\[\notag
	\xymatrix{
		\A \amalg \A \ar[d]\ar[r]& \A \times [1] \\
		\cX \amalg \cY &
	}
\]
in $\Cat$, where the horizontal arrow is the `cylinder' embedding when $\A \amalg \A$ is identified with $\A^{\{0,1\}} = \Cat(\{0,1\}, \A)$. Show that $\bsmat F \\ G \esmat$ is the collage of $\cX$ and $\cY$ along the profunctor $\proP^G \ccirc \proP_F : \cX \pto \cY$.
\item Given profunctors $\A \overset{\proP}\pto \B \overset{\proQ}\pto \C$ consider the categories $\A\uplus_\proP\B$ and $\B\uplus_\proQ\C$. Describe the pushout
\[\notag
	\xymatrix{
		\B \ar[r]\ar[d]& \A\uplus_\proP\B \ar[d]\\
		\B\uplus_\proQ\C \ar[r] & \ar@{}[ul]|(.2)\ulcorner \cate H
	}
\]
in $\Cat$. Is there a relation between $\cate H$ and the collage $\A \uplus \C$ along $\proQ\ccirc \proP$?
\item Isbell duality (\ref{isbella-duella}) can be regarded as an adjunction between the categories
\[\notag
	\Dist(\uno, \C)^\opp \leftrightarrows \Dist(\C, \uno).
\]
where $\uno$ is the terminal category. Is it possible to extend this result to an adjunction $\Dist(\D, \C)^\opp \leftrightarrows \Dist(\C, \D)$? (hint: yes; use right Kan extensions in $\Dist$).
\item \label{forexactness} Show that there is a canonical isomorphism
\[\notag
	\Cat(\cX, \cY)(FG, HK) \cong \Dist(\cW, \cZ)(\proP_G\ccirc \proP^K, \proP^F\ccirc \proP_H)
\]
for each square
\[\notag
	\xymatrix{
	\cX \ar[r]^G \ar[d]_K \ar@{}[dr]|{\Swarrow\alpha}& \cZ\ar[d]^F \\
	\cW \ar[r]_H & \cY
	}
\]
filled by a 2\hyp{}cell $\alpha$. This map sends $\alpha$ to a 2\hyp{}cell $\alpha^\sharp$ in $\Dist$. Does this equivalence restrict to strongly commutative squares, giving strongly commutative squares in $\Dist$?
\item \label{ex:promonad} \index{Monad!--- in $\Dist$}\index{Promonad|see{Monad}} A \emph{promonad} is a monad $T : \A \pto \A$ over an object o $\Dist$; this means that there are maps $T \diamond T \To T$ and $\hom \To T$ fitting in diagrams similar to \ref{def:monad}.
\begin{enumtag}{pm}
	\item Regard a set $A$ as a discrete category; show that every promonad $T : A \pto A$ determines and is determined by a category structure for $A$, \ie that a promonad on a discrete category amounts exactly as a category having $A$ as set of objects.
	\item \awful Show a similar result for nondiscrete categories: more precisely, observe that a promonad on $\A$ in $\Dist$ corresponds to a monad on $\Cat(\A^\opp,\Set)$ whose underlying endofunctor preserves colimits. Show that every category $\B$ that admits an identity\hyp{}on\hyp{}objects functor $F : \A \to \B$ induces such a monad; show that every colimit preserving monad between presheaf categories induces an identity\hyp{}on\hyp{}object functor $\A \to \B$.
\end{enumtag}
\item \label{ex5:promonads_go_to_cmc}Show that the equivalence
\[\notag
	\text{LAdj}(\Cat(\cX^\opp,\Set), \Cat(\cX^\opp,\Set)) \cong \Dist(\cX,\cX)\]
of \ref{yext_are_good} is monoidal, i.e. colimit\hyp{}preserving monads $T$ on $\Cat(\cX^\opp,\Set)$ go to promonads $\tau_T : \cX \pto \cX$.
\item \label{ex:traced} Gather from Chapter 1 the results you need to prove that the trace operator $\tr_{XY}^U$ in \ref{def:traced} really define a traced monoidal structure:
\begin{enumtag}{tm}
	\item naturality in $X$ and $Y$: for every $f:X\otimes U\to Y\otimes U$ and $g:X'\to X$,
	\[\notag
		\tr^U_{X',Y}(f \circ (g\otimes \id_U)) = \tr^U_{X,Y}(f) \circ g,
	\] and for every $f:X\otimes U\to Y\otimes U$ and $g:Y\to Y'$,
	\[\notag
		\tr^U_{X,Y'}((g\otimes \id_U) \circ f) = g \circ \tr^U_{X,Y}(f)
	\]
	\item dinaturality in $U$: for every $f:X\otimes U\to Y\otimes U'$ and $g:U'\to U$
	\[\notag
		\tr^U_{X,Y}((\id_Y\otimes g) \circ f)=\tr^{U'}_{X,Y}(f \circ (\id_X\otimes g))
	\]
	\item two vanishing conditions: for every $f:X \otimes I \to Y \otimes I$, (with $\rho_X \colon X\otimes I\cong X$ being the right unitor),
	\[\notag \tr^I_{X,Y}(f)=\rho_Y \circ f \circ \rho_X^{-1},\]
	and for every $f:X\otimes U\otimes V\to Y\otimes U\otimes V$
	\[\notag
		\tr^U_{X,Y}(\tr^V_{X\otimes U,Y\otimes U}(f)) = \tr^{U\otimes V}_{X,Y}(f)
	\]
	\item superposing: for every $f:X\otimes U\to Y\otimes U$ and $g:W\to Z$,
	\[\notag
		g\otimes \tr^U_{X,Y}(f)=\tr^U_{W\otimes X,Z\otimes Y}(g\otimes f)
	\]
	\item yanking: $\tr^X_{X,X}(\gamma_{X,X})=\id_X$
	(where $\gamma$ is the symmetry of the monoidal category).
\end{enumtag}
\item \label{es8_1} Prove equations \ref{promonoshit} using associativity and unitality for $\fkP$.
\item \label{es8_2} \index{Category!promonoidal ---}\index{Promonoidal category}\index{Convolution product}
Let $* : \Cat(\C,\Set)\times\Cat(\C,\Set) \to \Cat(\C,\Set)$ be a monoidal structure with monoidal unit $u : \C \to \Set$; show that the assignment
\[
	P(A,B;C) := (\coyon(A) * \coyon(B))(C)\qquad JA := uA
\]
\index{_aaa_yoncontra@$\yon$}
\index{_aaa_yoncov@$\coyon$}
is a promonoidal structure on $\C$, regarded as an object of $\Dist$. An elegant result of Day shows that this sets up a bijection between the ways in which $\Cat(\C,\Set)$ is a (pseudo)monoid in $\Cat$, and the ways in which $\C$ is a (pseudo)monoid in $\Dist$: prove it.
\item \label{es8_3} \index{Category!promonoidal ---}\index{Promonoidal category}\index{Convolution product}Outline the promonoidal structure $\fkP$ giving the Day convolution described in \ref{day}. If $\C$ is any small category, we define $P(A,B;C) = \C(A,C)\times \C(B,C)$ and $J$ to be the terminal functor $\C\to\Sets$. Outline the convolution product on $\Cat(\C, \Set)$, called the \emph{Cauchy convolution}, obtained from this promonoidal structure.
\item \label{es8_4} Is the composition of two kernels (see \ref{mulker}) again a kernel? Define the category of multiplicative kernels $\ker(\A,\C)\subset \Dist(\A,\C)$.
\item \label{es8_5} Show that a profunctor $K : \A\pto \C$ is a multiplicative kernel if and only if the cocontinuous functor $\Lan_\yon K=\hat K : [\A, \Sets]\to [\C, \Sets]$ corresponding to $\bar K : \A \to \VCat(\C,\V)$ under the construction in \ref{alternative} is monoidal with respect to the convolution monoidal structure on both $[\A, \Sets]_{\fkP}$ and $ [\C, \Sets]_{\fkQ}$.

Describe the isomorphisms $k_1, k_2$ when $\fkP$ is Day convolution.
\item \label{es8_6} Show that a functor $F : (\A,\otimes_{\A}, I) \to (\C,\otimes_\C, J)$ between monoidal categories is strong monoidal if and only if $\proP^F=\hom(F,1)$ is a multiplicative kernel.

Dually, show that for $\A,\C$ promonoidal, $F : \C\to \A$ preserves convolution on $[\A, \Sets]_{\fkP}, [\C, \Sets]_{\fkQ}$ precisely if $\proP_F=\hom(1, F)$ is a multiplicative kernel.
\item Show the following properties of the $K$-Fourier transform:\index{Parseval identity}
\begin{itemize}
	\item There is the canonical isomorphism
	      \[\notag
		      \hat K(f) \cong \int^A K(A,\firstblank)\times f(A)
	      \]
	\item $\hat K$ preserves the convolution monoidal structure (this is the \emph{Parseval identity} for the Fourier transform);
	\item $\hat K$ has a right adjoint defined by
	      \[\notag
		      \check{K}(g) \cong \int_x [K(\firstblank, X), g(X)].
	      \]
\end{itemize}
\item \awful Prove the following statements:
\begin{itemize}
	\item There is a monad $\tilde S$ on the category of profunctors, such that the following square is commutative:
	      \[ \notag \xymatrix{
			      \tCat\ar[r]^{ S}\ar[d] & \tCat \ar[d]\\
			      \Dist \ar[r]_{\tilde S} & \Dist
		      } \]
	      where $S$ sends a category $\A$ to the \emph{free monoidal category} on $\A$ and $\tCat \to \Dist$ is the canonical embedding of \ref{embare_i_profi} (on which cells it is contravariant?).
	\item Prove that A $\tilde S$-algebra is the same thing as a monad in the bicategory of $\tilde S$-profunctors
	\item Show that the following conditions are equivalent, for a profunctor $\proP : \A\pto \B$ between promonoidal categories
	      $ (\A,\fkP, J_\A),(\B,\fkQ, J_\B)$:
	      \begin{enumtag}{pa}
		      \item \label{pa:uno} $\proP$ is a pseudo-$\tilde S$-algebra morphism;
		      \item \label{pa:due} The cocontinuous left adjoint associated to $\proP$, $\hat \proP : [\B^\opp,\Set] \to [\A^\opp,\Set]$ is strong monoidal with respect to the convolution monoidal product on presheaf categories;
		      \item \label{pa:tre} $\proP$ is endowed with arrows
		      \[
			      \notag
			      \xymatrix@R=1mm{
			      \proP^{A_1}_{B_1}\times \proP^{A_2}_{B_2} \times \fkQ^{B_1B_2}_B \ar[r]^-{\gamma} & \int^A \fkP^{A_1A_2}_A \times \proP^A_B \\
			      (\proP^\opp)^{B_1}_{A_1} \times (\proP^\opp)^{B_2}_{A_2} \times \fkP^{A_1A_2}_A \ar[r]^-{\sigma} & \int^B (\proP^\opp)^B_A \times \fkQ^{B_1B_2}_B  \\
			      p^A_B \times {J_\A}_A \ar[r]^-{\delta} & {J_\B}_B
			      }
		      \]
		      or more precisely, $\sigma_{B_1B_2 [A_1A_2];A}, \sigma_{B_1B_2 [A_1A_2];A}, \delta_{[A]B}$, exhibiting universal cowedges in the bracketed variables (i.e., the maps induced on coends are isomorphisms), and natural in the others.
	      \end{enumtag}
	\item Assume the promonoidal structures $\fkP, \fkQ$ on $\A,\B$ are representable; then, the conditions above are in turn equivalent to the following: both mates $p^\lhd : \A \to \psh{\B}$ che $p^\rhd : \B \to [\A,\Set]^\opp$ are strong monoidal with respect to convolution on their codomains.
\end{itemize}
\end{exercises}

\chapter{Operads}\label{sec:opd}
\begin{abstract}
	We introduce the theory of \emph{operads} employing co\fshyp{}end calculus; the material is entirely classical and draws from \cite{Kelly2005a} and equally classical sources. A (symmetric) operad is a collection $O(n)$ of objects of a monoidal category whose objects are natural numbers, and endowed with maps
	\[\notag\textstyle
		O(n_1)\otimes \dots\otimes O(n_k) \otimes O(k) \to O\big(\sum n_i\big)
	\]
	satisfying suitable axioms of associativity and compatibility with a natural action of the symmetric group on each component; each $O(n)$ models a set of generalised $n$\hyp{}ary operations serving to describe in a neat and intrinsic way the equipment of a `set' with `structure'. The notion of operad lies at the very core of modern approaches to universal and categorical algebra.

	Operads are monoid objects in the presheaf category of representations of the groupoid of natural numbers; drawing from \cite{curiennone} we provide a standard characterisation of operads as the coKleisli category of a comonad on the bicategory $\Dist$; more precisely, an operad is an object in the coKleisli category of a comonad generated by the presheaf construction $\bsP$ and by the `free symmetric monoidal category' functor $S$; this allows for plenty of generalisations to other kinds of operads, by changing the role of $S$ in the same formal argument.
\end{abstract}
\epigraph{
	The sixth [is] the method of returning the letters to their prime\hyp{}material state and giving them form in accordance with the power of wisdom that confers form. [\dots\unkern] Regarding this method, it is stated in the \emph{Sefer Yet\d{z}irah}: `Twenty\hyp{}two cardinal letters; He engraved them and hewed them and weighed them and permuted and combined them and formed by their means the souls of all formed beings.'
}
{A. Abulafia --- \emph{O\d{z}ar Eden Ganuz}, quoted in \cite{book802037}}
\section{Introduction}
Operads are mathematical structures of manifold nature: they appear in algebra, topology, algebraic geometry, logic, and in each of these settings they model the notion of `set endowed with operations' providing an extremely powerful conceptual tool to categorise the old discipline known as universal algebra.

Operads were introduced by P\@. May in his \cite{may1972geometry} in order to solve a purely algebraic\hyp{}topological problem: topologists are often interested to classify spaces $Y$ which are homotopy equivalent to a loop space $\Omega X$; every such space carries the structure of a group up to homotopy. But they are much more interested in spaces $Y\simeq \Omega^n X$ for higher $n$, as such spaces carry a structure of \emph{$n$\hyp{}fold commutative} group. All the way up to infinity, there are spaces $Y$ that arise as infinitely many looped $X$'s, \ie spaces $Y$ such that $Y\simeq \Omega X$ for an $X$ which is $\Omega X'$, for an $X'$ which is $\Omega X''$\dots: these are called \emph{infinite loop spaces}, and they behave like abelian groups.

The notion of operad, introduced in \cite{may1972geometry}, offers a way to recognise infinite loop spaces among all spaces, as they are \emph{algebras} for a suitable operad. The reader is invited to consult \cite{adams1978infinite} for more information; Adams' book is one of the nicest introductions to the topic.

Shortly put, an operad is a family of spaces $O(n)$, one for each natural number $n$, subject to suitable axioms; one of which is that for every $k\in\N$, and every $k$-tuple of numbers $n_1,\dots,n_k$, there is a map
\[\gamma_{k,\vec n} : O(k)\times O(n_1)\times \dots\times O(n_k) \to O(n_1+\dots+n_k).\]
Since their very introduction it has been clear that operads are \emph{monoid\hyp{}like objects} in some category of functors, and that the maps above behave like multiplications of some sort: this is the reason why they can naturally act on other objects, and why the \emph{algebras} for an operad are so important;\footnote{This principle echoes Mt 7:20: \emph{Wherefore by their fruits ye shall know them}: the subject of the sentence are, of course, monoids; and the fruits are the monoid actions other object can carry.} making this analogy a precise statement, using the power of co\fshyp{}end calculus, is the content of a seminal paper by Max Kelly \cite{Kelly2005a}, that the present chapter follows extremely closely.

Co/end calculus is a perfect bookkeeping tool in otherwise extremely involved combinatorial arguments involving quotients of sets of $n$\hyp{}tuples by the action of a symmetric group.

Unfortunately, a thorough introduction to the theory of operads exceeds the aims of the present chapter: beginners (the author of the present note is undoubtedly among them) may feel rather disoriented when approaching any book on the subject, because algebraists might feel baffled by a geometric approach, and geometers might feel the same way reading about algebraic operads. So, it's extremely difficult to advise a single, comprehensive reference. Among classical textbooks, we can't help but mention \cite{may1972geometry}, and more recent monographies like \cite{loday2012algebraic, markl2007operads} written respectively from the algebraist's and topologist's\fshyp{}geometer's point of view. Among less classical and yet extremely valid points of view, the author profited a lot from a lucid, and unfortunately still unfinished, online draft \cite{Trimbled} written by T\@. Trimble.

We shall provide a glance on the use of operads in algebra in §\ref{some_adv} below; the exposition is terse, and maybe somewhat hasty, but keeps to a minimum the cognitive overload and tries to employ ideas and notation from the previous section.
\subsection{Local conventions}\label{prelim_notata} \leavevmode
\begin{itemize}\index{_aaa_P@$\P$}
	\item We will denote $\P$ the \emph{groupoid of natural numbers}, \ie the category having objects the nonempty sets $\{1,\dots, n\}$ (denoted as $n$ for short, assuming that $0 = \varnothing$) where $\P(m,n)=\varnothing$ if $n\neq m$ and $S_n$ (the group of bijections of $n$\hyp{}element sets) if $n=m$. It is evident that the groupoid $\P$ is the disjoint union of symmetric groups $\coprod_{n\ge 0} \text{Aut}(n) = \coprod_{n\ge 0} S_n$.\footnote{A subtlety that will never be mentioned in the discussion is that we blur the distinction between $\P$ and the \emph{free $\V$-category} $\overline{\P}^{\V}$ on $\P$, defined having the same objects, and where the hom-object between $n$and $m$ is the $\P(n,m)$-fold tensor of the monoidal unit $I\in\V$; to motivate the name `free', think of the case when $\V=\Ab$ is the category of abelian groups. Exercise \ref{operad_1} gives the proper definition.}
	\item like elsewhere, $\V$ is a fixed (symmetric) monoidal closed category, having all (weighted) co\fshyp{}limits needed to cast the relevant co\fshyp{}ends.
	\item We make a moderate use of $\lambda$\hyp{}notation: a function $x\mapsto Fx$ will be denoted as $\lambda x. Fx$ as if it were a $\lambda$\hyp{}term; the usual rules of $\alpha$\hyp{}conversion and $\beta$\hyp{}reduction straightforwardly apply.
	\item When needed, we freely employ the `Einstein notation' defined in \ref{einstein}; this will allow to maintain in a single line a few involved computations.
\end{itemize}
\begin{remark}\label{suppressed}
	A fundamental rule of Einstein convention is the following, really akin to the one for tensor operations: variables over which we integrate are always subscript\hyp{}superscript pairs. Moreover, monoidal structures symbols $\otimes$ are suppressed when this does not create ambiguity; thus, for example, the convolution \eqref{convognuscio} of two presheaves in Einstein notation has the following form:
	\[\notag
		F\ast G = \lambda X .\int^{C{C'}} \C^{C\otimes {C'}}_X F_C G_{C'}.
	\]
\end{remark}
\begin{remark*}
	We record that $\P$ has a natural choice of symmetric monoidal structure, with tensor the sum of natural numbers
	\[(n,m)\mapsto n+m = n \sqcup m;\]
	the action on arrows is given by $(\sigma,\tau)\mapsto \sigma+\tau$ defined acting as $\sigma$ on the set $\{1,\dots, m\}$ and as $\tau$ on the set $\{m+1,\dots, m+n\}$ (these permutations are called \emph{shuffles}).
\end{remark*}
\section{The convolution product}
\index{Convolution product}
The convolution product can be thought as a categorification of the convolution of regular functions: let $G$ be a topological group, and $C(G)$ the set of `regular', \ie continuous, complex\hyp{}valued functions $f : G \to \mathbb C$. Then, the set $C(G)$ can be endowed with a \emph{convolution} product, given by the integral $(f,g)(x) = \int_G f(xy^{-1})g(y)dy$ (not a coend!), once a suitable left\hyp{}invariant measure has been chosen on $G$.

In a similar fashion, if $\C$ is a \emph{monoidal} category, we can endow the category of functors $F : \C\to\V$ with a monoidal structure, which is in general different from the pointwise one. This is called \emph{convolution product} of functors.
\begin{proposition}[Day convolution]\label{day}
	\index{Convolution!Day ---}
	\index{Category!monoidal ---}
	\index{Product!convolution ---}\index{Convolution}
	\index{_aaa_ast@$\ast$}
	Let $\C$ be a symmetric monoidal $\V$\hyp{}category with monoidal product $\oplus$; the functor category $\VCat(\C,\V)$ is itself a symmetric monoidal category, (and in fact a cosmos if $\V$ is such) with respect to the monoidal structure given by \emph{Day convolution}: given $F,G\in\VCat(\C,\V)$ we define
	\[\label{convognuscio}
		F\ast G := \int^{CC'}\C(C\oplus C',\firstblank)\cdot FC\otimes GC'
	\]
	where we recall (see \ref{tenscotens}) that $X\cdot V$ for $X\in\Sets,V\in \V$ is the \emph{copower} (or \emph{tensor}) $X\cdot V$ such that
	\[
		\V(X\cdot V,W)\cong \Sets(X, \V(V,W)).
	\]
\end{proposition}
\index{Einstein notation}
\begin{proof*}
	We have to show that this really defines a monoidal structure:
	\begin{itemize}
		\item Associativity follows from the associativity of the tensor product on $\C$ and the ninja Yoneda lemma \ref{ninjayo}:\footnote{See \ref{suppressed} above for the way in which we employ Einstein notation; here and elsewhere, it is also harmless to suppress the distinction between monoidal products in $\V$ and $\V$\hyp{}tensors, since once the infix symbol has been removed to become mere juxtaposition, the two operations behave similarly; in all cases it can be easily devised which operation is which with a `dimensionality check'.}
		      \begin{align*}
			      F\ast(G\ast H) & = \lambda X . \int^{AB}\C^{A\oplus B}_X F_A (G\ast H)_B                           \\
			                     & \cong \lambda X . \int^{AB}\C^{A\oplus B}_X \int^{CD}\C^{C\oplus D}_B F_A G_C H_D \\
			                     & \cong \lambda X . \int^{ABCD}\C^{A\oplus B}_X \C^{C\oplus D}_B F_A G_C H_D        \\
			                     & \cong \lambda X . \int^{ACD}\C^{A\oplus(C\oplus D)}_X F_A G_C H_D
		      \end{align*}
		      A similar computation shows that
		      \[
			      (F\ast G)\ast H \cong \lambda X . \int^{ACD}\C^{(A\oplus C)\oplus D}_X F_A G_C H_D.
		      \]
		\item (Right) unitality : we show that $J=\yon_0=\C(0,\firstblank)$ plays the role of monoidal unit for the convolution $\ast$, if $0$ is the monoidal unit for $\oplus$: again, the ninja Yoneda lemma yields
		      \begin{align*}
			      F\ast J & \cong \lambda X . \int^{CD}\C^{C\oplus D}_X F_C J_D    \\
			              & \cong \lambda X . \int^{CD}\C^{C\oplus D}_X \C^0_D F_C \\
			              & \cong \lambda X . \int^C\C^{C\oplus 0}_X F_C\cong F.
		      \end{align*}
		      Similarly, we obtain left unitality.
	\end{itemize}
\end{proof*}
\index{Simplicial!--- subdivision}
\begin{example}[Subdivision and joins as convolutions]
	Compare Example \ref{ex.sub} and the definition of join of augmented\footnote{The category $\bDelta$ would have ordinal sum as monoidal structure, but it lacks an initial object $[-1] = \varnothing$ as monoidal unit; if we add such an object, we get a category $\bDelta_+$, and an \emph{augmented} simplicial set is a presheaf on $\bDelta_+$; the category of augmented simplicial sets is denoted $\sSet_+$. There is a triple of adjoints induced by the inclusion $i : \bDelta \subset \bDelta_+$ and linking the categories of simplicial and augmented simplicial sets.} simplicial sets given in \cite{Joy}: given $X, Y \in\sSet_+$ we define
	\[
		X \star Y = \int^{p,q} X_p \times Y_q \times \bDelta(\firstblank, p\oplus q)
	\]
	where $\oplus$ is the \emph{ordinal sum} operation.\index{Ordinal sum}
\end{example}
\begin{proposition}
	The convolution product of $F,G : \C \to \V$ has the universal property of the following left Kan extension:
	\[
		\vcenter{\xymatrix{
		\C \times \C \ar[d]_\oplus\ar[r]^{F\times G} & \V \times \V \ar[r]^\otimes & \V \\
		\C \ar@{.>}@/_1pc/[urr]_{F* G}&
		}}
	\]
\end{proposition}
\begin{proof}
	Just recognise that the coend expression given in \eqref{convognuscio} coincides with the coend formula for $\Lan_\oplus(F\times G)$ given in \ref{kanend}.
\end{proof}
Finding an explicit expression for the unit $\eta$ of this left extension is the scope of exercise \ref{ex:convo_al_lan}
\begin{remark}\label{its_closed_guys}
	The category $\VCat(\C,\V)$ is left and right closed: the internal hom $\llbracket G,H\rrbracket$ (or rather the functor $\llbracket G,\firstblank\rrbracket$ which is right adjoint to $\firstblank\ast G$) is given by
	\[
		\llbracket G,H\rrbracket := \lambda X. \int_C [GC,H(C\oplus X)]
	\]
	where $[\firstblank,\secondblank]$ is the internal hom in $\V$.
\end{remark}
(When $\C$ is not symmetric monoidal, we shall distinguish between a \emph{right} internal hom and a \emph{left} internal hom.)
\begin{proof*}
	We can compute directly:
	\begin{align*}
		\VCat(\C,\V)(F\ast G,H) & \cong \int_C \V \big( (F\ast G)C,HC\big)                                         \\
		                        & \cong \int_C {\V}\left( \int^{AB}\C^{A\oplus B}_C F_A G_B,\;H_C \right)          \\
		                        & \cong \int_{ABC}{\V}\big( \C^{A\oplus B}_C F_A G_B,\;H_C \big)                   \\
		                        & \cong \int_{ABC} {\V} \big(F_A, [\C^{A\oplus B}_C G_B,H_C]\big)                  \\
		                        & \cong \int_{ABC} {\V} \big(F_A, \big[G_B,[\C^{A\oplus B}_C,H_C]\big]\big)        \\
		                        & \cong \int_{AB} {\V} \big(F_A, \big[G_B, \int_C [\C^{A\oplus B}_C,H_C]\big]\big) \\
		                        & \cong \int_{AB} {\V} \big(F_A,[G_B,H_{A\oplus B}]\big)                           \\
		                        & \cong \int_A {\V}\left(F_A, \int_B[G_B,H_{A\oplus B}]\right)                     \\
		                        & \cong \int_A {\V}(F_A,\llbracket G,H\rrbracket_A)                                \\
		                        & \cong \VCat(\C,\V)\big(F,\llbracket G,H\rrbracket\big).
	\end{align*}
\end{proof*}
\begin{remark}
	We shall now specialise the above result to the particular case when $\C=\P$ with the sum monoidal structure of \ref{prelim_notata}; this means that $\VCat(\P,\V)$ is monoidal closed if we define
	\begin{align*}
		\displaystyle (F\ast G)_k                & := \int^{mn}\P(m+n,k)\cdot F_m \otimes G_n \\
		\displaystyle \llbracket F,G\rrbracket_k & := \int_n [F_n, G_{n+k}]
	\end{align*}
	In particular, we have a formula for the iterated convolution of $F_1,\dots, F_n\in [\C, \V]$, which we will make frequent use of:
	\[\label{useful_nfold}
		F_1\ast\dots\ast F_n = \int^{k_1,\dots, k_n} \P\big(\textstyle\sum k_i,\firstblank\big)\cdot F_1 k_1 \otimes \dots \otimes F_nk_n.
	\]
	Representing the tuple $k_1,\dots,k_n$ as a vector $\vec k$, we can also write \eqref{useful_nfold} in a more compact fashion as
	\[F_1\ast\dots\ast F_n = \int^{\vec k} \P(\sum \vec k, \firstblank)\cdot F\vec k.\notag\]
\end{remark}
\section{Substitution product and operads}
\index{Category!monoidal ---}
\index{Substitution product}
\index{Product!substitution ---}
\index{Operad}
\index{_aaa_odot@$\odot$}
The gist of the definition of a $\V$\hyp{}operad relies on the possibility to endow $\VCat(\P,\V)$ with an additional monoidal structure, defined by means of the Day convolution: this is called the \emph{substitution} product.
\begin{definition}[Substitution product on \protect{$\VCat(\P,\V)$}]\label{subbo_pro}
	Let $F,G\in \VCat(\P,\V)$. Define
	\[
		F\odot G:=\int^m Fm\otimes G^{\ast m},
	\]
	where $G^{\ast m} := G\ast\dots \ast G$ is the iterated convolution.
\end{definition}
Associativity exploits the following
\begin{lemma}
	There exists a natural equivalence $(F\odot G)^{\ast m}\cong F^{\ast m}\odot G$.
\end{lemma}
\begin{proof}
	It's a formal manipulation based on \eqref{useful_nfold} and the Yoneda lemma: the vector $\vec n$ denotes the $m$\hyp{}tuple of integers $(n_1,\dots,n_m)$ so that \eqref{useful_nfold} becomes $\bigast_{i=1}^n F_i \cong \int^{\vec k} \P\big(\sum k_i,\firstblank\Big)\cdot \bigotimes_{i=1}^n F_i k_i$.
	\begin{align*}
		(F\odot G)^{\ast m} & =\int^{\vec n} \P\Big(\sum n_i,\firstblank\Big)\cdot (F\odot G)n_1\otimes\dots\otimes (F\odot G)n_m                                                     \\
		                    & \cong \int^{\vec n,\vec k}\P\Big(\sum n_i,\firstblank\Big)\cdot Fk_1\otimes G^{\ast k_1}n_1\otimes\dots\otimes Fk_m\otimes G^{\ast k_m}n_m              \\
		                    & \cong \int^{\vec n,\vec k}Fk_1\otimes\dots\otimes Fk_m\otimes \P\Big(\sum n_i,\firstblank\Big)\cdot G^{\ast k_1}n_1\otimes\dots\otimes  G^{\ast k_m}n_m \\
		                    & \cong \int^{\vec k} Fk_1\otimes\dots\otimes Fk_m\otimes \big(G^{\ast k_1}\ast \dots\ast G^{\ast k_m}\big)                                               \\
		                    & \cong \int^{\vec k} Fk_1\otimes\dots\otimes Fk_m\otimes G^{\ast \sum k_i}                                                                               \\
		\ref{ninjayo}       & \cong \int^{\vec k,r}\P\Big(\sum k_i,r\Big)\otimes Fk_1\otimes\dots\otimes Fk_m\otimes G^{\ast r}                                                       \\
		                    & \cong \int^r F^{\ast m}r\ast G^{\ast r} = F^{\ast m}\odot G
	\end{align*}
	(Yoneda Lemma is used in the form
	\[G^{\ast n}\cong \int^r\P(n,t)\cdot Gt = \P(n,\firstblank)\odot G, \]
	because $(n,G)\mapsto G^{\ast n}$ is easily seen to be a bifunctor $\P \times [\P, \V]\to [\P, \V]$).
\end{proof}
Associativity of the substitution product now follows at once: we have
\begin{align*}
	F\odot(G\odot H) & =\lambda k. \int^m Fm\otimes (G\odot H)^{\ast m}k                   \\
	                 & \cong \lambda k. \int^m Fm\otimes (G^{\ast m}\odot H)k              \\
	                 & \cong \lambda k. \int^{m,l}Fm\otimes G^{\ast m}l\otimes H^{\ast l}k \\
	                 & \cong \lambda k. \int^l (F\odot G)l\otimes H^{\ast l}k              \\
	                 & = (F\odot G)\odot H.
\end{align*}
A unit object for the $\odot$\hyp{}product is $J=\P(1,\firstblank)\cdot I$; indeed $J(1)=I$, $J(n)=\varnothing_\V$ for any $n\neq 1$ and the ninja Yoneda lemma applies on both sides to show unitality rules:
\begin{itemize}
	\item On the left one has
	      \[
		      J\odot F = \int^m Jm\otimes F^{\ast m}=\int^m \P(1,m)\cdot F^{\ast m} \cong F^{\ast 1}=F.
	      \]
	\item On the right, $G\odot J\cong G$ once noticed that $J^{\ast m}\cong \P(m,\firstblank)\cdot I$ since
	      \begin{align*}
		      J^{\ast m}    & = \int^{\vec n}\P\Big(\sum n_i,\firstblank\Big)\cdot \P(1,n_1)\cdot\dots\cdot \P(1,n_m)\cdot I \\
		      \ref{ninjayo} & \cong \P(1+\dots+1,\firstblank)\cdot I                                                         \\
		                    & = \P(m,\firstblank)\cdot I
	      \end{align*}
	      because of Fubini rule \ref{fubozzo} and the Yoneda lemma, which says
	      \begin{multline}
		      \int^{\vec n}\P\left(n_1 + \cdots + n_m,\firstblank\right)\cdot \P(1,n_i)\cong \\ \cong \P(n_1+\dots+n_{i-1}+1+n_{i+1}+\dots + n_m,\firstblank),
	      \end{multline}
	      for any $1\le i\le m$. One has then
	      \[
		      G\odot J = \int^m Gm\otimes J^{\ast m}\cong \int^m Gm\otimes \P(m,\firstblank)\cdot I\cong G.
	      \]
\end{itemize}
\index{Category!closed ---}
The substitution product is highly non commutative; as an example take two representable presheaves and compose them in different order with respect to $\odot$. It is however closed on one side:
\begin{theorem}
	The $\odot$\hyp{}monoidal structure is \emph{left closed} in the sense that each $\firstblank\odot G$ has a right adjoint, but not right closed: not each $F\odot\firstblank$ has a right adjoint.
\end{theorem}
\begin{proof}
	To prove left closure, we compute:
	\begin{align*}
		\VCat(\C,\V)(F\odot G,H) & \cong \VCat(\C,\V)\Big( \int^m Fm\otimes G^{\ast m}, H \Big) \\
		                         & \cong \int_k \V\Big( \int^m Fm\otimes G^{\ast m}, H \Big)    \\
		                         & \cong \int_{km}\V(Fm, [G^{\ast m}k,Hk])                      \\
		                         & \cong \int_m \V(Fm, \int_k[G^{\ast m}k,Hk])                  \\
		                         & = \VCat(F, \{G,H\})
	\end{align*}
	if we define $\{G,H\}m=\int_k[G^{\ast m}k,Hk]$. Hence the functor $(\firstblank)\odot G$ has a right adjoint for any $G$.

	\index{Category!promonoidal ---}The functor $F\odot (\firstblank)$ can't have such an adjoint (Incidentally, this shows also that the substitution product can't come from a convolution product with respect to a \emph{promonoidal structure} in the sense of \ref{promonoshit}). We leave the reader to find a counterexample (find a colimit that is not preserved by $F\odot\firstblank$).
\end{proof}
\index{Operad}
\begin{definition}\label{opd_defn}
	An \emph{operad} in $\V$ consists of a $\odot$\hyp{}monoid object in $\VCat(\P,\V)$.

	\index{Unit!--- of a monad} In more explicit terms, an operad is a functor $T\in \VCat(\P,\V)$ endowed with a natural transformation called \emph{multiplication}, $\mu : T\odot T\to T$ and a \emph{unit} $\eta : J\to T$ such that
	\[\vcenter{\xymatrix{
				T\odot T \odot T \ar[d]_{\mu\odot T}\ar[r]^{T\odot\mu} & T\odot T\ar[d]^\mu & J\odot T\ar[dr]_\cong \ar[r]^{\eta\odot T} & T\odot T \ar[d]^\mu & \ar[l]_{T\odot\eta} T\odot J\ar[dl]^\cong \\
				T\odot T \ar[r]_\mu & T && T & \\
			}}\]
	are commutative diagrams.
\end{definition}
\index{Operad}
\begin{remark}\label{opd_explicit}
	Unraveling \ref{opd_defn}, we can notice that an operad in $\V$ consists of
	\begin{itemize}
		\item a natural transformation $\eta : J\To T$  that amounts to a map $\eta_1 : I\to T(1)$, since $J(1)=I, J(n)=\varnothing$ for $n\neq 1$;
		\item a natural transformation $\mu : T\odot T\To T$ that, in view of the universal property of the two coends involved, amounts to a cowedge
		      \[\vcenter{\xymatrix@C=2cm{
			      Tm \otimes \P(\vec n ,k)\cdot Tn_1\otimes\cdots\otimes Tn_m \ar[r]^-\tau & Tk
			      }}\]
		      for any $m,n_1,\dots, n_m,k\in\N$, natural in $k$ and the $n_i$ and such that the following diagram commutes:
		      \[\vcenter{\xymatrix@C=1.5cm{
			      Tm \otimes \P(\sum n_i ,k)\cdot T^\otimes \vec n\ar[r]^{\sigma^*}\ar[d]_{\sigma_*} & Tm \otimes \P(\sum n_{\sigma i} ,k)\cdot T^\otimes \vec n\ar[d]\\
			      Tm \otimes \P(\sum n_i ,k)\cdot T^{\otimes\sigma} \vec n \ar[r]& Tk
			      }}\]
		      (the notation is self\hyp{}evident) for every morphism $\sigma\in \P$. This is equivalent to a transformation
		      \[\vcenter{\xymatrix@C=1.5cm{
			      Tm \otimes Tn_1\otimes\cdots\otimes Tn_m \ar[r]^{\hat\tau} & [\P(\sum n_i, \firstblank), T(\firstblank)]
			      }}\]
		      (considering the $n_i$ fixed and the first functor constant in $k$) \ie, by the Yoneda Lemma a natural transformation from $Tm \otimes Tn_1\otimes\cdots\otimes Tn_m$ to $T(\sum n_i)$.
	\end{itemize}
\end{remark}
\begin{example}[Examples of operads]\label{eo}\leavevmode
	\begin{enumtag}{eo}
		\item \label{eo:uno} For any $F\in \VCat(\P,\V)$ the object $\{F,F\}$ is an operad whose multiplication is the adjunct of the arrow
		\[\vcenter{\xymatrix@C=2cm{
			\{F,F\}\odot \{F,F\}\odot F \ar[r]^-{\{F,F\}\odot \epsilon}& \{F,F\}\odot F \ar[r]^-{\text{ev}}& F
			}}\]
		($\epsilon$ is the counit of the $\firstblank\odot F \dashv \{F,\firstblank\}$ adjunction) and whose unit is the adjunct of the isomorphism $J\odot F\cong F$. This is called the \emph{endomorphism operad}.
		\index{Operad!--- from a monad}
		\item \label{eo:due} (see \cite[2.2.11]{Leinster2004}) Let $T : \Set \to \Set$ be a monad; then, the collection $(Tn\mid n\in \Set)$ (\ie the restriction of $T$ on \emph{finite} sets) defines an operad if we take
		\begin{itemize}
			\item the unit $\eta_1 : 1 \to T1$ to be the component of the unit of $T$ at the singleton set;
			\item the multiplication maps (see \ref{opd_explicit} above)
			      \[\gamma_{m, \vec n} : Tm \times Tn_1\times\cdots \times Tn_m \to T\left(\sum n_i\right)\]
			      via the following rule: let $n = \sum n_i$, and let us consider the image under $T$ of the coproduct inclusion $u_i : n_i \to n = n_1\sqcup \dots\sqcup n_m$, \ie the map $Tu_i : Tn_i \to Tn$; this defines a map $\prod_{i=1}^m Tn_i \to \prod_{i=1}^m Tn = (Tn)^m = \Set(m, Tn)$, that can be post\hyp{}composed with the action of $T$ on arrows, and with the multiplication of the monad, obtaining
			      \[Tn_1\times \cdots\times Tn_m \to \Set(m, Tn) \to \Set(Tm, TTn) \xto{\Set(Tm, \mu_n)} \Set(Tm, Tn);\]
			      the transpose of this map is the desired $\gamma_{m,\vec n}$.

			      As a motivation for this, we recall that in an algebraic theory the free algebra $T(A)$ is the set of all terms (inductively) built applying the operations of the theory $T$ to the `variables' in $A$; in this interpretation, there is a substitution map $\gamma_{m,\vec n}$ defined by sending a term $t\in Tm$ in the variables $x_1,\dots, x_m$, and $m$ terms $g_1,\dots, g_m$ in the variables $(x_1^j,\dots, x_{n_j}^j)$ with $j = 1,\dots, m$; $\gamma_{m,\vec n}(t; g_1,\dots, g_m)$ is then defined as $t[x_j := g_j]$ (a moderate amount of knowledge of $\lambda$\hyp{}calculus will allow to appreciate the notation).
		\end{itemize}
		\item \label{eo:ter} As a particular case, the monad `free commutative $k$\hyp{}algebra' yields that the family of polynomial algebras $(k[X_1,\dots, X_n] \mid n\in \mathbb N)$ defines an operad if the unit $k \to k[X]$ chooses the multiplicative identity, and the multiplication (if $n=\sum n_i$)
		\[k[T_1,\dots, T_m]\otimes k[X_1^1,\dots,X^1_{n_1}]\otimes\dots\otimes k[X^m_1,\dots, X^m_{n_m}] \to k[X_1,\dots,X_n]\]
		is defined by the rule
		\[(q; g_1,\dots,g_m)\mapsto q(g_1,\dots,g_m)\]
		(evaluation of a $m$\hyp{}variable polynomial into an $m$\hyp{}tuple of polynomials in $n_1,n_2,\dots,n_m$ variables each). This is evidently an associative and unital composition law, compatible with the action of the permutation group on $n$ elements.
	\end{enumtag}
\end{example}
As already said, every attempt to propose a self\hyp{}contained treatment of operads in a single chapter would turn into a complete failure: the theory is too large to be reduced to a bunch of not\hyp{}too\hyp{}technical concepts. In some sense, this is to be expected: an operad is nothing but a multicategory (see \cite{Leinster2004}) with a single object; the resulting theory must then be at least as expressive as the theory of categories.

Thus, we contempt ourselves to record a few of the less\hyp{}than\hyp{}elementary results (not without a certain dose of cherry\hyp{}picking towards the ones that can be easily expressed using co\fshyp{}end calculus); this should be sufficient to clarify how operads are a fundamental tool to describe `stuff endowed with $n$\hyp{}ary operations' that are coherent with the environment.
\subsection{Substitution and algebraic theories}
\index{Algebraic theory}
Once one familiarises with their definition, they might find that operad-like structures are way more common than expected. The scope of the present subsection is to convey an intuitive explanation for their ubiquity, while at the same time showing, as it should be, that most of the gadgets devised to axiomatise ``structures borne by operations and properties thereof'' form equivalent categories. The nontrivial connection between combinatorics and the theory of structures defined by operations and equations is still far from being perfectly understood, despite the enormous effort spent in studying its main features.

The subsection builds on various classical presentations of the subject, such as \cite{lawvere1963functorial,hyland2007category,lack2011notions}, and contains nothing essentially new.

Let $\Fin$ be the category of finite sets and functions; as we have already seen, the opposite category $\Fin^\opp$ exhibits the universal property of the free category with products generated by the point; namely, if $\C$ is any category with finite products, the subcategory of product\hyp{}preserving functors $F : \Fin^\opp \to \C$ is equivalent to $\C$:
\[
	\C \cong \Cat(*,\C)\cong \Cat_\times(\Fin^\opp,\C).
\]
In more concrete terms, this means that a product-preserving functor $p : \Fin^\opp \to \C$ is uniquely determined by the image of $[1]\in\Fin$, because it must send every other $[n]$ to $p[n]=p[1]^n$.

\index{Functor!identity on objects ---}
An \emph{identity on objects} functor (``idonob'' functor, for short) is just a functor $F : \C \to \D$ whose function on objects $F_o : \C_o \to \D_o$ is the identity of the set/class $\C_o = \D_o$.
\begin{definition}\label{opd_equiv_1}
	\index{Lawvere!--- theory}
	\index{Lawvere}
	We call a \emph{Lawvere theory} an idonob functor $p : \Fin^\opp \to \LL$, that strictly preserves products.
\end{definition}
The category $\BF{Law}$ of Lawvere theories is thus defined as having objects the pairs $(p,\LL)$ as above, and morphisms $(p,\LL)\to (p',\LL')$ the functors $h : \LL \to \LL'$ such that $h\circ p=p'$; note that this says that also $h$ is bijective on objects and product preserving.

When regarded as subcategory of the coslice $\Fin^\opp/\Cat$, the category of Lawvere theory seems to be a 2-category; the request that a Lawvere theory is bijective on objects, however, forces a natural transformation $\alpha : h\To h'$ to be the identity, thus proving that $\BF{Law}$ is a \emph{2-discrete} 2-category.
\index{Monad!finitary ---}
\index{Monad}
\index{Finitary monad}
\begin{definition}[Finitary monad]\label{opd_equiv_2}
	A monad (see \ref{def:monad}) on $\Set$ is \emph{finitary} if it preserves colimits over a filtered category (see \ref{def:filtcat}); finitary monads form a full subcategory of the category of monads, where a morphism of monads has been defined in \ref{rmk_monad_mor}.
\end{definition}
\index{Functor!finitary ---}
\begin{definition}\label{opd_equiv_3}
	Let us consider the category $\Cat(\Fin,\Set)$ of all functors $F : \Fin \to\Set$, endowed with its cartesian monoidal structure; by mimicking the construction of operads in \ref{subbo_pro} and \ref{opd_defn}, we can define a \emph{substitution} monoidal structure on $\Cat(\Fin,\Set)$ as follows (evidently, $\Fin(n,m)$ is a much larger set of morphisms than $\P(n,m)$, so the equaliser defining the following coend performs a much smaller quotient):
	\[F \odot G : n \mapsto \int^m Fm\times G^mn\]
	where the $m^\text{th}$ power of $G$, $G^mn = Gn\times \dots \times Gn$, is the product of $Gn$ with itself repeated $m$ times, and it coincides with the $m$-fold convolution of $G$ with itself, when both $\Fin$ and $\Set$ are given their cartesian monoidal structure.

	The monoids with respect to this monoidal structure are called \emph{cartesian operads} or \emph{clones} in \cite{curiennone} and in our \ref{def_clone} below.

	When the obvious notion of homomorphism of internal monoids is specialised to this particular case, we obtain a category $\BF{Clo}$ of clones.
\end{definition}\index{Clone}
\begin{lemma}\label{finita}
	There is an equivalence between the category of functors $\Fin \to \Set$ and the category $\Cat_\omega(\Set,\Set)$ of finitary (see \ref{accepre}) endofunctors of $\Set$; this equivalence is monoidal, when $\Cat(\Fin,\Set)$ is endowed with the substitution product defined above.
\end{lemma}
\begin{proof}
	The equivalence $\Cat(\Fin,\Set)\cong \Cat_\omega(\Set,\Set)$ is induced by the adjunction $\Lan_J \dashv J^*$, if $J : \Fin\subset \Set$, and the category of finitary endofunctors of $\Set$ is clearly monoidal with respect to composition.

	It is a general fact that given an equivalence of categories $F : \C \leftrightarrows \D : G$, if $(\D,\otimes,I)$ is monoidal we can define a monoidal structure on $\C$ that turns the adjunction $(F,G)$ into a monoidal equivalence, by setting $C \diamond C' := G(FC \otimes FC')$ and $J := GI$.

	It remains to show that the composition  on $\Cat_\omega(\Set,\Set)$, when transported to $\Cat(\Fin,\Set)$ along the equivalence, coincides with the substitution product, \ie that
	\[\Lan_JF\circ \Lan_JG\cong\Lan_J(F\odot G)\]
	for every $F,G: \Fin\to\Set $ if $J : \Fin \subset \Set$ is the inclusion, and vice-versa that
	\[UJ\odot VJ \cong U\circ V\circ J\]
	if $U,V : \Set \to \Set$ are finitary endofunctors.

	We can easily show both these equations using coend calculus:
	\begin{align*}
		\Lan_J(F\odot G)A & \cong \int^n \Set(n,A) \times (F \odot G)(n)       \\
		                  & \cong \int^n \Set(n,A) \times \int^m Fm\times G^mn \\
		                  & \cong \int^{nm} \Set(n,A) \times Fm\times G^mn     \\
		                  & \cong \int^m Fm \times \Set(m,\Lan_JGA)            \\
		                  & \cong Lan_JF(\Lan_JGA),
	\end{align*}
	where the last passage is motivated by the isomorphism
	\[
		\int^n \Set(n,A) \times G^mn\cong \Lan_J(\Set(m,G\firstblank))(A)\cong \Set(m, \Lan_JGA);
	\]
	conversely,
	\begin{align*}
		(UJ\odot VJ)n & \cong \int^m UJm\times (Vn)^m       \\
		              & \cong \int^m UJm \times \Set(m, Vn) \\
		              & \cong \int^m Um \times \Set(Jm, Vn) \\
		              & \cong \Lan_J(UJ)(Vn) = UVJn;
	\end{align*}
	This concludes the proof.
\end{proof}
\index{Clone}
\begin{corollary}\label{this_gia_fatto}
	There is an equivalence between the category of finitary monads $T : \Set \to \Set$ and the category $\BF{Clo}$ of clones.
\end{corollary}
\index{Functor!cmc ---}
\index{Monad!cmc ---}
\index{cmc functor}
\index{cmc monad}
\begin{definition}\label{opd_equiv_4}
	Let $T : \Cat(\Fin,\Set) \to \Cat(\Fin,\Set)$ be a monad; we call it a \emph{cmc functor} if it preserves colimits and finite products (which, in this case, coincide with Day convolution as already said). We denote $\BF{CMC}(\Fin,\Set)$ the category of cmc monads on $\Cat(\Fin,\Set)$ and monad morphisms (see \ref{rmk_monad_mor}).
\end{definition}
The reader can find a more detailed account of this definition in \ref{ciemmeci} below.
\begin{definition}\label{opd_equiv_5}
	Recall that the category $\Dist(\Fin^\opp,\Fin^\opp)$ is monoidal with respect to composition of profunctors; a monoid internal to it is a \emph{promonad}, \ie a profunctor $\mathfrak{t} : \Fin^\opp\pto \Fin^\opp$ that in addition is a monoid with respect to profunctor composition.
\end{definition}
All these seemingly diverse structures turn out, instead, to be the same:
\begin{theorem}
	There is an equivalence between the categories of \ref{opd_equiv_1}, \ref{opd_equiv_2}, \ref{opd_equiv_3}, \ref{opd_equiv_4}, \ref{opd_equiv_5}.
\end{theorem}
The equivalence between \ref{opd_equiv_2} and \ref{opd_equiv_3} follows from \ref{finita} and in particular from the fact that the equivalence is monoidal: it has been already said in \ref{this_gia_fatto}.

The equivalence between \ref{opd_equiv_4} and \ref{opd_equiv_5} follows from the fact that the equivalence in \ref{finita} can be promoted to a monoidal equivalence as well, when both hom-categories are endowed with the composition monoidal structure; in our case, this boils down to the statement that a promonad $\mathfrak{t}$ correspond to a cocontinuous monad $T$ on $\Cat(\Fin,\Set)$ under the equivalence of \ref{finita} and it is the content of Exercise \ref{ex5:promonads_go_to_cmc}.

The remaining implications need an intricate series of definitions and preliminary results, thus their proof occupies the following two subsections.
\subsubsection{Algebraic theories are finitary monads}
\index{Lawvere theory}
\index{Monad!finitary ---}
We will prove the equivalence between \ref{opd_equiv_1} and \ref{opd_equiv_2} building mutually functors in opposite directions between the category $\BF{Law}$ and the category of finitary monads on $\Set$.
\begin{remark}
	A preliminary remark is in order: the notion of Lawvere theory abstracts the notion of algebraic theory in the sense that one can define a category $\LL_M$ with the following properties:
	\begin{itemize}
		\item $\LL_M$ has objects the natural numbers $[0],[1],\dots$, and the operation $[n],[m] \mapsto [n+m]$ is a functor $\LL_M \times \LL_M \to \LL_M$ endowing $\LL_M$ with (strictly associative) finite products.
		\item $\LL_M$ contains morphisms $\mu : [2]\to [1], \eta : [0] \to [1]$ subject to equations given by the commutative diagrams
		      \[
			      \vcenter{
			      \xymatrix{
			      [2+1] = [1+2] \ar[r]^-{1+\mu}\ar[d]_{\mu+1} & [2] \ar[d]^\mu & [0+1]=[1+0]\ar@{=}[dr]\ar[r]^-{1+\eta}\ar[d]_{\eta+1}& [1+1]\ar[d]^\mu\\
			      [2]\ar[r]_\mu & [1] & [1+1] \ar[r]_\mu & [1]
			      }}
		      \]
		\item $\LL_M$ has a category of \emph{models}, \ie functors $F : \LL_M \to \Set$ such that $F[n]=F[1]^n$, and thus the functions $F\mu : A^2\to A$, $F\eta : A^0 = *\to A$ endow the set $A=F[1]$ with an associative, unital binary operation; in short, $A$ is a \emph{monoid}. A natural transformation between models turns out to be just a homomorphism between the monoids.
	\end{itemize}
	$\LL_M$ can thus be thought as the ``theory'' whose models $F : \LL_M \to \Set$ are monoids; it is a category harbouring the most abstract shape a monoid can possibly have, and a concrete realisation of such an abstraction inside set theory turns out to be just a functor $F : \LL_M \to \Set$.
\end{remark}
\begin{remark}
	At this point, the reader might want to embark in an easy exercise: are natural transformations $F \To G : \LL_M \to \Set$ of models really in bijection with monoid homomorphisms? Meaning: shouldn't we ask for a natural transformation to preserve the cartesian product itself? It turns out that this is not needed, in that such an $\alpha : F \To G$ \emph{must} preserve products, in the sense that the map $\alpha_{[n]} : F[1]^n \to G[1]^n$ must be $\alpha_{[1]}^n$.
\end{remark}
Given a general Lawvere theory $p : \Fin^\opp \to \LL$ we can form the category of its models as the functors
\[F : \LL \to \Set\]
that preserve products; such models are uniquely determined by their action on the object $[1]\in\LL$, in the sense that if $F$ is such a model, $F[n]\cong F[1]^n$ for every $[n]\in\Fin^\opp$, so we can characterise the subcategory $\Mod(p,\LL)\subset \Cat(\LL,\Set)$ as the upper left corner in the pullback
\[\label{modichcio}
	\vcenter{\xymatrix{
	\Mod(p,\LL) \pb \ar[d]_U\ar[r]^-j & \Cat(\LL,\Set)\ar[d]^{p^*} \\
	\Set \ar[r]_-{N_J} & \Cat(\Fin^\opp,\Set)
	}}
\]
where $N_J : A\mapsto (\lambda n.\Set(n,A))$ is the nerve (see \ref{nervereal}) associated to the functor $J : \Fin\subset\Set$, and $U$ the functor that evaluates a model $F : \LL \to \Set$ on the object $[1]$.
\begin{remark}
	The functor $U$ commutes with all limits and with filtered colimits.
\end{remark}
\begin{proof}
	Every limit of a diagram $(F_i)$ of product-preserving functors is still product-preserving, since limits commute with limits.\footnote{More formally: if $\bsmat n \\ m \esmat$ is the functor $\{0,1\}\to\C$ choosing two objects $n,m\in \LL$, then $[m]\times [n]\cong \lim_{\{0,1\}}\bsmat n \\ m \esmat$ and then one has $$\textstyle \lim_{\I}\lim_{\{0,1\}}F_i\bsmat n\\m\esmat \cong \lim_{\{0,1\}}\lim_\I F_i\bsmat n\\m\esmat,$$ by virtue of the Fubini rule for limits.}
	Now, the ``evaluation at $[1]$'' functor $U$ is easily seen to be isomorphic to the functor
	\[F \mapsto \Cat(\LL,\Set)(\yon[1],F)\]
	by virtue of the Yoneda lemma; thus, it preserves all limits (because every covariant representable $\C(X,\firstblank)$ does), and all filtered colimits, since the representables are finitely presentable (see \ref{accepre}) objects of $\Mod(p,\LL)$.
\end{proof}
The functor $U$ does not, however, preserve all colimits: for example, it doesn't always preserve coproducts or initial objects (consider, for example, the category of monoids $\Mod(\LL_M)$).\footnote{More formally, the representable object $y[1]$ on $[1]$ is tiny (\ie $\hom(y[1],\firstblank)$ is cocontinuous) as an object of $\Cat(\LL,\Set)$ by virtue of the Yoneda lemma, \emph{but not} as an object of $\Mod(p,\LL)$, as there are colimits in $\Mod(p,\LL)$ that are constructed differently from the way they are constructed in $\Cat(\LL,\Set)$.}
\begin{lemma}
	The functor $U$ admits a left adjoint.
\end{lemma}
\begin{proof}
	Given a set $A$, consider the functor $N_J(A) : \Fin^\opp\to \Set$, and the functor $p : \Fin^\opp\to \LL$ defining the theory; we shall prove that $A\mapsto \Lan_p(N_J(A))$ in the diagram
	\[\vcenter{\xymatrix{
				&\Fin \ar[dl]_p \ar[dr]^{N_JA}& \\
				\LL \ar[rr]_{\Lan_pN_JA}&& \Set
			}}\]
	defines a functor that is a left adjoint to $U$: to do this, we call $F_pA := \int^n \LL(m,n)\times A^n$ and prove that there is an adjunction $\Mod(p,\LL)(F_pA,H)\cong \Set(A,H[1])$ for every model $H$ of $p$. In order to see that, one can expand the definition of the left hand side:
	\begin{align*}
		\Mod(p,\LL)(F_pA,H) & \cong \int_n \Set(F_pAn, Hn)                                \\
		                    & \cong \int_n \Set\left(\int^m \LL(n,m)\times A^m, Hn\right) \\
		                    & \cong \int_{nm} \Set\left( \LL(n,m)\times A^m, Hn\right)    \\
		                    & \cong\int_{nm} \Set(A^m,\Set(\LL(n,m),Hn))                  \\
		                    & \cong \int_m \Set(A^m, H[1]^m) \cong \Set(A,H[1])
	\end{align*}
	where in the last steps we used the ninja Yoneda lemma twice.
\end{proof}
So, given a Lawvere theory $(p, \LL)$ we can associate a monad to it, precisely the monad induced by the adjunction $F_p\dashv U$, \ie to the functor $UF_p : \Set \to \Set$.
\begin{proposition}
	The monad obtained in this way is finitary.
\end{proposition}
\begin{proof}
	The functor $U$ commutes with finitely filtered colimits, and $F_p$ commutes with all colimits because it is a left adjoint; thus, the composition $UF_p$ is finitary.
\end{proof}
Finally, let's see that the adjunction $F_p \dashv U$ is \emph{monadic} in the sense of \ref{monadic_fu}: this means that the category of $UF_p$-algebras coincides up to equivalence with the category of $(p,\LL)$-modules.

It is rather easy to prove that $U$ is conservative having in mind the pullback in \eqref{modichcio}; it remains to show that $U$ creates the coequalisers of $U$-split pairs, in order to fulfil all the requests of \ref{monadici}; this can be checked directly, and in fact it is a general result about pulling back a monadic functor along another functor; we leave the details as an exercise in \ref{pull_is_mona}; note that it is not possible to remove the assumption that the lower horizontal functor in the pullback \eqref{modichcio} is fully faithful.

We now construct a correspondence in the opposite direction: given a finitary monad $T$ on $\Set$ we consider the composition
\[\xymatrix{\Fin \ar[r]^J & \Set \ar[r]^{F^T} & \Alg(T)}\]
where the functor $F^T$ is the free functor of \ref{frealg}, and its factorisation as an idonob functor $p : \Fin \to \LL^\opp$ followed by a fully faithful functor $\LL^\opp\to \Alg(T)$. The functor $p$ now is a Lawvere theory, and its category of models coincides with the category of $T$-algebras (see \ref{alg_for_a_mona}).

These two correspondences in opposite directions extend to functors as follows. On the category of finitary monads we take only \emph{restrained} morphisms in the sense of \ref{restrained_mormonad}.
\begin{itemize}
	\item Given a morphism of Lawvere theories $h : (p,\LL) \to (q,\M)$, \ie a commutative triangle
	      \[\vcenter{\xymatrix{
				      &\Fin^\opp \ar[dr]^q\ar[dl]_p& \\
				      \LL \ar[rr]_h && \M
			      }}\]
	      Let's say $p$ gives rise to the monad $UF$ and $q$ to the monad $U'F'$; then, we define a natural transformation
	      \[\lambda : UF \xto{UF * \eta'} UFU'F' \to UFU\bar k F' \xto{U * \epsilon * \bar k F'} U\bar kF' = U'F'\]
	      using the unit of $F'\dashv U'$, the counit of $F\dashv U$ and the functor $\bar h : \Mod(q,\M) \to \Mod(p,\LL)$ obtained from the universal property of the pullback, from $k^* : \Cat(\M,\Set)\to \Cat(\LL,\Set)$.

	      It is now a matter of using the zig-zag identities of the two adjunctions, to show that $\lambda$ satisfies the axioms of \ref{rmk_monad_mor}; we leave this as an unenlightening exercise for the reader.
	\item Given finitary monads $T,S$ on $\Set$, and a restrained morphism of monads $\lambda : S \To T$, we shall obtain a morphism of factorisations to fill in the center of the diagram
	      \[\label{quello}\vcenter{\xymatrix{
				      \Fin \ar[dr]_q\ar@{=}[rrr]\ar[d]_J &&& \Fin \ar[dl]^p\ar[d]^J \\
				      \Set \ar[d]_{F^S}& \M \ar[dl]&\LL\ar[dr]&\Set\ar[d]^{F^T} \\
				      \Alg(S) \ar[rrr]_{\lambda_\sharp}&&& \Alg(T).
			      }}\]
	      It is evident that $\lambda$ induces a morphism between $T$- and $S$-algebras
	      \[\lambda^\sharp : \Alg(S) \to \Alg(T)\]
	      defined sending a $T$-algebra $\var{TA}{A}$ into $SA \xto{\lambda_A} TA \to A$; the axioms satisfied by a monad morphism now entail that this composition is a $S$-algebra. However, the outer rectangle in \eqref{quello} does not commute, because $\lambda^\sharp$ does not restrict to a morphisms of \emph{free} algebras; in fact, the triangle
	      \[\vcenter{\xymatrix{
				      &\Set\ar[dr]^{F^T}\ar[dl]_{F^S}&\\
				      \Alg(S) \ar[rr]_{\lambda_\sharp} && \Alg(T)
			      }}\]
	      has no reason to commute.

	      Instead, \emph{if} $\lambda^\sharp$ had a left adjoint $\lambda_!$, it would be true that $\lambda_! \circ F^T = F^S$, and we could use it to define a morphism between the factorisations $\LL,\M$ in \eqref{quello} above, and thus a morphism of theories.

	      It turns out that $\lambda^\sharp$ is indeed a right adjoint: both $\Alg(S),\Alg(T)$ are complete and accessible categories, $\lambda^\sharp$ commutes with limits and filtered colimits because $S,T$ are finitary monads, and thus for the adjoint functor theorem \cite[2.45]{Adamek1994} it must admit a left adjoint $\lambda_!$. An explicit description of such left adjoint as a sequential colimit that exploits finite accessibility of $\Alg(S),\Alg(T)$ is left as Exercise \ref{the_soa_at_work}.
\end{itemize}
\subsubsection{Substitution monoids are cocontinuous monads}
To conclude the section, we prove the equivalence between \ref{opd_equiv_3} and \ref{opd_equiv_4}. This is a sensibly slicker argument, thanks to coend calculus!

Let $T$ be an $\odot$-monoid; it is a general fact that tensoring with an internal monoid in a monoidal category gives a monad $\firstblank\odot T : \Cat(\Fin,\Set) \to \Cat(\Fin,\Set)$; this gives a way to construct a monad out of $T$; this correspondence is of course functorial, as $\odot$ is a bifunctor. The resulting monad is furthermore cocontinuous, since the substitution monoidal structure is left closed (see \ref{its_closed_guys}). We leave to the reader the proof that the monad is convolution preserving.

On the other hand, given a cocontinuous monad $S$ on $\Cat(\Fin,\Set)$, we evaluate $S$ on the representable $J=y[1]$, \ie on the substitution monoidal unit, and it is now a matter of coend calculus to show that
\[S(A) \cong A\odot SJ,\]
so that every cocontinuous monad arises this way:
\begin{align*}
	(A \odot SJ) & = \int^m Am \times Sy[1]^m                                         \\
	             & \cong \int^m Am \times \Set(m, S(y[1])\firstblank)                 \\
	             & \cong \int^m S\left(Am \times \Set(m, y[1]\firstblank)\right)      \\
	             & \cong S\left(\int^m Am \times \Set(m, \firstblank)\right) = (SA)n.
\end{align*}
This concludes the proof of the equivalence between cmc monads and substitution monoids.
\section{Some more advanced results}\label{some_adv}
We expand a little bit more on the theory of operads, also putting in a broader perspective the results we have already encountered in the previous section (notably, the notion of cartesian operad or \emph{clone}).

We begin with a theorem linking operads and monads: every $\V$\hyp{}operad induces a monad on $\V$.
\index{Monad!--- from an operad}
\begin{theorem}
	Let $T : \P \to \V$ be an operad; then we can define a monad on $\V$ by the rule
	\[\label{opd_to_monad}
		M_T : A \mapsto \int^{n\in\P} Tn\otimes A^n.
	\]
\end{theorem}
\begin{proof}
	We shall first define unit and multiplication:
	\begin{itemize}
		\item To define the unit, we employ Exercise \ref{opd_adjoints}: there is a functor $\Psi : \V \to \VCat(\P, \V)$, precisely the left adjoint to evaluation at $0\in\P$, that acts as follows: $\Psi(A)$ is the functor that sends $0$ to $A$, and all other $n\ge 1$ to the initial object of $\V$. Unwinding this definition, we see that $A\cong \V(I,\firstblank)\odot \Psi A$ and that $T\odot \Psi A = M_TA$. This means that we can define the components of the unit as
		      \[
			      C \cong \V(I,\firstblank)\odot \Psi C \xto{\eta \odot \Psi C} T\odot \Psi C = M_T(C)
		      \]
		\item To define the multiplication, let us follow the string of equivalences
		      \begin{align*}
			      \int^m Tm \otimes (M_T(C))^m & \cong \int^m Tm \otimes \left[\int^n Tn \otimes C^n\right]^m                                     \\
			                                   & \cong \int^m Tm \otimes  \int^{\vec n} Tn_1\otimes \cdots\otimes Tn_m \otimes C^{n_1+\cdots+n_m} \\
			                                   & \cong \int^{m,\vec n}Tm \otimes Tn_1\otimes \cdots\otimes Tn_m \otimes C^{n_1+\cdots+n_m}
		      \end{align*}
		      Now, using the multiplication of the operad we get that this object maps canonically to $\int^k Tk\otimes C^k$.
	\end{itemize}
	With a certain amount of work, associativity and unitality for $T$ now prove the associativity and unitality for the endofunctor $M_T$, thus showing that it is a monad (see \ref{def:monad}).
\end{proof}
\begin{remark}
	The correspondence described in \ref{eo}.\ref{eo:due} can be promoted to a functor from the category of monads to the category of $Set$\hyp{}operads; it pairs with the correspondence described above, yielding an adjunction between the category of monads and the category of operads. This functor is however not an equivalence, because there are non\hyp{}isomorphic operads that generate the same monad (see \cite{leinster2006operads}).
\end{remark}
Our analysis continues with the main theorem in \cite{curiennone}:
\begin{quote}
	There exist 2\hyp{}comonads on the bicategory $\Dist$ [\emph{author: see our \ref{profdef}}] whose coKleisli category has objects respectively operads, symmetric operads and clones.
\end{quote}\index{Clone}
We begin by introducing a special class of operads of particular interest in universal algebra, called \emph{clones}. In short, a clone is an operad with respect to a cartesian monoidal structure.
\begin{definition}\label{ciemmeci}
	A category $\C$ is called cartesian monoidally cocomplete (\emph{cmc} for short) if it admits finite products, all colimits, and each functor $A \times \firstblank : \C \to \C$ preserves these colimits.
\end{definition}
\index{_aaa_Fin@$\Fin$}
Recall from Exercise \ref{monoidaux} that we denote $\Fin_*$ the category $1/\Set_{<\omega}$ of pointed finite sets; in the following, $\Fin$ will denote the category of unpointed finite sets.
\begin{proposition}\label{lemmone}
	Let $\C$ be a cmc category; then are equivalences of categories
	\[
		\C
		\cong \Cat_\times(\Fin^\opp, \C)
		\cong \Cat_{\times,!}(\Cat(\Fin, \Set), \C)
	\]
	of $\C$ with the category of product preserving functors $\Fin \to \C$, and with the category of colimit\hyp{} and product\hyp{}preserving functors $\Fin \to \Set$ (we call them \emph{cmc functors} for short: see \ref{ciemmeci} for the motivation behind this notation), induced by evaluation at the inclusion functor $\Fin\subset \Set$ (this is obviously colimit- and product\hyp{}preserving).
\end{proposition}
\begin{proof}
	The second equivalence is just the universal property of the Yoneda embedding described in \ref{yext_are_good} and specialised to this context.

	As for the first equivalence, an object $C\in\C$ goes to the functor $W : [n]\mapsto C^n$. Each such functor is uniquely determined by its action on $[1]$, since $[n] \cong [1] \times \cdots\times [1]$ in $\Fin^\opp$ (or equivalently, $[n] = [1] \sqcup\cdots\sqcup [1]$ in $\Fin$).
\end{proof}
We can actually trace the image of an object $C$ into the category of cmc functors: $C$ goes first to the functor $[n]\mapsto C^n$, and then to the functor that sends $F : \Fin\to \Set$ to its $W$\hyp{}weighted colimit
\[
	\int^{[n]\in\Fin} Fn \otimes C^n
\]
where $\otimes$ is the usual tensor of a cocomplete category over $\Set$ (see \ref{tenscotens}).

An immediate corollary of our \ref{lemmone} is that we have an equivalence of categories
\[
	\Cat_{\times,!}(\Cat(\Fin,\Set), \Cat(\Fin,\Set))\cong \Cat(\Fin,\Set)
\]
between the cmc\hyp{}endofunctors of $\Cat(\Fin,\Set)$ and itself. In fact, something more general is true (see \cite{Trimbled}): the monoidal product on endofunctors given by composition transfers across the equivalence to a monoidal product on $\Cat(\Fin,\Set)$, and this monoidal product is exactly the nonsymmetric substitution
\[
	F \odot G = \int^{[n]\in\Fin} Fn \otimes G^n
\]
where $G^n = G \times \cdots\times G$.

As a consequence, we can define clones:
\begin{definition}\label{def_clone}
	A \emph{cartesian operad}, also called a \emph{clone}, is a monoid object in $(\Cat(\Fin,\Set), \odot)$.
\end{definition}
Now, the main theorem in \cite{curiennone} that characterizes operads as the coKleisli category of a certain comonad on $\Dist$ involves the choice of a monad on $\Cat$, among a set of three:
\begin{itemize}
	\item the \emph{free monoidal category} $M(\A)$: given a category $\A$, the free strict monoidal category $M(\A)$ has objects the finite sequences of objects of $\A$, and the hom\hyp{}set between two tuples $\underline A, \underline A'$ is defined to be
	      \[
		      M(\A)(\underline A, \underline A') =
		      \begin{cases}
			      \prod_{i=1}^n \A(A_i, A_i') & \text{if } \ell(\underline{A}) = \ell(\underline{A}') \\
			      \varnothing                 & \text{otherwise}
		      \end{cases}
	      \]
	      where $\ell : M(\A)_o = \coprod_{n\ge 0} \A_o^n \to \N$ sends a tuple $\underline A = (A_1,\dots, A_n)$ to its \emph{length} $n$.
	\item the \emph{free symmetric monoidal category} $S(\A)$: given a category $\A$, the objects of $S(\A)$ are finite tuples $\vec A =(A_1,\dots, A_n)$ of objects of $\A$; morphisms between two tuples $\vec A$ and $\vec b$ exist only if the tuples have the same length; the monoidal product $\uplus$ is given by juxtaposition of tuples. Morphisms $f \uplus g$ are defined accordingly.

	      The category $S(\A)$ thus splits as a disjoint union of categories $S^n(\A)$ whose objects are $n$\hyp{}tuples of objects $\vec A =(A_1,\dots, A_n)$ and the hom\hyp{}objects are
	      \[
		      S^n(\A)(\vec X,\vec Y)  := \bigsqcup_{\sigma\in S_n} \A(X_1,Y_{\sigma 1})\times\dots\times\A(X_n,Y_{\sigma n}).
	      \]
	      The symmetry of the monoidal structure has components $\vec X\uplus \vec Y \to \vec Y \uplus \vec X$, the shuffle permutation swapping the first $m$\hyp{}elements of the first sequence with the last $n$\hyp{}elements of the second. The unit is the empty sequence $()$ and $S^0(\A) = \{()\}$ is the terminal category. The inclusion functor $\A \hookrightarrow S(\A)$ takes an object $x$ to the one\hyp{}element sequence $(x)\in S^1(\A)$.
	\item the \emph{free cartesian category} $K(\A)$, defined adjoining to a category $\A$ all finite products; of course, $K$ is a `sub\hyp{}monad' of $S$, in that every cartesian category is (nonstrictly) symmetric monoidal.
\end{itemize}
each functor $T \in \{M,S,K\}$ induces a monad $\bsP\circ T$, and by self\hyp{}duality a \emph{comonad}, on the bicategory of profunctors, whose coKleisli category is
\begin{itemize}
	\item the category of operads, if $T = M$;
	\item the category of symmetric operads, if $T = S$;
	\item the category of clones, as defined in \ref{def_clone}, if $T = K$.
\end{itemize}
\index{Clone}
The whole proof boils down to the existence of a (pseudo)\footnote{We decide to hide the coherence constraints of the presheaf construction functor $\bsP : A \mapsto \Cat(A^\opp,\Set)$, like \cite{curiennone} does; this is a hairy matter that we have no room or reason to expand as it would deserve; the interested reader can consult \cite{fiore2016relative}. We also set aside the problem given by the fact that strictly speaking, $\bsP$ is not an endofunctor: \cite{fiore2016relative} addresses this matter in the best possible way. See also \cite{altenkirch2010monads, yosegi}.}distributive law (see \ref{distribbio}) $\gamma : T\circ \bsP \To \bsP \circ T$, whose definition and well\hyp{}posedness can be proved by means of co\fshyp{}end calculus.
\begin{lemma}
	Let $T$ be any one of the three monads $M,S,K$; then there is a distributive law $\gamma : T\circ \bsP \To \bsP \circ T$ between $T$ and the presheaf construction $\bsP$.
\end{lemma}
\begin{proof}
	We define $\gamma$ on components, as \cite[§9]{curiennone} does: let $T$ be for example the monad $M$  and let us take $\V=\Set$ (in all other cases, a similar discussion can be carried over almost unchanged). We define the $\A$\hyp{}component of $\gamma$ to be
	\begin{gather}
		\gamma_A : M(\bsP \A) \to \bsP(M\A) \notag \\
		(F_1,\dots,F_n) \longmapsto \int^{\underline A} F_1A_i \times \cdots \times F_nA_n\times M(\A)(\firstblank, \underline{A})
	\end{gather}
	where $\underline A = A_1,\dots, A_n$ (recall that by its very definition $\gamma$ must turn `tuples of presheaves' into `presheaves over tuples' in each degree $n\ge 0$).

	We shall now show that the constraints listed in our \ref{distribbio} hold, thus defining a distributive law.
	\begin{itemize}
		\item the unit constraint of \ref{distribbio} means that the triangle
		      \[
			      \vcenter{\xymatrix{
			      M\bsP\ar[rr]^\gamma && \bsP M \\
			      & \bsP \ar[ur]_{\bsP\eta^{(M)}} \ar[ul]^{\eta^{(M)}\bsP}
			      }}
		      \]
		      commutes; this is evident, as the composition $\gamma\circ \eta^{(M)}\bsP$ sends a presheaf $F \in\bsP \A$ into $\int^A FA \times M(\A)(\firstblank, A)$, \ie into the left Kan extension of $F$ along the unit $\eta^{(M)}_{\A} : \A \to M\A$, \ie exactly into the image of $F$ along $\bsP\eta^{(M)}$.
		\item the second constraint of \ref{distribbio} involves the diagram
		      \[\label{multa}
			      \vcenter{\xymatrix{
			      MM\bsP  \ar[r]^{M\gamma}\ar[d]_{\mu^{(M)}\bsP }& M\bsP M\ar[r]^{\gamma M} & \bsP MM \ar[d]^{\bsP \mu^{(M)}}\\
			      M\bsP \ar[rr]_\gamma && \bsP M
			      }}
		      \]
		      Establishing its commutativity is a bit harder: following \cite[§9]{curiennone} very closely, we will keep track of the position of the diagram we are in, as the coend computation proceeds.

		      First, let's fix a component $\A$: the generic element of upper left corner of \eqref{multa} is thus a tuple of tuples of presheaves, and has the form
		      \[\label{upperleft}
			      \lambda \underline C.
			      \int^{\underline{D_1},\dots, \underline{D_n}} F_1^1 D_1^1 \times \cdots \times F_n^{i_n}D_n^{i_n} \times M(\A)\Big(\underline C, \underline{(D_1,\dots,D_n)}\Big)
		      \]
		      where the tuple of tuples $F_1^1 D_1^1 \times \cdots \times F_n^{i_n}D_n^{i_n}$ has not been flattened, and the tuple of tuples $\underline{\underline{D}} =  \underline{(\underline{D_1},\dots,\underline{D_n})}$ has been flattened to a long tuple $\underline D$ using the multiplication of the monad; this is sent to the lower left corner, to an identically written object where said tuple has been flattened by $\mu^{(M)}$ (formally, we consider the product $F_1^1 D_1^1 \times \cdots \times F_n^{i_n}D_n^{i_n}$ parenthesized in two different, but equivalent, ways). Thus, the object
		      \[
			      \lambda \underline C.
			      \int^{\underline{D_1},\dots, \underline{D_n}} F_1^1 D_1^1 \times \cdots \times F_n^{i_n}D_n^{i_n} \times M(\A)\Big(\underline C, \underline{(D_1,\dots,D_n)}\Big)
			      \notag
			      \vcenter{\begin{tikzpicture}[scale=.4]
					      \draw[gray!30] (0,0) rectangle (2,1);
					      \fill (0,0) circle (4pt);
					      \draw (2,0) circle (4pt);
					      \draw (0,1) circle (4pt);
					      \draw (1,1) circle (4pt);
					      \draw (2,1) circle (4pt);
				      \end{tikzpicture}}\]
		      lies in the lower left corner of the diagram. Applying the distributive law $\gamma$ now yields the functor
		      \begin{multline}
			      \lambda \underline C.\int^{\underline{A_1},\dots,\underline{A_n}}
			      \int^{\underline{D_1},\dots, \underline{D_n}} F_1^1 D_1^1 \times \cdots \times F_n^{i_n}D_n^{i_n} \times \\ \times MM(\A)(\underline{\underline{A}}, \underline{\underline{D}})\times M(\A)(\underline C, (A_1^1,\dots, A_p^{k_p}))
			      \label{lowerright}
		      \end{multline}
		      We shall now transport the object of \eqref{upperleft} to the upper right corner, by application of $\gamma M \circ M\gamma$: the result is equal to
		      \begin{multline*}
			      \lambda (\underline{A_1},\dots, \underline{A_p}).
			      \int^{\underline{B_1},\dots,\underline{B_n}} \int^{\underline{D_1},\dots,\underline{D_n}} F_1^1 D_1^1 \times \cdots \times F_n^{i_n}D_n^{i_n} \times \\
			      \times M(\A)(\underline{B_1},\underline{D_1}) \times \cdots \times M(\A)(\underline{B_n},\underline{D_n}) \times MM(\A)(\underline{\underline{A}},\underline{\underline{B}})
		      \end{multline*}
		      This is in turn equal to
		      \begin{multline*}
			      \lambda (\underline{A_1},\dots, \underline{A_p}).
			      \int^{\underline{D_1},\dots,\underline{D_n}} F_1^1 D_1^1 \times \cdots \times F_n^{i_n}D_n^{i_n} \times \\
			      \times \int^{\underline{B_1},\dots,\underline{B_n}} M(\A)(\underline{B_1},\underline{D_1}) \times \cdots \times M(\A)(\underline{B_n},\underline{D_n}) \times MM(\A)(\underline{\underline{A}},\underline{\underline{B}})
		      \end{multline*}
		      and by application of the ninja Yoneda lemma \ref{ninjayo} this is equal to
		      \[
			      \lambda (\underline{A_1},\dots, \underline{A_p}).\int^{\underline{D_1},\dots,\underline{D_n}}
			      F_1^1 D_1^1 \times \cdots \times F_n^{i_n}D_n^{i_n} \times
			      MM(\A)(\underline{\underline{A}},\underline{\underline{D}})
			      \notag
			      \vcenter{\begin{tikzpicture}[scale=.4]
					      \draw[gray!30] (0,0) rectangle (2,1);
					      \draw (0,0) circle (4pt);
					      \draw (2,0) circle (4pt);
					      \draw (0,1) circle (4pt);
					      \draw (1,1) circle (4pt);
					      \fill (2,1) circle (4pt);
				      \end{tikzpicture}}
		      \]
		      Flattening, \ie applying the functor $\bsP \mu^{(M)}$, falls in the lower right corner, and results in the same functor as in \eqref{lowerright}.
	\end{itemize}
	This concludes the proof.
\end{proof}
\begin{theorem}\label{curien_main}
	The monad $\bsP \circ T$ can be turned into a comonad under the equivalence $\Dist\cong \Dist^\opp$, and the coKleisli category of such a comonad corresponds to the category of operads if $T=M$, the category of symmetric operads if $T=S$, and the category of clones if $T=K$.
\end{theorem}
This means, more in detail, that
\begin{itemize}
	\item a functor $F : \C \to \V$ is a free $(\bsP\circ M)$\hyp{}coalgebra if and only if it is an operad (and similarly for symmetric operads, and clones if $\V$ is cartesian);
	\item the composition in the coKleisli category of $\bsP\circ T$ coincides with the substitution product for $T=M$, with the symmetric substitution for $T=S$, and with the cartesian substitution if $T=K$.
\end{itemize}
\begin{exercises}
\item \label{operad_1} Let $\A$ be an ordinary (\ie, $\Set$-enriched) category, and let $\V$ be a cosmos; define the \emph{free} $\V$-category $\overline{\A}^{\V}$ on $\A$ as the $\V$-category with the same objects of $\A$, and where
\[\notag \overline{\A}^{\V}(A,A') := \A(A,A')\cdot I\]
is a $\A(A,A')$-fold tensor of the monoidal unit $I$ of $\V$ (\ie, a coproduct of as many copies of $I$ as there are elements in $\A(A,A')$).
\begin{itemize}
	\item Prove that there is an isomorphism of categories $\overline{(\A\times \B)}^{\V}\cong \overline{\A}^{\V} \times \overline{\B}^{\V}$; deduce that if $\A$ is a monoidal category, $\overline{\A}^{\V}$ is a monoidal $\V$-category.
	\item Does the `free' $\V$-category $\overline{\A}^{\V}$ deserve its name? Is is true that $\V$-enriched functors
	      \[\notag F : \overline{\A}^{\V} \to \C\]
	      to a $\V$-category $\C$ correspond to unenriched functors
	      \[\notag  F_0 : \A \to |\C| \]
	      to the underlying ordinary category of $\C$?
\end{itemize}
\item \label{ex:convo_al_lan} Show that the convolution product on $[\C, \V]$ results as the following left Kan extension:
\[\notag
	\xymatrix{
	\C \times \C \ar[d]_\odot \ar[r]^{F\times G} & \V \times \V \ar[r]^(.6){\otimes_\V} & \V \\
	\C \ar@{.>}@/_1pc/[urr]
	}
\]
\item Show that the Isbell duality
\[\notag\vcenter{\xymatrix@C=1.5cm{
	\VCat(\C,\V)^\opp\ar@{}[r]|\perp \ar@<5pt>[r]^{\isbO} &\ar@<5pt>[l]^{\Spec} \VCat(\C^\opp, \V)
	}}\notag\]
of \ref{isbella-duella} is a pair of adjoint strong monoidal functors, when (the category $\C$ is monoidal and) their domains are endowed with the convolution product.
\item Let $\Set_*$ be the category of pointed sets; such category is monoidal with respect to the \emph{smash product}, where a pointed set $(X,x_0)$ and a pointed set $(Y,y_0)$ are smashed in the set $X\land Y$ defined as the pushout
\[\notag
	\xymatrix{
		1 \ar[r]^{x_0} \ar[d]_{y_0}& X \ar[d]\\
		Y \ar[r]& X\land Y
	}
\]
Let $\Seq$ be the discrete category whose objects are natural numbers $\{0,1,2,\dots\}$. Endow the category $\Cat(\Seq,\Set_*)$ with the convolution product, and show that
\[\notag
	F \ast G = \lambda n . \coprod_{i+j=n} Fi \land Gj
\]
\item Describe the convolution product in the category $\Cat(\N, \Set_*)$, where $\N = \{0<1<2<\dots\}$, and in the category $\Cat(\cN, \Set_*)$, where $\cN$ is the monoid $(\N, \cdot)$ regarded as a category with a single object (the category of pointed sets is always endowed with the smash product). Compare the latter with the convolution product on the category $\Cat(\cN, \Set)$ (cartesian monoidal structure).
\item \label{opd_adjoints} Define two functors
\begin{gather*}
	\Phi : \VCat(\P,\V)\to \V \qquad \text{(evaluation at $0$)}\\
	\Psi : \V \to \VCat(\P,\V)\qquad \text{(the left adjoint to } \Phi)
\end{gather*}
Prove that $\Psi(A\oplus B)\cong \Psi A\otimes \Psi B$ and $\Phi\circ\Psi\cong 1$; finally, for every object $V\in \V$, if $\yon V$ is the representable functor on $V$, then $\VCat(\P,\V)(V\ast F,G)\cong \V(V,[F,G])$, where
\[\notag
	[F,G] := \int_n \V(Fn, Gn)
\]
has been defined in \ref{cosmuclosed}.
\item Fill in the details of the proof that \eqref{opd_to_monad} defines a monad, by sowing that the associativity and unitality of an operad yield associativity and unitality of the monad $M_T$.
\item   Let $\A, \B$ be small categories; an \emph{$S$\hyp{}profunctor} $ P : \A \overset{S}{\pto} \B$ is a profunctor $ P : \A \pto S\B$, i.e. a functor $ P : S\B^\opp \times \A \to \Set$, where $S$ is the free symmetric monoidal category functor; given $ F : \cX \overset{S}{\pto} \cY$ and $ G : \cY \overset{S}{\pto} \cZ$, we define the \emph{composition} $G\circ F : \cX \overset{S}{\pto} \cZ$ as the coend
\[\notag
	(G\circ F)[\underline{Z}; X] = \int^{\vec Y} G^e[\underline{Z}; \vec Y]\otimes F[\vec Y; X]
\]
where
\[\notag
	G^e[\underline{Z}; \vec Y] := \int^{\underline{Z}_1,\dots,\underline{Z}_n}   G[\underline{Z}_1,Y_1]\times \dots G[\underline{Z}_n,Y_n]\times S(\cZ)(\underline Y, \underline{Z}_1\uplus \dots\uplus \underline{Z}_n)
\]
where $\uplus$ is concatenation.

Show that this is indeed an associative composition rule for a bicategory of $S$\hyp{}profunctors; find the identity 1\hyp{}cell of an object $\A$.
\item Fill in the details in the proof of \ref{curien_main}, and in particular prove that the composition in the coKleisli category of $\bsP\circ T$ coincides with substitution products.
\begin{itemize}
	\item Show that the monad $\bsP \circ T$, for $T \in \{M,S,K\}$ induces a monad $\hat T$ on $\Dist$, and thus a comonad $W$ by posing $W\A := T(\A^\opp)^\opp$. (Hint: use the self\hyp{}duality of \ref{dualiseur}).
	\item Show that $\hat T$ has the following expression on 1\hyp{}cells: given $\proP : \A \pto \B$,
	      \[\notag
		      \hat T \proP : (\underline A,\underline B') \mapsto \int^{\underline B} \proP(A_1,B_1) \times \cdots\times \proP(A_n, B_n) \times T(\B)(\underline B, \underline B')
	      \]
	\item Show that given a 2\hyp{}comonad $W$ on $\Dist$, if we denote $\Dist_W$ its coKleisli bicategory, the coKleisli composition is defined as soon as every map has a \emph{coKleisli lifting}: the composition $g\bullet f = g \circ Wf \circ \sigma_C$ of $f : TC \to C'$ and $g : TC' \to C''$ equals the composition $g\circ R(f)$ in the diagram
	      \[\notag
		      \vcenter{\xymatrix{
		      \drtwocell<\omit>{}& TC' \ar[r]^g \ar[d]^{\epsilon_{C'}}& C''\\
		      TC\ar@/^1pc/[ur]^{Rf}\ar[r]_f & C'
		      }}
	      \]
	      where $Rf$ is the right lifting of $f$ along $\epsilon_{C'} : TC' \to C'$.
	\item Use this, and Exercise \ref{riftuzzo}, to prove that the coKleisli composition coincides with the substitution product.
\end{itemize}
\item \label{pull_is_mona} Let
\[\notag\vcenter{\xymatrix{
			\A\ar[d]_{U'} \ar[r]^{F'}& \B\ar[d]^U \\
			\C \ar[r]_F & \D
		}}\]
be a strict pullback square of categories and functors, and assume that $F$ is fully faithful; prove that if $U$ is monadic, and $U'$ has a left adjoint, then $U'$ is monadic too.
\item Show that if $F : \E \to \A$ is a discrete fibration (see \ref{def:dfib}) which is moreover a strong monoidal functor, then the associated presheaf $\hat F : \A^\opp\to \Set$ is a monoid with respect to the Day convolution product (see \ref{day}) on $\Cat(\A^\opp,\Set)$. Show that this sets up a bijection under the equivalence of \ref{thm:equconfib}.
\item \awful Let $\Cat(\Set,\Set)_s$ be the category of \emph{small} functors\footnote{An endofunctor of $\Set$ is called small if it results as the left Kan extension of a functor $\A\to\Set$, where $\A\subseteq\Set$ is a small subcategory, along said inclusion. This restriction is needed in order for $\Cat(\Set,\Set)$ to (exist and to) form a locally small category.} $F : \Set \to\Set$ and let $F,G$ be two comonads; show that the Day convolution $F * G$ is itself a comonad. (hint: there is a neat argument that involves the theory of \emph{duoidal} categories, see \cite[§8.1]{garner2016commutativity})
\item \awful Is it possible to dualise the Day convolution of \ref{day} to define a \emph{Day involution} operation on $\Cat(\A^\opp,\Set)$? For a suitable operation $A,B\mapsto \langle A,B\rangle$ and presheaves $F,G : \A^\opp \to \Set$ we shall define
\[\notag
	\lceil F,F\rfloor := \int_A \A(\firstblank, \langle A,B\rangle)\pitchfork [FA,GB]
\]
How does this operation behave?
\item \awful Prove that the operad associated to a monad on $\Set$ (see \ref{eo}.\ref{eo:due}) really is an operad, by showing the commutativity of the following diagrams:
\[\notag
	\vcenter{\xymatrix{
			Tm \times T1 \ar[r]^-\gamma \times\dots\times T1 & Tm & T1 \times T1 \ar[r]^\gamma  & T1 \\
			Tm \times 1 \times\dots\times 1 \ar[r]_-\sim \ar[u]^{Tm\times \eta\times\dots\times \eta} & Tm\ar@{=}[u] & 1\times T1 \ar[u]^{\eta\times T1} \ar[r]_-\sim & T1\ar@{=}[u]
		}} \\
\]
where the lower horizontal isomorphism comes from the unitor of the cartesian monoidal structure on $\Set$. (Hint: you might want to recall the commutativities given by naturality of the monad structure, \ref{def:monad}, as well as naturality of the counit $\epsilon$ of the cartesian closed structure of $\Set$, \ref{ex:adjoints}.\ref{ad:cin}, as well as the fact that $\epsilon$ is a cowedge.)

Good luck proving associativity: the axiom is stated as follows. Given a positive integer $m$, positive integers $p_1,\dots,p_m$, and a matrix of positive integers $\bsmat q_{1,1} &\dots & q_{1,p_1} \\ \vdots & \ddots & \vdots \\ q_{m,1} &\dots & q_{m,p_m} \esmat$, consider the diagram
\[\notag
	\xymatrix@d@R=-.75cm{
	Tm \times T\vec p \times T\vec q_1 \times \cdots\times T\vec q_m \ar[r]^-{\id\times\gamma_{p_1,\vec q_1}\times\cdots\times \gamma_{p_m,\vec q_m}}\ar[d]_-{\gamma_{m,\vec p}\times \id} & Tm \times T\left(\sum_{j=1}^{p_1} q_{1,j}\right)\times\cdots\times T\left(\sum_{j=1}^{p_m} q_{m,j}\right)\ar[d]^-{\gamma_{m, \vec{\vec q}}}\\
	T\left(\sum_{i=1}^m p_i\right) \ar[r]_{\gamma_{\Sigma p_i,\vec q}} & T\left(\sum_{i=1}^m\sum_{j=1}^{p_i} q_{i,j}\right)
	}
\]
(there is an implicit identification between reshuffled copies of the factors in the upper left cartesian product, and we shorten the product $Tq_{j,1}\times\cdots\times Tq_{j,p_j}$ as $T\vec q_j$.) The associativity axiom for $\gamma_{\firstblank, \secondblank}$ amounts to the request that this diagram commutes.
\item \label{the_soa_at_work} (the small object argument at work) Let $\lambda : T \To S$ be a morphism $S \to T$ between two finitary monads (see \ref{accepre}) on the category of sets; let $\lambda^\sharp : \Alg(S) \to \Alg(T)$ be the functor defined sending a $S$-algebra $(A,a:SA\to A)$ into the $T$-algebra $(A, a\circ \lambda_A)$.

Show that $\lambda^\sharp$ has a left adjoint $\lambda_!$:
\begin{itemize}
	\item Define the arrow $(m_0, t_0)$ as the pushout
	      \[\notag
		      \vcenter{\xymatrix{
				      TA \ar[r]^{\lambda_A}\ar[d]_a & SA \ar[d]^{m_0}\\
				      A \ar[r]_{t_0}& P_0
			      }}
	      \] and by induction, define $(m_{i+1}, t_{i+1})$ to be the pushout of the pair of arrows $St_i$ and $m_i$:
	      \[\notag\vcenter{\xymatrix{
				      TA \po\ar[r]^{\lambda_A}\ar[d]_a & SA \po\ar[d]^{m_0}\ar[r]^{St_0} & SP_0\po\ar[d]^{m_1}\ar[r]^{St_1} & SP_1 \ar[d]\ar[r]^{St_2}& \dots\ar[r]  & Q\ar[d]^{\alpha=m_\infty}\\
				      A \ar[r]_{t_0}& P_0\ar[r]_{t_1} & P_1\ar[r]_{t_2} & P_2\ar[r]_{t_3} & \dots\ar[r] & P_\infty
			      }}\]
	      and let $\alpha : Q \to P_\infty$ be its colimit.
	\item Show that $SP_\infty\cong Q$; show that $\alpha$ endows $P_\infty$ with the structure of a $S$-algebra.
	\item Show that $\lambda_!$, defined sending $(A,a)$ to $(SP_\infty,\alpha_a)$ is a functor. Show that there is an adjunction $\lambda_!\dashv \lambda^\sharp$.
\end{itemize}
\item \upeyes Read \cite{gen-operads} and take it as an outstandingly useful mental exercise to generalise the topics of this section.
\end{exercises}

\def\p{\mathrm{p}}
\chapter{Higher dimensional co\fshyp{}ends}\label{sec:higher}
\begin{abstract}
	The present chapter studies co\fshyp{}end calculus in higher category theory; the basic theory of 2\hyp{}dimensional co\fshyp{}end calculus is introduced and explored in fair completeness, as well as `homotopy coherent' versions of co\fshyp{}ends (in model categories), $(\infty,1)$\hyp{}categorical co\fshyp{}ends, and co\fshyp{}ends inside a Groth(endieck) derivator \cite{groth2013derivators}. The focus is on showing how co\fshyp{}end calculus is a paradigm that can be exported in whatever formal context to do category theory, instead than a mere set of theorems.
\end{abstract}
\epigraph{I Pitagorici raccontano che l'uso di insegnare geometria iniziò così: a uno dei membri della setta, perdute le proprie sostanze in un disastro, fu fatta la concessione di trarre guadagno dall'insegnamento della geometria.}{M. Timpanaro Cardini --- \emph{Vite dei Pitagorici}}
We progressively raise the dimension we work in, starting with the theory of co\fshyp{}ends in 2\hyp{}categories. We begin the next section by recalling the bare minimum of 2\hyp{}dimensional category theory needed to appreciate the discussion; the definition of bicategory and 2\hyp{}category (=strict bicategory) is understood as in \ref{bicat}.
\begin{notat}
	We freely employ the theory sketched in §\ref{2categories} and in particular \ref{pseudocolax}: as always when dealing with higher dimensional cells and their compositions, there are several `flavours' in which one can weaken strict commutativity. Besides this strictness (where every diagram commutes with an implicit identity 2\hyp{}cell filling it), there is a notion of \emph{strong} commutativity and universality, where filling 2\hyp{}cells are requested to be invertible, and \emph{lax} commutativity/universality, where 2\hyp{}cells are possibly non\hyp{}invertible.
\end{notat}
\section{2\hyp{}dimensional coends}
A \emph{lax functor} $F : \C \laxto \D$ between two 2\hyp{}categories $\C,\D$ behaves like a functor, up to the fact that there is a non\hyp{}invertible 2\hyp{}cell linking the composition of the images $Ff, Fg$ and the image of a composition $F(gf)$ (see \ref{colaxe} for the precise definition).

\index{_aaa_Catl@$\Cat_l$}
A lax natural transformation $\alpha : F \To_l G$ is a family of 1\hyp{}cells $\alpha_C : FC \to GC$ coupled with 2\hyp{}cells $\alpha_f$, filling the diagrams
\[
	\vcenter{\xymatrix{
	FC \ar[d]_{\alpha_C}\ar[r]^{Ff}\drtwocell<\omit>{\alpha_f} & FC' \ar[d]^{\alpha_{C'}} \\
	GC \ar[r]_{Gf} & GC'.
	}}
\]
these 2\hyp{}cells are subject to suitable coherence conditions linking $\alpha_f,\alpha_g$ and $\alpha_{gf}$ for composable 1\hyp{}cells $C \xto{f} C' \xto{g} C''$ and characterizing $\alpha_{\id_C}$. 

In the setting of 2\hyp{}categories, co\fshyp{}end calculus has a rather natural interpretation in terms of enriched category theory;\footnote{See \cite{kelly1974review,2catlimits}, for invaluably complete surveys on the matter; the reader should however be aware that 2\hyp{}category theory is \emph{not} completely subsumed by the theory of $\Cat$\hyp{}enriched categories.} even though we strive for clarity, the following discussion has little hope to be a self\hyp{}contained exposition, and instead heavily relies on the existing literature (see for example \cite{kelly1982basic,dubuc1970kan}).

\smallskip
The definition of \emph{lax co\fshyp{}end} of a 2\hyp{}functor $S$ is given in terms of the notion of a \emph{lax wedge} $\omega : B\din S$, and it is the most general and less symmetric one can give: asking for the diagrams to be filled by isomorphisms, it specialises to the notion of \emph{strong} co\fshyp{}end, and strict co\fshyp{}end (of course, we can also define \emph{oplax} co\fshyp{}ends by reversing all 2\hyp{}cells); the present subsection is designed in such a way to reduce to the strong and strict cases as particular examples.

This said, we shall warn the reader that 2\hyp{}category theory exists in many dialects and follows slightly different notational conventions (minor differences that appear innocuous to the expert, but they can become annoying when studying a part of category theory still lacking a truly comprehensive monograph).

We try to follow an auto\hyp{}explicative notation based on our appendix A and on canonical references like \cite{Kelly2005b,2catlimits}, but we feel free to diverge from it from time to time.

The material on co\fshyp{}lax co\fshyp{}ends in the rest of the present section comes in its entirety from \cite{bozapalides1980some}: the original paper, as well as Bozapalides' PhD thesis \cite{bozapalides1976theories}, are very difficult to retrieve, and together they provide useful references and a starting point to develop 2\hyp{}dimensional coend calculus. We hope this survey can provide a useful and more accessible reference in the future, for a piece of elegant mathematics that is sorely lacking from the mainstream literature.

As usual in our Appendix on basic category theory, we denote $\A_o$ the class of objects of a category, or 2\hyp{}category $\A$.
\begin{definition}[Lax wedge]\label{laxwedge}\index{Co/wedge!Lax ---}
	Let
	\[S : \A^\opp\times \A \to \B\]
	be a strict 2\hyp{}functor between strict 2\hyp{}categories $\A,\B$. A \emph{lax wedge} $\omega$ for $S$ consists of
	\begin{itemize}
		\item a triple $\{B, \omega_o, \omega_h\}$, where $B\in \B_o$ is an object (called the \emph{tip} of the wedge);
		\item collections of 1\hyp{}cells $\big\{ \omega_A : B\to S(A,A)\big\}$, one for each $A\in \A_o$;
		\item 2\hyp{}cells $\big\{ \omega_f : S(A, f)\circ \omega_{A} \To S(f, A')\circ \omega_{A'} \big\}$, in a diagram
	\end{itemize}
	\[\label{prima}
		\vcenter{\xymatrix@C=1.5cm{
		B\drtwocell<\omit>{\quad\omega_f} \ar[r]^{\omega_A} \ar[d]_{\omega_{A'}} & S(A,A) \ar[d]^{S(A,f)}\\
		S(A',A')\ar[r]_{S(f,A')} & S(A,A') \\
		}}
	\]
	These data must fit together in such a way that the coherence axioms listed below, expressed by the commutativity of the following diagrams of 2\hyp{}cells, are satisfied:
	\begin{enumtag}{lw}
		\item The diagram of 2\hyp{}cells having faces
		\[
			\vcenter{\xymatrix@C=1.5cm{
			B \drtwocell<\omit>{<-1>\quad\omega_f}\ar[r]^{\omega_A}\ar[d]_{\omega_{A'}}& S(A,A)\ar[d]^{S(A,f)} & B \drtwocell<\omit>{\kern-12mm\omega_{f'}}\ar[r]^{\omega_A}\ar[d]_{\omega_{A'}}& S(A,A)\dtwocell^{\qquad S(A,f)}_{S(A,f')\qquad}{} \\
			S(A',A') \rtwocell^{S(f',A')}_{S(f,A)}{} & S(A,A') & S(A',A')\ar[r]_{S(f,A')} & S(A,A') \\
			}}
		\]
		is commutative for any $\lambda : f\To f'$, \ie the equation
		\[\notag
			\omega_{f'}\circ (S(A, \lambda)*\omega_A) = (S(\lambda, A') * \omega_{A'})\circ \omega_f
		\]
		holds.
		\item For each pair $A\xto{f}A'\xto{f'}A''$ of composable arrows in $\A$, the diagram of 2\hyp{}cells
		\[
			\vcenter{\xymatrix@R=4mm@C=4mm{
			& B \ar@{}[dddl]|{\Nwarrow\omega_{f'}}
			\ar@{}[dddr]|{\Swarrow\omega_f}\ar[dr]^{\omega_A}\ar[dl]_{\omega_{A''}}\ar[dd]_{\omega_{A'}} & && B\ddtwocell<\omit>{\omega_{f'f}}\ar[dr]^{\omega_A}\ar[dl]_{\omega_{A''}} & \\
			S(A'',A'') \ar[dd]_{S(f',A'')} && S(A,A)\ar@{}[ddr]|=\ar[dd] & S(A'',A'')\ar[dddr]|(.6){S(f'f,A'')}\ar[dd]  && S(A,A) \ar[dddl]_(.3){S(A,f'f)}\ar[dd]^{S(A,f)}\\
			& S(A',A')\ar@{}[dd]|= \ar[dr]\ar[dl]& & & \\
			S(A',A'')\ar[dr]_{S(f,A'')} && S(A,A') \ar[dl]^{S(A,f')} & S(A',A'')\ar[dr]_{S(f,A'')} && S(A,A') \ar[dl]^{S(A,f')}\\
			& S(A,A'') & & & S(A,A'') &
			}}
		\]
		is commutative, \ie the equation
		\[\notag
			(S(f, A'') * \omega_{f'})\circ (S(A, f') * \omega_f) = \omega_{f'f}
		\]
		holds.
		\item For each object $A\in\A$, $\omega_{\id_A} = \id_{\omega_A}$.
	\end{enumtag}
\end{definition}
\begin{notat}
	We use the compact notation $\omega : B\din S$ for a lax wedge with domain constant at the object $B$; this is evidently reminiscent of our \ref{extranatural}.
\end{notat}
\begin{definition}[Modification]\label{modification}\index{Modification}
	A \emph{modification} $\Theta : \omega\Rrightarrow\sigma$ between two lax wedges $\omega,\sigma : B\din S$ for $S : \A^\opp\times \A \to \B$ consists of a collection of 2\hyp{}cells $\Theta_A : \omega_A\To \sigma_A$ indexed by the objects of $\A$ such that the diagram of 2\hyp{}cells
	\[
		\vcenter
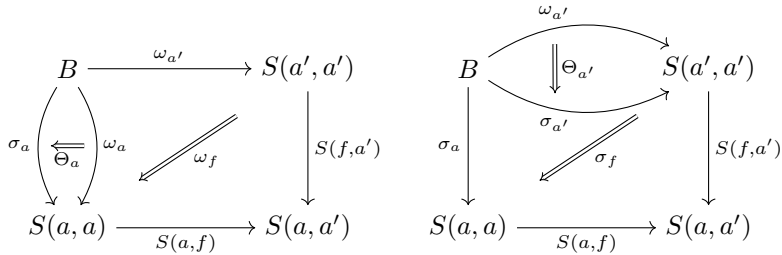
{
		\xymatrix@R=1.4cm@C=1.4cm{
		B \drtwocell<\omit>{\quad\omega_f}\ar[r]^{\omega_{A'}} \dtwocell^{\omega_A}_{\sigma_A}{\Theta_A}& S(A',A')\ar[d]^{S(f,A')} & B\drtwocell<\omit>{\quad\sigma_f} \rtwocell^{\omega_{A'}}_{\sigma_{A'}}{\quad\Theta_{A'}}\ar[d]_{\sigma_A} & S(A',A') \ar[d]^{S(f,A')}\\
		S(A,A) \ar[r]_{S(A,f)}& S(A,A') & S(A,A) \ar[r]_{S(A,f)}& S(A,A')
		}}
	\]
	is commutative, \ie
	\[\notag
		(S(A,f) * \Theta_A)\circ \omega_f = \sigma_f \circ (S(f,A') * \omega_{A'})
	\]
\end{definition}
The above definition of a modification is modeled on the definition of modification between (lax) natural transformations of functors.
\begin{remark}
	There is a more general definition for a modification $\Theta : \omega\Rrightarrow\sigma$ between lax wedges having different domains, say $\{B, \omega\}$ and $\{B', \sigma\}$: it consists of a morphism $\varphi : B\to B'$ and a 2\hyp{}cell $\lambda_A : \sigma_A \circ\varphi \To \omega_A$ such that
	\[
		(\sigma_f * \varphi)\circ (S(A,f) * m_A) = (S(f, A') * m_{A'})\circ \omega_f
	\]
	(draw the corresponding diagram of 2\hyp{}cells). We are not interested in this alternative definition, and this will not be investigated further. We thus take \ref{modification} as our working definition without further mention.

	Such definition entails that the set $\mathrm{LWd}(B, S)$ of lax wedges $B\din S$ is a category having morphisms precisely the modifications $\Theta : \omega\Rrightarrow \sigma$ with the same tip, and the correspondence $\beta_S : B\mapsto \text{LWd}(B,S)$ is functorial. The definition of \emph{lax end} for $S$ relies on the representability of this 2\hyp{}functor.
\end{remark}
\begin{definition}[Lax end of $S$]\index{Co/end!lax ---}
	Let $S : \A^\opp\times \A \to \B$ be a 2\hyp{}functor; an object of $\B$ is called the \emph{lax end} of $S$, and denoted $\twoint_A S(A,A)$, if there is a terminal lax wedge $\twoint_A S(A,A)\din  S$, \ie such that for any other lax wedge $\sigma : B'\to S$ there exists a unique 1\hyp{}cell $u : B'\to B$ between the tips of the wedges such that
	\[
		\omega_A \circ u = \sigma_A,\qquad
		\omega_f * u = \sigma_f.
	\]
	(This pair of equations is conveniently depicted by the diagram of 2\hyp{}cells
	\[\label{ultima}
		\vcenter
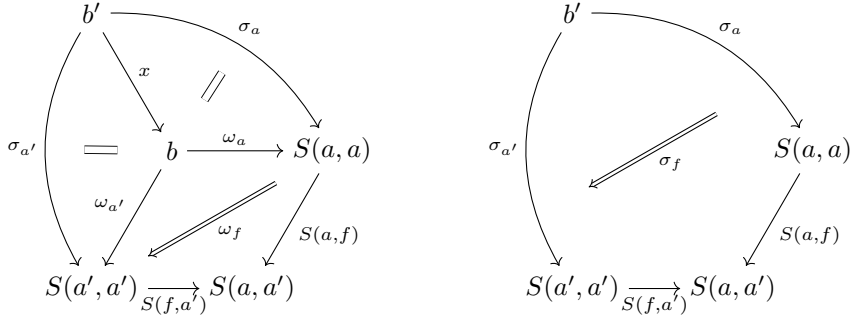
{\xymatrix{
		B' \drruppertwocell^{\sigma_A}_{}{<1>\quad\id}
		\ar[dr]_u
		\ddrlowertwocell^{}_{\sigma_{A'}}{<-1>\id}
		&& & B'\ddrrtwocell<\omit>{\quad\sigma_f} \ar@/^1.5pc/[drr]^{\sigma_A}  \ar@/_1.5pc/[ddr]_{\sigma_{A'}}\\
		& B\drtwocell<\omit>{\quad\omega_f} \ar[r]^{\omega_A}\ar[d]_{\omega_{A'}} & S(A,A)\ar@{}[drr]|(.3)=\ar[d]^{S(A,f)}  & && S(A,A)\ar[d]^{S(A,f)}\\
		& S(A',A') \ar[r]_{S(f,A')} & S(A,A') && S(A',A')\ar[r]_{S(f,A')} & S(A,A')
		}}
	\]
	that we require to commute). Moreover, every modification $\Theta : \sigma\Rrightarrow \sigma'$ induces a unique 2\hyp{}cell $\lambda : u\To u'$ ($u'$ is the arrow induced by $\sigma'$) in such a way that
	$\lambda * \omega_A = \Theta_A$.

	This equivalence sets up an isomorphism of categories between lax wedges with tip $B$ and the category $\B(B, \twoint_A S(A,A))$ (if needed, this can be relaxed to an \emph{equivalence} of categories).
\end{definition}
\begin{remark}\label{laxendnotat}
	We denote the lax end of $S$ as a `squared integral'
	\[
		\twoint_A S(A,A).
	\]
	This notation has a meaning: in the world of $n$\hyp{}categories, the $n$\hyp{}co\fshyp{}end operation should be depicted by an integral symbol (with suitable super- or subscripts) overlapped by an $2^n$\hyp{}agon; in this way, the polygon for a 2\hyp{}end has the correct number of edges since it is denoted as a `square\hyp{}over\hyp{}integral' symbol $\twoint$, and the circle being a polygon with an infinite number of sides, the notation is consistent for $\infty$\hyp{}co\fshyp{}ends (see \ref{infend}, where we denoted the co\fshyp{}end as an $\infint$ symbol).
\end{remark}
\begin{remark}
	Reversing the direction of 1\hyp{}cells yields the notion of \emph{lax coend}; but be careful! There is an additional dimension that can be reversed, \ie the direction of 2\hyp{}cells. Doing so, we obtain the notion of \emph{oplax} end and coend. We will rarely need to invoke oplax cells; of course every statement involving a lax widget can be properly dualised to get an oplax one.

	In compliance with the above \ref{laxendnotat}, we denote the lax end of $S$ as a squared integral with a superscript,
	\[
		\twoint^A S(A,A).
	\]
\end{remark}
Imposing stronger conditions on the 2\hyp{}cells filling diagrams \eqref{prima}--\eqref{ultima} we obtain the notion of strong coend (if each component $\omega_f$ of a wedge is an isomorphism for every arrow $f : A\to A'$) and strict coend (if each $\omega_f$ is the identity). Of course, a similar terminology applies to those cases, as well as the co\fshyp{}end calculus we are going to develop in the next section.
\subsection{Lax co\fshyp{}end calculus}
Co/end calculus and its rules remain true and expressive in the lax setting: several results proved in the previous chapters can be suitably `laxified', justifying the intuition of lax co\fshyp{}ends as the right 2\hyp{}categorical generalisation of strict co\fshyp{}ends.

We collect the most notable examples of this phenomenon in the rest of the section; the content of \ref{laxyonedaninja} surely deserves a special mention, as well as other remarks chosen to convey a sense of continuity and analogy. In \ref{laxyonedaninja} we prove that the lax counterpart of the ninja Yoneda lemma \ref{ninjayo} provides a reflection (using coends) and a reflection (using ends) of the category of strong presheaves into the category of lax presheaves.
\begin{example}\label{commaobj}
	The \emph{comma objects} $(f/g)$ \cite{Gray} of a 2\hyp{}category $\B$ can be identified with the lax end of functors $[1]^\opp\times[1]\to \B$ choosing the two 1\hyp{}cells $f,g$; this is a perfect analogy with Exercise \ref{ispull}, also in view of the characterisation of the comma object $(f/g)$ as a lax limit.
\end{example}
\begin{example}\label{laxnat}\index{Co/end!natural transformations as ---s}\index{Natural transformation!---s as ends}
	If $F,G : \A\to \B$ are 2\hyp{}functors, then the lax end of the functor
	\[
		\B(F,G) : \A^\opp\times \A\to \Cat
	\]
	characterises \emph{lax natural transformations between lax functors $F,G : \A\to \B$}:
	\[\label{this_laxnat}
		\twoint_A \B(FA,GA)\cong \twoCat_l(\A, \B)(F,G)
	\]
	(see \cite{Gray} for more information).
\end{example}
\begin{proof}
	We abbreviate $ \twoCat_l(\A, \B)(F,G)$ as $\hom(F,G)$. The reader will notice that the argument is fairly elementary and echoes our proof of \ref{naturalu}.

	A lax wedge $\tau$ for the 2\hyp{}functor $(A,A')\mapsto \hom(FA, GA')$ amounts to a diagram
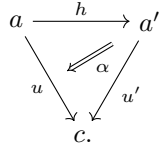
	\[
		\vcenter{\xymatrix{
		\E\ar[r]^{\tau_A}\ar[d]_{\tau_{A'}}\drtwocell<\omit>{\,\,\,\tau_f} & \hom(FA,GA) \ar[d]^{Gf\circ\firstblank} \\
		\hom(FA', GA') \ar[r]_{\firstblank\circ Ff}& \hom(FA, GA')
		}}
	\]
	filled by a 2\hyp{}cell $\tau_f : Gf \circ \tau_A \To \tau_{A'}\circ Ff$. Similarly to what happens in \ref{naturalu}, each of the functors $\tau_A : \E \to \hom(FA, GA)$ sends an object $E\in \E$ into an object $\tau_A(E) : FA \to GA$ such that
	\[
		G(f)\circ \tau_A(E) \overset{\tau_{f,E}}\Longrightarrow \tau_{A'}(E)\circ F(f)
	\]
	is a 2\hyp{}cell filling the square above. This is precisely the lax naturality condition needed to show that the correspondence $E\mapsto \{\tau_AE\}_{A\in \A}$ factors through $\twoCat_l(\A, \B)(F,G)$ (a moderate amount of work is needed to check that the correct coherence properties hold for $\tau_{f,E}$, but we leave such minor chore to the willing reader).
\end{proof}
Lax natural transformations $\eta : F\laxto G$, described as the lax end above, can also be characterised as \emph{lax limits} in the enriched sense: this motivates the search for a more general description of lax co\fshyp{}ends as lax co\fshyp{}limits, analogue to our \ref{coends_as_colims}, where instead of strict co\fshyp{}equalisers we use \cite{2catlimits}'s notion of \emph{co\fshyp{}inserter} (see \ref{inserters}).
\begin{definition}[Lax co\fshyp{}limit]\label{laxlimo}
	Let $F : \A \to \B$ be a 2\hyp{}functor; a \emph{lax limit} for $F$ consists of a family of 1\hyp{}cells $p_A : \llim F \to FA$ from an object $\llim F$, together with 2\hyp{}cells $\pi_f : Ff\circ p_A \To p_{A'}$ for every $f : A \to A'$, such that the pair of families $(p_A, \pi_f)$ is terminal with this property.

	This means that every other family of 1\hyp{}cells $u_A : X \to FA$ and 2\hyp{}cells $\upsilon_f : Ff\circ u_A \To u_{A'}$ factors uniquely through a 1\hyp{}cell $\bar u : X \to \llim F$, as in the diagram
	\[
		\vcenter{\xymatrix{
		&& FA \ar[dd]^{Ff}\\
		X \ar@/^1.5pc/[urr]^{u_A}\ar@/_1.5pc/[drr]_{u_{A'}}\ar@{.>}[r]^{\bar u}& \llim F \ar@{}[r]|(.6){\Swarrow \pi_f}\ar[ur]^{p_A}\ar[dr]_{p_{A'}}& \\
		&& FA'
		}} \quad = \upsilon_f
	\]
	(this takes into account the 1\hyp{}dimensional universal property of $\llim F$) and moreover, given a 2\hyp{}cell $\mu : p_A \bar u \To p_A \bar v$, there is a unique 2\hyp{}cell $\bar \mu$ such that
	\[
		\vcenter{\xymatrix@R=5mm@C=5mm{
		& \llim F \ar[dr]^{p_A}& \\
		X\rrtwocell<\omit>{\mu} \ar[ur]^{\bar u}\ar[dr]_{\bar v}&& FA \\
		& \llim F\ar[ur]_{p_A}
		}} \quad = p_A * \bar \mu
	\]
\end{definition}
Reversing the direction of 1\hyp{}cells yields the notion of \emph{lax colimit}; but be careful! There is an additional dimension that can be reversed, \ie the direction of 2\hyp{}cells. Doing so, we obtain the notion of \emph{oplax} limit and colimit. Again, we employ without further mention each of these names whenever needed.
\begin{proposition}[Co/ends of mute functors]
	Suppose that the 2\hyp{}functor $S : \A^\opp\times\A \to \B$ is mute in the contravariant variable, \ie that there is a factorisation $S = S'\circ p : \A^\opp\times\A \xto{p} \A \xto{S'} \B$; then the lax end of $S$ exists if and only if the lax limit of $S$ does; moreover, the two objects are canonically isomorphic, sharing the same universal property:
	\[
		\twoint_A S(A,A)\cong \llim S'
	\]
	This means that every lax co\fshyp{}limit is a degenerate form of lax co\fshyp{}end.
\end{proposition}
\begin{example}
	As a particularly simple example of this, if $\A$ is locally discrete (informally, this means that it can be identified with a 1\hyp{}category) and if the functor $c_B : \A\to \B$ is constant on a single object $B$, \ie $c_B(A) \equiv B$ for each $A\in \A$, and $c_B(f)\equiv \id_B$, then $\twoint_A c_B$ is canonically identified with the \emph{cotensor} of $B$ by $\A$ and is denoted $\A\pitchfork B$.
\end{example}
\begin{proof}
	It is easily seen that a lax wedge $B'\to \twoint c_B$ corresponds to a family of maps $B'\to B$ indexed by the objects of $\A$; the wedge condition now entails that the correspondence $A\mapsto f_A : B'\to B$ can be lifted to a functor $\A \to \B(B',B)$. This concludes the proof of an isomorphism
	\[\Cat(\A, \B(B',B))\cong \B\Big(B', \twoint_{\B} c_B\Big)\qedhere\]
\end{proof}
\begin{theorem}[Parametric lax ends]
	Whenever we have a functor $F : \A^\opp\times \A\times \B\to \C$, and the lax end $\twoint_A F(B, A,A)$ exists for every object $B\in \B$, then $B\mapsto \twoint_A F(B, A,A)$ extends to a 2\hyp{}functor $\B\to \C$ which has the universal property of the lax end of its mate $\hat F : \A^\opp\times \A \to \C^{\B}$ under an obvious cartesian closed adjunction.
\end{theorem}
We omit the proof (one way to proceed is to show how a functor $F : \A^\opp\times \A\times \B\to \C$ transposes to $\hat F :  \B\to \twoCat(\A^\opp\times \A,\C)$ once the functoriality of $\twoint : \twoCat(\A^\opp\times \A,\C) \to \C$ has been proved); we leave the reader enjoy this exercise, letting them discover how strict the resulting functor $\twoint_A F(\firstblank,A,A)$ shall be.

We also offer without proof the statement the 2\hyp{}dimensional analogue of \emph{Fubini theorem} in \ref{fubozzo}. This is a much more daunting exercise, but the challenge is merely notational, not conceptual.
\begin{theorem}[Fubini rule for lax co\fshyp{}ends]\label{laxfubini}\index{Fubini rule!--- for lax co\fshyp{}ends}
	If one among the following lax ends exists, then so does the others, and the three are canonically isomorphic:
	\[
		\twoint_{B,C}T(B,C,B,C)\cong \twoint_B\left(\twoint_C T(B,C,B,C)\right) \cong \twoint_C\left(\twoint_B T(B,C,B,C)\right)
	\]
	In a similar fashion,
	\[
		\twoint^{B,C}T(B,C,B,C)\cong \twoint^B\left(\twoint^C T(B,C,B,C)\right) \cong \twoint^C\left(\twoint^B T(B,C,B,C)\right)
	\]
\end{theorem}
\begin{corollary}[Fubini rule for lax co\fshyp{}limits]
	Let $T : \B\times \C\to \D$ be a 2\hyp{}functor; then
	\begin{enumtag}{lc}
		\item lax limits commute: we have canonical isomorphisms in $\D$
		\[
			\llim_{B\in \B}
			\llim_{C\in \C} T(B,C)
			\cong
			\llim_{C\in \C}
			\llim_{B\in \B}  T(B,C);
		\]
		\item lax colimits commute: we have canonical isomorphisms in $\D$
		\[
			\lcolim_{B\in \B}
			\lcolim_{C\in \C} T(B,C)
			\cong
			\lcolim_{C\in \C}
			\lcolim_{B\in \B}  T(B,C).
		\]
	\end{enumtag}
\end{corollary}
The hom functor must preserve lax ends:
\begin{theorem}\label{homcommuteslaxcoend}
	Let $S : \A \to \B$ be a 2\hyp{}functor such that $\twoint_A S(A,A)$ exists in $\B$; then we have a canonical isomorphism of categories
	\[\label{laxcomm}\B\Big(B, \twoint_A S(A,A)\Big) \cong \twoint_A \B(B,S(A,A)).\]
	A completely dual statement involves the lax coend of $S :\A\to \B$: there is a canonical isomorphism of categories
	\[\label{laxcommdue}\B\Big(\twoint^A S(A,A),B\Big) \cong \twoint_A \B(S(A,A),B).\]
\end{theorem}
This can be seen also as an alternative definition: the lax co\fshyp{}end of $S$ exists if and only if for every object $B\in\B$ the functor $\A^\opp\times \A \xto{S} \B \xto{\hom_B} \Cat$ ($\hom_B$ is understood here as a covariant or contravariant representable over $B\in\B$) has a lax co\fshyp{}end; indeed if this is the case, there is the isomorphism above.
\subsection{The lax ninja Yoneda lemma}\label{laxyonedaninja}\index{Yoneda lemma!lax ---}
In the world of 2\hyp{}categories and lax morphisms, the ninja Yoneda lemma \ref{ninjayo} acquires a peculiar flavour, since it is equivalent to the co\fshyp{}reflectivity of strict presheaves $F : \C^\opp\to \Cat$ inside lax presheaves $F' : \C^\opp\laxto \Cat$.

More precisely, there is a diagram of adjoint 2\hyp{}functors \index{_aaa_diesis@$(\firstblank)^\diesis$}\index{_aaa_diesis'@$(\firstblank)^\bemo$}
\[\vcenter{\xymatrix{
			\Cat(\C^\opp, \Cat) \ar[r]&\ar@<6pt>[l]^{(\firstblank)^\bemo}\ar@<-6pt>[l]_{(\firstblank)^\diesis} \Cat_l(\C^\opp, \Cat)\\
		}}\]
where the central arrow is the inclusion. This means that for each strict 2\hyp{}functor $H\in \twoCat(\C^\opp, \Cat)$ there are two natural isomorphisms
\begin{gather}
	\twoCat(\C^\opp, \Cat)(H, F^\bemo) \cong \twoCat_l(\C^\opp, \Cat)(H,F),\notag \\
	\twoCat(\C^\opp, \Cat)(F^\diesis, H)\cong \twoCat_l(\C^\opp, \Cat)(F,H).
\end{gather}
The proof is completely formal. To prove the isomorphism, we will show that the functors $F^\bemo$ and $F^\diesis$ defined by means of the lax coends
\[
	F^\diesis \cong \twoint^A \C(\firstblank,A)\times FA \qquad\qquad
	F^\bemo \cong \twoint_A \C(A,\firstblank)\pitchfork FA.
\]
have the above property.
\begin{proof}
	The proof exploits \ref{laxnat} as well as the evident fact that strict co\fshyp{}ends are particular cases of lax co\fshyp{}ends, and (as a consequence) the fact that any number of co\fshyp{}ends and lax co\fshyp{}ends commute with each other by virtue of the lax Fubini rule, plus the preservation of lax co\fshyp{}ends by the hom functor \ref{homcommuteslaxcoend}, and the strict ninja Yoneda lemma \ref{ninjayo}.

	Now, let us denote $F^\bemo = \twoint_A \C(A,\firstblank)\pitchfork FA$; we have that
	\begin{align*}
		\Cat(\C^\opp,\Cat)(H, F^\bemo) & = \int_C \Cat(HC, F^\bemo C)                                    \\
		                               & = \int_C \Cat \left(HC, \twoint_A \Cat(\C(A,C), FA) \right)     \\
		                               & \cong \int_C \twoint_A \Cat(HC, \Cat(\C(A,C), FA))              \\
		                               & \cong \twoint_A \int_C \Cat(HC\times \C(A,C), FA)               \\
		                               & \cong \twoint_A \Cat\left( \int^C HC \times \C(A,C), FA \right) \\
		                               & \cong \twoint_A \Cat(HA, FA)                                    \\
		\eqref{this_laxnat}            & = \twoCat_l(\C^\opp, \Cat)(H,F).
	\end{align*}
	The proof that $\Cat(\C^\opp, \Cat)(F^\diesis, H)\cong \twoCat_l(\C^\opp, \Cat)(F,H)$ is similarly formal, thus we leave it to the reader.
\end{proof}
Lax coend defining $F^\diesis$ can take as argument a very simple functor, and spit out a not-so-simple one: an instructive example is given by the $(\firstblank)^\diesis$ of the terminal functor:
\begin{example}
	Let $\I : \C^\opp\to \Cat$ be the constant functor sending $C\in\C$ into the terminal category, regarded as a lax functor. Then the strict functor $\I^\diesis : \C^\opp\to \Cat$ is the lax coend
	\[
		\I^\diesis C \cong \twoint^A \I A\times \C(A,C)\cong \twoint^A \C(A,C).
	\]
	Hence the category $\twoint^A \C(A,C)$ coincides with the lax colimit of the strict presheaf $\C^\opp\to \Cat$, $A\mapsto \C(A,C)$, which is \cite[p. 171]{Street19} the \emph{lax slice category} $\C/\!\!/C$ of commutative diagrams of 2\hyp{}cells
	\[\vcenter{\xymatrix{
				A \ar[rr]^h \ar[dr]_u \drrtwocell<\omit>{<-1>\alpha}&& A'\ar[dl]^{u'} \\
				& C &
			}}\]
\end{example}
Another striking instance of this phenomenon, by which quite complex shapes can arise as lax co\fshyp{}limits of simpler diagrams, is given by the twisted arrow category of $\A$, whose opposite category is isomorphic to the lax colimit of the `slice category' diagram.
\begin{proposition}[The twisted arrow category as a lax colimit]
	Let $\A$ be a category; then it is possible to characterise the (opposite of the) twisted arrow category of $\A$, as defined in \ref{twisted}, as the lax colimit of the diagram $\A \to \Cat : A\mapsto \A/A$, \ie as the lax coend $\twoint^A \A/A$.
\end{proposition}
\begin{proof}
	In order to prove the statement, we shall first find a good candidate for a universal cocone. Such cocone $q : \A/\firstblank\to \tw(\A)^\opp$ is defined on objects sending every $f : X \to A$ into itself, and a morphism
	\[
		\vcenter{\xymatrix{
				X \ar[rr]^u \ar[dr]_f && X'' \ar[dl]^{f'} \\
				& A &
			}}
	\]
	again into itself, regarding the triangle as a square with identity bottom:
	\[
		\vcenter{\xymatrix{
				X\ar[r]^u \ar[d]_f & X'\ar[d]^{f'} \\
				A \ar@{=}[r]& A
			}}
	\]
	(as an aside comment, this definition makes obvious that the functor we are defining is conservative; we shall however get an indirect proof of invertibility).

	Let us denote $t_* : \A/A \to \A/A'$ the functor induced by $t :A \to A'$ by post\hyp{}composition. Then, we get a tautological laxity cell $q_t : q_A \To q_{A'}\circ t_*$, having components
	\[
		\vcenter{\xymatrix@R=5mm{
		X \ar[dd]_f \ar@{=}[r]& X \ar[d]^f\\
		& A \ar[d]^t \\
		A \ar[r]_t & A'
		}}
	\]
	These components glue together defining a lax cocone $q : \A/\firstblank\to \tw(\A)$.

	We shall now prove that this is universal: in order to do so, we have to verify that the diagram
	\[ \vcenter{\xymatrix@R=5mm{
		\A/A \ar[dd]_{t_*} \ar[dr]^{q_A} \ar@/^1.5pc/[drr]^{q_A} && \\
		& \tw(\A)^\opp \ar@{.>}[r]^u \ar@{}[l]|(.6){\Swarrow q_t} & \K \\
		\A/A'\ar[ur]_{q_{A'}} \ar@/_1.5pc/[urr]_{q_{A'}} &&
		}} \]
	satisfies the (dual of the) conditions in \ref{laxlimo}.
	\begin{itemize}
		\item Each time we are given a lax cocone $h : \A/\firstblank \to \tw(\A)^\opp$, there is a \emph{unique} induced 1\hyp{}cell $u : \tw(\A)^\opp\to \K$ such that $u * p = h$.
		\item Every 2\hyp{}cell $\sigma : u\To u'$ such that the horizontal composition $\sigma \boxminus q_t$ in
		      \[ \label{iddu} \vcenter{\xymatrix@R=5mm@C=1.25cm{
			      \A/A \ar[dd]_{t_*} \ar[dr]^{q_A} && \\
			      & \tw(\A)^\opp  \ar@{}[l]|(.6){\Swarrow q_t} \rtwocell^u_{u'}{\sigma}& \K \\
			      \A/A'\ar[ur]_{q_{A'}} &&
			      }} \]
		      is invertible for every $t :A \to A'$, is itself invertible.
	\end{itemize}
	Given $h$ as above, since for every $t :A \to A'$ we have
	\[
		h_{A'}(t\circ f) = h_A(f)
	\]
	we can define $u : \tw(\A)^\opp\to \K$ sending the object $f : X \to A$ and the morphism $(s,t) : f \to f'$ respectively to
	\begin{itemize}
		\item $u\!\var{X}{A} = h_A\!\var{X}{A} \in\K$;
		\item $u\bsmat X & \overset{s}\leftarrow & Y \\ \downarrow && \downarrow \\ A & \underset{t}\to & B  \esmat = h_A \!\var{X}{A} \leftarrow h_B\!\var{Y}{B}$ obtained as the composition
		      \[ h_B(t) = h_B(tfs) \to h_A(fs) \to h_Af \]
		      where in the first step we use the laxity cell $h_t : h_B(tfs) \to h_A(fs)$ and in the second we apply $h_A$ to the morphism $f \to fs$ in $\A/A$.
	\end{itemize}
	Using these definitions, it is rather easy to see that $u *p =h$ as components of suitable modifications, since
	\[ u(q_t) = u\bsmat X & = & X \\ \downarrow && \downarrow \\ A &\underset{t}\to& A'\esmat \]
	or in other words to the action of the laxity cell $q_t : u(tf)\to u(f)$ composed with $q_{\id_A} : q_A(f) \to q_A(f)$.

	It is equally easy to see that every morphism in $\tw(\A)^\opp$ can be expressed as a morphism in the image of some $q_A$, and as a consequence, $u$ is uniquely determined by the properties we just checked.

	It remains to check the last condition; however, thanks to the definition of $q_t$, if the diagram \eqref{iddu} satisfies the property that the horizontal composition $\sigma \boxminus q_t$ is invertible for a 2\hyp{}cell $\sigma : u\To u'$, and each component of the composition $\sigma \boxminus q_{\id_A}$ is invertible, then each component $\sigma_{(X\to A)} : u\!\var{X}{A} \to u'\!\var{X}{A}$ is invertible as well.
\end{proof}
The following result is a partial analogue of \ref{ends-are-weighted}: in order to characterize a lax co\fshyp{}end as a weighted co\fshyp{}limit, we have to `twist' the hom functor using the $\diesis$ construction of \ref{laxyonedaninja}.
\begin{proposition}\label{lax.is.wcolim}
	\cite[§2]{bozapalides1980some} There is a canonical isomorphism between the lax end of a 2\hyp{}functor $T : \C^\opp\times \C\to \B$ and the limit of $T$ weighted by the bifunctor $\C((\firstblank)^\diesis,\secondblank)$, \ie
	\[
		\wlim{\C((\firstblank)^\diesis,\secondblank)}T \,\cong\, \twoint_A T(A,A)
	\]
	where $\C((\firstblank)^\diesis,\secondblank) : (C, C')\mapsto \C(\firstblank,C')^\diesis(C)= \twoint^A \C(C,A)\times \C(A, C')$. A dual statement holds for lax coends.
\end{proposition}
\begin{proof*}
	The proof is completely formal, by virtue of the results established so far: we can compute
	\begin{align*}
		\B\left( B, \wlim{\C((\firstblank)^\diesis,\secondblank)}T \right) & = \wlim{\C((\firstblank)^\diesis,\secondblank)} \B(B, T) \tag{see \ref{homcommuteswei}} \\
		(\ref{wlimcoends})                                                 & \cong \int_{C,D} \B(B, T(C,D))^{\C(C^\diesis,D)}                                        \\
		                                                                   & = \int_{C,D} \B(B, T(C,D))^{\twoint^A \C(C,A)\times \C(A, D)}                           \\
		                                                                   & \cong \twoint_A \int_{C,D} \left(\B(B, T(C,D))^{\C(A,D)} \right)^{\C(C,A)}              \\
		(\ref{laxfubini})                                                  & \cong \twoint_A \int_C \Cat\left( \C(C,A), \int_D \B(B, T(C,D))^{\C(A,D)} \right)       \\
		                                                                   & \cong \twoint_A \int_C \Cat\left( \C(C,A), \B(B, T(C,A)) \right)                        \\
		                                                                   & \cong \twoint_A \B(B, T(A,A)) = \B\left(B, \twoint_A T(A,A)\right).
	\end{align*}
\end{proof*}
Let $S : \A^\opp\to \Cat$, $T : \A\to \B$ be two functors and suppose $\B$ has $\Cat$\hyp{}tensors; then the lax coend of the 2\hyp{}functor $\A^\opp\times \A\xto{S\times T} \Cat\times\B\xto{\otimes} \B$ is the \emph{tensor product} of $S$ and $T$, denoted
\[
	S\mathop{\overline{\otimes}} T =: \twoint^A Sa\otimes Ta.
\]
\subsection{$2$\hyp{}profunctors, lax Kan extensions}
The present section is intended to provide an analogue of profunctor theory for lax co\fshyp{}ends; it can be safely skipped at first reading.

The proof of \ref{lax.is.wcolim} above suggests that $\C((\firstblank)^\diesis,\secondblank)$ is the \emph{lax composition} of two representable profunctors. The intuition that the theory of Chapter 5 can be suitably adapted to the lax context is correct (unfortunately though, we have few applications in mind for such a theory: exercise \ref{ex:laxprofunctors} below tries to suggests a few lines of investigation).

Shortly, a \emph{2\hyp{}profunctor} $\proP : \A\leadsto \B$ is a 2\hyp{}functor $\proP : \B^\opp\times \A\to \Cat$. Lax coends provide a weak composition rule for 2\hyp{}profunctors: more precisely, let
\[
	\A \overset{\proP}{\pto}\B \overset{\proQ}{\pto} \C
\]
be a couple of composable 2\hyp{}profunctors, namely two 2\hyp{}functors $\proP : \B^\opp\times \A\to \Cat$ and $\proQ : \C^\opp\times \B\to \Cat$; then the composition $\proQ\diamond \proP$ is defined by the coend
\[
	\proQ\diamond_l \proP(C,A) = \twoint^B \proP(B,A)\times \proQ(C,B)
\]
The compatibility between lax colimits and products ensures that the expected associativity holds up to a canonical identification:
\[
	(\proH \diamond_l \proQ)\diamond_l\proP \cong \proH \diamond_l (\proQ\diamond_l \proP)
\]
for any three $\A\overset{\proP}\pto\B\overset{\proQ}\pto\C\overset{\proH}\pto\D$.

A much more useful laxification involves Kan extensions: let $T : \A\to \B$ and $F : \A\to \C$ be two 2\hyp{}functors.
\begin{definition}
	We call \emph{left lax Kan extension} of $F $ along $T$ a 2\hyp{}functor $\lLan_TF  : \B\to \C$ endowed with a lax natural transformation $\eta : F  \To_l \lLan_TF \circ T$ (a \emph{unit}) such that, for every pair $S : \B\to \C$ and $\lambda : F  \To_l S\circ T$ (respectively, a 2\hyp{}functor and a lax natural transformation) there exists a unique $\zeta : \lLan_TF  \To_l S$ such that
	\[
		(\zeta * T)\circ \alpha = \lambda
	\]
	and moreover, if $\Sigma : \lambda \Rrightarrow \lambda'$ is a modification between lax natural transformations, there is a unique modification $\Omega : \zeta\Rrightarrow \zeta'$ (where $\zeta$ is induced by $\lambda$, and $\zeta'$ by $\lambda'$) between $\Cat$\hyp{}natural transformations such that $(\Omega * T)\circ \alpha = \Sigma$.

	This entire list of conditions can be expressed by means of the isomorphism
	\[
		\twoCat(\A, \C)\big( F , S\circ T \big)\cong \twoCat(\B,\C)\big( \lLan_TF , S \big)
	\]
	which is natural in $S$.
\end{definition}
\begin{example}
	Let $1$ be the terminal 2\hyp{}category. Then the left lax Kan extension of a 2\hyp{}functor $F  : \A\to\C$ along the terminal 2\hyp{}functor $\A\to 1$ is the lax colimit of $F $.
\end{example}
\begin{remark}
	We can obtain different flavours of lax Kan extension by reversing the directions of $\alpha, \lambda, \zeta\dots$ and imposing invertibility.
	The example above, as well as the following theorem, shows that the choice of $\Cat$\hyp{}natural transformations instead of lax natural transformations is the right choice (see also \cite{bozapalides1980some} for a dual statement):
\end{remark}
\begin{theorem}\label{pseudolan}
	In the same notation as above, assume $\C$ admits tensors for hom categories $\B(TA, B)\otimes F A'$ for each $A, A'\in \A$, and that the lax coend
	\[
		\twoint^A \B(TA, B)\otimes F A
	\]
	exists; then the lax Kan extension of $F$ along $T$ exists too, and it is canonically isomorphic to the coend above.
\end{theorem}
We can mimic also Exercise \ref{closed.via.coends} to obtain a lax analogue of it:
\begin{proposition}
	Let $\twoCat_l(\C,\Cat)(U,V)$ denote the category of lax natural transformations between two 2\hyp{}functors $U,V : \C \to \Cat$. Then
	\[
		\twoCat_l(\C,\Cat)(F\times G,H)\cong \twoCat_l(\C,\Cat)(F, H^G),
	\]
	where $H^G(x) = \twoCat_l(\C,\Cat)(\yon(A)\times G, H) = \twoint_{\,Y}\Set(\hom(Y,X)\times GY, HY)$
\end{proposition}
\begin{proof*}
	Every step can be motivated by results in the present section:
	\begin{align*}
		\twoCat_l(\C,\Cat)(F, H^G) & = \twoint_X \Sets(FX, \twoCat_l(\C,\Cat)(\yon(A)\times G, H))                    \\
		                           & \cong \twoint_X \twoint_Y \Sets(FX, \Set(\C(Y,X)\times GY, HY))                  \\
		                           & \cong \twoint_Y\Set\Big(\Big(\twoint^X FX \times \C(Y,X)\Big)\times GY, HY \Big) \\
		                           & \cong \twoint_Y\Set\big( FY\times GY, HY \big)                                   \\
		                           & =\twoCat_l(\C,\Cat)(F\times G,H).
	\end{align*}
\end{proof*}
\section{Coends in homotopy theory}
\index{Co/end!homotopy ---}
Higher category theory is nowadays living a Renaissance, thanks to a massive collaboration of several people drawing from various fields of research, and cooperating to re\hyp{}analyze every feature of category theory inside, or in terms of, the topos of simplicial sets.

The purpose of the present section is to study what this `homotopification' process does to co\fshyp{}end calculus.

The urge to keep this chapter self\hyp{}contained, force us to take for granted a certain acquaintance with model categories, simplicial categories, $\infty$\hyp{}categories \emph{à la} Joyal-Lurie, dg\hyp{}categories\dots{} Each of these theories is vast and constitutes enough material for a dedicated monograph. Nonetheless, we strive to offer the best intuition we can, in the little space we have.

We start by presenting the theory of \emph{homotopy co\fshyp{}ends} \ref{coends-in-model} in model category theory; the coend functor $\int^\C : \Cat(\C^\opp\times\C,\D) \to \D$ admits a `derived' counterpart $\int^\C_{\mathbb L}$ (see \cite{Isaacson}) that preserves weak equivalences, in the same sense the colimit functor does.

Subsequently, in \ref{coend-in-qcat}, we will explore the theory of \emph{quasicategorical} co\fshyp{}end calculus (providing a proof of the Fubini rule for $\infty$\hyp{}coends: this is not strictly new material, but brings together the pieces present in \cite{gepner2015lax}).

Subsequently, in \ref{coend-in-ssetcat}, we address the study of \emph{simplicially coherent} co\fshyp{}ends (\ie enriched co\fshyp{}ends in $\sSet$\hyp{}categories), and the definition of a co\fshyp{}end in a \emph{derivator} (this is nothing more than a paragraph; some additional results are presented as exercises \ref{ex8:derivcoend}, \ref{ex8:derivfubini}, and \ref{ex8:derivweighlim}).

The discussion closes the circle about co\fshyp{}end calculus in each of the most common models for higher category theory (model categories, simplicially enriched categories, simplicial sets, derivators). The reader will notice that most of these flavours of co\fshyp{}end calculus is far from being a full\hyp{}fledged theory. This testifies how much work there is still to be done in the field!

As authors, we hope to have convinced our reader that the endeavor to complete the theory is worth trying. Given its utility, co\fshyp{}end calculus shall reach a status of well\hyp{}understood and standardised tool for category theorists, in such a way all its users, no matter what is their origin or purpose, can profit from its conceptual simplicity.

The age to set up and use a full\hyp{}fledged \emph{homotopy coherent co\fshyp{}end calculus} is of our readers.
\begin{remark}
	We set aside a rather important question, that is \emph{model dependence}: the models to study the theory of $(\infty,1)$\hyp{}categories form a complicated web of equivalences. Uniqueness results \cite{Barwick2011,toennone} or even synthetic approaches \cite{riehlverity-1,riehlverity-4,Riehl2015a} are nowadays very trendy as they offer extremely powerful inter\hyp{}model comparison theorems, but things can become rather hairy in explicitly proving that (say) an homotopy co\fshyp{}end corresponds to the same notion of an $\infty$\hyp{}co\fshyp{}end, if we move from one model to the other.

	This is a subtle issue, derailing us from our primary objective; we thus choose to bluntly put it aside, but we maintain at least an agnostic point of view towards the matter: as no model is privileged, we shall glance at them all.
\end{remark}
\subsection{Coends in model categories}\label{coends-in-model}
A \emph{model category} is a category $\C$ having all small co\fshyp{}limits (or, in the original definition of Quillen, all finite co\fshyp{}limits), which is endowed with three distinguished classes of morphisms, the \emph{weak equivalences}, the \emph{fibrations} and the \emph{cofibrations}, interacting in such a way that the following properties are satisfied:
\begin{itemize}
	\item The initial category $\C[\wk^{-1}]$ where all weak equivalences become invertible admits a presentation as a category with the same objects, which is a hom\hyp{}wise quotient of $\C$: more precisely, there is an equivalence relation $R_{XY}$ on each $\C(X,Y)$ such that $\C[\wk^{-1}](X,Y) \cong \C(X,Y)/R_{XY}$; we call $\C[\wk^{-1}]$ the \emph{homotopy category} of $\C$ with respect to $\wk$.
	\item The pairs $(\wk\cap\cof, \fib)$ and $(\cof,\wk\cap\fib)$ form two \emph{weak factorisation systems} on $\C$ (simply put, this means that there is a way to factor every $f :  \to Y$ in $\C$ as a composition $X \xto{\cof} E \xto{\wk\cap\fib} Y$ and as a composition $X \xto{\wk\cap\cof} E' \xto{\fib} Y$).
\end{itemize}
\begin{definition*}
	Let $\C$ be a model category; we say that an object $A\in\C$ is \emph{cofibrant} if its initial arrow $\varnothing \to A$ is a cofibration. Dually, we say that an object $X$ is \emph{fibrant} if its terminal arrow $X \to *$ is a fibration.
\end{definition*}
One of the most important parts of model category theory is the study of \emph{homotopy co\fshyp{}limits}. To put it shortly, the vastness of homotopy theory and homological algebra arises from the fact that the colimit functor $\colim$ is rather ill\hyp{}behaved in terms of its interaction with weak equivalences. The story goes as follows.

It turns out that many common functors between model categories do not send weak equivalences to weak equivalences; so their behaviour must be corrected, either restricting their domain to subcategories of suitably nice objects, or by changing the shape of the functors themselves; an illustrious example of such common but ill\hyp{}behaved functor is the colimit $\varinjlim_{\cate J} : \C^\J\to \C$ over a diagram of shape $\J$: given a class of weak equivalences $\cW$ then every diagram category $\C^{\J}$ acquires an analogous structure $\cW^\J$, where $\eta : F\To G$ is in $\cW^\J$ if and only if each component $\eta_j : Fj\to Gj$ is itself a weak equivalence.

It turns out that the image of such a natural transformation $\eta : F\To G$ under the colimit functor, $\varinjlim\eta : \varinjlim F\to \varinjlim G$ is rarely a weak equivalence in $\C$.\footnote{A minimal example of this goes as follows: take $\cate J$ to be the generic span $1\leftarrow 0\to 2$ and the functor $F : \J \to \Spc$ sending it to $\{*\} \leftarrow S^{n-1}\to \{*\}$ ($S^k$ is the $k$\hyp{}dimensional sphere); the colimit of $F$ is the one\hyp{}point space $\{*\}$. We can replace $F$ with the diagram $\tilde F : D^2 \leftarrow S^{n-1}\to D^2$, and since disks are contractible there is a homotopy equivalence $\tilde F \To F$; unfortunately, the induced arrow $\varinjlim \tilde F = S^2\to *$ is not a weak equivalence (because $S^2$ is not contractible).}

One of the main tenets of homotopy theory is, nevertheless, that it doesn't matter if we replace an object of a model category with another, as soon as the two become (controllably) isomorphic in the homotopy category. There is a spark of hope, then, that the category of functors $\C^\J$ contains a better\hyp{}behaved representative for the functor $\varinjlim$, and that the two are linked by some sort of natural weak equivalence.

That's what homotopy colimits are: they provide \emph{deformations} $\hocolim$ of $\varinjlim$ that preserves pointwise weak equivalences, and are linked by object\hyp{}wise weak equivalence $\hocolim \To \colim$.
\begin{remark}\label{coend_is_quillen}
	Let $\boxtimes : \A \times \B \to \C$ be the `tensor' part of a thc situation (see Exercise \ref{tiaccaci}), and let us assume that it is a left Quillen functor (see \cite[1.3.1]{Hov}); let $\J$ be a Reedy category (\cite[5.2.1]{Hov}). Then the coend functor
	\[
		\int^{\J} : \Cat(\J^\opp,\A) \times \Cat(\J,\B) \to \C
	\]
	is a left Quillen bifunctor if we regard the functor categories $\Cat(\J^\opp,\A)$ and $\Cat(\J,\B)$ endowed with the Reedy model structure.

	The coend functor remains left Quillen even if $\J$ is not Reedy, but we have to assume the categories $\A,\B,\C$ are all combinatorial,\footnote{A model category is \emph{combinatorial} if it is locally presentable \cite[1.17]{Adamek1994} and cofibrantly generated, \ie if the cofibrations are generated as the weak orthogonal (see \ref{accepre}) of a small set $I$.} and we put the projective model structure on $\Cat(\J,\B)$, and the injective model structure on $\Cat(\J^\opp,\A)$.\footnote{It's not worth to enter the details, but see \cite{Hov,Hirschhorn2003} for more information: the \emph{projective} model structure on $\Cat(\J,\C)$ is determined by the fibrations of $\C$, and the \emph{injective} model structure is determined by cofibrations of $\C$, in a suitable sense.}
\end{remark}
In the following, $\C$ will be a small category, and $F : \C \to \D$, $G : \C^\opp\to \Set$ will be functors; when needed, we freely employ the notation of Chapter 4 on weighted co\fshyp{}limits. The category $\D$ will admit the co\fshyp{}limits allowing the object we define to exist.
\index{Bar construction}
\index{Cobar construction|see{Bar construction}}
\begin{definition}[Bar and cobar complexes]\leavevmode
	\begin{itemize}
		\item The bar complex of a pair of functors $F,G$ is the simplicial object in $\D$ $B(G,\C,F)_\bullet$ whose set of $n$\hyp{}simplices is
		      \[
			      B(G,\C,F)_n := \coprod_{C_0,\dots,C_n} (N(\C)_n\times GC_n)\otimes FC_0
		      \]
		      where $\otimes : \Set\times \D \to \D$ is the tensor functor of \ref{tenscotens}. More explicitly, $B(G,\C,F)_n$ is the disjoint union over $(n+1)$\hyp{}tuples of objects of $\C$ of sets
		      \[\big(GC_n\times\C(C_0,C_1)\times\dots\times \C(C_{n-1}, C_n)\big)\otimes FC_n.\]
		      Dually,
		\item the cobar complex of a pair of functors $F,G$ is the cosimplicial object in $\D$ $C(G,\C,F)^\bullet$ whose set of $n$\hyp{}simplices is
		      \[
			      C^n(G,\C,F) := \prod_{C_0,\dots,C_n} (N(\C)_n\times GC_n)\pitchfork FC_0
		      \]
		      where $\pitchfork : \Set^\opp\times \D \to \D$ is the cotensor functor of \ref{tenscotens}. More explicitly, $C(G,\C,F)^n$ is the product over $(n+1)$\hyp{}tuples of objects of $\C$ of sets
		      \[\big(GC_n\times\C(C_0,C_1)\times\dots\times \C(C_{n-1}, C_n)\big)\pitchfork FC_n.\]
	\end{itemize}
\end{definition}
We leave as Exercise \ref{barciobar} the check of some elementary properties of these objects; in particular, we will constantly exploit the fact that $B(G,\C,F)$ and $C(G,\C,F)$ are functorial in $F,G$ (with appropriate variance).

The proof of the following statement is conducted in the cited reference. The bar and cobar constructions allow to reduce the computation of every (weighted) co\fshyp{}limit to the computation of a certain (weighted) co\fshyp{}limit over the simplex category $\bDelta$ (these colimits and limits are called \emph{diagonalisation} and \emph{totalisation} respectively, and will return in \ref{totalizia} and \ref{diagonalizia} to define simplicially coherent co\fshyp{}ends).
\begin{theorem}[\protect{\cite[4.4.2]{riehl2014categorical}}]
	Let $F : \J \to \D$ be a functor, and $W : \J \to \Set$ be a weight (covariant, if we compute a limit; contravariant, if we compute a colimit in the formula below). Then we have canonical isomorphisms
	\begin{align}
		\wlim{W}F    & \cong \colim_{\bDelta} B(G,\J,F)_\bullet\notag \\
		             & \cong \int^{n\in\bDelta} B(G,\J,F)^n           \\
		\wcolim{W} F & \cong \lim_{\bDelta} C(G,\J,F)^\bullet \notag  \\
		             & \cong \int_{n\in\bDelta} C(G,\J,F)^n.
	\end{align}
\end{theorem}
The same result holds (and it is much more interesting) over most (but not all: the proof relies on the presence of an \emph{enriched Grothendieck construction}) bases of enrichment; the most important instance of such an enrichment base is that of simplicial sets.

Moreover, the bar and cobar constructions provide \emph{replacement} functors for $\colim$ and $\lim$ respectively, and as a consequence they provide models for the homotopy colimit and limit functors:
\index{_aaa_yoncontra@$\yon$}
\index{_aaa_yoncov@$\coyon$}
in the following, if $F: \J \to \D$ is a functor and $\D$ is a simplicial model category, we denote $B(\J,\J,F)_\bullet$ the simplicial set $[n]\mapsto B(\yon J, \J, F)_n$, and dually $C(\J,\J,F)^n$ is the cosimplicial set $[n]\mapsto C(\coyon J, \J, F)^n$. (The statement found in \cite{riehl2014categorical} is more general than the one we introduce here; its essence is, however, unchanged.)
\begin{theorem}[\protect{\cite[5.1.3]{riehl2014categorical}}]
	The left derived functor of $\colim : \Cat(\J,\D)\to \D$ and the right derived functor of $\lim : \Cat(\J,\D)\to \D$ are computed respectively as
	\[\mathbb{L}\colim F \cong \diag \big(B(\J,\J,\tilde F)\big) \qquad
		\mathbb{R}\lim F \cong \tot\big(C(\J,\J,\hat F)\big)\]
	where $\tilde F$ is a cofibrant replacement, and $\hat F$ is a fibrant replacement for the diagram $F$, in suitable model structures on the diagram category $\Cat(\J,\D)$
	\ie as the co\fshyp{}ends
	\begin{gather}
		\mathbb{L}\colim F \cong \int^n \bDelta(\firstblank, [n])\otimes B(\Delta[n],\J,\tilde F) \notag\\
		\mathbb{R}\lim F \cong \int_n \bDelta(\firstblank, [n])\pitchfork C(\Delta[n],\J,\hat F)
	\end{gather}
\end{theorem}
The other way to compute homotopy co\fshyp{}limits involves a resolution of the weight: the equivalence between the two approaches was fist proved in \cite{Gamb}.

Recall \ref{busfokane}: we can express the Bousfield\hyp{}Kan construction (\cite{Bousfield1972}) for the homotopy co\fshyp{}limit functor using co\fshyp{}end calculus.
\begin{theorem*}
	Let $F \colon \J\to \C$ be a diagram in a model category $(\C, \wk,\cof,\fib)$ which is tensored and cotensored (see \ref{tenscotens}) over simplicial sets. Then the homotopy limit $\holim F$ of $F$ can be computed as the end
	\[
		\int_J N(\J/J)\pitchfork F(J),
	\]
	where $N$ is the categorical nerve of \ref{catnerve}; in the same notation, the homotopy colimit $\hocolim F$ of $F$ can be computed as the coend
	\[
		\int^J N(\J/J) \otimes F(J).
	\]
\end{theorem*}
\begin{remark}\label{busfokane-explained}
	These two universal objects are the weighted co\fshyp{}limit of $F$ with the nerve of a slice category as weight; the idea behind this characterisation is that $\colim$ is a weighted colimit over the terminal weight. The problem is that usually the constant terminal weight won't be a cofibrant object in $\Cat(\J,\sSet)$; thus, when we want to pass to the homotopy correct version of $\colim$ we shall replace the weight $W$ with a homotopy equivalent, but cofibrant, diagram $\tilde W$.

	The Bousfield\hyp{}Kan formula arises precisely in this process: $N(J/\J)$ and $N(\J/J)$ are both contractible categories, and they are linked to $N(*)$ (the nerve of the terminal category) by an homotopy equivalence induced by the terminal functor. The simplicial presheaves $N(\J/\firstblank),N(\firstblank/\J)$ can thus be thought as proper \emph{replacements} for the terminal functor.
\end{remark}
\subsection{Simplicially coherent co\fshyp{}ends}\label{coend-in-ssetcat}\index{Co/end!simplicially coherent ---}
All the material in the following subsection comes from \cite{cordier1997homotopy}. We begin the exposition establishing a convenient notation and a series of useful short\hyp{}hands to adapt the discussion to our choice of notation.

We strive to keep this introduction equally self\hyp{}contained and simple, but the reader shall be warned that
\begin{itemize}
	\item there is a sheer amount of unavoidable sins of omissions in the exposition, essentially due to our inability to master the topic in its entirety; moreover, the price we pay to obtain a self\hyp{}contained exposition is that we deliberately ignore most of the subtleties of the combinatorics of simplicial sets. The blame is on us if the reader feels that our exposition is gawky or incomplete.
	\item Since the times when \cite{cordier1997homotopy} was published, newer and more systematic approaches to a similar topic were developed; among many, the reader shall take the exceptionally clear \cite{riehl2014categorical,shulman}. These references reduce the construction of a simplicially coherent co\fshyp{}end to an application of the `unreasonably effective' co\fshyp{}bar construction \cite[4]{riehl2014categorical}, and it can be proved using \cite[21.4]{shulman} that the coherent co\fshyp{}end of $T$ (see \ref{cohcoend}) results as a suitable \emph{derived weighted co\fshyp{}limit} of the functor $T$. This last remark is important, in particular in light of the fact that it can be exported to define an homotopy coherent co/end calculus in every 2-category $\VCat$ of categories enriched over a monoidal model category (see \cite[Ch. 4]{Hov} and \cite{Ber}).\footnote{This remark deserves better expansion, and yet \emph{hanc marginis exiguitas non caperet}: in short, the category $\VCat(\C^\opp\times\D,\V)$ can be endowed with a model structure that allows to compute cofibrant resolutions $\delta\hom$ for the identity profunctor $\hom_\C$ when $\C=\D$; the coherent co/end of an endo-profunctor $T : \C\pto\C$ is then computed as the homotopy co/limit of $T$ weighted by the resolved weight $\delta\hom$, in the exact same way as the incoherent co/end is the co/limit weighted by $\hom$ (see \ref{ends-are-weighted}). At the moment of writing this book, this appears to be an interesting subject of further investigation.}
\end{itemize}
It is our sincere hope that this does not affect the outreach of this elegant and neglected piece of Mathematics, and our clumsy attempt to popularise an account of \cite{cordier1997homotopy} has to be seen as a reverent act of outreach of the branch of algebraic topology called \emph{categorical homotopy theory}.
\paragraph{\bf Local notation.} All categories $\A,\B,\dots$ appearing in the present subsection are enriched over the category $\sSet = \Cat(\bDelta^\opp , \Set)$. All such categories possess the co\fshyp{}tensors (see \ref{tenscotens}) needed to state definitions and perform computations.

Moreover, the tensor, internal hom, and cotensor functors assemble into a thc situation (see \ref{tiaccaci}) $(\otimes, \hom, \pitchfork)$ where $\otimes : \sSet \times \A \to \B$ determines the variance of the other two functors. A useful shorthand to denote the functor $\pitchfork(K,A) = K \pitchfork A$ is $A^K$. We feel free to employ such exponential notation when needed (especially when it is necessary to save space or invoke its behaviour, similar to the one of an internal hom).
Let $\B$ be a simplicially enriched category. A \emph{simplicial\hyp{}cosimplicial object} in $\B$ is a functor $Y : \bDelta^\opp \times \bDelta \to \B$.
\begin{definition}[Totalisation]\label{totalizia}\index{Totalisation}\index{_aaa_tot@$\tot$}
	Given such a simplicial\hyp{}cosimplicial object, we define its \emph{totalisation} $\tot(Y)$ as the end
	\[
		\int_{n\in\bDelta} \bDelta[n]\pitchfork Y_n
	\]
	(note that it is a cosimplicial object $m \mapsto \int_{n\in\bDelta} \bDelta[n]\pitchfork Y^n_m$). The totalisation of $Y$ is also denoted with the shorthand $\bDelta^\bullet\pitchfork Y$ or similar.
\end{definition}
Dually, a \emph{bisimplicial object} in $\B$ is a functor $X : \bDelta^\opp\times \bDelta^\opp \to \B$.
\begin{definition}[Diagonalisation]\label{diagonalizia}\index{Diagonalisation}\index{_aaa_diag@$\diag$}
	Given a bisimplicial object $X : \bDelta^\opp\times \bDelta^\opp \to \B$ we define the \emph{diagonalisation} $\diag(X)$ of $X$ to be the coend
	\[
		\int^{n\in\bDelta} \bDelta[n]\otimes X_n
	\]
	(note that it is a simplicial object $m\mapsto \int^{n\in\bDelta} \bDelta[n]\otimes X_n^m$) The diagonalisation of $X$ is also denoted with the shorthand $\bDelta^\bullet \otimes X$ or similar.
\end{definition}
Note that as a consequence of the ninja Yoneda lemma \ref{ninjayo}, the diagonalisation $\bDelta^\bullet \otimes X$ is the simplicial object $n\mapsto X_n^n$ (often denoted simply $X_{nn}$).
\begin{notat}[Chain co\fshyp{}product]\index{Chain co\fshyp{}product}
	Let $\A \in \Cat_{\bDelta}$; we shall denote as $\vec X_n = (X_0,\dots, X_n)$ the `generic $n$\hyp{}tuple of objects' in $\A$; given two other objects $A,B\in\A$, we define a bisimplicial set $\amalg\A[A| \vec X_\bullet | B]_\bullet$ whose simplicial set of $n$\hyp{}simplices is
	\index{_aaa_Aab@${\amalg\A[a\mid \firstblank \mid B]_\bullet}$}
	\[
		\amalg\A[A| \vec X_n | B]_\bullet :=
		\coprod_{X_0,\dots, X_n \in \A}
		\A(A,X_0)_\bullet\times \A(X_0, X_1)_\bullet \times \cdots \times \A(X_n,B)_\bullet.
		\footnote{It is useful to extend this notation in a straightforward way: $\A[A|\vec X|B]$ denotes the product $\A(A,X_0)\times \A(X_0, X_1) \times \cdots \times \A(X_n,B)$, and $\Pi \A[A|\vec X|B]$, $\A[\vec X]$, $\Pi \A[\vec X]$, $\amalg \A[\vec X]$ are defined similarly. Note that $\amalg\A[A| \vec X_n | B]_\bullet$ does not depend on $\vec X_n$ since the coproduct is quantified over all such $\vec X_n$'s.}
	\]
	Faces and degeneracies are induced, respectively, by composition and identity\hyp{}insertion (see Exercise \ref{ex8:Y-of-T}).

	Finally we define the simplicial set $\delta\A(A,B)$ to be the diagonalisation
	\[
		\diag (\amalg\A[A|\vec X_\bullet|B]_\bullet) = n\mapsto \amalg\A[A|\vec X_n|B]_n.
	\]
	Couched as a coend, the object $\delta\A(A,B)$ is written
	\begin{align*}
		\delta\A(A,B) & \cong \int^{n\in\bDelta} \bDelta[n]\times \amalg\A[A| \vec X_n | B]                                                                   \\
		              & = \int^{n\in\bDelta} \bDelta[n]\times \coprod_{X_0,\dots, X_n \in \A} \A(A,X_0)\times \A(X_0, X_1) \times \cdots \times \A(X_n,B)     \\
		              & \cong \int^{n\in\bDelta} \coprod_{X_0,\dots, X_n \in \A} \bDelta[n]\times \A(A,X_0)\times \A(X_0, X_1) \times \cdots \times \A(X_n,B)
	\end{align*}
\end{notat}
\begin{example}
	If we consider $\A$ to be trivially enriched over $\sSet$ (\ie as a \emph{discrete} simplicial category, where each $\A(A,A')$ is a discrete simplicial set), then the object $\A[A|\vec X_n|B]_\bullet$ coincides with the nerve of the `double slice' category $A/\A /B$ of arrows `under $A$ and above $B$'.
\end{example}
\begin{definition}[Simplicially coherent co\fshyp{}end]\label{cohcoend}\index{Co/end!simplicially coherent ---}
	Let $T : \A^\opp\times \A \to \B$ be a $\sSet$\hyp{}functor. We define
	\index{_aaa_intint@$\oint$}
	\begin{gather}
		\oint_A T(A,A) := \int_{A',A''} \delta \A(A', A'')\pitchfork T(A',A'') \notag\\
		\oint^A T(A,A) := \int^{A',A''}\delta \A(A',A'') \otimes T(A',A'')
	\end{gather}
	to be respectively the \emph{simplicially coherent end} and \emph{coend} of $T$.
\end{definition}
Expanding these definitions, we see that
\begin{gather*}
	\oint_A T(A,A) \cong \int_{A',A'',n} \bDelta[n]\times \amalg\A[A| \vec X_n | B]  \pitchfork T(A',A'') \notag\\
	\oint^A T(A,A) \cong \int^{A',A'',n} \bDelta[n]\times \amalg\A[A| \vec X_n | B]   \otimes T(A',A'')
\end{gather*}
(At this point, the reader will surely understand why performing even elementary computations with this sort of objects compels us to establish a compact notation.)
\begin{remark}
	In a few words, the definition of a simplicially coherent co\fshyp{}end mimics the classical construction, and in particular the characterisation of a co\fshyp{}end as a hom weighted co\fshyp{}limit (see \ref{ends-are-weighted} and \ref{weicolims}.\ref{coends_are_hom_weighted}), replacing hom with co\fshyp{}tensors for a `fattened up' mapping space $\A[A|\vec X|B]$.
	The `deformation' perspective is very useful, since we write that $\oint_A T$ corresponds to the end
	\[
		\int_{(A',A'')\in \A^\opp\times \A} \hom(A',A'') \pitchfork T(A',A'')
	\]
	where we applied a suitable `deformation' (or `resolution', or `replacement') functor $\delta$ to the hom functor $\A(\firstblank,\secondblank)$, seen as the identity profunctor (Remark \ref{profundefs}). This point of view is rather fruitfully explored in \cite{nashphd}, in the particular case of dg-$\Cat$ (see \ref{dg_stuff} below to draw a connection between simplicially enriched and dg\hyp{}categories, and for a precise definition of dg-category).

	This perspective is of great importance to encompass coherent co\fshyp{}ends into a general theory `compatible' with a model structure on $\VCat$, for some monoidal model $\V$ and the Bousfield-Kan model structure on $\VCat$.
\end{remark}
We now embark in the study of what shall be called \emph{simplicially coherent coend calculus}. Classical co\fshyp{}end calculus consists of the triptych Fubini - Yoneda - Kan; we shall now reproduce these steps in full detail.

The authors of \cite{cordier1997homotopy} succeed in the quite ambitious task to rewrite the most important pieces of classical category theory in this coherent model describing co\fshyp{}limits, mapping spaces, the Yoneda lemma, and Kan extensions. The aim of the rest of this subsection is to sketch some of these original definitions, hopefully helping one of the alternative approaches to $(\infty,1)$\hyp{}category theory to escape oblivion.
\begin{definition}[The functors $\Y$ and $\W$]\index{_aaa_Y@$\Y$}\index{_aaa_W@$\W$}\label{Y_and_W}
	Let $T : \A^\opp\times \A \to \B$ be a functor; we define $\Y(T)^\bullet$ to be the cosimplicial object (in $\B$)
	\[
		\Y(T)^n := \prod_{\vec X_n = (X_0,\dots, X_n)} \A[\vec X_n]\pitchfork T(X_0, X_n)
	\]
	where $\A[\vec X_n] = \A(X_0,X_1) \times \dots \times \A(X_{n-1},X_n)$.

	Dually, given the same $T$, we define $\W(T)_n$ to be the simplicial object (in $\B$)
	\[
		\W(T)_n := \coprod_{\vec X_n = (X_0,\dots, X_n)} \A[\vec X_n] \otimes T(X_0, X_n).
	\]
\end{definition}
\index{_aaa_intint@$\oint$}
\begin{proposition}
	Let $T : \A^\opp\times \A \to \B$ be a $\sSet$\hyp{}functor. Then there is a canonical isomorphism
	\[
		\oint_A T(A,A) \cong \tot(\Y(T)^\bullet)
	\]
\end{proposition}
\begin{proof*}
	We make heavy use of the ninja Yoneda lemma \ref{ninjayo} in its enriched form, \ie of the fact that given a $\sSet$\hyp{}functor $F : \A \to \sSet$ we have a canonical isomorphism
	\[
		\int_X \A(X,B) \pitchfork F(X) \cong F(B)
	\]
	and the fact that $K \pitchfork (H \pitchfork A)\cong (K\otimes H)\pitchfork A$, naturally in all arguments.

	Now, let $\vec X_n = (X_0,\dots, X_n)$ be a generic tuple of objects of $\A$: to save some space we employ the exponential notation $A^K$ to denote the cotensor $K\pitchfork A$.
	\begin{align*}
		\oint_A T(A,A) & := \int_{A',A''} T(A',A'')^{\delta\A(A',A'')}                                                                               \\
		               & \cong \int_{A',A''} 		 T(A',A'')^{\int_n \bDelta[n]\times \amalg \A[A'|\vec X| A'']}                                          \\
		               & \cong \int_{A',A'',n}		 T(A',A'')^{\bDelta[n]\times \amalg \A[A'|\vec X| A'']}                                                \\
		               & \cong \int_{A',A'',n} \left( T(A',A'')^{\A(X_n,A'')} \right)^{\bDelta[n]\times \amalg \A[A'|\vec X]}                        \\
		               & \cong \int_{A',A'',n} \prod_{X_0,\dots, X_n} \left( T(A',A'')^{\A(X_n,A'')} \right)^{\bDelta[n]\times \A[A'|\vec X]}        \\
		               & \cong \int_{A',n} \prod_{X_0,\dots, X_n} \left( \int_{A''} T(A',A'')^{\A(X_n,A'')} \right)^{\bDelta[n]\times \A[A'|\vec X]} \\
		               & \cong \int_{A',n} \prod_{X_0,\dots, X_n} T(A',X_n)^{\bDelta[n]\times \A[A'|\vec X]}                                         \\
		               & \cong \int_{n} \prod_{X_0,\dots, X_n} \left( \int_{A'} T(A',X_n)^{\A(A',X_0)} \right)^{\bDelta[n]\times \A[\vec X]}         \\
		               & \cong \int_{n} \Big(\prod_{X_0,\dots, X_n} T(X_0,X_n)^{\A[\vec X_n]}\Big)^{\bDelta[n]} \cong \tot(\Y(T)).
	\end{align*}
\end{proof*}
For the sake of completeness, we notice that the universal wedge testifying that $\oint_A T(A,A) \cong \tot(\Y(T))$ is induced by the morphisms
\[\vcenter{\xymatrix@R=4mm{
	{\displaystyle \oint_A T(A,A) = \int_{A',A''} \delta\A(A',A'')\pitchfork T(A', A'')} \ar[d]^\wr\\
	{\displaystyle \int_{A',A'',n}\big(\Delta[n]\times \amalg \A[A'|\vec X| A'']\big)\pitchfork T(A',A'')} \ar[d]\\
	{\displaystyle \Delta[n]\pitchfork \Y(T)^n.}
	}}\]
Prove the dual statement as an exercise (it is of vital importance that you establish a good notation):
\index{_aaa_intint@$\oint$}
\begin{proposition}
	Let $T : \A^\opp\times \A \to \B$ be a $\sSet$\hyp{}functor. Then there is a canonical isomorphism
	\[
		\oint^A T(A,A) \cong \diag(\W(T)_\bullet)
	\]
\end{proposition}
\begin{remark}
	The homotopy coherent co\fshyp{}ends admit `comparison' maps to the classical co\fshyp{}ends, induced by the fact that the `fattened hom' $\delta\A(\firstblank,\secondblank)$ has canonical maps to/from the plain hom $\A(\firstblank,\secondblank)$; this is part of a general tenet, where homotopically correct objects result as a \emph{deformation} of classical ones, and the deformations maps in/out of the plain object.

	The comparison map $\oint T(A,A) \to \int T(A,A)$ arises, here, as an homotopy equivalence between the simplicial set $\A(A,B)$ (seen as bisimplicial, and constant in one direction) and the bisimplicial set $\delta \A(A,B) = \diag \A[A|\bullet|B]_\bullet$: this is \cite[p. 15]{cordier1997homotopy}.

	Note that it is possible to write an explicit contracting homotopy between the two objects. The map
	\[\textstyle
		d_0 : \coprod_{X_0} \A(A, X_0)\times \A(X_0,B) \to \A(A,B)
	\]
	given by composition has an homotopy inverse given by
	\[
		s_{-1} : \A(A,B) \to \A(A, A)\times \A(A,B) : g\mapsto (\id_A, g).
	\]
	Indeed, the composition $d_0 s_{-1}$ is the identity on $\A(A,B)$, whereas the composition $s_{-1} d_0$ admits is homotopic to the identity on $\delta \A(A,B)$ (we use the same name for the maps $d_0, s_{-1}$ and the induced maps $\bar d_0 : \delta \A(A,B) \to \A(A,B)$, induced by the universal property, and $\bar s_{-1} : \A(A,B) \to \delta \A(A,B)$).

	There is an important difference between these two maps, though: while $d_0$ is natural in both arguments, $s_1$ is natural in $B$ but not in $A$. This has an immediate drawback: while $d_0$ can be obtained canonically, as the universal arrow associated to a certain natural isomorphism (see (\ref{dizzero}) below), $s_{-1}$ can't (the best we can do is to characterise the natural argument of $s_{-1}$ via \cite[Example 2, p. 16]{cordier1997homotopy}).
\end{remark}
\subsubsection{Simplicially coherent natural transformations}
\index{Natural transformation!simplicially coherent ---}
As we have seen in \ref{naturalu}, the set of natural transformations between two functors $F,G : \C \to \D$ coincides with the end $\int_X \D(FX, GX)$, and (see \ref{laxnat}) the category of lax natural transformations between two 2\hyp{}functors coincides with the lax end $\sqint_X \D(FX, GX)$. It comes as no surprise, then, that the following characterisation of \emph{homotopy coherent} natural transformations between two simplicial functors hold:
\begin{definition}[Coherent natural transformations]\index{Co/end!natural transformations as ---s}\index{Natural transformation!---s as ends}
	Let $F,G : \C \to \D$ be two simplicial functors; then the simplicial set of \emph{coherent transformations} between $F$ and $G$ is defined as
	\[
		\Cat_{\bDelta}(\C,\D)(\!( F,G )\!) := \oint_A \D(FA, GA).
	\]
\end{definition}
We define the following operations of \emph{coherent tensoring} and \emph{cotensoring} a simplicial functor with a representable:
\def\opitchfork{\mathop{\overline{\pitchfork}}}
\def\uotimes{\mathop{\underline{\otimes}}}
\begin{definition}[Mean tensor and cotensor]
	Let $F : \A \to \B$, $G : \A \to \sSet$, $H : \A^\opp \to \sSet$.
	We define $G \opitchfork F$, $H \uotimes F$ respectively as
	\[
		G \opitchfork F := \oint_A GA\pitchfork FA \qquad\qquad
		H \uotimes F := \oint^A HA \otimes FA.
	\]
\end{definition}
This yields the notion of \emph{standard resolution} of a simplicial functor:
\begin{definition}[Standard resolutions]\label{standard_res}
	Let $F : \A \to \B$ be a simplicial functor; we define
	\begin{align*}
		\overline{F}A  & := \A(A,\firstblank)\opitchfork F = \oint_X \A(A,X)\pitchfork FX \\
		\underline{F}A & := \A(\firstblank,A)\uotimes F = \oint^X \A(X,A)\otimes FX.
	\end{align*}
\end{definition}
\begin{example}\label{thehoms}
	We specialise the above definition to compute the functors
	$\overline{\hom}(A,\firstblank)$ and $\underline{\hom}(A,\firstblank)$: in particular we concentrate on the second case, since the first is completely dual.
	\begin{align*}
		\underline{\hom}(A,B) & = \oint^A \A(A,X) \times \A(X,B)                                           \\
		                      & \cong \int^{XY} \delta \A(X,Y) \times \A(A,X)\times \A(Y,B)                \\
		                      & \cong \int^{XY} \A(A,X)\times \delta\A(X,Y) \times \A(Y,B)                 \\
		                      & \cong \int^{XYn} \A[A|\tilde{X}_n|B]\times\bDelta[n] \cong \delta \A(A,B).
	\end{align*}
\end{example}
We leave as an easy exercise in co\fshyp{}end calculus the proof of the following result (see Exercise \ref{ex8:cohnat}), which shows that the standard resolutions $\underline{F}, \overline{F}$ of $F$ `absorb the coherence':
\begin{proposition}\label{absorb}
	Let $F,G : \C \to \D$ be two simplicial functors; then there are canonical isomorphisms
	\[
		[F, \overline G] \cong \Cat_{\bDelta}(\C,\D)(\!( F,G )\!) \cong [\underline F,G].
	\]
\end{proposition}
This result has a number of pleasant consequences: the simplicially coherent setting is powerful enough to retrieve several classical constructions.
\begin{itemize}
	\item Example \ref{thehoms} above shows that $\underline{\hom}(A,\firstblank)(b)\cong \delta \A(A,B)$; this entails that there is an isomorphism
	      \[\label{dizzero}
		      \Cat_{\bDelta}(\A,\sSet)(\delta \A(A,\firstblank), \A(A, \secondblank)) \cong
		      \Cat_{\bDelta}(\A,\sSet)(\!( \A(A,\firstblank ),\A(A,\secondblank))\!)
	      \]
	      and it is a matter of verifying some additional nonsense to see that the $\sSet$\hyp{}natural transformation corresponding to the identity coherent transformation is precisely $d_0$.
	\item The map $d_0$ defines additional universal maps $\eta_F , \eta^F$ which `resolve' a functor $F : \A \to \B$ whenever $\underline{F}, \overline{F}$ exist (it is sufficient that $\B$ admits all the relevant co\fshyp{}limits to perform the construction of $\underline{F}, \overline{F}$). From the chain of isomorphisms
	      \begin{align*}
		      \eta^F : \overline{F}B  & = \oint_A \A(B,A)\pitchfork FA                                        \\
		                              & \cong \int_{A',A''} \delta\A(A',A'')\pitchfork\A(B,A'')\pitchfork FA' \\
		                              & \leftarrow \int_{A',A''} \A(A',A'')\pitchfork\A(B,A'')\pitchfork FA'  \\
		      (\ref{ninjayo})         & \cong FB;                                                             \\
		      \eta_F : \underline{F}B & = \oint^A FA \otimes \A(A,B)                                          \\
		                              & \cong \int^{A',A''} FA' \otimes \A(A'',B)\delta\A(A',A'')             \\
		                              & \to \int^{A',A''} FA' \otimes \A(A'',B)\A(A',A'')                     \\
		                              & \cong FB;
	      \end{align*}
	      we obtain natural transformations corresponding to suitable coherent identities under the isomorphism of \ref{absorb}.
	\item The maps $\eta_F , \eta^F$ behave like resolutions: \cite[3.4]{cordier1997homotopy} shows that they are level\hyp{}wise homotopy equivalences (meaning that $\eta_F : FA \to \overline{F}A$ induces homotopy equivalences of simplicial sets $\B(B, FA)\xto{(\eta_F)_*} \B(B, \overline{F}A)$ for each $B$, naturally in $B$).\footnote{We decide to skip the proof of this proposition, as it is quite long, technical, and even though it relies on co\fshyp{}end calculus it doesn't add much to the present discussion.}
\end{itemize}
\subsubsection{Simplicially coherent Kan extensions}\index{Kan extension!simplicially coherent ---}
The universal property of a Kan extension is inherently 2\hyp{}dimensional: uniqueness is stated at the level of 2\hyp{}cells, and any sensible generalisation of it to the higher world involves a `space' of 2\hyp{}cells between 1\hyp{}cells.

This entails that any reasonable definition of a (left or right) Kan extension ultimately relies on a nice definition for a space of coherent natural transformations between functors, which has been the subject of the previous subsection.

There are, nevertheless, several subtleties as there are many choices available for a definition: in the words of \cite{cordier1997homotopy},
\begin{quote}
	Clearly one can replace natural transformations by \emph{coherent} ones [in the definition of a Kan extension], but should isomorphism be replaced by homotopy equivalence, should they be natural, in which direction should this go\dots ?
\end{quote}
Solving this problem can be tricky; one of the reasons is that simplicial combinatorics captures really well the behaviour of $(\infty,1)$\hyp{}categories, whereas any satisfactory model for homotopy coherent Kan extensions shall speak about $(\infty,2)$\hyp{}categories.

In order to define coherent Kan extensions for $\B  \xot{G} \A\xto{F} \C$ we ask the isomorphisms
\begin{gather}
	\Cat_{\bDelta}(\B,\C)(\!( H,\hoRan_G F)\!)\cong \Cat_{\bDelta}(\A,\C)(\!( HG,F )\!)\label{kan:chistu}\\
	\Cat_{\bDelta}(\B,\C)(\!(\hoLan_G F,K)\!)\cong \Cat_{\bDelta}(\A,\C)(\!( F,KG )\!)\label{kan:chiddu}
\end{gather}
to hold at the level of \emph{coherent} transformations. This can be achieved as follows:
\begin{definition}[Coherent Kan extensions]\label{cohkan}
	Let $F : \A \to \C$ and $G : \A \to \B$ be a span of simplicial functors; we define
	\begin{gather*}
		\hoRan_G F(\firstblank) = \oint_A \B(\firstblank,GA)\pitchfork FA\\
		\hoLan_G F (\firstblank) = \oint^A \B(GA,\firstblank)\otimes FA
	\end{gather*}
\end{definition}
Proving that the isomorphisms \eqref{kan:chistu} and \eqref{kan:chiddu} hold follows from an easy computation with the explicit form of the coherent co/ends above. We leave to the reader to either solve this as an exercise in \ref{ex8:cohkan}, or to consult \cite{cordier1997homotopy}.
\paragraph{A glance at dg\hyp{}coends}\label{dg_stuff}
The Dold-Kan correspondence (see \ref{doldekanne}) establishes an equivalence of categories between $\Ch_{\ge}(\Z)$, chain complexes of abelian groups concentrated in positive degree, and simplicial objects in the category of abelian groups, \ie functors $G : \bDelta^\opp\to \Mod(\Z)$. The equivalence of categories is generated by a cosimplicial object
\[
	DK : \bDelta \to \Ch_{\ge}(\Z)
\]
and the universal property of the Yoneda embedding now yields an adjunction
\[
	\Lan_\yon(DK) : \sSet \leftrightarrows \Ch_{\ge}(\Z) : \Lan_{DK}\yon
\]
It turns out that this is an equivalence of categories.

(This result can be restated in fair more generality, but in this paragraph we stick to the $R$\hyp{}linear case.)

Dold-Kan equivalence induces an equivalence of 2\hyp{}categories between simplicial\hyp{}abelian\hyp{}group enriched categories on one side, and categories enriched in (positive) chain complexes on the other, or suitable \emph{dg\hyp{}categories} (\ie categories enriched in chain complexes) for short.

The theory of dg\hyp{}categories is deeply rooted in homological algebra and finds applications in algebraic geometry \cite{kuznetsov2014categorical,nashphd}: in order to study nicer versions of derived categories, one can attach a \emph{derived dg\hyp{}category} $\underline{\sD}(X)$ to a space/scheme/abelian category, and such a category is way better behaved than the `incoherent' derived category $\sD(X)$ (where the cone construction of \ref{cone-is-a-colim} is not functorial).

There is of course a link between the two objects; the incoherent category can be recovered as the homotopy category of $\underline{\sD}(X)$ as follows: each cohomology functor $H^n$ extend degree\hyp{}wise to functors $\underline{H}^n : \text{dg-}\Cat \to \Ab\text{-}\Cat$ (=categories enriched on abelian groups, also called \emph{preadditive}), and there is an isomorphism $\sD(X) = \underline{H}^0(\underline{\sD}(X))$.

In order to approach derived algebraic geometry with the tools of enriched category theory, it might be interesting to restrict coherent co\fshyp{}end calculus to $[\bDelta^\opp,\Ab]$\hyp{}categories, and think about the result as dg\hyp{}categories using the Dold-Kan equivalence.

This is an enticing application of co\fshyp{}end calculus in an homotopical/homological setting, and many questions arise naturally from the expressive power of co\fshyp{}end calculus. For example: if $\A$ is any dg\hyp{}category its identity profunctor $\A \pto \A$ is a functor $\A^\opp\boxtimes \A \to \Ch(\Z)$, so that the coherent end
\[
	\oint_A \A(A,A)
\] \ie the object of derived natural transformations of the identity functor $\id_{\A}$, recovers the \emph{Hochschild complex} of $\A$. Then, if $\A$ is an associative algebra regarded as a one\hyp{}object dg\hyp{}category concentrated in degree zero, the object $H^n(\int_* A)$ is the \emph{Hochschild cohomology} of $A$, understood in the classical sense of, say, \cite[Ch. 11]{pierce1982associative}.

We don't have enough space to expand on this interesting topic, but one can find that the bicategory of dg\hyp{}profunctors, and in particular of the \emph{endoprofunctors} of a single object $\cX$ gives rise to plenty of derived invariants of a dg\hyp{}category, but we leave to the interested reader the endeavor of re\hyp{}reading the paper \cite{kuznetsov2014categorical} wearing appropriate co\fshyp{}end\hyp{}goggles (among many, that paper seems the most liable to a co\fshyp{}end\hyp{}theoretic reformulation).

The same approach can be carried over in the more general setting of a category enriched over a monoidal model category; a perfect starting point accounting for the state of the art on the matter is Shulman's \cite{shulman}.
\section{Co/ends in quasicategories}\label{coend-in-qcat}
\index{Co/end!$\infty$- ---}
As a rule of thumb, the translation procedure from category to $\infty$\hyp{}category theory is based on the following meta\hyp{}principle: first you rephrase the old definition in a `simplicially meaningful' way, so that the $\infty$\hyp{}categorical definition specialises to the old one for quasicategories $N(\C)$ which arise as nerves of categories. Then you forget about the original gadget and keep the simplicial one; this turns out to be the right definition.

The first victim of this procedure is the twisted arrow category \ref{twisted} of an $\infty$\hyp{}category.
\begin{definition}[Twisted arrow $\infty$\hyp{}category]
	\index{Category!twisted arrows $\infty$- ---}
	Let $\varepsilon : \bDelta\to \bDelta$ be the functor $[n]\mapsto [n]\star [n]^\opp$, where $\star$ is the join of simplicial sets \cite{Joy,ehlers2008ordinal}, and the \emph{opposite} of a simplicial set is defined in \cite[\S 6.19]{rezk2017stuff}. Let $\C$ be a $\infty$\hyp{}category; the twisted arrow category $\tw(\C)$ is defined to be the simplicial set $\varepsilon^*\C$, where $\varepsilon^*=\firstblank\circ \varepsilon^\opp : \sSet \to \sSet$ is the induced precomposition functor. More explicitly, and consequently, the $n$\hyp{}simplices of $\tw(\C)$ are characterised by the relation
	\[
		\tw(\C)_n \cong \sSet(\bDelta[n], \tw(\C))\cong \sSet(\bDelta[n]\star \bDelta[n]^\opp, \C).
	\]
\end{definition}
The most important feature of the twisted arrow category is that it admits a fibration over $\C^\opp\times\C$ (part of its essential properties can be deduced from this); the machinery of left and right fibrations exposed in \cite[2.0.0.3]{HTT} gives that
\begin{itemize}
	\item There is a canonical simplicial map $\Sigma : \tw(\C) \to \C^\opp\times\C$ (induced by the two join inclusions $\bDelta[\firstblank], \bDelta[\firstblank]^\opp \to \bDelta[\firstblank]\star \bDelta[\firstblank]^\opp$);
	\item This $\infty$\hyp{}functor is a right fibration in the sense of \cite[2.0.0.3]{HTT}.
\end{itemize}
\begin{remark}
	It is rather easy to see that the above definition is reasonable: a 0\hyp{}simplex in $\tw(\C)$ is an edge $f : \bDelta[1]\to \C$, and a 1\hyp{}simplex of $\tw(\C)$ is a 3\hyp{}simplex thereof, that we can depict as a pair of edges $(u,v)$, such that the square having twisted edges
	\[
		\xymatrix{
			\ar[d]_f & \ar[l]_u \ar[d]^{f'}\\
			\ar[r]_v &
		}
	\]
	commutes. This suggest (as it must be) that the definition of $\tw(\C)$ for a $\infty$\hyp{}category specialises to the 1\hyp{}dimensional one and adds higher\hyp{}dimensional information to it.
\end{remark}
\begin{definition}\label{infend}
	\index{Co/end!$\infty$- ---}
	Let $\C,\D$ be two $\infty$\hyp{}categories; the $\infty$-\emph{co}/\emph{end} of a simplicial map $F : \C^\opp\times \C\to \D$ is the co\fshyp{}limit $\oint F$ of the composition
	\[
		\tw(\C)\xto{\Sigma} \C^\opp\times \C \xto{F} \D
	\]
\end{definition}
The main interest of the authors in \cite{gepner2015lax} is to formulate an analogue of \ref{fibelem}, which characterises the Grothendieck construction of a $\Cat$\hyp{}valued functor as a particular weighted colimit (see \ref{elts-as-coend}).

It is rather easy to formulate such an analogue definition: this appears as \cite[ 2.8]{gepner2015lax}.
\begin{definition}[op/lax colimit of $F$]\index{Colimit!op/lax ---}
	Let $F : \C\to \Cat_\infty$ be a functor between $\infty$\hyp{}categories. We define
	\begin{itemize}
		\item the \emph{slice fibration} for $\C \in\BF{QCat}$ to be the functor of quasicategories $\chi_\C : \C \to \BF{QCat}$ sending $C\in\C$ to $\C/C$, and dually the \emph{coslice fibration} to be $\chi^\C : \C \to \BF{QCat} : C \mapsto C/\C$;
		\item the \emph{lax colimit} of $F$ to be the coend
		      \[
			      \infint^C C/\C\times FC;
		      \]
		\item the \emph{oplax colimit} of $F$ to be the coend
		      \[
			      \infint^C \C/C\times FC.
		      \]
	\end{itemize}
\end{definition}
The Grothendieck construction associated to $F$, discussed in \cite{HTT} with the formalism of un/straightening functors results precisely as the oplax colimit of $F$. This is concordant with our \ref{elts-as-coend} and \ref{its.another.nerve}.
\subsubsection{A Fubini rule for $\infty$\hyp{}coends}\label{fubi_for_infty}
\index{Fubini rule!--- for $\infty$\hyp{}coends}
The following section establishes the analogue of \ref{fubozzo} for $\infty$\hyp{}coends; the present proof, as it stands, currently appears as \cite{infub}.

We freely employ the terminology on $\infty$\hyp{}category theory recalled in §\ref{higher_caz}. In particular we denote $\BF{Kan}$ the category of \emph{Kan complexes} \ie simplicial sets that lie in the orthogonal of all horn inclusions $\Lambda^n_k \hookrightarrow \Delta[n]$, for $0\le k\le n$ and $n\ge 0$; the \emph{nerve} functor $N : \Cat \to \sSet$ establishes a Quillen equivalence between the category of categories and the category of simplicial sets; fibrant objects in the latter category are the \emph{$\infty$\hyp{}categories} of \cite{HTT}.
\begin{lemma}\label{adjunccc}
	Let $\C$ be a small $\infty$\hyp{}category, and $\D$ be a presentable $\infty$\hyp{}category; then $\D$ is tensored and cotensored over $\catS = N(\BF{Kan})$ (the \emph{$\infty$\hyp{}category of spaces}). This entails that there is a two\hyp{}variable adjunction
	\[
		\D^\opp \times \D \xto{\Map_\D} \catS
		\qquad
		\catS \times \D \xto{\otimes} \D
		\qquad
		\catS^\opp\times\D \xto{\pitchfork}\D
	\]
	such that
	\[ \D(X\otimes D,D')\cong \catS(X, \Map_{\D}(D,D'))\cong \D(D, X\pitchfork D') \]
\end{lemma}
From the existence of these isomorphisms it is clear that
\begin{align}
	V\otimes(W\otimes D) \cong W\otimes(V\otimes D) \cong (V\times W)\otimes D\label{iso:fst} \\
	V\pitchfork(W\pitchfork D) \cong W\pitchfork(V\pitchfork D)\cong (V\times W)\pitchfork D\label{iso:snd}
\end{align}
\begin{lemma}\label{end-is-functor}
	Let $F : \C^\opp\times\C\to \D$ be a $\infty$\hyp{}functor and $\C,\D$ $\infty$\hyp{}categories as in the assumptions of \ref{adjunccc}. Then
	\begin{itemize}
		\item $F\mapsto \int^C F$ is functorial, and it is a left adjoint;
		\item $F\mapsto \int_C F$ is functorial, and it is a right adjoint.
	\end{itemize}
\end{lemma}
\begin{proof}
	We only prove the first statement for coends; the other one is dual.

	Since $\int^C F = \colim_{\tw(\C)}(F\circ \Sigma) = \colim_{\tw(\C)} \circ \Sigma^*(F)$ results as a composition of $\infty$\hyp{}functors, it is clearly functorial; then
	\[\textstyle
		\int^C :
		\xymatrix@C=2cm{
		[\C^\opp\times\C,\D] \ar@<5pt>[r]^-{\Sigma^*} &\ar@<5pt>@{.>}[l]^-{\Ran_\Sigma} \ar@{}[l]|-\perp[\tw(\C),\D] \ar@<5pt>[r]^-{\colim_{\tw(\C)}} &\ar@<5pt>@{.>}[l]^-{c} \ar@{}[l]|-\perp\D\\
		}
	\]
	is a left adjoint because it is a composition of left adjoints ($c = t^*$ is the constant functor inverse image of the terminal map $\tw(\C)\to *$).

	Dually, the left adjoint to $\int_C $ is given by $\Lan_\Sigma \circ c(D)$.
\end{proof}
Now, the Fubini rule asserts that when the domain of a functor $F : \A^\opp\times\A\to \D$ is of the form $(\C\times\E)^\opp\times(\C\times\E)$, then the co\fshyp{}ends of $F$ can be computed as `iterated integrals'
\begin{gather}
	\int^{(C,E)} F \cong \int^{CE} F \cong \int^{EC} F\label{coends}\\
	\int_{(C,E)} F \cong \int_{CE} F \cong \int_{EC} F\label{ends}
\end{gather}
These identifications hide a slight abuse of notation, that is worth to make explicit in order to avoid confusion: thanks to \ref{end-is-functor} the three objects of \eqref{coends} can be thought as images of $F$ along certain functors, and the Fubini rule asserts that they are linked by canonical isomorphisms; we can easily turn these functors into having the same type by means of the cartesian closed structure of $\sSet$ (somewhat sloppily, we denote the internal hom of $\E,\D$ as $[\E,\D]$):
\[
	\vcenter{\xymatrix@C=0cm{
	[\C^\opp\times\C\times\E^\opp\times\E,\D] \ar@{=}[d]& [\C^\opp\times\C\times\E^\opp\times\E,\D] \ar@{=}[d]&
	[\C^\opp\times\C\times\E^\opp\times\E,\D] \ar@{=}[d]\\
	[\C^\opp\times\C, [\E^\opp\times\E,\D]] \ar[d]_{[\C^\opp\times\C,\int^E]}&
	[\E^\opp\times\E, [\C^\opp\times\C,\D]] \ar[d]^{[\E^\opp\times\E,\int^C]}&
	[(\C\times\E)^\opp\times(\C\times\E),\D]\ar[dd]^{\int^{(C,E)}} \\
	[\C^\opp\times\C,\D] \ar[d]_{\int^C}&
	[\E^\opp\times \E,\D] \ar[d]^{\int^E}& \\
	\D & \D & \D
	}}
\]
(of course, we can provide similar definitions for the iterated end functor).

Once that this has been clarified, we can deduce the isomorphisms \eqref{coends} and \eqref{ends} from the fact that the three functors $\int^{CE},\int^{EC},\int^{(C,E)}$ have right adjoints isomorphic to each other, and hence they must be isomorphic themselves.

This argument evidently mimics the one given in \ref{fubozzo}.
\begin{proposition}
	The functor $R = \Ran_\Sigma(c(\_))$ acts `cotensoring with mapping space':
	\[D\mapsto \Big((C,C')\mapsto \Map_{\C}(C,C')\pitchfork D\Big)\]
	Dually, the functor $L = \Lan_\Sigma(c(\_))$ acts `tensoring with mapping space':
	\[D\mapsto \Big((C,C')\mapsto \Map_{\C}(C,C')\otimes D\Big)\]
\end{proposition}
\begin{proof}
	We only prove the first statement about $R$; the other one is dual.

	It turns out that the statement relies on very well\hyp{}known features of simplicial categories: we move to that setting, lacking an equally simple and conceptual quasicategorical proof.

	Recall from \ref{cohcoend} that translating this result into the simplicially\hyp{}enriched setting, we get the \emph{coherent end} of a simplicial functor $F : \A\to \underline{\B}$ of $\BF{Kan}$\hyp{}tensored and cotensored simplicial categories as the end
	\index{_aaa_intint@$\oint$}
	\[\oint^A F(A,A) := \int^{A',A''}\delta \A[A'|A''] \otimes F(A',A'')\]
	defined in \ref{cohcoend}; the object $\delta\A[X|Y]$ is a cofibrant resolution of the hom functor, regarded as the identity $\hom : \A\pto\A$. In view of this, and since the weighted colimit functor
	\[ {} \firstblank \otimes\secondblank : [\cX^\opp,\BF{Kan}]\times[\cX,\B]\to \B : (G,F)\mapsto \wcolim{G}{F}\]
	is functorial in its $F$ argument, is also a left adjoint with right adjoint the cotensoring with the weight:
	\[
		\D(W\otimes F, D) \cong [\cX,\D](F, W\pitchfork D).
	\]
	Since $\oint^A F(A,A)\cong \delta\A \otimes F$ is the weighted colimit with $\delta\A$ as a weight, it turns out that there is an adjunction
	\[\textstyle
		\oint^A \cong \delta\A \otimes\firstblank \quad\dashv\quad \lambda D.\lambda XY.\delta\A[X|Y] \pitchfork D
	\]
	and by the uniqueness of adjoint functors (of course valid also in this setting) we obtain that $\Ran_\Sigma \circ c(D)\cong \delta\A \pitchfork D$.
\end{proof}
The Fubini rule now follows from uniqueness of adjoints: in diagram
\[
	\vcenter{\xymatrix@R=3mm{
	\lambda F . \int^C\kern-.35em\int^E F \ar@{-|}[r] & \lambda D. \lambda CC'.\lambda EE'.\Map(C,C')\pitchfork \Big(\Map(E,E')\pitchfork D\Big)\ar@{=}[d]^\wr\\
	\lambda F . \int^E\kern-.35em\int^C F \ar@{-|}[r] & \lambda D. \lambda EE'.\lambda CC'.\Map(E,E')\pitchfork \Big(\Map(C,C')\pitchfork D\Big)\ar@{=}[d]^\wr\\
	\lambda F . \int^{(C,E)} F \ar@{-|}[r] & \lambda D.\lambda CEC'E' \Big(\Map(C,C')\times\Map(E,E')\Big)\pitchfork D\ar@{=}[d]^\wr\\
	& \Map\big((C,E),(C',E')\big)\pitchfork D
	}}
\]
the vertical isomorphisms on the right are justified by \eqref{iso:snd}. A completely analogous argument, using \eqref{iso:fst} instead, and the left adjoints given by tensoring with the derived mapping space, gives the Fubini rule for \eqref{ends}.
\section{Co/ends in a derivator}\index{Co/end!--- in a derivator}
The theory of derivators provides a purely 2\hyp{}categorical model for higher category theory, where all the coherence information is encoded in conditions that are imposed on suitable diagrams of 2\hyp{}cells.

The theory was invented by A. Grothendieck in order to address the many shortcomings of triangulated categories, a categorical structure naturally arising in stable homotopy theory.

Here we only sketch some of the basic definitions needed to pave the way to \ref{coendinder} below.
\begin{definition}[The 2\hyp{}category of prederivators]\label{preders} We define a 2-category $\mathsf{PDer}$ such that
	\begin{itemize}
		\item an object of $\mathsf{PDer}$, called \emph{prederivator}, is a strict 2\hyp{}functor $\sD : \caat^\opp\to \tCat$;
		\item a \emph{morphism of prederivators} is a pseudonatural transformation between pseudofunctors, $\eta : \sD\To \sD'$;
		\item a \emph{2\hyp{}cell} between morphisms of prederivators is a modification (see \ref{modification}) $\Theta : \eta\Rrightarrow \eta'$ between pseudonatural transformations.
	\end{itemize}
\end{definition}
\begin{notation}
	Along this paragraph we employ a local notation which is specific of the literature in derivator theory: the terminal category is often denoted $e$ (perhaps French for `\emph{e}nsemble avec un seul \emph{é}lément'); small categories are in Roman uppercase $I,J,K\dots$ and functors are in Roman lowecase $u : I \to J$. If $u : I \to J$, then its image along a prederivator $\sD$ is denoted $u^* : \sD(J) \to \sD(I)$. A morphism of prederivators has components $F_I : \sD(I)\to \sD'(I)$ and it is specified by families of invertible 2\hyp{}cells $\gamma_u : F_I\circ u^* \To u^* \circ F_J$ subject to suitable coherence conditions.
\end{notation}
The notion of derivator is a refinement of \ref{preders}, originally motivated by the desire to provide a satisfactory axiomatisation for triangulated categories --and more generally, homotopy categories of model categories-- that only appeals to 2\hyp{}categorical language.

A \emph{derivator} is then a prederivator that satisfies the following additional conditions (we adopt the same labeling convention of \cite{groth2013derivators}):
\begin{enumerate}[label=$\text{Der}\oldstylenums{\arabic*})$, ref=$(\text{Der}\oldstylenums{\arabic*})$]
	\item \label{derax:1} The functor $\sD(I\sqcup J)\to \sD(I)\times\sD(J)$ obtained from the canonical inclusions $i_I : I\to I\sqcup J\leftarrow J : i_J$ is an equivalence.
	\item \label{derax:2} Each object $j : e\to J$ induces a family of functors $\sD(J)\xto{j^*} \sD(e)$; we ask that this family \emph{jointly reflects isomorphisms}, \ie a morphism $f\in \sD(J)$ is invertible if and only if each $j^*f$ is invertible in $\sD(e)$.
	\item \label{derax:3} Each functor $u^* : \sD(J)\to \sD(I)$ induced by $u : I\to J$ admits both a left adjoint $u_!$ and a right adjoint $u_*$. These functors are called, respectively, the \emph{homotopy left Kan extension} and \emph{homotopy right Kan extension} along $u$.
	\item \label{derax:4} Given a functor $u : J\to K$, and the two squares in the left column below, there exist two squares in $\tCat$, in the right column below, induced by the colax pullbacks defining the slice and coslice categories:
	      \begin{align*}
		      \vcenter{\xymatrix{
		      J_{\!/k} \ar[r]^t\ar[d]_p                                  & e\ar[d]^k                         \\
		      J \ar[r]_u                                                 & K
				      \ultwocell<\omit>{\varpi}
			      }}
		      \quad
		      \rightsquigarrow
		      \quad
		      \vcenter{\xymatrix{
		      \sD(J_{\!/k})\ar[r]^{t_*}                                  & \sD(e)\dltwocell<\omit>{\varpi_*} \\
		      \sD(J)\ar[u]^{p^*} \ar[r]_{u_*}                            & \sD(K)\ar[u]_{k^*}
		      }}                                                                                             \\
		      \vcenter{\xymatrix{
		      J_{k/} \drtwocell<\omit>{\varpi'} \ar[r]^t\ar[d]_p         & e\ar[d]^k                         \\
		      J \ar[r]_u                                                 & K
			      }}
		      \quad
		      \rightsquigarrow
		      \quad
		      \vcenter{\xymatrix{
		      \sD(J_{k/})\ar[r]^{t_!}                                    & \sD(e)                            \\
		      \sD(J)\urtwocell<\omit>{\varpi_!}\ar[u]^{p^*} \ar[r]_{u_!} & \sD(K).\ar[u]_{k^*}
			      }}
	      \end{align*}
	      We ask that these squares are filled by invertible 2\hyp{}cells $\varpi'_! : t_!p^*\To k^* u_!$, and $\varpi_* : k^*u_* \To t_* p^*$.
\end{enumerate}
\begin{remark}
	Taken all together, the axioms of derivator are meant to ensure that we can build a theory which is expressive enough for applications: te fundamental idea is that we can do category theory `over $\caat$' as a family of large categories $\sD(J)$ contravariantly depending on functors $u : I \to J$.

	Let's take a deeper look at the axioms: \ref{derax:1} asks $\sD$ to act independently on different (finite) connected components; \ref{derax:2} asks that being an isomorphism in $\sD(J)$ is a `local' notion; the other two axioms are meant to express the fact that we can compute left and right Kan extensions for every functor $u : I\to J$ \ref{derax:3}, and that these extensions are pointwise \ref{derax:4}. More precisely, axiom \ref{derax:4} (see \cite[1.10]{groth2013derivators}) states that these Kan extensions can always be computed with a \emph{pointwise} formula; since, according to our \ref{kan_are_wei}, Kan extensions can be identified with weighted co\fshyp{}limits on representable weights, and more precisely because all the following concepts are equivalent:
	\begin{itemize}
		\item the left Kan extension of $F : \A\to \C$ along $G : \A\to\B$ computed in $B$;
		\item the weighted colimit of $F$ with respect to the representable $\hom(G\firstblank,B)$;
		\item the conical co\fshyp{}limit of $F$ over the category of elements of $\hom(G\firstblank,B)$;
		\item the conical colimit of the diagram $(G/B)\to \A\xto{F}\B$ over the comma category of $G$ and $B$.
	\end{itemize}
	we can try to express co\fshyp{}end calculus using pointwise Kan extensions in a derivator. We exploit two basic facts: first, according to \ref{is.a.colim} co\fshyp{}ends are colimits over twisted arrow categories (or suitable opposite thereof); second, co\fshyp{}limits are Kan extensions along terminal arrows $\p : K \to *$.
\end{remark}
With these remarks in mind, let $\tw(K)$ be the twisted arrow category of $K$ (see \ref{twisted}) or equivalently the category of elements (see \ref{eltsf}) of $\hom_K$, for a small category $K$; then there exists a functor $\Sigma_K = (t,s) : \tw(K)\to K^\opp\times K$ (see \ref{fibelem}).
\begin{definition}[Homotopy coend in a derivator]\label{coendinder}
	Let $\sD$ be a derivator, and $K\in\caat$ a category. The \emph{homotopy coend}
	\[ \oint^K : \sD(J\times K^\opp\times K)\to \sD(J)\]
	is defined as the pseudonatural transformation with components obtained from the composition
	\[
		\oint^{K,[J]} : \sD(J\times K^\opp\times K) \xto{\Sigma_K^*} \sD(J\times \tw(K))\xto{\p_!} \sD(J)
	\]
\end{definition}
\begin{remark}
	If we rephrase the above definition in terms of the \emph{shifted derivator} $\sD(J|\firstblank) : \caat^\opp \to \tCat$ of $\sD$, \ie the functor $I\overset{\firstblank\times J}\longmapsto I\times J \overset{\sD(\firstblank)}\mapsto \sD(I\times J)$, the homotopy coend $\oint^K$ defines a morphism between the shifted derivators $\sD(K^\opp\times K|\firstblank)\to \sD$. (Pseudonaturality follows from a simple pasting rule between pseudonatural transformations: filling in the details is straightforward.)
\end{remark}
\begin{remark}[Homotopy ends as homotopy limits]
	We can state the dual notion of homotopy \emph{end} in $\sD$: it is enough to replace $\p_!$ with the right adjoint $\p_*$ computing limits instead of colimits in the definition above (the twisted arrow category shall be replace by a suitable opposite thereof): in components,
	\[
		\int_{K,[J]} : \sD(J\times K^\opp\times K) \xto{J\times\Sigma_K^*} \sD(J\times \tw(K))\xto{\p_*} \sD(J)
	\]
\end{remark}
\begin{lemma}
	If $F : \sD\to \sD'$ is a morphism of derivators, there is a canonical `comparison' morphism
	\[\textstyle
		\varsigma : \int^K \circ F \to F\circ \int^K
	\]
	obtained as the composition
	\[\label{chistu}\vcenter{\xymatrix@C=4cm{
		\sD(J\times K\times K^\opp)
		\ar[r]|{\p_! \circ F_{J\times L\times \tw(K)}\circ (t,s)^*}
		\ar@/^2.5pc/[r]^{\p_! \circ (t,s)^*\circ F_{J\times L\times K^\opp\times K}}
		\ar@/_2.5pc/[r]_{F_{J\times L}\circ \p_!\circ (t,s)^*}
		\rtwocell<\omit>{<3>\wr}
		\rtwocell<\omit>{<-3>\wr}
		& \sD'(J)
		}}\]
	where the second morphism results from the pasting of 2\hyp{}cells:
	\[\label{chiddu}\vcenter{\xymatrix{
		\sD'(e) \ar@/_1pc/@{-}[dr]\drtwocell<\omit>{<-1>\epsilon}& \ar[l]\sD'(J)\drtwocell<\omit>{} & \ar[l]\sD(J) \ar@/^1pc/@{-}[dr]\drtwocell<\omit>{<1>\eta}\\
		& \sD'(e) \ar[u]& \ar[l]\sD(e)\ar[u] &\ar[l] \sD(J)
		}}\]
\end{lemma}
\begin{proof}
	The morphism $F : \sD\to \sD'$ has components $F_I : \sD(I)\to \sD'(I)$ and there is a 2\hyp{}cell filling the central square in \eqref{chiddu}. This allows for \eqref{chistu} to be a well\hyp{}posed definition. The rest is an easy check of pseudonaturality conditions, following from conditions on the data defining $\varsigma$.
\end{proof}
A derivator morphism $F$ \emph{preserves homotopy coends} (\ie the above 2\hyp{}cell is invertible) if (and only if? Think about it: is it still possible to describe a co\fshyp{}limit as a certain co\fshyp{}end) it preserves colimits, or more generally left homotopy Kan extensions. Dually, $F$ \emph{preserves homotopy ends} if (and only if?) it preserves right homotopy Kan extensions.
\begin{exercises}
\item \label{ex:laxprofunctors} Study the category of 2\hyp{}profunctors with the composition of 1\hyp{}cells given by the lax coend
\[\notag
	\proP \diamond \proQ (A,C) = \twoint^B \proP(A,B)\times \proQ(B,C)
\]
\begin{itemize}
	\item the 3\hyp{}category $\twoCat$ can be embedded into $\text{2-}\Dist$ in various ways, using $F \mapsto \hom(F,1)$, $F\mapsto \hom(F^\diesis, 1)$, etc.; what is the strictest of these embeddings? Are they fully faithful in a suitable sense?
	\item Is there an adjunction $\hom(F,1)\dashv \hom(1,F)$? Or rather an adjunction $\hom(F^\diesis,1) \dashv \hom(1,F^\bemo)$?
	\item Does $\text{2-}\Dist$ have a dualiser in the sense of \ref{dualiseur}? Does it have products, coproducts?
\end{itemize}
\item \label{ex8:tw-as-a-laxlim} A \emph{lax colimit} for a diagram $F : \J\to \K$ in a 2\hyp{}category $\K$ is an object $L$ with a \emph{lax cocone} $\{FJ \to L \mid J \in\J\}$ satisfying a suitable universal property (state it, mimicking --the dual of-- \ref{laxwedge}).
\item \label{ex8:laxcoends} Dualise \ref{laxwedge}: state the definition of lax \emph{cowedge} $S\din d$ for a 2\hyp{}functor $S : \A^\opp\times \A\to \B$; state the definition of lax coend for $S$ as an \emph{initial} cowedge, the representing object of the functor $d\mapsto \text{LCwd}(S,d)$.
\item \label{barciobar} Prove that the bar and cobar complexes define functors
\[\notag
	B(\firstblank,\C,F)\quad B(G,\C,\firstblank)\quad C(\firstblank,\C,F)\quad C(G,\C,\firstblank).\]
Deduce that there exists a canonical morphism $B(G,\C,F)\to N\C$ in the category of simplicial sets.

Prove that if $F,G : \C\to \Set$ are parallel functors, there is an isomorphism
\[\notag B(G,\C,F)_\bullet^\opp \cong B(F,\C^\opp,G)_\bullet.\]
Prove that the bar complex is isomorphic to the nerve of a certain category of elements (see \ref{eltsf}). Can you provide an intuition about this identification? Dualise this statement.
\item \label{ex8:Y-of-T} Define co\fshyp{}face and co\fshyp{}degeneracy maps for the co\fshyp{}simplicial objects $\Y(T)$ and $\W(T)$ of \ref{Y_and_W} (hint: there is an isomorphism
\[\notag
	\tau : T(X_0,X_n)^{\Pi \A[\vec X]}\cong \Big(T(X_0,X_n)^{\A(X_0,X_1)}\Big)^{\A[X_1|\vec Y|X_n]_{n-1}},
\]
and you want to assemble a map $\Y(T)^{n-1}\to \Y(T)^n$ from its components $\Pi \A[\vec X] \pitchfork T(X_1,X_n) \to \Y(T)^n$; this defines $d^0$. The map $d^n$ is defined via an isomorphism $\sigma$ and a similar argument).
\item \label{ex8:cohfubini} \righteyes Prove the Fubini theorem for simplicially coherent co\fshyp{}ends: given a functor $T : \A^\opp\times \A\times \B^\opp\times \B\to\C$, then
\begin{align*}
	\oint^A\left( \oint^B T(A,A,B,B)\right) & \cong	\oint^{(A,B)\in\A\times\B} T(A,B,A,B)    \\
	                                        & \cong \oint^B\left( \oint^A T(A,A,B,B)\right)
\end{align*}
(hint: it is a simple theorem about the relation between $\delta(\A\times\B)$ and $\delta\A\times\delta\B$ and about \ref{fubozzo}).
\item \label{ex8:cohnat} Prove that $[F, \overline G] \cong \Cat_{\bDelta}(\!( F,G )\!) \cong [\underline F,G]$, using  \ref{cohcoend} and a formal argument.
\item \label{ex8:cohkan} Prove that the isomorphisms \eqref{kan:chistu} and \eqref{kan:chiddu} hold, and thus define coherent Kan extensions for simplicial functors.
\item Prove that the standard resolutions of \ref{standard_res} `absorb coherence' in coherent Kan extensions, showing that
\[\notag
	\Ran_G \overline F \cong \hoRan_G F \qquad\qquad \Lan_G \underline F \cong \hoLan_G F .
\]
(a preliminary lemma: prove that $\overline F(\firstblank)\cong \int_A FA^{\delta\A(\firstblank,A)}$).
\item Find an expression for $\hoRan_HF$, given a cospan of functors $\Cat_{\bDelta}(\C,\sSet) \xot{H} \A \xto{F}\Cat_{\bDelta}(\B,\sSet)$.
\item \label{ex8:derivcoend} \awful Prove that $\int^K : \sD(K^\opp\times K|-)\to \sD$ defines a morphism of derivators (you can either prove that a functor $u : K\to L$ induces a morphism between the shifted derivators $\sD(L|-)\to \sD(J|-)$, or prefer an explicit argument --both ways are considerably long).
\item Prove that `coends in a derivator are pointwise', \ie that given an arrow $j : e\to J$ there is a canonical isomorphism $j^*\big(\int^K X\big)\cong \int^K j^*X$ for each $X\in \sD(J\times K^\opp\times K)$.
\item \label{ex8:derivfubini} State and prove the Fubini theorem for homotopy coends in $\sD$: the diagram
\[\notag\vcenter{\xymatrix@C=3cm{
	\sD(J\times K^\oo\times L^\oo) \ar[r]^{(J\times \Sigma_K\times L^\oo)^*}\ar[dddr]_{\Sigma_{K\times L}}
	\ddrtwocell<\omit>{\alpha}
	\ddrtwocell<\omit>{<10>\beta}
	\ar[d]_{\Sigma_L^*}& \sD(J\times \tw(K)\times L^\oo) \ar[d]^{\p_!^K}\\
	\sD(J\times K^\oo\times \tw(L))\ar[d]_{\p_!^L} & \sD(J\times L^\oo)\ar[d]^{\Sigma_L^*} \\
	\sD(J\times K^\oo) \ar[d]_{\Sigma_K^*}& \sD(J\times \tw(L))\ar[d]^{\p_!^L}\\
	\sD(J\times \tw(K)) \ar[r]_{\p_!^K}& \sD(J)
	}}
\]
commutes for canonically determined 2\hyp{}cells $\alpha$ and $\beta$.
\item \label{ex8:derivweighlim} We define a \emph{bimorphism} between three derivators to be a family of functors
\[\notag B_{IJ} : \sD(I)\times \mathbb{E}(J) \to \mathbb{F}(I\times J)\]
in $\tCat$ endowed with 2\hyp{}cells $\gamma_{u_1,u_2}$ filling the diagrams
\[\notag
	\vcenter{\xymatrix{
			\sD(J_1) \times \mathbb{E}(J_2)
			\drtwocell<\omit>{\gamma}
			\ar[r]\ar[d] & \mathbb{F}(J_1\times J_2)\ar[d]\\
			\sD(I_1) \times \mathbb{E}(I_2) \ar[r]& \mathbb{F}(I_1\times I_2)
		}}
\]
These $\gamma_{u_1,u_2}$ are subject to certain coherence conditions: for every pair $(\alpha_1,\alpha_2)$ of natural transformation $\alpha_\epsilon : u_\epsilon \To v_\epsilon : I_\epsilon \to J_\epsilon$ and $\epsilon = 1,2$ we have that the two diagrams
\begin{gather*}
	\vcenter{\xymatrix{
	(u_1\times u_2)^\ast \circ (v_1\times v_2)^\ast\circ  B\ar[r]^-\gamma \ar@/_1.0pc/[dr]_-\gamma&(u_1\times u_2)^\ast\circ B\circ (v_1^\ast\times v_2^\ast)\ar[d]^\gamma\\
	&  B\circ (u_1^\ast\times u_2^\ast)\circ (v_1^\ast\times v_2^\ast)
	}}\\
	\vcenter{\xymatrix{
			(u_1\times u_2)^\ast\circ B\ar[r]\ar[d]_\gamma&(u_1'\times u_2')^\ast\circ B\ar[d]^\gamma\\
			B \circ(u_1^\ast\times u_2^\ast)\ar[r]& B\circ (u_1'^\ast \times u_2'^\ast)
		}}
\end{gather*}
commute. Write suitable diagrams of 2\hyp{}cells expressing these commutativities less tersely. Show that there is a category $\sfDer_b((\sD,\mathbb E), \mathbb F)$ of bimorphisms of derivators, and that there is an equivalence
\[\notag\sfDer_b((\sD,\mathbb E), \mathbb F)\cong \sfDer(\sD\times \mathbb E,\mathbb F)\]
Is there a way to prove this result using the fact that the right side is a pseudo\hyp{}end?
\item \awful Let $\sD^{\Sets}$ be the \emph{discrete} derivator sending a small category $J$ into the functor category $\Set^J = \Cat(J,\Set)$. State a definition and an existence theorem for \emph{weighted colimits} in a derivator $\sD$: given a small category $I\in\caat$ and a bimorphism $\boxplus : (\sD^{\Sets},\sD)\to \sD(-|I)$, we define the colimit of $X\in \sD(J)$, weighted by $W\in \sD^{\Sets}(J^\opp)$ as the coend (in $\sD$) $\int^J W\boxplus X$, \ie as the image of the pair $(W,X)$ under the composition
\[\notag
	\sD^{\Sets}(J^\opp)\times \sD(J) \xto{\boxplus_{J^\opp,J}} \sD(J^\opp\times J|I) \xto{\int^{J,[I]}}\sD(I).
\]
\end{exercises}

\chapter{Addenda et Paralipomena}
\begin{abstract}
	The present chapter collects a few results that are complementary to the main discussion; they are very recent development of co\fshyp{}end calculus that find applications in algebra, analysis (in the intention of the author of \cite{day2011monoidal}, even to quantum physics!) and computer science (especially functional programming). Our focus here is not on proofs; instead, we offer a birdseye view on material not developed enough to deserve a separate chapter, and yet too complex to be included as a terse example elsewhere; the chapter strives for a more relaxed tone, and tries to open a window on, so to speak, what is currently an `open research problem' that employs co\fshyp{}end calculus.
\end{abstract}
\epigraph{El dorso estaba numerado con ocho cifras. Llevaba una pequeña ilustración, como es de uso en los diccionarios: un ancla dibujada a la pluma, como por la torpe mano de un niño.	Fue entonces que el desconocido me dijo: -Mírela bien. Ya no la verá nunca más.}{J.L. Borges --- \emph{El libro de arena}}
\section{Fourier theory}\label{sec:promono}\index{Category!promonoidal ---}
According to our §\ref{struct_of_Prof}, the bicategory $\Dist$ is monoidal with respect to the pseudo\hyp{}cartesian structure (similarly, every $\Dist(\V)$ inherits a symmetric monoidal structure from a symmetric monoidal structure on $\VCat$).

This means that we can consider internal monoids in $\Dist$: objects endowed with maps
\[
	\xymatrix{
	\M\times\M \ar[r]|-@{|}^-{\mathfrak{m}} & \M & 1 \ar[r]|-@{|}^-{\mathfrak{i}} & \M
	}
\]
in $\Dist$ that witness the fact that $\M$ is an internal (pseudo)monoid.

Such internal monoids take the name of \emph{promonoidal categories}.

Informally speaking, a promonoidal category is what we obtain if we replace every occurrence of the word \emph{functor} with the word \emph{profunctor} in the definition of monoidal category (of course, taking care of the coherence conditions imposed by the weak 2\hyp{}category structure of $\Dist$).

More precisely, we can give the following definition:
\begin{definition}[Promonoidal structure]\index{Promonoidal category}
	Let $\C$ be a category. A \emph{promonoidal structure} consists of a tuple
	\[\fkP = (\C, P, J,\alpha,\lambda,\rho)\]
	where
	\begin{enumtag}{pm}
		\item $\C$ is a category endowed with
		\item a bi\hyp{}profunctor $P : \C\times \C \pto \C$ (the monoidal \emph{multiplication}) and
		\item a profunctor $J : 1\pto \C$ (the monoidal unit), such that the following two diagrams
		\[\notag
			\xymatrix@C=2cm{
			\C \drtwocell<\omit>{\alpha} \times \C \times \C \ar[d]|-@{|}_{\hom\times P}\ar[r]|-@{|}^{P\times\hom}& \C \times\C \ar[d]|-@{|}^{P}\\
			\C \times \C \ar[r]|-@{|}_{P}& \C
			}\qquad
			\xymatrix@C=2cm{
			\C \drtwocell<\omit>{\rho} \ar@/_1pc/@{=}[dr]_\hom\ar[r]|-@{|}^{J\times \hom} & \C \times \C \ar[d]|-@{|}^{P} & \C \ar[l]|-@{|}_{\hom\times J} \ar@/^1pc/@{=}[dl]^\hom \\
			& \C \urtwocell<\omit>{\lambda}&
			}
		\]
		are filled by the indicated 2\hyp{}cells,
		\item 2-cells, respectively called the \emph{associator}
		\[\label{pms:uno}\alpha : P \diamond (P\times \hom) \cong P\diamond(\hom\times P)\]
		and the \emph{left} and \emph{right unitors}
		\[\label{pms:due}\lambda : P\diamond(\hom\times J)\cong \hom \quad \rho : P\diamond(J\times\hom)\cong\hom\]
	\end{enumtag}
\end{definition}
\begin{remark}
	Coend calculus allows to turn \eqref{pms:uno} and \eqref{pms:due} into diagrammatic relations:
	\index{Associator!promonoidal ---}
	\begin{itemize}
		\item The associator amounts to an isomorphism linking the two sets below (note that each component $\alpha^{ABC}_{D}$ has four arguments, three contravariant and one covariant, whereas $P$ has components $P^{AB}_C = P(A,B;C)$ as a functor $\C^\opp\times\C^\opp\times \C \to \Set$).
		      \begin{align*}
			      (P\diamond (\hom\times P))_{ABC;D} & =\int^{XY} P_D^{XY}H_A^X P_Y^{BC}                     \\
			                                         & \cong \int^YZ\Big(\int^X P_D^{XY}H_A^X \Big) P_Y^{BC}
			      \cong \int^Z P_D^{AY}P_Y^{BC}                                                              \\
			      (P\diamond(P\times\hom))_{ABC;D}   & \cong \int^{XY} P_D^{XY} H_Y^C P_X^{AB}
			      \cong \int^Z P_X^{AB} P_D^{XC}.
		      \end{align*}
		\item The left unit axiom is equivalent to the isomorphism between the functor
		      \begin{align*}
			      (A,B) & \mapsto \int^{YZ} J_Z H^A_Y P^{YZ}_B              \\
			            & \cong \int^Z J_Z\Big( \int^Y H^A_Y P^{YZ}_B \Big) \\
			            & \cong \int^Z J_Z P^{AZ}_B
		      \end{align*}
		      and the hom functor $(A,B)\mapsto \C(A,B)$.
	\end{itemize}
\end{remark}
The most interesting feature of promonoidal structure in categories is that they correspond bijectively with monoidal structures on the category of functors $\Cat(\C, \Sets)$, framing the construction of \emph{Day convolution} \ref{day} in its maximal generality.
\begin{proposition}\label{promonoshit}\index{Convolution product}
	Let $\fkP = (P, J,\alpha, \rho,\lambda)$ be a promonoidal structure on the category $\C$; then we can define a $\fkP$-\emph{convolution} monoidal structure on the category $\Cat(\C, \Set)$, via
	\begin{gather}
		[F\ast_{\fkP} G]C = \int^{AB} P(A,B;C)\times FA\times GB\\
		J_{\fkP} = J
	\end{gather}
	and this turns out to be a monoidal structure on $\Cat(\C, \Set)$. We denote the monoidal structure $(\Cat(\C, \Set), \ast_{\fkP}, J_{\fkP})$ shortly as $[\C, \Set]_{\fkP}$.
\end{proposition}
The same definition, changing the cartesian structure with the monoidal structure of $\V$, yields a notion of $\fkP$\hyp{}convolution on $\VCat(\C,\V)$ for a $\V$\hyp{}category $\C$.
\begin{definition}
	A functor $\Phi : [\A, \Sets]_{\fkP} \to [\B, \Sets]_{\fkQ}$ is said to \emph{preserve the convolution product} if the obvious isomorphisms hold in $[\B, \Sets]_{\fkQ}$:
	\begin{itemize}
		\item $\Phi(F\ast_{\fkP} G) \cong \Phi(F)\ast_{\fkQ}\Phi(G)$;
		\item $\Phi(J_{\fkP}) = J_{\fkQ}$.
	\end{itemize}
\end{definition}
\begin{remark}
	It is observed in \cite{imkelly} that for a monoidal $\A$ the category of presheaves $[\A^\opp,\V]$ endowed with the convolution monoidal structure is the \emph{free monoidal cocompletion} of $\A$, having in $\BF{Mon}$ (monoidal categories, monoidal functors and monoidal natural transformations) the same universal property that $[\A^\opp,\V]$ has in $\Cat$.
\end{remark}
There is a bijection between the promonoidal structures on $\C$, and monoidal structure on $\Cat(\C,\Set)$; this is the content of Exercise \ref{es8_2}.
\subsection{Fourier transforms via coends}
\index{Multiplicative kernel|see{Fourier transform}}
\index{Fourier transform}
For the rest of the section, $\V$ is assumed to be a complete and cocomplete *\hyp{}autonomous symmetric monoidal closed category.
\begin{definition}\label{mulker}
	Let $\A ,\C$ be two promonoidal categories (thus implicitly regarded as objects of $\Dist$) with promonoidal structures $\fkP$ and $\fkQ$ respectively; a \emph{multiplicative kernel} from $\A$ to $\C$ consists of a profunctor $K : \A\pto \C$ endowed with two natural isomorphisms
	\begin{enumtag}{k}
		\item $\displaystyle \int^{YZ}K^A_Y K^B_Z P^{YZ}_X\cong \int^C K^C_X P^{AB}_C$;
		\item $\displaystyle \int^C K^C_X J_C\cong J_X$.
	\end{enumtag}
	These isomorphisms say that $K$ mimicks the behaviour of the hom functor (in fact, the hom functor $\hom_\A$ is the \emph{identity} multiplicative kernel $\A\pto\A$: the isomorphisms above follow from \ref{ninjayo}).

	We define a \emph{multiplicative} natural transformation $\alpha : K\to H$ between two kernels as a 2\hyp{}cell in $\Dist$ commuting with the structural isomorphisms given in \ref{mulker}. This, together with the fact that multiplicative kernels compose, yields a category of kernels $\ker(\A,\C)$.
\end{definition}
\begin{definition}
	Let $K : \A\pto\C$ be a multiplicative kernel between promonoidal categories; we define the $K$-\emph{Fourier transform} $f\mapsto \hat K(f) : \C\to \Sets$, obtained as the image of $f : \A\to\Sets$ under the left Kan extension $\Lan_\yon K : [\A, \Sets]\to [\C, \Sets]$. Explicitly, this is the coend
	\[
		\F_K(f) : X\mapsto \int^A K(A,X)\otimes fA.
	\]
\end{definition}
We can also define the dual Fourier transform:
\[
	\F^\lor(g) : Y\mapsto \int_A [K(A,X), gA]
\]
and find the relation $\F^\lor_K(g) \cong \F_K(g^*)^*$.

The following results are easily proved using standard co\fshyp{}end calculus:
\begin{proposition}
	Let $K : \A \pto \cate X$ be a multiplicative kernel, and let $\A$ be a promonoidal category; then
	\begin{enumtag}{mk}
		\item $\F_K$ preserves the upper $\mathfrak P_\A$\hyp{}convolution of presheaves $f,g$, defined as
		\[f\mathrel{\overline{\ast}} g = \int^{AA'} fA\otimes gA'\otimes P(A,A',\firstblank);\] dually,
		\item $\F_K^\lor$ preserves the lower $\mathfrak P_\A$\hyp{}convolution of presheaves $f,g$, defined as
		\[f\mathrel{\underline{\ast}} g = \int_{AA'} \Big(fA^*\otimes (gA')^*\otimes P(A,A',\firstblank)\Big)^*;\]
	\end{enumtag}
\end{proposition}
Observe that (as stated in \cite{day2011monoidal}), both the upper and lower convolution product yield associative and unital monoidal structures on the functor category $\VCat(\A,\V)$; the upper product preserves $\V$\hyp{}colimits in each variable, while the lower product preserves $\V$\hyp{}limits
in each variable.

\index{Convolution!--- product}
\index{Convolution!upper and lower ---}
The  lower and upper convolution transform into each other under the equivalence of $\V$\hyp{}categories
\[
	\VCat(\A,\V)^\opp\cong \VCat(\A^\opp,\V^\opp);
\]
this means that under the above equivalence $(f \mathrel{\overline{\ast}} g)^* \cong f^* \mathrel{\underline{\ast}} g^*$.
\begin{theorem}
	Let $\V$ be a *\hyp{}autonomous monoidal base; then we can define the pairing $\VCat(\A,\V)\times \VCat(\A,\V) \to \V$ as the twisted form of functor tensor product (as defined in \ref{tenso_pro_offunc})
	\[
		\langle f,g\rangle = \int^A fA^*\otimes gA
	\]
	\index{Parseval formula|see{Fourier transform}}
	\index{Fourier transform}
	If $K$ is a kernel such that the Fourier transform $\Lan_\yon K$ is fully faithful, we have an analogue of \emph{Parseval formula}:
	\[
		\langle f,g\rangle \cong \langle \F_K(f),\F_K(g)\rangle.
	\]
\end{theorem}
\index{Combinatorial species}
Fourier theory is linked to the theory of Joyal's \emph{combinatorial species}
(see \ref{analuo}): let $E : \A\pto \cate X$ be a multiplicative kernel; if we write
$E(A,X):= X^A$ (without any reference to a tensor operation between $A$ and
$X$), the $E$-Fourier transform can be expressed as an $\A$\hyp{}indexed formal power
series as follows:
\begin{align}
	\F_E(f) & = \int^A fA\otimes X^A \notag        \\
	        & \cong \sum_{A\in\A} fA\otimes_\A X^A
\end{align}
(it is understood that $f : \A \to \V$ is a fixed combinatorial species.) This can be made precise as follows: let $\A$ be the (free $\V$\hyp{}category on the) permutation category of
\ref{prelim_notata}; then $E(n,X) := X^{\otimes n}=X\otimes\dots\otimes X$ is a
multiplicative kernel and an $E$\hyp{}analytic functor results as the left Kan
extension
\begin{align}
	FX & =\int^n f(n)\otimes X^{\otimes n} \notag                       \\
	   & \cong \sum_{n\in\P} f(n)\otimes_{\text{Sym}(n)} X^{\otimes n}.
\end{align}
Given two combinatorial species $f,g : \A \to \V$ and the associated analytic
functors $F,G$, the convolution $F \ast G$ is again an analytic functor, and its
generating combinatorial species is the upper convolution product of $f,g$ (with
respect to the implicit promonoidal structure of $\A$):
\begin{align}
	F \ast G(X) & \cong \int^{AB} FA\otimes GB\otimes p(A,B;X) \notag                              \\
	            & \cong \int^{ABUV} fU\otimes E(U,A)\otimes gV\otimes E(V,B)\otimes p(A,B;X)\notag \\
	\star       & \cong \int^{UVC} fU\otimes gV\otimes p(U,V;C)\otimes E(C,X) \notag               \\
	            & \cong \int^C (f\mathrel{\overline{\ast}} g)(C)\otimes E(C,X).
\end{align}
(note that in ($\star$) we used the fact that $E$ is a multiplicative kernel.)
\section{Tambara theory}
\index{Tambara module}
\index{Profunctor}
\begin{definition}
	Let $\C$ be a monoidal category with monoidal unit $I$. A (\emph{left}) dio \emph{Tambara module} on $\C$ consists of:
	\begin{itemize}
		\item a profunctor
		      $P : \C^\opp \times \C \rightarrow \Set$;
		\item a family of functions
		      $\tau_{A}(X, Y) : P(X, Y) \longrightarrow P(A \otimes X, A \otimes Y)$
		      natural in $X, Y$ and a wedge in $A$, satisfying the two equations:
		      \begin{gather}
			      \vcenter{\xymatrix{
			      P(X, Y) \ar[rr]^{\tau_I(X,Y)}\ar@{=}[dr]&& P(I\otimes X, I\otimes Y) \ar[dl]^{P(l_X^{-1},l_Y)}\\
			      & P(X,Y)
			      }}
			      \notag\\
			      \vcenter{\xymatrix{
			      P(X,Y) \ar[rr]^{\tau_{A'}(X, Y)}
			      \ar[dr]_{\tau_{A \otimes A'}(X, Y)\quad}
			      && P(A'\otimes X, A'\otimes Y) \ar[dl]^{\quad\tau_{A}(A' \otimes X, A' \otimes Y)} \\
			      &  P(A \otimes A' \otimes X, A \otimes A' \otimes Y)&
			      }}
		      \end{gather}
	\end{itemize}
	The notion of \emph{right} Tambara module is given in a similar fashion, using maps $\nu_A(X,Y) : P(X, Y) \longrightarrow P(X \otimes A, Y \otimes A)$ satisfying the relations
	\begin{enumtag}{lt}
		\item $P(r_X, r_X^{-1})\circ \nu_A(X,Y) = \id_{P(X,Y)}$;
		\item $\nu_A(X\otimes A',Y\otimes A')\circ \nu_{A'}(X,Y) = \nu_{A'\otimes A}(X,Y)$.
	\end{enumtag}
	The definition can be given for profunctors enriched over any other monoidal base different from $\Set$; for example, in the original work by Tambara, \cite{tambara2006distributors}, categories are assumed to be enriched over vector spaces. In the paper \cite{pastro2008doubles} the enrichment base is completely arbitrary (\ie it is just a symmetric monoidal closed category).
\end{definition}
\begin{definition}
	Define the category $\Tamb(\C)$ whose:
	\begin{itemize}
		\item objects are Tambara modules $(P, \tau)$ consisting of a
		      profunctor $P : \C^\opp \times \C \rightarrow \Set$
		      and Tambara structures $\tau_{A}(X, Y)$.
		\item morphisms $(P, \tau) \rightarrow (Q, \sigma)$
		      are natural transformations
		      $\gamma : P \Rightarrow Q$ such that for all $A, X, Y$ the following
		      diagram commutes:
		      \[ \vcenter{\xymatrix@C=1.5cm{
			      P(X, Y)
			      \ar[r]^-{\tau_{A}(X, Y)}
			      \ar[d]_{\gamma_{(X, Y)}}
			      & P(A \otimes X, A \otimes Y)
			      \ar[d]^{\gamma_{(A \otimes X, A \otimes Y)}}
			      \\
			      Q(X, Y)
			      \ar[r]_-{\sigma_{A}(X, Y)}
			      & Q(A \otimes X, A \otimes Y)
			      }}
		      \]
	\end{itemize}
	There is a functor to the category of endo\hyp{}profunctors on $\C$,
	\[
		\iota : \Tamb(\C) \longrightarrow \Dist(\C,\C)
	\]
	which forgets the Tambara structure.
\end{definition}
The codomain of $\iota$ is monoidal with respect to composition of 1\hyp{}cells, as every hom\hyp{}category of endomorphisms: it turns out that Tambara modules can be composed, and that $\iota$ is strong monoidal with respect to this monoidal structure.
\begin{remark}
	The category $\Tamb(\C)$ has a monoidal structure whose:
	\begin{itemize}
		\item \emph{unit} is the hom\hyp{}functor
		      $\hom_\C : \C^\opp \times \C \rightarrow \Set$
		      which has a canonically associated Tambara structure:
		      \[
			      \C(X, Y) \longrightarrow \C(A \otimes X, A \otimes Y)
		      \]
		\item The profunctor composition of
		      $(P, \tau)$ and $(Q, \sigma)$ given by the coend
		      \[
			      (P \diamond Q)(X, Y) = \int^{Z} P(X, Z) \times Q(Z, Y)
		      \]
		      has a Tambara structure
		      $(P \diamond Q)(X, Y) \rightarrow (P \diamond Q)(A \otimes X, A \otimes Y)$
		      induced by the maps
		      \[
			      \tau_{A} \times \sigma_{A} : P(X, Z) \times Q(Z, Y)
			      \rightarrow P(A \otimes X, A \otimes Z) \times Q(A \otimes Z, A \otimes Y)
		      \]
		      using the universal property of the coend.
	\end{itemize}
	This makes the functor $\iota : \Tamb(\C) \rightarrow \Dist(\C,\C)$
	strong monoidal.
\end{remark}
\begin{proposition}
	The forgetful functor $\iota : \Tamb(\C) \rightarrow \Dist(\C,\C)$
	forms part of an adjoint triple:
	\[
		\vcenter{\xymatrix{
				\Tamb(\C) \ar[r]|\iota & \ar@<6pt>[l]^\varphi  \ar@<-6pt>[l]_\theta \Dist(\C,\C)
			}}
	\]
\end{proposition}
\begin{itemize}
	\item The left adjoint $\phi : \Dist(\C,\C) \rightarrow \Tamb(\C)$ constructs
	      the \emph{free Tambara module} from a profunctor. This is given by the formula
	      \[
		      \phi_{P}(X, Y)
		      = \int^{C, U, V} \C(X, C \otimes U) \times \C(C \otimes V, Y) \times P(U, V)
	      \]
	      with Tambara module structure given by
	      \[
		      \vcenter{\xymatrix{
				      \C(X, C \otimes U) \times \C(C \otimes V, Y) \times P(U, V)
				      \ar[d] \\
				      \C(A \otimes X, A \otimes C \otimes U)
				      \times \C(A \otimes C \otimes V, A \otimes Y) \times P(U, V)
			      }}
	      \]
	      together with the coprojection $q_{A \otimes C}$,
	      using the universal property of the coend.
	\item The right adjoint $\theta : \Dist(\C,\C) \rightarrow \Tamb(\C)$
	      constructs the \emph{cofree Tambara module} from a profunctor. This is given by the formula
	      \[
		      \theta_{P}(X, Y) = \int_{C} P(C \otimes X, C \otimes Y)
	      \]
	      with Tambara module structure given by
	      $\theta_{P}(X, Y) \rightarrow \theta_{P}(A \otimes X, A \otimes Y)$ is induced
	      by the projection functions,
	      \begin{gather}
		      p_{C \otimes A} : \int_{C} P(C \otimes X, C \otimes Y)
		      \rightarrow P(C \otimes A \otimes X, C \otimes A \otimes Y)\notag\\
		      p_{C} : \int_{C} P(C \otimes A \otimes X, C \otimes A \otimes Y)
		      \rightarrow P(C \otimes A \otimes X, C \otimes A \otimes Y)
	      \end{gather}
	      using the universal property of the end.
\end{itemize}
The proof of the following proposition (a `recognition principle' for Tambara modules) goes by inspection
using the definition of coalgebra: such a map is determined by
\begin{itemize}
	\item An object $P$ of $\Dist(\C,\C)$ given by
	      $P : \C^\opp \times \C \rightarrow \Set$;
	\item A structure map given by a natural transformation
	      $\tau : P \Rightarrow \theta_{P}$ in $\Dist(\C,\C)(P,\theta_P)$ whose components
	      $\tau(X, Y) : P(X, Y) \rightarrow \theta_{P}(X, Y)$ are given by:
	      \[
		      \tau(X, Y) : P(X, Y) \longrightarrow \int_{A} P(A \otimes X, A \otimes Y)
	      \]
	\item By the universal property of the end at codomain, the structure map is
	      determined by a wedge in $A$,
	      \[
		      \tau_{A}(X, Y) :
		      P(X, Y) \longrightarrow P(A \otimes X, A \otimes Y)
	      \]
	      that is moreover natural in $X,Y$.
\end{itemize}
\begin{proposition}
	The adjunction $\iota \dashv \theta$ yields a comonad
	\[
		\Theta : \Dist(\C,\C) \rightarrow \Dist(\C,\C)
	\]
	whose category of coalgebras is
	isomorphic to $\Tamb(\C)$. Dually, the adjunction $\phi \dashv \iota$ yields a monad
	$\Phi : \Dist(\C,\C) \rightarrow \Dist(\C,\C)$ whose category of algebras is
	isomorphic to $\Tamb(\C)$. Moreover, there is an
	adjunction
	\[
		\xymatrix{
			\Dist(\C,\C) \ar@<4pt>[r]^\Phi\ar@{}[r]|\perp & \ar@<4pt>[l]^\Theta\Dist(\C,\C)
		}
	\]
	between the resulting monad $\Phi = \iota\circ\varphi$ and comonad $\Theta = \iota\circ\theta$ on $[\C^\opp \times \C, \Set]$.
\end{proposition}
Profunctors $\A \pto \B$ can equivalently be described as left adjoints $\Cat(\B^\opp,\Set) \to \Cat(\A^\opp,\Set)$; thus we obtain that
\begin{corollary}
	The left adjoint
	\[
		\Phi : \Cat(\C^\opp \times \C, \Set)
		\to \Cat(\C^\opp \times \C, \Set)
	\]
	is equivalent to the \emph{endo\hyp{}profunctor} $\check\Phi : \C^\opp\times\C\pto\C^\opp\times\C$ whose action on objects is given by the coend:
	\[
		\check{\Phi}(X, Y, U, V)
		= \int^{C} \C(X, C \otimes U) \times \C(C \otimes V, Y)
	\]
\end{corollary}
This endo\hyp{}profunctor
$\check{\Phi} : \C^\opp \times \C \pto \C^\opp \times \C$
actually underlies a \emph{promonad} (see Exercise \ref{ex:promonad}) in the bicategory $\Dist$.
From the formal theory of monads \cite{Street1972} it is known that the bicategory $\Dist$ admits the Kleisli
construction for promonads, so we can ask what is the Kleisli category of $\check\Phi$: such category
is called the (\emph{left}) \emph{double} of the monoidal category $\C$ and it is denoted $\BF{Db}(\C)$; it has the same objects as
$\C^\opp \times \C$, and hom\hyp{}sets defined by the coend
\[
	\BF{Db}(\C)\big( (X, Y), (U, V) \big)
	= \int^{C} \C(X, C \otimes U) \times \C(C \otimes V, Y)
\]
This formula provides the foundation for all of \emph{profunctor optics} (see \cite{barto,noi,pickering2017profunctor}).
\cite{pastro2008doubles} proves that there is an equivalence of categories:
\[
	\Tamb(\C) \simeq \Cat(\BF{Db}(\C), \Set)
\]
\begin{exercises}
\item \label{es8_1} Prove equations \ref{promonoshit} using associativity and unitality for $\fkP$.
\item \label{es8_2} \index{Category!promonoidal ---}\index{Promonoidal category}\index{Convolution product}
Let $* : \Cat(\C,\Set)\times\Cat(\C,\Set) \to \Cat(\C,\Set)$ be a monoidal structure with monoidal unit $u : \C \to \Set$; show that the assignment
\[
	P(A,B;C) := (\coyon(A) * \coyon(B))(C)\qquad JA := uA
\]
is a promonoidal structure on $\C$, regarded as an object of $\Dist$. An elegant result of Day shows that this sets up a bijection between the ways in which $\Cat(\C,\Set)$ is a (pseudo)monoid in $\Cat$, and the ways in which $\C$ is a (pseudo)monoid in $\Dist$: prove it.
\item \label{es8_3} \index{Category!promonoidal ---}\index{Promonoidal category}\index{Convolution product}Outline the promonoidal structure $\fkP$ giving the Day convolution described in \ref{day}. If $\C$ is any small category, we define $P(A,B;C) = \C(A,C)\times \C(B,C)$ and $J$ to be the terminal functor $\C\to\Sets$. Outline the convolution product on $\Cat(\C, \Set)$, called the \emph{Cauchy convolution}, obtained from this promonoidal structure.
\item \label{es8_4} Is the composition of two kernels (see \ref{mulker}) again a kernel? Define the category of multiplicative kernels $\ker(\A,\C)\subset \Dist(\A,\C)$.
\item \label{es8_5} Show that a profunctor $K : \A\pto \C$ is a multiplicative kernel if and only if the cocontinuous functor $\Lan_\yon K=\hat K : [\A, \Sets]\to [\C, \Sets]$ corresponding to $\bar K : \A \to \VCat(\C,\V)$ under the construction in \ref{alternative} is monoidal with respect to the convolution monoidal structure on both $[\A, \Sets]_{\fkP}$ and $ [\C, \Sets]_{\fkQ}$.

Describe the isomorphisms $k_1, k_2$ when $\fkP$ is Day convolution.
\item \label{es8_6} Show that a functor $F : (\A,\otimes_{\A}, I) \to (\C,\otimes_\C, J)$ between monoidal categories is strong monoidal, \ie
\begin{itemize}
	\item $F(A\otimes B)\cong FA\otimes FB$;
	\item $FI\cong J$
\end{itemize}
naturally in $A,B$ if and only if $\proP^F=\hom(F,1)$ is a multiplicative kernel.

Dually, show that for $\A,\C$ promonoidal, $F : \C\to \A$ preserves convolution on $[\A, \Sets]_{\fkP}, [\C, \Sets]_{\fkQ}$ precisely if $\proP_F=\hom(1, F)$ is a multiplicative kernel.
\item Show the following properties of the $K$-Fourier transform:\index{Parseval identity}
\begin{itemize}
	\item There is the canonical isomorphism
	      \[\notag
		      \hat K(f) \cong \int^A K(A,\firstblank)\times f(A)
	      \]
	\item $\hat K$ preserves the convolution monoidal structure (this is the \emph{Parseval identity} for the Fourier transform);
	\item $\hat K$ has a right adjoint defined by
	      \[\notag
		      \check{K}(g) \cong \int_x [K(\firstblank, X), g(X)].
	      \]
\end{itemize}

\end{exercises}

\appendix
\chapter{Review of category theory}\label{chap:review}
\begin{abstract}
	The scope of the present appendix is to recall the bare minimum of category theory that we employ along the book. This is also meant to fix our notation beyond what is already done in the introduction. Even though we assume the typical reader of this book is already acquainted with basic category theory, we will still indulge in a certain desire of self\hyp{}containment: as a rule of thumb, capital results like \ref{lem:the-real-yoda} or \ref{thm:yoda-is-dense} are proved in full detail; most of the other proofs are barely sketched, and some of the marginal results are not proved. It is in fact unrealistic to aim at such a big target as providing a complete account of basic category theory in this appendix; the reader not feeling at ease while consulting the present section is warmly invited to parallel it with more classical references as \cite{McL, riehlcontext,leinster2014basic,simmons2011introduction}.
\end{abstract}
\epigraph{A major explanation for the cognitive advantages of diagrams is \emph{computational offloading}: They shift some of the processing burden from the cognitive system to the perceptual system, which is faster and frees up scarce cognitive resources for other tasks.}{Daniel L. Moody --- \cite{moody2009physics}}
\section{Categories and functors}\label{cat_and_fun}
In simple terms, a category is a structure capable to abstract a number of working assumptions of everyday mathematics:
\begin{itemize}
	\item All objects of a given `kind' can be collected in a class;
	\item such objects form coherent conglomerates, allowing for relations between structures to form;
	\item far from being rare, these relational conglomerates are pretty common and arise at every corner of pure and applied mathematics.
\end{itemize}
Of course, it is not trivial at all to define what is the `kind' of a structure; nor it is easy to define properly what is the process leading to the formation of `homomorphisms', understood as a map preserving the structure (why a `homomorphism of topological spaces' is a continuous function and not an open map?). As always, only experience gives the right answer to such questions; but category theory is of great help in building and strengthening some sort of sixth sense for which Mathematics among the many is `the right one'.

Lacking the sufficient authority to break the unspoken rule that a mathematician should talk about explicit theorems or concrete examples, and not about speculations, we shall refrain from this kind of ramblings into philosophy of Mathematics. The problem of what is category theory, and what's it for, and what's there to do in it, has however been addressed to some extent. A much better informed opinion than ours can be found in \cite{kromer2007tool, marquis2010category}:\footnote{This shows how category theory fits into a solid track of prior philosophical and mathematical research; in a certain sense it is the pinnacle of such research, and the result of its declination in the field of pure mathematics: category theory is what structuralism becomes when it is merged with mathematical craftmanship.} we warmly invite the interested reader to refer to these sources, but it is our humble opinion that a too deep study of the philosophical ground that made category theory possible quickly turns out to be counterproductive, when is not backed up by at least an elementary knowledge of his techniques.\index{Category!--- theory and structuralism}
\begin{definition}[Category]\label{def:cate}
	\index{Category}
	\index{Object}
	\index{Morphism}
	A \emph{(locally small) category} $\C$ consists of the following data:
	\begin{enumtag}{c}
		\item \label{c:uno} A class $\C_o$ whose elements are termed \emph{objects}, usually denoted with Latin letters like $A,B,\dots$;
		\item \label{c:due} A collection of sets $\C(A,B)$, indexed by the pairs $A,B\in \C_o$, whose elements are termed \emph{morphisms} or \emph{arrows} (see \ref{why-arrows} below) with \emph{domain} $A$ and \emph{codomain} $B$;
		\item An associative\footnote{If $\C$ is a class of sets, we say that a family of functions $\{f_{XYZ} : X\times Y\to Z\}$ indexed by the elements $X,Y,Z \in \C$ is \emph{associative} if
			\[\notag
				f_{WZU}(w, f_{XYZ}(x,y)) = f_{ZYU}(f_{WXZ}(w,x),y)
			\]
			for every tuple $X,Y,Z,U,W$ and elements for which this is meaningful. When $f_{XYZ} = \circ$ is the composition map of a category this translates, of course, into the familiar associativity property $u\circ (v\circ w) = (u\circ v)\circ w$.} \emph{composition law}
		\[
			\circ = \circ_{\C,ABC} : \C(B,C)\times\C(A,B) \to \C(A,C) : (g,f)\mapsto g\circ f
		\]
		defined for any triple of objects $A,B,C$. The composition $\circ(g,f)$ is always denoted as an infix operator, $g\circ f := \circ(g,f)$;
		\item \label{c:tre} for every object $A\in  \C_o$ there is an arrow $\id_A\in \C(A,A)$ such that for every $A,B\in \C_o$ and $f:A\to B$ we have $f\circ \id_A=f=\id_B\circ f$.
	\end{enumtag}
\end{definition}
\begin{remark}[On composition]
	From time to time, the composition $g \circ f$ in a category may be denoted by similar `monoid\hyp{}like' infix operations as $g \cdot f$ or $g \bullet f$, or even mere juxtaposition $gf$. What will \emph{never} happen is that we denote function application and morphism composition with an infix semicolon as in $f \mathrel{;} g$.
\end{remark}
\begin{remark}[On arrows]\label{why-arrows}\index{_aaa_hom@$\hom$}
	The fact that for every $f\in \C(A,B)$ we call $A$ the \emph{domain} of $f$ and $B$ the \emph{codomain} of $f$ suggests how a morphism can be pictorially represented as an arrow $f:A\to B$ `traveling', so to speak, from the domain to the codomain. Alternative notation for the set of arrows $A\to B$ are: $\hom_\C(A,B)$, $\hom(A,B)$ (when the category $\C$ is understood from the context), or more rarely $[A,B]$.
\end{remark}
\begin{remark}[On morphism application]
	There are few illustrious exceptions to the tradition of accumulating function symbols to the left, when denoting the composition
	\[f_n \circ f_{n-1}\circ\cdots\circ f_1\]
	of a tuple of morphisms $C_{i-1} \xto{f_i} C_i$. We stick to the most common notation that function application is on the left, without further mention; in this sense, a composition $A \xto{f}B\xto{g} C$ is written $g\circ f$.

	The fact that composition is associative in a category $\C$ makes every such tuple of compositions well\hyp{}defined, and we will refer to the arrow $f_n \circ \cdots\circ f_1$, where the $f_i$ are composed according to an arbitrary parenthesisation, as \emph{the} composition of the tuple $(f_n,\dots, f_1)$.
\end{remark}
\begin{definition}[Functor]\label{def:funct}
	\index{Functor}
	Let $\C$ and $\D$ be two categories; we define a \emph{functor} $F :\C\to \D$ as a pair $(F_0, F_1)$ consisting of the following data:
	\begin{enumtag}{f}
		\item \label{f:uno} $F_0$ is a function $ \C_o\to  \D_o$ sending an object $C\in \C_o$ to an object $FC \in \D_o$;
		\item \label{f:due} $F_1$ is a family of functions $F_{AB} : \C(A,B)\to \D(FA,FB)$, one for each pair of objects $A,B\in\C_o$, sending each arrow $f:A\to B$ into an arrow $Ff:FA\to FB$, and such that:
		\begin{itemize}
			\item $F_{AA} (\id_A)=\id_{F A}$;
			\item $F_{AC} (g\circ_\C f)=F_{BC} (g)\circ_\D F_{AB} (f)$.
		\end{itemize}
	\end{enumtag}
\end{definition}
\begin{remark}
	Every family of arrows $F_{AB}$ like in \ref{def:funct}.\ref{f:uno},\ref{f:due} will be said to satisfy a \emph{functoriality property}; from now on, we will always denote both the action on objects and on arrows of $F$ with the same symbol $F$: so, $F$ sends an object $A$ into $FA$ and an arrow $f :A \to B$ into an arrow $Ff : FA \to FB$.
\end{remark}
\begin{definition}[Subcategory]\label{def:subcat}
	\index{Category!sub---}\index{Subcategory}
	Let $\C$ be a category. A \emph{subcategory} $\catS$ of $\C$ is a category defined by the following conditions
	\begin{enumtag}{sc}
		\item The objects of $\catS$ form a sub\hyp{}class of the class of objects of $\C$;
		\item For every $A,B\in\catS_o$ there is an injective function $\catS(A,B)\subseteq \C(A,B)$.
	\end{enumtag}
	If in the second condition above the inclusion is in fact an equality, or a bijection, the subcategory is called \emph{full}.
\end{definition}
\begin{definition}[Isomorphism]
	\index{Morphism!isomorphism}
	Let $\C$ be a category, and $f:A\to B$ one of its morphisms. We say that $f$ is an \emph{isomorphism} (or an \emph{invertible morphism}) if there exists a morphism $g:B\to A$ going in the opposite direction, such that $f\circ g=\id_B$  and $g\circ f=\id_A$. When $A=B$ we call an arrow $f : A \to A$ an \emph{endomorphism}, and an invertible arrow an \emph{automorphism} of $A$.
\end{definition}
\begin{remark}
	As it happens for groups, if an inverse of $f : A \to B$ exists, then it is unique: if there are two such inverses,
	\[g = g\circ \id = g\circ f\circ g' = g'\notag\]
	using associativity of the composition and the definition of inverses.\footnote{Note that this is actually showing something stronger: if $f$ has a left inverse $g$, and a right inverse $g'$, then $f$ is invertible, and $g=g'$.} We can thus call an inverse of $f$ \emph{the} inverse $f^{-1}$.
\end{remark}
\begin{example}[Examples of categories]\label{es:categ}\leavevmode
	\begin{enumtag}{c}
		\item \index{Category!monoids as ---s}
		\index{Category!empty ---}
		\index{Empty category}\label{catex:cin} Let's rule out all edge examples in a single item: the \emph{empty} category, having no objects and morphisms, satisfies all axioms of \ref{def:cate}. So does the \emph{singleton} category, having a single object $*$ and a single morphism, its identity $\id_*$. Every set $A$ can be regarded as a category $A^\delta$, having as objects the elements of $A$, and where $A^\delta(x,y)=\varnothing$ if $x\neq y$, and where there is a unique arrow $x\to x$, which must be the identity $\id_x$. This is called the \emph{discrete} category on the set $A$. In a dual fashion, every set $A$ can be regarded as a category $A^\chi$, having as objects the elements of $A$, and where $A^\chi(x,y)$ has \emph{exactly} one element for each pair $(x,y)\in A\times A$; this is called the \emph{chaotic} category on the set $A$.
		\item \index{Category!--- of sets}\label{catex:uno} The collection of sets, and functions between sets is a category $\Set$. The set $\Set(A,B)$ is the set of all functions $f :A \to B$, seen as the subset of the power\hyp{}set of $A\times B$ of those relations that are functions (\ie such that for every $a\in A$ such that $(a,b)\in f$ there is a unique such $b\in B$).
		\item \index{Category!--- of structured sets}
		The above category contains as (nonfull) subcategories those of \emph{structured sets}, having objects the sets endowed with operations like groups, rings, vector spaces, and morphisms the \emph{homomorphisms} of these structures (homomorphisms of groups, rings, linear maps of vector spaces, \dots). These categories are denoted with evocative terms like $\BF{Grp}$ (groups and their homomorphisms), $\Ab$ (\emph{abelian} groups and their homomorphisms), $\Mod(K)$ (vector spaces and linear maps), or more generally $\Mod(R)$, $\BF{Gph}$ (graphs\footnote{A \emph{graph} consists of a pair of sets $E,V$ (\emph{edges} and \emph{vertices}) such that $E\subseteq V\times V$; a graph homomorphism is a pair of functions $f_E, f_V$ between the sets of edges and vertices of two graphs. A graph $\underline{G}=(E,V)$ is \emph{directed} if each element $(e_0, e_1)$ in $E$ is an ordered pair; in this case $e_0$ is the \emph{source} or \emph{domain} of the edge and $e_1$ is its \emph{target} or \emph{codomain}. Of course, a category can be regarded as a (possibly $\mho^+$-)graph of a particular kind.} and their homomorphisms)\dots
		\item Moreover, there is a category of \emph{topological spaces}, whose objects are pairs $(X,\tau)$ ($\tau$ is a topology on the set $X$) and the morphisms $f : X\to Y$ are continuous functions.
		\item \index{Category!posets as ---s}\label{catex:due} Every partially ordered set (\emph{po}set) $P$ endowed with a relation `$\le$' can be seen as a category $\cate P$. Indeed, the objects of $P$ are its elements, and there is an arrow $x\to y$ if and only if $x\le y$ in $P$ (in all other cases, the set $P(x,y)$ is empty).
		\item \label{catex:tre} Every monoid $M$ can be regarded as a category $\cate M$ with a single object $*_M$, and where the set of arrows $m : *_M \to *_M$ is the set of elements of $M$; indeed, in this case, the monoid axioms and the axioms in \ref{def:cate} translate perfectly into each other. The group of invertible elements of $M$ identifies to the automorphism group of the unique object $*_M$. Given this, every group $G$ can be regarded as a category with a single object, such that every morphism is invertible.
		\item \label{catex:qtr} Given a directed graph $\Gamma$, we can consider the \emph{free category} generated by the graph: its class of objects is the same collection of vertices of $\Gamma$; given two vertices $U,V$ the set of morphisms $f : U\to V$ is the set of paths $U \to W_1\to\dots\to W_n \to V$ of length $n$; composition of morphisms is concatenation of paths; the identity morphism is the (unique) empty path of length $0$.
	\end{enumtag}
\end{example}
\begin{definition}\label{def:monoepi}
	\index{Morphism!monomorphism}
	\index{Morphism!epimorphism}
	Let $\C$ be a category; we say that an arrow $f : X \to Y$ is
	\begin{enumtag}{me}
		\item a \emph{monomorphism} if every two parallel arrows $u,v : A\rightrightarrows X$ such that $f \circ u = f \circ v$ are equal;
		\item an \emph{epimorphism} if every two parallel arrows $u,v : Y\rightrightarrows B$ such that $u \circ f = v \circ f$ are equal.
	\end{enumtag}
	These conditions are in turn respectively equivalent to the following ones, stated in terms of the hom\hyp{}sets of $\C$: the arrow $f : X \to Y$ is
	\begin{enumtag}{em}
		\item a monomorphism if for every object $A$ the function $\C(A,f) : \C(A,X)\to \C(A,Y) : u \mapsto f \circ u$ is injective;
		\item an epimorphism if for every object $B$ the function $\C(f,B) : \C(Y,B)\to \C(X,B) : u \mapsto u \circ f$ is injective.
	\end{enumtag}
\end{definition}
Arrangements of objects and arrows in a category are called \emph{diagrams}; to some extent, category theory is the art of making diagrams \emph{commute}, \ie the art of proving that two paths $X\to A_1\to A_2 \to\cdots\to A_n \to Y$ and $X \to B_1\to\cdots\to B_m \to Y$ result in the same arrow when they are fully composed.

We attempt at the difficult task of providing a precise formalisation of what is a commutative diagram; it is with a certain surprise that we noticed how even reliable sources as \cite{McL, acc} fail to provide more than a bland intuition for such a fundamental notion.

We begin recording an easy and informal remark.
\begin{remark}[Diagrams and their commutation]\leavevmode
	Depicting morphisms as arrows allows to draw regions of a given category $\C$ as parts of a (possibly non planar) graph; we call a \emph{diagram} such a region in $\C$, the graph whose vertices are objects of $\C$ and whose edges are morphisms of suitable domains and codomains. For example, we can consider the diagram
	\[
		\vcenter{\xymatrix{
				&&X_2\ar[drr]^v&& \\
				X_1\ar[dr]_g\ar[urr]^u\ar[rrrr]|(.375)\hole_(.7)h&&&&X_3\\
				& X_0 \ar[uur]_(.3)q\ar[rr]_p&& X_4 \ar[ur]_k&
			}}
	\]
	The presence of a composition rule in $\C$ entails that we can meaningfully compose \emph{paths} $[u_0,\dots, u_n]$ of morphisms of $\C$. In particular, we can consider diagrams having distinct paths between a fixed source and a fixed `sink' (say, in the diagram above, we can consider two different paths $\fkP = [k,p,g]$ and $\fkQ = [v,u]$); both paths go from $X_1$ to $X_3$, and we can ask the two compositions $\circ[k,p,g] = k\circ p\circ g$ and $v\circ u$ to be the same arrow $X_1\to X_3$; we say that a diagram \emph{commutes at $\fkP,\fkQ$} if this is the case; we say that a diagram \emph{commutes} (without mention of $\fkP, \fkQ$) if it commutes for every choice of paths for which this is meaningful.
\end{remark}
Searching a formalisation of this intuitive pictorial idea leads to the following:
\index{Diagram}
\begin{definition}
	A \emph{diagram} is a map of directed graphs (`digraphs') $D : J \to \C$ where $J$ is a digraph and $\C$ is the digraph underlying a category.\footnote{Every small category has an underlying graph, obtained keeping objects and arrows and forgetting all compositions; there is of course a category of graphs, and regarding a category as a graph is another example of forgetful functor. Of course, making this precise means that the collection of categories and functors form a category on its own right.} Such a diagram $D$ \emph{commutes} if for every pair of parallel edges $f,g : i\rightrightarrows j$ in $J$ one has $Df = Dg$.
\end{definition}
\begin{definition}[Full, faithful, conservative functors]\label{fullfaithconsfun}
	\index{Functor!full, faithful, conservative}
	Let $F  : \C \to \D$ be a functor between two categories.
	\begin{enumtag}{ffc}
		\item \label{fu} $F$ is called \emph{full} if each $F_{XY} : \C(X,Y)\to \D(FX,FY)$ in \ref{def:funct}.\ref{f:due}  is surjective;
		\item \label{fai} $F$ is called \emph{faithful} if each $F_{XY} : \C(X,Y)\to \D(FX,FY)$ in \ref{def:funct}.\ref{f:due} is injective;
		\item \label{fufai} $F$ is called \emph{fully faithful} if it is full and faithful;
		\item \label{cons} $F$ is called \emph{conservative} if whenever an arrow $Fv$ is an isomorphism in $\D$, then the arrow $v$ is already an isomorphism in $\C$.\footnote{Note that the fact that $F$ preserves invertibility of $v$ is a consequence of the functoriality conditions in \ref{def:funct}; one also often says that a functor $F$ is conservative if, besides preserving them, it \emph{reflects} isomorphisms. A conservative functor is thus such that $v$ is invertible if and only if $Fv$ is invertible.}
	\end{enumtag}
\end{definition}
\begin{remark}\label{pleine}
	\index{Category!full sub---}\index{Full subcategory}
	A subcategory $\catS\subseteq \C$ is full if and only if the inclusion functor $\catS\hookrightarrow\C$ is full in the sense of \ref{fullfaithconsfun}.\ref{fu} above.
\end{remark}
\section{Natural transformations}
Category theory was born in order to find the correct definition of a natural transformation.

The original motivation for this search lied in algebraic topology: there are two well-known ways to attach an homotopy invariant algebraic structure to a topological space $X$:
\begin{itemize}
	\item the homotopy groups $\pi_n(X,x)$, obtained as homotopy classes of pointed maps $S^n \to X$; these are abelian groups if $n\ge 2$, and they are notoriously from\hyp{}difficult\hyp{}to\hyp{}impossible to compute;
	\item the homology groups $H_n(X,\mathbb Z)$ with integer coefficients, a family of abelian groups way easier to compute, and very well\hyp{}behaved, but not as expressive and comprehensively descriptive as homotopy groups.
\end{itemize}
The two constructions are tightly linked: the homology group $H_n(X,\mathbb Z)$ is obtained as the $n$-th homology of the chain complex whose groups of $n$-simplices are free on the set of all continuous functions $s_n : D^n \to X$ ($D^n$ is the $n$-dimensional ball, identified up to an obvious homeomorphism with the topological $n$-simplex of \ref{realizia}); since $H_n(S^n,\mathbb Z)\cong \mathbb Z\langle u\rangle$, the homotopy a continuous function $f : S^n \to X$ determines a unique element $f_*(u):= H_n(f)(u) \in H_n(X)$, and this defines a function $\pi_n(X,x)\to H_n(X,\mathbb Z)$ by sending $f$ to $f_*u$.

Now, a capital observation is in order: this construction is compatible with the functoriality of $\pi_n, H_n$ in the following sense. If $g : X \to Y$ is a continuous function of spaces, then the square
\[\label{ureviccio}\vcenter{\xymatrix{
			\pi_n(X) \ar[r]^h \ar[d]_{g_*} & H_n(X) \ar[d]^{g_*}\\
			\pi_n(Y) \ar[r]_h & H_n(Y)
		}}\]
is commutative, so that the maps $h_X : \pi_n(X,x)\to H(X,\mathbb Z)$ `coherently' or `naturally' vary according to the action of $\pi_n,H_n$ on morphisms (since a category is just a bunch of objects linked by morphisms, it is `natural' to ask functors to preserve composition, and to families of morphisms $\alpha_X : FX \to GX$ between two functors to vary accordingly, in the same sense of \eqref{ureviccio}).

This leads directly to the notion of \emph{natural transformation}.
\begin{definition}[Natural transformation]\label{def:trasnat}
	\index{Transformation!natural ---}
	Let $F , G :\C\rightrightarrows \D$ be functors between two categories; a \emph{natural transformation} $\tau:F \To  G $ consists of a family of morphisms $\tau_X:FX\to  GX$, one for each object $X\in \C_o$, called the \emph{components} of the transformation, such that for every morphism $f:X\to Y$ the diagram
	\[\label{natusquare}
		\vcenter{\xymatrix{
				F X \ar[r]^{F (f)}\ar[d]_{\tau_X}& F Y \ar[d]^{\tau_Y}\\
				G X \ar[r]_{ G (f)}&  G Y
			}}
	\]
	commutes, \ie we have the equation $\tau_Y\circ F (f)= G (f)\circ \tau_X$.
\end{definition}
\begin{definition}[Natural equivalence]\label{def:equnat}
	\index{Natural equivalence}
	\index{Equivalence!---  of functors}
	Let $F , G :\C\rightrightarrows\D$ be two parallel functors. A natural transformation $\tau:F \To G $ such that every component is an isomorphism in $\D$ is called a \emph{natural equivalence} or (less often) an \emph{isomorphism of functors}, or a \emph{functorial isomorphism}.

	Note that if $\tau:F \To G$ is a natural transformation, and each component $\tau_C : FC \to GC$ is an isomorphism in $\D$, then the family of maps $\tau_C^{-1}$ is also natural.
\end{definition}
\begin{definition}[Whiskering]\label{comporizz}
	\index{Whiskering}
	Let $\tau : F \To G$ be a natural transformation between functors $F , G : \C\rightrightarrows\D$; given a third functor $H : \mathcal H\to\C$, we define the natural transformations
	\[
		\tau * H : FH\To  GH : (\tau * H)_C = \tau_{ HC} : FHC\to GHC
	\]
	Similarly, given a functor $ K : \D\to\K$ we can define the natural transformation
	\[
		K * \tau : KF \To KG : ( K *\tau)_C= K(\tau_C) : KFC\to KGC
	\]
\end{definition}
It is clear that the two equations
\begin{gather}
	(\tau * H) * H' = \tau * ( H\circ H')\\
	K' * ( K *  \tau)=( K'\circ K) * \tau
\end{gather}
hold true for every $\tau, H, H', K, K'$ that make them meaningful. It is equally obvious that $(K*\tau)*H=K*(\tau *H)$. This operation is called \emph{whiskering}, since it acts `drawing whiskers' on the left and on the right of $\tau$:
\[
	\vcenter{\xymatrix{
	\cate H \ar[r]^H & \C \rtwocell^F_G{\tau} & \D \ar[r]^K & \cate K
	}}
\]
Finally, we make precise the idea that categories and functors form a category in their own right:
\begin{remark}[The category of functors]
	Let $\C,\D$ be two \emph{small} categories; then the functors $F : \C\to \D$ form the object class of a category, whose morphisms $\alpha : F\To G$ are the natural transformations. Two natural transformations can be joined component\hyp{}wise, in such a way that if $F  \overset{\alpha}\To G \overset{\beta}\To H$, we have that the components of the composite transformation $\beta\circ\alpha$ are
	\[
		(\beta\circ\alpha)_C=\beta_C\circ\alpha_C.
	\]
	The naturality square is of course
	\[
		\vcenter{
		\xymatrix{
		FC\ar[r]^{\alpha_C}\ar[d]_{Ff} & GC\ar[r]^{\beta_C}\ar[d]_{Gf} & HC \ar[d]^{Hf}\\
		FC' \ar[r]_{\alpha_{C'}}& GC' \ar[r]_{\beta_{C'}}& HC'
		}
		}
	\]
	\index{_aaa_Cat@$\Cat$}
	The identity in $\Cat(\C,\D)(F ,F)$ is the natural transformation having object\hyp{}wise identity components. This is called the \emph{vertical} composition of natural transformations, because $\alpha,\beta$ can be arranged in a diagram
	\[\vcenter{\xymatrix{
				\C \ruppertwocell^{F}{\alpha} \ar[r]\rlowertwocell_{H}{\beta}& \D
			}}\]
	and the natural transformation $\beta\circ\alpha$ is the result of stacking $\alpha,\beta$ one on top of the other.
\end{remark}
\index{Horizontal composition}
\index{_aaa_boxminus@$\boxminus$}\index{Natural transformation!composition of ---s}
Now, we shall show that there is another way to compose the natural transformations $\beta$ and $\alpha$ called \emph{horizontal composition}: let $\A, \B,\C$ be three categories, and $F,G,H,K$ functors arranged as follows:
\[
	\vcenter{\xymatrix{\A \rtwocell^{F}_{G}{\alpha} & \B \rtwocell^{H}_{K}{\beta} & \C }}
\]
\begin{definition}
	Given categories $\A, \B,\C$, functors $F , G : \A\to \B$ and $ H, K : \B\to \C$ we can define the \emph{horizontal composition} of natural transformations $\alpha : F \To  G $ and $\beta : H\to K$ to be the natural transformation $\beta \boxminus\alpha$, whose components are defined thanks to the fact that
	$(\beta *  G )\circ ( H * \alpha) = ( K * \alpha)\circ (\beta * F )$ (as it is easy to check):
	\[(\alpha\boxminus\beta)_X = \beta_{ G X}\circ H(\alpha_X) = K(\alpha_X)\circ \beta_{F X}.\]
\end{definition}
Applying the definition, we can show that the horizontal composition is (well defined and) associative, in the sense of \ref{def:cate}.(3).
\index{Natural transformation!composition of ---s}
\begin{definition}[Categorical equivalence]\label{def:equcat}
	\index{Equivalence!--- of categories}
	Let $\C,\D$ be two categories. An \emph{equivalence of categories} $(F,G,\xi,\eta)$ between $\C$ and $\D$ consists of two functors $F :\C\to \D, G :\D\to \C$ endowed with two natural equivalences $\xi : F \circ  G \to \id_\C$ e $\eta: G \circ F \To \id_\D$ between the two compositions of $F $ and $G$ and the respective identity functors of $\C$ and $\D$ (see \ref{def:funct}.\ref{f:uno}). If this is the case, the categories $\C$ and $\D$ are called \emph{equivalent}; we write $\C\cong \D$.
\end{definition}
\begin{remark}\label{adji_are_contrattible}
	It is quite common to denote an equivalence with the only functor $F :\C\to \D$, and to say that it is an equivalence if there exists a functor $G$ in the opposite direction, and natural transformations $\xi,\eta$ such that the tuple $(F,G,\xi,\eta)$ is an equivalence; this is customary and harmless since (as we will observe in \ref{adjbasta}) such a $G$ is unique up to natural isomorphism, provided it exists.
\end{remark}
The richness of category theory is to some extent due to the fact that the correct way to assert `sameness' for two categories is the above definition of equivalence, and not the stricter notion of isomorphism:
\begin{definition}[Isomorphism of categories]\label{def:isocat}
	\index{Isomorphism of categories}
	We call two categories $\C,\D$ \emph{isomorphic} if there is an invertible functor (called an \emph{isomorphism} between the two categories) $F :\C\to \D$; this means that there exists $G :\D\to \C$ with the property that $F \circ  G =\id_\D$ and $ G \circ F =\id_\C$.
\end{definition}
In a certain deep sense, equivalences and isomorphisms of categories stand in the same relation as \emph{homotopy equivalences} and \emph{homeomorphisms} of topological spaces (this path leads directly to higher category theory. Another example is the following: observe that in \ref{adji_are_contrattible} we are saying that the `space' of inverses to $F : \C \to \D$ is either empty or contractible).

\index{Natural transformation!composition of ---s}

\begin{proposition}\label{equicrit}
	Let $\C,\D$ two categories and $F  : \C\to \D$ a functor. Then, $F$ is (part of) an equivalence of categories if and only if it is fully faithful and \emph{essentially surjective}, \ie if every object in $\D$ is isomorphic to some $FC$ lying in the image of $F$.
\end{proposition}
\subsection{Duality and slices}\label{duality_and_slices}
As for every other algebraic structure, there is plenty of ways we can obtain new categories out of given ones.
\begin{itemize}
	\item First of all, we observe that the shape of the axioms in \ref{def:cate} entails that the structure obtained taking the same object of a category $\C$, but where `all arrows have been reversed' remains a category $\C^\opp$. This is called the \emph{dual}, or \emph{opposite} category of $\C$.
	\item Then, we notice how given an object $C$ of $\C$, the class of all arrows with fixed co\fshyp{}domain $C$ becomes a category, called the \emph{co\fshyp{}slice} of $\C$.
	\item Finally, as for every other algebraic structure, two categories can be arranged in a \emph{cartesian product} and a \emph{disjoint union}: these last two examples will appear as part of a fairly more general theory of \emph{co\fshyp{}limits} in \ref{lim_and_colims}.
\end{itemize}
\index{Category!dual|see{opposite}}
\index{Category!opposite ---}\index{Opposite category}
The first procedure we introduce is \emph{dualisation}: the \emph{dual}, or \emph{opposite} category of $\C$ is made by the same objects, and the same class of arrows, but we have interchanged domain and codomain. If an arrow is represented pictorially, as $f : A\to B$, then its namesake in $\C^\opp$ is $f^\opp : B\to A$.

To be more precise, there exists a functorial correspondence of $\Cat$ into itself called \emph{duality involution} and denoted $(\firstblank)^\opp : \Cat \to \Cat$ that has the following properties: it is the identity on objects, and it is the `swapping' involution on arrows, as soon as morphisms $f :A\to B$ are identified with triples $(A, B, f)$;
\[
	\vcenter{\xymatrix{A\ar[r]^f \ar[dr]_{g\circ f}& B\ar[d]^g\\ &C}}
	\qquad
	\overset{(\firstblank)^\opp}\longmapsto
	\vcenter{\xymatrix{A & B \ar[l]_{ f^\opp}\\ & C\ar[u]_{g^\opp}\ar[ul]^{f^\opp\circ_\opp g^\opp}}}
\]
\begin{remark}\label{duality-from-acc}
	A commutative triangle in $\C$ as in the left diagram gets modified as in the right diagram by the duality involution; the object $\C^\opp$ thus defined satisfies all the axioms of a category, as stated in \ref{def:cate}. The reason why this construction is interesting is that every assertion made in $\C$ has a `companion' in $\C^\opp$: in the words of \cite{acc},
	\begin{quote}
		The concept of category is well balanced, which allows an economical and useful duality.
		Thus in category theory the `two for the price of one' principle holds: every concept is
		two concepts, and every result is two results.
	\end{quote}
	In short, the situation goes as follows: in order to define a category, axioms regarding certain indefinite notions (object, morphism, domain, codomain\dots); now, every statement $\varphi$ made in the language of category theory, involving solely relations between objects of $\C$ and the notions of object, arrow, domain, codomain, identity, composition is valid in said language if and only if the statement $\varphi^\opp$ obtained from $\varphi$ substituting each occurrence of `$g\circ f$' with `$f^\opp\circ g^\opp$', and every occurrence of \emph{domain} (resp., \emph{codomain}) with \emph{codomain} (resp. \emph{domain}).
\end{remark}
Of course, this does not mean that a statement $\varphi$ about a commutative diagram is true in $\C$ if and only if $\varphi^\opp$ is true in $\C$! Exercise \ref{set-aint-autodual} will show that (for example) the opposite category of $\Set$ is not equivalent to $\Set$, exhibiting a property that $\Set^\opp$ does enjoy, while $\Set$ does not.
\begin{definition}[Contravariant Functor]\label{def:contfun}
	\index{Functor}\index{Functor!contravariant ---}
	Let $\C,\D$ be two categories. A \emph{contravariant functor} is a functor $F :\C^\opp\to \D$ (in the sense of \ref{def:funct}); more explciitly, a contravariant functor amounts to the same data of \ref{def:funct}, where the second condition in \ref{f:due} is replaced by $F (g\circ f)=F f\circ F g$, for every pair of composable morphisms $f,g$.
\end{definition}
Functors in the sense of \ref{def:funct} are called \emph{covariant}.
\begin{remark}
	If we call $(-)^\opp:\C\mapsto \C^\opp$ the duality involution, its functoriality amounts to say that every functor $F :\C\to \D$ induces a functor $F^\opp : \C^\opp\to \D^\opp$ (so a contravariant functor $F : \C\to \D^\opp$):
	\[
		F ^\opp (g\circ f)=F ((g\circ f)^\opp)=F (g^\opp\circ f^\opp)=F ^\opp f\circ F ^\opp g
	\]
	Similarly, to every natural transformation $\tau : F\To G$ is associated a natural transformation $\tau^\opp  : G^\opp \To F^\opp $ (write down the relevant diagram and check that the components of $\tau$ are indeed reversed).
\end{remark}
\begin{definition}[Slice categories]\label{defn:slice}
	\index{Slice category|see{Category}}
	\index{Category!slice ---}\index{Slice category}
	\index{_aaa_C/C@$\C/C$}
	Let $\C$ be a category and $C\in \C$ be an object of $\C$; we define
	\begin{itemize}
		\item the `slice' category $\C/C$ of \emph{arrows over} $C$ having class of objects the arrows with codomain $C$, and morphisms between $f :C'\to C$ and $g : C''\to C$ the arrows $h : C'\to C''$ such that $g\circ h=f$.
		\item the `co\hyp{}slice' category $C/\C$ of \emph{arrows under} $C$ having class of objects the arrows with domain $C$, and morphisms between $f :C\to C'$ and $g : C\to C''$ the arrows $h : C'\to C''$ such that $h\circ f=g$.
	\end{itemize}
\end{definition}
Note that $C/\C\cong (\C^\opp/C)^\opp$ and similarly $\C/C = (C/\C^\opp)^\opp$.
\begin{definition}[Comma categories]\label{def:comma}
	\index{Comma category}\index{Category!comma ---}
	The \emph{comma category} of a diagram $\catS \xto{F} \C \xot{G}\cT$ of functors is the category having objects the tuples $(S\in\cate S,T\in\cate T,\phi \in \C(FS, GT))$, and morphisms the pairs $(u:S\to S', v:T\to T')$ such that the square
	\[
		\vcenter{\xymatrix{
				FS \ar[r]^{Fu}\ar[d]_\varphi & FS'\ar[d]^{\varphi'} \\
				GT \ar[r]_{Gv}& GT'
			}}
	\]
	commutes. This category is denoted $(F\downarrow G)$ or $(F/G)$.
\end{definition}
The reader is invited to study a few properties of the comma category $(F/G)$ in Exercise \ref{da_comma}; we now turn our attention to the \emph{product} and \emph{coproduct} of categories.
\begin{definition}[Product of categories]\label{prodocat}
	\index{Category!product ---}\index{Product of categories}
	Let $\C,\D$ be two categories; we define the \emph{product category} $\C\times \D$ to be the following structure:
	\begin{itemize}
		\item it has objects the product of the classes of objects $\C_o\times \D_o$;
		\item it has as set of morphisms $(C,D)\to (C',D')$ the cartesian product $\C(C,C')\times \D(D,D')$.
	\end{itemize}
	All the axioms of a category are satisfied, as compositions and identities are defined factor\hyp{}wise and thus act independently one from the other. Of course, the definition inductively extends to the product of any finite number of categories $\C_1\times \cdots\times \C_n$.
\end{definition}
Dually, we can define
\begin{definition}[Coproduct of categories]\label{coprodocat}
	\index{Category!coproduct ---}\index{Coproduct of categories}
	Let $\C,\D$ be two categories; we define the \emph{coproduct category} $\C\coprod \D$ to be the following structure:
	\begin{itemize}
		\item it has objects the set\hyp{}theoretical disjoint union of the classes of objects $\C_o\times \D_o$;
		\item it has as set of morphisms  $X\to Y$ the set $\C(X,Y)$ if $X,Y\in\C$, $\D(X,Y)$ if $X,Y\in\D$, and the empty set otherwise.
	\end{itemize}
	All the axioms of a category are satisfied, as compositions and identities are defined summand\hyp{}wise and thus act independently one from the other. Of course, the definition inductively extends to the coproduct of any finite number of categories $\C_1\amalg \cdots\amalg \C_n$.
\end{definition}
Observe tat the following property holds true for the product $\C\times \D$:
\begin{quote}
	Let $\E$ be a category and $\C \xot{F} \E \xto{G} \D$ be two functors; then there exists a unique functor $\E \to \C\times \D$ such that the compositions $\E \to \C\times \D \xto{P} \C$ and $\E \to \C\times \D \xto{Q}\D$ with the projections on the factors of the product correspond respectively to $F$ and $G$.
\end{quote}
Observe that $\C\coprod\D$ enjoys the same property in the dual category $\Cat^\opp$, or rather, the dual property in $\Cat$:
\begin{quote}
	Let $\E$ be a category and $\C \xto{F} \E \xot{G} \D$ be two functors; then there exists a unique functor $\C\coprod \D \to \E$ such that the compositions $\C \to \C\coprod \D \to \E$ and $\D \to \C\coprod \D \to \E$ with the embeddings into the coproduct correspond respectively to $F$ and $G$.
\end{quote}
The scope of the following section is to formalise the intuition that these two construction are `universal' and `dual' to each other.
\section{Limits and colimits}\label{lim_and_colims}\index{Limit}
Right after having introduced the definition of an algebraic structure, it is often shown that given one or more such structures one can build others from it; one can for example assemble the \emph{product} of two groups, sets or topological spaces.

More than often the objects that are built in this way are characterised by some sort of uniqueness; for example, there is only one way to assemble two vector spaces $V,W$ into a third space $V\oplus W$ that contains `no more than $V,W$'; this vector space is made by pairs of vectors $(v,w)$, and the vector sum $(v,w)+(v',w')$, as well as scalar multiplication $\alpha(v,w)$, are done component\hyp{}wise. Moreover, there is a diagram
\[
	\vcenter{\xymatrix{ V \ar[r]^-{i_V} & V\oplus W & W\ar[l]_-{i_W}}}
\]
such that \emph{every other} similar diagram $V \xto{f} U \xot{g} W$ factors through the first: there is a unique $[f/\!\!/g] : V\oplus W \to U$ that coincides with $f,g$ when restricted to the summands, \ie such that $[f/\!\!/g]\circ i_V = f$ and $[f/\!\!/g]\circ i_W = g$.

Category theory provides a very neat way to organise these data and explain what are the essential features of this phenomenon, through the theory of co\fshyp{}limits.

We assume that our reader is already familiar with elementary algebra, and in particular that they are familiar with the simple idea that algebraic structures of a same kind (e.g., groups or vector spaces) can be assembled together.
\begin{definition}[Initial and terminal objects]\label{def:inifin}\index{Object!initial and terminal}
	Let $\C$ be a category.
	\begin{itemize}
		\item An object $\varnothing\in  \C_o$ is called \emph{initial} if for every other object $C\in  \C_o$ there is a single morphism $i_C : \varnothing\to C$.
		\item Dually, an object $1$ is called \emph{final} o \emph{terminal} if for every object $C\in  \C_o$ there is a unique morphism $t_C : C\to 1$.
	\end{itemize}
	(Note the substantial difference between `there is at most one morphism' and `there is exactly one morphism'!)
\end{definition}
\begin{remark}
	As a consequence of the definition, if $\varnothing$ is an initial object in $\C$, then there is a single arrow $i_\varnothing : \varnothing\to \varnothing$, the identity of $\varnothing$. Similarly, if $1$ is terminal, there is a unique arrow $t_1 : 1\to 1$, the identity of $1$; if $\C$ has both an initial and a terminal object, then there is a unique arrow $z : \varnothing\to 1$; we say that $\C$ has a \emph{zero} object if $z$ is an isomorphism.
\end{remark}
The simple proof of the following statement will enlighten the nature of the notion of universal property. An initial object $\varnothing\in\C$ enjoys what is called a \emph{universal property}:
\begin{remark}
	Let $\C$ be a category with an initial object $\varnothing$. If $\varnothing'$ is another initial object, then there is a unique isomorphism $\varnothing\cong \varnothing'$.
\end{remark}
\begin{proof}
	Assume that there are two initial objects $\varnothing, \varnothing'$; then, by the respective universal properties of $\varnothing$ and $\varnothing'$, there is a unique arrow $u : \varnothing \to \varnothing'$, and similarly a unique arrow $v : \varnothing'\to \varnothing$. The compositions $v\circ u$ and $u\circ v$ must be the identities of $\varnothing$ and $\varnothing'$ respectively, thus showing that $u,v$ are mutually inverse isomorphisms.
\end{proof}
The notion of \emph{universal property} arises to generalise this phenomenon; if $\J$ is a category, a functor $D : \J\to \D$ can be thought as a diagram (`of shape $\J$') representing the category $\J$ in its codomain; for example if $\J$ is the category having three objects $0,1,2$ and whose nonidentity morphisms are
\[
	\vcenter{\xymatrix{ &0\ar[dr]\ar[dl]&\\
			1 && 2}}
\] then a diagram of shape $\J$ is simply a triple of objects $D_0,D_1,D_2$ linked by similar morphisms $D_0\to D_1,D_0\to D_2$.
The key observation of the theory of co\fshyp{}limits is that to each diagram $D$ one can associate a category of \emph{cones} (and dually, a category of \emph{cocones}), obtained as suitable extensions of $D$ to a `cone category' $\tilde \J \supset \J$, adding to $\J$ a new initial (or terminal) object.

Extensions of this kind can be organised in categories $\J^\rhd$ if we add a terminal object, and $\J^\lhd$ if we add an initial object, of co/cones over $\J$.\footnote{The reader will surely appreciate the origin of the name: if $X$ is a topological space, the \emph{cone} of $X$ is the space obtained adding a distinguished point $\infty$, disconnected from $X$ and then adding a path from $\infty$ to each point of $X$.}

This terse procedure directly leads to the following
\begin{definition*}
	The \emph{limit} of a diagram $D$ is (whenever it exists) the \emph{terminal} object in the category of its \emph{cones}, and dually the colimit of $D$ is the \emph{initial} object in its category of \emph{cocones}.
\end{definition*}
Of course, we will not leave the reader alone deciphering this mysterious axiomatics; let us start to make this construction precise by defining the
\index{Object!initial and terminal}
\begin{definition}[Cone completions of $\J$]\label{conecompletion}\index{Cone!completion}\index{_aaa_Jlhd@$\J^\lhd, \J^\rhd$}
	Let $\J$ be a small category; we denote $\J^\rhd$ the category obtained adding to $\J$ a single terminal object $\infty$; more in detail, $\J^\rhd$ has objects $\J_o \cup\{\infty\}$, where $\infty\notin \J$, and it is defined by
	\begin{align*}
		{\J^\rhd}(J,J')     & = \J(J,J') \\
		{\J^\rhd}(J,\infty) & = \{*\}
	\end{align*}
	and it is empty otherwise. This category is called the \emph{right cone} of $\J$.

	Dually, we define a category $\J^\lhd$, the \emph{left cone} of $\J$, as the category obtained adding to $J$ a single \emph{initial} object $-\infty$; this means that ${\J^\lhd}(J,J')=\J(J,J')$, ${\J^\lhd}(-\infty, J)=\{*\}$, and it is empty otherwise.
\end{definition}

\begin{definition}[Cone of a diagram]\label{daskone}\index{Diagram!Cone of a ---}
	Let $\J$ be a small category, $\C$ a category, and $D : \J \to \C$ a functor; all along this section, an idiosyncratic way to refer to $D$ will be as a \emph{diagram of shape $\J$}. We call a \emph{cone for $D$} any extension of the diagram $D$ to the left cone category of $\J$ defined in \ref{conecompletion}, so that the diagram
	\[\label{cone-as-lift}
		\vcenter{\xymatrix{
				\J \ar[r]^D \ar[d]_{i_\lhd} & \C \\
				\J^\lhd\ar@{.>}[ur]_{\bar D}
			}}
	\]
	commutes.
\end{definition}
Every such extension is thus forced to coincide with $D$ on all objects in $\J\subseteq \J^\lhd$; the value of $\bar D$ on $-\infty$ is called the base of the cone; dually, the value of an extension of $D$ to $\J^\rhd$ coincides with $D$ on $\J\subseteq \J^\rhd$, and $\bar D(\infty)$ is called the \emph{tip} of the cone.

There is of course a similar definition of a \emph{cocone} for $D$: it is an extension of the diagram $D$ to the right cone category of $\J$ so that the diagram
\[\label{cone-as-lift2}
	\vcenter{\xymatrix{
			\J \ar[r]^D \ar[d]_{i_\rhd} & \C \\
			\J^\rhd\ar@{.>}[ur]_{\bar D}
		}}
\]
commutes.
\begin{remark}\leavevmode
	\begin{itemize}
		\item The class of cones for $D$ forms a category $\Cn(D)$, whose morphisms are the natural transformations $\alpha : D'\To D'' : \J \to \C$ such that the right whiskering of $\alpha$ with $i : \J \to \J^\lhd$ coincides with the identity natural transformation of $D$; this means that a morphism $\alpha$ of this sort is a natural transformation such that
		      \[
			      \vcenter{\xymatrix{
					      &\J\ar[dl]_i\ar[dr]^D & \\
					      \J^\lhd \rrtwocell^{D'}_{D''}{\alpha}&& \C
				      }}	\quad = \id_D
		      \]
		      as a 2-cell
		      $D\To D$.	\item Dually, the class of cocones for $D$ forms a category $\Ccn(D)$, whose morphisms are the natural transformations $\alpha : D'\To D'' : \J \to \C$ such that the right whiskering of $\alpha$ with $i : \J \to \J^\rhd$ coincides with the identity natural transformation of $D$; this means that a morphism $\alpha$ of this sort is a natural transformation such that
		      \[
			      \vcenter{\xymatrix{
					      &\J\ar[dl]_i\ar[dr]^D & \\
					      \J^\rhd \rrtwocell^{D'}_{D''}{\alpha}&& \C
				      }}	\quad = \id_D
		      \]
		      as a 2-cell $D\To D$.
	\end{itemize}
\end{remark}
\begin{definition}[Colimit, limit]\label{colimlim}\index{Limit and colimit}
	The \emph{limit} of a diagram $D : \J\to \C$ is the terminal object denoted `$\lim_\J D$' in the category of cones for $D$; dually, the \emph{colimit} of $D$ is the initial object denoted `$\colim_\J D$' in the category of cocones for $D$.
\end{definition}
\begin{proposition} \label{lims_iff_pis_and_eqs}
	Let $\C$ be a category; the following conditions are equivalent:
	\begin{enumtag}{l}
		\item $\C$ has all limits of shape $\J$;
		\item $\C$ has products indexed by every set $\J$, and equalisers;
		\item $\C$ has a terminal object, and pullback of every co-span $X \to Z \leftarrow Y$.
	\end{enumtag}
	Dually, the following conditions are equivalent:
	\begin{enumtag}{c}
		\item $\C$ has all colimits of shape $\J$;
		\item $\C$ has coproducts indexed by every set $\J$, and coequalisers;
		\item $\C$ has an initial object, and pushout of every span $X \leftarrow Z \to Y$.
	\end{enumtag}
\end{proposition}
\begin{proof}
	See \cite[2.8.2]{Bor1} (and a suitable dual statement). We only record that the limit of a diagram $D$ has the same universal property of the equaliser of the pair
	\[
		u,v :\prod_A DA \rightrightarrows \prod_{f : A \to B}DB
	\]
	where the arrows $(u,v)$ are defined as
	\begin{gather*}
		\pi_{(f : A \to B)}\circ u =\pi_B\\
		\pi_{(f : A \to B)}\circ v = (Df)\pi_A
	\end{gather*}
	if $\pi_X : \prod_A DA \to DX$ denotes the projection at $X$ coordinate. The terminal cone exhibiting $\lim D$ is of course the composition $\lim D \to \prod_A DA \to DA$.

	Dually, the coequaliser of
	\[
		(u', v') : \coprod_{f : A \to B} DA \rightrightarrows \coprod_A DA
	\]
	for a similar pair of arrows $(u',v')$ has the same universal property of $\colim D$.
\end{proof}
\begin{definition}[Preservation of co\fshyp{}limits]\label{preservation}\index{Limit!preservation of ---}
	Let $\J,\C,\D$ be categories, $\J$ small, $D : \J \to \C$ be a diagram, and $F : \C \to \D$ a functor. Assume that the limit $\lim_\J D $ exists in $\C$; then, applying the functor $F$ to the terminal cone $\lambda : \J^\lhd \to \C$ of $D$, we get a cone $F * \lambda$ for the composed diagram $\J \xto{D}\C \xto{F} \D$. We say that $F$ \emph{preserves} the limit of $D$ if $F * \lambda$ is the limit of $F\circ D$, \ie if
	\[\textstyle \lim_\J(F\circ D) \cong F(\lim_\J D).\]
\end{definition}
\begin{definition}[Co/complete category]\label{being-complete}\index{Limit!Co/complete category}
	\index{Category!co\fshyp{}complete ---}\index{Co/complete category}
	A category $\C$ \emph{has all limits of shape $\J$} if every diagram $D : \J\to \C$ has a limit; dually, we define a category having all colimits of shape $\J$. A category is said to \emph{have all $\kappa$\hyp{}co\fshyp{}limits} if it has co\fshyp{}limits of shape $\J$ for every category $\J$ with less than $\kappa$ objects. A category is said to have \emph{all small co\fshyp{}limits}, or simply \emph{all co\fshyp{}limits} (but the implicit smallness request on $\J$ is needed!) if every diagram $D: J \to \C$ with small domain has a co\fshyp{}limit.\footnote{Every category admitting `too large' co\fshyp{}limits must be a (large) poset; this is a theorem by P. Freyd, and constitutes the main reason why the smallness assumption on $\J$ can not be dropped.}
\end{definition}
In view of \ref{lims_iff_pis_and_eqs} above, a category $\C$ has $\kappa$\hyp{}co\fshyp{}limits if and only if it has co\fshyp{}equalisers and co\fshyp{}products of every family of less than $\kappa$ objects.

\begin{proposition}
	The category $\Cat$ of small categories and functors has all small limits and colimits.
\end{proposition}
\begin{definition}\label{def:join}\index{Category!Join of two ---s}\index{Join operation}
	Let $\C,\D$ be two categories; we define the \emph{join} of $\C$ and $\D$ to be the category having objects the disjoint union $\C_o \coprod \D_o$ of $\C_o$ and $\D_o$, and where morphisms from $X$ to $Y$ are defined as follows:
	\[
		\begin{cases}
			\C(X,Y) & \text{ if } X,Y\in \C      \\
			\D(X,Y) & \text{ if } X,Y\in \D      \\
			\{*\}   & \text{ if } X\in\C,Y\in \D
		\end{cases}
	\]
	where $\{*\}$ is a singleton set. In every other case, there are no morphisms $X \to Y$.
\end{definition}
\begin{remark}
	The join operation gives the category $\Cat$ a monoidal structure, having the empty category as monoidal unit; we shall observe that the \index{Associator!monoidal ---}associator is easily determined, but the structure is highly nonsymmetric.

	As already observed, in the terminology of Exercise \ref{collage} the join of $\C$ and $\D$ coincides with the \emph{collage} of $\C$ and $\D$ along the terminal profunctor $* : \C \pto \D$, \ie with the category of elements of the functor $ * : \C^\opp\times \D\to \Set$ constant in the singleton set. The join operation can thus be regarded as a particular case of an operation in the bicategory of profunctors.
\end{remark}
\begin{definition}[\protect{\cite[3.3]{riehlcontext}} Creation of co\fshyp{}limits] \index{Co/limit!creation of ---s}\label{creation}
	Let $F : \C \to \D$ be a functor; we say that $F$ \emph{creates} the co\fshyp{}limit of a diagram $D : \J \to \C$ if for every co\fshyp{}limit co\fshyp{}cone for $FD$ in $\D$ there exists a co\fshyp{}limit co\fshyp{}cone for $D$ in $\C$, and every co\fshyp{}cone for $D$ that is sent to a co\fshyp{}limit co\fshyp{}cone by $F$ is itself a co\fshyp{}limit co\fshyp{}cone.
\end{definition}
The above definition admits a simple restriction in case $D$ is the empty diagram: a functor $F : \C \to \D$
\begin{itemize}
	\item preserves the initial object if $F(\varnothing_\C)$ is an initial object in $\D$;
	\item reflects the initial object if the fact that $F(A)$ is initial in $\D$ entails that $A$ is initial in $\C$;
	\item creates the initial object if the existence of an initial object in $\D$ of the form $FA$ entails the existence of an initial object in $\C$, and every object such that $FA=\varnothing$ is itself initial.
\end{itemize}
Of course, the word `initial' can be replaced with `terminal', obtaining preservation, reflection and creation of terminal objects.

The notion of preserving, reflecting and creating an initial object is sufficient to capture \ref{creation} by virtue of the same argument in \ref{adj-preserve}: $F$ creates the colimit of $D$ if and only if the functor $F_* : \Ccn(D)\to \Ccn(FD)$ creates the initial object. A simple dualisation yields that $F$ creates the limit of $D$ if and only if the functor $F_* : \Cn(D)\to \Cn(FD)$ creates the terminal object.

\section{Adjunctions}\label{adjunctions}
The notion of adjunction lies at the very core of category theory.
\begin{definition}[Adjunction]\label{def:adj}\index{Adjunction}\index{Functor!adjoint ---s}
	Let $\C,\D$ be two categories. We define an \emph{adjunction} between $\C,\D$ to be a pair of functors $F :\C\to \D$, $ G :\D\to \C$ endowed with a collection of bijections $\phi_{CD}$,
	\[
		\label{eqn:adj}
		\phi_{CD} : \D(FC,D)\cong \C(C, G D)
	\]
	one for each $C\in \C_o$, $D\in  \D_o$, natural in both arguments $C,D$.
\end{definition}
We denote the presence of an adjunction $(F,G)$ between $\C$ and $\D$ writing
\[F :\C\leftrightarrows \D: G \]
and we say that $F$ is a \emph{left adjoint} to $G$, or that $G$ is a \emph{right adjoint} to $F$, or that $(F,G)$ is an adjoint pair, or that $F,G$ are mutually adjoint, etc.

Given an adjunction as in \ref{def:adj}, we write \index{_aaa_dashv@$\dashv$} $F \dashv  G$ to shortly denote this (highly nonsymmetric) relation among $F$ and $G$.
Once we will have introduced the Yoneda lemma in \ref{lem:the-real-yoda}, we will be able to instantly deduce the following uniqueness result:
\begin{proposition}\label{adjbasta}
	If a functor $F  : \C\to \D$ has a left adjoint $ G : \D\to \C$, such adjoint is unique up to a unique natural isomorphism; in other words, if we have adjunctions $F\dashv G$ and $F\dashv G'$ then there exists a \emph{unique} natural isomorphism $\tau : G\cong G'$.
\end{proposition}
Of course, a similar result holds for the uniqueness of a left adjoint $F$.
\begin{example}[Examples of adjunctions]\label{ex:adjoints}\leavevmode
	\begin{enumtag}{ad}
		\item \label{ad:uno} Let $f : M \to N$ be a morphism of monoids; according to \ref{es:categ}.\ref{catex:tre} we can regard it as a functor between one\hyp{}object categories; it is easy to see that such a functor has a left adjoint $g : N \to M$ if and only if it has a right adjoint, if and only if it is an equivalence of categories, if and only if it is an isomorphism of monoids.
		\item \label{ad:due} Let $P$ be a poset regarded as a category: an adjunction $f : P \leftrightarrows Q :g$ consists of a pair of monotone functions $f,g$ such that $fp \le q$ if and only if $p \le gq$; these pairs of monotone mappings are called \emph{Galois connections} \index{Galois connection} (the name is motivated by exercise \ref{ex:galua}).

		\item \index{Reflective subcategory}\index{Category!reflective sub---}\label{ad:qtr} More generally, let $i : \cate S \hookrightarrow \C$ be the embedding of a full subcategory. We say that $\cate S$ is \emph{reflective} into $\C$, or a \emph{reflective subcategory} of $\C$, if $i$ admits a left adjoint $L : \C \to \cate S$. Almost all familiar subcategories of algebraic structures (magmas, semigroups, monoids, groups, abelian groups, $R$\hyp{}modules, \dots) fit into an adjunction with the category of sets, but almost none of these categories is reflective; this is because the forgetful functor $U : \C \to \Set$ is rarely full. (Difficult exercise: generalise the construction of unit and counit maps to an abstract algebraic structure.)
		\item \label{ad:cin} Given three sets $A,B,C$ there exists a bijection between the set $\Set(A\times B,C)$ of functions $f : A\times B\to C$ and the set $\Set(A,C^B)$ of functions from $A$ to $\Set(B,C)$; given $f : A\times B\to C$, we can define a \emph{transposed} (or \emph{curried}) function $\hat f : A\to C^B$ sending $a\in A$ into the function $f_a=f(a,\firstblank) : B\to C$. Currying a function $f$ into $\hat f$ is a bijection, and the inverse sends a function $g : A\to C^B$ into $\tilde g : A\times B \to C$ via an \emph{uncurrying} operation.\index{Currying}\index{Uncurrying|see{Currying}}

		If for every $B\in \Sets$ we define $(\firstblank)\times B : \Sets\to \Sets : A\mapsto A\times B$ and $ (\firstblank)^B : \Sets\to \Sets : C\mapsto C^B$ we get two functors (defined accordingly on arrows) that form an adjunction $\firstblank\times B \dashv (\firstblank)^B$.

		Every category $\C$ having finite products and exhibiting the same adjunction
		\[\C(A\times B, C) \cong \C(A, C^B)\cong \C(B, C^A)\]
		is called \emph{cartesian closed}; the category of sets, and the category of categories are both cartesian closed; no category with a zero object can be cartesian closed without being the terminal category $\boldsymbol *$; the category of all topological spaces and continuous maps is notoriously not cartesian closed, but many subcategories of `nice' spaces are.\index{Category!closed ---}
		\item \index{Object!initial and terminal} \label{ad:sei} Let $\C$ be a category; then $\C$ has a terminal object (see \ref{def:inifin}) if and only if the unique functor $t_\C : \C \to \boldsymbol{*}$ has a right adjoint; dually, $\C$ has an initial object if and only if $t_\C$ has a left adjoint.
	\end{enumtag}
\end{example}
\subsection{Unit and counit, triangle identities}\label{unit_counit}\index{Adjunction!unit and counit}\index{Unit!--- of an adjunction}\index{Counit!--- of an adjunction}
From every adjoint pair $(F,G)$ we can obtain a pair of natural transformations, the \emph{unit} and the \emph{counit} of the adjunction: if in \eqref{eqn:adj} we put $D=FC$ we have
\[
	\phi_{CFC}:\D(FC,FC)\to \C(C, G FC)
\]
Now, the domain $\D(FC,FC)$ is nonempty (by \ref{def:cate} it must contain at least the identity of $FC$) thus $\varphi$ is not the empty function, and it is well defined the image $\eta_C=\phi_{CFC}(\id_{FC})$ of the identity under $\varphi$; $\eta_C$ is the component at $C$ of a natural transformation $\eta:\id_\C\to  G \circ F $ called the \emph{unit} of the adjunction. Indeed, given $h:C'\to C$ we have that
\begin{align}
	\label{eqn:aggiu1}
	GF (h)\circ \eta_{C'} & = GF (h)\circ \phi_{C'FC'}(\id_{FC'})\notag  \\
	                      & =\phi\circ F (h)\circ \id_{FC'} \notag       \\
	                      & =\phi_{C'FC}\circ \id_{FC}\circ F (h) \notag \\
	                      & = \phi_{CFC}(\id_{FC})\circ h \notag         \\
	                      & = \eta_C\circ h.
\end{align}
This is equivalent to the joint commutativity of both parts in the diagram
\[\vcenter{
	\xymatrix{
	\D(FC',FC)\ar[rr]^{\D(FC',F (h))}
	\ar[d]_{\phi_{C'FC'}}
	&& \D(FC',FC)\ar[d]_{\phi_{C'FC}} && \D(FC,FC)
	\ar[d]_{\phi_{CFC}}\ar[ll]_{\D(F (h),FC)}\\
	\C(C', GFC')\ar[rr]_{\C(C', G (F (h)))} && \C(C', G (FC)) && \C(C, GFC)
	\ar[ll]^{\C(h, GFC)}
	}}
\]
Dually, if in \ref{eqn:adj} we put $C= GD$ we have
\[
	\phi_{ G DD} : \D(F GD,D)\to \C( G D, G D)
\]
(or equivalently $\phi_{ G DD}^{-1} : \C( G D, G D)\to\D(F ( G D),D)$). Again, $\C( G D, G D)$ is nonempty, thus we can define $\epsilon_D = \phi_{ G DD}^{-1}(\id_{ G D})$; this arrow is the component at $D$ of a natural transformation $\epsilon :F \circ G \to \id_\D$ called the \emph{counit} of the adjunction, in such a way that for every $k:D\to D'$ the diagram
\[
	\vcenter{\xymatrix{
			F ( G D)\ar[d]_{F ( G (k))}\ar[r]^(.6){\epsilon_D} & D \ar[d]^k\\
			F ( G (D')) \ar[r]_(.6){\epsilon_{D'}}& D'
		}}
\]
commutes.
\begin{notation}\label{adjnotaz}
	A compact way to denote all the information of an adjunction is
	\[
		F :\C\adjunct{\epsilon}{\eta}\D : G
	\]
	(to be parsed as `there is an adjunction between the functors $F$ and $G$, $F : \C\to \D$ and $G : \D \to \C$, with counit $\epsilon$ and unit $\eta$'): one of the first important results about adjunctions is that we can characterise them solely through their unit and counit instead of using the adjunction isomorphisms in \ref{def:adj}.
\end{notation}
\begin{proposition}[Zig-zag identities]
	\label{zigozago}
	\index{Adjunction!triangle identities}
	\index{Unit!--- of an adjunction}
	\index{Counit!--- of an adjunction}
	\index{Adjunction!Triangle identities|see{Adjunction!triangle identities}}
	Let $F  : \C\leftrightarrows \D :  G $ be an adjunction, having  $\eta$ and $\epsilon$ as unit and counit; using the whiskering operation of \ref{comporizz} we have
	\begin{gather}
		( G * \epsilon)\circ (\eta *  G ) = \id_G \notag\\
		(\epsilon * F )\circ (F  * \eta) = \id_F  \label{adjident}
	\end{gather}
	($\id_F$ denotes, here, the identity natural transformation of $F$ into itself, and similarly for $\id_G$).
\end{proposition}
\begin{proposition}[Adjoint preserve co\fshyp{}limits]\leavevmode\label{adj-preserve}
	Let $F\dashv G$ be an adjunction, and let $F: \C\to \D$ be the left adjoint. Then
	\begin{itemize}
		\item $F$ preserves all colimits that exist in $\C$, in the sense of \ref{preservation};
		\item dually, the right adjoint $G$ preserves all limits that exist in $\D$, in the sense of (the dual) of \ref{preservation}.
	\end{itemize}
\end{proposition}
\begin{theorem}
	Let $\C$ be a category, $F  : \C\dashv\D :  G $ a pair of adjoint functors. Then, the following conditions are equivalent:
	\begin{enumtag}{eq}
		\item $F , G $ are both fully faithful (see \ref{fullfaithconsfun}.\ref{fufai});
		\item Te unit $\eta$ and the counit $\epsilon$ of the adjunction are both natural equivalences (see \ref{def:equnat});
		\item $F $ is an equivalence of categories between $\C,\D$, whose inverse is $G$ (see \ref{def:equcat}).
	\end{enumtag}
	If one of these conditions holds true, then there is also an adjunction
	\[G : \D\adjunct{\eta^{-1}}{\epsilon^{-1}} \C : F .\]
\end{theorem}
\section{The Yoneda lemma}\label{yoneda_lemma}\index{Yoneda lemma}
\subsection{Presheaves}
\begin{definition}[Category of presheaves]\label{def:psh}
	A \emph{presheaf} is a functor $F : \C^\opp\to \Set$; the category of presheaves has morphisms the natural transformations between functors, as defined in \ref{def:trasnat}, \ie the families of maps $\{\alpha_A : FA \to GA\}$ such that the square of functions
	\[
		\vcenter{\xymatrix{
				FA \ar[r]^{\alpha_A} & GA \\
				FB \ar[u]^{Ff}\ar[r]_{\alpha_B}& GB\ar[u]_{Gf}
			}}
	\]
	commutes for every arrow $f : A\to B$ in $\C$.
\end{definition}

\begin{definition}[Representable presheaf]\label{def:reprepsh}
	Every object $X\in \C$ defines a presheaf obtained sending $A\in\C$ into the set of morphisms $u : A\to X$ in $\C$. This defines a functor, called the presheaf \emph{associated} to $X$, thanks to the associative property of composition in $\C$: the action of $\C(\firstblank,X)$ on morphisms is given by the function that sends $f : X\to Y$ in the natural transformation $\C(\firstblank,X)\To \C(\firstblank,Y)$ having components
	\[f_{*,A} : \C(A,X) \to \C(A,Y) : u \mapsto f\circ u\]
	A completely analogous definition can be given for a functor $X\mapsto \C(X,\firstblank)$, only in this case the functor is contravariant. We call a presheaf $F$ \emph{representable} if it is isomorphic to $\C(\firstblank,X)$ or $\C(X,\firstblank)$ (depending on its variance), for some object $X\in\C$.
\end{definition}
\subsection{The Yoneda lemma}
Yoneda lemma is one of the tautologies on which our understanding of reality is built.
\index{Yoneda lemma|(}
\index{_aaa_yoncontra@$\yon$}
\index{_aaa_yoncov@$\coyon$}
\begin{lemma}[Yoneda lemma]\label{lem:the-real-yoda}
	Let $\yon = \yon_{\C} : \C \to \Cat(\C^\opp,\Sets)$ be the functor sending $X\in\C$ to the presheaf associated to $X$; then, for every $F\in\Cat(\C^\opp,\Sets)$, there exists a bijection between the set of natural transformations $\yon X\To F$ and the set $FX$. This bijection is moreover natural in the object $X$.
\end{lemma}
\begin{proof}
	The proof is elementary, in the algebraic sense of the word: we do the only possible thing with the data we have, and it works.

	We define a function of sets
	\[
		Y : [\C^\opp, \Sets](\yon X,F) \to FX \label{yoneda-compo}
	\]
	and we show that it is a bijection.
	Given a natural transformation $\alpha : \yon X \To F$ we can consider its $X$\hyp{}component $\alpha_X : \C(X,X)\to FX$; as such, it is a function of sets, and it is not equal to the empty function (because $\C(X,X)$ contains at least one element called $\id_X$, thanks to axiom \ref{c:tre} of \ref{def:cate}).

	Thus we can define $Y(\alpha) := \alpha_X(\id_X)\in FX$. We shall show that
	\begin{itemize}
		\item $Y$ is injective. Assume that $\alpha_X(\id_X) = \beta_X(\id_X)$ for a pair of natural transformations $\alpha,\beta : \yon X\To F$. The naturality request then entails that for every $u : A \to X$ one has
		      \begin{align*}
			      \alpha_A(u) & = \alpha_A(\id_X\circ u)  \\
			                  & =Fu \circ \alpha_X(\id_X) \\
			                  & =Fu \circ \beta_X(\id_X)  \\
			                  & =\beta_A(u)
		      \end{align*}
		      (if there are no maps $u : A \to X$, then $\alpha_A,\beta_A$ coincide vacuously, as they both are the unique empty function $\varnothing \to FX$).
		\item $Y$ is surjective. Given an element $x \in FX$ we can define, for every $A\in \C$ and every $u : A \to X$ the function $\alpha_A^x : \C(A,X) \to FA : u\mapsto Fu(x)$; this is well\hyp{}defined, because $Fu$ is a function $FX\to FA$, given that $F$ is contravariant. Now, we shall show that this is the $A$\hyp{}component of a natural transformation $\yon X\To F$ (once this has been proved, the fact that $Y(\alpha^x)=x$ is obvious in view of the definition of $\alpha^x$).

		      The naturality of this correspondence follows from the fact that
		      \[Ff \circ \alpha_B^x(v) = \alpha_A^x (v\circ f)\]
		      for every $f :A\to B$, \ie from the fact that $Ff(Fv(x)) = F(vf)(x)$, true because $F$ is a contravariant functor.
	\end{itemize}
	This concludes the proof.
\end{proof}
\index{Yoneda lemma|)}
There are many ways to read this result.
\begin{itemize}
	\item A natural transformation $\yon X\to F$ is uniquely determined by the value that its $X$\hyp{}component $\alpha_X : \C(X,X)\to FX$ takes on $\id_X$.
	\item The identity arrow $\id_X$ can be thought as the \emph{universal element} witnessing the representability of a functor; thanks to the Yoneda lemma every element $t \in FX$ induces a unique natural transformation $\hat t_X : \yon X\To F$, and we call $t$ a \emph{universal element} if $\hat t_X$ is invertible. The Yoneda lemma can thus be thought as the statement that $\id_X\in \C(X,X)$ is  a universal element.
\end{itemize}
This leads to an immediate corollary:
\begin{corollary}\label{yon_is_ff}
	The Yoneda map is a fully faithful functor $\yon : \C \to \Cat(\C^\opp,\Sets)$.
\end{corollary}
As a consequence, we call $\yon$ the \emph{Yoneda embedding} \index{Yoneda!--- embedding} of $\C$ in $\Cat(\C^\opp,\Sets)$.
\begin{proof}
	We can apply \ref{lem:the-real-yoda} to the special case of a representable $F$: if $F = \hom(\firstblank,A)$ for an object $A\in\C$, we have
	\[
		\Cat(\C^\opp,\Sets)(\yon X, \yon A)\cong \yon AX = \C(X,A)
	\]
	According to the way in which the bijection $Y$ is defined above, in this case it coincides with the action of the functor $\yon$ on arrows, and this allows to conclude. The proof that $Y_X$ defines the $X$\hyp{}component of a natural transformation $\yon\To F$ is an easy exercise for the reader stated in \ref{exe:yon-is-nat}.
\end{proof}
\index{_aaa_yoncontra@$\yon$}
\index{_aaa_yoncov@$\coyon$}
\begin{remark}
	There is of course an analogue of the Yoneda lemma for \emph{covariant} functors $G : \C \to \Set$; given such $G$, there is a natural bijection
	\[ \Cat(\C, \Set)(\coyon_\C A, G)\cong GA \]
	given by evaluating at the universal element $\id_A$ for the representable $\coyon_\C A = \C(A,\firstblank)$; a similar argument as the one in \ref{yon_is_ff} yields that $\coyon_\C : \C^\opp\to \Cat(\C, \Set)$ is fully faithful. The same happens for the next result, that can be seen as another capital result of basic category theory: the essential image of $\yon_\C$ generates all the category $[\C^\opp, \Sets]$ under colimits.
\end{remark}
\begin{lemma}\label{its-a-cone}
	The natural transformation $\hat t_X$ is the $X$\hyp{}component of a cocone for the diagram
	\[
		\textstyle \elts{F}{\C} \xto{\Sigma} \C \xto{\yon} \Cat(\C^\opp,\Sets)
	\]
	This means that for every $f \in \C(X,Y)$ which is a morphism in $\elts{F}{\C}$ between objects $(X, x\in FX)$ and $(Y,y\in FY)$, \ie such that $Ff(y)=x$, the diagram
	\[
		\vcenter{\xymatrix{
				& F & \\
				\yon X\ar[rr]_{f_*}\ar[ur]^{\hat x} && \yon Y \ar[ul]_{\hat y}
			}}
	\]
	commutes.
\end{lemma}
\begin{theorem}[Density of the Yoneda embedding]\label{thm:yoda-is-dense}
	\index{Yoneda!density of the --- embedding}
	Let $ F: \C^\opp\to \Sets$ be a presheaf; then $F$ is canonically isomorphic to the colimit of the diagram
	\[\textstyle
		\elts{\C}{F} \xto{\Sigma_F} \C \xto{\yon} \Cat(\C^\opp,\Sets)
	\]
	or in other words, the presheaf $F$ is isomorphic to $\varinjlim_{(X,x)\in \elts{\C}{F}}\yon X$.
\end{theorem}

\begin{remark}
	The density theorem can be rephrased using coend calculus, as we know from \ref{section:weight}:  every presheaf is isomorphic to the weighted colimit of the Yoneda embedding, and it is its own weight. Alternatively, the ninja Yoneda lemma of chapter \ref{sec:tre} constitutes an elegant rewriting of the density theorem: the isomorphism given by the coend
	\[
		\int^A FA\times\C(\firstblank,A)\cong F(\firstblank)
	\]
	describes explicitly the way in which the presheaf $F$ is a colimit of all representable presheaves $\C(\firstblank,A)$.
\end{remark}
\subsection{Alternative looks on the Yoneda lemma}\label{alt_yoneda}
\index{Grothendieck!construction}
\begin{definition}\label{eltsf}\index{Category!--- of elements}\index{Category of elements}
	Let $\C$ be an ordinary category, and let $W : \C\to \Sets$ be a functor; the \emph{category of elements} $\elts{\C}{W}$ of $W$ is the category having objects the pairs $(C\in\C, u\in WC)$, and morphisms $(C,u)\to (C',v)$ those $f\in \C(C,C')$ such that $W(f)(u)=v$.
\end{definition}
\begin{notat}
	The notation `$\elts{\C}{W}$' for the category of elements of $W$ is borrowed from \cite{Graya}. Other references call it $\int W$ (it is obvious why we can't stick to this more compact notation) or $\text{El}(W)$.
\end{notat}
\begin{proposition}\label{eltsf:char}
	The category $\elts{\C}{W}$ defined in \ref{eltsf} can be equivalently characterised as each of the following objects:
	\begin{enumerate}[label=$\roman*$)]
		\item The category which results from the pullback
		      \[
			      \xymatrix{
				      \elts{\C}{W}\ar[r]\ar[d] \pb & \Sets_* \ar[d]^U \\
				      \C \ar[r]_W & \Sets
			      }
		      \]
		      where $U : \Sets_*\to\Sets$ is the forgetful functor which sends a pointed set to its underlying set;
		      \index{Comma category}\index{Category!comma ---}
		\item The comma category $(*\downarrow W)$ of the cospan $\{*\}\to \Sets \xot{W} \C$, where $\{*\}\to \Sets$ chooses the terminal object of $\Sets$;
		\item The opposite of the comma category $(\coyon_{\C}\downarrow \lceil W\rceil)$, where $\lceil W\rceil : \{*\}\to [\C, \Sets]$ is the \emph{name} of the functor $W$, \ie the unique functor choosing the presheaf $W\in [\C, \Sets]$:
		      \[
			      \vcenter{
			      \xymatrix{
			      (\elts{\C}{W})^\opp \ar[d] \ar[r] \drtwocell<\omit>{} & {*} \ar[d]^{\lceil W\rceil}\\
			      \C^\opp \ar[r]_{\coyon_{\C}} & \Cat(\C, \Set)
			      }
			      }
		      \]
	\end{enumerate}
\end{proposition}
\begin{proof}
	The proof that these categories are all canonically isomorphic to $\elts{\C}{W}$ is an exercise in Yoneda lemma and universal properties that we leave to the reader.
\end{proof}
\begin{definition}[Discrete opfibration]
	\label{def:dfib}
	\index{Discrete fibration}
	A \emph{discrete opfibration} of categories is a functor $G : \E \to \C$ with the property that for every object $E\in\E$ and every arrow $p : C\to GE$ in $\C$ there is a unique $q : E'\to E$ `over $p$', \ie such that $Gq=p$.
\end{definition}
With a straightforward definition of morphism between discrete opfibrations, we can define the category $\DFib(\C)$ of discrete opfibrations over $\C$. The nomenclature here comes from the fact that there exists a dual notion of discrete \emph{fibration}: it is a functor $G : \E \to \C$ such that for every $E\in \E$ and $p : GE \to C$, there exists a unique $q : E \to E'$ such that $Gq=p$. Of course, there is a category of discrete opfibrations over $\C$. As we will see later on, discrete \emph{op}fibrations determine contravariant functors, while fibrations determine covariant ones out of $\C$.
\index{_aaa_DFib@$\DFib(\C)$}
\begin{proposition}\label{fibelem}\index{Category!--- of elements}\index{Category of elements}
	The category of elements $\elts{\C}{W}$ of a functor $W : \C\to \Sets$ comes equipped with a canonical discrete fibration to the domain of $W$, which we denote $\Sigma : \elts{\C}{W}\to \C$, defined forgetting the distinguished element $u\in Wc$.
\end{proposition}
Now, let $S$ be any set; it is well\hyp{}known (see for example \cite{mac1992sheaves}) that the category of functors $S \to \Sets$ (viewing the set $S$ as a discrete category) is equivalent to the category $\Sets/S$ of \emph{functions over $S$}; this seemingly innocuous result is a particular instance of what is called \emph{Grothendieck construction}. If possible, the Grothendieck construction is equally important than the Yoneda lemma, as it clarifies the way in which categories are intrinsically geometric entities; a presheaf $F :\C^\opp\to \Sets$ is equivalent to some sort of generalised space `lying over' the category $\C$, in a similar way the total space of a fiber bundle lies over its base.
\index{Category!--- of elements}\index{Category of elements}
\index{_aaa_CintF@$\elts{\C}{F}$}
Now, we can consider the \emph{category of elements} \ref{eltsf} of a presheaf $F : \C^\opp\to \Sets$; this sets up a functor from $\Cat(\C^\opp,\Sets)$ to the category of discrete opfibrations over $\C$: the Grothendieck construction asserts that this is an equivalence of categories, as defined in \ref{def:equcat}.
\begin{theorem}\label{thm:equconfib}
	There is an equivalence of categories
	\[
		\Cat(\C^\opp,\Sets) \to \DFib(\C)
	\]
	defined by the correspondence sending $F\in\Cat(\C^\opp,\Sets)$ to its \emph{fibration of elements} \index{Category!--- of elements}\index{Category of elements} $\Sigma_F : \elts{\C}{F} \to \C$.
\end{theorem}
\begin{proof*}
	First of all, let's show that $\Sigma_F$ is a discrete opfibration; then, we shall show that there is a correspondence in the opposite direction assigning to a discrete opfibration $G : \E \to \C$ a presheaf of sets $p^G : \C^\opp\to \Sets$. To this end, given a functor $G$ having codomain $\C$, we define a correspondence that sends an object $X\in\C$ into the category $G^\leftarrow X$, where the notion of fiber is understood as follows: we take the objects $E\in\E$ such that $GE=X$, and morphisms $E\to E'$ such that $Gf=\id_X$.

	It is easy to notice that the fiber of $G$ over $X$ is a discrete category if $G$ is a discrete opfibration; this means that each fiber of $G$ over $X$ can be regarded as a set (incidentally, this motivates the \emph{discrete} in the name); it is equally easy to see that a morphism $f : X\to X'$ induces a function of sets between the fibers $G^\leftarrow X'$ and $G^\leftarrow X$ (see Figure \ref{fig:discfib} below); this allows us to show that the correspondence with domain $\DFib(\C)$ that sends $G$ into $p^G : X \mapsto G^\leftarrow X$; the above argument entails that $p^G$ is a presheaf of sets on $\C$ (the discrete opfibration condition entails that it is a contravariant functor).\index{_aaa_DFib@$\DFib(\C)$}

	It is a straightforward check that the composition of the two correspondences is the identity in both directions:
	\begin{itemize}
		\item Given a discrete opfibration $G :\E \to \C$, the category $\elts{\C}{p^G}$ receives a natural morphism of fibrations from $G$, in a commutative triangle
		      \[
			      \vcenter{\xymatrix{
			      \elts{\C}{p^G}\ar[dr]_{\Sigma_{p^G}} && \ar[ll]^{\eta_G} \E \ar[dl]^G \\
			      & \C &
			      }}
		      \]
		      this is the unit of an adjunction $(\firstblank)^G \dashv \elts{\C}{\firstblank}$, and it is quite easy to show that $\eta_G$ is an isomorphism over $\C$ (\ie an isomorphism of categories over $\C$).
		\item given a presheaf $F$, and given the definition of $p^{\elts{\C}{F}}$, the identity maps work as components of a natural transformation $p^{\elts{\C}{F}} \To F$.
	\end{itemize}
\end{proof*}
\begin{center}
	\begin{figure}[t]
		\begin{tikzpicture}
			\fill (0,0) circle (1pt) node[below] (X) {$X$};
			\fill (2,-.5) circle (1pt) node[below] (X') {$X'$};
			\fill[gray!20] (X) -- (-.35,4) -- (.35,4.25) node[above, black] {$G^\leftarrow X$} -- cycle;
			\fill[gray!20] (X') -- (2-0.35,4-.5) -- (2+.35,3.75) node[black, above] {$G^\leftarrow X'$} -- cycle;
			\fill (2,2) circle (1pt) node[above] (e) {$E$};
			\fill (0,2.75) circle (1pt) node[above] (e') {$E'$};
			\draw[dashed,->, >=stealth] (e) to[bend right] (e');
		\end{tikzpicture}
		\caption{A discrete opfibration $G : \E \to \C$ induces a presheaf, sending $X\in\C$ to the fiber $G^\leftarrow X$: a morphism $f :  X\to X'$ induces a function of sets $p^GX'\to p^GX$ sending $E\in p^GX'$ to the unique $E'$ such that $GE'=X$.}
		\label{fig:discfib}
	\end{figure}
\end{center}
With \ref{thm:equconfib}  in our hand, we shall now attempt to answer the very natural question: what does the Yoneda lemma become, if a presheaf $F :\C^\opp\to\Sets$ is regarded as a discrete opfibration over $\C$? First of all, the proof of the following result is an immediate consequence of the definition of $\elts{\C}{\C(\firstblank, C)}$:
\begin{lemma}\label{lem:fibofrepre}
	The fibration of elements $\Sigma_{\C(\firstblank,C)}=\Sigma_C$ of a representable presheaf coincides with the slice category $\C/C$ of arrows over $C$, defined in \ref{defn:slice}; the functor $\Sigma_C$ coincide with the `src' functor sending $u : A\to C$ into $A$.
\end{lemma}
We are now ready to investigate the form of \ref{lem:the-real-yoda} in terms of fibrations:
\begin{proposition}[Yoneda lemma, the geometric way]\label{geoyoda}
	\index{Yoneda lemma!geometric ---}
	Let $G : \E\to \C$ be a discrete opfibration over a small category $\C$; let $X\in\C$ an object; then, there is a bijection between the set of functors $H : \C/X\to \E$ such that $G\circ H = \mathrm{src}$, \ie of discrete opfibration maps from the element opfibration of $\yon X$ and $G$, and the set $G^\leftarrow X$.
\end{proposition}
\begin{proof}
	The opfibration of elements of a representable functor $\yon X$ is the category $\C/X$ of arrows having codomain $X$; the functor $\mathrm{src} : \C/X\to \C$ being the opfibration $\Sigma_C$ in the notation of \ref{lem:fibofrepre}. A morphism of discrete opfibrations is a functor $H : \C/X\to \C$ such that $GH(\id_X) = X$, so $H\mapsto H(\id_X)$ defines a function $\DFib(\C) \to G^\leftarrow X$. This assignment determines the desired bijection, as it is easy to see using a similar argument of \ref{lem:the-real-yoda}.
\end{proof}
\section{Monoidal categories and monads}\label{moncat_and_monads}
\begin{definition}[Monoid in a monoidal category]\label{def:monoid-in-cat}\index{Monoid!internal ---}
	\index{Object!monoid ---}
	Let $\C$ be a monoidal category with monoidal structure $\otimes$ and monoidal unit $I$; we define a \emph{monoid} in $\C$ to be an object $M$ endowed with maps $m : M\otimes M \to M$ and $u : I\to M$ such that the diagrams
	\[
		\vcenter{\xymatrix{
				M\otimes M\otimes M \ar[r]^{M\otimes m}\ar[d]_{m\otimes M}& M\otimes M \ar[d]^m & M\otimes I \ar@{=}[dr]\ar[r]^{M\otimes u}& M\otimes M \ar[d]_m & I\otimes M\ar[l]_{u\otimes M}\ar@{=}[dl]\\
				M\otimes M \ar[r]_m & M & & M
			}}
	\]
	commute; these testify the associativity and unitality property of $(m,u)$.
\end{definition}
\begin{remark}\index{Monad!---s as monoids}
	The category $[\C,\C]$ of endofunctors of a category $\C$ has a natural choice of a monoidal structure given by composition, whose monoidal unit is the identity functor $\id_\C$; it is a very useful exercise to verify the axioms of monoidal category one by one; all coherence morphisms are in fact identity, so this is an example of a \emph{strict} monoidal structure (the composition of functors is strictly associative, and $\id_\C$ is a strict unit). Of course, this monoidal structure is highly nonsymmetric.
\end{remark}
\begin{definition}[Monad]\label{def:monad}\index{Monad}\index{Unit!--- of a monad}
	Let $\C$ be a category; a \emph{monad} on $\C$ consists of an endofunctor $T : \C\to \C$ endowed with two natural transformations
	\begin{itemize}
		\item $\mu : T\circ T\To T$, the \emph{multiplication} of the monad, and
		\item $\eta : \id_\C \To T$, the \emph{unit} of the monad,
	\end{itemize}
	such that the following axioms are satisfied:
	\begin{itemize}
		\item the multiplication is associative, \ie the diagram
		      \[
			      \vcenter{\xymatrix{
					      T\circ T\circ T\ar[r]^{T *\mu}\ar[d]_{\mu *T} & T\circ T \ar[d]^\mu\\
					      T\circ T \ar[r]_\mu & T
				      }	}
		      \]
		      is commutative, \ie the equality of natural transformations $\mu\circ (\mu * T) = \mu \circ (T * \mu)$ holds;
		\item the multiplication has the transformation $\eta$ as unit, \ie the diagram
		      \[
			      \vcenter{\xymatrix{
					      T \ar[r]^{\eta *T}\ar@{=}[dr]& T\circ T \ar[d]_\mu & T\ar[l]_{T*\eta} \ar@{=}[dl]\\
					      & T &
				      }	}
		      \]
		      is commutative, \ie the equality of natural transformations $\mu\circ (\eta *T)=\mu\circ (T * \eta)= \id_T$ holds.
	\end{itemize}
\end{definition}
\begin{proposition}\index{Monad!---s from adjunctions}\index{Unit!--- of a monad}
	Let $F\adjunct{\epsilon}{\eta} G$ be an adjunction between two categories; say $F:\C\to \D$; then the composition $GF$ is the endofunctor of a monad on $\C$; the multiplication map is given by the whiskering $G *\epsilon * F : GFGF \To GF$, and the unit of the monad coincides with the unit of the adjunction, $\eta : \id_\C \To GF$.
\end{proposition}
The correspondence sending an adjunction $F\adjunct{\epsilon}{\eta} G$ into a monad $(T,G * \epsilon *F, \eta)$ is a functor towards a suitable category of monads; but it is predictably highly nonbijective, and in fact there is an entire \emph{category} $\sfSpli(T)$ of adjunctions inducing the same monad.
The category $\sfSpli(T)$ is in general quite difficult to describe; we know, however, that it always has a terminal and an initial object.
\index{Object!initial and terminal}
\begin{definition}\label{alg_for_a_mona}\index{Monad!algebra for a ---}\index{Category!Eilenberg-Moore ---}
	Let $T : \C\to \C$ be a monad; we define a \emph{$T$\hyp{}algebra} as a pair $(A,a)$, where $A\in\C$ and $a : TA\to A$ is a morphism (the \emph{algebra map}) in $\C$, such that the following properties hold
	\begin{itemize}
		\item compatibility with the multiplication, namely the equation $a\circ Ta = a \circ \mu_A$;
		\item compatibility with the unit, namely the equation $a\circ \eta_A=\id_A$.
	\end{itemize}
	A \emph{morphism of $T$\hyp{}algebras} is a morphism $f : A\to B$ that commutes with the algebra maps $a,b$ of the objects $A,B$, namely such that the square
	\[
		\vcenter{\xymatrix{
				TA \ar[d]_a \ar[r]^{Tf}& TB\ar[d]^b \\
				A \ar[r]_f & B
			}}
	\]
	commutes.
\end{definition}
\begin{proposition}\label{em_is_terminal}
	The category of $T$\hyp{}algebras, so defined, has the universal property of a terminal object in $\sfSpli(T)$.
\end{proposition}
\index{_aaa_@$\Alg(T)$}

\begin{definition}\label{frealg}\index{Monad!free algebra for a ---}\index{Category!Kleisli ---}
	Every pair $(TX, \mu_X)$ is, trivially, a $T$\hyp{}algebra; these objects define a (full) subcategory of \emph{free $T$\hyp{}algebras}.
\end{definition}
\begin{proposition}\label{kl_is_initial}
	The category of free $T$\hyp{}algebras has the universal property of an initial object in $\sfSpli(T)$.
\end{proposition}
\index{_aaa_Kl@$\Kl(T)$}
The universal properties \ref{em_is_terminal} and \ref{kl_is_initial} of $\Kl(T)$ and $\Alg(T)$ yields a unique functor $\Kl(T) \to \Alg(T)$, corresponding to the inclusion of free algebras into all algebras. Exercise \ref{kleisli-is-kl} makes this characterisation precise.

More in general, the universal property of $\Alg(T)$ as a terminal object yields a unique functor for every other adjunction $F\dashv G$ splitting a given monad:
\begin{definition}\label{comparimona}
	Given a pair of adjoint functors $(F\dashv G)$ in $\sfSpli(T)$ we call \emph{comparison functor} associated to the adjunction the functor $K$ fitting in the diagram
	\[
		\vcenter{\xymatrix{
				& \C \ar@<-3pt>[dr]_{F^T}\ar@<-3pt>[dl]_{F}& \\
				\D \ar@<-3pt>[ur]_{G}\ar@{.>}[rr]_K && \Alg(T)\ar@<-3pt>[ul]_{G^T}
			}}
	\]
\end{definition}
\begin{definition}\label{monadic_fu}
	An adjunction $(F\dashv G)$ is \emph{monadic} if the associated comparison functor $K$ is an equivalence.
\end{definition}
\emph{Beck's theorem} gives a necessary and sufficient condition to recognise if a functor is the right adjoint of a monadic adjunction.
\begin{definition}
	Let $U : \C \to \D$ be a functor; we say that a diagram
	\[	\xymatrix{
			X \ar@<-4pt>[r]_g \ar@<4pt>[r]^f & Y \ar[r]^h & Z
		}\]
	is a \emph{split coequaliser} if there exist a fourth and a fifth arrows $s : Z \to Y$ and $t : Y \to X$ such that $h\circ s = \id_Z$, $g\circ t = s\circ h$ and $f\circ t = 1_Y$.
\end{definition}
Note that if a diagram like the one above is a split coequaliser, then $h$ is forced to be the coequaliser of the pair $f,g$, so that $h$ is uniquely determined up to isomorphism by the pair $f,g$; note also that every functor $F : \D \to \E$ sends a split coequaliser in $\D$ into a split coequaliser in $\E$.

We say that a pair of arrows $f,g$ in $\C$ extends to a split coequaliser if there is an $h$ such that the diagram $X \underset{g}{\overset{f}\rightrightarrows} Y \xto{h} Z$ is a split coequaliser; given the above remark, this is a well-defined notion.
\begin{definition}\label{monadici}
	Let $U : \C \to \D$ be a functor; then $U$ is monadic if and only if admits a left adjoint and
	\begin{enumtag}{mf}
		\item $U$ is conservative (\ie $Uf$ is an isomorphism if and only if $f$ is);
		\item $\D$ has, and $U$ preserves, the coequalisers of \emph{$U$-split pairs}, \ie those parallel pairs of morphisms in $\C$,
		\[
			\xymatrix{
				X \ar@<-4pt>[r]_g \ar@<4pt>[r]^f & Y
			}
		\]
		such that the pair $Uf,Ug$ has a split coequaliser in $\D$.
	\end{enumtag}
\end{definition}
\begin{definition}\label{def:idemmona}\index{Monad!idempotent}
	An \emph{idempotent monad} is a monad $(T,\mu,\eta)$ such that the multiplication $\mu  : TT\To T$ is an isomorphism. If $T$ is an idempotent monad, each of these properties is true (and equivalent to the request of idempotency):
	\begin{enumtag}{im}
		\item the natural transformations $T*\eta$ and $\eta *T$ are (invertible and) equal;
		\item If $A$ is a $T$\hyp{}algebra, then there is a unique algebra map $TA\to A$, and it is invertible;
		\item The category of $T$\hyp{}algebras embeds into the category $\C$ via its forgetful functor $U : \Alg(T)\to \C$, which is thus full and faithful;
		\item Every adjunction $(F,G)\in\sfSpli(T)$ has invertible counit.
	\end{enumtag}
\end{definition}
Being endofunctors of $\C$, any two monads $T,S$ on $\C$ can be composed, but the composition $ST$ is often not a monad. In order for such functor to have a multiplication $\mu : STST\To ST$ we must provide an `intertwining' operator $\lambda : TS \To ST$.
\begin{definition}\index{Distributive law}\index{Monad!distributive law}\label{distribbio}
	A \emph{distributive law} between two monads consists of a 2\hyp{}cell $\lambda : TS \To ST$ such that the following commutativities hold:
	\begin{enumtag}{dl}
		\item $\lambda\circ (\eta^{(T)} * S) = S * \eta^{(T)}$;
		\item $(S * \mu^{(T)}) \circ (\lambda * T) \circ (T * \lambda) = \lambda \circ (\mu^{(T)} * S)$.
	\end{enumtag}
	\begin{enumtag}{dr}
		\item $\lambda\circ (T*\eta^{(S)}) = \eta^{(S)} * T$;
		\item $\lambda \circ (T * \mu^{(S)}) = (\mu^{(S)} * T)\circ (S * \lambda)\circ (\lambda * S)$.
	\end{enumtag}
	These conditions are expressed by the commutativity of diagrams
	\begin{gather}
		\vcenter{\xymatrix{
		TS\ar[rr]^\lambda && ST & TTS \ar[r]^{T\lambda}\ar[d]_{\mu^{(T)}S}& TST\ar[r]^{\lambda T} & STT \ar[d]^{S\mu^{(T)}}\\
		& S \ar[ur]_{S\eta^{(T)}} \ar[ul]^{\eta^{(T)}S}&& TS\ar[rr]_\lambda && ST
		}}\\
		\vcenter{\xymatrix{
		TS\ar[rr]^\lambda && ST & TSS \ar[r]^{\lambda S}\ar[d]_{T\mu^{(S)}}& STS\ar[r]^{S\lambda} & SST \ar[d]^{\mu^{(S)}T}\\
		& T \ar[ur]_{\eta^{(S)}T} \ar[ul]^{T\eta^{(S)}}&& TS\ar[rr]_\lambda && ST
		}}
	\end{gather}
\end{definition}
\begin{remark}\label{rmk_monad_mor}\index{Monad}\index{Monad!morphism of ---s}\index{Morphism!--- of monads}
	The notion of distributive law can be seen as just a particular instance of a \emph{monad morphism}: in particular, as a certain kind of endomorphism $(S,\lambda)$ of the monad $T$; we now provide the reader with the precise definition of such notion. Let $S,T$ be two monads, respectively $S : \C \to \C$ and $T : \D \to \D$; a \emph{monad morphism} $(X,\lambda) : S \to T$ consists of a pair $X : \C \to \D$ and $\lambda : TX \To XS$ such that the following diagrams are commutative:

	\[
		\vcenter{
		\xymatrix{
		&\ar[dr]^{X\eta^{(T)}}\ar[dl]_{\eta^{(S)} X} X && SSX \ar[r]^{S\lambda} \ar[d]_{\mu^{(S)} X} & SXT \ar[r]^{\lambda T} & XTT \ar[d]^{X\mu^{(T)}}\\
		SX \ar[rr]_\lambda && XT & SX \ar[rr]_\lambda && XT
		}
		}
	\]
	Moreover, a 2-cell between two parallel monad morphisms $(X,\lambda), (Y, \sigma) : (S,\C) \to (T,\D)$ consists of a natural transformation $\nu : X \To Y$ such that the square
	\[\vcenter{\xymatrix{
				SX \ar[d]_\lambda \ar[r]^{S\nu} & SY\ar[d]^\sigma \\
				XT \ar[r]_{\nu T}& YT
			}}\]
	is commutative.

	Now, let $\C$ be a 2-category. A distributive law between two monads can be characterised as a monad in the 2-category of monads on $\C$, monad morphisms and monad 2-cells (the reader is invited to make this statement precise as an exercise: how do monad morphisms compose? What is the object on which a distributive law is a monad on?).
\end{remark}
\index{Monad}\index{Restrained morphism of monads|see{Monad!morphism of ---s}}
\begin{remark}\label{restrained_mormonad}
	Some authors prefer to restrict to a more rigid definition of monad morphism, where the endofunctor $X$ of \ref{rmk_monad_mor} above is the identity: we call such special monad morphisms \emph{restrained}.
\end{remark}
\index{Monad!co---}
\begin{remark}\label{comona}
	There is a dual theory of \emph{comonads}, these are comonoids in the monoidal category of endofunctors $[\C,\C]$; outlining the basic definitions is left as an exercise in \ref{ex:commies}.
\end{remark}
\section{2\hyp{}categories}\label{2categories}\index{Category!enriched ---}
\begin{definition}\label{enrichcat}
	A \emph{category enriched over the monoidal base $\V$}, or briefly a \emph{$\V$\hyp{}category} $\underline{\A}$, consists of
	\begin{enumtag}{mc}
		\item a class of objects $\underline{\A}_o$;
		\item an object $\underline{\A}[A,B]\in\V$ for each pair of objects $A,B\in\A_0$;
		\item a family of \emph{composition maps} $c_{ABC} : \underline{\A}[A,B]\otimes \underline{\A}[B,C]\to \underline{\A}[A,C]$, one for each triple of objects $A,B,C\in\underline{\A}_o$;
		\item a family of \emph{identity arrows} $i_A : I \to \underline{\A}[A,A]$, one for each object $A\in \underline{\A}_o$.
	\end{enumtag}
	These data satisfy the following axioms:
	\begin{enumtag}{a}
		\item Composition is associative, where associativity is defined via the associator of $\V$: the diagram
		\[{\scriptstyle
			\vcenter{\xymatrix@C=5mm{
			\big(\underline{\A}(C,D)\otimes \underline{\A}(B,C)\big)\otimes \underline{\A}(A,B)\ar[d]_{c_{BCD}\otimes 1}\ar[rr]^a && \underline{\A}(C,D)\otimes \big(\underline{\A}(B,C)\otimes \underline{\A}(A,B)\big) \ar[d]^{1\otimes c_{BCD}}\\
			\underline{\A}(B,D)\otimes \underline{\A}(A,B)\ar[dr]_{c_{ABD}} && \underline{\A}(C,D)\otimes\underline{\A}(A,C)\ar[dl]^{c_{ACD}}\\
			& \underline{\A}(A,D) &
			} }}
		\]
		is commutative.
		\item Composition has the identities $i_A$ as neutral elements, \ie the two diagrams
		\[\vcenter{\xymatrix@C=1.4cm{
			\underline{\A}(B,B)\otimes\underline{\A}(A,B) \ar[r]^-{c_{ABB}}& \underline{\A}(A,B) &\ar[l]_-{c_{AAB}} \underline{\A}(A,B)\otimes\underline{\A}(A,A)\\
			I\otimes\underline{\A}(A,B)\ar[u]^{i_B\otimes 1} \ar[ur]&&\underline{\A}(A,B)\otimes I\ar[u]_{1\otimes i_A} \ar[ul]}}\]
		commute.
	\end{enumtag}
\end{definition}
\begin{definition}[$\V$\hyp{}functor and $\V$\hyp{}natural transformation]\label{enrifun}\index{Functor!enriched ---}
	If $\underline{\A}$ and $\underline{\B}$ are two $\V$\hyp{}categories, a \emph{$\V$\hyp{}functor} $ F : \underline{\A} \to \underline{\B}$ consists of
	\begin{itemize}
		\item A function $F_o : \A_o \to \B_o$;
		\item A family of $\V$\hyp{}morphisms $F_{AA'} : \A(A,A') \to \B(FA,FA')$, one for each pair $A,A'\in\A_o$;
	\end{itemize}
	These data are such that the following diagrams commute:
	\begin{gather}
		\xymatrix@C=2cm{
		\A(A,A')\otimes \A(A',A'') \ar[d]_c\ar[r]^-{F_{AA'}\otimes F_{A'A''}} & \B(FA, FA')\otimes \B(FA', FA'') \ar[d]^c \\
		\A(A, A'') \ar[r]_{F_{AA''}} & \B(FA, FA'')
		}\notag\\
		\vcenter{\xymatrix{
				\B(FA,FA) &\ar[l] \ar[d]I \\
				&\A(A,A)\ar[ul]
			}}
	\end{gather}
	A $\V$\hyp{}natural transformation $\alpha$ between two $\V$\hyp{}functor $F,G : \underline{\A}\to \underline{\B}$ consists of a family of maps $\alpha_A : I \to \B(FA,GA)$ such that the following diagrams commute:
	\[\vcenter{\xymatrix@d{
		\A(A, A') \ar[r]^\rho_\wr \ar[d]_\lambda^\sim & \A(A, A') \otimes I \ar[r]^{G_{AA'}\otimes \alpha_A} & \B(G A, G A') \otimes \B(F A, G A) \ar[d]^c \\
		I \otimes \A(A, A') \ar[r]_-{\alpha_{A'}\otimes F_{AA'}} & \B(F A', G A') \otimes \B(F A, F A')\ar[r]_-c & \B(F A, G A')
		}}\]
	These diagrams express the fact that $\alpha$ is `natural' in the sense that $\alpha_A \circ Ff = Gf\circ \alpha_{A'}$ for every $f :A\to A'$; of course, there is no such thing as $f : A \to A'$ here, because $\A(A,A')$ does not have `elements' \emph{strictu senso}.
\end{definition}
\begin{definition}\index{Category!2-category}\label{duecatte}
	A \emph{$2$\hyp{}category} is a $\Cat$\hyp{}enriched category, with the tensor product given by the product of categories and the terminal category as unit. Similarly, a 2\hyp{}functor is a $\Cat$\hyp{}functor, and a 2\hyp{}natural transformation is a $\Cat$\hyp{}natural transformation. We explicitly spell out the definition of a 2\hyp{}category, leaving the definition of functor to the reader, once they will have understood how to (easily) argue by analogy.
\end{definition}
Spelling out the definition above, a $2$\hyp{}category $\C$ consist of
\begin{enumtag}{ma}
	\item a class of objects $ \C_o$,
	\item for any pair $X,Y\in C$ a small category $\C(X,Y)$,
	\item for any triple $X,Y,Z\in C$ a functor $\mu\colon\C(X,Y)\times\C(Y,Z)\to\C(X,Z)$ called composition law,
	\item a unit functor $\I\to\C(X,X)$, that is to say an object $\text{id}_X\in\C(X,X)$ for every object $X\in\C$.
\end{enumtag}
Furthermore, these data are subject to the following conditions
\begin{enumtag}{ma}
	\item for every $X\in\C$, the functors
	\begin{gather*}
		\mu(\firstblank,\text{id}_Y)\colon\C(X,Y)\to\C(X,Y)\\
		\mu(\text{id}_X,\firstblank)\colon\C(X,Y)\to\C(X,Y)
	\end{gather*}
	are the identity functors,
	\item for every $X,Y,Z,W\in\C$, the diagram
	\[
		\vcenter{\xymatrix{
		\C(X,Y)\times\C(Y,Z)\times\C(Z,W) \ar[r]^-{\mu\times\id}\ar[d]_{\id\times\mu} & \C(X,Z)\times\C(Z,W)\ar[d]^{\mu}\\
		\C(X,Y)\times\C(Y,W)\ar[r]_{\mu} & \C(X,W)
		}}
	\]
	commutes.
\end{enumtag}
The objects of $\C$ are called $0$\hyp{}cells, the objects of $\C(X,Y)$ are called $1$\hyp{}cells and the morphisms of $\C(X,Y)$ are called $2$\hyp{}cells. The notations are the same as for categories, functors and natural transformations in $\Cat$, which is the prototypical example of $2$\hyp{}category.
\begin{remark}
	The vertical composition of $2$\hyp{}cells is defined by means of the composition law of the hom\hyp{}categories, while the functor $\mu$ recovers the composition of $1$\hyp{}cells and the horizontal composition of $2$\hyp{}cells. Moreover, the functoriality of $\mu$ can be used to prove the interchange law for $2$\hyp{}cells.
\end{remark}
\begin{definition}
	A $1$\hyp{}cell $f\colon X\to Y$ inside a $2$\hyp{}category $\C$ is said to be an \emph{equivalence} if there exists another $1$\hyp{}cell $g\colon Y \to X$ together with two invertible $2$\hyp{}cells $1_X \To gf$ and $fg \To 1_Y$.
\end{definition}
Naturally, the next step is to describe how the notion of enriched functor specialises to this case.
\begin{remark}\label{pseudocolax} The previous is also known as \emph{strict} $2$\hyp{}functor, to distinguish it from other weak versions of $2$\hyp{}functors between $2$\hyp{}categories. Indeed, a 2\hyp{}functor sends identities to identities and respect compositions, just as an ordinary functor does, with an extra action on $2$\hyp{}cells. Nevertheless, it makes sense to ask for a weak version of the coherences, whose diagrams commute only up to a $2$\hyp{}cell. If these $2$\hyp{}cells are invertible we get a \emph{pseudofunctor}, otherwise we have a \emph{lax} or \emph{colax functor} (depending on the direction of the $2$\hyp{}cell). Further details can be found in \cite[\S 7.5]{Bor1}, and in \ref{colaxe} below. We say that a co\fshyp{}lax functor is \emph{normal} if it preserves identities strictly.
\end{remark}
\begin{example}
	The simplest example of $2$\hyp{}functors is the $2$\hyp{}dimensional analogue of an hom\hyp{}functor. For instance, the correspondence $\C(X,\firstblank)\colon\C\to\Cat$, for a fixed $X\in \C_o$ defines a 2\hyp{}functor. The action of this functor is really simple:
	\begin{itemize}
		\item[(i)] it sends every $0$\hyp{}cell $Y$ to the small category $\C(X,Y)$,
		\item[(ii)] it maps every $1$\hyp{}cell $f\colon Y\to Z$ to the functor
		      \begin{align*}
			      f\circ\firstblank\colon\C(X,Y) & \to\C(X,Z)         \\
			      g                              & \mapsto f\circ g   \\
			      (\gamma\colon g\To g')         & \mapsto 1_f*\gamma
		      \end{align*}
		      where $*$ is the horizontal composition of $2$\hyp{}cells.
		\item[(iii)] and finally it sends every $2$\hyp{}cell $\alpha\colon f \To g$ to the horizontal post\hyp{}composition $\alpha*\firstblank$.
	\end{itemize}
	The contravariant case $\C(\firstblank,Y)\colon\C^\opp\to\Cat$ is completely analogous.
\end{example}
\begin{definition}
	Let $F,G\colon\C\to\D$ be $2$\hyp{}functors between $2$\hyp{}categories. A $2$\hyp{}natural transformation $\alpha\colon F\To G$ is the datum of a 1\hyp{}cell \[
		\alpha_C\colon FC\to GC
	\] for every $C\in\C$, in such a way that the following diagram
	\[ \vcenter{\xymatrix{
		\C(C,C') \ar[d]_{G_{CC'}}\ar[r]^{F_{CC'}} & \D(FC,FC') \ar[d]^{\alpha_{C'}\circ\firstblank}
		\\
		\D(GC,GC')\ar[r]_{\firstblank\circ\alpha_{C}}
		& \D(FC,GC')
		}}\]
	commutes.
\end{definition}
\noindent In Remark \ref{pseudocolax} we said that $2$\hyp{}functors have weaker counterparts, namely pseudofunctors and lax functors. So it happens for $2$\hyp{}natural transformations, which in turn can be weakened into \emph{pseudonatural transformations} and \emph{lax natural transformations}. For the sake of simplicity, we prefer to give the definitions in the strict case. An explicit definition can be found, again, in \cite[\S 7.5]{Bor1}.
\begin{definition}
	Let $F, G\colon\C\rightarrow\D$ be $2$\hyp{}functors and $\alpha,\beta\colon F\To G$ 2\hyp{}natural transformations.\index{Modification} A \emph{modification} $\Xi\colon\alpha\Rrightarrow\beta$ is a family of $2$\hyp{}cells \[
		\Xi_C\colon\alpha_C\To\beta_C
	\]
	such that for any two $1$\hyp{}cells $f,g\colon C\rightrightarrows C'$ and any 2\hyp{}cell $\gamma\colon f\To g$,  we have that
	\[
		\Xi_{C'}*F\gamma=G\gamma*\Xi_C
	\]
\end{definition}
The definition above applies, basically unchanged, also to pseudonatural and lax natural transformations. By its very definition, a modification is a kind of `morphism in dimension 3'.
\begin{theorem}\label{2-yoneda}
	Let $\C$ be a $2$\hyp{}category, $F\colon\C^\opp\to\Cat$ a $2$\hyp{}functor and $C$ a $0$\hyp{}cell in $\C$, then there exist an isomorphism of categories
	\[
		[\C^\opp,\Cat](\C(C,\firstblank),F)\cong FC
	\]
	where $[\C^\opp,\Cat](\C(C,\firstblank),F)$ is the category whose objects are the $2$\hyp{}natural transformations $\C(C,\firstblank)\To F$ and with the modifications between those $2$\hyp{}natural transformations as morphisms.
\end{theorem}
\index{Bicategory}
\begin{definition}[Bicategory]\label{bicat}\index{Category!bicategory}\index{Natural transformation!composition of ---s}
	A \emph{(locally small) bicategory} $\B$ consists of the following data:
	\begin{enumtag}{bc}
		\item \label{bica:uno} A class $\B_o$ of \emph{objects}, denoted with Latin letters like $A,B,\dots$;
		\item \label{bica:due} A collection of (small) categories $\B(A,B)$, one for each $A,B\in \B_o$, whose objects are called \emph{1\hyp{}cells} or \emph{arrows} with \emph{domain} $A$ and \emph{codomain} $B$, and whose morphisms $\alpha : f \To g$ are called \emph{2\hyp{}cells} or \emph{transformations} with domain $f$ and codomain $g$; the composition law $\circ$ in $\B(A,B)$ is called \emph{vertical composition} of 2\hyp{}cells;
		\item A \emph{horizontal composition} of 2\hyp{}cells
		\[
			\boxminus_{\B,ABC} : \B(B,C)\times\B(A,B) \to \B(A,C) : (g,f)\mapsto g\boxminus f
		\]
		defined for any triple of objects $A,B,C$. This is a family of functors between hom\hyp{}categories;
		\item \label{bica:tre} for every object $A\in  \B_o$ there is an arrow $\id_A\in \B(A,A)$ such that for every $A,B\in \C_o$ and $f:A\to B$ we have $f\boxminus \id_A=f=\id_B\boxminus f$.
	\end{enumtag}
	To this basic structure we add\index{Associator!monoidal ---}
	\begin{enumtag}{bs}
		\item a family of invertible maps $\alpha_{fgh} : (f \boxminus g) \boxminus h \cong f \boxminus (g \boxminus h)$ natural in all its arguments $f,g,h$, that taken together form the \emph{associator} isomorphisms;
		\item a family of invertible maps $\lambda_f  : \id_B \boxminus f \cong f$ and $\varrho_f : f \boxminus \id_A \cong f$ natural in its component $f : A \to B$, that taken together form the \emph{left unitor} and \emph{right unitor} isomorphisms.
	\end{enumtag}
	And these data are subject to the following axioms:
	\index{Horizontal composition}
	\index{_aaa_boxminus@$\boxminus$}
	\begin{enumtag}{ba}
		\item For every quadruple of 1\hyp{}cells $f,g,h,k$ we have that the diagram
		\[
			\vcenter{\xymatrix{
			((f\boxminus g)\boxminus h)\boxminus k \ar[d]_{\alpha_{f,g,h}\boxminus k}\ar[r]^{\alpha_{fg,h,k}} & (f\boxminus g)\boxminus (h\boxminus k) \ar[r]^{\alpha_{f,g,hk}} & f\boxminus (g\boxminus (h\boxminus k))\\
			(f\boxminus (g\boxminus h))\boxminus k \ar[rr]_{\alpha_{f,gh,k}} && f\boxminus ((g\boxminus h)\boxminus k)\ar[u]_{f\boxminus \alpha_{g,h,k}}
			}}
		\]
		commutes.
		\item For every pair of composable 1\hyp{}cells $f,g$,
		\[
			\vcenter{\xymatrix{
			(f \boxminus \id_A)\boxminus g\ar[dr]_{\varrho_f\boxminus g}\ar[rr]^{a_{A,\id_A,g}} && f\boxminus(\id_A\boxminus \, g)\ar[dl]^{f\boxminus \lambda_g}\\
			& f\boxminus g
			}}
		\]
		commutes.
	\end{enumtag}
\end{definition}
\index{Bicategory}
\begin{definition}[Pseudofunctor, co\fshyp{}lax functor]\label{colaxe}\index{Functor!pseudo---}
	Let $\B,\C$ be two bicategories; a \emph{pseudofunctor} consists of
	\begin{enumtag}{pf}
		\item A function $F_o : \B_o \to \C_o$;
		\item A family of functors $F_{AB} : \B(A,B) \to \C(FA, FB)$;
		\item An invertible 2\hyp{}cell $\mu_{fg} : Ff \circ Fg \To F(fg)$ for each $A \xto{g}B\xto{f} C$, natural in $f$ (with respect to vertical composition) and an invertible 2\hyp{}cell $\eta : \eta_f : \id_{FA} \To F(\id_A)$, also natural in $f$.
	\end{enumtag}
	These data are subject to the following commutativity conditions for every 1\hyp{}cell $A \to B$:
	\[\label{psedofu}
		{\scriptstyle \vcenter{\xymatrix{
		Ff\circ \id_A \ar[r]^{\varrho_{Ff}}\ar[d]_{Ff * \eta} & Ff\ar[d]^{F(\varrho_f)} & \id_B \circ Ff\ar[d]_{\eta * Ff} \ar[r]^{\lambda_{Ff}}& Ff\ar[d]^{F(\lambda_f)}\\
		Ff \circ F(\id_A)\ar[r]_{\mu_{f,\id_A}} & F(f \circ \id_A) & F(\id_B)\circ Ff\ar[r]_{\mu_{\id_B,f}} & F(\id_B \circ f)\\
		(Ff\circ Fg) \circ Fh \ar[rrr]^{\alpha_{Ff,Fg,Fh}}\ar[d]_{\mu_{fg} * Fh} &&& Ff\circ (Fg\circ Fh)\ar[d]^{Ff * \mu_{gh}}\\
		F(fg)\circ Fh \ar[d]_{\mu_{fg} * Fh} &&& Ff\circ F(gh)\ar[d]^{\mu_{f,gh}}\\
		F((fg)h) \ar[rrr]_{F \alpha_{fgh}}&&& F(f(gh))
		}}}
	\]
	(we denote invariably $\alpha,\lambda,\varrho$ the associator and unitor of $\B,\C$).

	A \emph{lax} functor is defined by the same data, but both the 2\hyp{}cells $\mu : Ff \circ Fg \To F(fg)$ and $\eta : \id_{FA} \To F(\id_A)$ can be noninvertible; the same coherence diagrams \eqref{psedofu} hold. A \emph{colax} functor reverses the direction of the cells $\mu,\eta$, and the commutativity of \eqref{psedofu} changes accordingly.\index{Functor!lax ---}
\end{definition}

\section{Higher categories}\label{higher_caz}
\begin{definition}[The simplex category]\label{deltacat}\index{Category!simplex ---}\index{_aaa_delta@$\bDelta$}
	The \emph{simplex category}, denoted $\bDelta$, is defined as the category having
	\begin{itemize}
		\item objects the nonempty finite sets that are totally ordered: the typical object of $\bDelta$ is denoted $[n]=\{0<\cdots<n\}$, in this way, $[0]$ is the terminal object of $\bDelta$;
		\item morphisms $[n]\to [m]$ are the order preserving functions.
	\end{itemize}
\end{definition}
\begin{remark}\label{faces_and_deg}
	If in $\bDelta$ we consider
	\begin{itemize}
		\item the $n+1$ injective functions $\delta_{n,k}\colon [n]\to [n+1]$ whose image misses $k$ (called \emph{co-faces});
		\item the $n+1$ surjective functions $\sigma_{n,k}\colon [n]\to [n-1]$ assuming the value $k$ twice (called \emph{co\hyp{}degeneracies});
	\end{itemize}
	we obtain that every $f\colon [m]\to [n]$ can be written as a composition
	\[
		f = \delta_{n_1, k_1}\circ \dots\circ \delta_{n_r, k_r}\circ \sigma_{m_1, h_1}\circ\dots \circ \sigma_{m_s, h_s}
	\]
	for some indices $n_i, m_j$ and $0\le k_i\le n_i, 0\le h_j\le m_j$.
\end{remark}
\begin{remark}
	The functions $\delta_{n,k}$ and $\sigma_{n,k}$ satisfy the \index{Simplicial!--- identities} \emph{cosimplicial identities}, namely they fit into commutative diagrams
	\begin{gather}\label{cosi}
		\xymatrix{
		[n-1] \ar@{}[dr]|{i<j}\ar[r]^{\delta_i}\ar[d]_{\delta_{j-1}} & [n] \ar[d]_{\delta_j}\\
		[n] \ar[r]_{\delta_i} & [n+1]
		}
		\xymatrix{
		[n-1] \ar@{}[dr]|{i<j}\ar[r]^{\delta_i}\ar[d]_{\sigma_j} & [n]\ar[d]^{\sigma^j}\\
		[n-2] \ar[r]_{\delta_i}& [n-1]
		}
		\xymatrix{
		[n-1] \ar[r]^{\delta_j}\ar[d]_{\delta_{j+1}}& [n]\ar[d]^{\sigma_j}\\
		[n] \ar[r]_{\sigma_j}& [n-1]
		}\notag\\
		\xymatrix{
		[n-1] \ar@{}[dr]|{i>j+1}\ar[r]^{\delta_i}\ar[d]_{\sigma_j}& [n]\ar[d]^{\sigma_j}\\
		[n-2] \ar[r]_{\delta_{i-1}} & [n-1]
		}
		\xymatrix{
		[n+1] \ar@{}[dr]|{i\le j}\ar[r]^{\sigma_i}\ar[d]_{\sigma_i}& [n]\ar[d]^{\sigma_j}\\
		[n] \ar[r]_{\sigma_{j+1}}& [n-1]
		}
	\end{gather}
\end{remark}
Along the discussion, we consider $\bDelta$ as a full subcategory of $\Cat$, via the identification
\[
	\{0<\dots<n\} = \{0\to 1\to\dots\to n\}.
\]
This secretly defines an identity\hyp{}on\hyp{}objects functor $\iota : \bDelta \to \Cat$.
\begin{definition}[The category of simplicial sets]\index{Simplicial!--- set}\index{Category!simplex ---}\index{_aaa_delta@$\bDelta$}
	The category of \emph{simplicial sets} is defined as the category $[\bDelta^\opp,\Sets]$ of presheaves on $\bDelta$.
\end{definition}
A simplicial set can equivalently be specified by a collection of sets (a \emph{graded set}) $\{X_n\mid n\ge 0\}$ with maps
\[
	d_{n,i}\colon X_n \to X_{n-1} \qquad s_{n,j}\colon X_n \to X_{n+1}, \qquad 0\le i,j\le n
\]
(respectively called the \emph{faces} and \emph{degeneracies} of $X$) satisfying the dual identities of \eqref{cosi}, called the \emph{simplicial identities}
\[\begin{cases}
		d_id_j = d_{j-1}d_i   & i<j    \\
		d_i s_j = s_{j-1}d_i  & i<j    \\
		d_j s_j = \id = d_{j+1}s_j     \\
		d_i s_j = s_j d_{i-1} & i>j+1  \\
		s_i s_j = s_{j+1}s_i  & i\le j
	\end{cases}\]
In the following discussion, we will freely employ, and without further mention, this identification. The elements of $X_n$ are called \emph{$n$\hyp{}simplices} of $X$.

The Yoneda embedding $\yon_{\bDelta}\colon \bDelta \to \sSet$ sends every object $\bDelta$ in the \emph{representable} simplicial set $\yon[n] = \Delta(\firstblank,[n])$, acting on objects and morphisms of $\bDelta$ in the expectable way. The usual notation for the representable $\yon[n]$ on $[n]$ is $\Delta[n]$.
\begin{itemize}
	\item Show that the $m$\hyp{}simplices of $\Delta[n]$ are in bijection with the tuples $(a_0,\dots, a_m)$ of elements of $\{0,\dots,n\}$, such that $a_0\le a_1 \le\dots \le a_m$;
	\item Given a simplicial set $X$ and one of its $n$\hyp{}simplices $x\in X_n$ is called \emph{non degenerate} if it can't be expressed as degeneracy of some lower dimensional simplex $y\in X_k$. Show that the non degenerate $m$\hyp{}simplices of $\Delta[n]$ are in bijection with the subsets of $\{0,\dots, n\}$ having cardinality $m+1$.
\end{itemize}
\begin{definition}[Boundaries and horns]\index{Simplicial!--- horns and boundaries}\label{cornibordi}
	We define the \emph{boundary} of $\Delta[n]$ as the union
	\[
		\partial\Delta[n] = \bigcup_{i=0}^n d_{n,i}(\Delta[n-1]),
	\]
	and we define the \emph{$k^\text{th}$ $n$\hyp{}dimensional horn} as the union
	\[
		\Lambda^k[n] = \bigcup_{i\neq k} d_{n,i}(\Delta[n-1]).
	\]
\end{definition}
\section{Miscellaneous definitions}
The already hard endeavour to write a self\hyp{}contained introduction to category theory is made harder by the vastness of the subject. As a consequence, we must necessarily leave out important fragments of theory that often constitute research areas in their own right; these sub\hyp{}theories are here only touched in one or two lines of the book, in so tiny space that it is impossible to do them justice.

While referring the interested reader to the customary sources, we employ this last section of the book to give a rapid glance to the category theory we left out of it.

\index{Ordinal number}\index{Cardinal number}
For us, an \emph{ordinal number} will be any well\hyp{}ordered set, and a \emph{cardinal number} is any ordinal which is not in bijection with a smaller ordinal. Every set $X$ has a unique \emph{cardinality}, \ie a cardinal $\kappa$ with a bijection $\kappa \cong X$ such that there are no bijections from a smaller ordinal. We freely employ results that depend on the axiom of choice when needed. A cardinal $\kappa$ is \emph{regular} if no set of cardinality $\kappa$ is the union of fewer than $\kappa$ sets of cardinality less than $\kappa$; all cardinals in the following subsection are assumed regular without further mention.
\begin{definition}[Filtered category]\label{def:filtcat}
	\index{Category!filtered ---}\index{Filtered category}
	Let $\kappa$ be a cardinal; we say that a category $\A$ is $\kappa$\hyp{}\emph{filtered} if for every category $\J\in\Cat_{<\kappa}$ with less than $\kappa$ objects, $\A$ is injective with respect to the cone completion $\J\to \J^\rhd$; this means that every diagram
	\[
		\vcenter{\xymatrix{
				\J\ar[d]\ar[r]^D & \A \\
				\J^\rhd\ar@{.>}[ur]_{\bar D}
			}}
	\]
	has a dotted filler $\bar D : \J^\rhd \to \A$.
\end{definition}
We say that a category $\C$ admits filtered colimits if for every filtered category $\A$ and every diagram $D : \A \to \C$, the colimit $\colim D$ exists as an object of $\C$. Of course, whenever an ordinal $\alpha$ is regarded as a category, it is a filtered category, so a category that admits all $\kappa$\hyp{}filtered colimits admits all colimits of chains
\[
	C_0 \to C_1 \to \cdots \to C_\alpha \to\cdots
\]
with less than $\kappa$ terms. A useful, completely elementary result is that the existence of colimits over all ordinals less than $\kappa$ implies the existence of $\kappa$\hyp{}filtered colimits; this relies on the fact that every filtered category $\A$ admits a cofinal functor from an ordinal $\alpha_\A$.
\begin{definition}\label{accepre}\index{Category!accessible ---}\index{Accessible category}\index{Category!presentable ---}\index{Presentable category}
	Let $\C$ be a category;
	\begin{itemize}
		\item We say that $\C$ is \emph{$\kappa$\hyp{}accessible} if it admits $\kappa$\hyp{}filtered colimits, and if it has a \emph{small} subcategory $\catS\subset \A$ of $\kappa$\hyp{}presentable objects such that every $A\in\A$ is a $\kappa$\hyp{}filtered colimit of objects in $\catS$.
		\item We say that $\C$ is \emph{(locally) $\kappa$\hyp{}presentable} if it is accessible and cocomplete.
	\end{itemize}
	There is a 2-category whose objects are $\kappa$-accessible categories, whose morphisms are functors that preserve $\kappa$-filtered colimits (also called $\kappa$-ary functors, or functors of rank $\kappa$), and all natural transformations between these. In particular we call \emph{finitary} the functors of rank $\omega$.

	The theory of presentable and accessible categories is a cornerstone of \emph{categorical logic}, \ie of the translation of model theory into the language of category theory.

	Accessible and presentable categories admit \emph{representation theorems}:
	\begin{itemize}
		\item A category $\C$ is accessible if and only if it is equivalent to the ind\hyp{}completion $\text{Ind}_\kappa(\catS)$ of a small category, \ie to the completion of a small category $\catS$ under  $\kappa$\hyp{}filtered colimits;
		\item A category $\C$ is presentable if and only if it is a full reflective subcategory of a category of presheaves $i : \C \to \Cat(\catS^\opp,\Set)$, such that the embedding functor $i$ commutes with $\kappa$\hyp{}filtered colimts.
	\end{itemize}
\end{definition}
All categories of usual algebraic structures are (finitely) accessible, and they are locally (finitely) presentable as soon as they are cocomplete; an example of a category which is $\aleph_1$\hyp{}presentable but not $\aleph_0$\hyp{}presentable: the category of metric spaces and short maps.

We now glance at \emph{topos theory}:
\begin{definition}\label{eletop}\index{Elementary topos}\index{Topos}\index{Category!elementary topos}\index{Category!cartesian closed ---}\index{Cartesian closed category}
	An \emph{elementary topos} is a category $\E$
	\begin{itemize}
		\item which is \emph{cartesian closed}, \ie each functor $\firstblank\times A$ has a right adjoint $[A, \firstblank]$;
		\item having a \emph{subobject classifier}, \ie an object $\Omega\in\E$ such that the functor $\text{Sub} : \E^\opp\to \Set$ sending $A$ into the set of isomorphism classes of monomorphisms $\var{U}{A}$ is representable by the object $\Omega$.
		      \index{Subobject classifier}
	\end{itemize}
	The natural bijection $\E(A,\Omega)\cong\text{Sub}(A)$ is obtained pulling back the monomorphism $U\subseteq A$ along a \emph{universal arrow} $t : 1\to \Omega$, as in the diagram
	\[
		\vcenter{\xymatrix{
				U \pb\ar[r]\ar[d]& 1\ar[d]^t \\
				A \ar[r]_{\chi_U}& \Omega
			}}
	\]
	so, the bijection is induced by the map $\var{U}{A}\mapsto \chi_U$.
\end{definition}
\begin{definition}\label{grotop}\index{Topos!Grothendieck ---}
	A \emph{Grothendieck topos} is an elementary topos that, in addition, is locally finitely presentable.
\end{definition}
Whenever we spoke about sheaves on a topological space or a Grothendieck site, we wer secretly talking about topos theory; the notion of Grothendieck topos is intimately connected with co\fshyp{}end calculus, as we have seen all along chapter 3, and especially in \ref{giraudo}.

In fact, Giraud theorem gives a proof for the difficult implication of the following \emph{recognition principle} for Grothendieck toposes:
\begin{theorem}
	Let $\E$ be a category; then $\E$ is a Grothendieck topos if and only if it is a left exact reflection of a category $\Cat(\A^\opp,\Set)$ of presheaves on a small category $\A$.
\end{theorem}
(recall that a \emph{left exact reflection} of $\C$ is a reflective subcategory $\cate{R}\hookrightarrow \C$ such that the reflector $r : \C \to \cate{R}$ preserves finite limits. It is a reasonably easy exercise to prove that a left exact reflection of a Grothendieck topos is again a Grothendieck topos; Giraud proved that all Grothendieck toposes arise this way.)

Next, we mention the existence of \emph{abelian categories}. \index{Abelian category}\index{Category!abelian ---}

Albeit tangential to the co\fshyp{}end calculus exposed in this book, the notion of abelian category is historically relevant: the definition of co\fshyp{}integration was given by Yoneda \cite{Yoneda} in the setting of module categories, and these constitute the main motivating example for the abstract definition of abelian category.
\begin{definition}\label{Abeliume}
	A category $\A$ is called \emph{abelian} if it satisfies the following list of axioms:
	\begin{enumtag}{a}
		\item \label{a:uno} it is enriched \ref{enrichcat} over abelian groups, \ie every $\A(A,A')$ is an abelian group and the composition operation is $\mathbb Z$\hyp{}bilinear, thus can be represented as an abelian group homomorphism $\A(A,A')\otimes \A(A',A'') \to \A(A,A'')$;
		\item \label{a:due} $\A$ has all finite limits and all finite colimits;
		\item \label{a:ter} in $\A$, every monomorphism is a kernel, and every epimorphism is a cokernel (this means that if $m : A\to B$ is a monomorphism, then it appears in an equaliser diagram $A \xto{m} B \underset{0}{\overset{f}\rightrightarrows} C$, and dually, if $e : B\to C$ is an epimorphism, then it appears in a coequaliser diagram $A \underset{0}{\overset{g}\rightrightarrows} B \xto{e} C$).
	\end{enumtag}
\end{definition}
\begin{remark}
	The notion of abelian category admits many equivalent definitions; so, \ref{Abeliume} above is not the only possible way to define it. An alternative presentation of the axioms disassembles \ref{a:due} into
	\begin{enumtag}{a2}
		\item \label{a2:uno} $\A$ has finite products and finite coproducts;
		\item \label{a2:due} Every morphism $f :A\to B$ in $\A$ has a kernel and a cokernel.
	\end{enumtag}
	A category satisfying only \ref{a:uno} is called \emph{preadditive} (note that this entails that $\A$ has a zero object if and only if it has an initial object, if and only if it has a terminal object); a category satisfying \ref{a:uno} and \ref{a2:uno} is called \emph{additive}; a category satisfying \ref{a:uno}, \ref{a2:uno} and \ref{a2:due} is called \emph{preabelian}.
	\index{Preadditive category}\index{Category!Preadditive ---}
	\index{Additive category}\index{Category!additive ---}
	\index{Preabelian category}\index{Category!preabelian ---}

	A merit of \ref{Abeliume} above is that all axioms are visibly auto\hyp{}dual, \ie $\A$ satisfies \ref{a:uno}--\ref{a:ter} if and only if the opposite category $\A^\opp$ satisfies \ref{a:uno}--\ref{a:ter}.
\end{remark}

Along this book, there are few explicit mentions of abelian categories: see for example \ref{doldekanne}, or the discussion in \ref{dg_stuff}; derived categories appearing in \ref{dg_stuff} are seldom abelian, because they lack finite co\fshyp{}limits.

\begin{exercises}
\item \label{inchia} Show that for every fixed cardinal number $\kappa$, there exists an abelian group of cardinality $\kappa$. If $\kappa$ is infinite, is it true that there is a field of cardinality $\kappa$ (of course, it's false if $\kappa$ is finite and not a power $p^n$ of a prime number $p$)?
\item Let $\C$ be a category admitting an initial and a terminal object; prove that the initial object is the colimit of the unique empty diagram $\boldsymbol{\varnothing} \to \C$ from the empty category, using \ref{colimlim}, and dually that a terminal object is the limit of the same empty diagram.
\item Prove that a category $\C$ has a terminal object if and only if the unique functor $\C \to *$ to the terminal category has a right adjoint; dually, a category $\C$ has an initial object if and only if the unique functor $\C \to *$ to the terminal category has a left adjoint.
\item \label{set-aint-autodual} Prove that the category $\Set^\opp$, \ie the opposite category of sets and functions, is not equivalent to $\Set$.
\index{Comma category}\index{Category!comma ---}
\item \label{da_comma} Consider the comma category of \ref{def:comma};\index{Category!slice ---}
\begin{itemize}
	\item Let $\id_\C$ be the identity functor: characterise the category $(\id_{\C}/G)$; let $C$ be the functor $\J \to \C$ assuming the value $C$ constantly: characterise the category $(C/\id_\C)$.
	\item What is the relation between $(F/G)$ and $(G/F)$?
	\item The \emph{iso\hyp{}comma} category of the functors $F$ and $G$ is given by a similar definition of \ref{def:comma}; the only difference is that we only consider arrows $Fs\to Gt$ that are invertible; we denote it as $[F/G]$. What is the relation, now, between $[F/G]$ and $[G/F]$? Are they isomorphic?
	\item Show that there are two functors $\catS \xot{p_S} (F/G) \xto{p_T} \cT$ such that the square
	      \[\notag
		      \vcenter{\xymatrix{
				      (F/G)\ar[r]^{p_S}\ar[d]_{p_T} & \catS \ar[d]^F\\
				      \cT \ar[r]_G & \C
			      }}
	      \]
	      commutes.
	\item Consider the functor $J_1 : (F/G) \xto{p_S} \catS \xto{F} \C$. Determine the category $(J_1/G)$. Does the composition of the projections above in the triangle $\star$ equals the projection $(J_1/G)\to \cT$ (in other words, does the triangle commute)?
	      \[\notag
		      \vcenter{\xymatrix{
		      (J_1/G) \ar[r]^{q_{(F/G)}}\ar[dr]_{q_T}& (F/G) \ar@{}[dl]|(.3)\star \ar[d]_{p_T}\ar[r]^{p_S} & \catS \ar[d]^F \\
		      & \cT \ar[r]_G & \C
		      }}
	      \]
	\item Does the category $(F/G)$ have a universal property?
\end{itemize}
\item \index{Sierpi\'nski space} Let {\sierpo} be the Sierpiński space, where a two\hyp{}point set $J=\{a,b\}$ has topology $\{\varnothing, \{a\}, J\}$. Let \codis be the codiscrete space, where $\{a,b\}$ has trivial topology $\{\varnothing, J\}$, and \dis the discrete space where $J$ has discrete topology $\{\varnothing, \{a\}, \{b\}, J\}$. Show that
\begin{itemize}
	\item the functor $O:\Spc^\opp \to \Set$ that sends a topological space into its set of open subsets is representable, and that {\sierpo} is its representing object;
	\item the functor $D : \Spc^\opp \to \Set$ that sends a topological space $X$ into its set of \emph{disconnections}, \ie the set of pairs $(U,V)$ such that $U\cup V = X$ and $U\cap V=\varnothing$, is representable by the discrete space \dis.
	\item the functor $S : \Spc^\opp\to\Set$ that sends a topological space into its set of subspaces is representable by the codiscrete space \codis.
\end{itemize}
Show that the set of natural transformations $[\Spc^\opp,\Set](O,S) $ has exactly four elements.
\item \label{cones-are-cocones-in-cop} Show that a cocone $\bar D : \J^\rhd \to \C$ for a diagram $D : \J \to \C$ is exactly a cone $\bar D^{\opp} : \J^\lhd \to \C^\opp$ for the opposite functor $D^\opp : \J^\opp \to \C^\opp$. Prove it directly, but then notice that \ref{conecompletion} and in particular \eqref{cone-as-lift} give a slick argument to conclude.
\item Let $f : [n]\to [m]$ be a map in $\bDelta$ (see \ref{deltacat}). We regard $\bDelta$ as a subcategory of $\Cat$ in the obvious way. Show that
\begin{itemize}
	\item \index{Category!simplex ---}\index{_aaa_delta@$\bDelta$} A morphism $f\colon [n]\to [m]$ in $\bDelta$ has a left adjoint if and only if $f(n)=m$; in such a situation the adjoint $f_L$ sends $i$ into $f_L(i)=\min\{j\in[n]\mid f(j)\ge i\}$.
	\item Dually, a morphism $f\colon [n]\to [m]$ in $\bDelta$ has a right adjoint if and only if $f(0)=0$; in such a situation, the adjoint $f_R$ sends $i$ into $f_R(i)=\max\{j\in[n]\mid f(j)\le i\}$.
\end{itemize}
Deduce that for every $j\in[n]$ the function $\sigma_{n,j}$ has both a left and a right adjoint, and in particular
\[\notag
	\delta_{n,j+1} \dashv \sigma_{n-1,j}\dashv \delta_{n,j} \colon
	\vcenter{\xymatrix@C=2cm{
	[n] \ar[r]|{\sigma_{n-1,j}}& \ar@<9pt>[l]\ar@<-9pt>[l] [n-1]
	}}
\]
\item \label{ex:galua} Let $F\hookrightarrow E$ be a homomorphism of rings between fields; it is thus a monomorphism. We define
\[ \notag
	Fix(E|F) = \{\sigma : E \xto{\sim} E \mid \sigma|_F = \id_F\}\\
\]
as the set of automorphisms of $E$ that become the identity map when restricted to $F$, and
\[ \notag	Ext(E|F) = \{ K \mid F \le K\le E\} \]
the set of intermediate extensions between $F$ and $E$. We define two functions $\fki : Fix(E|F) \to Ext(E|F)$ and $\fkj : Ext(E|F) \to Fix(E|F)$ sending respectively $\sigma\in Fix(E|F)$ in the intermediate field $\{a\in E\mid \sigma(a)=a\}$ and the intermediate field $F\le K\le E$ in the group of those $E$\hyp{}automorphisms that become the identity when restricted to $K$. Show that $\fki,\fkj$ set up an adjunction (thus the name \emph{Galois connections} \index{Galois connection} for these adjunctions) when the posets $Fix(E|F)$ and $Ext(E|F)$ are regarded as categories; who is the right adjoint, and who is the left adjoint? Who is the unit, and who is the counit?
\item \label{ex:pizero} Show that the limit of a constant diagram $\Delta_\J A : \J \to \A$ is a product of copies of $A$ indexed by the set $\pi_0\J$ of connected components of the domain $\J$: prove that there is a bijection
\[\notag
	\Cn(X,A)\cong \Set(\pi_0\J, \A(X,A))
\]
and conclude, using \eqref{coten}.
\item \label{ex:orthoequivs} Show that the following conditions are equivalent in a category $\C$ with a terminal object:
\begin{itemize}
	\item The terminal arrow $\var{X}{1}$ is right orthogonal to $f :A\to B$;
	\item The functor $\C(\firstblank,X)$ sends $f$ to an isomorphism;
	\item Every arrow $A\to X$ in $\C$ has a unique extension to $B\to X$ along $f$.
\end{itemize}
\item \label{exe:yon-is-nat} Show that in the same notation of \eqref{yoneda-compo}, $Y_X$ is the $X$\hyp{}component of a natural transformation $\Cat(\C^\opp,\Sets)(\yon,F)\To F$, that is natural in both its arguments; in other words, show that for every arrow $X\to X'$ the square
\[\notag
	\vcenter{\xymatrix{
	\Cat(\C^\opp,\Sets)(\yon X,F) \ar[r]^-Y & FX \\
	\Cat(\C^\opp,\Sets)(\yon(X'),F)\ar[u]^{[\C^\opp,\Set](\yon(f),F)}\ar[r]_-Y & FX' \ar[u]_{Ff}
	}}
\]
is commutative.

Similarly, for every natural transformation $\tau : F \To G$ the square
\[\notag
	\vcenter{\xymatrix{
			\Cat(\C^\opp,\Sets)(\yon X,F) \ar[d]_{\Cat(\C^\opp,\Sets)(\yon X,\tau)} \ar[r]^-Y & FX\ar[d]^{\tau_X} \\
			\Cat(\C^\opp,\Sets)(\yon X,G)\ar[r]_-Y & GX
		}}
\]
is commutative.
\item \index{Category!Kleisli ---}\label{kleisli-is-kl} Denote $\C_T$ the category defined in \ref{frealg}. Define the category $\Kl(T)$ as follows: it has the same objects of $\C$, and $\Kl(T)(X,Y) := \C(X, TY)$. Composition of $ X \xto{f} TY$ with $Y\xto{g} TZ$ is defined by the rule
\[\notag
	g \bullet_T f := \big( X \xto{f} TY \xto{Tg} TTZ \xto{\mu_Z} TZ \big)
\]
\begin{itemize}
	\item Prove that this really defines a category, if the identity map of an object $A$ is the unit $\eta_A : A \to TA$ (prove that the composition is associative and that $\eta_B\bullet_T f = f = f \bullet_T \eta_A$ for every $f : A\to TB$).
	\item Defines a functor $ W : \Kl(T) \to \C_T$ that acts as $T$ on objects, and such that $Wf = \mu_B\circ Tf$. Prove that $W$, so defined, is fully faithful, and that its essential image is made by free $T$\hyp{}algebras.
\end{itemize}
\item Show that each square
\[\notag
	\vcenter{\xymatrix{
	[n+1] \ar@{}[dr]|{i\le j}\ar[r]^{\sigma_i}\ar[d]_{\sigma_i}& [n]\ar[d]^{\sigma_j}\\
	[n] \ar[r]_{\sigma_{j+1}}& [n-1]
	}}\]
is an absolute pushout, and that every square
\[\notag
	\vcenter{\xymatrix{
	[n-1] \ar@{}[dr]|{i<j}\ar[r]^{\delta_i}\ar[d]_{\delta_{j-1}} & [n] \ar[d]_{\delta_j}\\
	[n] \ar[r]_{\delta_i} & [n+1]
	}}\]
is an absolute pullback.
\item \label{ex:commies} A comonoid in a monoidal category is an object $M \in\C$ endowed with maps $c : M\to M\otimes M$ and $e : M\to I$ such that the following diagrams commute
\[\notag
	{\scriptstyle\vcenter{\xymatrix{
	M\otimes M\otimes M \ar@{<-}[r]^{M\otimes c}\ar@{<-}[d]_{c\otimes M}& M\otimes M \ar@{<-}[d]^c & M\otimes I \ar@{<-}@{=}[dr]\ar@{<-}[r]^{M\otimes e}& M\otimes M \ar@{<-}[d]_m & I\otimes M\ar@{<-}[l]_{e\otimes M}\ar@{<-}@{=}[dl]\\
	M\otimes M \ar@{<-}[r]_c & M & & M
	}}}
\]
witnessing \emph{coassociativity} and \emph{counitality} properties for $(c,e)$. 	A \emph{comonad} on a category $\C$ is now a functor $S : \C \to \C$ which is a comonoid with respect to the monoidal structure $([\C,\C],\circ)$; a
\begin{itemize}
	\item Write down the axioms of \emph{coassociativity} ad counitality for $S$;
	\item Show that if $F\adjunct{\epsilon}{\eta} G$ is an adjunction then $S=FG$ is a comonad with comultiplication $F *\eta * G$ and counit $\epsilon$;
	\item Define a category of \emph{coalgebras} for $S$, and of \emph{free coalgebras} for $S$; show that they enjoy suitable universal properties.
\end{itemize}
\item Let
\[\notag
	\xymatrix{
		\C \ar[r]|F & \D \ar@<1em>[l]^L \ar@<-1em>[l]_R
	}
\] be a triple of adjoint functors, meaning that $L\dashv F$ and $F\dashv R$; show that $L$ is a fully faithful functor if and only if $R$ is a fully faithful functor.
\item \index{Monoidal category}\cite{segal1974a} Given a monoidal category $(\V, \otimes)$, show that there exists a category $\V^\otimes$ defined as follows:
\begin{itemize}
	\item the objects of $\V^\otimes$ are $n$\hyp{}tuples of objects in $\V$, denoted $[C_1,\dots,C_n]$ (this follows the convention that if $n=0$, the tuple is empty);
	\item morphisms $[C_1,\dots, C_n]\to [D_1,\dots,D_m]$ are defined as pairs $(\alpha,\{f_j\})$ where $\alpha$ is a partial function $[n]\to [m]$ having domain $S_\alpha$ and
	      \[\textstyle\notag
		      \big\{f_j : \bigotimes_{\{i\mid \alpha(i)=j\}} C_i\to D_j \mid 1\le j\le m\big\}
	      \]
	      is a family of morphisms in $\V$.
	\item Define the composition law of two morphisms $(\alpha,f), (\beta,g)$.
	      \index{_aaa_Fin@$\Fin$}
	\item \label{monoidaux} Let $\Fin_*$ be the category $*/\Set_{<\omega}$ of pointed finite sets. Show that there is a functor $p : \V^\otimes\to \Fin_*$ sending $[C_1,\dots,C_n]$ into $[n]_*$; show that $p$ is an \emph{opfibration}: for every object $\overline C= [C_1,\dots,C_n]\in\V^\otimes$ and every arrow $f : p(\overline C) \to [m]_*$ in $\Fin_*$ there is an arrow $(\theta_f,\bar f) : \overline C \to \overline D = [D_1,\dots, D_m]$ such that $p(\theta_f,\bar f)=f$, and such that the composition with $\bar f$ induces a bijection for every $\overline E = [E_1,\dots,E_d]$, as follows:
	      \[
		      \V^\otimes(\overline D, \overline E) \cong \V^\otimes(\overline C, \overline E) \times_{\Fin_*([n]_*, [d]_*)}\Fin_*([m]_*, [d]_*).
	      \]
	\item Show that if we denote $\V^\otimes_n$ the fiber of $[n]_*$ along $p$, then the functor $p$ induces a functor $\V^\otimes_m \to \V^\otimes_n$ among the fibers, for every $f : [m]_* \to [n]_*$ in $\Fin_*$.
	\item Show that $\V^\otimes_0 \cong \{0\}$, $\V^\otimes_1 \cong \V$, and that more in general $\V^\otimes_n \cong \V\times\dots\times\V$ ($n$ times).
	\item Show that the correspondence sending $\V$ into $\V^\otimes$ is functorial in $\V$; are strong monoidal functors enough?
	\item Show that $\Fin_* \cong \{0\}^\otimes$ with respect to the unique monoidal structure that exists on the terminal category $\{0\}$; show that $p$ is the functor induced by the unique functor $\V \to \{0\}$.
\end{itemize}
\end{exercises}

\chapter{}

\begin{adjustbox}{width=\textwidth}
	\begingroup
	\renewcommand{\arraystretch}{3.2}
	\newcommand{\extra}{.8em}
	\begin{tabular}{l@{\hspace{-3em}}c}
		\multicolumn{2}{c}{\bf Table of notable integrals}                                                                                                                    \\\toprule
		\multirow{2}{*}{Fubini rule}
		 & \hspace{-3em}$\displaystyle \int^C\int^D T(C,D,C,D)\cong \int^D\int^C T(C,D,C,D)\cong \int^{(C,D)} T(C,C,D,D)$                                                     \\
		 & \hspace{-3em}$\displaystyle \forall n\ge 2,\sigma \in \text{Sym}(n),\quad \int^{C_{\sigma 1}}\kern-1em\dots\int^{C_{\sigma n}} T \cong \int^{(C_1,\dots, C_n)} T $ \\[\extra]\midrule
		Natural transformations
		 & $\displaystyle \Cat(\C,\D)(F,G)\cong \int_C\D(FC,GC)$                                                                                                              \\[\extra]\midrule
		Tensor product of functors
		 & $\displaystyle G\boxtimes F = \int^C GC\otimes FC$                                                                                                                 \\[\extra]\midrule
		\multirow{2}{*}{Yoneda lemma}
		 & $\displaystyle FX \cong \int^A [X,A]\times FA \qquad FX \cong \int_A \Set([A,X], FA)$                                                                              \\
		 & $\displaystyle GX \cong \int^A [A,X]\times GA \qquad GX \cong \int_A \Set([X,A], GA)$                                                                              \\[\extra]\midrule
	\end{tabular}
	\endgroup
\end{adjustbox}

\begin{adjustbox}{width=\textwidth}
	\begingroup
	\renewcommand{\arraystretch}{3.2}
	\newcommand{\extra}{.8em}
	\begin{tabular}{l@{\hspace{-2em}}c}
		\multicolumn{2}{c}{\bf Table of notable integrals}                                                                                           \\\toprule
		\multirow{2}{*}{Kan extensions: $\displaystyle \C \xot{G} \A \xto{F} \B$}
		                 & $\displaystyle \Lan_GF(C) \cong \int^A \C(GA, C)\otimes FA$                                                               \\
		                 & $\displaystyle \Ran_GF(A) \cong \int_A \C(C, GA)\pitchfork FA $                                                           \\[\extra]\midrule
		density comonad: & $\displaystyle T_F \cong \int_A \C(\firstblank,FA)\pitchfork FA$                                                          \\[\extra]
		codensity monad: & $\displaystyle S^F \cong \int^A \C(FA,\firstblank)\otimes FA$                                                             \\[\extra]\midrule
		\multirow{2}{*}{Weighted co\fshyp{}limits}
		                 & $\displaystyle \wcolim{W} F \cong \int^A WA \otimes FA$                                                                   \\
		                 & $\displaystyle \wlim{W} F \cong \int_A WA\pitchfork FA$                                                                   \\[\extra]\bottomrule
		\multirow{2}{*}{Profunctor theory}
		                 & $\displaystyle \proP \bullet \proQ = \int^X \proP(\firstblank, X)\otimes \proQ(X,\secondblank)$                           \\
		                 & $\displaystyle \Ran_\proP\proQ = \int_A\hom(\proP_A,\proQ_A) \qquad \Rift_{\proP}{\proQ} = \int_A \hom(\proH_A, \proP_A)$ \\[\extra]\bottomrule
	\end{tabular}
	\endgroup
\end{adjustbox}

\begin{adjustbox}{width=\textwidth}
	\begingroup
	\renewcommand{\arraystretch}{3.2}
	\newcommand{\extra}{.8em}
	\begin{tabular}{l@{\hspace{-2em}}c}
		\multicolumn{2}{c}{\bf Table of notable integrals}                                                                                                              \\\toprule
		\multirow{3}{*}{Operads}
		                & $\displaystyle F\ast G := \int^{XY}\C(X\oplus Y,\firstblank)\otimes FX\otimes GY$                                                             \\
		                & $\displaystyle F\odot G:=\int^m Fm\otimes G^{\ast m}_{\qquad\substack{G^{\ast m} = G \ast\cdots\ast G                                         \\ m \text{ times}}}$ \\
		                & $\displaystyle \frac{F\diamond G \to H}{F \to \{G,H\}} \qquad \{G,H\}m=\int_k[G^{\ast m}k,Hk]$                                                \\[\extra]\midrule
		\multirow{3}{*}{Simplicial coends}
		                & $\displaystyle \oint_A T(A,A) := \int_{A',A''} \delta \A(A', A'')\pitchfork T(A',A'')$                                                        \\
		                & $\displaystyle \oint^A T(A,A) := \int^{A',A''}\delta \A(A',A'') \otimes T(A',A'')$                                                            \\
		                & $\displaystyle \delta\A(A',A'') = \int^{n\in\bDelta} \coprod_{X_0,\dots, X_n \in \A} \Delta[n]\times \A(A,X_0)\times \cdots \times \A(X_n,B)$ \\[\extra]\midrule
		\multirow{3}{*}{Promonoidal structures}
		                & $\displaystyle [F\ast_{\mathfrak P} G]C = \int^{AB} P(A,B;C)\times FA\times GB$                                                               \\
		                & $\displaystyle \int^Z P_D^{AY}P_Y^{BC} \cong \int^Z P_X^{AB} P_D^{XC}$                                                                        \\
		                & $\displaystyle \int^{YZ} J_Z H^A_Y P^{YZ}_B \int^z J_Z\big( \int^Y H^A_Y P^{YZ}_B \big) \cong \int^Z J_Z P^{AZ}_B \cong \hom(A,B)$            \\[\extra]\midrule
		Tambara module: & $\displaystyle \phi_{P}(X, Y)
		= \int^{C, U, V} \C(X, C \otimes U) \times \C(C \otimes V, Y) \times P(U, V)$                                                                                   \\[\extra]\bottomrule
	\end{tabular}
	\endgroup
\end{adjustbox}

\backmatter

\addtocontents{toc}{\vspace{\baselineskip}}
\renewcommand{\refname}{Bibliography}
\bibliography{allofthem}{}
\label{refs}
\bibliographystyle{amsalpha}
\cleardoublepage

\printindex
\end{document}